\pdfoutput=1
\documentclass[letterpaper, oneside, reqno]{amsart}
\usepackage[margin=1in]{geometry}

\usepackage{amsmath,amssymb,mathtools}
\usepackage[foot]{amsaddr}

\usepackage{placeins}

\usepackage{amssymb} 
\usepackage{stmaryrd}
\usepackage{graphicx}
\usepackage[colorlinks=true, pdfstartview=FitV, linkcolor=blue, citecolor=Green, urlcolor=WildStrawberry, linktoc=page]{hyperref} 

\usepackage{array}
\usepackage{ragged2e}

\usepackage{bm}

\usepackage{dsfont}
\usepackage[T1]{fontenc}
\usepackage[utf8]{inputenc}

\DeclareFontFamily{OT1}{rsfs}{}
\DeclareFontShape{OT1}{rsfs}{n}{it}{<-> rsfs10}{}
\DeclareMathAlphabet{\mathscr}{OT1}{rsfs}{n}{it}
\usepackage{mathrsfs}
\usepackage{MnSymbol}
\usepackage{microtype} 
\usepackage{enumitem}
\usepackage[dvipsnames]{xcolor}

\usepackage[normalem]{ulem}

\usepackage{booktabs} 

\usepackage{comment}

\usepackage{braket}

\usepackage{etoolbox}
\newtoggle{focs}
\toggletrue{focs}
\newcommand{\iffocs}[2]{\iftoggle{focs}{#1}{#2}}

\definecolor{darkgreen}{rgb}{0,0.5,0}
\definecolor{darkblue}{rgb}{0,0,0.7}
\definecolor{darkred}{rgb}{0.9,0.1,0.1}

\newcommand{\js}[1]{{\color{ForestGreen}{[JS: #1]}}}

\newlength{\bibitemsep}
\let\oldthebibliography\thebibliography
\renewcommand\thebibliography[1]{%
  \oldthebibliography{#1}%
  \setlength{\parskip}{\bibitemsep}%
  \setlength{\itemsep}{-7pt}%
}

\newtheoremstyle{break}%
{}{}%
{\itshape}{}%
{\bfseries}{.\vphantom{$p_{p_{p_p}}$}}%
{\newline}
{\thmname{#1}\thmnumber{ #2}\thmnote{\ \,\textmd{(#3)}}}

\theoremstyle{break}

\newtheorem{proposition}{Proposition}
\newtheorem{theorem}[proposition]{Theorem}
\newtheorem{lemma}[proposition]{Lemma}
\newtheorem{corollary}[proposition]{Corollary}

\theoremstyle{remark}
\newtheorem{remark}[proposition]{Remark}

\theoremstyle{definition}
\newtheorem{definition}[proposition]{Definition}

\newtheorem{fact}[proposition]{Fact}

\newtheorem*{lemma*}{Lemma}

\newcommand{\vocab}[1]{\emph{#1}}

\numberwithin{equation}{section}
\numberwithin{proposition}{section}
\numberwithin{figure}{section}
\numberwithin{table}{section}

\newcommand{\Z}{\mathbb{Z}}
\newcommand{\N}{\mathbb{N}}

\newcommand{\R}{\mathbb{R}}

\newcommand{\calN}{\mathcal N}

\renewcommand{\P}{\mathop{{}\mathbb{P}}}
\renewcommand{\Pr}{\P}
\newcommand{\Cov}{\mathop{{}\boldsymbol{\mathrm{Cov}}}}
\newcommand{\E}{\mathop{{}\mathbb{E}}}

\newcommand{\hQ}{\widehat{Q}}
\renewcommand{\hm}{\widehat{m}}
\newcommand{\hx}{\widehat{x}}
\newcommand{\hy}{\widehat{y}}
\newcommand{\hw}{\widehat{w}}

\newcommand{\hp}{\widehat{p}}
\newcommand{\hf}{\widehat{f}}
\newcommand{\tx}{\widetilde{x}}
\newcommand{\ty}{\widetilde{y}}
\newcommand{\tw}{\widetilde{w}}

\renewcommand{\le}{\leqslant}
\renewcommand{\ge}{\geqslant}
\renewcommand{\leq}{\leqslant}

\renewcommand{\subset}{\subseteq}
\renewcommand{\bar}{\overline}

\renewcommand{\tilde}{\widetilde}
\newcommand{\td}{\widetilde}

\renewcommand{\hat}{\widehat}

\newcommand{\subeq}{\subseteq}

\newcommand{\al}{\alpha}
\newcommand{\be}{\beta}
\newcommand{\ga}{\gamma}
\newcommand{\de}{\delta}
\newcommand{\ka}{\kappa}
\newcommand{\lm}{\lambda}
\newcommand{\La}{\Lambda}
\newcommand{\ph}{\varphi}
\newcommand{\De}{\Delta}
\newcommand{\ep}{\varepsilon}
\newcommand{\eps}{\varepsilon}
\newcommand{\si}{\sigma}
\newcommand{\Si}{\Sigma}
\newcommand{\om}{\omega}
\newcommand{\Om}{\Omega}
\newcommand{\te}{\theta}
\newcommand{\Te}{\Theta}
\newcommand{\rh}{\rho}

\newcommand{\rec}{\rho_{\mathrm{EC}}}

\newcommand{\mg}[0]{m}
\newcommand{\ub}[2]{\underbrace{#1}_{#2}}

\newcommand{\pl}{\partial}
\newcommand{\dd}[2]{\frac{d #1}{d #2}}

\newcommand{\ddd}[1]{\frac{d}{d #1}}

\newcommand{\pat}[1]{(\textup{#1})}

\DeclareMathOperator{\dist}{\mathcal{L}}
\DeclareMathOperator{\diam}{diam}
\DeclareMathOperator{\tr}{tr}
\DeclareMathOperator{\supp}{supp}
\DeclareMathOperator{\sign}{sign}

\newcommand{\Tr}{\mathsf{Tr}}

\newcommand{\KL}{\operatorname{KL}}
\newcommand{\TV}{\operatorname{TV}}
\newcommand{\GOE}{\operatorname{GOE}}

\newenvironment{e*}{\begin{equation*}}{\end{equation*}\ignorespacesafterend}
\newcommand{\norm}[1]{\left\lVert{#1}\right\rVert}

\renewcommand{\Tr}{\mathsf{Tr}}

\newcommand{\TAP}{{\mathrm{TAP}}}
\newcommand{\FT}{\calF_{\,\mathrm{TAP}}}

\newcommand{\pol}[0]{^{\mathrm{pol}}}
\newcommand{\polar}[0]{\operatorname{polar}}
\newcommand{\homg}[0]{^{\mathrm{hom}}}

\newcommand{\ip}[2]{\langle#1, #2\rangle}
\newcommand{\iprod}[1]{\langle#1\rangle}
\newcommand{\Iprod}[1]{\left\langle#1\right\rangle}

\newcommand{\an}[1]{\left\langle#1\right\rangle}
\newcommand{\fl}[1]{\left\lfloor#1\right\rfloor}
\newcommand{\ce}[1]{\left\lceil#1\right\rceil}

\newcommand{\sT}{\mathsf{T}}

\newcommand{\Lip}{\mathsf{Lip}}
\newcommand{\ot}{\otimes}

\newcommand{\ab}[1]{\left| {#1} \right|}
\newcommand{\ba}[1]{\left[ {#1} \right]}
\newcommand{\bc}[1]{\left\{ {#1} \right\}}
\newcommand{\pa}[1]{\left( {#1} \right)}
\newcommand{\ve}[1]{\left\Vert {#1}\right\Vert}

\newcommand{\fc}[2]{\frac{#1}{#2}}
\newcommand{\rc}[1]{\frac{1}{#1}}
\newcommand{\sfc}[2]{\sqrt{\frac{#1}{#2}}}
\newcommand{\pf}[2]{\pa{\frac{#1}{#2}}}
\newcommand{\prc}[1]{\pa{\frac{1}{#1}}}

\newcommand{\sumo}[2]{\sum_{#1=1}^{#2}}
\newcommand{\sumz}[2]{\sum_{#1=0}^{#2}}
\newcommand{\prodo}[2]{\prod_{#1=1}^{#2}}

\newcommand{\gd}[0]{\nabla}
\newcommand{\bs}[0]{\setminus}

\DeclareMathOperator{\sech}{sech}

\newcommand{\Id}[0]{I}
\DeclareMathOperator{\id}{id}

\DeclareMathOperator{\im}{Im}

\DeclareMathOperator{\Spec}{Spec}

\DeclareMathOperator{\emp}{emp}
\DeclareMathOperator{\osc}{osc}

\renewcommand{\tr}{\operatorname{tr}}
\DeclareMathOperator{\diag}{\mathsf{diag}}

\renewcommand{\norm}[1]{\left\lVert{#1}\right\rVert}

\newcommand{\sop}{_\mathsf{op}}
\newcommand{\opnorm}[1]{\ensuremath{\left\lVert #1 \right\rVert\sop}}
\newcommand{\ED}{E_{\mathcal{D}_n}}
\newcommand{\lpnorm}[2][2]{\ensuremath{\left\lVert {#2} \right\rVert_{#1}}}
\newcommand{\schnorm}[2][2]{\ensuremath{\left\lVert {#2} \right\rVert_{L^{#1}}}}
\newcommand{\Lpnorm}[2][2]{\ensuremath{\left\lVert {#2} \right\rVert_{L^{#1}}}}
\newcommand{\lipnorm}[1]{\ensuremath{\left\lVert {#1} \right\rVert_{\Lip}}}
\newcommand{\frenorm}[2][]{\ensuremath{\left\lVert {#2} \right\rVert_{\mathrm{Fr\acute{e}}{#1}}}}

\newcommand{\poly}{\mathsf{poly}}

\newcommand{\Var}{\mathop{{}\boldsymbol{\mathrm{Var}}}}
\newcommand{\Ent}{\mathop{{}\boldsymbol{\mathrm{Ent}}}}

\newcommand{\coltwo}[2]{
\begin{pmatrix}
{#1}\\
{#2}
\end{pmatrix}}

\renewcommand{\le}{\leqslant}
\renewcommand{\leq}{\leqslant}
\renewcommand{\ge}{\geqslant}

\usepackage{prettyref}
\newcommand{\savehyperref}[2]{\texorpdfstring{\hyperref[#1]{#2}}{#2}}

\protected\def\verythinspace{%
  \ifmmode
    \mskip0.5\thinmuskip
  \else
    \ifhmode
      \kern0.083em
    \fi
  \fi
}
\newrefformat{eq}{\savehyperref{#1}{\textup{(\ref*{#1})}}}
\newrefformat{e}{\savehyperref{#1}{\textup{(\ref*{#1})}}}
\newrefformat{ineq}{\savehyperref{#1}{\textup{(\ref*{#1})}}}
\newrefformat{eqn}{\savehyperref{#1}{\textup{(\ref*{#1})}}}
\newrefformat{l}{\savehyperref{#1}{Lemma~\ref*{#1}}}
\newrefformat{lem}{\savehyperref{#1}{Lemma~\ref*{#1}}}
\newrefformat{def}{\savehyperref{#1}{Definition~\ref*{#1}}}
\newrefformat{d}{\savehyperref{#1}{Definition~\ref*{#1}}}
\newrefformat{t}{\savehyperref{#1}{Theorem~\ref*{#1}}}
\newrefformat{thm}{\savehyperref{#1}{Theorem~\ref*{#1}}}
\newrefformat{cor}{\savehyperref{#1}{Corollary~\ref*{#1}}}
\newrefformat{c}{\savehyperref{#1}{Corollary~\ref*{#1}}}
\newrefformat{cha}{\savehyperref{#1}{Chapter~\ref*{#1}}}
\newrefformat{sec}{\savehyperref{#1}{\S\verythinspace\ref*{#1}}}
\newrefformat{s}{\savehyperref{#1}{\S\verythinspace\ref*{#1}}}
\newrefformat{subsec}{\savehyperref{#1}{\S\verythinspace\ref*{#1}}}
\newrefformat{app}{\savehyperref{#1}{\S\verythinspace\ref*{#1}}}
\newrefformat{tab}{\savehyperref{#1}{Table~\ref*{#1}}}
\newrefformat{fig}{\savehyperref{#1}{Figure~\ref*{#1}}}
\newrefformat{hyp}{\savehyperref{#1}{Hypothesis~\ref*{#1}}}
\newrefformat{alg}{\savehyperref{#1}{Algorithm~\ref*{#1}}}
\newrefformat{a}{\savehyperref{#1}{Algorithm~\ref*{#1}}}
\newrefformat{rem}{\savehyperref{#1}{Remark~\ref*{#1}}}
\newrefformat{item}{\savehyperref{#1}{Item~\ref*{#1}}}
\newrefformat{step}{\savehyperref{#1}{step~\ref*{#1}}}
\newrefformat{conj}{\savehyperref{#1}{Conjecture~\ref*{#1}}}
\newrefformat{fact}{\savehyperref{#1}{Fact~\ref*{#1}}}
\newrefformat{p}{\savehyperref{#1}{Proposition~\ref*{#1}}}
\newrefformat{prop}{\savehyperref{#1}{Proposition~\ref*{#1}}}
\newrefformat{prob}{\savehyperref{#1}{Problem~\ref*{#1}}}
\newrefformat{claim}{\savehyperref{#1}{Claim~\ref*{#1}}}
\newrefformat{clm}{\savehyperref{#1}{Claim~\ref*{#1}}}
\newrefformat{relax}{\savehyperref{#1}{Relaxation~\ref*{#1}}}
\newrefformat{rem}{\savehyperref{#1}{Remark~\ref*{#1}}}
\newrefformat{red}{\savehyperref{#1}{Reduction~\ref*{#1}}}
\newrefformat{part}{\savehyperref{#1}{Part~\ref*{#1}}}
\newrefformat{ex}{\savehyperref{#1}{Exercise~\ref*{#1}}}
\newrefformat{property}{\savehyperref{#1}{Property~\ref*{#1}}}
\newrefformat{type}{\savehyperref{#1}{Type~\ref*{#1}}}
\newrefformat{eg}{\savehyperref{#1}{Example~\ref*{#1}}}
\newrefformat{obs}{\savehyperref{#1}{Observation~\ref*{#1}}}
\newrefformat{que}{\savehyperref{#1}{Question~\ref*{#1}}}
\newrefformat{cond}{\savehyperref{#1}{Condition~\ref*{#1}}}
\newrefformat{ass}{\savehyperref{#1}{Assumption~\ref*{#1}}}
\newrefformat{not}{\savehyperref{#1}{Notation~\ref*{#1}}}
\newrefformat{cond}{\savehyperref{#1}{Condition~\ref*{#1}}}
\newcommand{\Sref}[1]{\hyperref[#1]{\S\ref*{#1}}}

\let\pref=\prettyref
\let\Cref=\prettyref

\renewcommand{\eps}{\varepsilon}

\newcommand{\calA}{\mathcal A}
\newcommand{\calB}{\mathcal B}

\newcommand{\calD}{\mathcal D}
\newcommand{\calE}{\mathcal E}
\newcommand{\calF}{\mathcal F}

\newcommand{\calH}{\mathcal H}

\newcommand{\calL}{\mathcal L}
\newcommand{\calM}{\mathcal M}
\renewcommand{\calN}{\mathcal N}
\newcommand{\calP}{\mathcal P}
\newcommand{\calQ}{\mathcal Q}
\newcommand{\calR}{\mathcal R}

\newcommand{\calW}{\mathcal W}

\newcommand{\sL}{\mathscr L}

\newcommand{\sE}{\mathscr E}

\newcommand{\bbD}{\mathbb D}

\newcommand{\bbR}{\mathbb R}

\newcommand{\bbH}{\mathbb H}

\newcommand{\one}[0]{\mathds{1}}
\renewcommand{\set}[2]{\left\{{#1}:{#2}\right\}}

\renewcommand{\R}{\mathbb R}
\newcommand{\C}{\mathbb C}
\renewcommand{\N}{\mathbb N}
\renewcommand{\Z}{\mathbb Z}

\newcommand{\la}{\langle}
\newcommand{\ra}{\rangle}

\newcommand{\iy}{\infty}

\newcommand{\slice}[2]{\binom{[#1]}{#2}}
\newcommand{\wdg}[2]{\binom{[#1]}{\le #2}}
\newcommand{\ball}[2]{B_{#1}(#2)}
\newcommand{\bball}[2]{\partial B_{#1}(#2)}

\newcommand{\tmix}[0]{t_{\mathrm{mix}}}

\newcommand{\CJET}[1]{C_{\mathrm{JE},T,#1}}

\usepackage{needspace}
\newlength{\ppartneed}
\newcommand{\ppart}[1]{%
  \par\Needspace{\ppartneed}\addvspace{\smallskipamount}%
  \noindent\textit{#1.}\hspace{0.5em}\ignorespaces%
}

\usepackage{algorithm}
\usepackage{algpseudocode}
\algnewcommand\algorithmicinput{\textbf{Input: }}
\algnewcommand\INPUT{\State\algorithmicinput}
\algnewcommand\algorithmicinitialize{\textbf{Initialize: }}
\algnewcommand\INIT{\State\algorithmicinitialize}
\algnewcommand\algorithmicrun{\textbf{Run: }}
\algnewcommand\RUN{\State\algorithmicrun}
\algnewcommand\algorithmicupdate{\textbf{Update: }}
\algnewcommand\UPDATE{\State\algorithmicupdate}
\algnewcommand\algorithmicset{\textbf{Set: }}
\algnewcommand\SET{\State\algorithmicset}
\algnewcommand\algorithmicquery{\textbf{Query: }}
\algnewcommand\QUERY{\State\algorithmicquery}
\algnewcommand\algorithmicoutput{\textbf{Output: }}
\algnewcommand\OUTPUT{\State\algorithmicoutput}

\makeatletter
\newcommand\appendix@section[1]{%
  \refstepcounter{section}%
  \orig@section*{\@Alph\c@section.\texorpdfstring{\,\,\,\;}{}#1}
}
\let\orig@section\section
\g@addto@macro\appendix{\let\section\appendix@section}
\makeatother

\renewcommand{\paragraph}[1]{\medskip\noindent{\bf #1{.}}}

\makeatletter
\newcommand{\saveequation}[2]{
  #2 \label{#1}
  \protected@write\@mainaux{}{\string\SAVEEQUATION{#1}{\unexpanded{\unexpanded{#2}}}}%
}
\newcommand{\savetagequation}[3]{
  #3 \label{#1} \tag{#2}
  \protected@write\@mainaux{}{\string\SAVEEQUATION{#1}{\unexpanded{\unexpanded{#3}}}}%
}
\newcommand{\SAVEEQUATION}[2]{%
  \global\@namedef{SAVEDEQUATION@#1}{#2}%
}
\newcommand{\repeatequation}[1]{%
  \ifcsname SAVEDEQUATION@#1\endcsname
    \@nameuse{SAVEDEQUATION@#1}\tag{\ref{#1}}%
  \else
    ?? \notag
  \fi
}
\makeatother

\setlist[itemize]{topsep=-4pt, partopsep=2pt}

\usepackage{setspace}

\usepackage{thmtools}
\declaretheoremstyle[%
  spaceabove=-2pt,%
  spacebelow=6pt,%
  headfont=\normalfont\itshape,%
  postheadspace=1em,%
  qed=\qedsymbol%
]{mystyle} 
\declaretheorem[name={Proof},style=mystyle,unnumbered,
]{prf}

\togglefalse{focs}

\begin{document}

\author{Ewan Davies}
\address{Department of Computer Science, Colorado State University, USA}
\email{\href{mailto:ewan.davies@colostate.edu}{ewan.davies@colostate.edu}}

\author{Holden Lee}
\address{Department of Applied Mathematics and Statistics, Johns Hopkins University, USA}
\email{\href{mailto:hlee283@jhu.edu}{hlee283@jhu.edu}}

\author{Juspreet Singh Sandhu}
\address{Department of Computer Science, Colorado State University, USA}
\email{\href{mailto:jsinghsa@ucsc.edu}{js.sandhu@colostate.edu}}

\author{Jonathan Shi}
\address{Chipletics Inc, Redmond WA, USA}
\email{\href{mailto:jshi@cs.cornell.edu}{jshi@cs.cornell.edu}}

\title[Potential Hessian Ascent III: sampling for the Sherrington--Kirkpatrick Model]{Potential Hessian Ascent III:\\ Sampling the Sherrington--Kirkpatrick Model at $\beta < 1/2$}

\begin{abstract}
\small
\noindent We give a polynomial-time algorithm to sample from the Gibbs measure of the Sherrington--Kirkpatrick model with negligible total-variation distance (TVD) error up to inverse temperature $\beta < 1/2$.
Prior work obtained TVD error guarantees only up to $\beta\approx 0.295$, while results covering the entire replica-symmetric regime $\beta < 1$ gave guarantees only in Wasserstein distance. 
\\ 

\noindent Our approach demonstrates that the same potential Hessian ascent previously developed for optimization also functions as a sampling algorithm by implementing algorithmic stochastic localization at high temperature.
By estimating the covariance of the tilted Gibbs distribution via Gaussian integration by parts, overlap concentration, and precise cavity estimates, we show that a Hessian-ascent process achieves an $O(1)$ Wasserstein error guarantee for finite-time localization, improving on the previous $o(n)$.
A careful comparison of stochastic localization with the Hessian ascent process and a free probability argument controlling the diagonal sub-algebra of the Hessian improves this to $O(1)$ in KL divergence.
We then use Jarzynski's equality with rejection sampling, along with entropy contraction on the time-$T$ localized distribution, to refine the error to $o(1)$ in TVD up to a constant time $T$ and to complete the sampling with the polarized walk.


\end{abstract}

\maketitle

\thispagestyle{empty}
\vspace{-5mm}
\renewcommand{\baselinestretch}{0.9}\normalsize
{
  \hypersetup{linkcolor=Red}
  \setcounter{tocdepth}{1}
  \tableofcontents
}
\renewcommand{\baselinestretch}{1.0}\normalsize

%
%
%
%
%
%

\newpage 
\pagenumbering{arabic}

\section{Introduction}

We give a polynomial-time algorithm to sample from the Gibbs measure of the Sherrington--Kirkpatrick (SK) model at inverse temperatures $\beta < 1/2$ with total variation error $o_n(1)$. 
The (random) Gibbs measure for the $n$-dimensional SK model is a probability measure on $\{-1,1\}^n$ given by
\[ \mu_{\beta A}(x) = \frac{1}{Z(\beta A)} e^{\frac{1}{2}\iprod{x,\beta A x}},\]
where $ A$ is a random $n\times n$ matrix from the Gaussian orthogonal ensemble and the normalizing constant
\[ Z(\beta A) = \sum_{x\in\{-1,1\}^n}e^{\frac{1}{2}\iprod{x,\beta A x}} \]
is the partition function. 
This means that $A$ is a symmetric matrix whose diagonal entries $A_{i,i}$ are independently distributed as $\calN(0,2/n)$ and whose off-diagonal entries $A_{i,j} = A_{j,i}$ for $i \ne j$ are independently distributed according to $\calN(0,1/n)$, where $\calN(\mu,\sigma^2)$ is the Gaussian distribution with mean $\mu$ and variance $\sigma^2$. 
The parameter $\beta>0$ represents inverse temperature, and with this parametrization, the \emph{replica-symmetric} (or high-temperature) regime corresponds to $\beta<1$. 
It is widely expected that a polynomial-time sampling algorithm exists throughout the replica-symmetric regime.
Our argument establishes this up to $\beta<1/2$.

\begin{theorem}[Informal version]
\label{t:informal}
    Let $0 <\be<1/2$. 
    Given any fixed $\de>0$, for all $n \ge n_0(\beta,\delta)$, with probability $1-\de$ over $A$, \pref{alg:main} gives a sample from a distribution that is $o_n(1)$ in TVD from $\mu_{\be A}$ in polynomial time.
\end{theorem}

\subsection{Spin glass models and sampling}
The SK model is the canonical example of a spin glass model: first introduced in statistical physics to study systems with dense, disordered interactions,
it has since become a central object in probability, with a Gibbs measure that exhibits intricate and rich structure, including \emph{replica symmetry breaking} at low temperatures~\cite{sherrington1975solvable, parisi1980order, parisi1980sequence, Gue01, guerra2002thermodynamic, Gue03, talagrand2006parisi, talagrand2010mean, talagrand2010mean2, panchenko2013parisi, panchenko2013sherrington, auffinger2015parisi, chen2023generalized, subag2024free}.
Interest in spin glass models in computer science is motivated by their utility as a model of noise in unsupervised learning, signal recovery, or hypothesis testing problems~\cite{montanari2014statistical,ma2016polynomial,hopkins2018statistical,aubin2019spiked,ashourian2023application}, their deep relationship to average-case constraint satisfaction problems~\cite{dembo2017extremal,panchenko2018k,gamarnik2021overlap,chou2022limitations,huang2022computational,jones2023random}, and their challenge as a testbed for algorithms working in high-dimensional and highly non-convex settings~\iffocs{\cite{subag2021following, montanari2021optimization,sandhu2024sum,jekel2024pha}}{\cite{subag2021following, montanari2021optimization,ivkov2024semidefinite,sandhu2024sum,jekel2024pha,jekel2025pha2}}, including noisy quantum algorithms~\cite{farhi2022quantum,basso2022quantum,boulebnane2025evidence}.

As random Gibbs measures on the hypercube, spin glass models are also a natural testing ground for sampling algorithms: the target measure is the stationary distribution of natural local Markov chain dynamics, such as Glauber dynamics, and a central question is whether these dynamics mix rapidly.
Average-case guarantees in these settings are well-motivated, with Bayesian posterior sampling as an example application where the target distribution is itself random, and with worst-case analyses often failing in regimes where algorithms actually succeed.
The SK model presents a prototypical challenge for sampling from these random measures induced by dense interactions~\cite[Open Problem 15]{bandeira2025randomstrasse101}.

One line of work proves functional inequalities for the Gibbs measure or transformations of it~\cite{eldan2022spectral,anari2022entropic}, as this immediately implies fast mixing of the associated Markov chain. 
For the SK model, this currently yields high-accuracy total-variation guarantees up to $\beta < 0.295$ \cite{anari2024trickle}.
These results rely on simple spectral properties of $A$, suggesting that further progress requires more fine-grained understanding of the Gibbs measure and the underlying random matrix.

A complementary line of work is \emph{Algorithmic Stochastic Localization} (ASL)~\cite{el2022sampling,huang2024sampling}.
Stochastic localization (SL)~\cite{Eld13} (equivalently, a diffusion model~\cite{Mon23}) is a measure-valued stochastic process that randomly evolves the original distribution $\mu_0$ to progressively more localized distributions $\mu_t$ while maintaining the Bayesian posterior property that sampling a random $\mu_t$ and then sampling $\sigma \sim \mu_t$ is equivalent to sampling from the original distribution $\sigma \sim \mu_0$.
ASL is the algorithmization of SL by replacing its computationally intractable component---namely the magnetization, or equivalently, the score---with efficiently computable surrogates.
Using Approximate Message Passing (AMP), prior work achieves sampling throughout the replica-symmetric regime $\beta < 1$ with $o(n)$ Wasserstein error~\cite{el2022sampling,celentano2024sudakov}.
We improve the error guarantee to $o_n(1)$ total variation error when $\beta < 1/2$ with a Hessian-based approach.

\subsection{Sampling strategy}
In SL, the \emph{tilt} or \emph{linear field} $y_t$ defines $\mu_t$ by the reweighting 
\[
    \mu_t(\sigma) \propto \exp\left(\iprod{\sigma,y_t}\right)\mu_0(\sigma)\,,
\] 
and the \emph{magnetization} $m_t$ is the average \[m_t := \sum_{\sigma} \sigma \,\mu_t(\sigma)\]
of the reweighted/conditional distribution.
The evolution of SL can be described in two equivalent ways that involve $m_t$: the usual stochastic differential equation (SDE) definition is
\begin{equation}
\label{eq:SL}
\tag{SL}
    dy_t = m_t\,dt + dB_t\,, 
\end{equation}
where the $m_t\,dt$ term gradually dominates as $t$ gets large.
A change of variables yields an equivalent SDE that we call Hessian Dynamics\footnote{The name reflects the general equivalence between the covariance of the Gibbs measure and the Hessian of the free energy with respect to the linear field~\cite[Proposition 3.1(a)]{wainwright2008graphical}. It also reflects that HD belongs to the same Hessian-based formalism as Potential Hessian Ascent~\iffocs{\cite{jekel2024pha}}{\cite{jekel2024pha,jekel2025pha2}}.} (HD):
\begin{equation}
\label{eq:HD}
\tag{HD}
    dm_t = \Cov(\mu_t)
    \,dB_t\,.
\end{equation}

Neither of these two SDEs are algorithmic due to the difficulty of estimating $m_t$ from $y_t$ or $\Cov_{\sigma \sim \mu_t}(\sigma)$ from $m_t$.
But mathematical spin glass theory gives us a third relationship between $m_t$ and $y_t$ by positing that at all times, they should jointly extremize an explicit function $\calF_{\TAP}(m,y)$ called the TAP free energy~\cite{TAP77, chen2023generalized}.
This \emph{TAP stationarity principle} is only approximately true, but the relationship it specifies is efficiently computable.

We obtain algorithmic versions of SL and HD by using the TAP stationarity principle to compute surrogates for the non-algorithmic $m_t$ in \pref{eq:SL} and the non-algorithmic covariance $\Cov(\mu_t)$ in \pref{eq:HD}. The former yields ASL-TAP and the latter Potential Hessian Dynamics (PHD). Although these arise from the SL and HD viewpoints respectively, both end up with the same diffusion coefficient for $dm_t$, namely the inverse TAP Hessian $(\nabla_m^2\calF_{\TAP}(m_t,y_t))^{-1}$; they differ only in drift (see \pref{s:locsamphess}).

The first major step of the analysis shows that the surrogate covariance $(\nabla_m^2\calF_{\TAP}(m_t,y_t))^{-1}$ is within $O(1)$ of the true covariance $\Cov(\mu_t)$ in Frobenius norm distance in expectation.
The argument is involved, using cavity interpolation estimates, overlap concentration, and various delicate symmetry arguments.
The result implies that $m_t$ in PHD stays within $O(1)$ Wasserstein distance of HD for constant time, improving on the $o(n)$ Wasserstein guarantee in~\cite{el2022sampling}.

The second step compares ASL-TAP to PHD-TAP to show that $y_t$ in ASL-TAP stays within $O(1)$ KL divergence error of SL for constant time.
ASL-TAP and PHD-TAP have the same diffusion terms, and their drift difference is determined by the diagonal entries of the square of the surrogate covariance, which we control by a free probability analysis.

Once this $O(1)$ KL divergence error is achieved, we boost it to $o(1)$ using rejection sampling.
To accomplish this, Jarzynski’s equality~\cite{Jar97} provides an explicit algorithmic density function by comparing ASL-TAP to a TAP-based scaffold of explicit density functions.

If we could run ASL to $T=\infty$, the localized measure would be a point mass. Since our analysis only permits $T=O(1)$, we stop at a large constant time $T=T(\be)$ and sample from $\mu_{T}$ using an appropriate Markov chain on the hypercube---the polarized walk---after showing entropy contraction holds for $\mu_T$ at large enough $T$.

\subsection{Technical challenges and contributions}
The analysis of our algorithm requires several major innovations or extensions of techniques beyond what is established in the literature, which we describe here.

\paragraph{Rejection sampling with Jarzynski's equality (\pref{s:sde})}
Jarzynski's equality (JE) is a nonequilibrium identity from statistical mechanics that yields importance weights for SDE-based sampling~\cite{Jar97,albergo2024nets}.
We use it to reweight ASL-TAP toward a scaffold of explicit densities, enabling rejection sampling.
To our knowledge, this is the first use of Jarzynski's equality in a theoretical sampling problem, showing that ASL can be corrected by reweighting probabilities instead of displacing mass,
and it may have broader applicability in algorithmic contexts.

\paragraph{Covariance estimate from overlap moments (\pref{s:est-cov})}
We show that the covariance of the Gibbs distribution of the planted SK model with plant-dependent linear field on the Nishimori line, and consequently also the SL-tilted SK model~\cite{el2022sampling}, is well-approximated by the inverse of the Hessian of the TAP free energy.
The main difficulty here is that, unlike the well-studied case of the SK model with independent external field, the tilt $y_t$ is itself a random function of the disorder $A$. 

El Alaoui and Gaitonde established this approximation for the standard SK model without plant, tilt, or field by decomposing the error into functions of overlaps~\cite{el2024bounds}, then invoking known cavity interpolation
~\cite[\S 1.6]{talagrand2010mean} and overlap concentration~\cite[\S 11]{talagrand2010mean2} estimates.
We generalize and extend their work by performing the decomposition for the planted model with plant-dependent linear field. Since no existing cavity estimates sufficiently covered the planted model, we also \emph{de novo} establish various analytic estimates made available by the cavity interpolation for the planted model (\pref{sec:cavity-interpolations-trace-identities}).

\paragraph{Regularity of ASL-TAP and PHD (\pref{sec:alg-properties})}
On a high-probability event in $A$, we show that the coefficients of ASL-TAP and PHD are uniformly bounded and Lipschitz. These regularity estimates control the SDE comparison error that goes into the $O(1)$ KL divergence bound, and they also imply concentration of Lipschitz functionals of the algorithmic trajectory, which is used to control the sample complexity of the Jarzynski reweighting.


\paragraph{Control of the diagonal through free probability (\pref{app:deformed-wigner-resolvent})}
We use a free probability argument to control the drift mismatch between ASL-TAP and PHD-TAP, since this mismatch is determined by the diagonal entries of a squared resolvent matrix.
The analysis is similar to a previous characterization for Potential Hessian Ascent~\cite{jekel2024pha}, though here the degree of control required is tighter, and singularities may arise because we cannot afford to regularize with an imaginary component in the resolvent.
The new ingredients are a chaining argument in place of an epsilon net, a cutoff to avoid real-axis singularities, and sharper edge and derivative estimates for the squared resolvent.



\paragraph{Entropy contraction for localized distributions on wedges (\pref{s:localized})}
The last technical challenge comes in sampling from the time-$T$ localized distribution $\mu_T$ for large enough $T$.
For this, we run the polarized walk restricted to a 
wedge 
of the hypercube\footnote{For technical reasons, the polarized walk is more convenient to analyze than Glauber dynamics.}. We show that this mixes rapidly by proving that 
the polarized walk satisfies entropy contraction with respect to such measures.
This is accomplished
through the framework of localization schemes~\cite{CE25}
by extending and combining Eldan, Kohler and Zeitouni's~\cite{eldan2022spectral} \emph{needle decomposition} technique into a 2-stage decomposition, with a bound on the covariance of the decomposed measures, itself established via a local limit theorem for Bernoulli random variables.

\paragraph{Hessian Ascent for sampling}
Potential Hessian Ascent~\cite{jekel2024pha} was introduced as an optimization algorithm for the SK model Hamiltonian at $\beta > 1$, with the covariance of the algorithmic increments chosen based on a resolvent of the Hessian of the generalized TAP free energy~\cite{chen2023generalized}. Here we show that at $\beta < 1/2$, the same Hessian-resolvent formalism yields a sampling algorithm: up to a small rank-one correction, PHD is the limit of PHA as the imaginary component in the resolvent and the step size both go to 0.
This demonstrates that Hessian Ascent is both an optimization algorithm and a sampling algorithm, and in fact implements a denoising diffusion model~\cite{Mon23}.


\paragraph{TVD sampling without a functional inequality}
Our method achieves $o_n(1)$ total-variation sampling under weaker hypotheses than current functional-inequality approaches appear to require. In particular, we do not assume a high-probability operator-norm bound on the covariance along the localization path~\cite{CE25,huang2025weak}---an estimate that is not presently known throughout the replica-symmetric regime---and instead combine $O(1)$ covariance control with Jarzynski correction and a final localized mixing argument.

\subsection{Prior work}
There are two main approaches to sampling from spin glass models. One is to prove functional inequalities (e.g.\ Poincaré or modified log-Sobolev inequalities), typically via stochastic localization, and deduce rapid mixing of Glauber or Langevin dynamics.
It can be sufficient to show \emph{weak} functional inequalities, which essentially establish the efficacy of these dynamics given a warm start or suitable annealing.
The other is to algorithmize stochastic localization itself. Our method combines these perspectives: we use ASL up to a large finite time, and then prove a functional inequality for the localized distribution.

\paragraph{Functional inequalities and fast mixing}
Eldan, Koehler and Zeitouni~\cite{eldan2022spectral}, taking inspiration from Bauerschmidt and Bodineau~\cite{bauerschmidt2019very},
were the first to give nontrivial sampling guarantees (i.e.\ to within small total variation distance) for the SK model up to $\be<1/4$ by proving a Poincaré inequality (PI). 
This is a corollary of their result which holds for any Ising model with a PSD interaction matrix $J$ satisfying $\opnorm{J}<1$, and they introduce stochastic localization as a tool to show functional inequalities. 
A limitation of their work is that it is a ``worst-case'' analysis over $\opnorm{J}$. 
The result is tight for the Curie--Weiss model, whereas for the SK model we expect the randomness inherent to the interactions to allow sampling at higher temperatures. 

Anari et al.~\cite{anari2022entropic} improved the result to a modified log-Sobolev inequality.
Chen and Eldan~\cite{CE25} established the general framework of localization schemes and streamlined the methods so that they rely only on bounding the covariance along the stochastic localization process.
In later work, Anari, Koehler, and Vuong~\cite{anari2024trickle} proved a log-Sobolev inequality that applies for the SK model up to $\be\approx 0.295$ through a more carefully designed scheme taking advantage of the symmetry of the spectrum of $A$. This is the best result to date for high-accuracy sampling for the SK model. 
Mikulincer and Sohn~\cite{mikulincer2025fastmixingisingmodels} show a modified log-Sobolev inequality for Ising models with a negative spectral outlier. All aforementioned works, when applied to the SK model, are limited by their use of simple properties of the spectrum of $A$ as a black box.

In another direction, Huang et al.~\cite{huang2025weak} showed a \emph{weak} Poincaré inequality for \emph{spherical} $p$-spin models by bounding the covariance along stochastic localization with high probability; this shows efficient sampling for \emph{annealed} Langevin dynamics (up to a so-called SL threshold in the inverse temperature that lies below the shattering threshold). An earlier work of Adhikari et al.~\cite{adhikari2024spectral} proved a Poincaré inequality for $p$-spin models on the hypercube with a different approach based on induction, giving the first results for $p>2$, though the inverse temperature threshold is weaker

Sampling with negligible TV error up to the conjectured threshold is known for the simpler Continuous Random Energy Model (CREM) on the hypercube~\cite{lee2024sampling}.

\paragraph{Algorithmic stochastic localization}
El Alaoui, Montanari and Sellke~\cite{el2022sampling} introduced algorithmic stochastic localization (ASL) to sample from the SK measure. By using AMP and natural gradient descent on (a custom version of) the TAP free energy to estimate the magnetization, they were the first to achieve nontrivial sampling guarantees throughout the high-temperature regime (initially for $\be<1/2$ and extended to $\be<1$ by Celentano~\cite{celentano2024sudakov}). The error guarantees are much weaker than TV distance, and they achieve a squared 2-Wasserstein distance of $o(n)$.

Huang, Montanari and Pham~\cite{huang2024sampling} showed that ASL obtains $o(1)$ TV error on spherical $p$-spin models by using a higher-order TAP correction to reduce aspects of the $p$-spin problem to the exactly solvable spherical 2-spin model. In contrast, for the Ising case, the SK model is itself already a 2-spin model. 
Wang et al.~\cite{WCZY25} also studied ASL for models supported on the hypercube but their analysis via Dobrushin uniqueness requires prohibitively small $\beta$ in the context of the SK model.

\begin{table}[t]

\newcolumntype{L}[1]{>{\RaggedRight\arraybackslash}m{#1}}
\centering
\small
\setlength{\tabcolsep}{12pt}
\renewcommand{\arraystretch}{1.48}
\begin{tabular}{@{}L{3.7cm}L{2.45cm}L{2.7cm}L{1.5cm}L{2.4cm}@{}}
\toprule
Work & Model & Method & Guarantee & Regime ($\beta<\cdots$) \\
\midrule
\multicolumn{5}{@{}l}{\emph{Functional inequalities / fast mixing}\vspace{0.4em}}\\
Eldan, Koehler, and Zeitouni~\cite{eldan2022spectral}; \linebreak Anari et al.~\cite{anari2022entropic}
& SK
& PI / MLSI via SL
& $\varepsilon$-TVD
& 
$1/4$ \\

Anari, Koehler, and Vuong~\cite{anari2024trickle}
& SK
& MLSI via SL
& $\varepsilon$-TVD
& 
$\approx 0.295$ \\

Huang et al.~\cite{huang2025weak}
& spherical $p$-spin
& weak PI via SL; \linebreak annealed Langevin
& $\varepsilon$-TVD
& SL threshold \\

Adhikari et al.~\cite{adhikari2024spectral}
& hypercube $p$-spin
& PI
& $\varepsilon$-TVD
& ${}<{}$ SL threshold \\

\midrule
\multicolumn{5}{@{}l}{\emph{Algorithmic stochastic localization}\vspace{0.3em}}\\
El Alaoui, Montanari, and Sellke~\cite{el2022sampling}; \linebreak Celentano~\cite{celentano2024sudakov}
& SK
& ASL + AMP
& $o(n)$-$W_2^2$
& 
initially $1/2$, later $1$ \\

Huang, Montanari, and Pham~\cite{huang2024sampling}
& spherical $p$-spin
& corrected ASL
& $o(1)$-TVD
& SL threshold \\

\midrule
\textbf{This work}
& SK
& ASL-TAP / PHD + JE
& $o(1)$-TVD
& $1/2$ \\
\bottomrule
\end{tabular}
\vspace{2.5mm}
\caption{Selected sampling guarantees for spin-glass models. All guarantees are polynomial-time. On the hypercube, functional-inequality--based results sample in $\ep$-TVD in $O(\log(1/\ep))$ time (noting just the dependence on $\ep$). We obtain $\ep$-TVD in $e^{O(1/\ep)}$ time.}
\label{tab:sampling-guarantees}
\end{table}

The (A)SL process is a reparameterization of a standard denoising diffusion model~\cite{Mon23}, and these models have been independently discovered by the machine learning community for use in generative modeling. 
A long line of work has established efficient sampling under access to an approximate score function~\cite{chen2022sampling,chen2023improved,benton2023nearly}. 

\paragraph{Jarzynski's equality} Jarzynski's equality \cite{Jar97} is an identity in statistical mechanics based on the Feynman-Kac formula \cite{kac1949distributions} (see also \cite[\S8.2]{oksendal2003stochastic}). It gives an SDE-ODE system to calculate changes in the partition function $Z_t = \int e^{-U_t(x)}dx$ as a potential $U_t$ evolves. It has recently seen wide application in sampling and in machine learning \cite{vargas2023transport,albergo2024nets}, with theoretical guarantees recently given in~\cite{guo2025complexity}.

\paragraph{Interpolations, overlaps and covariances}
The analysis of our algorithm requires understanding the overlap array of independent samples from the Gibbs measure of the planted SK model with random linear fields given by SL. 
Historically, the ``quadratic'' Guerra interpolation~\cite{Gue01,Gue03} and the cavity interpolation~\cite{aizenman2003extended, panchenko2013sherrington} have both played significant roles in reasoning about such overlap arrays. 
El-Alaoui and Gaitonde~\cite{el2024bounds} used such tools to show a bound on the operator norm of the covariance matrix for the Gibbs measure of the standard SK model (resolving a conjecture of Talagrand).  

Related Bayes-optimal/Nishimori work includes adaptive interpolation formulas of Barbier and Macris, the Franz--Parisi viewpoint of El Alaoui and Krzakala, general overlap-matrix concentration results of Barbier, and fluctuation refinements for overlap arrays on the Nishimori line due to Camilli, Contucci, and Mingione~\cite{barbier2019adaptive,el2018estimation,barbier2021overlap,camilli2023central}. These works provide together provide overlap concentration in planted models, but fall short of directly supplying the precise moment estimates needed for various terms that show up in the covariance estimate using the planted-SK model with the SL random external field.

\paragraph{Potential Hessian ascent, free probability, and control of the diagonal}
Subag pioneered the Hessian Ascent algorithm---following the top eigenvector of the Hessian of the objective function---to optimize spherical spin glass models~\cite{subag2021following}, which inspired both Montanari's AMP-based optimization work~\cite{montanari2021optimization} and the later Potential Hessian Ascent framework for the SK model~\cite{jekel2024pha,jekel2025pha2}. In the Ising/hypercube setting, PHA adds a coordinate-wise \emph{potential function} to the objective function (recovering the generalized TAP free energy~\cite{chen2023generalized} in the case of the SK model), showing that this adapts the algorithm to the geometry of the hypercube.
Resolvent bounds and free-probability methods then obtain control of the diagonal entries of the Hessian, which determines how individual coordinates evolve.
We use the same basic mechanism here in a sampling setting.

\subsection{Organization of paper}
\pref{s:locsamphess} introduces stochastic localization, Hessian dynamics, the algorithmic ASL-TAP/PHD-TAP framework, Jarzynski reweighting, and the reduction to sampling from the time-$T$ localized distribution. 
\pref{sec:proof-overview} in particular summarizes the main ideas and the structure of the remaining sections as they fit into the main argument.
\pref{s:sde} develops the quantitative SDE analysis and the comparison between the ideal and algorithmic processes. \pref{s:est-cov} and \pref{sec:cavity-interpolations-trace-identities} prove the main covariance estimate for tilted and planted SK measures using Gaussian integration by parts, overlap concentration, and cavity arguments. \pref{sec:alg-properties} establishes the regularity, Lipschitz, and diagonal-control estimates for the algorithmic coefficients, relying on the resolvent and free-probability analysis developed in \pref{app:deformed-wigner-resolvent}. Then \pref{s:localized} proves concentration of the time-$T$ localized measure on wedges and the functional inequalities needed for the final sampling step. 
The final section \pref{s:main-proof} combines these ingredients to prove the main theorem.

\subsection{Notation}
 $\calP(\R^m)$ refers to the set of all probability measures on the measure space $\left(\R^m, \calB(\R^m)\right)$.
 $\dist(x)$ is the law of $x$ when $x$ is a random variable.
 $d_{W,p}(X,Y)$ denotes the Wasserstein-$p$ distance between two distributions $X,Y \in \calP(\R^m)$. We denote the empirical measure of a vector by $\emp: \R^n \to \calP(\R)$, so that $\emp(v)$ is given by $x \mapsto \frac{1}{n}\sum_{i \in [n]} \delta(x - v_i)$, where $\delta$ is the Dirac delta measure.
 When $\Omega$ is an event in a probability space, $\one_{\Omega}$ is the indicator function for $\Omega$. $\Pr\Omega$ is the probability of $\Omega$.

$M_n(\bbD)$ refers to the set of $n \times n$ matrices whose entries are in $\bbD$. 
Similarly, $\mathcal{D}_n(\bbD)$ denotes the $n \times n$ diagonal matrices whose entries are in $\bbD$.
The map $\diag:\R^n\to \calD_n(\R)$ takes a vector to the diagonal matrix whose diagonal entries are given by that vector, and $\diag^*:\R^{n\times n}\to \R^n$ takes a matrix to the vector of its diagonal entries.
$E_{\calD_n}\left[X\right]$ is the diagonal matrix that has the same diagonal entries as $X \in M_n(\bbD)$, so that $E_{\calD_n}\left[X\right] = \diag(\diag^*(X))$.\footnote{
In operator theory, $E_{\calD_n}\left[X\right]$ is known as the non-commutative conditional expectation onto the diagonal subalgebra when $M_n(\bbD)$ is regarded as the base set of a non-commutative probability space, see \pref{sec:free-prob} for details.}
The all-ones matrix will be denoted as $J := \mathbf{1}_n\mathbf{1}_n^\sT$, where $\mathbf{1}_n \in \R^n$ is the all-ones vector.
The notation $\circ$ denotes the Hadamard (entrywise) product, and the algebra of $\calD_n$ is closed under this operation.
When $X$ is a function of other variables, $dX$ is its Frech\'et derivative given in \pref{def:frechet}.

Given $X \in M_n(\R)$, we define $\Spec(X) := \{\lambda \in \R: X - \lambda\Id_n\text{ is not invertible}\}$ to be the spectrum of $X$ (later we will use the same notation for more general operators).
By $\tr_n$ we denote the normalized trace so that $\tr_n[X] := \frac{1}{n}\sum_{i \in [n]}X_{i,i}$ when $X \in M_n(\R)$. $\Tr$ denotes the un-normalized trace.
The norms $\schnorm[p]{X}$ are normalized Schatten norms $\schnorm[p]{X} := (\tr_n[|X|^p])^{1/p}$, where $|X| := (XX^{\sT})^{1/2}$.
However, the vector $\ell^p$ norms $\lpnorm[p]{v} := (\sum_{i \in [n]} |v_i|^p)^{1/p}$ are not normalized for $v \in \R^n$.
We also use $\norm{X}_F$ to refer to the non-normalized Frobenius norm $\norm{X}_F := (\sum_{i,j\in [n]}X_{i,j}{}^2){}^{1/2}$, and $\opnorm{X}$ is the spectral/operator norm.
$\lipnorm{f}$ is the Lipschitz norm of a function $f$.
$\frenorm{dX}$ is given in \pref{def:frechet}.

In \pref{s:est-cov}, given a vector $\sigma \in \R^n$, $m(\sigma) := \frac{1}{n}\sum_{i=1}^n \sigma_i$ denotes its magnetization. The normalized inner-product between two (independent) samples form the Gibbs measure is denoted as $R_{1,2} := \frac{1}{n}\an{\sigma^1,\sigma^2}$. The notation $\left\langle f(a_1,\dots,a_m)\right\rangle_s$ denotes an expectation with respect to the Gibbs measure $G_{n,\beta,s}^{\ot m}$ at time $s$ where the $s$-dependent Hamiltonian may be deduced from context. When the Gibbs itself is random, then $\nu_s(f) := \E\left[\left\langle f(a_1,\dots,a_m)\right\rangle_s\right]$. For an ideal process $\bullet_t$, we denote the continuous algorithmic process as $\hat{\bullet}_t$ and the discretized algorithmic process as $\td{\bullet}_t$. In \pref{s:est-cov} and \pref{sec:cavity-interpolations-trace-identities}, we use the notation $\mathsf{span}\{a_1,\dots,a_k\}$ to denote a finite linear combination of terms made from $a_1,\dots,a_k$ with real coefficients that depend on $\beta$ and fixed integers $k$ (but not $n$). The terms $a_i$ will (often) be of the form $\nu(\cdot)$. The difference between the $i$-th spin of two configurations $\sigma^a$ and $\sigma^b$ is denoted $\Delta^{ab}_i$.

Let $d(x,y)$ denote the  Hamming distance between two vectors. Define $B_k(x_0),\,\pl B_k(x_0)\subset \{\pm 1\}^n$ and $\binom{[n]}{\le k},\,\binom{[n]}{k} \subset \{0,1\}^n$ as follows:
\begin{align*}
    B_k(x_0) &= \set{x\in \{\pm 1\}^n}{d(x,x_0)\le k} & \binom{[n]}{\le k} &= \set{x\in \{0,1\}^n}{\ve{x}_1\le k}\,,\\
    \pl B_k(x_0) &= \set{x\in \{\pm 1\}^n}{d(x,x_0)= k} & \binom{[n]}{k} &= \set{x\in \{0,1\}^n}{\ve{x}_1= k}.
\end{align*}
The objects $B_k$ and $\binom{[n]}{\le k}$ are termed ``wedges'', while the objects $\pl B_k$ and $\binom{[n]}{k}$ are termed ``slices'' of the hypercube. We will identify $\{0,1\}^n$ with $\mathcal P([n])$ (the power set of $[n]$) through $x\leftrightarrow \set{i\in [n]}{x_i=1}$.

Lastly, standard notation (and results) from probability theory, spin-glass theory, random matrix theory and Itô calculus are used in multiple places. See~\cite{anderson2010introduction, talagrand2010mean,talagrand2010mean2,karatzas2014brownian,durrett2019probability} for reference.  

\section{Localization, Hessians and sampling}\label{s:locsamphess}

One of our main tools is the stochastic localization framework of Eldan~\cite{Eld13,Eld20} and Chen and Eldan~\cite{CE25}. 
Given a sample space $\Omega$, a \emph{localization process} is a stochastic process $(\nu_t)_{t\ge 0}$ where $\nu_t$ is a probability measure on $\Omega$ and for all events $A\subset\Omega$ the process $t\to \nu_t(A)$ is a martingale such that $\lim_{t\to\infty}\nu_t(A)\in\{0,1\}$. 
A \emph{localization scheme} is then a mapping which assigns to each probability measure $\nu$ on $\Omega$ a localization process $(\nu_t)_{t\ge 0}$ with $\nu_0=\nu$.
An important consequence of the localization property is that $\nu_t$ must approach a Dirac measure as time progresses. Hence, a localization scheme is a method for interpolating from the given $\nu$ to a random Dirac measure via a measure-valued martingale. 
We call the measure $\nu_T$ obtained from $\nu$ by running a localization scheme for time $T$ a \emph{$T$-localized measure}.


\emph{Algorithmic stochastic localization} hinges on the observation that 
if one can accurately and efficiently simulate the localization process, then one can obtain an approximate sampling algorithm for the measure~\cite{el2022sampling,EMS25,huang2024sampling}, even without proving a functional inequality.

The localization scheme that we implement in our algorithm is the ``linear tilt localization'' that is commonly referred to as ``stochastic localization''. 
For convenience, in the following we assume the measure $\nu$ is supported on a subset of the hypercube. 
A concrete way to understand the measure-valued process is via a process $(y_t)_{t\ge 0}$ that generates a tilt.
We can think of $\si^*\sim \nu$ and define the tilt $y_t = t\si^* + B_t$ where $B_t$ is a standard Brownian motion. 
Then the \emph{linear-tilt localization} of $\nu$ is given by $\nu_t(\si) \propto \nu(\si)e^{\an{y_t,\si}}$.
Since $y_t/t\to \si^* \sim \nu$ as $t\to\infty$ almost surely, this description makes it clear that the process localizes in the sense that $\nu_t\Rightarrow \delta_{\si^*}$ as $t\to\infty$. It is also straightforward to check the required martingale property.

Let $m_t = \E_{\nu_t}[\si]$ be the magnetization of the measure $\nu_t$, and note that $m_t\in[-1,1]^n$. 
Another perspective on the linear-tilt localization is the fact that the tilt defined above is the solution to the following SDE,
\begin{align}
    \savetagequation{e:SL}{SL}{dy_t &= m_t dt + dB_t,& y_0&=0}.
\end{align}
We remark that this localization scheme is equivalent to the diffusion process
\begin{align*}
    dx_t &= dB_t,& x_0&\sim\nu,
\end{align*}
under the substitution $y_t = tx_{1/t}$.
It follows conceptually that to sample from the measure $\nu$ it suffices to simulate~\eqref{e:SL}. 
It is not at all clear, however, how to efficiently compute or approximate the drift term $m_t$ corresponding to the magnetization.
For more details on these assertions, see e.g.\ the paper of Chen and Eldan~\cite{CE25} introducing localization schemes, the survey of Montanari~\cite{Mon23} on the connection between linear-tilt localization and diffusion, and the citations therein.


We introduce Hessian Dynamics (HD) by means of a change of coordinates using Itô's lemma (\pref{lem:ito-local}).
From the definition of $m_t$ and the SDE~\eqref{e:SL} we obtain (see \pref{s:sde} for details)
\begin{align}
    \savetagequation{e:HD}{HD}{dm_t &= Q_tdB_t,& m_0&=\E_{x_0 \sim \nu}[\sigma]}.
\end{align}
where $Q_t=\Cov(\nu_t)$ is the covariance matrix.
By the localization property of~\eqref{e:SL}, as $t\to\infty$ we also have $m_t\to \si^*\sim \nu$, and hence to sample from $\nu$ it is also sufficient to simulate~\eqref{e:HD}.
The name ``Hessian dynamics'' comes from the general equivalence between the covariance of the Gibbs distribution and the Hessian of the free energy~\cite[Proposition 3.1(a)]{wainwright2008graphical}.

In principle, one could approximate~\eqref{e:SL} or~\eqref{e:HD} for constant time $T$, round the resulting magnetization $m_T$ to the cube to obtain a sample $\si_0$ (or equivalently take $\si_0=\sign(y_T)$), and hope this is good enough. 
We improve the accuracy of the sampler by running the polarized walk restricted to a ball around the point $\si_0$ for for $T'$ steps, establishing the effectiveness of this step via entropy contraction.
As explained above, 
this is insufficient for us to prove $o_n(1)$ TV distance error guarantees. The final idea we implement in the algorithm is to correct the distribution of $y_T$ with a round of rejection sampling using 
importance weights computed during the simulation of~\eqref{e:SL}.
Using $\hat{\,\cdot\,}$ to denote an algorithmic version of some quantity and $\td{\,\cdot\,}$ to denote a discretized algorithmic version, this results in the following outline for our algorithm. 
We give the full details in the following subsections.

\begin{algorithm}[!ht]
\caption{Informal description of our algorithm}
\label{alg:informal}
\begin{algorithmic}[1]
\Repeat{}
\State Initialize $y_0=0$.
\State Approximate~\eqref{e:SL} for time $T=O(1)$ to obtain a tilt $\td y_T$, while computing importance weights $\td w_T$.
\State Accept $\td y_T$ with some probability computed from $\td w_T$.
\Until{accept} \Comment{Now $\td y_T$ is approximately from the time-$T$ distribution of \eqref{e:SL}.}
\State Initialize $\si_0=\sign(\td y_T)$.
\State Run the polarized walk (\Cref{d:pw}) restricted to a wedge $\{\si : d_H(\si_0,\si)\le k\}$ for $T'$ steps and output the final point.
\end{algorithmic}
\end{algorithm}

\subsection{Algorithmic processes}

Fixing $\be>0$, we augment the SK model with an external field $y\in\mathbb{R}^n$ and define
\[ \mu_{\beta A,y}(\si) = \frac{1}{Z(\beta A,y)} e^{\frac{1}{2}\iprod{\si,\beta A \si} + \iprod{y,\si}}.\]
Note that we generalize the partition function $Z$ to the case of an external field as necessary. 
We denote expectations with respect to the measure $\mu$ with $\an\cdot$, so that a function $g:\{\pm 1\}^n\to \mathbb{R}$ has expectation
\[
\an{g(\sigma)}_{\beta A, y} := \E_{\mu_{\beta A,y}}[g(\sigma)].
\]
Where $\beta A$ is clear from context, we omit it from subscripts and the arguments of the function $Z$. For the rest of this section we suppose an arbitrary fixed $\beta A$.

Algorithmic versions of~\eqref{e:SL} and~\eqref{e:HD} arise if one replaces the quantities involved in the SDEs with approximations that are carefully chosen to be efficiently computable. 
For the SK model we can use the so-called TAP free energy of Thouless, Anderson and Palmer~\cite{TAP77} for this purpose.
Note that we retain the continuous-time perspective for now and separately consider discretization at the end.

For a partition function $Z$, let the free energy be $\calF = \log Z$ (note that this definition may differ by multiplicative terms from standard definitions). 
For the SK model with fixed $\beta A$, and varying external field $y$, we consider $\calF$ a function of $y$. 
By a direct computation, that the gradient $\partial_y\calF$ is the magnetization $\an{\si}_y$, and the Hessian $\nabla^2_y \calF$ is the covariance matrix of the Gibbs measure $\mu_{\be A, y}$.

Writing $\calF^*(m)$ for the the Legendre transform of the free energy, which is valid as $\calF(y)$ is convex, we have
\[ \calF^*(m) := \sup_y \Big\{\iprod{m,y}-\calF(y)\Big\}. \]
A little calculus and the fact $\partial_y\calF=\an{\si}$ shows that for fixed $m$, the value of $y$ which achieves the supremum is precisely one such that $\an{\si}_y=m$.
An important mathematical property of the Legendre transform is that the Hessians $\nabla^2_m\calF^*(m)$ and $ \nabla^2_y\calF(y)$ are inverses of each other. 
Thus, one might take the inverse of the Hessian of an approximation to $\calF^*(m)$ to derive an approximation of the covariance $\Cov(\mu_{\be A, y})$. 
This is one of the key insights of algorithmic approaches to the SK model, where the TAP free energy provides an approximation of $\calF^*(m)$.

The ASL process of El Alaoui, Montanari and Sellke~\cite{el2022sampling} replaces the term $m_t$ in~\eqref{e:SL} with an approximation $\hat m_t$ computed with AMP (and natural gradient descent to ensure certain regularities), and the computation is designed such that $\hat m_t$ is close to a stationary point of a custom TAP free energy. 
Our perspective, which sidesteps AMP, is to change coordinates and replace the covariance $Q_t$ in~\eqref{e:HD} with an approximation obtained from the (inverse of the) Hessian of the classic TAP free energy.
The accuracy of this approximation and the connection to the operator-theoretic notion of a resolvent has been the subject of several recent works~\cite{ABSY21,BSXY24,el2024bounds}.

The classic TAP free energy is given by 
\begin{align}
\label{e:FTAP}
\calF_{\TAP}(m,y) 
&:= -\fc{\be}{2}\iprod{m,A m} + \iprod{y,m} - \sumo in h(m_i) - \rc 4 n\be^2 \pa{1-\rc n \norm{m}_2^2}^2, \qquad\text{where}\\
\label{e:h}
h(x) &:= -\fc{1+x}{2}\log \pf{1+x}2 - \fc{1-x}{2}\log \pf{1-x}2,
\end{align}
and for fixed $y$ one can think of this as an approximation to the Legendre transform of the free energy of the system. 
Standard computations and properties of the Legendre transform give that the covariance $Q_t = \Cov(\mu_{y_t})$ we see in~\eqref{e:HD} is equivalently the Hessian $\nabla_y^2 \calF(y)$ (evaluated at $y=y_t$) or the inverse of the Hessian $\nabla_m^2 \calF^*(m)$ evaluated at $m=m_t = \an{\si}_{y_t}$ and $y=y_t$. 
Then it is natural to define an algorithmic version of Hessian Dynamics~\eqref{e:HD}, which we call Potential Hessian Dynamics (PHD), where one replaces the true covariance $Q_t$ with the inverse of the Hessian of $\calF_{\TAP}$. 
The use of the word ``potential'' refers to the role of the binary entropy $h(x)$ as a \emph{potential function} that adapts the process to the geometry of the hypercube, as in previous works on Hessian ascent~\iffocs{\cite{jekel2024pha}}{\cite{jekel2024pha,jekel2025pha2}}.

Here we see the parallel formalism:~\eqref{e:SL} considers the evolution of the external field $y_t$ and its algorithmic counterpart is obtained by replacing the ``true'' drift term $m_t$ with an approximation that can be computed efficiently given $y_t$ (note that the noise term is trivial). 
Motivated by the coordinate change inherent to the Legendre transform, \eqref{e:HD} considers the evolution of the magnetization $m_t$ and its algorithmic counterpart is obtained by replacing the ``true'' noise term $Q_t$ with something that can be computed efficiently given $m_t$ (note that the drift term is trivial).

\subsection{An unattainable trio}\label{s:alg-process}

Recall that we have fixed $\beta A$ and consider the ideal processes~\eqref{e:SL} and~\eqref{e:HD} for the SK measure with external field $y_t$
\begin{align}
    \repeatequation{e:SL},\\
    \repeatequation{e:HD},\\
    \mg_t &= \an{\si}_{y_t},& Q_t &= \Cov(\mu_{y_t}).\nonumber
\end{align}

Assuming that we have access to the external field, we obtain algorithmic versions of these processes by replacing the terms $m_t$ and $Q_t$---which involve averages over the true Gibbs measure $\mu_{y_t}$---with terms that are easier to compute. 
The TAP free energy gives us both: for fixed $y$ we expect the true magnetization $m$ to be close to a stationary point of $\calF_{\TAP}$ (i.e.\ to satisfy the classic ``TAP equations'') and the true covariance $Q$ to be close to $(\nabla^2 \calF_{\TAP})^{-1}$. 
The accuracy of the TAP equations for the magnetization at high temperature is covered, e.g., in Talagrand's monograph~\cite{talagrand2010mean}. For the covariance, the TAP approximation has been proved accurate in the case of no external field~\cite{el2024bounds,ABSY21} and at high enough temperatures~\cite{BXY25a}. 
Since the fields introduced by~\eqref{e:SL} depend on the random interaction matrix $\beta A$, which is not exactly a classic setting for mathematical studies of the SK model, none of these results apply in our setting and we must establish analogues where necessary. 

Inspired by the TAP free energy, given the external field $\hat y_t$ we want $\hat m_t$ to be a stationary point of $\calF_{\TAP}$, and we can approximate the covariance with $\hat Q(\hat m_t)$ where $\hat Q : [-1,1]^n\to\R^{n\times n}$ is the operator-valued function
\begin{equation}
\savetagequation{e:hQ}{Q\textsubscript{TAP}}{\hQ(m) = \left(-\beta A + \diag\prc{1-m^2} + \beta^2\left(1-\fc{\norm{m}_2^2}{n}\right)\Id_n - \fc{2\beta^2}{n} mm^{\sT}\right)^{-1}}.
\end{equation}
That is, we want the following system of equations:
\begin{align}
    \savetagequation{e:ASL}{ASL}{d\hat{y}_t &= \hat m_t dt+ dB_t,& \hat y_0 &=0},\\
    \savetagequation{e:PHD}{PHD}{d\hat m_t &= \hat Q(\hat m_t)dB_t,& \hat m_0 &=0},\\
    \savetagequation{e:TAP}{TAP}{0 &= \gd_m \calF_\TAP(\hat m_t,\hat y_t)}.
\end{align}
See \eqref{e:TAP-expanded} for an explicit expression for $\gd_m \FT$.
The intuition behind the first two equations is clear from the ideal processes, and~\eqref{e:TAP} plays the role of the Gibbs average relating $m_t$ and $y_t$ for the algorithmic analogues $\hat\mg_t$ and $\hat y_t$. 
For use later, we define $\hat\mg(y)$ to be the function that, given an external field $y$, returns an magnetization $\hat\mg$ such that $0 = \gd_m \calF_\TAP(\hat m,y)$; and we define $\hat y(m')$ similarly to be the function that, given a magnetization $m'$, returns a tilt $\hat y$ such that  $0 = \gd_m \calF_\TAP(m',\hat y)$. These functions are used as the algorithmic analogue of computing $\an{\si_y}$ and solving for the field $y$ such that $m'=\an{\si}_y$ respectively.

The trio of equations above form an overdetermined system, hence we can choose to satisfy a pair of the equations and derive the consequences for the third. 
The simplest option is to just take the ASL and PHD processes as written:
\begin{equation}
\savetagequation{e:ASL-PHD}{ASL-PHD}{\begin{cases}
\begin{aligned}
        d\hat y_t &= \hat \mg_tdt+dB_t\\
        d\hat \mg_t &= \hat Q_t(\hat \mg_t)dB_t
\end{aligned}
    \end{cases}}
\end{equation}
This, however, does not enforce any kind of consistency between $\hat y_t$ and $\hat \mg_t$ that parallels the ideal situation $m_t = \an{\si}_{y_t}$. 
It turns out to be more useful to enforce that $(\hat \mg_t, \hat y_t)$ satisfy~\eqref{e:TAP}, choose one of ASL and PHD, and use Itô's lemma to derive a version of the other. For ASL, from the calculations in \pref{s:sde}, the resulting system is
\begin{align}
\repeatequation{e:ASL-TAP}
\end{align}
where with $D(m) = \diag\prc{1-m^2}$ the function
\begin{align}\saveequation{e:Itomag}{f(\mg) &= \left(E_{\calD_n}\left[\hQ(\mg)^2\right]D(\mg)^2 - \frac{\beta^2}{n}\left(\Tr\left[\hat{Q}^2(\mg)\right]\Id_n +2\hat{Q}^2(\mg)\right)\right)\mg}
\end{align}
represents a modified magnetization term. 
We show that at high temperatures, the difference between the ideal process given by \eqref{e:SL}, \eqref{e:HD} and the algorithmic process \eqref{e:ASL-TAP} is small. 
For the analysis, we first prove that the true covariance $Q_t$ and the approximation $\hat Q(m_t)$ are $O(1)$ close on average under the ideal process; integration (via Gr\"onwall's inequality) then shows the same for the magnetization $\mg_t$ and $\hat\mg(y_t)$.

However, $O(1)$ closeness is not enough to obtain negligible ($o_n(1)$) total variation distance. One approach, taken by \cite{huang2024sampling} for the spherical $p$-spin model, is to add carefully-chosen higher order correction terms to the TAP free energy. We instead take an alternate route: correct the output of ASL-TAP via rejection sampling with importance weights derived from Jarzynski's equality.


\subsection{Rejection sampling with importance weights from Jarzynski's equality} 
\label{s:iw}

Rejection sampling is a method of sampling a distribution $p(x)$ by choosing an easier-to-sample proposal distribution $q(x)$ such that $p(x)\le C q(x)$, drawing a sample from $q(x)$, and then accepting it with probability $\fc{p(x)}{Cq(x)}$. This works if the bound is satisfied pointwise; if it is sampled with high probability over $x\sim q$, then we incur some additional error. 

However, we do not have access to the probability density of the sample from \eqref{e:ASL} to correct it to \eqref{e:SL} by rejection sampling. Instead, we will instead first target a certain distribution $\rh_T$ whose density is explicit and hence acts as a scaffold; Jarzynski's equality gives us the importance weights to sample from this. We then $\rh_T$ as a proposal distribution, to sample from the SL density. 
Technically, since we only obtain weighted samples from $\rh_T$, we incorporate those weights into the acceptance ratio.

Jarzynski's equality is a method of designing a stochastic differential equation to compute an unbiased estimator for a partition function given an annealing sequence of distributions. 
It can also be used as an asymptotically correct sampler for a measure $\rh_t\propto e^{-U_t}$ for some weights $U_t$, or as an estimate of an expectation under $\rh_t$. 
We take the statement below from Albergo and Vanden-Eijnden~\cite[Proposition 3]{albergo2024nets}.
\begin{theorem}[Jarzynski's equality~\cite{vaikuntanathan2008escorted,vargas2023transport,albergo2024nets}]\label{t:je}
Let $\rh_t \propto e^{-U_t}$ be a sequence of probability densities where $U_t$ is twice differentiable with respect to $(x,t) \in [0,T]\times \R^n$.
Let $b_t$ be continuously differentiable, and suppose that uniformly for $t\in [0,T]$, $\gd U_t$ and $b_t$ are Lipschitz in $x$. Consider the coupled SDE-ODE system
    \begin{align*}
        dx_t & = \ba{-\rc2\gd U_t(x_t) + b_t(x_t)}dt + dB_t, &x_0&\sim \rh_0\\
        dw_t &= \ba{\gd \cdot b_t(x_t) - \an{\gd U_t(x_t), b_t(x_t)} - \pl_t U_t(x_t)}dt, & w_0&=0.
    \end{align*}
    Then $\fc{Z_t}{Z_0} = \fc{\int_{\R^n} e^{-U_t(x)}dx}{\int_{\R^n} e^{-U_0(x)}dx} = \E[e^{w_t}]
    $
    and for any bounded measurable function $h$,
    \[
\int_{\R^n} h(x)\rh_t(x) dx = \fc{\E[e^{w_t}h(x_t)]}{\E[e^{w_t}]}.
    \]
\end{theorem}

We can use Jarzynski's equality with either \eqref{e:ASL-TAP} or \eqref{e:PHD-TAP} to maintain importance weights $\hat w_t$ designed to correct the distribution of $\hat y_t$ to follow a distribution $\rh_t$ that satisfies 
\begin{align}
\saveequation{e:rho-t}{\rh_t(y) &\propto e^{\FT(\hat\mg(y),y) - \fc{\ve{y}^2}{2t}}},
\end{align} 
where we recall that $\hat\mg(y)$ is the magnetization $\hat\mg$ such that $0 = \gd_m \calF_\TAP(\hat m,y)$.
Recall that the function $f$ defined in~\eqref{e:Itomag} gives the ``modified'' magnetization one gets from Itô's lemma in the derivation of ASL-TAP and PHD-TAP.
Augmenting ASL-TAP with importance weights, we obtain
\begin{align}
\repeatequation{e:asl-tap-je}
\end{align}
PHD-TAP-JE can be similarly derived; see \pref{s:asl-phd-tap}.
Our algorithm computes a standard (Euler-Maruyama) discretization of the system~\eqref{e:asl-tap-je} up to some constant time $T=T(\beta)$ and uses the importance weights to ``correct'' the resulting density to be $\rho_T$. 
Now that we have a known density, we can further correct towards a distribution on tilts arising from the ideal SL process with rejection sampling.
This lets us sample a tilt $\td y_T$ and round it to via cube via $\si_0 = \sign(\td y_T) \in\{-1,1\}^n$. 
We then run the polarized walk restricted to a Hamming ball of radius $k=\ep(\be) n$ around $\si_0$ and output the resulting point on the hypercube.
\label{s:T-loc}

It remains to show that (1) we can estimate the density $p_T$ of the random tilt $y_T$ obtained by running SL for time $T$, and (2) we can sample from the localized distribution $\mu_{y_T}$, both with good probability over $y_T\sim p_T$. 
The standard reduction from partition function estimation to sampling via simulated annealing means it suffices to prove (2).
Technically, to avoid erroneously accepting tilt samples with low probability under $p_T$, we show that with high probability over $A$ that we can obtain an \emph{underestimate} of $p_T(y)$ for all $y$, since we can always sample $\mu_{y_T}$ restricted to a small neighborhood (wedge) around the rounded $y_T$, and it is only with good probability over $y_T\sim p_T$ (rather than $\rh_T$) that the neighborhood contains most of the mass.

Efficient sampling from the $T$-localized distribution is shown in \pref{s:localized} via proving entropy contraction for the polarized walk.



\subsection{Algorithm} 

We display the full algorithm in \pref{alg:main}.  
First, we run ASL-TAP-JE 
as derived in \pref{s:alg-process}--\ref{s:iw}; we separate this out into \pref{alg:asl-ta-je-dre}.
We conduct rejection sampling as described in \pref{s:T-loc} to turn our our weighted sample for $\rh_T$ into an approximate sample for $\mu_T$. Finally, we sample from the localized distribution $\mu_{y_T}$ using the polarized walk restricted to a neighborhood.

\begin{algorithm}[!ht]
\caption{SK Sampler}
\label{alg:main}
\begin{algorithmic}[1]
\INPUT matrix $A$, inverse temperature $\be$, total time $T$ which is a multiple of step size $\eta$, total steps of polarized walk $T_{\mathrm{PW}}$, wedge size $k$.
\OUTPUT Approximate sample from Gibbs measure $\mu_A$.
\State \textbf{Part 1:} Approximately sample from $\mu_T$. 
\State Run \pref{a:ars} (Approximate Rejection Sampler) 
with subroutine for sampling and density ratio estimation given by \pref{a:asl-ta-je-dre} (ASL-TAP-JE with density ratio estimation) to obtain $\td y_T$.
\State \textbf{Part 2:} Approximately sample from localized distribution $\mu_{\beta A,\td y_T}$.
\State Initialize $\si_0 = \sign(\td y_T)$.
\State Run the polarized walk (\Cref{d:pw}) restricted to $\set{\si}{d_H(\si_0,\si)\le k}$ for $T_{\mathrm{PW}}$ steps and output sample.
\end{algorithmic}
\end{algorithm}

\begin{algorithm}[!ht]
\caption{(ASL-TAP-JE) with density ratio estimation}
\begin{algorithmic}[1]
\INPUT matrix $A$, inverse temperature $\be$, total time $T$, step size $\eta$, number of steps $S$ for mirror descent, error parameter $\ep$
\State Let $\td y_0=0$, $\td \mg_0=0$, $\td x_0=0$, and $\tw_0=0$. 
\For{$t\in \{0,\eta,\ldots, T-\eta\}$}
    \State Draw $\xi_t \sim \calN(0,I_n)$. 
    \State $\td y_{t+\eta} = \td y_t + \eta \td\mg_t + \sqrt{\eta}\xi_t$.
    \State $\td x_{t+\eta}^{(0)} = \td x_t +\eta \ba{D(\td\mg_t)\hQ(\td\mg_t) (\td \mg_t - f(\td \mg_t)) + \ED\pa{\hQ(\td\mg_t)^2} D(\td \mg_t)^2\td \mg_t} + \sqrt{\eta}D(\td \mg_t)\hQ(\td \mg_t) \xi_t$.
    \For{$s\in \{0,\ldots, S-1\}$} \Comment{Solve \eqref{e:TAP} for $\mg$ when $y=\td y_{t+\eta}$ using mirror descent}
    \label{line:asl-tap-je-mirror-descent}
        \State $\td x_{t+\eta}^{(s+1)} = \td x_{t+\eta}^{(s)} - \lm\gd_m \FT(\tanh(\td x_{t+\eta}^{(s)}),\td y_{t+\eta})$ where $\gd_m \FT$ is given in \eqref{e:TAP-expanded} and $\lm$ is as in \ref{d:Q-reg}.
    \EndFor
    \State $\td x_{t+\eta} = \td x_{t+\eta}^{(S)}$.
    \State $\td m_{t+\eta} = \tanh(\td x_{t+\eta})$. 
    \State $\td w_{t+\eta}=\td w_t + \fc{\eta}2 \ba{\Tr(\hQ(\td\mg_t))+ \ve{\td\mg_t}^2}$.
\EndFor
\State Let $\hat Z_{\be A, \td y_T}$ be the output of \pref{alg:sa} (Simulated annealing with partition function estimation) with temperature ladder $\be_\ell = \fc{(\ell-1)\be}{4n}$ ($1\le \ell \le 4n+1$), $N=\Te\pf{n}{\ep^2}$, $R=\Om\pa{\log\prc{\ep}}$, initial partition function $\sumz ik \binom ni$, with algorithm the polarized walk (\Cref{d:pw}) initialized at $\sign(\td y_T)$ and restricted to $\ball{\ep(\be) n}{\sign(\td y_T)}$ with $\Te\pa{n\log\pf n\ep}$ steps ($\ep(\be)$ is given by \pref{t:ls-sk}).
\OUTPUT $\pa{\td y_T , e^{\tw_T}\fc{\hat Z_{\be A, \td y_T}}{e^{\FT(\td \mg_T,\td y_T)}}}$
\end{algorithmic}
\label{a:asl-ta-je-dre}
\label{alg:asl-ta-je-dre}
\end{algorithm}

\subsection{Desiderata}
\label{s:desiderata}
To analyze the algorithm sketched above, it suffices to establish the following desiderata for $\be<1/2$. We will prove (d1) in \pref{s:est-cov} and then (d0) and (d2)--(d5) in \pref{sec:alg-properties}, and then show how they will give us guarantees for our sampling via the stochastic processes described in \pref{s:sde}. We will prove (d6) in \pref{s:localized} to show sampling from the localized distribution.

Below, constants are allowed to depend on $\be$ but not $n$, and $T$ is a constant which depends on $\be$. Recall that $\hQ(\mg)$ is defined as in \eqref{e:hQ} and define $Q(m)$ such that $Q(m)$ is the covariance matrix $\Cov(\mu_{\beta A,y})$ for the unique external field $y$ such that $\an{\si}_{\beta A,y}=m$.
Recall also that $D(\mg) = \diag\left(\frac{1}{1-\mg^2}\right)$.
Then both $\hQ$ and $Q$ depend implicitly on $\be A$. We would like the following to hold with high probability over $A$. (The expectations below are with respect to the various random processes, for fixed $A$.)
\begin{enumerate}[itemsep=0.7em,label={(d\arabic*)},start=0]
\item \label{d:SDE-existence-uniqueness-stay-in-cube} \textbf{(SDE existence, uniqueness, and non-escape, \pref{cor:alg-desiderata})}
    The algorithmic SDEs \eqref{e:PHD} and \eqref{e:asl-tap-je} have unique strong solutions such that $\hat{m}_t \in (-1,1)^n$ for all $t < \infty$ almost surely.
    \item \label{d:Q-error}
    \textbf{(Error of covariance estimate, \pref{thm:covar-estimate})} 
    For all $t\in [0,T]$, for $\mg_t$ drawn from HD:
    \[\E \ve{\hQ(\mg_t)^{-1}Q(\mg_t) - \Id_n}_F^2 \le \ep(t)^2,\]
    where $\ep(t)$ is a constant depending only on $\be$ and $t$ and bounded on $[0,T]$.
    \item \label{d:Q-reg}
    \textbf{(Convexity and second-order regularity of $\FT$, \pref{cor:alg-desiderata})} There are constants $\lambda,\Lambda>0$ such that 
    for all $\mg\in (-1,1)^n$,
    \[
        \lambda D(\mg)^{-1} \preceq \hQ(\mg) \preceq \Lambda D(\mg)^{-1}.
    \]
    \item \label{d:drift-error}
    \textbf{(Drift error in PHD-TAP, \pref{cor:alg-desiderata})}
    For all $t\in [0,T]$, for $\mg_t$ drawn from HD, 
    \[
    \E\ve{f(\mg_t)-m_t}_2^2 \le \ep_{\mathrm{drift}}(t)^2,\]
    where $\ep_{\mathrm{drift}}(t)$ is a constant depending only on $\beta$ and $t$ and bounded on $[0,T]$.
    \item \label{d:Q-Lip}
    \textbf{(Lipschitzness of diffusion and drift terms, \pref{cor:alg-desiderata})}
    \begin{enumerate}
        \item There exists a constant $L$  such that for all $\mg, w\in (-1,1)^n$, 
    \[
    \ve{\hQ(\mg) - \hQ(w)}_F \le L\ve{\mg-w}_2.
    \]
        \item
        There exists a constant $L_f$  such that for all $\mg, w\in (-1,1)^n$, with $f$ defined as in \eqref{e:Itomag},
    \[
    \ve{f(\mg) - f(w)}_2 \le L_f\ve{\mg - w}_2.
    \]
    \end{enumerate}
    \item \textbf{(Bounded mgf of JE weights, \pref{cor:alg-desiderata})}
    \label{d:JE-warm-start}
    Let $(\hat y_t, \hat \mg_t,\hat w_t)$ evolve as in \eqref{e:asl-tap-je}. 
    Then for all $\lambda \in \R$, there is some finite $C_{\mathrm{JE},T,\lambda} > 0$ so that
    \[ \E
    \left[\exp\left(\lambda\hat{w}_T - \E
    [\lambda\hat{w}_T]\right)\right] \le C_{\mathrm{JE},T,\lambda}. \]
    
    \item \label{d:lsi}
    \textbf{(Entropy contraction for localized distribution, \Cref{t:ls-sk})} With high probability over $A$, with $1-e^{-\Om(n)}$ probability over $y_T$ from the ASL process, the localized distribution $\mu_{A, y_T}$ is concentrated on a wedge $\ball{\ep n}{\sign(y_T)}$, and $\mu_{A, y_T}|_{\ball{\ep n}{\sign(y_T)}}$ satisfies a entropy contraction with constant $\ge\fc{\rh_\be}{n}$.
\end{enumerate}
These desiderata  encapsulate the main technical challenges. The main sampling result then follows from these via  general arguments.
In particular, \ref{d:SDE-existence-uniqueness-stay-in-cube}, \ref{d:Q-error}, and \ref{d:Q-Lip}(a) are sufficient to conclude $O(1)$ Wasserstein-2 error up to a constant time $T$, by a comparison of PHD with HD using Gr\"onwall's inequality (\pref{cor:wasserstein-bound-asl}). 
Adding \ref{d:Q-reg}, \ref{d:drift-error}, and \ref{d:Q-Lip}(b) improves this to a $O(1)$ error in KL divergence up to time $T$, by comparison of ASL-TAP with SL via PHD (\pref{c:O1-KL}) and mirror descent to solve the TAP equation (\pref{s:mirror-tap}). 
Finally, adding \ref{d:JE-warm-start} and \ref{d:lsi} enables sampling $\mu_0$ with $o(1)$ total variation error via rejection sampling and the polarized walk (\pref{s:main-proof}).

\begin{remark}\hspace{1em}
    \begin{enumerate}[itemsep=0.7em]
        \item By ``with high probability over $A$'', we mean that for $\de>0$ arbitrarily small, with probability $1-\de$ over the draw of $A$, the desiderata hold with constants depending on $\de$. 
        Hence, for the desiderata involving expectations, \ref{d:Q-error} and \ref{d:drift-error}, by Markov's inequality it suffices to prove an analogous bound on the expectation over $A$ as well, restricted to a ``good'' set. For instance, for \ref{d:Q-error}, it suffices to show for some $A$-measurable good set $G$ such that $\Pr(A\in G)\ge 1-o(1)$ that 
            \[
\E_{A,\mg_t} \one_G \ve{\hQ(A,\mg_t)^{-1}Q(A,\mg_t) - I_n}_F^2 =O(1)
    \]
    (where $\mg_t$ is drawn from HD with interaction matrix $\beta A$).
    \item Note \ref{d:Q-error} and \ref{d:drift-error} are with respect to an ideal process (HD), while \ref{d:JE-warm-start} is with respect to the algorithmic process (ASL-TAP-JE). 
    Conditions with respect to the ideal process are used as input to standard SDE error bounds, e.g.\ Girsanov's theorem.
    \item Two high-probability events on $A$ suffice to establish \ref{d:SDE-existence-uniqueness-stay-in-cube} and \ref{d:Q-reg} through \ref{d:JE-warm-start} for ASL-TAP-JE: one simple bound on $\opnorm{A}$ and one capturing control of the diagonal entries of $\hat{Q}^2$ through free probability arguments. These events are described in \pref{sec:alg-events}.
    \item Note \ref{d:Q-reg} and \ref{d:Q-Lip} will not hold for all  $\be < 1$, because there is a $\beta < 1$ for which the Hessian $\hat{Q}(m)$ of the TAP free energy fails to be PSD at all points $m \in (-1,1)^n$~\cite{gufler2023concavity}.
    \end{enumerate}
\end{remark}

\begin{corollary}[Expected operator norm bound on covariance of HD-tilted measure]
\label{c:Eop}
    With high probability over $A$, for each $t\in [0,T]$, there is a constant $K(t)>0$ such that when $\mu_t=\mu_{\beta A,y_t}$ is the SK measure tilted by running (SL) or (HD) for time $t$, $\E\ba{ \opnorm{\Cov(\mu_t)}^2} \le K(t)^2$.
\end{corollary}
\begin{prf}
Note that $\Cov(\mu_t) = Q(\mg_t)$.
    We have by desiderata~\ref{d:Q-error} and~\ref{d:Q-reg},
    \[
\La^{-1} \E\opnorm{Q(\mg_t)}^2 
\le \E \ba{\opnorm{\hQ^{-1}(\mg_t) Q(\mg_t) }^2}
\le \E \ba{\pa{1+\ve{\hQ^{-1}(\mg_t) Q(\mg_t) - I_n}_F}^2}\le 2+2\ep(t)^2. \qedhere   \]
\end{prf}

\subsection{Proof overview}
\label{sec:proof-overview}
In this subsection, we overview the design and structure of the proofs of the desiderata (and their interdependencies).

\paragraph{SDE error bounds and rejection sampling with JE (\pref{s:sde})} 
In \pref{s:asl-phd-tap}--\ref{s:je}, we derive the equations \eqref{e:asl-tap-je} and \eqref{e:phd-tap-je}, where the weights are such that we obtain an importance weighted sample from the scaffolding distribution $\rh_T(y)\propto e^{\FT(\hat \mg(y),y) - \fc{\ve{y}^2}{2T}}$. This is a natural choice of scaffolding distribution because it would equal the probability density of $y_T$ if $e^{\FT(\hat \mg(y),y)}$ were equal to $Z(\be A, y)$.

\ppart{Wasserstein bounds and integrating to get magnetization error (\pref{s:wasserstein})} A Gr\"onwall argument gives Wasserstein-2 error between the laws of 2 SDE's given Lipschitzness of the approximate process and closeness of the mean and diffusion terms under the ideal process (\pref{l:W-gw}). Applying this to \eqref{e:PHD-TAP}, we obtain a Wasserstein bound and that under \eqref{e:SL}, $\E\ba{\ve{\hat \mg(y_t)-\mg_t}^2} = O(1)$ (\pref{l:m-error}). This allows obtaining $O(1)$ KL divergence between \eqref{e:SL} and \eqref{e:ASL} up to time $T$.

\ppart{Rejection sampling with JE (\pref{s:je-error})}
We adapt rejection sampling to the case when only the unnormalized ratio between the target distribution $q$ and the proposal distribution $p$ is (approximately) known and has upper tails bounded, by first drawing a few samples and taking an appropriate quantile of the samples as normalization (\pref{l:ars}). In our setting, our weights are the product of $e^{\hw_T}$ from Jarzynski's equality and an estimate of the ratio $\dd{p_T}{\rh_T}$, where $p_T=\dist(y_T)$. The density $p_T$ is a partition function that can be approximated from efficient sampling of the tilted distribution $\mu_{\be A, y_T}$, with high probability under $y_T$ (shown in \pref{s:localized}). This is not necessarily true for $y_T$ from the proposal distribution, but because we always have mixing on a wedge, this can only result in an underestimate, which is tolerable. \pref{l:je-rs} is the main result which tabulates all sources of error: drift error, tails of the Jarzynski weights $\hw_T$, discretization error for the samples and weights (which are furnished in \pref{s:disc}), and ratio error. The main theorem is proved in \pref{s:main-proof} by plugging in all these sources of error, bounded through the desiderata. We additionally consider the error from approximating the solution to the TAP equation in \pref{s:mirror-tap}: the theory of mirror descent gives linear convergence since $\FT$ is relatively smooth and strongly convex in $m$ with respect to $D(\mg)$.

\paragraph{Estimate of the covariance (\pref{s:est-cov} and \pref{sec:cavity-interpolations-trace-identities})}\label{sec:tech-overview-covariance} The objective of these two sections is to justify the TAP prediction that, for the tilted Gibbs law \(G_{A,y_t,\beta}\), the covariance matrix
\[
P:=P(A,y_t):=\Cov(G_{A,y_t,\beta})
\]
is well approximated by the inverse Hessian of the TAP free energy at the Gibbs magnetization \(\mg_t:=\an{\sigma}\) (note that $P$ is referred to as $Q_t$ in other sections). Equivalently, if
\[
\hat Q^{-1}:=
\beta^2\Bigl(1-\frac{\|\mg_t\|_2^2}{n}\Bigr)\Id_n
-\beta A
+D(\mg_t)
-\frac{2\beta^2}{n}\mg_t\mg_t^\sT,
\qquad
D(\mg_t):=\diag\Bigl(\frac{1}{1-\mg_{t,1}^2},\dots,\frac{1}{1-\mg_{t,n}^2}\Bigr),
\]
then the main claim proved in this section (\pref{thm:covar-estimate}) is
\[
\E_{A,B_t}\norm{\hat Q^{-1}P-\Id_n}_F^2=O_{\beta,t}(1).
\]
The proof is carried out under the planted representation \(x_0\sim G_{A,\beta}\), \(y_t=t x_0+B_t\), and contiguity~\cite[\S4]{el2022sampling} transfers the estimate back to the original tilted model.

The argument has two parts: a decomposition of the Frobenius norm error \(\norm{\hat Q^{-1}P-\Id_n}_F^2\) into functions of overlaps using Gaussian integration by parts and various algebraic simplifications using replica identities, followed by a local cavity interpolation computation identifying the exact values of the non-negligible contributions from certain subsets of those overlap functions.

\subsubsection*{Cavity interpolation for coefficients of overlaps (\pref{sec:cavity-interpolations-trace-identities})} A cavity interpolation isolates one spin, say \(\sigma_n=\varepsilon\), writes \(\sigma=(\rho,\varepsilon)\) with \(\rho\in\{\pm1\}^{n-1}\), and interpolates the cavity field seen by \(\varepsilon\) (that is, the effective local field in the conditional law of \(\varepsilon\) given the other spins) between the true field and an auxiliary Gaussian field with matched mean and variance (see \eqref{eq:smartpath-corrected} in \pref{sec:cavity-interpolation}; cf. \pref{def:cavity-interpolation-planted-sk}). In the planted model, the cavity field contains both an SK contribution and a Curie--Weiss drift, so both pieces must be interpolated simultaneously. At the decoupled endpoint \(s=0\), the last spin factorizes from the bulk and sees the one-dimensional replica-symmetric field (\pref{lem:fact-s0})
\[
Y^*:=y_{t,n}+\beta^2 m^*+\beta\sqrt{q^*}\,Z,
\]
where the $y_t$ is the SL tilt and $Z \sim \calN(0,1)$ is a standard Gaussian. The extra drift term $\beta^2 m^*$ comes from the Curie--Weiss component and is absent in the standard SK cavity method (see \pref{lem:deriv}).
This field goes into the constants in \eqref{eq:rho2rho4} that occur later:
\[
\rho_2:=\E[\sech^2(Y^*)],
\qquad
\rho_4:=\E[\sech^4(Y^*)].
\]

Following Talagrand’s rectangular-sum method from~\cite[\S~1.8]{talagrand2010mean}, the cavity computation is organized around the antisymmetric rectangular sums
\[
f:=R_{1,3}-R_{1,4}-R_{2,3}+R_{2,4},
\qquad
\hat f:=\hat R_{1,3}-\hat R_{1,4}-\hat R_{2,3}+\hat R_{2,4},
\]
where
\[
\hat\sigma_i^\ell:=\frac{\sigma_i^\ell-\langle \sigma_i\rangle}{\sqrt{1-\langle \sigma_i\rangle^2}},
\qquad
\hat R_{\ell,\ell'}:=\frac1n\sum_{i=1}^n \hat\sigma_i^\ell \hat\sigma_i^{\ell'}.
\]
The point of this choice is that \(f\) and \(\hat f\) change sign under the swap \(3\leftrightarrow4\), while replica-symmetric bulk quantities remain invariant; this kills non-replica fluctuations and leaves a single stability denominator \(1-\beta^2\rho_4\) (\pref{lem:orthog-hatR}). Writing \(\nu(\cdot):=\E\an{\cdot}\), the key second moments are
\[
U_n:=\nu(f^2),\qquad V_n:=\nu(f\hat f),\qquad W_n^{(2)}:=\nu(f\,f^{(2)}),
\]
where \(f^{(2)}\) is the \(D^2\)-weighted rectangular sum that later enters \(\Tr(PD^2P)\). The cavity recursion gives
\[
(1-\beta^2\rho_4)U_n=\frac{4\rho_4}{n}+O(n^{-2}),
\qquad
(1-\beta^2\rho_4)V_n=\frac{4\rho_2}{n}+O(n^{-2}),
\]
and
\[
W_n^{(2)}=\frac{4}{n}+\beta^2U_n+O(n^{-2})
=\frac{4}{n(1-\beta^2\rho_4)}+O(n^{-2}).
\]
These formulas are proved in \pref{prop:D17} and \pref{prop:D18} and supply one set of the exact coefficients needed in the final trace expansion. The other critical exact estimate in the cavity interpolation comes in the form of the proof that
\[
    \E[c\Tr[PDP]] = \E[c^2\Tr[P^2]] + \E[c\Tr[P]] + O(1)\, ,
\]
where $c = \frac{\beta^2}{n}\Tr[P]$ (\pref{lem:cTrPDP-identity}). This expansion provides the second set of exact coefficients needed in the final trace expansion, and the fact that this form should exist is ``guessed'' from the behavior of the algorithmic counterpart $\hat{Q}$ to $P$\footnote{~As $\hat{Q}$ is a resolvent evaluated at $c$, one can use analytic subordination to evaluate its derivatives at $c$ and take their projections onto the diagonal sub-algebra $\calD_n$. Doing this reveals that $\hat{Q}$ satisfies a similar identity, and leads to the natural guess that $P$ should as well (if they are good proxies for each other).}. Armed with the exact estimates for $W^{(2)}_n$ and $\E[c\Tr[PDP]]$, one can arrange for exact cancellations between certain terms in $P$ that ``behave like'' the diagonal $D$ in $\hat{Q}$ and another set that comes from applying Gaussian integration-by-parts to ``mixed'' terms that interact with the $\mathsf{GOE}(n)$ matrix $A$.

After (identifying and) arranging for the cancellation of the heavy terms, the task of tremendous complexity that remains is to show every other term is of $O(1)$ (cf.~\pref{lem:remainder-bound-ea}). Doing this is an involved and technical exercise, requiring very subtle book-keeping which simultaneously uses the linear-algebraic properties of $P$ in conjunction with the cancellations in bulk terms (such as $f$) afforded by good concentration and cavity estimates. 

The decisive estimates that control a large set of the remainder terms come from the observation that linear and cubic cavity terms are one order smaller than a naive H\"older bound would suggest. If \(G^-\) is a quadratic bulk observable in the \((n-1)\)-spin cavity system, then uniformly along the cavity path (\pref{lem:mean-bias-On-1} \& \pref{prop:D15}),
\[
\nu_s(f^-G^-)=O(n^{-2}),
\qquad
\nu_s(\hat f^-G^-)=O(n^{-2}).
\]
The Hölder estimate would only give \(O(n^{-3/2})\) (cf.~\pref{lem:holder-bookkeeping}) because \(f^-\) is \(O(n^{-1/2})\) and \(G^-\) is \(O(n^{-1})\). That weaker bound would be fatal: these terms appear with an \(n^2\) prefactor in the final trace expansion, so \(O(n^{-3/2})\) would leave an uncontrolled \(O(n^{1/2})\) remainder. The extra factor \(n^{-1/2}\) comes from a second-order Taylor expansion along the cavity path and the exact antisymmetry of the rectangular sum in \pref{lem:D16-corrected}, while the denominator \(1-\beta^2\rho_4\) being bounded away from zero ensures that the resulting self-consistent estimate closes without losing this cancellation. A similar use of the exact control of the Gibbs derivative available during the cavity interpolation~(\pref{lem:deriv}) sets up a self-consistent linear system to prove that $\nu_s(R_{ab}-q^*)$ and $\nu_s(M_a - m^*)$ are also $O(1/n)$~(cf.~\pref{lem:mean-bias-On-1}) and not just $O(n^{-1/2})$ as suggested by naive concentration estimates.

Just as the requirement that $\beta^2\rho_4 < 1$ is critical in the cubic estimates, so is the requirement that $\rho(\beta^2 \mathbf{M}) < 1$ for the linear fluctuation estimate (where $\mathbf{M} \in \R^{2\times 2}$ is an operator that governs the ``stability'' of a certain linear system). For $\beta < 1/2$, both these conditions are true (cf.~\pref{lem:beta2M-spr-bound}).

\subsubsection*{Covariance error estimate} With cavity asymptotics in hand and concentration for bulk quantities sourced from \cite[\S 3]{alaoui2017finite}, the covariance theorem reduces to an exact trace expansion. One decomposes
\[
\E\norm{\hat Q^{-1}P-\Id_n}_F^2=\mathcal E_D+\mathcal E_A + \calE_R,
\]
where \(\mathcal E_D\) collects the diagonal/TAP terms handled in \pref{lem:diag-frob-term}, \(\mathcal E_A\) collects the terms involving the disorder matrix \(A\) handled in \pref{lem:EA-term}, and \(\mathcal E_R\) collects terms involving the rank-$1$ plant which are handled in \pref{lem:rank-1-O1}. The trace-to-replica identities in \pref{lem:trace-to-replica} and \pref{lem:diag-frob-term}, which include
\[
\Tr[P^2]=\frac{n^2}{4}\an{f^2},
\qquad
\Tr[PDP]=\frac{n^2}{4}\an{f\hat f},
\qquad
\Tr[PD^2P]=\frac{n^2}{4}\an{f\,f^{(2)}},
\]
convert \(\mathcal E_D\) into a linear combination of \(U_n\), \(V_n\), \(W_n^{(2)}\), \(\Tr(P)\), and the rank-one correction coming from \(\frac{2\beta^2}{n}\mg_t\mg_t^\sT\).

The term \(\mathcal E_A\) is treated by Gaussian integration by parts in the GOE entries in \eqref{eq:GOE-ibp} and \eqref{eq:gibbs-deriv-basic}, giving $\E[A_{ij}F] = \frac{1}{n}\E[\partial_{ij}F]$ for any well-behaved function $F$ of $A$, where $\partial_{ij}$ is the partial derivative with respect to $A_{ij}$. The basic derivative identity is
\[
\partial_{ij}\an{F}
=
\beta\Bigl(\an{F\,\sigma_i\sigma_j}-\an{F}\an{\sigma_i\sigma_j}\Bigr).
\]
Applying this to \(P_{ab}\) and expanding the centered four-spin term yields the exact decomposition in \pref{lem:dP-pairing}:
\[
\partial_{ij}P_{ab}
=
\beta\bigl(P_{ai}P_{bj}+P_{aj}P_{bi}\bigr)
+
\beta\bigl(\mg_{t,i}\,T_{abj}+\mg_{t,j}\,T_{abi}+\Gamma_{abij}\bigr),
\]
where \(T_{abj}\) is a centered three-point response term and \(\Gamma_{abij}\) is a centered four-point cumulant. The first bracket produces the leading replica traces; every other term carries an extra centered fluctuation and is, therefore, lower order once the \(O(n^{-2})\) cubic cavity bounds from \pref{prop:D15} and the one-site moment bound $\E[D_{jj}^p]\le C_{p,\beta,t}$ of \pref{lem:D-moments} are utilized. Some of the other terms carry ``repeated indices'' along with factors of $D$ that requires slightly more delicate control, and this often requires the use of a ``diagonally-weighted'' cubic lemma, shown in in~\pref{lem:D-oneleg-O(n)} --  this lemma is also proved via a second-order Taylor expansion, with initial terms controlled via exact decoupling lemmata (cf. \pref{lem:fact-s0},~\pref{lem:rect-scaled} \&~\pref{lem:rect-doublyscaled}) in conjunction with~\pref{lem:D-moments}, and second-derivative terms controlled via symmetry and H\"older arguments.

After substituting the explicit asymptotics of \(U_n\), \(V_n\), and \(W_n^{(2)}\) into \(\mathcal E_D\) and \(\mathcal E_A\), the order-\(n\) contributions cancel exactly. Once the cancellations involving $(U_n,V_n,W^{(2)}_n)$, $\E[c^2\Tr[P^2]]$, and $\beta^2\E[\Tr[P^2]]$ are made, the remaining contributions are all \(O(1)\): the one-site terms $\sigma_i$ weighted by the entries of $D$ are controlled by the uniform moments of \(D_{jj}\), the rank-one magnetization terms are controlled by the overlap and magnetization moment bounds, and the $T_{abj}$ and $\Gamma_{abij}$ terms in \pref{lem:dP-pairing} are controlled because their \(n^2\) prefactors are exactly balanced by the \(O(n^{-2})\) cubic cavity estimates. Consequently,
\[
\E_{A,B_t}\norm{\hat Q^{-1}P-\Id_n}_F^2=O_{\beta,t}(1).
\]

\paragraph{Properties of the algorithmic process (\pref{sec:alg-properties})}
The analysis in \pref{sec:alg-properties} reduces to one random-matrix input from \pref{app:deformed-wigner-resolvent}: a uniform law for the diagonal of a real-axis squared resolvent.  Once that input is available, the rest of the argument is deterministic resolvent calculus.  Throughout, \(\gamma:=(1-2\beta)/2>0\) is fixed, but to keep the roadmap readable we display only the powers of \(n\); all suppressed constants depend polynomially on \(\beta\) and \(\gamma^{-1}\).  

Note that in keeping with free probability conventions, in these two sections we use the non-commutative $L^2$ norm $\schnorm{\cdot}$ on linear operators, which differs from the Frobenius norm $\norm{\cdot}_F$ by a factor of $1/\sqrt{n}$ when the linear operator is a matrix. We also use the normalized trace $\tr_n$.

\subsubsection*{Deterministic resolvent control from a spectral gap}

The deterministic part of the argument is developed in \pref{sec:alg-basic-bounds-lipschitz}.  For \(m\in(-1,1)^n\), the relevant coefficients are
\[
D(m)=(\Id_n-\diag(m)^2)^{-1},
\qquad
a(m)=\beta^2\tr_n(D(m)^{-1}),
\qquad
\bar Q(m)=\bigl(a(m)\Id_n-\beta A+D(m)\bigr)^{-1},
\]
together with the rank-one-corrected resolvent \(\hat Q(m)\) appearing in ASL--TAP.  The basic high-probability event is
\[
\Omega_\gamma^{(A)}:=\{\beta\opnorm{A}\le 1-\gamma\}.
\]
On this event, every matrix \(a\Id_n-\beta A+D\) with \(a\in[0,1]\) and \(D\succeq \Id_n\) is uniformly invertible, so \pref{lem:alg-resolvent-bounds} and \pref{lem:alg-resolvent-lipschitz} give
\[
\opnorm{\hat Q(m)}=O(1),\qquad
\opnorm{D(m)\hat Q(m)}=O(1),\qquad
\opnorm{D(m)\hat Q(m)^2D(m)}=O(1),
\]
together with Lipschitz bounds (as measured with $\schnorm{\cdot}$)
\[
\operatorname{Lip}(\hat Q)=O(n^{-1/2}),\qquad
\operatorname{Lip}(\hat Q^2)=O(n^{-1/2}),\qquad
\operatorname{Lip}(D\hat Q^2D)=O(n^{-1/2}).
\]
Thus, after \(\Omega_\gamma^{(A)}\) is imposed, the only nontrivial quantity left is the diagonal of \(\hat Q(m)^2\), encoded by
\[
\delta_{\diag}(m)
=
E_{\calD_n}\!\left[\hat Q(m)^2\right]D(m)^2
-
\Bigl(1+\beta^2\tr_n(\hat Q(m)^2)\Bigr)\Id_n
\]
from \pref{eq:delta-diag}.  This is the term that enters the It\^o correction in ASL--TAP and controls its drift mismatch with PHD (observe that $\delta_{\diag}(m)$ differs from $m - f(m)$ by a small $\frac{2\beta^2}{n}\hat{Q}^2(m)m$). 

\subsubsection*{Free-probability prediction for the diagonal}

The target for the diagonal law is explained in \pref{sec:free-prob}. 
Essentially, free probability provides a way to characterize sums and products with GOE matrices by a limiting \emph{semicircular operator} $S$ as $n \to \infty$.

The operator \(S\) is an infinite-dimensional object serving as the limiting case of $A$, and (following Biane~\cite{biane1997free}) analytic subordination in \pref{lem:subordination-semicircular} produces an analytic map \(F\) such that
\[
E_{\calD_n}\!\left[(z\Id_n-\beta S + D)^{-1}\right]=(F(z)\Id_n + D)^{-1}.
\]
Differentiating this identity at a real point \(a\) predicts
\[
E_{\calD_n}\!\left[D(a\Id_n-\beta S + D)^{-2}D\right]
=
F'(a)\,D(F(a)\Id_n+D)^{-2}D.
\]
At the value used in the algorithm $a=\beta^2\tr_n(D^{-1})$, \pref{lem:subordination-semicircular} gives
\[
F(a)=0,
\qquad
F'(a)=\frac{1}{1-\beta^2\tr_n(D^{-2})},
\]
so the free model collapses to the scalar matrix
\[
E_{\calD_n}\!\left[D(a\Id_n-\beta S + D)^{-2}D\right]
=
\frac{\Id_n}{1-\beta^2\tr_n(D^{-2})}.
\]
This is the only consequence of free probability later used in \pref{sec:alg-properties}: after this specialization, the main task is simply to show that the finite-\(n\) diagonal is uniformly close to this scalar matrix. 

\subsubsection*{Finite-\(n\) comparison with the freely independent limit}
The main finite-\(n\) statement is \pref{cor:diagonal-Q2-controlled-in-l2}.  Writing
\[
Y(a,D):=E_{\calD_n}\!\left[D(a\Id_n-\beta A + D)^{-2}D\right],
\]
that result states that with high probability,
\begin{equation}
\label{eq:free-prob-diagonal-goal}
\sup_{a\in[0,1]}\sup_{D\succeq \Id_n}
\schnorm{Y(a,D)-F'(a)D(F(a)\Id_n+D)^{-2}D}
\le
O(n^{-1/2}).
\end{equation}
Specializing to \(a=\beta^2\tr_n(D^{-1})\) produces the event \(\Omega^{(\mathrm{free})}\) used in \pref{sec:alg-diagonal-control}.  The proof combines the concentration result in \pref{thm:chaining-square-no-sqrtlog} with the expectation comparison of \pref{prop:gamma-regime-square-expectation-estimate}.

The expectation comparison follows the same overall template as the modified-Hessian analysis in~\cite[Section 4]{jekel2024pha}\iffocs{}{ and~\cite[Section 4]{jekel2025pha2}}, but the observable is different.  Here one needs the real-axis squared
\[
D(a\Id_n-\beta A + D)^{-2}D
\]
rather than a regularized imaginary-part resolvent.
We approach the squared resolvent through the limit
\[
-\frac1b\,\Im\bigl[D(a\Id_n +ib\Id_n-\beta A + D)^{-1}D\bigr]\;\to\; D(a\Id_n - \beta A + D)^{-2}D\qquad (b\downarrow0),
\]
To control the diagonal of this matrix, we perform Gaussian interpolation on the trace inner product between $DR_k(a+bi,t)D$ and an arbitrary diagonal matrix $T$, where
\[
W_k(t)=\sqrt{1-t}\,(A\otimes \Id_k)+\sqrt t\,B_k,
\qquad
R_k(z,t)=\bigl(z\Id_{nk}-\beta W_k(t)+D\otimes \Id_k\bigr)^{-1}.
\]
We take $k \to \infty$ to approach the freely independent limit at $t=1$, whereas at $t=0$ it stays equal to the quantity of interest $\iprod{T,\, Y(a,D)}_{\tr_n}$.
As in \iffocs{\cite{jekel2024pha}}{\cite{jekel2024pha,jekel2025pha2}}, we bound the $t$-derivative of this interpolated observable by applying Gaussian integration by parts, after which we see explicit single-trace terms of order \(n^{-1}\) and also a covariance term that is shown to be of order \(n^{-3/2}\) by a Poincar\'e inequality.  Thus the expectation-vs.-free-model comparison is strictly smaller than the final fluctuation scale:
\[
\Iprod{T,\;\E_A[Y(a,D)]-F'(a)D(F(a)\Id_n+D)^{-2}D}_{\tr_n}
\;\le\;
O(n^{-1})\schnorm{T}+O(n^{-3/2})\schnorm{T}.
\]

Unlike in \iffocs{\cite{jekel2024pha}}{\cite{jekel2024pha,jekel2025pha2}}, this interpolation encounters a singularity if $a\Id_n - \beta A + D \not\succ 0$.
We therefore multiply the interpolated observable by a smooth cutoff function $\chi_k(A,B_k)$ to make it identically 0 in the rare event when $\opnorm{W_k(t)}$ can be large enough to potentially cause a singularity.
This introduces an exponentially small error, though it requires a free probability characterization of the edges of the spectrum of $W_k(t)$ in \pref{cor:AIk+Bk-final-bound}.

\subsubsection*{Chaining bound over all $D$}
The preceding argument controlled $\E_A \iprod{Y(a,D),T}$ for diagonal test matrices $T$, but to obtain \pref{eq:free-prob-diagonal-goal}, we must obtain a high-probability uniform bound over all choices of $D$ and $T$.
Writing \(D=W^{-1}\), we study the scalar process
\[
\Phi_{a,W,T}(A)
=
\iprod{T,\, Y(a,D)}_{\tr_n},
\qquad
\schnorm{T}=1.
\]
A direct computation, proved in \pref{lem:mixed-lipschitz-gradient}, gives
\[
\nabla_A\Phi_{a,W,T}(A)
=
\beta\Bigl(
X^{-2}DTD\,X^{-1}
+
X^{-1}DTD\,X^{-2}
\Bigr),
\qquad
X=a\Id_n-\beta A + D.
\]
This shows that \(\Phi_{a,W,T}\) is Lipschitz in \(A\) and also Lipschitz in the parameter triple \((a,W,T)\), with constants independent of \(n\).  The parameter family has effective dimension \(2n+1\): one scalar parameter \(a\), \(n\) diagonal degrees of freedom in \(W\), and \(n-1\) in the diagonal test matrix \(T\).  Pointwise GOE concentration is already of order \(n^{-1}\), but a naive \(\varepsilon\)-net over this family (as used in \cite{jekel2024pha}) would have size \(\exp(\Theta(n\log n))\), introducing an extra \(\sqrt{\log n}\) and yielding only a bound of \(O(\sqrt{\log n}/\sqrt n)\).  Generic chaining (in the form of Dudley-style entropy bounds) removes this logarithmic loss and upgrades the pointwise \(n^{-1}\) concentration to the sharp uniform \(n^{-1/2}\) bound in \pref{thm:chaining-square-no-sqrtlog}.  

\subsubsection*{The TAP resolvent and consequences for ASL--TAP}
The argument so far applies to \(\bar Q(m)\), which omits the rank-1 term in $\hat{Q}(m)$.  The transfer to the actual TAP resolvent is isolated in \pref{sec:alg-diagonal-error-rank-1}: \pref{lem:alg-resolvent-bounds-DQ2-hQ2D} shows that
\[
\schnorm{
E_{\calD_n}\!\left[\hat Q(m)^2-\bar Q(m)^2\right]D(m)^2
}
\le 
O(n^{-1/2}).
\]
Hence the same diagonal law holds for \(\hat Q(m)\) at the same scale.  Combined with \pref{lem:alg-diag-basic} and \pref{lem:alg-resolvent-bounds-delta-frob}, this yields
\[
\sup_{m\in(-1,1)^n}\schnorm{\delta_{\diag}(m)}\le O(n^{-1/2}),
\qquad
\operatorname{Lip}(\delta_{\diag}) \le O(n^{-1/2}).
\]

The drift comparison in \pref{sec:alg-phd-asl-tap-closeness} is then deterministic.  If \(f(m)\) denotes the It\^o correction map from that section, \pref{lem:alg-resolvent-bounds-m-f} and \pref{lem:alg-resolvent-lipschitz-drift} imply
\[
\lpnorm{m-f(m)} \le O(1),
\qquad
\operatorname{Lip}\bigl(m\mapsto \hat Q(m)(m-f(m))\bigr) \le O(1).
\]
Together with the \(O(1)\) operator bounds and \(O(n^{-1/2})\) Lipschitz bounds on \(\hat Q\), this yields existence and uniqueness of the SDE, Wasserstein-2 robustness to path perturbations, nonescape from \((-1,1)^n\) in \pref{lem:stay-in-cube}, and sub-Gaussian concentration of Lipschitz path functionals (including the Jarzynski weights) in \pref{thm:mSDE-path-subgaussian}.

\paragraph{Entropy contraction for Ising models on wedges (\pref{s:localized})} 
After running \eqref{e:asl-tap-je} up to large constant time $T$, we wish to sample from the localized distribution $\mu_T = \mu_{\be A, y_T}$ using the polarized walk, which is enabled by the main theorem of this section (\Cref{t:ls-sk}):
\begin{enumerate}
\item For the stochastic localization process \(y_T=T\si^*+B_T\), for large enough constant \(T\), the posterior \(\mu_{\be A,y_T}\) is exponentially concentrated on a small Hamming ball (wedge) $\ball{\ep n}{\si^*}$ or $\ball{\ep n}{\sign(y_T)}$.
    \item If \(\Omega=\ball{\ep n}{x_0}\) for sufficiently small $\ep$, then with high probability over \(A\sim \GOE(n)\), the restricted measure $\mu_{\be A,h}\big|_{\Omega}$
satisfies 
entropy contraction with constant $\ge \fc{\rh_\be}{n}$, 
with
$\rh_\be$ depending only on $\be$.
\end{enumerate}
The intuition is that for large $T$, $y_T = T\si^* + B_T$ 
has a large signal to noise ratio, so that the mass of the posterior $\mu_T$ is concentrated around $\si^*$. We know that for any Ising model $\mu_{\be A}$ on the hypercube, a log-Sobolev inequality is satisfied when $\lm_{\max}(\be A)-\lm_{\min}(\be A)<1$ \cite{anari2022entropic}. One way to think of this is $\lm_{\max}(\be A)-\lm_{\min}(\be A)<\rc{2\sqrt n} \diam(\{\pm 1\}^n)$. One might expect that if we restrict to a small $\ep n$-wedge, which has smaller diameter, that this allows sampling when $A$ has larger spectral width. The condition we will use to ensure this is that $A_{S\times S}$ is small for sets of coordinates of size $O(\ep n)$, which is satisfied whp when we sample $A\sim \GOE(n)$ as for the SK model.

Huang, Montanari and Pham~\cite{huang2024sampling} prove a functional inequality for the localized distribution for the spherical $p$-spin model by showing log-concavity in a chart; since our distribution is on the hypercube we have to work harder to tame the geometry. 

\ppart{Two-stage decomposition}
We make a two-stage decomposition into product measures restricted to the wedge $\Om$, and then make a key estimate of the covariance of these measures. The two stages are as follows.
\begin{enumerate}
    \item The needle decomposition of \cite{eldan2022spectral} (\Cref{t:needle}), given a test function $\ph$, decomposes an Ising model $\mu_{A}$ (where $A$ is PSD) into rank-1 Ising models \(\mu_{xx^{\sT},v}\) with $xx^{\sT}\preceq A$, which preserve the expectation under $\ph$. Therefore, it suffices to prove functional inequalities for the rank-1 models which arise.

    \noindent Since this requires a PSD interaction matrix, we replace \(\be A\) by \(\be A+\ga I\). This does not change the Gibbs measure on the hypercube, because \(\an{\si,I_n\si}=n\) is constant.
    \item The Hubbard--Stratonovich transform (\Cref{t:hs}) writes each such rank-1 model as a 1-dimensional mixture of product measures,
\[
\mu_{\lm xx^{\sT},v}\big|_{\Om}
=
\int p_\lm(u)\,\mu_{u\sqrt{\lm}\,x+v}\big|_{\Om}\,du .
\]
\end{enumerate}
By the localization framework \cite{CE25} applied to the rank-1 models, it suffices to control the covariance of $\mu_{\lm xx^{\sT},v}\big|_{\Om}$. 
By variance decomposition, since the variance of $\mu_{u\sqrt{\lm}\,x+v}\big|_{\Om}$ is easily controlled, it suffices to control the variance of the 1-dimensional distribution $p_\lm$. We summarize the two-stage decomposition in \pref{lem:LS-2-stage}. 

This would follow if $p_\lm$ were log-concave, so we calculate
\[
-(\log p_\lm)''(u)
=
1-\Var_{\mu_{u\sqrt{\lm}x+v}\mid_\Om}\bigl(\an{\sqrt{\lm}x,\si}\bigr).
\]

\ppart{Entropy contraction for $\be<\rc{2\sqrt 2}$ (\pref{s:lsi-1})} We need to bound the covariance of a restricted product measure. For convenience, we bound the covariance when restricted to $\Om$, a slice or wedge of size $k$ in $\{0,1\}^n$. The basic bound, \Cref{l:prod-dist-cov}, is that it can be bounded in terms of its diagonal entries,
\begin{align}\label{e:cov-2-diag}
\Cov(\mu_h|_\Omega)&\preceq 2\diag\bigl(p_i(1-p_i)\bigr),
\qquad p_i=\E_{\mu_h|_\Omega}\si_i,
\end{align}
proved using negative dependence and a Gershgorin argument. 
This gives the key inequality \eqref{e:var-by-AS} for $\Om=\ball{k}{x_0}$, $xx^{\sT}\preceq A$,
\[
\Var_{\mu_h|_{\Om}}(\ip{x}{\si}) \le 2 \max_{|S|=4k} \opnorm{A_{S\times S}}^2.
\]

The gain from restricting to a wedge is that only \(O(\ep n)\) coordinates can fluctuate. Accordingly, the relevant interaction strength is not the full norm of \(A\), but the norms of small principal submatrices. A union bound together with the GOE operator-norm tail estimate shows that with high probability,
\[
\sup_{|S|\le 4\ep n}\opnorm{A_{S\times S}}
\lesssim \sqrt{h(4\ep)},
\]
where \(h\) is the binary entropy; see \Cref{l:AS}. 

Now recall that we actually replace $A$ by $\be A + \ga I$ to make it PSD; we must take $\ga>2\be$. The contribution of $\be A$ to the operator norm of small submatrices is small, but there is an unavoidable contribution of $\ga>2\be$. In the worst case, the covariance of rank-1 Ising models can be bounded only when the norm is at most 1 (consider the Curie--Weiss model), so this places a limit of $\be<\rc2$. Essentially, considering small submatrices can squash the positive side of the spectrum of $A$, enabling an improvement of a factor of 2 from \cite{eldan2022spectral}.

However, the covariance bound \eqref{e:cov-2-diag} is a factor of 2 loose from what it would be for an unrestricted product distribution, so this only yields entropy contraction for $\be< \rc{2\sqrt 2}$ (\pref{l:ls-sk-wedge}). This is exactly the obstruction to attaining $\be<\rc 2$.

\ppart{Attaining $\be<\rc{2}$ using a refined covariance bound (\pref{s:lsi-2})} Although \eqref{e:cov-2-diag} is not true without the factor 2, we note that the covariance is indeed smaller when the total uncertainty $m=\sum_i p_i(1-p_i)$ between the coordinates is large: \pref{cor:cov-hi-under-uncert} gives that 2 can be replaced by $1+\de$ when $m\ge m(\de)$. This is a corollary of the main result (\pref{l:cov-tight-off-diag}) that as $m\to \iy$,
the off-diagonal covariances are asymptotically rank-one: 
\[
\Cov_{\mu_h|_\Omega}(\si_i,\si_j)
=
-\frac{p_i(1-p_i)p_j(1-p_j)}{m}(1+o(1))
\]
for $\Om$ a slice of $\{0,1\}^n$, and the same formula holds for wedges up to an additional factor \(c\in[0,1]\). 

The proof proceeds through an asymptotic cross-ratio formula for slice probabilities (\pref{l:cross-ratio}) to obtain the asymptotics for the slice. The computation for the wedge is done through conditioning on slices and summing, using a local limit theorem (\pref{t:llt}) to control slice probabilities and carefully restricting to windows of size $\om(\sqrt m)$.

One still has to handle the complementary regime where \(m\) is small. In this case, only \(O(1)\) coordinates can have marginals bounded away from \(0\) and \(1\). The structural \Cref{l:exist-small-uncertain} shows that if the variance of \(\an{\si,x}\) is still too large, then there must exist coordinates \(i,j\) with \(x_i<0<x_j\) whose marginals are both non-extremal. Under the one-dimensional tilt \(h\mapsto h+ux\), however, such a pair can remain simultaneously non-extremal only for \(u\) in an interval of length \(O(1)\); see \Cref{l:var-large-exist-sets}. Consequently, the HS density \(p_\lm\) is strongly log-concave except on a set of \(u\)'s of bounded measure, while everywhere it still has a crude lower curvature bound. The 1-dimensional \Cref{l:mostly-lc-var} then implies that \(p_\lm\) has bounded variance, which is exactly what the two-stage decomposition needs.

\ppart{Concentration of localized measure (\pref{s:hi-prob-wedge})}
The remaining ingredient is the concentration of the stochastic-localization posterior. If \(y_t=tx_0+B_t\), then each coordinate of \(\sign(y_t)\) differs from the corresponding coordinate of \(x_0\) with probability \(\Phi(-\sqrt t)\), independently across coordinates. 

Therefore, using a Chernoff bound, taking 
\(T=O(\log(1/\ep))\), and interpreting the localized measure as a posterior, we obtain concentration on a wedge of radius \(\ep n\) (\pref{l:conc-x0}), and choosing \(\ep=\ep(\be)\) small enough to satisfy the deterministic wedge condition completes the proof of \Cref{t:ls-sk}.

Finally, once the restricted measure on the wedge has entropy contraction, efficient sampling and partition function estimation follow by standard arguments (\pref{s:samp-Z}).

\section{Combining ASL, PHD, and JE: A quantitative analysis}
\label{s:sde}

We first introduce the many SDEs that we consider (\pref{s:asl-phd-tap}) and Jarzynski's equality (\pref{s:je}). Then we introduce the Wasserstein error bound using Gr\"onwall (\pref{s:wasserstein}), which we apply to bounding the error in the estimated magnetization in ASL. Finally, we give the guarantees of Jarzynski's equality (\pref{s:je-error}). This will give us sampling guarantees under our desiderata, provided we can also sample from localized distributions after some constant time.

Because we often apply It\^o's lemma to functions that are only defined on $(-1,1)^n$, we use the following local version of It\^o's lemma along with the almost surely infinite escape time of ASL-TAP (\pref{lem:stay-in-cube}).
\begin{lemma}[{Local It\^o's lemma, see e.g. \cite[Theorem 3.4.3]{lawler2010stochastic}}]
\label{lem:ito-local}
Let $\mathbb{D} \subseteq \R^n$ be an open set and let $t_0 < t_1$ be real numbers.
    Suppose the It\^o process $X_t \in \R^n$ satisfies $X_{t_0} \in \mathbb{D}$ and also satisfies the SDE 
    \[dX_t = \mu_tdt +\sigma_tdB_t\]
    with $\mu_t \in \R^n$ and $\sigma_t \in \R^{n \times n}$ when $t \in (t_0, t_1)$ and $X_t \in \mathbb{D}$.
    Suppose that $f_t(x)$ is a vector-valued function that is continuously differentiable in $t \in [t_0, t_1]$ and continuously twice-differentiable in $x \in \mathbb{D}$.
    Let $T := \inf\{t \in (t_0, t_1):  X_t \not\in \mathbb{D}\}$.
    Then for $t \in (t_0, T)$, $f_t(X_t)$ is an It\^o process satisfying the SDE
    \[ df_t(X_t) = [\partial_tf_t](X_t)dt + \nabla_xf_t(X_t)^{\sT}\mu_tdt  + \nabla_xf_t(X_t)^{\sT}\sigma_tdB_t + \frac{1}{2}\nabla_x^2f_t(X_t)\pmb{:} \sigma_t\sigma_t^{\sT} dt, \]
    where $\pmb{:}$ is the tensor contraction $(A \pmb{:} M)_k := \sum_{i,j} A_{i,j,k}M_{i,j}$ and the first two coordinates of $\nabla_x^2f_t(X_t)$ are the ones corresponding to the $\nabla$s.
\end{lemma}

\subsection{ASL, PHD and TAP}
\label{s:asl-phd-tap}
First, we recall the ideal processes stochastic localization (SL) and Hessian dynamics (HD) where we assume a fixed interaction matrix $\beta A$ throughout, so the subscripts on $\mu$ and $\an\cdot$ refer to the tilt (external field).
\begin{align}
    \repeatequation{e:SL},\\
    \repeatequation{e:HD},\\
    \mg_t &= \an{\si}_{y_t}, & Q_t &= \Cov(\mu_{y_t}).\label{e:mg}
\end{align}
\begin{lemma}
    Defining stochastic localization by \eqref{e:SL} and \eqref{e:mg}, the magnetization $\mg_t$ evolves according to \eqref{e:HD}.
\end{lemma}
\begin{prf}
    By Itô's formula, 
    \[
d\mg_t = f(y_t) dt + D\mg(y_t)\, dB_t
    \]
    for some function $f$. Note that $D\mg(y_t) = \Cov(\mu_{y_t}) = Q_t$ as required. 
    Rather than compute $f$ explicitly, we note that $\mg_t$ is a martingale so $f=0$. 
\end{prf}
Now we consider the algorithmic version of these equations. 
Recall 
\begin{equation} 
\repeatequation{e:hQ},
\end{equation}
and consider
\begin{align}
    \repeatequation{e:ASL},\\
    \repeatequation{e:PHD},\\
    \repeatequation{e:TAP},
\end{align}
where $\FT$ is given by \eqref{e:FTAP}. Expanding the last equation gives (with $\mg=\hat\mg_t$, $y=\hy_t$)
\begin{align}
\label{e:TAP-expanded}
    0 &= 
    -\be A m + y + \ub{\rc 2\log \pf{1+\mg}{1-\mg}}{\tanh^{-1}(\mg)} 
    + \be^2 \pa{1-\rc n \ve{m}_2^2}m,
\end{align}
where operations on $m$ are defined componentwise.
As explained before, these equations cannot all be satisfied at the same time but we can choose a pair to satisfy. If we choose \eqref{e:TAP} and one of the SDEs, we can write the evolution for the other variable by using Itô's lemma. 
Recall the definition
\begin{align}
    \repeatequation{e:Itomag}
\end{align}
where $D(m) = \diag\prc{1-m^2}$. 
Then the two options give 
\begin{align}
\savetagequation{e:ASL-TAP}{ASL-TAP}{&\begin{cases}\begin{aligned}
        d\hat y_t &= \hat \mg_tdt+dB_t\\
        d\hat \mg_t &= \hQ(\hat \mg_t)\pa{\hat \mg_t - f(\hat \mg_t)} dt + \hQ(\hat \mg_t)dB_t
    \end{aligned}\end{cases}}\\
\savetagequation{e:PHD-TAP}{PHD-TAP}{&\begin{cases}\begin{aligned}
        d\hat \mg_t &= \hQ(\hat \mg_t)dB_t\\
        d\hat y_t &= f(\hat \mg_t)dt + dB_t 
    \end{aligned}\end{cases}}
\end{align}
which we derive below.

We also note, for the purposes of solving \eqref{e:TAP},
the Hessian of $\FT$:
\begin{align}
    \label{e:Hess-TAP}
    \gd_{\mg}^2 \FT(\mg,y)
    &= -\be A + D(\mg) + \be^2\pa{1-\fc{\ve{m}_2^2}{n}}\Id_n - \fc{2\be^2}{n} \mg\mg^{\sT} = 
    \hQ(\mg)^{-1}
\end{align}
By \ref{d:Q-reg}, $\gd_{\mg}^2 \FT(\mg,y)= \hQ(\mg)^{-1}\succeq \La^{-1}D(\mg) \succeq \La^{-1}\Id_n$ is strongly convex for $\be<\rc 2$. Thus, if \eqref{e:TAP} has a solution, then it is unique; moreover, as long as $\FT(m,y)$ attains minimum in $(-1,1)^n$, the minimum is the solution to \eqref{e:TAP}. Noting $(\gd\FT)_i\to \pm \iy$ as $x_i\to \pm 1$, we conclude the minimum must be attained in $(-1,1)^n$ and equal the unique solution to \eqref{e:TAP}. In \pref{s:mirror-tap}, we show how the solution can be efficiently found using the fixed point iteration of mirror descent. 

\begin{remark}
    Standard results in stochastic analysis imply that the Lipschitz characterization of the dual SDE coefficients in \pref{l:dual-phd} and \pref{lem:stay-in-cube} also suffices to estimate $\hat{m}_t$ by an Euler-Maruyama discretization of \eqref{eq:d-ASL-TAP} instead of by mirror descent.
    We have opted for mirror descent due to the convenience of geometric convergence.
\end{remark}

\paragraph{Derivation of ASL/PHA-TAP} 

Differentiating \eqref{e:TAP} with Itô's formula gives, using \eqref{e:Hess-TAP},
\begin{align}
    dy 
    &= 
    \hat Q(\mg)d\mg 
    + \rc 2 \sumo in
    \Tr\ba{\pa{\gd_\mg^2 \pa{\tanh^{-1}(\mg_i) + \be^2 \pa{1-\fc{\ve{\mg}^2}{n}}\mg_i}}\hat Q(\mg)^2}
    e_i dt.
\end{align}
Now 
\begin{align*}
    \gd_{\mg}^2 \pa{\tanh^{-1}(\mg_i) + \be^2 \pa{1-\fc{\ve{\mg}^2}{n}}\mg_i}
    &= \fc{2\mg_i}{(1-\mg_i^2)^2}e_i^{\ot 2} - \fc{2\be^2\mg_i}{n}\Id_n - \fc{4\be^2\mg_i}{n}e_i^{\ot 2} 
    - \sum_{j\ne i}\fc{2\be^2 \mg_j}{n}(e_i\ot e_j+e_j\ot e_i)
\end{align*}
so the Itô term is
\begin{align*}
    \sumo in \fc{\mg_i}{(1-\mg_i^2)^2}(\hat Q(\mg)^2)_{ii} e_i
    -\fc{\be^2}n \sumo jn \ve{\hat Q_j(\mg)}^2\mg 
    - \fc{2\be^2}n\hat Q(\mg)^2\mg.
\end{align*}
This is exactly $f(\mg)$. Hence enforcing \eqref{e:TAP} during the time evolution gives
\begin{align*}
    d\hat y_t &= \hat Q(\hat \mg_t)^{-1} d\hat \mg_t + f(\hat \mg_t) dt.
\end{align*}
Assuming \eqref{e:PHD}, we obtain \eqref{e:PHD-TAP}. Assuming \eqref{e:ASL}, we get
\[
\hat \mg_t dt + dB_t = \hQ^{-1}(\hat \mg_t)d\hat \mg_t + f(\hat \mg_t)dt,
\]
and solving for $d\hat \mg_t$ gives \eqref{e:ASL-TAP}. 

\subsection{Jarzynski's equality}
\label{s:je}

Recall that we use Jarzynski's equality with either \eqref{e:ASL-TAP} or \eqref{e:PHD-TAP} to maintain importance weights so that the distribution of $y_t$ follows $\rho_t$ satisfying
\begin{align*}
\repeatequation{e:rho-t}.
\end{align*}
This is a natural choice of scaffolding distribution because letting $p_t = \dist(y_t)$ where $y_t$ is as in \eqref{e:SL}, noting that $y_t = t\si + \sqrt t\xi$ for $\si\sim \mu_{\be A}$, we expect $e^{\FT(\hat \mg (y), y)}$ to be an approximation of $Z(\be A,y)$ so that
\allowdisplaybreaks
\begin{align*}
    p_t(y) &\propto \sum_{\si\in \{\pm 1\}^n} e^{\rc 2 \ip{\si}{\be A\si}}e^{-\fc{\ve{y-t\si}_2^2}{2t}}= \sum_{\si\in \{\pm1\}^n} e^{\rc 2 \ip{\si}{\be A\si} + \an{y,\si}} e^{-\fc{\ve{y}_2^2}{2t}}\\
    &= Z(\be A,y) e^{-\fc{\ve{y}_2^2}{2t}} 
    \approx e^{\FT(\hat \mg(y),y)-\fc{\ve{y}^2}{2t}} \propto \rh_t(y).
\end{align*}
By \Cref{t:je}, for any $\rh_t$-integrable function $h$ we will have
$
\int_{\R^n} h(x)\rh_t(x)dx = \fc{\E[e^{w_t}h(y_t)]}{\E[e^{w_t}]}
$
in both cases.\footnote{Note that $\rh_0$ is singular so \Cref{t:je} does not directly apply. See the end of this subsection for the formal limiting argument.} For ASL, we claimed the resulting system is
\allowdisplaybreaks
\begin{align}
\savetagequation{e:asl-tap-je}{ASL-TAP-JE}{
    &\begin{cases}\begin{aligned}
        d\hat y_t &= \hat\mg_tdt+dB_t\\
        d\hat \mg_t &= \hQ(\hat\mg_t)\pa{\hat\mg_t - f(\hat\mg_t)} dt + \hQ(\hat\mg_t)dB_t\\
        d\hat w_t &= \rc 2 \ba{\Tr(\hQ(\hat\mg_t))+ \ve{\hat\mg_t}^2}dt.
    \end{aligned}\end{cases}}
\end{align}
We will let $\om(\mg)= \rc 2 \ba{\Tr(\hQ(\mg))+ \ve{\mg}^2}$, so we can write the last equation as $d\hat w_t = \om(\hat\mg_t)dt$. 
For PHD, write $f_{\text{err}}(\mg)= f(\mg)-\mg$; the resulting system is
\allowdisplaybreaks
\begin{align}
\savetagequation{e:phd-tap-je}{PHD-TAP-JE}{
    &\begin{cases}\begin{aligned}
        d\hat \mg_t &= \hQ(\hat\mg_t)dB_t\\
        d\hat y_t &= f(\hat\mg_t) dt + dB_t \\
        d\hat w_t &= \bc{\rc 2 \ba{\Tr(\hQ(\hat\mg_t)) + \ve{\hat\mg_t}^2}
        + \gd_y \cdot f_{\text{err}}(\hat\mg(y_t))
        + \an{\hat\mg_t - \fc{\hat y_t}t, f_{\text{err}}(\hat\mg_t)}}
        dt.
    \end{aligned}\end{cases}}
\end{align}
We derive these systems of equations below.

\paragraph{Derivation of ASL/PHD-TAP-JE}
Note
\[
\FT(\hat\mg(y),y) = 
\inf_{\mg} \an{y,\mg} - g(\mg) 
=: g_*(y)
\quad \text{where}\quad 
g(\mg) = \fc{\be}2\an{\mg, A\mg} + \sumo in h(\mg_i) + \fc{\be^2n}{4}\pa{1-\rc n\ve{\mg}^2}^2,
\]
where we denote $\hat \mg(y)$ as the argmin in the optimization problem, i.e., that solves $y=\gd g(\mg)$, so $\hat \mg(y)=(\gd g)^{-1}(y)$.
Then because $\hat \mg(y) = (\gd g)^{-1}(y)$, 
\begin{align}
\nonumber
    \gd_y (\FT(\hat\mg(y),y)) 
    &= \gd_y 
    g_*(y)
    = \gd_y (\an{y,\hat\mg(y)} - g(\hat\mg(y))) \\
    &= \hat\mg(y) + D_y\hat \mg(y) \cdot  (y-\gd g(\hat\mg(y))) = \hat\mg(y)
\label{e:gd-y-FTAP}
\end{align}
and
\begin{align*}
    \gd_y^2 (\FT(\hat\mg(y),y)) 
    &= D_y ((\gd_y g)^{-1}(y))
    = (D_y (\gd_y g)(\hat\mg(y)))^{-1}
    = \gd_y^2 g(\hat\mg(y))^{-1} = D_y\hat\mg(y) = \hQ(y)\\
    \De_y (\FT(\hat\mg(y),y))  
    &= \gd\cdot \hat\mg(y) = \Tr[\hQ(\hat\mg(y))].
\end{align*}
We now compute the $dw_t$ equation in 
\Cref{t:je} for $\rh_t\propto e^{-U_t}$ where 
$-U_t(y)=\FT(\hat\mg(y),y) - \fc{\ve{y}^2}{2t}$. 
We calculate
\begin{align*}
    -\gd U_t(y) &= \hat\mg(y) - \fc{y}{t}
\end{align*}
For \eqref{e:ASL-TAP}, we take $b_t(y) = \rc2 \pa{\hat\mg(y) + \fc yt}$ and compute
\begin{align}
\label{e:dw-calc-1}
    \gd \cdot b_t(y) - \an{\gd U_t(y), b_t(y)} - \pl_t U_t(y)
    &= \ba{\rc2 \Tr(\hQ(\hat \mg(y))) + \fc n{2t}}
    + \an{\hat\mg(y)-\fc yt, \rc2\pa{\hat\mg(y) + \fc{y}{t}}}
    + \fc{\ve{y}^2}{2t^2}.
\end{align}
Thus (dropping the part of the expression that doesn't depend on $y$) we obtain that we can take
\[dw_t = \rc2\ba{\Tr(\hQ(\hat\mg_t))+ \ve{\hat\mg_t}^2}dt,\]
giving \eqref{e:asl-tap-je}.

For \eqref{e:phd-tap-je}, we note that to match \Cref{t:je}, we take
\[
b_t(y) = \rc2 \pa{\hat\mg(y) + \fc yt} + \pa{f(\hat\mg(y))-\hat \mg(y)}.
\]
The last term will give rise to an extra term in $\gd \cdot b_t(y)$ and $\an{\gd U_t(y),b_t(y)}$, showing \eqref{e:phd-tap-je}.

\paragraph{Limiting argument} 
Although we have existence and uniqueness of \eqref{e:asl-tap-je} for all $t\ge 0$, formally, we can only apply \Cref{t:je} on $[s,t]$ for $s>0$. 
Let $h$ be bounded and measurable. 
If we initialize \eqref{e:asl-tap-je} with $\hx_s\sim \rh_s$, $\hw_s=0$, then \Cref{t:je} gives us
\begin{align}
\label{e:je-limit-1}
    \fc{\E_{\hx_s\sim \rh_s} e^{\hw_t} h(\hx_t)}{\E_{\hx_s\sim \rh_s} e^{\hw_t}} = \E_{x\sim \rh_t} h(x).
\end{align}
Note that $\om(m)$ is bounded; hence $|\hw_s| \le s \max|\om| \to 0$ as $s\to 0$.
Starting at $x_0=0$, we have as $s\to 0$ that
\begin{align}
\label{e:je-limit-2}
    \fc{\E_{\hx_0=0} e^{\hw_t} h(\hx_t)}{\E_{\hx_0=0} e^{\hw_t}}
    & \sim 
    \fc{\E_{\hx_0=0} e^{\hw_t - \hw_s} h(\hx_t)}{\E_{\hx_0=0} e^{\hw_t - \hw_s}}
    = \fc{\E_{\hx_s\sim \hp_s} e^{\hw_t - \hw_s} h(\hx_t)}{\E_{\hx_s\sim \hp_s} e^{\hw_t - \hw_s}}
    = \fc{\E_{\hx_s\sim \hp_s, \hw_s = 0} e^{\hw_t} h(\hx_t)}{\E_{\hx_s\sim \hp_s, \hw_s = 0} e^{\hw_t}}
\end{align}
where $\hp_s$ be the distribution of $\hx_s$ started at $x_0=0$.
It remains to note that, because the drift term $\hm_t$ of the diffusion is uniformly bounded,
\begin{align*}
    \TV(\hp_s, \rh_s)
    \le 
    \TV(\hp_s, \calN(0,s)) + \TV (\calN(0,s), \rh_s) \to 0
\end{align*}
as $s\to 0$, so the difference between the numerators and denominators of \eqref{e:je-limit-1} and \eqref{e:je-limit-2} goes to 0 as $s\to 0$, as hence we obtain $\fc{\E_{\hx_0=0} e^{\hw_t} h(\hx_t)}{\E_{\hx_0=0} e^{\hw_t}} = \E_{x\sim \rh_t} h(x)$, as needed.


\subsection{Wasserstein bound, with application to TAP error}
\label{s:wasserstein}
We first give a Wasserstein bound for general SDE's with different drift and diffusion terms. 
Consider 2 SDE's
\begin{align}\label{e:sdes}
\begin{split}
    dX_t &= f_t(X_t)dt + \Si_t(X_t)dB_t\\
    d\hat X_t &= \hat f_t(\hat{X}_t)dt + \hat \Si_t (\hat X_t)dB_t
    \end{split}
\end{align}
where $f, \hat f:\R\times \R^n\to \R^n$ and $\Si, \hat \Si:\R \times \R^n\to \R^{n\times n}$ are continuous in $t$ and $x$. Suppose that $X_0=\hat X_0$. 
We will think of $X_t$ as an ideal process and $\hat X_t$ as an approximate process. 

\begin{lemma}[Gr\"onwall bound for Wasserstein error of SDE's]\label{l:W-gw}
Suppose that in \eqref{e:sdes}, $f,\hat f, \Sigma, \hat \Sigma$ are Lipschitz in $x\in \R^n$, uniformly over $[0,T]$. Suppose also that for $t\in [0,T]$,
\begin{enumerate}
    \item 
    $\E \ve{f_t(X_t)-\hat f_t(X_t)}^2\le \ep_f(t)^2$ and 
    $\E \ve{\Si_t(X_t)-\hat \Si_t(X_t)}_F^2\le \ep_\Si(t)^2$, where $X_t$ evolves according to the SDE \eqref{e:sdes}.
    \item 
    $\hat f_t$ is $L_f$-Lipschitz and 
    $\hat \Si_t$ is $L_\Si$-Lipschitz (with respect to the $\ve{\cdot}_2$ and $\ve{\cdot}_F$ norms).
\end{enumerate}
    Then for $t\in [0,T]$, 
    with $X_t$ and $\hat X_t$ driven by the same Brownian motion in \eqref{e:sdes} (i.e. synchronous coupling), 
    \[
W_2(\mathcal L(X_t),\mathcal L(\hat X_t))^2
\le 
\E \ba{\ve{X_t -\hat X_t}^2} \le 
\int_0^t 
\pa{\fc{\sqrt 2}{L_f} \ep_f(s)^2 + 2\ep_\Si(s)^2}
e^{(2\sqrt 2 L_f + 2L_\Si^2)(t-s)}ds.
    \]
\end{lemma}
Note that (1) the expectations are with respect to the distribution of the ideal process, and (2) the Lipschitz condition is on the approximations $\hat f_t$ and $\hat \Si_t$, not on $f_t$ and $\Si_t$.
\begin{prf}
    Let $\De_t:= X_t-\hat X_t$. 
    Note 
    \[
d\De_t = (f_t(X_t) - \hat f_t(\hat X_t))dt
+ (\Si_t(X_t) - \hat \Si_t(\hat X_t))dB_t.
    \]
    By Itô's formula, 
    \begin{align*}
d\ve{\De_t}^2 &= 
2\an{\De_t, (f_t(X_t) - \hat f_t(\hat X_t))dt
+ (\Si_t(X_t) - \hat \Si_t(\hat X_t))dB_t}
+ \ve{\Si_t(X_t) - \hat \Si_t(\hat X_t)}_F^2dt\\
\implies
\ddd t\E \ve{\De_t}^2
&= 2 \E \an{\De_t, f_t(X_t) - \hat f_t(\hat X_t)} + 
\E \ve{\Si_t(X_t) - \hat \Si_t(\hat X_t)}_F^2\\
&\le c\E \ve{\De_t}^2 + \rc{c}\E \ve{f_t(X_t)-\hat f_t(\hat X_t)}^2 + \E \ve{\Si_t(X_t) - \hat \Si_t(\hat X_t)}_F^2
    \end{align*}
    for any $c>0$.
Now
\begin{align*}
    \ve{\hat f_t(\hat x) - f_t(x)}^2 &\le 2 \ve{\hat f_t(\hat x) - \hat f_t(x)}^2 + 
    2 \ve{\hat f_t(x) - f_t(x)}^2
    \le 2L_f^2 \ve{\hat x - x}^2 + 2 
    \ve{\hat f_t(x) - f_t(x)}^2
    \\
    \ve{\hat \Si_t(\hat x) - \Si_t(x)}_F^2
    &\le 2\ve{\hat \Si_t(\hat x) - \hat \Si_t(x)}_F^2 + 2\ve{\hat \Si_t(x) -  \Si_t(x)}_F^2
        \le 
        2L_{\Si}^2\ve{\hat x - x}^2 + 2 
        \ve{\hat \Si_t(x) -  \Si_t(x)}_F^2.
    \end{align*}
    Plugging in gives
    \begin{align*}
        \ddd t\E \ve{\De_t}^2
        &\le 
        \pa{c+\fc{2L_f^2}{c} + 2L_\Si^2} \E \ve{\De_t}^2 + \fc{2}{c} \E \ve{\hat f_t(X_t) - f_t(X_t)}^2 + 2 \E \ve{\hat \Si_t(X_t) -  \Si_t(X_t)}_F^2\\
        &\le 
        \pa{c+\fc{2L_f^2}{c} + 2L_\Si^2} \E \ve{\De_t}^2 + \fc 2c \ep_f(t)^2 + 2 \ep_\Si(t)^2\\
        &= 
        \pa{2\sqrt 2 L_f + 2L_\Si^2} \E \ve{\De_t}^2 + \fc{\sqrt 2}{L_f} \ep_f(t)^2 + 2\ep_\Si(t)^2 \text{ by taking $c=\sqrt 2 L_f$.}
    \end{align*}
A Gr\"onwall inequality argument then gives
\[
\E \ve{\De_t}^2 \le \int_0^t 
\pa{\fc{\sqrt 2}{L_f} \ep_f(s)^2 + 2\ep_\Si(s)^2}
e^{(2\sqrt 2 L_f + 2L_\Si^2)(t-s)}ds.
\]
\end{prf}

We apply this to the ideal coupled SL-HD process and approximate coupled ASL-PHD process,
\begin{align}
\label{e:sl-hd-asl-phd}
    &\begin{cases}
        dy_t = \mg_t dt + dB_t&\\
        d\mg_t = Q(\mg_t) dB_t&
    \end{cases}
    &&\begin{cases}
        d\hat y_t = \hat \mg_t dt + dB_t&\\
        d\hat \mg_t = \hat Q(\hat \mg_t) dB_t,&
    \end{cases}
\end{align}
which can be written as 
\begin{align}
\label{e:sl-hd-asl-phd-block}
    d\coltwo{y_t}{\mg_t} &= \coltwo{\mg_t}{0} dt + \coltwo{I}{Q(\mg_t)} dB_t &
    d\coltwo{\hat y_t}{\hat \mg_t} &= \coltwo{\hat \mg_t}{0} dt + \coltwo{I}{\hat Q(\hat \mg_t)} dB_t.
\end{align}
This allows us to obtain a Wasserstein bound for ASL-PHD, using just \ref{d:SDE-existence-uniqueness-stay-in-cube}, \ref{d:Q-error}, and \ref{d:Q-Lip}(a).
\begin{corollary}[$O(1)$ Wasserstein bound for ASL-PHD]
    \label{cor:wasserstein-bound-asl}
    Suppose that in \eqref{e:sl-hd-asl-phd}, $y_0=\hat y_0$ and $m_0 = \hat{m}_0$ and the following hold:
    \begin{enumerate}
        \item The solution to \eqref{e:sl-hd-asl-phd-block} exists and satisfies $\hat{m}_t \in (-1,1)^n$ for all $t$ almost surely.
        \item $\E\ve{Q(\mg_t) - \hat Q(\mg_t)}_F^2 \le \ep(t)^2$.
        \item $\hat Q$ is $L$-Lipschitz in $(-1,1)^n$ with respect to $\ve{\cdot}_F$. 
    \end{enumerate}
    Then
    \[
W_2(\mathcal L(y_t),\mathcal L(\hat y_t))^2
\le 
W_2(\mathcal L(y_t,\mg_t),\mathcal L(\hat y_t,\hat \mg_t))^2
\le \int_0^t 
2\ep(s)^2
e^{(2\sqrt 2 + 2L_\Si^2)(t-s)}ds.
    \]
\end{corollary}
\begin{prf}
    Apply \Cref{l:W-gw} to \eqref{e:sl-hd-asl-phd-block} with $L_f=1$, $L_\Si=L$, $\ep_f(t)=0$, and $\ep_\Si(t)=\ep(t)$.
\end{prf}

Note, however, that this does not allow us to obtain KL divergence bounds via Girsanov's Theorem, because in ASL-PHD, $\hat m_t$ is not a function of $\hat y_t$. For KL divergence bounds, we will need $L^2$ error for ASL-TAP as compared to the ideal process.

\paragraph{Integrating PHD to get ASL/score error}
\label{sec:pha-to-score}
\begin{lemma}[Error in estimated magnetization/score]
\label{l:m-error}
Suppose that \ref{d:SDE-existence-uniqueness-stay-in-cube}, \ref{d:Q-error}, \ref{d:Q-reg}, \ref{d:drift-error}, and \ref{d:Q-Lip} hold to time $T$.
    Let $y_t$ denote the SL process, and $\mg_t=\an{\si}_{y_t}$ be the true magnetization for tilt $y_t$. Let $\hat \mg(y)$ denote the magnetization obtained by solving the TAP equation. Then for each $t\in [0,T]$ there is a constant $C(t)$ such that 
    \[
\E\ba{\ve{\hat \mg(y_t)-\mg_t}^2}\le C(t).
    \]
\end{lemma}
\begin{prf}
    Consider the synchronously coupled SDEs
\begin{align*}
    (\textup{PHD-TAP})\quad &
    \begin{cases}
        d\hm_t = \hat Q(\hm_t)dB_t\\
        d\hy_t = f(\hm_t)dt + dB_t
    \end{cases} \\
    (\textup{Ideal})\quad &
    \begin{cases}
        d\mg_t= Q(\mg_t)dB_t = \Cov(\mu_{y_t}) dB_t\\
        dy_t= \mg_tdt+dB_t
        .
    \end{cases}
\end{align*}
From the desiderata, we have that 
\begin{align}
\label{e:mg-primal-error-2}
    \E\left[\ve{\hat Q(\mg_t) -Q(\mg_t)}_F^2\right] &=O(1)\\
\label{e:mg-primal-error-4}
    \E\left[\ve{f(\mg_t)-\mg_t}^2\right]&=O(1)
\end{align}
where
\eqref{e:mg-primal-error-2} follows from \ref{d:Q-error} and \eqref{e:mg-primal-error-4} follows from \ref{d:drift-error}. 
By \ref{d:Q-Lip}, we also see
$f(\mg)$ and $\hat Q(\mg)$ are $O(1)$-Lipschitz in $(-1,1)^n$. 

Thus, with \ref{d:SDE-existence-uniqueness-stay-in-cube} ensuring non-escape from $(-1,1)^n$, we can invoke \Cref{l:W-gw} to get for some constant $C_1(t)$ that
\[
\E[\lpnorm{(\hm_t,\hy_t) - (\mg_t,y_t)}^2] \le e^{Ct}-1=:C_1(t).
\]
Note $\hm(y)$ has Lipschitz constant $O(1)$ as its Jacobian is $\hat Q(\hm(y))$, which is bounded in operator norm by \ref{d:Q-reg}.
Therefore (note $\hm_t=\hm(\hy_t)$),
\begin{align*}
    \E\ba{\ve{\hm(y_t)-\mg_t}^2}
    &\lesssim 
    \E[\ve{\mg_t-\hm_t}^2] + \E[\ve{\hm(\hy_t)-\hm(y_t)}^2] \\
    &\lesssim \E[\ve{\mg_t-\hm_t}^2] + \E[\ve{\hy_t-y_t}^2] 
    \lesssim C_1(t). \qedhere
\end{align*}

\end{prf}

A standard Girsanov argument (see e.g. \cite[\S 4.4]{chewi2024log}) then bounds the KL divergence between the SL and ASL-TAP.
\begin{corollary}[$O(1)$ KL error for ASL-TAP]
\label{c:O1-KL}
Suppose that \ref{d:SDE-existence-uniqueness-stay-in-cube}, \ref{d:Q-error}, \ref{d:Q-reg}, \ref{d:drift-error}, and \ref{d:Q-Lip} hold to time $T$. Let $y_t$ and $\hy_t$ follow \eqref{e:SL} and \eqref{e:ASL-TAP}, respectively.
Then 
\[
\KL(\dist(y_T)\|\dist(\hy_T)) = O(1).
\]
\end{corollary}
\begin{prf}
Let $P_T$ and $\hat P_T$ denote the path measures of $(y_t)_{t=0}^T$ and $(\hy_t)_{t=0}^T$. 
Since $\mg, \hat \mg$ are bounded, Novikov's condition holds: $\E \exp\pa{\rc 2 \int_0^T \ve{\hat \mg(y_t) -m(y_t)}^2dt}<\iy$. 
By the data processing inequality and Girsanov's Theorem,
\begin{align*}
    \KL(\dist(y_T)\|\dist(\hy_T))
    &\le 
    \KL(P_T\|\hat P_T)
    \le 
    \E_{P_T}
    \log \dd{P_T}{\hat P_T}\\
    &=
    -\E  
    \int_0^T \ip{\mg(y_t) - \hat\mg(y_t)}{dB_t} + \rc 2 
    \int_0^T \ve{\mg(y_t) - \hat\mg(y_t)}^2dt\\
    &= \rc 2 
    \int_0^T \ve{\mg(y_t) - \hat\mg(y_t)}^2 = O(1),
\end{align*}
    where the last equality follows from \pref{l:m-error}.
\end{prf}

\subsection{Error analysis for Jarzynski's equality}
\label{s:je-error}
First we state a standard guarantee for rejection sampling, but allowing for error in sampling from the proposal distribution,  and in computing the ratio. We also allow unbounded ratio of densities. Note that we do not assume the normalizing constant is known, so we first draw a set of samples, take an appropriate quantile, and then use it to rescale our draws for rejection sampling. 
Compared to \cite[Lemma 2]{fan2023improved}, we use a quantile, rather than a single sample as reference; this is necessary to avoid needing a two-sided bound for the density ratio.

\begin{algorithm}[!ht]
\caption{Approximate rejection sampler with unknown normalization 
}
\begin{algorithmic}[1]
\INPUT Oracle for sampling from $\td Q\approx Q$, (possibly randomized) function $\td R\approx R$ where $\dd{P}{Q} \propto R$, cutoff parameters $C_1$, $C_2$, failure probability $\de'$ 
\State Draw $N$ samples from $\td Q$, where $N=\Om(C_1^2\log(1/\de'))$ for appropriate constant
\State Let $R_0$ be the $p$th quantile of the samples, where $p=1-\rc{6C_1}$. 
\Repeat{}
    \State Draw $X\sim \td Q$.
    \State Draw $U\sim \mathsf{Uniform}([0,1])$.
\Until{$U\le \fc{\td R(X)}{C_3R_0}$, where $C_3=8C_2$.}
\OUTPUT $X$.
\end{algorithmic}
\label{a:ars}
\end{algorithm}

\newcommand{\esm}{\varepsilon_{\textup{sample}}}
\newcommand{\ewt}{\varepsilon_{\textup{weight}}}
\newcommand{\ert}{\varepsilon_{\textup{ratio}}}

\begin{lemma}[Approximate rejection sampling with unknown normalization]
\label{l:ars}
Suppose that \pref{a:ars} is run with $\td Q$, $\td R$, and parameters $C_1, C_2$  satisfying the following.
\begin{enumerate}
    \item \label{i:esm}
    (Sampling error) $\td Q$ is a distribution such that $\TV(\td Q, Q)\le \esm$.
    \item \label{i:ewt}
    (Tails of ratio)
    Let 
    \[
    C(\ep) = \inf_C \set{C}{\E_Q\pa{\dd PQ \wedge C}\ge 1-\ep}, 
    \]
    and suppose $C_1\ge C(1/2)\vee 1$, $C_2\ge C(\ewt)$. 
    \item \label{i:ert}
    (Ratio error) 
    Let $G_1$, $G_2$ be events such that 
    \begin{align*}
    G_1 &\subeq \bc{\fc{\td R}{R}\in [0,e^{\ert}]}\\
G:&= G_1\cap G_2 \subeq 
\bc{ \fc{\td R}{R} \in [e^{-\ert}, e^{\ert}] }
    \end{align*}
    and $\td Q(G_1^c)\le \de_1$, $P(G_2^c)\le \de_2$. (Note these probabilities can involve randomness in the algorithm.)
\end{enumerate}
Then letting $\hat P$ be the output distribution, with probability $\ge 1-\de'$ over the $N$ initial samples,
\[
\TV(\hat P, P) \le \ewt + \de_2 + 384C_1C_2(\de_1 + \ert + \esm) 
\]
and if 
$\ewt + \de_2 + 96C_1C_2(\de_1 + \ert + \esm)\le \rc2$, the acceptance probability is at least $\rc{192C_1C_2}$.
\end{lemma}
Note that in assumption \ref{i:ert}, we consider two different bad events separately. We exclude the estimated ratio being too large under the sampling distribution $\td Q$, while we only need to exclude the estimated ratio being too small under $P$. This will be more convenient for us than excluding it under $\td Q$, and this is sufficient because it is okay to underestimate the ratio in places where the target distribution $P$ has small mass anyways.
\begin{prf} We proceed in two steps.
\ppart{Step 1} We show that with probability at least $1-\de'$ that $R_0\in \E_Q[R]\cdot \ba{\rc 8, 12 C_1}$. 
For $0\le p\le 1$, let $q(p)$ denote the $p$th quantile of the distribution of $\dd PQ$ over $Q$.
Suppose $C_1$ is such that $C_1\ge 4[C(1/2)\vee 1]$. 
Let $L_1 = q\pa{1 - \rc{4C_1}}$. Then 
\allowdisplaybreaks
\begin{align*}
    \rc 2 &\le \E_Q\ba{\dd PQ \wedge C_1}\le L_1 \pa{1-\rc{4C_1}} + C_1 \rc{4C_1}\\
    \implies 
    L_1 &\ge \fc{1/4}{1-\rc{4C_1}} \ge \rc 8.
\end{align*}
Let $L_2 = q\pa{1-\rc{12C_1}}$. Then
\begin{align*}
    1 = \E_Q\ba{\dd PQ} \ge \rc{12C_1} L_2 \quad \implies \quad L_2 \le 12C_1.
\end{align*}
Now if $N=\Om\pa{C_1^2 \log(1/\de')}$ with appropriate constant, then taking the $1-\rc{6C_1}$ quantile, 
noting moreover that 
\[
\td Q
\pa{\fc{\td R}{R}
        \nin 
        [e^{-\ert}, e^{\ert}]}\le\de_1 + \esm\le \rc{24C_1},
\]
we obtain the result by Hoeffding's inequality applied to the random variables $\one_{X_i<q\pa{1-\rc{4C_1}}}$ and $\one_{X_i>q\pa{1-\rc{12C_1}}}$ for $X_i$ drawn iid from $\td Q$. 

\ppart{Step 2} We now analyze the rejection sampling. First, we consider the true ratio under the true proposal distribution:
For $C_3=8C_2$, $C_3R_0\in \E_Q[R]\cdot [C_2, 96C_1C_2]$. Note that $\fc{R}{\E_Q[R]} = \dd PQ$. Then 
\allowdisplaybreaks
\begin{align}
\nonumber
    \rc{C_2} &\ge
    \E_{Q}\ba{\fc{R}{C_3R_0}\wedge 1}
    \\
\label{e:ars-1}
    &= 
    \rc{C_3R_0}
    \E_{Q}[R] \E_Q\ba{\dd PQ \wedge\fc{C_3R_0}{\E_Q[R]}}\\
\label{e:ars-2}
    &\ge \rc{C_3R_0} 
    \E_Q[R] (1-\ewt)\ge \rc{96C_1C_2}(1-\ewt),
\end{align}
where we use $\fc{C_3R_0}{\E_Q R}\ge C_2\ge C(\ewt)$ in \eqref{e:ars-1}.
Therefore, for any measurable $A\subeq \Om$,
\begin{align}\label{e:S>C}
\E_Q[(R\wedge C_3R_0)\one_A]
\ge \E_Q[R\one_A] - \ewt \E_Q[R].
\end{align}
Next, we consider the approximate ratio under the approximate proposal distribution. 
Note that 
\allowdisplaybreaks
\begin{align}
\label{e:ars-hPA}
    \hat P(A) &= 
    \fc{\E_{\td Q} [(\td R\wedge C_3R_0) \one_{A}]}{\E_{\td Q} [\td R\wedge C_3R_0]}\\
\label{e:ars-PA}
    P(A) &= \fc{\E_{Q} [R \one_{A}]}{\E_{Q} [R]}.
\end{align}
Note also that 
\begin{align}
\label{e:ars-de2}
    \E_Q[(R\wedge C_3R_0)\one_A \one_{G_2^c}]
    &\le \E_Q [R\one_{G_2^c}] = \E_Q[R] P(G_2^c)\le \de_2 \E_Q[R]
\end{align}
    Using $\TV(\td Q, Q)\le \esm$, we have that the numerator in \eqref{e:ars-hPA} is
\allowdisplaybreaks
\begin{align*}
    &\E_{\td Q}[(\td R \wedge C_3R_0) \one_{A}]\\
    &\in_{\eqref{i:ert}} 
    \ba{ \E_{\td Q} (Re^{-\ert} \wedge C_3R_0) \one_{A}\one_G, 
    \E_{\td Q} (Re^{\ert} \wedge C_3R_0) \one_{A}\one_{G_1} + C_3R_0\de_1
    }\\
    &\subeq \ba{\E_{\td Q} (R\wedge C_3R_0) \one_{A}\one_G - C_3R_0\ert, 
    \E_{\td Q} (R\wedge C_3R_0) \one_{A}\one_{G_1} + C_3R_0(\de_1+\ert)}\\
    &\subeq_{\eqref{i:ert}}
    \ba{\E_{\td Q} (R\wedge C_3R_0) \one_{A}\one_{G_2} - C_3R_0(\de_1+\ert), 
    \E_{\td Q} (R\wedge C_3R_0) \one_{A} + C_3R_0(\de_1+\ert)}\\
    &\subeq_{\eqref{i:esm}}
     \ba{\E_{Q} (R\wedge C_3R_0) \one_{A} \one_{G_2}  - C_3R_0(\de_1 + \ert + \esm), 
     \E_{Q} (R\wedge C_3R_0) \one_{A} + 
     C_3R_0(\de_1 + \ert+\esm)}\\
     &\subeq_{\eqref{e:ars-de2}}
     \E_{Q} (R\wedge C_3R_0) \one_{A} + \ba{- \de_2 \E_Q[R] - C_3R_0(\de_1 + \ert + \esm), 
     C_3R_0(\de_1 + \ert+\esm)}\\
     &\subeq_{\eqref{e:S>C}}
      \E_{Q} R \one_{A}
     + \ba{ -(\ewt+\de_2) \E_{Q} [R] - C_3R_0(\de_1 + \ert + \esm), 
     C_3R_0(\de_1 + \ert+\esm)}\\
     &\subeq \E_{Q} R \one_{A}
     + \E_{Q} [R]\cdot \ba{ -\ewt - \de_2  - 96C_1C_2(\de_1 + \ert + \esm), 
      96C_1C_2(\de_1 + \ert+\esm)}.
\end{align*}
Hence 
\begin{align*}
    &\ab{\hat P(A) - P(A)} =
    \ab{\fc{\E_{\td Q} [(\td R\wedge C_3R_0) \one_{A}]}{\E_{\td Q} [\td R\wedge C_3R_0]}- \fc{\E_{Q} [R \one_{A}]}{\E_{Q} [R]}}\\
    &\le 
    \max \bc{
P(A) - \fc{P(A) -\ewt - \de_2 - 96C_1C_2(\de_1 + \ert + \esm)}{1+ 96C_1C_2(\de_1 + \ert + \esm)}
, \fc{P(A) +  96C_1C_2(\de_1 + \ert + \esm)}{1- 96C_1C_2(\de_1 + \ert + \esm)} - P(A),1
    }\\
    &\le \ewt + \de_2 + 384C(\de_1 + \ert + \esm),
\end{align*}
where the last inequality follows 
using $\fc{p+\ep_1}{1-\ep_2}\wedge 1\le p + 2\ep_1+2\ep_2$ if $p\le 1$ (split into cases based on whether $\ep_2\le \rc 2$).
Taking a maximum over measurable $A$ gives the TV bound. Finally, the acceptance probability equals
\begin{align*}
    \E_{\td Q} \ba{\fc{\td R}{C_3R_0} \wedge 1} 
    &= \rc{C_3R_0} \E_{\td Q} [\td R\wedge C_3R_0]
    \ge \rc{C_3R_0} \E_{Q}[R]\cdot \pa{1-\ewt - \de_2 - 96C_1C_2(\de_1 + \ert + \esm)}\\
    &\ge \rc{96C_1C_2}  \pa{1-\ewt - \de_2 - 96C_1C_2(\de_1 + \ert + \esm)}\ge \rc{192 C_1C_2}.
\end{align*}
where we use~\eqref{e:ars-2} in the last line.
\end{prf}

For Jarzynski's equality, we give TV sampling guarantees with a number of draws that is exponential in the total drift error. In addition to a bound on the drift error, we need an exponential upper tail and a lower tail for the log-weights. Note that these are over the algorithmic process. 

There are 3 distributions in our analysis: the target measure $ p_t$ of the ideal process, the measure $\hat p_t$ of the approximate process, and the scaffolding measure $\rh_t$ with JE targets. Samples from $\hat  p_T$ are reweighted by $e^{w_T}$ to approximate the measure $\rh_T$, because we have an explicit expression for $\rh_T$ but not $\hat  p_T$. Then we can conduct rejection sampling with ratio $\dd{ p_T}{\rh_T}$ (or an approximation of this ratio) to obtain (approximate) samples from $ p_T$.
\begin{lemma}[Jarzynski's equality with rejection sampling]
\label{l:je-rs}
    Suppose that $f_t$, $\hat f_t$ are Lipschitz in $x\in \R^n$ uniformly over $t\in [0,T]$, 
    \begin{align}
        & dx_t = f_t(x_t) dt + dB_t,&& x_0\sim  p_0\\
        \label{e:je-generic}
        &
        d\hx_t = \hf_t(\hx_t) dt+ dB_t 
        && 
        \hx_0\sim  p_0 
    \end{align}
    and define $\hw_t((\hx_s)_{s\in [0,T]})$ for $t\in [0,T]$ by 
    \begin{align}
    \label{e:je-generic-2}
        &d\hw_t = \om_t(\hx_t)dt, && \hw_0=0
    \end{align}
    where \eqref{e:je-generic} and \eqref{e:je-generic-2}
    is a SDE-ODE system arising from Jarzynski's equality applied to $\rh_t$, where $\rh_0= p_0$.  
    Let $ p_t$, $\hp_t$ be the distribution of $x_t$, $\hx_t$, respectively. 
    Below, $\hw_t=\hw_t((\hx_s)_{s\in [0,T]})$ unless specified otherwise.
    Let $Z_t = \E[e^{\hw_t}]$.
    Suppose
    \begin{enumerate}
        \item (Drift error)
        $\E\ba{\ve{\hf_t(x_t)-f_t(x_t)}^2}\le \ep(t)^2$.
        \item (Exponential upper tail for log-weights)
        $\E\ba{e^{\hw_T}\one_{e^{\hw_T}/Z_T\ge C_1}}\le \ep_1 Z_T$, for $\ep_1 = \fc{c_3\ep}{6 \exp\pa{\fc{3}{\ep}\pa{\fc 12 \int_0^T \ep(t)^2dt+1}}}$.
        \item (Lower tail for log-weights)
        $\P\ba{\fc{e^{\hw_T}}{Z_T}\le c_3}\le \ep_3$, for $\ep_3 = \fc{\ep}{6 \exp\pa{\fc{3}{\ep}\pa{\fc 12 \int_0^T \ep(t)^2dt+1}}}$.
    \end{enumerate}
    Then for $L=\exp\pa{\fc{3}{\ep}\pa{\fc 12 \int_0^T \ep(t)^2dt+1}}$, $\E\ba{
            \fc{e^{\hw_T}}{Z_T} \dd{ p_T}{\rh_T}(\hx_T)\wedge \fc{2C_1L}{c_3}}\ge 1-\ep$, and given query access to $\dd{ p_T}{\rh_T}$, we can sample from $ p_T$ with TV error $2\ep$ by $O\pf{C_1^2L^2\log\prc\ep}{c_3^2}$ simulations of \eqref{e:je-generic}, by choosing the parameters  $C_1'=C_2'=\fc{2C_1L}{c_3}$ and using the estimate 
            $e^{\hw_T} \dd{ p_T}{\rh_T}(\hx_T)$
            in \pref{a:ars}. 
Suppose moreover that $(\tx_t)_{t\in [0,T]}$ is a process and $\tw_T((x_t)_{t\in [0,T]})$ are weights such that 
\begin{enumerate}[resume]
    \item (Path error) \label{i:je-path-error}
    $\TV((\hx_t)_{t\in [0,T]}, (\tx_t)_{t\in [0,T]}) \le \ep_4$, where $\ep_4 = \fc{\ep}{384C_1^{\prime2}}$.
    \item (Weight error) \label{i:je-ewt}
    $\P\pa{|\tw_T((\tx_t)_{t\in [0,T]}) - \hw_T((\tx_t)_{t\in [0,T]})|>\fc{\ep_4}2}\le \fc{\ep_4}2$. 
    \item (Ratio error) \label{i:je-ert}
    $\td R_2(x)$ is a (randomized) estimate such that there is $R_2(x)$ with $\dd{ p_T}{\rh_T} \propto R_2$ with
\begin{align*}
(a) &&\forall x,\quad  \P\pa{\fc{\td R_2(x_T)}{R_2(x_T)}\nin [0, e^{\ep_4/2}] \Big| x_T=x} &\le \fc{\ep_4}2\\
(b)&&    \P\pa{\fc{\td R_2(x_T)}{R_2(x_T)}\nin [e^{-\ep_4/2}, e^{\ep_4/2}]}&\le \ep
\end{align*}
\end{enumerate}
Then we can sample from $ p_T$ with TV error $6\ep$ by $O\pf{C_1^2L^2\log\prc\ep}{c_3^2}$ simulations of \eqref{e:je-generic} using the estimate $e^{\tw_T((\td x_t)_{t\in [0,T]})} \td R_2(\tx_T)$
            in \pref{a:ars}.
\end{lemma}
Note for the ratio estimate, in addition for having a 2-sided bound with high probability over the target distribution, 
we need that for any point $x$, it is with high probability not much of an overestimate.
\begin{prf}
    Note first that from Jarzynski's equality,
    \[
\fc{\E[e^{\hw_t}|\hx_t]}{Z_t} = \dd{\rh_t}{\hp_t}(\hx_t).
    \]
    Define the good sets (1) the ratio between the ideal and approximate distribution is bounded, (2) the ratio between the scaffolding and approximate distribution is not too small, and (3) the weights are bounded,
    \begin{align*}
        G_1 &= \bc{\dd{ p_T}{\hp_T}\le L}\\
        G_2 &= \bc{\fc{\E[e^{\hw_T}|\hx_T]}{Z_T}=\dd{\rh_T}{\hp_T}\ge \fc{c_3}2}\\
        G_3 &= \bc{\fc{e^{\hw_T}}{Z_T}\le C_1}
    \end{align*}
    and let their complements be $B_1,B_2,B_3$, respectively. 
    Note $G_1,G_2$ are $x_T$-measureable and $G_3$ is $\hat w_T$-measurable. 
    Now
    \begin{align}
        \E\ba{\fc{e^{\hw_T}}{Z_T} \dd{ p_T}{\rh_T}(\hx_T) \one_{B_1\cup B_2 \cup B_3}}
        &\le 
        \E\ba{\fc{e^{\hw_T}}{Z_T}\dd{ p_T}{\rh_T}(\hx_T)\one_{B_1} }
        + \E\ba{\fc{e^{\hw_T}}{Z_T}\dd{ p_T}{\rh_t}(\hx_T)\one_{G_1\cap B_2} }
        + \E\ba{\fc{e^{\hw_T}}{Z_T}\dd{ p_T}{\rh_t}(\hx_T)\one_{G_1\cap G_2\cap B_3} }.
    \end{align}
    We bound each term separately. 
    Note that if $\mu\ll \nu$, then 
    \begin{align}
    \nonumber
        \KL(\mu\|\nu)& = \E_\mu \pa{\log \dd{\mu}{\nu}}\one_{\dd{\mu}{\nu}\ge 1} + \E_\mu\pa{\log \dd{\mu}{\nu}} \one_{\dd{\mu}{\nu}< 1} 
        \ge  \E_\mu \pa{\log \dd{\mu}{\nu}}\one_{\dd{\mu}{\nu}\ge 1} - 1
        \ge \mu\pa{\dd{\mu}{\nu}\ge L}(\log L)-1\\
        \implies
        \mu\pa{\dd{\mu}{\nu}\ge L} & \le \fc{\KL(\mu\|\nu) + 1}{\log L}.
        \label{e:kl-markov}
    \end{align}
    For the first term,
    \begin{align*}
        \E\ba{\fc{e^{\hw_T}}{Z_T}\dd{ p_t}{\rh_t}(\hx_T)\one_{B_1} }
        &= \E\ba{\E\ba{\fc{e^{\hw_T}}{Z_T}|\hx_T}\dd{ p_T}{\rh_T}(\hx_T)\one_{B_1} }\\
        &= \E \ba{\dd{\rh_T}{\hp_T}(\hx_T)\dd{ p_T}{\rh_T}(\hx_T)\one_{B_1}}
        = \E \ba{\dd{ p_T}{\hp_T}(\hx_T)\one_{B_1}}\\
        &= \E\ba{\one_{B_1}(x_T^*)}\le \P\ba{\dd{ p_T}{\hp_T}(x_T)>L}
        = \P\ba{\log\dd{ p_T}{\hp_T}(x_T)>\log L}\\
        &\le_{\eqref{e:kl-markov}} \fc{\KL( p_T\|\hp_T)+1}{\log L}
        \le \fc{\rc 2 \int_0^T \ep(t)^2dt + 1}{\log L}\le 
        \fc{\ep}3
    \end{align*}
    by assumption 1 and Girsanov's theorem and choosing $L=\exp\pa{\fc{3}{\ep}\pa{\fc 12 \int_0^T \ep(t)^2dt+1}}$. 
    Next, 
    note that by assumption 3, since $\fc{\E[e^{\hw_T}|\hx_T]}{Z_T}\le \fc{c_3}2$ implies 
    $\P\ba{e^{\hw_T} \ge c_3Z_T | \hx_T}\le \fc{\E[e^{\hw_T}|\hx_T]}{c_3 Z_T}\le  \rc2$ and hence
    $\P\ba{e^{\hw_T} \le c_3Z_T | \hx_T}\ge \rc2$, 
    \begin{align*}
        \ep_3 \ge \P\ba{e^{\hw_T} \le c_3Z_T}
        &= \E\ba{\P\ba{e^{\hw_T} \le c_3Z_T | \hx_T}}\\
        &\ge  \E\ba{\rc2\one_{\fc{\E[e^{\hw_T}|\hx_T]}{Z_T}\le \fc{c_3}2}}= \rc 2\P\ba{\fc{\E[e^{\hw_T}|\hx_T]}{Z_T}\le \fc{c_3}2}
        \ge \rc 2\P[B_2].
    \end{align*}
    Hence we can bound the second term
    \begin{align*}
        \E\ba{\fc{e^{\hw_T}}{Z_T}\dd{ p_T}{\rh_t}(\hx_T)\one_{G_1\cap B_2} }
        &= \E\ba{\E\ba{\fc{e^{\hw_T}}{Z_T}|\hx_T} \dd{ p_T}{\rh_T}(\hx_T)\one_{G_1\cap B_2}} = 
        \E\ba{\dd{ p_T}{\hp_T}(\hx_T)\one_{G_1\cap B_2}}\le L \P[B_2] \le 2L\ep_3.
    \end{align*}
    Finally, note that under $G_1\cap G_2$, 
    $\dd{ p_T}{\rh_T} = \dd{ p_T}{\hp_T}\big/ \dd{\rh_T}{\hp_T} \le \fc{2L}{c_3}$, so 
    for the third term, by assumption 2
    \begin{align*}
        \E\ba{\fc{e^{\hw_T}}{Z_T}\dd{ p_T}{\rh_t}(\hx_T)\one_{G_1\cap G_2\cap B_3} }
        &\le 
        \fc{2L}{c_3}\E\ba{\fc{e^{\hw_T}}{Z_T}\one_{B_3}}
        \le \fc{2L\ep_1}{c_3}.
    \end{align*}
    Putting everything together, 
    \begin{align*}
        \E\ba{\fc{e^{\hw_T}}{Z_T} \dd{ p_T}{\rh_T}(\hx_T) \one_{B_1\cup B_2 \cup B_3}}
        &\le \fc{\ep}3 + 2L\ep_3 + \fc{2L\ep_1}{c_3}=\ep
    \end{align*}
    by choosing $\ep_3 = \fc{\ep}{6L}$ and $\ep_1 = \fc{c_3\ep}{6L}$. Hence 
    \begin{align*}
        \E\ba{
            \fc{e^{\hw_T}}{Z_T} \dd{ p_T}{\rh_T}(\hx_T)\wedge \fc{2C_1L}{c_3}
        } \ge \E\ba{\fc{e^{\hw_T}}{Z_T} \dd{ p_T}{\rh_T}\one_{G_1\cap G_2 \cap G_3}}
        = \E\ba{\dd{ p_T}{\hat  p_T}(x_T)} - 
        \E\ba{\fc{e^{\hw_T}}{Z_T} \dd{ p_T}{\rh_T}\one_{B_1\cup B_2\cup B_3}}
        \ge 1-\ep.
    \end{align*}  
    The first part now follows from \Cref{l:ars}, with $\esm = \de = \ert=0$. Here, we are considering path measures $P$ and $Q$ where 
    $Q = \hat  p_{[0,T]}$, the measure of $(\hx_t)_{t\in [0,T]}$, and $P$ is such that $\dd PQ((x_t)_{t\in [0,T]}) = 
    \dd{ p_T}{\hat p_T}(x_T)$, and then taking the marginal distribution of $x_T$. Note we allow probability $\ep$ of failure (in the quantile estimate or rejection sampling taking too many trials).

    For the second part, let $\td Q = \td  p_{[0,T]}$ be the path measure of $(\td x_t)_{t\in [0,T]}$. 
    Let 
    $\td R = e^{\td w_T}\td R_2$, and $R = e^{\hw_T}R_2$. 
    Note that assumption \ref{i:je-path-error} gives assumption \ref{i:esm} of \Cref{l:ars} with $\esm=\ep_4$. Consider the events
    \begin{align*}
        G_1' &= \bc{\fc{\td R_2}{R_2}\in [0,e^{\ep_4/2}]} \cap 
        \bc{|\tw_T - \hw_T|\le \fc{\ep_4}2}
        \subeq \bc{\fc{\td R}{R}\in [0,e^{\ep_4}]}.
        \\
        G_2' &= \bc{\fc{\td R_2}{R_2}\in [e^{-\ep_4/2}, e^{\ep_4/2}]},
    \end{align*}
    and note $G_1'\cap G_2'\subeq \bc{\fc{\td R}{R}\in [e^{-\ep_4},e^{\ep_4}]}$. It remains to check assumption \ref{i:ert} of \Cref{l:ars}. We check
    \begin{align*}
        \td Q({G_1'}^c)
        &\le \P\pa{\fc{\td R_2(\tx_T)}{R_2(\td x_T)}\nin [0,e^{\ep_4/2}]}
        + \P\pa{|\tw_T - \hw_T|>\fc{\ep_4}2}
        \le_{\eqref{i:je-ewt}, \eqref{i:je-ert}(a)} \fc{\ep_4}2+\fc{\ep_4}2 = \ep_4\\
        P({G_1'}^c)
        &= \P \pa{\fc{\td R_2(\tx_T)}{R_2(\td x_T)}\nin [e^{-\ep_4/2},e^{\ep_4/2}]
        } \le_{\eqref{i:je-ert}(b)} \ep.
    \end{align*}
    Thus the assumption holds with $\ert=\de_1=\de_2=\ep_4$. The final TV error of the sampled distribution $\hat P$ is
    \[
\TV(\hat P, P)
\le \ewt + \de_2 + 384 C_1'C_2'(\de_1+\ert+ \esm)
\le \ep +\ep + 384 {C_1'}^2(\ep_4+\ep_4+\ep_4) = 5\ep.
    \]
    Allowing probability $\ep$ of failure then gives the result.
\end{prf}

Having a bounded exponential moment on $[-1,\lm_+>1]$ is sufficient to satisfy assumptions 2 and 3 of \Cref{l:je-rs}.
\begin{lemma} \label{l:wt-tail}
    Suppose that $W$ is a random variable such that for some $\lm_+>1$, for $\lm\in \{-1,1,\lm_+\}$, we have
    \[
\E[e^{\lm W - \E[\lm W]}] \le K_{\lm}. 
    \]
    Then letting $Z = \E[e^W]$, 
    \begin{align*}
       \E[e^W \one_{e^W/Z>C_1}]& \le \ep_1 Z&
       \text{when }C_1 &= \pf{K_{\lm_+}}{\ep_1}^{\rc{\lm_+-1}}\\
       \P\ba{\fc{e^W}{Z}\le c_3} &\le \ep_3 & \text{when }c_3 &= \fc{K_{-1}K_1}{\ep_3}.
    \end{align*}
\end{lemma}
\begin{prf}
    By Jensen's inequality and the mgf bound,
    \begin{align*}
        1\le \fc{Z}{e^{\E W}} \le K_1.
    \end{align*}
    Then
    \begin{align*}
        K_{\lm_+} e^{\lm_+\E W} &\ge 
        \E[e^{\lm_+W}] \ge 
        (ZC_1)^{\lm_+-1} \E[e^W \one_{e^W/Z\ge C_1}]\\
        \implies K_{\lm_+} C_1^{ - (\lm_+-1)} Z&\ge K_{\lm_+} C_1^{ - (\lm_+-1)}  e^{\lm_+\E W} Z^{-(\lm_+-1)} \ge \E[e^W\one_{e^W/Z\ge C_1}].
    \end{align*}
    Taking $C_1$ as defined gives the first inequality.

    For the second inequality, 
    \begin{align*}
        \P(e^W/Z\le c_3) = \P(e^{-W} \ge Zc_3^{-1}) 
        \le \fc{\E[e^{-W}]}{Z^{-1}c_3^{-1}}
        \le \fc{K_{-1}e^{-\E W}}{Z^{-1}c_3^{-1}}
        \le K_1K_{-1}c_3.
    \end{align*}
    Taking $c_3$ as defined gives the second inequality.
\end{prf}

\subsection{Discretization analysis}
\label{s:disc}
We need to consider the discretization of the SDE for the tilt $y_t$, as well as the ODE for the weights $w_t$. 
First, we use the standard result that the path measure of the interpolated discrete SDE is close to that of the actual SDE by Girsanov's theorem. Then for the weights, we conduct the analysis on the interpolated discrete SDE, so the analysis reduces to bounding Brownian motion.
This will show that conditions (4) and (5) in \Cref{l:je-rs} are satisfied.

\begin{lemma}[{SDE discretization bound, cf. \cite[Lemma 4.4.1]{chewi2024log}, \cite{dalalyan2012sparse}}]
\label{l:sde-disc}
Fix a stepsize $h>0$ and let $t^-=h\fl{\frac th}$.
Suppose that $f:\R^n\to \R^n$ is $L$-Lipschitz. 
    Consider the two SDE's in $\R^n$,
    \begin{align*}
        dy_t &= f(y_t) dt + dB_t\\
        d\ty_t &= f(\ty_{t^-}) dt + dB_t,
    \end{align*}
    where $y_0,\ty_0\sim \rh_0$.
    Let $T=Nh$, $N\in \N$. 
    Then
    \begin{align*}
\KL(\calL((\ty_t)_{t\in [0,T]})\| \calL((y_t)_{t\in [0,T]}))
&\le 
\fc{L^2h^3}{6} \sum_{k=0}^{N-1} \E\ve{f(\ty_{kh})}^2 + \fc{T L^2 n h}{2}\\
&\le 
TL^2\pa{\fc{h^2M^2}{6}+ \fc{n h}{2}}
    \end{align*}
where the last inequality holds if $\ve{f}_{\iy}\le M$. If $h\le \fc{\sqrt{3/2}\ep}{\sqrt T LM}\wedge \fc{\ep^2}{2nTL^2}$, then $\TV(\calL((\ty_t)_{t\in [0,T]}), \calL((y_t)_{t\in [0,T]}))\le \ep$.
\end{lemma}
\begin{prf}
    Note that while \cite[Lemma 4.4.1]{chewi2024log} specifically considers Langevin diffusion, where $f=\gd V$, the argument extends to any Lipschitz $f$. Since the diffusion term is $dB_t$ rather than $\sqrt 2 dB_t$, the error is scaled by $\rc2$.
    The last statement holds by Pinsker's inequality, $\TV(\mu,\nu)\le \sqrt{2\KL(\mu\|\nu)}$. 
\end{prf}

\begin{lemma}[Discretization for JE weights under interpolated discrete SDE]
\label{l:je-wt-disc}
    Fix a stepsize $h>0$ and let $t^-=h\fl{\frac th}$.
    Consider a discrete-time process
    \[
\ty_{t+h} = \ty_t + h\hf(\ty_t) + \sqrt h \xi_t , \quad \xi_t\sim \calN(0,\Id_n), \quad \text{for }t=kh, \,k\in \N_0.
    \]
    Consider the continuous-time interpolation of this process, 
    \[d\ty_t = \hf(\ty_{t^-}) dt + dB_t\]
    and consider
    \[
d\hw_t = \om(\ty_t)dt.
    \]
    Consider the discretization of this ODE,
    \[
    \tw_{t+h} = \tw_t + h\om(\ty_t), \quad \text{for }t=kh,\,k\in \N_0.
    \]
    Suppose $\ve{\hf}_\iy\le M$. 
    Then for some constant $C$, for $T=Nh$, $N\in \N$, we have
    \begin{align}
    \label{e:disc-wt-tail-bound}
    \P\pa{
|\tw_T - \hw_T|
\ge L_\om \pa{\fc{TMh}2 + \fc{Th^{1/2}n^{1/2}}{2^{1/2}} + \sqrt{Tn}hu}
    }\le e^{-\fc{u^2}{C}}
\end{align}
Thus for $h\le \fc{\ep}{2L_\om TM} \wedge \fc{\ep^2}{2L_\om^2 T^2 n} \wedge \fc{1}{4L_\om \sqrt{TnC \log \prc{\de}}}$, we have $\P\pa{
|\tw_T - \hw_T|
\ge\ep}\le \de$.
\end{lemma}

\begin{prf}
Note
\begin{align*}
    \ty_t &= \ty_{t^-} + (t-t^-) \hf(\ty_{t^-}) + (B_t - B_{t^-})\\
    \ve{\ty_t - \ty_{t^-}}&\le 
    (t-t^-) M + \ve{B_t - B_{t^-}}
\end{align*}
    We first bound
    \allowdisplaybreaks
    \begin{align*}
        |\tw_T - \hw_T| 
        &= \int_0^T |\om(\ty_t) - \om(\ty_{t^-})|dt\\
        &\le L_\om \int_0^T \ve{\ty_t - \ty_{t^-}}dt\\
        &\le L_\om \sumz k{N-1} \int_{kh}^{(k+1)h} (t-kh) M + \ve{B_t - B_{kh}}dt\\
        &\le L_\om \pa{\fc{TMh}2 + \sumz k{N-1} \int_{kh}^{(k+1)h}\ve{B_t - B_{kh}}dt}
    \end{align*}
The second term is a sum of $N$ iid copies of the random variable $\int_0^h \ve{B_t}dt$. Note by Jensen's inequality that 
\allowdisplaybreaks
\begin{align*}
    e^{\lm \int_0^h \ve{B_t}dt}
    &\le \rc h\int_0^h e^{\lm h \ve{B_t}}dt.
\end{align*}
Note that 
$\ve{\ve{B_t}}_{\psi_2}^2 = 
\ve{\ve{B_t}^2}_{\psi_1} \lesssim tn$ as $\ve{B_t}^2$ is $t$ times a $\chi_n^2$ random variable, which is a sum of $n$ iid random variables satisfying $\ve{X}_{\psi_1} = O(1)$.
Hence $\ve{\int_0^h \ve{B_t}dt - \E \int_0^h \ve{B_t}dt}_{\psi_2}^2\lesssim \ve{h\ve{B_h}}_{\psi_2}^2\lesssim h^3n$.
Note
\begin{align*}
    \E \int_0^h \ve{B_t}dt \le h^{1/2} \pa{\int_0^h \ve{B_t}^2dt}^{1/2}\le h^{1/2}\pa{\int_0^h tn\,dt} = \sfc{h^3n}2.
\end{align*}
Hence we obtain \eqref{e:disc-wt-tail-bound}.
\end{prf}

\section{Estimating the covariance of tilted Gibbs measures via the TAP Hessian}\label{s:est-cov}
Given a tilt $y_t$, assuming that $m_t$ is a stationary point of the TAP free energy $m\mapsto \calF_{\TAP}(\beta A,m,y_t)$, one expects that with high probability over $A$ and $y_t$ sampled from (SL), that
\[
    \Cov(\mu_{\beta A,y_t}) \approx  \left(\nabla^2_m \calF_{\,\mathsf{TAP}}\left(\beta A, m_t, y_t\right)\right)^{-1}\,.
\]
Differentiating explicitly, the Hessian of the TAP free energy (evaluated at magnetization $\mg_t$) is 
\[
    -\beta A - \frac{\beta^2}{n}x_0x_0^\sT + D(\mg_t) + \beta^2\left(1-\frac{\iprod{\mg_t,\mg_t}}{n}\right)\Id_n - \fc{2\beta^2}{n}\mg_t\mg_t^\sT\,,
\]
where $D(\mg)$ is the diagonal matrix with $D(\mg)_{ii}=\frac{1}{1-\mg_i^2}$.   

The true covariance matrix for the planted SK model with a SL tilt is given as $P(A,y_t): = \Cov(G_{A,y_t,\beta})$. Furthermore, the SL process runs with a SDE that satisfies the following equivalent rewrite,
\[
    x_0 \sim G_{A,\beta},\,\,\,\,\, y_t = tx_0 + B_t,\,\,B_t \sim \calN(0,\Id_n)\,.
\]
We are interested in bounding the following quantity
\[
    \E_{A,B_t}\left[\norm{\hat{Q}^{-1}P - \Id_n}_F^2\right] \le O_\beta(1)\,,
\]
for every fixed choice of $x_0 \in \{-1,1\}^n$. At that point, the contiguity between the planted model and the original model~\cite[\S2]{el2022sampling} yields that a $O_\beta(1)$ bound also holds for the original model under the ASL process. The main result of this section is a rigorous proof of the aforementioned desired bound.

\subsection{Gibbs measures, Gaussian processes, and interpolations} The results in this section will heavily use certain facts about the derivatives of Gibbs measures with quadratic interactions and Gaussian processes. These facts (given below) primarily memorize the fact that the derivative of moments of Gibbs measures leads to larger (explicit) moments, while the derivatives of Gaussian processes are evaluated using Stein's lemma.

The other critical fact that will be used (in various forms) in \pref{sec:cavity-interpolation} is that the time derivative of a time interpolated Gaussian processes can be expressed in terms of its second spatial derivatives (resembling the behavior of the heat equation PDE).

\paragraph{Basic analytic properties of Gibbs measures}
Below we state two elementary propositions. The first computes the derivative of the tilted Gibbs measure with respect to the interactions $A_{i,j}$, and the second extends this to obtain derivatives of the first two moments of the Gibbs measure.
\begin{proposition}[Derivatives of the tilted Gibbs measure]\label{prop:tilted-gibbs-derivatives}
    Let the tilted partition function be%
    \[
        Z_{A,\beta, x_0} := \sum_{\sigma \in \{-1,1\}^n} e^{\frac{\beta}{2}\sum_{i , j}A_{i,j}\sigma_i\sigma_j + \an{\sigma,y_t}}\,.
    \]
    Then, the derivative with respect to the interaction variable $A_{i,j}$ is given as
    \begin{align*}
        \partial_{A_{i,j}} G_{A,\beta,x_0}(\sigma) &= \partial_{A_{i,j}} \left(\frac{e^{H(\sigma)}}{Z_{A,\beta, x_0}}\right) = \frac{e^{H(\sigma)}\beta/2\left(Z_{A,\beta, x_0}\left(\sigma_i\sigma_j\right) - \sum_{\tau}\left(\tau_i\tau_j\right)e^{H(\tau)}\right)}{Z^2_{A,\beta, x_0}} \\
        &= \frac{\beta}{2} G_{A,\beta,x_0}(\sigma)\left(\sigma_i\sigma_j - \an{\tau_i\tau_j}_{G_{A,\beta,x_0}}\right)\,.
    \end{align*}
\end{proposition}
\begin{prf}
    Treat $A_{i,j}$ as a variable and apply product rule followed by an elementary factorization argument to rewrite the Gibbs average.
\end{prf}

\begin{proposition}[Derivatives of first and second moments under the symmetric convention]\label{prop:derivatives-gibbs-moments}
Let $A=A^{\sT}$
and consider the Gibbs measure $G_{A,\beta,y_t}$ on $\{\pm1\}^n$ with Hamiltonian
\[
H(\sigma)\;:=\;\frac{\beta}{2}\,\sigma^{\sT} A\sigma \;+\;\frac{\beta^2}{2n}\;+\;\langle y_t,\sigma\rangle,
\qquad 
G_{A,\beta,y_t}(\sigma)=\frac{e^{H(\sigma)}}{Z}.
\]
Write $\langle\cdot\rangle$ for expectation under $G_{A,\beta,y_t}$.  For $k<\ell$ let
$\partial_{k,\ell}:=\partial/\partial A_{k\ell}$ denote differentiation with respect to the
\emph{independent} upper-triangular GOE entry (so $A_{\ell k}=A_{k\ell}$), and let $\partial_{k,k}:=\partial/\partial A_{kk}$.

Define centered spins and centered ``edge'' observables by
\[
m_i:=\langle\sigma_i\rangle,\qquad 
\td\sigma_i:=\sigma_i-m_i,\qquad
\widetilde{\sigma_{ab}}:=\sigma_a\sigma_b-\langle\sigma_a\sigma_b\rangle.
\]
Then for every $k<\ell$, every $i,j\in[n]$, and every $k_i,k_\ell\in\{1,2\}$,
\[
\begin{aligned}
\partial_{k,\ell}\,\langle\sigma_i\rangle^{k_i}
&= (2\langle\sigma_i\rangle)^{k_i-1}\,\beta\,\big\langle \td\sigma_i\,\widetilde{\sigma_{k\ell}}\big\rangle,\\
\partial_{k,\ell}\,\langle\sigma_\ell\rangle^{k_\ell}
&= (2\langle\sigma_\ell\rangle)^{k_\ell-1}\,\beta\,\big\langle \td\sigma_\ell\,\widetilde{\sigma_{k\ell}}\big\rangle,\\
\partial_{k,\ell}\,\langle\sigma_i\sigma_j\rangle
&= \beta\,\big\langle \widetilde{\sigma_{ij}}\,\widetilde{\sigma_{k\ell}}\big\rangle,\\
\partial_{k,\ell}\,\langle\sigma_i\sigma_\ell\rangle
&= \beta\,\big\langle \widetilde{\sigma_{i\ell}}\,\widetilde{\sigma_{k\ell}}\big\rangle.
\end{aligned}
\]
Moreover, for every $k\in[n]$ and every observable $F(\sigma)$ one has $\partial_{k,k}\langle F\rangle=0$
(in particular $\partial_{k,k}\langle\sigma_i\rangle^{k_i}=0$ and $\partial_{k,k}\langle\sigma_i\sigma_j\rangle=0$).
\end{proposition}

\begin{proof}
For $k<\ell$, the symmetric convention gives
\[
\partial_{k,\ell}H(\sigma)=\beta\,\sigma_k\sigma_\ell.
\]
Hence for any observable $F$,
\[
\partial_{k,\ell}\langle F\rangle
=\big\langle (F-\langle F\rangle)\,\partial_{k,\ell}H\big\rangle
=\beta\,\big\langle (F-\langle F\rangle)\,\widetilde{\sigma_{k\ell}}\big\rangle.
\]
Taking $F=\sigma_i$ yields $\partial_{k,\ell}\langle\sigma_i\rangle=\beta\langle\td\sigma_i\,\widetilde{\sigma_{k\ell}}\rangle$,
and the case $\langle\sigma_i\rangle^{k_i}$ follows by the chain rule (noting that for $k_i\in\{1,2\}$,
$k_i\langle\sigma_i\rangle^{k_i-1}=(2\langle\sigma_i\rangle)^{k_i-1}$).  Taking $F=\sigma_i\sigma_j$
gives the stated centered formula for $\partial_{k,\ell}\langle\sigma_i\sigma_j\rangle$, and similarly for $F=\sigma_i\sigma_\ell$.
Finally, for diagonal entries $A_{kk}$ one has $\partial_{k,k}H(\sigma)=\frac{\beta}{2}\sigma_k^2=\frac{\beta}{2}$, which is constant in $\sigma$,
so $\partial_{k,k}\langle F\rangle=\langle(F-\langle F\rangle)\,\frac{\beta}{2}\rangle=0$ for all $F$.
\end{proof}

Extending the two propositions above allows for the computation of the $4$-tensor which gives the Hessian of the Gibbs covariance matrix for a planted (and tilted) SK model. This turns out to be of crucial importance in evaluating the symmetries that help bound the terms that show up in the proof of~\pref{thm:covar-estimate}.

We now state a proposition that explicitly writes the magnetization of a coordinate $i \in [n]$ under the Gibbs measure $G(A,tx_0 + B_t,\beta)$ for a \emph{fixed} instance of the disorder.

\begin{lemma}[Magnetization under the tilted SK model]\label{lem:magnetization-tilted-SK}
    For every $i \in [n]$, a \emph{fixed} choice of $A \in \R^{n \times n}$, and 
    \[
        y_t = tx_0 + B_t\, ,
    \]
    where $B_t \sim \calN(0,t)$ independently for every $t \in \R_{\ge 0}$, the following holds
    \[
        (m_t)_i := \an{\sigma_i}_{G(A,y_t,\beta)} = \an{\tanh\left(\frac{\beta}{\sqrt{n}}\sum_{j \ne i}A_{ij}\sigma_j + t(x_0)_i + (B_t)_i\right)} \in_{\text{a.s.}} (-1,1)\,.
    \]
\end{lemma}
\begin{proof}
    The proof proceeds by direct computation.
    \allowdisplaybreaks
    \begin{align*}
        &\an{\sigma_i}_{G(A,y_t,\beta)} = \sum_{\sigma}\sigma_i \frac{e^{\frac{\beta}{\sqrt{n}}\sum_{k < \ell}A_{k\ell}\sigma_k\sigma_\ell + t\iprod{x_0,\sigma} + \iprod{B_t,\sigma}}}{\sum_{\tau}e^{\frac{\beta}{\sqrt{n}}\sum_{k < \ell}A_{k\ell}\tau_k\tau_\ell + t\iprod{x_0,\tau} + \iprod{B_t,\tau}}} \\
        &= \sum_{\sigma_{\setminus\{i\}}}\sum_{\sigma_i}\sigma_i\frac{e^{\frac{\beta}{\sqrt{n}}\sum_{k < \ell, k \ne i, \ell \ne i}A_{k\ell}\sigma_k\sigma_\ell + t\iprod{(x_0)_{\setminus\{i\}},\sigma_{\setminus\left\{i\right\}}} + \iprod{(B_t)_{\setminus\left\{i\right\}},\sigma_{\setminus\left\{i\right\}}} + \sigma_i\sum_{j\ne i}A_{ij\sigma_j} + t(x_0)_i\sigma_i + (B_t)_i\sigma_i}}{\sum_{\tau_{\setminus\{i\}}}\sum_{\tau_i} e^{\frac{\beta}{\sqrt{n}}\sum_{k < \ell, k \ne i, \ell \ne i}A_{k\ell}\tau_k\tau_\ell + t\iprod{(x_0)_{\setminus\{i\}},\tau_{\setminus\left\{i\right\}}} + \iprod{(B_t)_{\setminus\left\{i\right\}},\tau_{\setminus\left\{i\right\}}} + \tau_i\sum_{j\ne i}A_{ij\tau_j} + t(x_0)_i\tau_i + (B_t)_i\tau_i}} \\
        &= \sum_{\sigma_{\setminus\left\{i\right\}}}\left(\frac{e^{\frac{\beta}{\sqrt{n}}\sum_{k < \ell, k \ne i, \ell \ne i}A_{k\ell}\sigma_k\sigma_\ell + t\iprod{(x_0)_{\setminus\{i\}},\sigma_{\setminus\left\{i\right\}}} + \iprod{(B_t)_{\setminus\left\{i\right\}},\sigma_{\setminus\left\{i\right\}}}}}{\sum_{\tau_{\setminus\left\{i\right\}}}e^{\frac{\beta}{\sqrt{n}}\sum_{k < \ell, k \ne i, \ell \ne i}A_{k\ell}\tau_k\tau\ell + t\iprod{(x_0)_{\setminus\{i\}},\tau_{\setminus\left\{i\right\}}} + \iprod{(B_t)_{\setminus\left\{i\right\}},\tau_{\setminus\left\{i\right\}}}}2\cosh\left(\frac{\beta}{\sqrt{n}}\sum_{i\ne j}A_{ij}\tau_j + t(x_0)_i + (B_t)_i\right)}\right) \\
        &\quad\quad\,\,\,\,\,\left(2\sinh\left(\frac{\beta}{\sqrt{n}}\sum_{i\ne j}A_{ij}\sigma_j + t(x_0)_i + (B_t)_i\right)\right) \\
        &= \an{\tanh\left(\frac{\beta}{\sqrt{n}}\sum_{j \ne i}A_{ij}\sigma_j + t(x_0)_i + (B_t)_i\right)}\, .\qedhere
    \end{align*}
\end{proof}

\paragraph{Gaussian and sub-Gaussian processes} \label{subsec:guassian-lemmata} To make explicit computations about the covariances under the Gibbs measure of the tilted model will require the following mild generalizations of Stein-type lemmata.

\begin{proposition}[Stein's lemma for uncentered Gaussian processes]\label{prop:steins-uncentered}
    Let $f: \R^n \to \R \in C^1$ be a function with slower than exponential growth.
    Then, for every $i \in [n]$,
    \[
        \E_{g \sim \calN(\mu,\Sigma)}\left[(g_i-\mu_i)f(g)\right]=\sum_{j=1}^n\Sigma_{ij}\E\left[\partial_{g_j}f(g)\right]\,.
    \]
\end{proposition}

\begin{proposition}[Stein's lemma with Gibbs averages~{\cite[Mild generalization of Lemma 1.1]{panchenko2013sherrington}}]\label{prop:steins-gibbs-uncentered}
    Let $f: \R \to \R \in C^1$ be a function with slower than exponential growth. Let $\{X(\sigma)\}_{\sigma \in \{-1,1\}^n}$ and $\{Y(\sigma)\}_{\sigma \in \{-1,1\}^n}$ be (possibly uncentered) Gaussian processes with a Gibbs measure $G_{n,\beta,Y}$ defined over the process $Y$ at inverse temperature $\beta$.

    Then,
    \begin{align*}
        &\E\left[\left\langle f\big(X(\sigma) - \E[X(\sigma)]\big)\right\rangle_{G_{n,\beta,Y}}\right] = \E\left[\left\langle f\big(\E\left[X(\sigma)Y(\sigma)\right] - \E\left[X(\sigma)Y(\sigma')\right]\big)\right\rangle_{G^{\ot 2}_{n,\beta,Y}}\right]\,. 
    \end{align*}
\end{proposition}

\paragraph{Guerra-type interpolations, overlaps and the cavity method}
To obtain overlap-concentration for the planted SK model at high-temperature, interpolating all spins \emph{simultaneously} via a Guerra interpolation is critical.

\begin{definition}[Guerra-type interpolation between two independent Gaussian processes]\label{def:guerra-interpolation}
    Let $\{H_0(\sigma)\}_{\sigma \in \Sigma}$ and $\{H_1(\sigma)\}_{\sigma \in \Sigma}$ be two independent Gaussian processes. Then, the Guerra-interpolated Gaussian process $\{H_t(\sigma)\}_{\sigma \in \Sigma}$ is given as
    \begin{equation}
        H_t(\sigma) = \sqrt{1-t}H_0(\sigma) + \sqrt{t}H_1(\sigma)\, ,        
    \end{equation}
    for $t \in [0,1]$.
\end{definition}

See~\cite[\S 5]{alaoui2017finite} for the concrete instantiation of a Guerra interpolation  for the planted SK model. Noting that $\{H_t(\sigma)\}_{\sigma \in \Sigma}$ is a Gaussian process for every $t \in [0,1]$ and using~\pref{prop:steins-gibbs-uncentered} leads to control of the time-derivative when interpolating between the free-energies of two Gaussian processes.

\begin{lemma}[Time-derivative of $C^2$ functions in a Guerra-type interpolation,~{\cite[Lemma 1.3.1]{talagrand2010mean}}]\label{lem:time-derivative-guerra-interpolation}
    Let $\{H_t(\sigma)\}_{\sigma \in \Sigma}$ be a Guerra-interpolated Gaussian process as defined in~\pref{def:guerra-interpolation}. Then, given a function $f: H_t(\cdot) \to \R$, the following holds
    \[
        \fc{d}{dt}\E\left[f(H_t(\cdot))\right] = \frac{1}{2}\sum_{\sigma,\tau \in \Sigma}\left(\E\left[H_1(\sigma)H_1(\tau)\right]-\E\left[H_0(\sigma)H_0(\tau)\right]\right)\E\left[\partial_{x_\sigma}\partial_{x_\tau}f(H_t(x))\right]\,.
    \]
\end{lemma}
\pref{lem:time-derivative-guerra-interpolation} is used to derive the free energy for the planted SK model in the replica-symmetric regime in \cite{alaoui2017finite}. The replica-symmetric regime for the planted SK model is characterized by magnetization and overlap concentration.

\begin{definition}[Overlap and magnetization moment generation functions]
    Let $\{H(\sigma)\}_{\sigma \in \Sigma}$ be a Gaussian process. Then, given a fixed $\beta > 0$ and any \emph{fixed} $k \in \Z_+$, the moments of the deviations of the magnetization and overlap parameters of this process are given as
    \begin{equation}\label{eq:magnetization-mgf}
        \E_H\left[\an{\left(m - R_{\sigma,\vec{1}}\right)^{2k}}\right]\, ,
    \end{equation}
    \begin{equation}\label{eq:overlap-mgf}
        \E_H\left[\an{\left(q - R_{1,2}\right)^{2k}}\right]\,.
    \end{equation}
\end{definition}

These quantities are bounded for the planted SK model when $\beta < 1$ via an analysis of some fixed-point equations, followed by an appeal to the Nishimori identity~(\pref{lem:nishimori-condition}) and a ``mapping'' of our model to the Guerra-interpolated model in~\cite[\S 3]{alaoui2017finite}.

Studying the precise limiting values of functions of the overlap arrays requires computing the \emph{exact} limits for some of its statistics. This is accomplished via a cavity interpolation.

\begin{definition}[Cavity interpolation for the planted SK model,~{Modification of \cite[\S2]{cadel2010sherrington}}]\label{def:cavity-interpolation-planted-sk}
    Consider the Hamiltonian of the planted model (with the functions of the last site isolated):
    \[
        \beta H(\sigma) = \sqrt{\frac{(n-1)}{n}}\beta H_{n-1}(\rho) + \frac{\beta^2}{2n} + \frac{\beta}{2}A_{nn} + \sigma_n\left(g(\rho) + (y_t)_n\right)\, , 
    \]
    where $\rho = \left(\sigma_1,\dots,\sigma_{n-1}\right)$, $(y_t)_n \sim \calN(t,\beta^2 t)$, $A_{nn}\sim\calN\left(0,\frac{2}{n}\right)$ and 
    \[
        g(\rho) := \frac{\beta^2}{n}\sum_{i=1}^{n-1}\sigma_i + \beta\sum_{i=1}^{n-1}A_{i,n}\sigma_i\,.
    \]
    The interpolation that smoothly decouples a $(n-1)$-dimensional configuration $\rho$ and a $n$-dimensional configuration $\sigma = \rho\cdot\sigma_n$ for the Hamiltonian $H_t(\sigma)$ with magnetization $m \in [-1,1]$ and overlap $q \in [0,1]$ is defined as
    \begin{align}
        \beta H_t(\sigma) := &\sqrt{\frac{(n-1)}{n}}\beta H_{n-1}(\rho) + \frac{\beta^2}{2n} + \frac{\beta}{2}A_{nn} +\left\{g_t(\rho) +(y_t)_n\right\}\sigma_n \nonumber\, ,
    \end{align}
    where
    \begin{equation}
        g_t(\rho) := \beta^2\left(\frac{t}{n}\sum_{i=1}^{n-1}\sigma_i+ (1-t)m\right)  + \beta\left(\sqrt{t}\sum_{i=1}^{n-1}A_{in}\sigma_i + \sqrt{1-t}\sqrt{q}\,h\right)\, ,
    \end{equation}
    for $t \in [0,1]$ and $h \sim \calN(0,1)$.
\end{definition}

This cavity interpolation in~\pref{def:cavity-interpolation-planted-sk} is used to extend the ideas in~\cite[\S1.6 \& \S1.8]{talagrand2010mean} to obtain various limiting estimates for the planted SK model (\pref{sec:cavity-interpolation}).

\subsection{Covariance estimation error via moments of overlaps of the planted Gibbs measure} The main theorem of this section is a bound on the Frobenius norm that measures the error between the ideal process, driven by the ideal covariance, and the algorithmic process, driven by the covariance obtained as the inverse Hessian of the TAP free energy for the (planted) SK model under SL tilt.

\begin{theorem}[Frobenius norm error between PHA \& ASL-tilted SK covariance matrices]\label{thm:covar-estimate}
        Let $D(x) = \diag\left(\frac{1}{1-x^2}\right)$ for any $x \in (-1,1)^n$ and 
        \[
            \hat{Q}(\mg_t) = \left(\beta^2\left(1-\frac{\norm{\mg_t}^2_2}{n}\right)\Id_n - \beta A - \frac{\beta^2}{n}x_0x_0^\sT + D(\mg_t) - \frac{2\beta^2}{n}\mg_t\mg_t^\sT\right)^{-1}\, ,
        \]
        where $t \in [0,T]$, and 
        \[
            \mg_t = \an{\sigma}\,.
        \] 
        Fix $x_0 \in \{-1,1\}^n$ and let 
        \[
            P=P(A,y_t) = \Cov(G_{A,y_t,\beta})\,.
        \]
        Then,
        \begin{equation}
            \E_{A,B_t}\norm{\hat{Q}^{-1}P - \Id_n}_F^2 = O_{\beta,t}(1)\,.
        \end{equation}
\end{theorem}
\begin{prf}
    Set $S= \beta^2(1-\norm{\mg_t}^2_2/n)\Id_n + D(\mg_t) - \frac{2\beta^2}{n}\mg_t\mg_t^\sT$. Then, by simple linear algebra,
    \begin{align*}
        &\E_{A,B_t}\norm{\hat Q^{-1} P -\Id_n}_F^2 = \E_{A,B_t}\left[\norm{\hat Q^{-1}P}_F^2\right]- 2\E_{A,B_t}\left[\Tr(\hat Q^{-1}P)\right] + n = \E\left[\Tr[P(\hQ^{-1})^2P]\right] - 2\E\left[\Tr\left[\hQ^{-1}P\right]\right] + n \\
        &=\E\left[\Tr\left[P\hQ^{-2}P\right]\right]  +2\beta\,\E\Tr[AP] + \frac{2\beta^2}{n}\E\Tr[x_0x_0^\sT P] -2\E\left[\Tr\left[SP\right]\right] + n \\
        &= \E\left[\Tr\left[PS^2P\right]\right] + \beta^2\E\left[\Tr\left[PA^2P\right]\right] -2\beta\E\left[\Tr\left[PSAP\right]\right] + \frac{\beta^4}{n^2}\E[\Tr[P(x_0x_0^\sT)^2P]] - \frac{2\beta^2}{n}\E[PS(x_0x_0^\sT)P] \\
        &+ \frac{2\beta^3}{n}\E\left[PA(x_0x_0^\sT)P\right] + 2\beta\,\E\Tr\left[AP\right] + \frac{2\beta^2}{n}\E\Tr[x_0x_0^\sT P] -2\E\left[\Tr\left[SP\right]\right] + n \\
        &= \underbrace{\E\left[\Tr\left[PS^2P\right]\right] -2\E\left[\Tr\left[SP\right]\right] + n}_{=\calE_D} + \underbrace{ \beta^2\E\left[\Tr\left[PA^2P\right]\right] -2\beta\E\left[\Tr\left[PSAP\right]\right] + 2\beta\,\E\Tr(AP)}_{\calE_{A}} \\
        & + \underbrace{\frac{\beta^4}{n^2}\E[\Tr[P(x_0x_0^\sT)^2P]] - \frac{2\beta^2}{n}\E[PS(x_0x_0^\sT)P] + \frac{2\beta^3}{n}\E\left[PA(x_0x_0^\sT)P\right] +  \frac{2\beta^2}{n}\E\Tr[x_0x_0^\sT P]}_{\calE_R}\, .
    \end{align*}
    $\calE_D$, $\calE_A$ and $\calE_R$ are reduced to Gibbs averages of products of overlaps between two replicas. However, for $\calE_A$, one needs various applications of Stein's lemma to further the simplification. After this, various cavity estimates help show that the $O(1/n)$ scaling terms \emph{exactly} cancel to obtain $\calE_A + \calE_D = O_{\beta,q^*}(1)$, while various ``remainder'' terms (higher-order moments and diagonally-weighted cumulants) can be handled via a mixture of cavity estimtates and concentration. The rank-$1$ terms make no important contribution as $\calE_R = O(1)$. The final bound follows from~\pref{cor:final-frobenius-bound}.
\end{prf}

These terms are individually simplified using replica identities for the Gibbs measures leading to various cancellations. The surviving terms reduce to understanding products of overlaps under the Gibbs measure for the planted SK model. Control of the even moments of the overlap deviations for the planted SK model is extracted from~\cite[\S 3]{alaoui2017finite}. This allows one to conclude that the deviations around the overlap are sufficiently small, and the net error in the estimation is $O_{\beta,t}(1)$ for the ``lower order'' terms. Arranging for the cancellations of the ``higher order'' terms requires precise limiting estimates for various mixed moments of the overlaps, which are established by using a modified cavity interpolation for the planted SK model in~\pref{sec:cavity-interpolation}. This approach can be seen as a substantial and technically involved generalization of the one adopted in~\cite[\S3]{el2024bounds}.

In the following proofs, we will use notions of ``rectangular sums'' ($f, \hat{f}$), rescaled-and-recentered spins ($\hat{\sigma}_i$), cubic ``bulk'' terms such as $\nu_1(f G^-)$ where $G^- \in \calQ$, and the limits for these denoted as $U_n, V_n$ and $W_n$, as well as various other objects (and the pertinent notation) introduced in~\pref{sec:cavity-interpolations-trace-identities}.

\begin{lemma}[Controlling the $\calE_D$ term]\label{lem:diag-frob-term}
    Let $G=G_{A,y_t,\beta}$ be the (tilted/planted) SK Gibbs measure at $t \in [0,T]$, and
    \[
        m_t := \an{\sigma}\in(-1,1)^n,
        \qquad 
        P:=\Cov(G),
        \qquad 
        P_{ij}=\an{\sigma_i\sigma_j}-\an{\sigma_i}\an{\sigma_j}\,.
    \]
    Then,
    \begin{align*}
        \calE_D &= \frac{n^2}{4}\left(W^{(2)}_n - \frac{4}{n}\right) + \frac{n^2}{4}\E\left[c^2\an{f^2}\right] + \frac{n^2}{2}\E\left[c\an{f\hat{f}}\right] - 2\E\left[c\Tr[P]\right] + \calR^V_D \\
        &= \beta^2\E\left[\Tr[P^2]\right] + \E\left[c^2\Tr[P^2]\right] + 2\left(\E\left[c\Tr[PDP]\right] - \E\left[c\Tr[P]\right]\right) + \calR^V_D + O(1) \\
        &= \beta^2 \Tr[P^2] + 3\mathbb{E}\left[c^2\Tr[P^2]\right] + \calR^V_D + O(1)\,.
    \end{align*}
    where $W^{(2)}_n$ is defined in \pref{sec:cavity-interpolation} and $\calR^D_V$ is defined in the proof. 
\end{lemma}
\begin{prf}
    Denote $m:=m_t=\an{\sigma}$ and define the rescaled and re-centered overlaps as $\hat{\sigma_i} := \frac{\sigma_i - m_i}{\sqrt{1-m_i^2}}$. Observe that
    \[
        P_{ij} = \an{\sigma_i\sigma_j} - \an{\sigma_i}\an{\sigma_j} = \frac{1}{2}\an{\left(\sigma^1_i-\sigma^2_i\right)\left(\sigma^1_j-\sigma^2_j\right)}\,.    
    \]
    Notate $V := \frac{2\beta^2}{n}\mg\mg^\sT$ and $c = \beta^2(1-\norm{\mg}^2_2/n)$. The goal is to separate the terms into those \emph{with} the rank-1 projector $V$ and those without them. The former will be shown to be small by using concentration of cubic ``bulk'' observables established in~\pref{prop:D15} and~\pref{lem:D16-corrected}. The latter will be estimated using replica identities and the CLT estimates developed in \pref{sec:cavity-interpolation} and made to cancel with some terms in $\calE_A$.
    
    Expanding the trace and using its linearity,
    \allowdisplaybreaks
    \begin{align*}
        \calE_D &= \E\left[\Tr[PS^2P]-2\Tr[SP] + n\right] = \E\left[\Tr[P(D + c\Id_n - V)^2P]\right] - 2\E\left[\Tr[(D + c\Id_n -V)P]\right] + n \\
        &= \E\left[\Tr[PD^2P]\right] + 2\E[c\Tr[PDP]] - \E[\Tr[P(VD + DV)P]] - 2\E[c\Tr[PVP]] + \E[c^2\Tr[P^2]] \\
        &\quad + \E[\Tr[PV^2P]] - 2\E[\Tr[DP]] - 2\E[c\Tr[P]] + 2\E[\Tr[VP]] + n \\
        &= \left(\E\left[\Tr[PD^2P]\right] - 2\E[Tr[DP]] + n\right) + \E\left[c^2\Tr[P^2]\right] + 2\E[c\left(\Tr[PDP]-\Tr[P]\right)] + \calR^V_{D}\, ,
    \end{align*}
    where $\calR^V_D$ contains all tracial terms that have at least one copy of the rank-$1$ projector $V$.
    \ppart{Evaluating $\Tr[DP]$}
    Note that $P_{ii}=1-m_i^2$ and $D(m)_{ii}=\frac{1}{1-m_i^2}$, which implies
    \[
        \Tr\left[D P\right]=\sum_{i=1}^n \frac{1}{1-m_i^2}(1-m_i^2)=n\,.
    \]
    \ppart{Evaluating $\Tr[P^2]$} Using the fact that $f = \frac{1}{n}\sum_{i=1}^n(\sigma^1_i - \sigma^2_i)(\sigma^3_i-\sigma^4_i)$, $\Tr[P^2] = \norm{P}^2_F$ and standard replica identities, one obtains the following algebraic simplification via~\pref{lem:trace-to-replica},~\eqref{eq:TrP2-f2} as
    \begin{align*}
        \Tr[P^2] &=\frac{n^2}{4}\an{f^2}\,.
    \end{align*}
    \ppart{Evaluating $\Tr[PDP]$} By a similar calculation using~\pref{lem:trace-to-replica}, one obtains
    \begin{align*}
        \Tr[PDP] &= \frac{n^2}{4}\an{f\hat{f}}\,.
    \end{align*}
    \ppart{Evaluating $\Tr[PD^2P]$} By yet another analogous calculation, this yields the ``doubly-scaled'' quadratic term,
    \[
        \Tr[PD^2P] = \frac{n^2}{4}\an{ff^{(2)}}\,.
    \]
    Taking expectation, and setting $W^{(2)}_n = \nu_1(f\,f^{(2)})$ along with an application of~\pref{prop:D18} gives the $\beta^2\E[\Tr[P^2]]$ in $\calE_D$. The final step is applying~\pref{lem:cTrPDP-identity} to the $2\E\left[c(\Tr[PDP] - \Tr[P])\right]$ term to obtain a factor of $2\E\left[c^2\Tr[P^2]\right]$. \qedhere
\end{prf}

We now dispense with the $\calR^V_D$ term first, by using a mixture of overlap concentration and control over the moments of certain quadratic bulk observables, showing that it is $O(1)$.

\begin{lemma}[$\calR^V_D = O(1)$]\label{lem:bound-d-remainder}
    Let $\calR^V_D$ be as defined in~\pref{lem:diag-frob-term}. Then,
    \[
        \calR^V_D = O(1)\,.
    \]
\end{lemma}
\begin{prf}
    From the proof of~\pref{lem:diag-frob-term}, note that
    \[
        \calR^V_D = \E\left[\Tr[PV^2P] - \Tr[P(VD + DV)P] -2c\Tr[PVP] + 2\Tr[VP]\right]\,.
    \]
    Since $V = \frac{2\beta^2}{n}mm^\sT$, it is immediate that
    \[
        \Tr[PVP] = \frac{2\beta^2}{n}\mg^\sT P^2\mg,\quad \Tr[VP] = \frac{2\beta^2}{n}\mg^\sT P\mg,\quad \Tr[PV^2P] = \frac{4\beta^4}{n^2}\norm{\mg}^2_2\mg^\sT P^2\mg\,. 
    \]
    Similarly,
    \[
        \Tr[PDVP] = \frac{2\beta^2}{n}\mg^\sT P^2(D\mg) = \Tr[PVDP]\,.
    \]
    The quadratic forms are scaled by powers of $1/n$ and $\norm{\mg}^2_2 \le n$, so it will suffice to show that
    \[
        \E\langle\mg, P\mg\rangle = O(n),\quad \E\langle \mg, P^2\mg\rangle = O(n),\quad \E\langle\mg, P^2(D\mg)\rangle = O(n)\,.  
    \]
    \ppart{$\E[\Tr[VP]] = O(1)$ via overlap-concentration} This is the easy term, and follows directly from overlap concentration after applying replica identities to rewrite the quadratic form as a straightforward quadratic polynomial of replica overlaps. Observe the following
    \allowdisplaybreaks
    \begin{align*}
        \Tr[VP] &= \frac{2\beta^2}{n}\E\left[ \mg^\sT P \mg\right] = \sum_{i,j}P_{i,g}\mg_i\mg_j = \sum_{i,j}\left(\an{\sigma_i\sigma_j}-\an{\sigma_i}\an{\sigma_j}\right)\an{\sigma_i}\an{\sigma_j} \\
        &= \frac{2\beta^2}{n}\sum_{i,j}\an{\sigma_i\sigma_j}\an{\sigma_i}\an{\sigma_j} - \sum_{i,j}\an{\sigma_i}^2\an{\sigma_j}^2 \\
        &= \frac{2\beta^2}{n}n^2\left(\an{R_{1,2}R_{1,3}} - \an{R_{1,2}}^2\right) = 2\beta^2n\left(\an{(R_{1,2}-q^*)(R_{1,3}-q^*)} - \an{(R_{1,2}-q^*)}^2\right)\,.
    \end{align*}
    Taking expectations $\E[\cdot]$, Cauchy-Schwarz, Jensen, and invoking overlap-concentration yields
    \begin{align*}
        2\beta^2 n\E\left[\left|\an{R_{1,2}R_{1,3}}- \an{R_{1,2}}^2\right|\right] &\le 2\beta^2n\E\left[\left|\an{R_{1,2}-q^*}\right|\left|\an{R_{1,3}-q^*}\right|\right] + 2\beta^2n\E\left[\an{(R_{1,2}-q^*)}^2\right] \\
        &\le 2\beta^2n\sqrt{\E\left[\an{(R_{1,2}-q^*)^2}\right]}\sqrt{\E\left[\an{(R_{1,3}-q^*)^2}\right]} + 2\beta^2n\E\left[\an{(R_{1,2}-q^*)^2}\right] \\
        &\le 2\beta^2 n\frac{C}{\sqrt{n}}\frac{C}{\sqrt{n}} + 2\beta^2 \frac{C}{n} = O(1)\,.
    \end{align*}
    This immediately implies that
    \[
        \E[\Tr[VP]] = O_\beta(1)\,.
    \]

    \ppart{$\E[\Tr[PVP]]$ and $\E[\Tr[PV^2P]]$ are $O_\beta(1)$ via~\pref{prop:D15}} We expand the trace and use replica identities to write the term as a product of centered overlaps, weighted by the ``rectangular sum'' $f$, which can shown to be $O(n^{-2})$ by~\pref{prop:D15}.
    \allowdisplaybreaks
    \begin{align*}
        \Tr[PVP] &= \frac{2\beta^2}{n}\mg^\sT P^2 \mg = \frac{2\beta^2}{n}\sum_{i,j,k}\mg_i P_{ik}P_{kj} \mg_j \\
        &= \frac{\beta^2}{2n}\sum_{i,j,k}\an{\sigma^5_i}\an{(\sigma^1_i - \sigma^2_i)(\sigma^1_k-\sigma^2_k)}\an{(\sigma^3_j-\sigma^4_j)(\sigma^3_k - \sigma^4_k)}\an{\sigma^6_j} \\
        &= \frac{\beta^2}{2n}\left(\sum_i\an{\sigma^5_i(\sigma^1_i-\sigma^2_i)}\right)\left(\sum_j\an{\sigma^6_j(\sigma^3_j - \sigma^4_j)}\right)\left(\sum_k\an{(\sigma^1_k-\sigma^2_k)(\sigma^3_k-\sigma^4_k)}\right) \\
        &= \frac{\beta^2 n^2}{2}\an{(R_{1,5}-R_{2,5})(R_{3,6} - R_{3,4})f} =  \frac{\beta^2 n^2}{2}\an{(\hat{R}_{1,5}-\hat{R}_{2,5})(\hat{R}_{3,6} - \hat{R}_{3,4})f}\\
        &= \frac{\beta^2 n^2}{2}\left(\an{f\hat{R}_{1,5}\hat{R}_{3,6}} - \an{f\hat{R}_{1,5}\hat{R}_{3,4}} - \an{f\hat{R}_{2,5}\hat{R}_{3,6}} + \an{f\hat{R}_{2,5}\hat{R}_{3,4}}\right)\,.
    \end{align*}
    Now, taking expectations $\E[\cdot]$ and applying~\pref{prop:D15} to each of the four terms above yields that
    \[
        \E\left[\Tr[PVP]\right] = \frac{\beta^2 n^2}{2}\left(\frac{C_1}{n^2} + \frac{C_2}{n^2} + \frac{C_3}{n^2} + \frac{C_4}{n^2}\right) = O_\beta(1)\,,
    \]
    since $\nu_1(f G^-) = Cn^{-2}$ for any quadratic bulk observable $G^-$ such as $\hat{R}^-_{\ell,\ell'}\hat{R}^-_{r,r'}$~\eqref{eq:fminus-replace}, and
    \allowdisplaybreaks
    \begin{align*}
        \nu(f\hat{R}_{\ell,\ell'}\hat{R}_{r,r'}) &= \underbrace{\nu(fG^-)}_{= Cn^{-2}\text{~\eqref{eq:fminus-replace}}} + \underbrace{\frac{1}{n^2}\nu(f(\eps^\ell\eps^{\ell'}\eps^r\eps^{r'}))}_{\le C'n^{-2}} + \frac{1}{n}\nu(f\eps^\ell\eps^{\ell'}\hat{R}^-_{r,r'}) + \frac{1}{n}\nu(f\eps^r\eps^{r'}\hat{R}^-_{\ell,\ell'})  \\
        &\le_{\text{Cauchy-Schwarz}} Cn^{-2} + \frac{2}{n}\sqrt{\nu(f^2)}\sqrt{\nu((\hat{R}^-_{\ell,\ell'})^2)}\\
        &=_{\text{\pref{lem:replicon-moments-from-mgf}\,+\,overlap-concentration}} Cn^{-2}\,.
    \end{align*}

    Similarly,
    \[
        \E[\Tr[PV^2P]] = \E\left[\frac{4\beta^4}{n^2}\norm{\mg}^2_2\braket{\mg | P^2\mg}\right] \le_{\norm{\mg}^2_2 \le n} \frac{4\beta^4}{n}\E\left[\Tr[PVP]\right] \le O_\beta(1)\,. 
    \]
    \ppart{$\E[\Tr[PDVP]] = \E[\Tr[PVDP]] = O_\beta(1)$ via \pref{prop:D15}} Since $\Tr[PVDP] = \Tr[PDVP] = \frac{2\beta^2}{n}\mg^\sT P^2 (D\mg)$, it suffices to bound the quadratic form $\E[\mg^\sT P^2(D\mg)] = \Theta(n)$. To this end, note that
    \allowdisplaybreaks
    \begin{align*}
        &\E\left[\mg^\sT P^2 (D\mg)\right] = \E\left[\sum_{i,j} P^2_{i,j}\mg_i D_{jj}\mg_j\right] = \sum_{i,j}\mg_i\mg_jD_{jj}\sum_k P_{ik}P_{kj} \\
        &= \frac{1}{4}\E\left[\sum_{ij}\mg_i\mg_jD_{jj}\frac{1}{4}\sum_k\left(\an{(\sigma^1_i-\sigma^2_i)(\sigma^1_k-\sigma^2_k)}\right)\left(\an{(\sigma^3_k-\sigma^4_k)(\sigma^3_j-\sigma^4_j)}\right)\right]\\
        &= \frac{1}{4}\E\left[\an{\left(\sum_i m_i(\sigma^1_i-\sigma^2_i)\right)\left(\sum_j m_j D_{jj}(\sigma^3_j - \sigma^4_j)\right)\left(\sum_k(\sigma^1_k-\sigma^2_k)(\sigma^3_k-\sigma^4_k)\right)}\right]\\
        &=\frac{n^2}{4}\E\left[\an{\left(\sum_j m_j D_{jj}(\sigma^3_j - \sigma^4_j)\right)\left(\frac{1}{n}\sum_i \sigma^1_i - m^*\right)\,f}\right] - \frac{n^2}{4}\E\left[\an{\left(\sum_j m_j D_{jj}(\sigma^3_j - \sigma^4_j)\right)\left(\frac{1}{n}\sum_i \sigma^2_i - m^*\right)\,f}\right]\,.
    \end{align*}
    Using the fact that $|m_j| \le 1$ uniformly for every $j \in [n]$, and noting that the three terms in the summand are either a bulk-deviation (like magnetization deviation), the rectangular sum ($f$), or a diagonally-weighted summand, means a direct invocation of~\pref{lem:D-oneleg-O(n)} shows the above is $O(n)$ which concludes the bound\footnote{~See~\pref{rem:boundedness-handling}.}. \qedhere
\end{prf}

Below, we bound the $\calA$ term using a similar approach. The computation and simplifications are longer and more tedious than those involved in controlling the $\calA$ term, but the changes are not fundamental.


\begin{lemma}[Controlling the $\mathcal E_A$ term and exact cancellations]\label{lem:EA-term}
Let $G=G_{A,y_t,\beta}$ be the Gibbs measure of the planted SK model with tilt $y_t$ at time $t\in[0,T]$, with
\[
m:=\langle\sigma\rangle\in(-1,1)^n,
\qquad
P:=\Cov(G),\qquad P_{ij}=\langle\sigma_i\sigma_j\rangle-\langle\sigma_i\rangle\langle\sigma_j\rangle.
\]
Let
\[
q:=\frac{\|m\|_2^2}{n},
\qquad
c:=\beta^2(1-q) = \frac{\beta^2}{n}\Tr[P],
\qquad
V:=\frac{2\beta^2}{n}mm^{\sT},
\qquad
S:=D+cI_n-V.
\]
Then,
\begin{equation}\label{eq:EA-master}
\mathcal E_A
= -\beta^2\E[\Tr[P^2]] - 3\E\left[c^2\Tr[P^2]\right] + O(1)
+\mathcal R_A^{V}
+\mathcal R_A^{\mathrm{resp}},
\end{equation}
where
\begin{itemize}[itemsep=0.25em]
\item $\mathcal R_A^{V}$ is the collection of all contributions containing at least one copy of $V$, and%
\item $\mathcal R_A^{\mathrm{resp}}$ is the collection of all ``response'' contributions coming from
Gaussian integration-by-parts in $A$ where at least one $A$--derivative hits either $D$, $c$,
or produces a non-diagonal spin-correlation response\,. 
\end{itemize}
\end{lemma}
\begin{proof}
    The proof follows by applying algebraic simplifications after applications of Stein's lemma~(\pref{prop:steins-uncentered}) and explicit evaluations of derivatives of the Gibbs density to obtain trace terms that exactly cancel those in $\mathcal{E}_D$, while collecting the remaining terms into another remainder $\calR^V_A$ which can be shown to be $O(1)$ via overlap/magnetization concentration, as well as the cubic and diagonally-weighted remainder estimates in \pref{sec:cavity-interpolation}.
    \ppart{Steins's lemma-type identities} Write the SK term in the planted SK model Hamiltonian as
    \[
        \frac{\beta}{2}\langle\sigma,A\sigma\rangle = \beta\sum_{1\le i<j\le n}A_{ij}\sigma_i\sigma_j +\frac{\beta}{2}\sum_{i=1}^n A_{ii}\,.
    \]
    The diagonal piece $\frac{\beta}{2}\sum_i A_{ii}$ is $\sigma$-independent, and so $P,m,D,c,V,S$ do \emph{not} depend on $\{A_{ii}\}_{i}$ (\pref{lem:diag-does-not-enter-gibbs}). All applications of~\pref{prop:steins-uncentered} are performed only over $\{A_{ij}\}_{i < j}$. For $i<j$, $A_{ij}\sim \calN(0,1/n)$, and so using~\pref{prop:steins-uncentered} gives
    \begin{equation}\label{eq:GOE-ibp}
    \E\big[A_{ij}\,F(A)\big]
    =
    \frac{1}{n}\,\E\big[\partial_{ij}F(A)\big]\,.
    \end{equation}
    For the quadratic term we apply~\pref{prop:steins-uncentered} in conjunction with the product rule to obtain
    \begin{equation}\label{eq:GOE-ibp-2}
    \E\big[A_{ij}A_{kl}\,F(A)\big]
    =
    \frac{1}{n}\mathbf 1_{\{(i,j)=(k,l)\}}\E[F(A)]
    +\frac{1}{n}\E\big[A_{kl}\,\partial_{ij}F(A)\big],
    \qquad (i<j,\;k<l).
    \end{equation}
    \ppart{Simplifying $\mathcal E_A$ and isolating the rank-$1$ part} Using $S=D+cI-V$ and trace cyclicity,
    \[
    \Tr[PSAP]=\Tr[P(D+cI-V)AP]=\Tr[AP^2D]+c\,\Tr[AP^2]-\Tr[AP^2V]\,.
    \]
    This immediately gives,
    \allowdisplaybreaks
    \begin{align}
    \mathcal E_A
    &=
    \underbrace{\beta^2\,\E\Tr[PA^2 P]
    -2\beta\,\E\Tr[AP^2D]
    -2\beta \,\E c\Tr[AP^2]
    +2\beta\,\E\Tr[AP]}_{:=\,\mathcal{E}^{(0)}_A}
    +\underbrace{2\beta\,\E\Tr[AP^2V]}_{:=\,\mathcal R_A^{V}}.
    \label{eq:EA-expand-1}
    \end{align}
    $\mathcal R_A^{V}$ will contain \emph{all} terms that contain at least one copy of $V$ after performing the Stein's lemma based expansions.
We first evaluate the ``$V$-free'' part, which gives
\[
\mathcal E_A^{(0)}
:=
\beta^2\,\E\Tr[PA^2P]
-2\beta\,\E\Tr[AP^2D]
-2\beta \,\E c\Tr[AP^2]
+2\beta\,\E\Tr[AP]\,.
\]
\ppart{Rewriting Gaussian trace terms as sums with symmetric coefficients}
For any $M \in M_n(\R)$ that is independent of the diagonal of $A$,
\[
\Tr[AM]=\sum_{i<j}A_{ij}\big(M_{ij}+M_{ji}\big).
\]
Since $P$ is symmetric, so is $P^2$, hence $P^2_{ij}=P^2_{ji}$.
Moreover $P^2D$ is not symmetric, but
\[ 
(P^2D)_{ij}+(P^2D)_{ji}=(P^2)_{ij}D_{jj}+(P^2)_{ij}D_{ii}=(P^2)_{ij}(D_{ii}+D_{jj}).
\]
Thus
\allowdisplaybreaks
\begin{align}
\Tr[AP] &= 2\sum_{i<j} A_{ij}P_{ij},\label{eq:AP-symm}\\
\Tr[AP^2] &= 2\sum_{i<j} A_{ij}(P^2)_{ij},\label{eq:AP2-symm}\\
\Tr[AP^2D] &= \sum_{i<j} A_{ij}(P^2)_{ij}(D_{ii}+D_{jj}).\label{eq:AP2D-symm} \\
\Tr[PA^2P]&= \Tr[AP^2A] = \sum_{i,j} A_{ij}(P^2A)_{ij}.
\end{align}
\ppart{Gibbs differentiation identities}
For $i<j$ and any bounded observable $F(\sigma)$,
\begin{equation}\label{eq:gibbs-deriv-basic}
\partial_{ij}\langle F\rangle
=
\beta\Big(\langle F\,\sigma_i\sigma_j\rangle-\langle F\rangle\langle\sigma_i\sigma_j\rangle\Big).
\end{equation}
In particular, with $m_k:=\langle\sigma_k\rangle$ and $P_{ab}:=\langle\sigma_a\sigma_b\rangle-m_am_b$, applying~\pref{prop:derivatives-gibbs-moments} gives
\allowdisplaybreaks
\begin{align}
\partial_{ij}m_a
&=
\beta\Big(\langle\sigma_a\sigma_i\sigma_j\rangle-m_a\langle\sigma_i\sigma_j\rangle\Big),
\label{eq:dm}\\
\partial_{ij}P_{ab}
&=
\beta\Big(\langle (\sigma_a-m_a)(\sigma_b-m_b)\sigma_i\sigma_j\rangle
-P_{ab}\langle\sigma_i\sigma_j\rangle\Big).
\label{eq:dP}
\end{align}
Similarly, using~\pref{prop:derivatives-gibbs-moments} on $D_{aa}=(1-m_a^2)^{-1}$ and $c = \beta^2(1-\norm{\mg}^2_2/n)$ gives
\begin{align}\label{eq:dD}
\partial_{ij}D_{aa}&=\frac{2m_a}{(1-m_a^2)^2}\,\partial_{ij}m_a=2m_aD_{aa}^2\,\partial_{ij}m_a\,, \\
\partial_{ij}c&=-\frac{2\beta^2}{n}\sum_{a=1}^n m_a\,\partial_{ij}m_a.
\end{align}
\ppart{Apply~\pref{prop:steins-uncentered} on $\calE^{(0)}_A$}
We apply \eqref{eq:GOE-ibp} and \eqref{eq:GOE-ibp-2} to \eqref{eq:EA-expand-1},
and we systematically split each resulting derivative into:
\begin{enumerate}[itemsep=0.2em]
\item A \emph{trace-polynomial} part, which is obtained by keeping only the contributions in which
the $A$-derivative produces $\partial_{ij}P_{ij}$-type or $\partial_{ij}(P^2)_{ij}$-type
contributions that close on traces,
\item A \emph{response remainder} part, which collects every other contribution -- those where the derivative hits $D$, $c$, or creates a non-diagonal spin-response such as
$\partial_{ij}P_{ik}$ with $k\notin\{i,j\}$, or $\partial_{ij}(PAP)_{ij}$ through $\partial_{ij}P$ with non-matching indices,
or those in which a second IBP leaves a surviving factor of $A$. As we will see, most of the higher-order remainder terms and trace powers that have $O(1)$ contribution come from these terms.
\end{enumerate}
We simplify each term individually. 
\ppart{The $2\beta\E[\Tr[AP]]$ term} Using \eqref{eq:AP-symm} and \eqref{eq:GOE-ibp},
\[
    \E\Tr[AP] = 2\sum_{i<j}\E[A_{ij}P_{ij}] =\frac{2}{n}\sum_{i<j}\E[\partial_{ij}P_{ij}].
\]
The contribution is obtained by applying \eqref{eq:dP} with $a=i$ and $b=j$ as
\[
    \partial_{ij}P_{ij} = \beta\Big(\langle(\sigma_i-m_i)(\sigma_j-m_j)\sigma_i\sigma_j\rangle -P_{ij}\langle\sigma_i\sigma_j\rangle\Big)
        = \beta\big(P_{ii}P_{jj}-P_{ij}^2\big),
\]
where the last equality follows from direct computation. Using the tracial expansion,
\begin{equation}\label{eq:EA-AP-main}
2\beta \E\Tr[AP] = \frac{4\beta^2}{n}\sum_{i<j}\E\big[P_{ii}P_{jj}-P_{ij}^2\big]\,.
\end{equation}
Note that
\[
    \sum_{i<j}P_{ii}P_{jj}=\frac12\Big(\Tr[P]^2-\sum_i P_{ii}^2\Big),
\qquad
\sum_{i<j}P_{ij}^2=\frac12\Big(\Tr[P^2]-\sum_i P_{ii}^2\Big).
\]
Subtracting and taking $\E[\cdot]$ gives
\[
    2\beta\E[\Tr[AP]] = \E\left[\frac{4\beta^2}{n}\sum_{i<j}(P_{ii}P_{jj}-P_{ij}^2)\right]=\frac{2\beta^2}{n}\left(\E\left[\Tr[P]^2\right]-\E\left[\Tr[P^2]\right]\right) = 2\E[c\Tr[P]] - \frac{2\beta^2}{n}\E[\Tr[P^2]]\,.
\]
The second term will be shown to be $O(1)$ in the proof of \pref{lem:remainder-bound-ea}.

\ppart{The $-2\beta\E\left(c\Tr[AP^2] + \Tr[AP^2D]\right)$ terms} Similar calculations allow for extraction of the ``heavy'' components. 
Using \eqref{eq:AP2D-symm} and \eqref{eq:GOE-ibp},
\[
\E\Tr[AP^2D]
=
\sum_{i<j}\E\Big[A_{ij}(P^2)_{ij}(D_{ii}+D_{jj})\Big]
=
\frac{1}{n}\sum_{i<j}\E\Big[\partial_{ij}\big((P^2)_{ij}(D_{ii}+D_{jj})\big)\Big].
\]
Using the product rule and summing yields,
\allowdisplaybreaks
\begin{align*}
    \E[\Tr[AP^2D]] &= \frac{1}{n}\sum_{i < j}\E[\partial_{ij}\big((P^2)_{ij}(D_{ii}+D_{jj})\big)] \\
    &= \frac{1}{n}\E\left[\sum_{i < j} (\partial_{ij}P^2_{ij})(D_{ii}+D_{jj}) \right] + \underbrace{\frac{1}{n}\E\left[(P^2)_{ij}\big(\partial_{ij}D_{ii}+\partial_{ij}D_{jj}\big)\right]}_{:=\,\mathrm{DW} \in \calR_A^{\mathrm{resp}}}\,.
\end{align*}
Writing $P^2_{ij}=\sum_{k}P_{ik}P_{kj}$, the main contribution is obtained by keeping only the $k=i$ and $k=j$ pieces
(which are the only ones in which $\partial_{ij}$ produces $\partial_{ij}P_{ij}$ rather than a genuinely non-diagonal response).
Specifically,
\[
    \partial_{ij}(P^2)_{ij} = \beta\Big(P_{ii}(P^2)_{jj}+P_{jj}(P^2)_{ii}+2P_{ij}(P^2)_{ij}\Big) + R^6_{i,j,ij}\,,
\]
where $R^6_{i,j,ij}$ once again collects the centered third-moment and fourth cumulant terms (\pref{lem:dP-pairing}).
Substituting gives,
\allowdisplaybreaks
\begin{align}
-2\beta\,\E\Tr(AP^2D)
&=
-\frac{2\beta^2}{n}\sum_{i<j}\E\Big[(D_{ii}+D_{jj})\Big(P_{ii}(P^2)_{jj}+P_{jj}(P^2)_{ii}+2P_{ij}(P^2)_{ij}\Big)\Big]
\label{eq:AP2D-main-sum}
\\
&\quad
-\frac{2\beta}{n}\sum_{i<j}\E\Big[(D_{ii}+D_{jj})\,R^6_{i,j,ij}\Big]
\;+\;\mathrm{DW}.\nonumber
\end{align}
We now rewrite the main sum in traces. Denote
\[
S_1:=\sum_{i<j}(D_{ii}+D_{jj})\big(P_{ii}(P^2)_{jj}+P_{jj}(P^2)_{ii}\big).
\]
Expand $S_1$ into the four relevant pieces as
\allowdisplaybreaks
\begin{align*}
S_1
&=\sum_{i<j}\Big(D_{ii}P_{ii}(P^2)_{jj}+D_{jj}P_{jj}(P^2)_{ii}\Big)
+\sum_{i<j}\Big(D_{ii}P_{jj}(P^2)_{ii}+D_{jj}P_{ii}(P^2)_{jj}\Big).
\end{align*}
Since $D_{ii}P_{ii}=1$,
\[
\sum_{i<j}\Big(D_{ii}P_{ii}(P^2)_{jj}+D_{jj}P_{jj}(P^2)_{ii}\Big)
=\sum_{i<j}\big((P^2)_{jj}+(P^2)_{ii}\big)
=(n-1)\Tr[P^2].
\]
For the “ratio” part, using symmetry under swapping $(i,j)$ gives
\begin{align*}
\sum_{i<j}\Big(D_{ii}P_{jj}(P^2)_{ii}+D_{jj}P_{ii}(P^2)_{jj}\Big)
&=\sum_{i\neq j}D_{ii}P_{jj}(P^2)_{ii}
=\sum_i D_{ii}(P^2)_{ii}\sum_{j\neq i}P_{jj}\\
&=\Tr[P]\sum_i D_{ii}(P^2)_{ii}-\sum_i D_{ii}P_{ii}(P^2)_{ii}\\
&=\Tr[P]\Tr[DP^2]-\Tr[P^2].
\end{align*}
Using $\Tr[DP^2]=\Tr[PDP]$, we get
\begin{equation}\label{eq:S1-traces}
S_1=(n-1)\Tr[P^2]+\Tr[P]\Tr[PDP]-\Tr[P^2]=(n-2)\Tr[P^2]+\Tr[P]\Tr[PDP].
\end{equation}
Similarly, let
\[
S_2:=\sum_{i<j}2(D_{ii}+D_{jj})P_{ij}(P^2)_{ij}.
\]
Again by symmetry,
\[
\sum_{i<j}(D_{ii}+D_{jj})P_{ij}(P^2)_{ij}
=\sum_{i\neq j}D_{ii}P_{ij}(P^2)_{ij}.
\]
Noting the following expansion for the cubic part,
\[
\Tr[DP^3]=\sum_{i,j}D_{ii}P_{ij}(P^2)_{ij}
=\sum_i D_{ii}P_{ii}(P^2)_{ii}+\sum_{i\neq j}D_{ii}P_{ij}(P^2)_{ij}
=\Tr[P^2]+\sum_{i\neq j}D_{ii}P_{ij}(P^2)_{ij}\,,
\]
gives the following final expression
\begin{equation}\label{eq:S2-traces}
S_2
=2\big(\Tr[DP^3]-\Tr[P^2]\big)
=2\Tr[DP^3]-2\Tr[P^2].
\end{equation}
Combining \eqref{eq:S1-traces} and \eqref{eq:S2-traces} yields
\begin{equation}\label{eq:S1S2-combined}
S_1+S_2
=(n-4)\Tr[P^2]+\Tr[P]\Tr[PDP]+2\Tr[DP^3].
\end{equation}
Inserting \eqref{eq:S1S2-combined} into \eqref{eq:AP2D-main-sum} yields that
\allowdisplaybreaks
\begin{align}
-2\beta\,\E\Tr[AP^2D]
&=
-\frac{2\beta^2}{n}\E\Big[(n-4)\Tr[P^2]+\Tr[P]\Tr[PDP]+2\Tr[DP^3]\Big]
\label{eq:AP2D-after-traces}
\\
&\quad
+\underbrace{\Big(
-\frac{2\beta}{n}\sum_{i<j}\E\Big[(D_{ii}+D_{jj})\,R^6_{i,j,ij}\Big]
\;+ \mathrm{DW}
\Big)}_{=:\ \mathcal R_{AP^2D}^{\mathrm{resp}}} \,.
\nonumber
\end{align}
To control the $-2\beta\E[c\Tr[AP^2]]$ term, note that
\[
    -2\beta\E[c\Tr[AP^2]] = -2\beta\E\left[c\sum_{i}(AP^2)_i\right] = -2\beta\E\left[c\sum_{i < j}A_{ij}\left((P^2)_{ij} + (P^2)_{ji}\right)\right] = -4\beta\sum_{i<j}\E\left[cA_{ij}(P^2)_{ij}\right]\,.
\]
Applying~\eqref{eq:GOE-ibp} gives
\[
    -2\beta\E[c\Tr[AP^2]] = -\frac{4\beta}{n}\sum_{i<j}\E\left[\partial_{ij}\left(c (P^2)_{ij}\right)\right]\,.
\]
At this point, the product rule followed by~\pref{lem:dP-pairing} along with some algebra and tracial identities yields that
\[
    -\frac{4\beta}{n}\E\left[c\,\partial_{ij} (P^2)_{ij}\right] = -4\E\left[c^2\Tr[P^2]\right] - \frac{4\beta^2}{n}\E\left[c\Tr[P^3]\right] + \frac{8\beta^2}{n}\E\left[c\sum_i P_{ii}(P^2)_{ii}\right] - R^{(3)}_1 - R^{(3)}_2\, ,
\]
where
\begin{align*}
    R^{(3)}_1 &= -\frac{4\beta^2}{n}\sum_{i<j}\E\left[c\left(\sum_k P_{kj}\left(m_iT_{ikj} + m_jT_{iki} + \Gamma_{ikij}\right)\right)\right] \\
    R^{(3)}_2 &= -\frac{4\beta^2}{n}\sum_{i<j}\E\left[c\left(\sum_k P_{ik}\left(m_iT_{kjj}+m_{j}T_{kji}+ \Gamma_{kjij}\right)\right)\right]\,,
\end{align*}
and
\[
    -\frac{4\beta}{n}\sum_{i<j}\E\left[(P^2)_{ij}\partial_{ij} c\right] = -\frac{8\beta^4}{n^2}\E\left[\sum_{i<j} (P^2)_{ij}(P^2)_{ji}\right] - R^{(3)}_3\, , 
\]
where
\begin{align*}
    R^{(3)}_3 = -\frac{4\beta^4}{n^2}\sum_{i<j}\E\left[(P^2)_{ij}\sum_k\left(m_iT_{kkj}  m_jT_{kki} + \Gamma_{kkij}\right)\right]\,.
\end{align*}
Putting these together yields the following final expression
\begin{align*}
    -2\beta\E\left[c\Tr[AP^2]\right] &= -4\E\left[c^2\Tr[P^2]\right] - \frac{4\beta^2}{n}\E\left[c\Tr[P^3]\right] + \frac{8\beta^2}{n}\E\left[c\left(\sum_i P_{ii}(P^2)_{ii}\right)\right]  - \frac{8\beta^4}{n^2}\E\left[\sum_{i<j}(P^2)_{ij}(P^2)_{ji}]\right] \\
    &- R^{(3)}_1 - R^{(3)}_2 - R^{(3)}_3\,. 
\end{align*}
Simple algebra along with the substitution that $c = \frac{\beta^2}{n}\Tr[P]$ yields the final expression as
\begin{align}\label{eq:AP2D-extracted}
&-2\beta\E\left(c\Tr[AP^2] + \Tr\left[AP^2D\right]\right)
=
-2\beta^2\,\E\Tr[P^2]\;-\;2\,\E[c\,\Tr(PDP)]- 4\E\left[c^2\Tr[P^2]\right] \nonumber \\
&\;-\;\underbrace{\frac{4\beta^2}{n}\E\Tr[DP^3]
\;+\;\frac{8\beta^2}{n}\E\Tr[P^2]
\;- \frac{4\beta^2}{n}\E\left[c\Tr[P^3]\right] + \frac{8\beta^2}{n}\E\left[c\left(\sum_i P_{ii}(P^2)_{ii}\right)\right] - \frac{8\beta^4}{n^2}\E\left[\sum_{i<j}(P^2)_{ij}(P^2)_{ji}\right]}_{:=\mathrm{R}_{t0}} \nonumber \\
&- R^{(1)} - R^{(2)} - R^{(3)} +\mathcal R_{AP^2D}^{\mathrm{resp}}\,.
\end{align}
All terms, excluding the first three, are shown to be $O(1)$. The last four terms are higher-order, and shown to be $O(1)$ in the proof of \pref{lem:remainder-bound-ea}, whereas the trace polynomial terms in $\mathrm{R}_{t0}$ are also $O(1)$ and this is proved at various places in the proof of \pref{lem:remainder-bound-ea} (see \pref{rem:trace-polys}). \\
At this point, we invoke~\pref{lem:cTrPDP-identity} to obtain that
\[
    -2\E\left[c\Tr[PDP]\right] = -2 \E\left[c\Tr[P]\right] - 2\E\left[c^2\Tr[P^2]\right] + O(1)\,.
\]
This yields
\[
    -2\beta\E\left(c\Tr[AP^2] + \Tr\left[AP^2D\right]\right) = -2\beta^2\,\E\Tr\left[P^2\right]-6\,\E\big[c^2\Tr\left[P^2\right]\big]-2\,\E\big[c\Tr\left[P\right]\big]+\text{ remainder terms}\,.
\]
\ppart{The $\beta^2\E[\Tr[PA^2P]]$ term} Using similar simplifications for tracial expansions,
\[
\beta^2 \E\left[\Tr[PA^2P]\right] = \beta^2\sum_{i,j}\E\left[A_{ij}(P^2A)_{ji}\right] = \beta^2 \E\left[\sum_{i < j}A_{ij}\left((P^2A)_{ij} + (P^2A)_{ji}\right) + \sum_{i=1}^n A_{ii}(P^2A)_{ii}\right]\,.
\]
The diagonal piece can be computed using the fact that $\Var(A_{ii})=2/n$ and every entry in $P$ is independent of the diagonal of $A$ via~\pref{lem:diag-does-not-enter-gibbs}. Specifically,
\begin{align*}
    \sum_i\E\left[A_{ii}(P^2A)_{ii}\right] &= \sum_i\sum_k\E\left[A_{ii}(P^2)_{ik}A_{ik}\right] = \sum_i \E\left[A^2_{ii}(P^2)_{ii}\right] \\
    &= \sum_i \E\left[A^2_{ii}\right]\E[(P^2)_{ii}] = \frac{2}{n}\E\left[\Tr[P^2]\right]\,.
\end{align*}
For the first term, we use~\eqref{eq:GOE-ibp} along with index algebra to simplify as
\begin{align*}
    \E\left[\sum_{i < j}A_{ij}\left((P^2A)_{ij} + (P^2A)_{ji}\right)\right] &= \frac1n\sum_{i<j}\E\Big[\partial_{ij}\big((P^2A)_{ji}+(P^2A)_{ij}\big)\Big] \\
    &= \frac1n\sum_{i<j}\E\left[(P^2)_{ii} + (P^2)_{jj} + \sum_k A_{kj}\partial_{ij}(P^2)_{ik} + \sum_k A_{ki}\partial_{ij}(P^2)_{jk} \right] \\
    &= \frac{1}{n}\sum_{i<j}\E\left[(P^2)_{ii} + (P^2)_{jj}\right] + \frac{1}{n}\sum_{i<j}\E\left[\sum_k A_{kj}\partial_{ij}(P^2)_{ik}+ \sum_k A_{ki}\partial_{ij}(P^2)_{jk} \right]\\
    &= \frac{n-1}{n}\E\left[\Tr[P^2]\right] + \frac{1}{n}\sum_{i<j}\E\left[\sum_k A_{kj}\partial_{ij}(P^2)_{ik}+ \sum_k A_{ki}\partial_{ij}(P^2)_{jk} \right]\,.
\end{align*}
Adding the previous simplification into the above equation yields
\[
    \beta^2\E\left[\Tr[PA^2P]\right] = \beta^2\E\left[\Tr[P^2]\right] + \frac{\beta^2}{n}\sum_{i<j}\E\left[\sum_k A_{kj}\partial_{ij}(P^2)_{ik} + \sum_k A_{ki}\partial_{ij}(P^2)_{jk}\right]+ \frac{\beta^2}{n}\E\left[\Tr[P^2]\right]
\]
We will now use~\pref{eq:dP-pairing} and~\eqref{eq:GOE-ibp-2} to simplify the remaining off-diagonal term further into a main term, and one that involves third centered moments and fourth cumulants (which will be part of a remainder term). First, using~\pref{eq:dP-pairing} gives
\[
    \partial_{ij}(P^2)_{ik} = \beta(P_{ii}(P^2)_{jk} + P_{ij}(P^2)_{ik} + (P^2)_{ii}P_{kj} + (P^2)_{ij}P_{ki}) + R^1_{i,j,ik}\, , 
\]
where $R^1_{i,j,ik}$ is a linear combination of terms involving the centered third-moments and fourth cumulants. Doing the same for the $\partial_{ij}(P^2)_{jk}$ term and summing over $k$ gives
\[
    \sum_k A_{kj}\partial_{ij}(P^2)_{ik} = \beta\sum_k A_{kj}\Big(
P_{ii}(P^2)_{jk}+P_{ij}(P^2)_{ik}+(P^2)_{ii}P_{kj}+(P^2)_{ij}P_{ki}\Big) + \sum_k R^{1}_{i,j,ik}\,.
\]
Collapsing the sums over $k$ for each summand gives
\[
     \sum_k A_{kj}\partial_{ij}(P^2)_{ik} = \beta\left( P_{ii}(AP^2)_{jj}+P_{ij}(AP^2)_{ji}+(P^2)_{ii}(AP)_{jj}+(P^2)_{ij}(AP)_{ji} \right) +  \sum_k R^{1}_{i,j,kj}\,.
\]
A similar calculation gives
\[
    \sum_k A_{ki}\partial_{ij}(P^2)_{jk} = \beta\left(P_{jj}(AP^2)_{ii}+P_{ij}(AP^2)_{ij}+(P^2)_{jj}(AP)_{ii}+(P^2)_{ij}(AP)_{ij}\right) + \sum_k R^{2}_{i,j,jk}\,.
\]
Note the following trace identities
\allowdisplaybreaks
\begin{align*}
\sum_{i<j}\Big(P_{ii}(AP^2)_{jj}+P_{jj}(AP^2)_{ii}\Big) &= \Tr[P]\Tr[AP^2]-\sum_i P_{ii}(AP^2)_{ii},\\
\sum_{i<j}P_{ij}\big((AP^2)_{ji}+(AP^2)_{ij}\big) &= \Tr[PAP^2]-\sum_i P_{ii}(AP^2)_{ii},\\
\sum_{i<j}\Big((P^2)_{ii}(AP)_{jj}+(P^2)_{jj}(AP)_{ii}\Big) &= \Tr[P^2]\Tr[AP]-\sum_i (P^2)_{ii}(AP)_{ii},\\
\sum_{i<j}(P^2)_{ij}\big((AP)_{ji}+(AP)_{ij}\big) &= \Tr[P^2AP]-\sum_i (P^2)_{ii}(AP)_{ii}\,,
\end{align*}
which we substitute into the in the above two terms and sum over $i < j$ to obtain
\allowdisplaybreaks
\begin{align*}
    &\frac{\beta^2}{n}\sum_{i<j}\E\left[\sum_k A_{kj}\partial_{ij}(P^2)_{ik} + \sum_k A_{ki}\partial_{ij}(P^2)_{jk}\right] = \\
    &\qquad \frac{\beta^3}{n}\E\left[\Tr[P]\Tr[AP^2] +\Tr[P^2]\Tr[AP]\right] + \underbrace{\left(\frac{\beta^3}{n}\E\left[2\Tr(PAP^2) -2\sum_i P_{ii}(AP^2)_{ii} -2\sum_i (P^2)_{ii}(AP)_{ii}\right]\right)}_{:= \calR^{lo}_{PA^2P}} \\
    &\qquad + \underbrace{\frac{\beta^2}{n}\E\left[\left(\sum_k R^{1}_{i,j,kj} + \sum_k R^{2}_{i,j,jk}\right)\right]}_{:= \calR^{ho}_{PA^2P}}\,. 
\end{align*}
We now apply~\pref{eq:GOE-ibp} to the first two terms in the expansion above (dropping diagonal terms, as $P$ is independent of the diagonal of $A)$ and collect the remaining three as remainder terms in $\calR^{AP}$. For the first term, using the product rule, \eqref{eq:GOE-ibp}, and~\pref{lem:dP-pairing} immediately gives
\allowdisplaybreaks
\begin{align*}
    &\E[\Tr[P]\Tr[AP^2]] = \frac{2}{n}\sum_{u<v}\E\big[\partial_{uv}\big(\Tr[P]\,(P^2)_{uv}\big)\big] = \frac{2}{n}\sum_{u<v}\E\big[(\partial_{uv}\Tr[P])(P^2)_{uv} + \Tr[P]\,\partial_{uv}(P^2)_{uv}\big] \\
    &=_{\text{\pref{lem:dP-pairing}}} \frac{2\beta}{n}\E\left[\Tr[P]^2\Tr[P^2]\right] + \underbrace{\frac{2\beta}{n}\E\left[Tr[P]\Tr[P^3] + \Tr[P^4] - 2\Tr[P]\sum_u P_{uu}(P^2)_{uu} - \sum_u (P^2_{uu})^2\right]}_{:=\mathrm{R}_{t1}} \\
    &\qquad + \underbrace{\frac{2\beta}{n}\sum_{u<v}\E\left[(P^2)_{uv}\left(m_u\sum_aT_{aav} + m_v\sum_a T_{aau}+\sum_a\Gamma_{aauv}\right)\right]}_{:= \mathrm{A}} \\
    &\qquad + \underbrace{\frac{2}{n}\sum_{u<v}\E\left[\Tr[P]\left(\sum_aP_{av}(m_uT_{uav} + m_vT_{uau}+\Gamma_{uauv}) + \sum_a P_{ua}(m_uT_{avv}+ + \Gamma_{avuv})\right)\right]}_{:= \mathrm{B}}\,.
 \end{align*}
 This immediately gives
 \[
    \frac{\beta^3}{n}\E[\Tr[P]\Tr[AP^2]] = 2\E[c^2\Tr[P^2]] + \frac{C(\beta)}{n}\mathrm{R}_{t1} + \frac{C(\beta)}{n}\mathrm{A} + \frac{C(\beta)}{n}\mathrm{B}\, ,
 \]
 where $\frac{C(\beta)}{n}\left(\mathrm{A} + \mathrm{B}\right)$ are of exactly the form of the terms $R^{(3)}_1, R^{(3)}_2, R^{(3)}_3$ defined in the proof of~\pref{lem:remainder-bound-ea}, where they are shown to be $O(1)$. For $\frac{C(\beta)}{n}\mathrm{R}_{t1}$, refer to \pref{rem:trace-polys}.
 
Similar calculations as the above also yield that
\allowdisplaybreaks
\begin{align*}
    &\E\left[\Tr[P^2]\Tr[AP]\right] = \frac{2}{n}\sum_{u<v}\E\left[\partial_{uv}(\Tr[P^2]P_{uv})\right] = \frac{2}{n}\sum_{u<v}\E\left[\Tr[P]\partial_{uv}P_{uv} + P_{uv}\partial_{uv}\Tr[P]\right] \\
    &=_{\text{\pref{lem:dP-pairing}}} \frac{\beta}{n}\E\left[\Tr[P^2]\Tr[P]^2\right] - \underbrace{\frac{\beta}{n}\E\left[\Tr[P^2]^2]\right] + \frac{4\beta}{n}\E\left[\Tr[P^4]\right] - \frac{4\beta}{n}\E\left[\sum_u P_{uu}(P^3)_{uu}\right]}_{:= \mathrm{R}_{t2}} \\
    &\qquad \underbrace{\frac{4\beta}{n}\sum_{u<v}\E\left[P_{uv}\sum_{a,b}P_{a,b}\left(m_uT_{abv} + m_vT_{abu}+\Gamma_{abuv}\right)\right]}_{:= \mathrm{C}}\,.
\end{align*}
This, again, immediately yields that
\[
    \frac{\beta^3}{n}\E\left[\Tr[P^2]\Tr[AP]\right] = \E\left[c^2\Tr[P^2]\right] + \frac{C(\beta}{n}\mathrm{R}_{t2} + \frac{C(\beta)}{n}\mathrm{C}\,, 
\]
where $\frac{C(\beta)}{n}\mathrm{C}$ has the same form as $\mathrm{HO}_9 + \mathrm{HO}_{10} + \mathrm{HO}_{11}$ (defined in \pref{lem:remainder-bound-ea}) and is $O(1)$ via the bounds proved on the latter in the proof of \pref{lem:remainder-bound-ea}, and $\frac{C(\beta)}{n}\mathrm{R}_{t2} = O(1)$ (see \pref{rem:trace-polys}). \\
Combining the above bounds gives the final expression and using the fact that $\frac{\beta^2}{n}\E[\Tr[P^2]] = O(1)$ as shown in \eqref{eq:c-p2-bound} gives
\[
    \beta^2\E\left[\Tr[PA^2P]\right] = 3\E\left[c^2\Tr[P^2]\right] + \beta^2E[\Tr[P^2]] + \calR^{lo}_{PA^2P} + \calR^{ho}_{PA^2P} + O(1)\,,
\]
and every type of remainder term in $ \calR^{lo}_{PA^2P}$ and $\calR^{ho}_{PA^2P}$ are captured in the set of remainders explicitly defined in the proof of \pref{lem:remainder-bound-ea}.

\ppart{Bounding the $2\beta\E[\Tr[AP^2V]]$ term} Showing that the rank-$1$ conjugated term is $O(1)$ is rather straightforward (though laborious) after applying \eqref{eq:GOE-ibp} followed by \pref{lem:dP-pairing} and doing some algebra after symmetrization of sums and the invocation of standard replica identities, as is done in the simplifications for $\E[\Tr[P^2]\Tr[AP]]$ and $\E[\Tr[P]\Tr[AP^2]]$. The details are delegated to \pref{lem:ea-rank-1-term}.

\ppart{Final bound for $\calE_A$} We now combine the bounds on all four terms, ignoring the $O(1)$ terms and higher-order remainder terms (that will be bounded in \pref{lem:remainder-bound-ea}) and obtain
\allowdisplaybreaks
\begin{align*}
    \calE_A &= \beta^2\,\E\Tr[PA^2 P]
    -2\beta\,\E\Tr[AP^2D]
    -2\beta \,\E c\Tr[AP^2]
    +2\beta\,\E\Tr[AP] + 2\beta\E\Tr[AP^2V] \\
    &= \left(3\E\left[c^2\Tr[P^2]\right] + \beta^2E[\Tr[P^2]] + \calR^{lo}_{PA^2P} + \calR^{ho}_{PA^2P} + O(1)\right) +\\
    & \left(-2\beta^2\,\E\Tr[P^2]-6\,\E\big[c^2\Tr[P^2]\big]-2\,\E\big[c\Tr(P)\big]+\text{ remainder terms}\right) + 2\E[c\Tr[P]] + O(1) \\
    &= -3\E\left[c^2\Tr[P^2]\right] -\beta^2\E\Tr[P^2] + \text{ all remainder terms } + O(1)\,.\qedhere
\end{align*}
\end{proof}

\begin{remark}\label{rem:trace-polys}
    Note that the trace polynomial terms in $\mathrm{R}_{t0}$ are explicitly bounded to be $O(1)$ in the proof of \pref{lem:remainder-bound-ea} as they appear in the explicit ``types'' of remainders associated to each term in $\calE_A$. $\mathrm{R}_{t1}$ and $\mathrm{R}_{t2}$ can be bounded to be $O(1)$ by using the bounds $\E[\Tr[P^2]] = O(n)$, $\E[\Tr[P^3]] = O(n)$ and $\E[\Tr[P^4]] = O(n^2)$ provided in \eqref{eq:c-p2-bound}, \eqref{eq:c-p3-bound} and \eqref{eq:c-p4-bound} in conjunction with the pointwise estimates that $0 \le P_{ii} \le 1$ and $|m_i| \le 1$. While this happens at many points during the (lengthy) proof of \pref{lem:remainder-bound-ea}, we collect some of these for clarity over here.
    \ppart{$\mathrm{R}_{t1}$ terms} Note the following bounds
    \begin{align*}
        \frac{C(\beta)}{n^2}\E\left[\Tr[P]\Tr[P^3]\right] &\le_{\Tr[P] \le n,\, P \succeq 0} \frac{C(\beta)}{n}\E\left[\Tr[P^3]\right] \le_{\text{\eqref{eq:c-p3-bound}}} O(1)\, ,\\
        \frac{C(\beta)}{n^2}\E\left[\Tr[P^4]\right] &=_{\text{\pref{eq:c-p4-bound}}} O(1)\,,\\
        \frac{C(\beta)}{n^2}\E\left[\Tr[P]\sum_u P_{uu}(P^2)_{uu}\right] &\le_{\Tr[P] \le n,\,|P_{uu}| \le 1} \frac{C(\beta)}{n}\E\left[\sum_u (P^2)_{uu}\right] =_{\text{\eqref{eq:c-p2-bound}}} O(1)\,,\\
        \frac{C(\beta)}{n^2}\E\left[\sum_u (P^2_{uu})^2\right] &\le_{(P^2)_{uu} = \sum_v (P_{uv})^2 \le Cn} \frac{C(\beta)}{n}\E\left[\Tr[P^2]\right] =_{\text{\eqref{eq:c-p2-bound}}} O(1)\,.
    \end{align*}

    \ppart{$\mathrm{R}_{t2}$ terms} Similarly, note the following set of complementary bounds
    \allowdisplaybreaks
    \begin{align*}
        \frac{C(\beta)}{n^2}\E\left[\Tr[P^2]^2\right] &= C(\beta)n^2 \nu(f^2_{1234}f^2_{5678}) =_{\text{CS\,+\,\pref{lem:replicon-moments-from-mgf}}} O(1)\, \\
        \frac{C(\beta)}{n^2}\E\left[\Tr[P^4]\right] &=_{\text{\eqref{eq:c-p4-bound}}} O(1)\, \\
        \frac{C(\beta)}{n^2}\E\left[\sum_u P_{uu}(P^3)_{uu}\right] &\le_{0 \le P_{uu} \le 1, P \succeq 0} \frac{C(\beta)}{n^2}\E\left[\Tr[P^3]\right] =_{\text{\eqref{eq:c-p3-bound}}} O(1)\,.
    \end{align*}
\end{remark}

\begin{lemma}[Bounding the remainder terms in $\calR^A$]\label{lem:remainder-bound-ea}
    The remainder terms in $\calE_A$ are $O(1)$.
\end{lemma}
\begin{proof}
    It is easily verified that the following are the types of remainders contributed by each of the 4 terms in $\calE_A$\footnote{~Note that, for instance, there are more remainder terms in $\beta^2\E[\Tr[PA^2P]]$ that come from centered third moments or fourth cumulants after applications of~\pref{lem:dP-pairing} in the proof of~\pref{lem:EA-term}, but they all have the same ``type'' with centered third moments with distinct indices or one repeated, and the same for the fourth cumulant terms. By symmmetry, therefore, it suffices to bound a single term of each type as this only loses $O(1)$ factors in the bound.}
    \allowdisplaybreaks
    \begin{align*}
        2\beta\E[\Tr[AP]] &\to \left\{-\frac{2\beta^2}{n}\E\Tr[P^2]\right\}\,, \\
        -2\beta\E[\Tr[AP^2D]] &\to \left\{\frac{8\beta^2}{n}\E\Tr[P^2], \frac{4\beta^2}{n}\E\Tr[DP^3], R_1,  R_2\right\}\,, \\
        -2\beta\E[c\Tr[AP^2]] &\to \left\{ \frac{4\beta^2}{n}\E\left[c\Tr[P^3]\right], \frac{8\beta^2}{n}\E\left[c\left(\sum_i P_{ii}(P^2)_{ii}\right)\right], \frac{8\beta^4}{n^2}\E\left[\sum_{i<j}(P^2)_{ij}(P^2)_{ji}\right], R^{(3)}_1, R^{(3)}_2, R^{(3)}_3\right\}\,, \\
        \beta^2\E[\Tr[PA^2P]] &\to \Bigg\{ \frac{\beta^4}{n}\sum_{i<j}\mathbb{E}\left[\sum_k A_{ki}m_i\sum_a P_{ja}T_{akj}\right], \frac{\beta^4}{n}\sum_{i < j}\mathbb{E}\left[\sum_k A_{ki}m_i\sum_b P_{bk}T_{jbj}\right], \frac{\beta^4}{n}\sum_{i<j}\mathbb{E}\left[\sum_k A_{ki}\sum_b P_{bk}\Gamma_{jbij}\right],\\
        &\frac{\beta^4}{n}\sum_{i<j}\mathbb{E}\left[\sum_k A_{ki}\sum_a P_{ja}\Gamma_{akij}\right], \frac{2\beta^3}{n}\E\left[\Tr[PAP^2]\right], \frac{2\beta^3}{n}\E\left[\sum_i P_{ii}(AP^2)_{ii}\right], \frac{2\beta^3}{n}\E\left[\sum_i (P^2)_{ii}(AP)_{ii}\right]\Bigg\}\,,
    \end{align*}
    where the the remainders containing centered third moments and fourth cumulants, and powers of diagonally weighted overlap-like terms, are specifically
    \allowdisplaybreaks
    \begin{align*}
        R_1 &= -\frac{2\beta}{n}\E\left[\sum_{i<j}(P^2)_{ij}\left(\partial_{ij}D_{ii} + \partial_{ij}D_{jj}\right)\right] \\
        R_2 &= -\frac{2\beta^2}{n}\E\left[\sum_{i<j}(D_{ii} + D_{jj})\left(\sum_kP_{kj}\left(m_iT_{ikj} + m_j T_{iki} + + \Gamma_{ikij}\right) + \sum_kP_{ik}(m_iT_{kjj} + m_jT_{kji} + \Gamma_{kjij} ) \right)\right]\,, \\
        R^{(3)}_1 &= -\frac{4\beta^2}{n}\sum_{i<j}\E\left[c\left(\sum_k P_{kj}\left(m_iT_{ikj} + m_jT_{iki} + \Gamma_{ikij}\right)\right)\right]\,, \\
        R^{(3)}_2 &= -\frac{4\beta^2}{n}\sum_{i<j}\E\left[c\left(\sum_k P_{ik}\left(m_iT_{kjj}+m_{j}T_{kji}+ \Gamma_{kjij}\right)\right)\right]\,, \\
        R^{(3)}_3 &= -\frac{4\beta^4}{n^2}\sum_{i<j}\E\left[(P^2)_{ij}\sum_k\left(m_iT_{kkj}  m_jT_{kki} + \Gamma_{kkij}\right)\right]\,.
    \end{align*}
    where $T_{abc}$ and $\Gamma_{abcd}$ are as defined in~\pref{lem:dP-pairing}. \\
    All of these remainder terms can be rewritten into linear combinations of cubic moments of rectangular sums $f$, weighted by quadratic terms in the overlap/magnetization, or can directly be rewritten in terms of third and fourth moments of rectangular sums that are diagonally weighted. Consequently, judicious applications of~\pref{prop:cubic-O4},~\pref{prop:D15} and~\pref{lem:D-oneleg-O(n)} suffice to show that the remainder terms are $O(1)$. 

    \ppart{The tracial terms} The easiest terms to bound are the tracial terms. Note by~\pref{eq:TrP2-f2}
    \[
        \frac{C(\beta)}{n}\E\Tr[P^2] = C(\beta)\, \nu(f^2)\,n \le_{\text{\pref{eq:replicon-2k-moment} with }k=1} C(\beta)\,. 
    \]
    Similarly, using the replica identities that $P_{ij} = \frac{1}{2}\an{(\sigma^1_i-\sigma^2_i)(\sigma^1_j-\sigma^2_j)}$ and $(P^2)_{ij} = \frac{n}{4}\an{(\sigma^1_i-\sigma^2_i)(\sigma^3_j-\sigma^4_j)f_{1234}}$, we have
    \allowdisplaybreaks
    \begin{align*}
        &\E\Tr[DP^3] = \E\Tr[PDP^2] = \E\left[\sum_{i,j}P_{ij}D_{jj}(P^2)_{ij}\right] \\
        &= \E\left[\sum_{i,j} \left(\frac{1}{2}\an{(\sigma^5_i-\sigma^6_i)(\sigma^5_j-\sigma^6_j)}\right)D_{jj}\left(\frac{n}{4}\an{(\sigma^1_i-\sigma^2_i)(\sigma^3_j-\sigma^4_j)f_{1234}}\right)\right] \\
        &= \frac{n}{8}\E\left[\sum_{i,j}D_{jj}\an{(\sigma^1_i-\sigma^2_i)(\sigma^3_j-\sigma^4_j)(\sigma^5_i-\sigma^6_i)(\sigma^5_j-\sigma^6_j)f_{1234}}\right] \\
        &= \frac{n^2}{8}\E\left[\sum_j D_{jj}\an{(\sigma^3_j-\sigma^4_j)(\sigma^5_j-\sigma^6_j)f_{1256}f_{1234}}\right] \\
        &= \frac{n^2}{8}\Bigg(\E\left[\sum_j D_{jj}\an{(\sigma^3_j-\sigma^4_j)(\sigma^5_j-\sigma^6_j)(R_{15}-q^*)f_{1234}}\right] + \E\left[\sum_j D_{jj}\an{(\sigma^3_j-\sigma^4_j)(\sigma^5_j-\sigma^6_j)(R_{26}-q^*)f_{1234}}\right] \\
        &\quad\quad\quad - \E\left[\sum_j D_{jj}\an{(\sigma^3_j-\sigma^4_j)(\sigma^5_j-\sigma^6_j)(R_{16}-q^*)f_{1234}}\right] - \E\left[\sum_j D_{jj}\an{(\sigma^3_j-\sigma^4_j)(\sigma^5_j-\sigma^6_j)(R_{25}-q^*)f_{1234}}\right]\Bigg)\,,
    \end{align*}
    where the last step uses the identity that $f_{1256} = \frac{1}{n}\sum_i a^{1256}_i = (R_{15}-q^*) + (R_{26}-q^*) - (R_{16}-q^*) - (R_{25}-q^*)$. At this point, observe that each of the four terms in a form compliant with~\pref{lem:D-oneleg-O(n)} with the addition of one more $(\sigma^3_j-\sigma^4_j)$ term, which is absolutely bounded by $2$. Invoking the lemma yields that
    \[
        \frac{C(\beta)}{n}\E\Tr[DP^3] = O(1)\,.
    \]
    Using the fact that $c = \frac{\beta^2}{n}\Tr[P] = \beta^2\an{1- R_{78}}$ and $\Tr[P^3]=\frac{n^3}{8}\an{f_{1234}f_{1256}f_{3456}}$ one obtains that
    \allowdisplaybreaks
    \begin{align}\label{eq:c-p3-bound}
        \frac{C(\beta)}{n}\E\left[c\Tr[P^3]\right] &= C(\beta)n^2\E\left[\an{f_{1234}f_{1256}f_{3456}}\right] - C(\beta) n^2 \E\left[\an{R_{78}f_{1234}f_{1256}f_{3456}}\right] \nonumber\\
        &\le_{(4,4,4,4)-\text{H\"older's,~\pref{lem:replicon-moments-from-mgf}}} C(\beta)n^2\,\mathsf{span}\left\{f_{1234}\,\Xi_{ab}\,\Xi_{cd}
        \right\} + O(1) \nonumber\\
        &=_{\text{\pref{prop:D15-full}}} O(1)\,,
    \end{align}
    where one uses the fact that $\nu(f\Xi_1\Xi_2) = \nu(f\Xi^-_1\Xi^-_2) + O(4)$. Similarly, a crude bound using the fact that $P_{ii} \le 1$ and $c = \beta^2(1-\an{R_{12}}) \le 2\beta^2$ implies that
    \allowdisplaybreaks
    \begin{align}\label{eq:c-p2-bound}
       \frac{C(\beta)}{n} \E\left[c\left(\sum_i P_{ii}(P^2)_{ii}\right)\right] &\le \frac{C(\beta)}{n}\E\left[c\sum_i (P^2)_{ii}\right] \le \frac{C'(\beta)}{n}\E\left[\Tr[P^2]\right] \nonumber \\
       &=_{\text{\pref{eq:TrP2-f2}}} C'(\beta)n\nu(f^2) =_{\text{\pref{lem:replicon-moments-from-mgf}}} O(1)\,. 
    \end{align}
    Using simple tracial algebra yields that $\sum_{i,j}(P^2)_{ij}(P^2)_{ji} = \Tr[P^4]$ and that $\sum_{i<j}(P^2)_{ij}P^2_{ji} \le \frac{1}{2}\Tr[P^4]$, which implies that
    \allowdisplaybreaks
    \begin{align}\label{eq:c-p4-bound}
        &\frac{C(\beta)}{n^2}\E\left[\sum_{i<j}(P^2)_{ij}(P^2)_{ji}\right] \le \frac{C'(\beta)}{n^2}\E\left[\Tr[P^4]\right] \nonumber \\
        &=C''(\beta)n^2\E\left[\an{f_{1234}f_{5678}f_{1278}f_{3456}}\right] \le_{(4,4,4,4)\text{-H\"older's}} C''(\beta)n^2 \nu(f^4_{1234})^{1/4}\nu(f^4_{5678})^{1/4}\nu(f^4_{1278})^{1/4}\nu(f^4_{3456})^{1/4} \nonumber \\
        &=_{\text{\pref{lem:replicon-moments-from-mgf}}} O(1)\,.
    \end{align}
    We now bound the tracial terms that involve a factor of $A$ and involve an application of~\pref{eq:GOE-ibp} before simplification with replica rewrites. Note using~\pref{lem:diag-does-not-enter-gibbs} and some standard algebra that
    \allowdisplaybreaks
    \begin{align*}
        \frac{2\beta^3}{n}\E\left[\Tr[PAP^2]\right] &=  \frac{2\beta^3}{n^2}\sum_{i<j}\E\left[\partial_{ij}P^3\right] \\
        &= \frac{2\beta^3}{n^2}\sum_{i<j}\E\left[\sum_{k,\ell}(\partial_{ij}P_{ik})P_{k\ell}P_{\ell j} + \sum_k P_{ik}(\partial_{ij}P_{k\ell})P_{\ell j} + \sum_kP_{ik}P_{k\ell}(\partial_{ij}P_{\ell j})\right] \\
        &=_{\text{\pref{lem:dP-pairing}}} \frac{2\beta^4}{n^2}\E\left[\sum_{i<j}(P_{ii}(P^3)_{jj} + P_{jj}(P^3)_{ii})\right] + \frac{2\beta^4}{n^2}\E\left[2\sum_{i<j}P_{ij}(P^3)_{ij}\right] + \mathrm{HO}^1 + \mathrm{HO}^2 + \mathrm{HO}^3\, ,
    \end{align*}
    where each term in $\mathrm{HO}_i$ is one of five types, namely
    \allowdisplaybreaks
    \begin{align*}
        \mathrm{HO}_1 &=\frac{2\beta^4}{n^2}\E\left[\sum_{i<j}m_i\sum_k(P^2)_{kj}T_{ikj}\right]\, , \\
        \mathrm{HO}_2 &=\frac{2\beta^4}{n^2}\E\left[\sum_{i<j}m_j\sum_k (P^2)_{kj}T_{iki}\right]\, , \\
        \mathrm{HO}_3 &=\frac{2\beta^4}{n^2}\E\left[\sum_{i<j}\sum_k (P^2)_{kj}\Gamma_{ikij}\right]\, , \\
        \mathrm{HO}_4 &=\frac{2\beta^4}{n^2}\E\left[\sum_{i<j}m_i\sum_{k,\ell} P_{ik}P_{\ell j}T_{k\ell j}\right]\, , \\
        \mathrm{HO}_5 &=\frac{2\beta^4}{n^2}\E\left[\sum_{i<j}\sum_{k,\ell} P_{ik}P_{\ell j}\Gamma_{k\ell ij}\right]\, .
    \end{align*}
    We will show, in the next part of the proof (where we bound the higher-order remainder terms involving centered third-moments and fourth cumulants) that each of these is $O(1)$. So, to obtain the final bound for $\frac{2\beta^3}{n}\E\left[\Tr[PAP^2]\right]$, we now bound the second-order parts. Simplifying further into trace terms and using the fact that $c = \frac{\beta^2}{n}\Tr[P]$
    \allowdisplaybreaks
    \begin{align*}
         \frac{2\beta^4}{n^2}\E\left[\sum_{i<j}(P_{ii}(P^3)_{jj} + P_{jj}(P^3)_{ii})\right] + \frac{2\beta^4}{n^2}\E\left[2\sum_{i<j}P_{ij}(P^3)_{ij}\right] &= \frac{2\beta^2}{n}\E\left[c\Tr[P^3]\right] + \frac{2\beta^4}{n^2}\E\left[\Tr[P^4]\right] - \frac{4\beta^4}{n^2}\E\left[\sum_i P_{ii}(P^3)_{ii}\right]\,.
    \end{align*}
    Note the first two terms are easily shown to be $O(1)$ by \eqref{eq:c-p3-bound} and \eqref{eq:c-p4-bound}. The final term can be handled by the observation that $\Tr[P^3] = \frac{n^3}{8}\nu(f_{1234}f_{3456}f_{1256})$ and $0 \preceq E_{\calD_n}[P] \le \Id_n$ and $0 \preceq P$ to obtain
    \allowdisplaybreaks
    \begin{align*}
        \frac{4\beta^4}{n^2}\E\left[\sum_i P_{ii}(P^3)_{ii}\right] \le \frac{4\beta^3}{n^2}\E[\Tr[P^3]] = C(\beta)n\nu(f_{1234}f_{3456}f_{1256}) \le O(n^{-1/2})\, ,
    \end{align*}
    by an application of $(2,4,4)$-H\"older's followed by an invocation of~\pref{lem:replicon-moments-from-mgf}. This implies that $\frac{2\beta^3}{n}\E\left[\Tr[PAP^2]\right] = O(1)$. \\
    We turn to the remaining two trace terms from $\beta^2\E[\Tr[PA^2P]]$. Note that using the fact that $|\partial_{ij}P_{ii}| \le 2\beta|(m^2_iP_{ij}-m_im_jP_{ii})| \le 8\beta$, we have
    \allowdisplaybreaks
    \begin{align*}
        &\frac{2\beta^3}{n}\E\left[\sum_i P_{ii}(AP^2)_{ii}\right] = \frac{2\beta^3}{n^2}\sum_{i<j}\E\left[\partial_{ij}\left((P^2)_{ij}\left(P_{ii}+P_{jj}\right)\right)\right] \\
        &= \frac{2\beta^3}{n^2}\sum_{i<j}\E\left[\left((P^2)_{ij}\left(\partial_{ij}P_{ii}+\partial_{ij}P_{jj}\right) + \left(P_{ii}+P_{jj}\right)\partial_{ij}(P^2)_{ij}\right)\right] \\
        &=_{\text{\pref{lem:dP-pairing}\,+\,symmetry}}\frac{32\beta^4}{n^2}\E\left[\sum_{i<j}(P^2)_{ij}\right] + \frac{C(\beta)}{n^2}\E\left[\sum_{i<j}(P_{ii}+P_{jj})\left(P_{ii}(P^2)_{jj} + P_{jj}(P^2)_{ii}\right) + 2P_{ij}(P^2)_{ij}\right] \\
        &\qquad\qquad\qquad\qquad\quad + \mathrm{HO}_6 + \mathrm{HO} _7 + \mathrm{HO}_8\,.
    \end{align*}
    At this point observe that because $P \succeq 0$, a Cauchy--Schwarz impliciation means that
    \[
        \frac{32\beta^4}{n^2}\E\left[\sum_{i<j}(P^2)_{ij}\right] \le \frac{C(\beta)}{n}\sqrt{\Tr[P^4]} = O(1)\, ,
    \]
    using \eqref{eq:c-p4-bound}. The other pairing term is handled similarly using the facts that $(P_{ii} + P_{jj})(P_{ii}(P^2)_{jj}) \le 2 (P^2)_{jj}$ and analogously on swapping the indices $i$ and $j$ for the pairing term. This gives
    \allowdisplaybreaks
    \begin{align*}
        &\frac{C(\beta)}{n^2}\E\left[\sum_{i<j}(P_{ii}+P_{jj})\left(P_{ii}(P^2)_{jj} + P_{jj}(P^2)_{ii}\right) + 2P_{ij}(P^2)_{ij}\right] \le \frac{C'(\beta)}{n}\E\left[\Tr[P^2]\right]  + \frac{C(\beta)}{n^2}\sum_{i<j}\E\left[2P_{ij}(P^2)_{ij}\right] \\
        &\le_{\text{\eqref{eq:c-p2-bound},\,CS}} O(1) + \frac{C(\beta)}{n^2}\sqrt{\E[\Tr[P^2]]}\sqrt{\E[\Tr[P^4]]} \le_{\text{\eqref{eq:c-p2-bound},\,\eqref{eq:c-p4-bound}}} O(1) + O(n^{-1/2}) = O(1)\,.   
    \end{align*}
    The higher-order terms have the following types,
    \begin{align*}
        \mathrm{HO}_6 &= \frac{2\beta^4}{n^2}\sum_{i<j}\E\left[(P_{ii} + P_{jj})m_i\sum_k P_{kj}T_{ikj}\right]\,, \\
        \mathrm{HO}_7 &= \frac{2\beta^4}{n^2}\sum_{i<j}\E\left[(P_{ii} + P_{jj})m_j\sum_k P_{kj}T_{iki}\right]\,, \\
        \mathrm{HO}_8 &= \frac{2\beta^4}{n^2}\sum_{i<j}\E\left[(P_{ii} + P_{jj})\sum_k P_{kj}\Gamma_{ikij}\right]\,.
    \end{align*}
    These higher-order terms are bounded in the next part of the proof and shown to be $O(1)$. This implies that $\frac{2\beta^3}{n}\E\left[\sum_i P_{ii}(AP^2)_{ii}\right] = O(1)$.
    For the second tracial term, observe that 
    \allowdisplaybreaks
    \begin{align*}
        \frac{2\beta^3}{n}\E\left[\sum_i (P^2)_{ii}(AP)_{ii}\right] &= \frac{2\beta^3}{n}\sum_{i<j}\E\left[A_{ij} P_{ij}\left((P^2)_{ii} + (P^2)_{jj}\right)\right] \\
        &=_{\eqref{eq:GOE-ibp}} \frac{2\beta^3}{n^2}\sum_{i<j}\E\left[(\partial_{ij}P_{ij})\left((P^2)_{ii} + (P^2)_{jj}\right)\right] + \frac{2\beta^3}{n^2}\sum_{i<j}\E\left[P_{ij}\left(\partial_{ij}(P^2)_{ii} + \partial_{ij}(P^2)_{jj}\right)\right]\,.
    \end{align*}
    Using the fact that $\partial_{ij}P_{ij} = \beta(P_{ii}P_{jj}-P^2_{ij})$ and that $|(P_{ii}P_{jj} - P^2_{ij})| \le 4$, the first term can be simplified as
    \allowdisplaybreaks
    \begin{align*}
        \frac{C(\beta)}{n^2}\left|\E\left[\sum_{i<j}(\partial_{ij}P_{ij})\left((P^2)_{ii} + (P^2)_{jj}\right)\right]\right| &\le \frac{C'(\beta)}{n^2}\left|\E\left[\sum_{i<j}\left((P^2)_{ii} + (P^2)_{jj}\right)\right]\right| \\
        &\le \frac{C'(\beta)(n-1)}{n^2}\E\left[\Tr[P^2]\right] = C'(\beta)(n-1)\nu(f^2) =_{\text{\pref{lem:replicon-moments-from-mgf}}} O(1)\,.
    \end{align*}
    To bound the second term, we use similar simplifications to obtain
    \[
        \frac{2\beta^3}{n^2}\sum_{i<j}\E\left[P_{ij}\left(\partial_{ij}(P^2)_{ii} + \partial_{ij}(P^2)_{jj}\right)\right] = \frac{4\beta^3}{n^2}\sum_{i<j}\E\left[P_{ij}\sum_k P_{ik}\partial_{ij}P_{ik}\right] + \frac{4\beta^3}{n^2}\sum_{i<j}\E\left[P_{ij}\sum_k P_{jk}\partial_{ij}P_{jk}\right]\,,
    \]
    and so, by symmetry in the indices under a $i \leftrightarrow j$ swap, it suffices to bound one term. So, note that using~\pref{lem:dP-pairing} we get that
    \allowdisplaybreaks
    \begin{align*}
        \frac{4\beta^3}{n^2}\sum_{i<j}\E\left[P_{ij}\sum_k P_{ik}\partial_{ij}P_{ik}\right] &= \frac{4\beta^4}{n^2}\sum_{i<j}\E\left[P_{ii}P_{ij}(P^2)_{ij}\right] + \frac{4\beta^4}{n^2}\sum_{i<j}\E\left[(P_{ij})^2(P^2)_{ii}\right] + \mathrm{HO}_6 + \mathrm{HO}_7 + \mathrm{HO}_8\, ,
    \end{align*}
    where
    \allowdisplaybreaks
    \begin{align*}
        \mathrm{HO}_9 &= \frac{4\beta^4}{n^2}\sum_{i<j}\E\left[m_iP_{ij}\sum_k P_{ik}T_{ikj}\right]\,,\\
        \mathrm{HO}_{10} &= \frac{4\beta^4}{n^2}\sum_{i<j}\E\left[m_jP_{ij}\sum_k P_{ik}T_{iki}\right]\,,\\
        \mathrm{HO}_{11} &= \frac{4\beta^4}{n^2}\sum_{i<j}\E\left[P_{ij}\sum_k P_{ik}\Gamma_{ikij}\right]\,.
    \end{align*}
    For the lower-order terms, observe that
    \begin{align*}
        \frac{4\beta^4}{n^2}\E\left[\sum_{i<j} (P^2)_{ii}(P_{ij})^2\right] &\le \frac{4\beta^4}{n^2}\E\left[\sum_i (P^2_{ii})\sum_{j}P_{ij}P_{ji}\right] = \frac{4\beta^4}{n^2}\E\left[\sum_i (P^2)^2_{ii}\right] \le \frac{C(\beta)}{n}\E\left[c\Tr[P^2]\right] \\
        &\le_{c \le \beta^2} \frac{C'(\beta)}{n}\E[\Tr[P^2]] \le C'(\beta) n \nu(f^2) \le_{\text{\pref{lem:replicon-moments-from-mgf}}} O(1)\,, 
    \end{align*}
    and
    \allowdisplaybreaks
    \begin{align*}
        \frac{4\beta^4}{n^2}\E\left[\sum_{i<j} (P_{ij})^2 (P^2)_{ii}\right] \le_{P = P^\sT} \frac{4\beta^4}{n^2}\E\left[\Tr[P^4]\right] \le C(\beta)n^2\nu(f_{1234}f_{5678}f_{1278}f_{3456}) = O(1)\, , 
    \end{align*}
    by \eqref{eq:c-p4-bound}. The higher-order terms are, again, shown to be $O(1)$ in the next part of the proof. Those bounds complete the proof and show that $\frac{2\beta^3}{n}\E\left[\sum_i (P^2)_{ii}(AP)_{ii}\right] = O(1)$.
    
    \ppart{The derivative terms} We now use the facts that $\partial_{ij}D_{ij}=2\beta D^2_{jj}m_j(m_i - m_j\an{\sigma_i\sigma_j})$ and $m_i-m_j\an{\sigma_i\sigma_j}=P_{jj}m_i - m_jP_{ij}$ in conjunction with replica identities and symmetrization arguments to bound terms such as $R_1$. Some elementary algebra immediately gives
    \allowdisplaybreaks
    \begin{align*}
    R_1 &= -\frac{4\beta^2}{n}\E\left[\sum_{j=1}^n m_j D^2_{jj}\sum_{i \ne j}(P^2)_{ij}\left(P_{jj}m_i - m_jP_{ij}\right)\right] \\
    &= -\frac{4\beta^2}{n}\,\E\Bigg[\sum_{j=1}^n m_jP_{jj}D_{jj}^2\sum_{i\neq j}m_i(P^2)_{ij}\Bigg] +\frac{4\beta^2}{n}\,\E\Bigg[\sum_{j=1}^n m_j^2D_{jj}^2\sum_{i\neq j}P_{ij}(P^2)_{ij}\Bigg].
    \end{align*}
    We collapse the inner sums with replica identities for $P_{ij}$ and $(P^2)_{ij}$ using the facts that 
    \[
        \sum_{i\ne j}\Delta\sigma^{k_1k_2}_i\Delta\sigma^{k_3k_4}_i = nf_{k_1k_2k_3k_4}\,,
    \]
    and
    \[
        \frac{1}{n}\sum_{i=1}^n \Delta\sigma^{k_1k_2}_i\sigma^{k_3}_i = R_{k_1k_3}-R_{k_2k_3}\, ,
    \]
    to simplify the term further as
    \allowdisplaybreaks
    \begin{align*}
        R_1 &= -\frac{4\beta^2}{n}\,\mathbb E\Bigg[\sum_{j=1}^n m_jP_{jj}D_{jj}^2\frac{n^2}{4}\Big\langle (R_{15}-R_{25})\,\Delta\sigma_j^{34}\,f_{1234}\Big\rangle\Bigg]\\
            &\quad +\frac{4\beta^2}{n}\,\mathbb E\Bigg[\sum_{j=1}^n m_j^2D_{jj}^2\cdot\frac{n^2}{8}\Big\langle\Delta\sigma_j^{56}\Delta\sigma_j^{34}\,f_{1234}\,f_{1256}\Big\rangle\Bigg] \\
            &= -\beta^2 n\E\left[\sum_j D_{jj}\left(m_j(\sigma^3_j-\sigma^4_j)\right)((R_{15}-q^*)-(R_{25}-q^*))f_{1234}\right] \\
            &\quad + \frac{\beta^2 n}{2}\E\left[\sum_j D^2_{jj} \left(m^2_j(\sigma^3_j-\sigma^4_j)(\sigma^5_j-\sigma^6_j)\left((R_{15}-q^*) + (R_{26}-q^*) - (R_{16}-q^*)\right)f_{1234}\right)\right]
    \end{align*}
    The first term is a finite linear combination of terms of the form $\E[\an{\sum_j D_{jj}a_nB_j\Xi f_{1234}}]$ where $B_j$ is an absolutely bounded function of the spin $j$ -- this term is bounded, therefore, by a direction application of~\pref{lem:D-oneleg-O(n)}. Note that the final term is a linear combination of terms of the form $\E[\an{\sum_j D^2_{jj} a_n B_j \Xi f_{1234}}]$ where $B_j$ is an absolutely bounded term in the $j$-th spin, and $\Xi$ is a centered bulk deviation. Similar to the bound for $\E[\Tr[DP^3]]$, one can invoke~\pref{lem:D-oneleg-O(n)} to bound this\footnote{~By~\pref{rem:boundedness-handling}, note that the extra power of $D_{jj}$ does not change the scaling of the bound, since $\E[\an{D^k_{jj}}] < \infty$ for every $k \in \N$ by~\pref{lem:D-moments}, and nor does any extra factor $B$ provided the $a_n$ is present in the term.}. 
    
    \ppart{The third-centered moment and fourth-cumulant terms} To control these terms we can use lemmata that bound quartic moments of rectangular sums that are weighted by diagonals and/or bulk deviation terms. To do this, we will use the fact that $(\sigma_j-m_j)^2 = P_{jj} - 2m_j(\sigma_j-m_j)$ which can be used to show the following simplifications for centered-third moment terms and fourth-cumulant terms when they have repeated indices
    \[
        T_{kjj} = -2m_jP_{kj},\qquad T_{iki} = -2m_iP_{ik},\qquad T_{kkj} = -2m_kP_{kj},\qquad T_{kki} = -2m_kP_{ki}\,,
    \]
    and
    \[
        \Gamma_{ikij} = -2P_{ik}P_{ij} -2m_iT_{ikj},\qquad \Gamma_{kjij} = -2P_{kj}P_{ij} -2m_jT_{kij},\qquad \Gamma_{kkij} =-2P_{ki}P_{kj} -2m_kT_{kij}\,.
    \]
    Using these we bound each term chronologically.
    \allowdisplaybreaks
    \begin{align*}
        &R_2 =  -\frac{2\beta^2}{n}\E\left[\sum_{i<j}(D_{ii} + D_{jj})\left(\sum_kP_{kj}\left(m_iT_{ikj} + m_j T_{iki} + + \Gamma_{ikij}\right) + \sum_kP_{ik}(m_iT_{kjj} + m_jT_{kji} + \Gamma_{kjij} ) \right)\right] \\
        &= -\frac{2\beta^2}{n}\E\left[\sum_{i<j}(D_{ii} + D_{jj})\left(\sum_k P_{kj}(-m_iT_{ikj} - 2m_im_jP_{ik} - 2P_{ik}P_{ij})  + P_{ik}(-m_jT_{kji} - 2m_im_jP_{kj} - 2P_{kj}P_{ij})\right)\right] \\
        &= \frac{4\beta^2}{n}\E\left[\sum_{i<j}(D_{ii}+D_{jj})\left(m_im_j(P^2)_{ij} + P_{ij}(P^2)_{ij}\right)\right] +\frac{2\beta^2}{n}\E\left[\sum_{i<j}(D_{ii}+D_{jj})\left(\sum_k (T_{ikj}m_iP_{kj} + T_{jki}m_jP_{ik})\right)\right] \\
        &= \frac{4\beta^2}{n}\E\left[\sum_{i\ne j}D_{ii}P_{ij}(P^2)_{ij}\right] + \frac{4\beta^2}{n}\E\left[\sum_{i \ne j}D_{ii}m_im_j(P^2)_{ij}\right] + \frac{2\beta^2}{n}\E\left[\sum_{i<j}(D_{ii}+D_{jj})\left(\sum_k (T_{ikj}m_iP_{kj} + T_{jki}m_jP_{ik})\right)\right] \\
        &= \frac{4\beta^2}{n}\E\left[\Tr[DP^3]- \Tr[P^2]\right] + \frac{4\beta^2}{n}\E\left[\langle m , P^2 (Dm) \rangle + O_\beta\left(\Tr[DP^2]\right)\right] \\
        &\qquad\qquad\qquad + \frac{2\beta^2}{n}\E\left[\sum_{i<j}(D_{ii}+D_{jj})\left(\sum_k (T_{ikj}m_iP_{kj} + T_{jki}m_jP_{ik})\right)\right]\,.
    \end{align*}
    We have already established earlier in the proof that the first two terms in the equality are $O(1)$. For the third term, note that $\frac{C(\beta)}{n}\E\left[\langle m, P^2(Dm)\rangle\right] = O(1)$ by the bound for the $\E[\Tr[PDVP]]$ term in the proof of~\pref{lem:bound-d-remainder}. To bound the fourth term, note that
    \[
        \frac{C(\beta)}{n}\E\left[\Tr[DP^2]\right] = \frac{C(\beta)}{n}\E\left[\sum_{i}D_{ii}(P^2)_{ii}\right] = C(\beta)\E\left[\sum_i \an{D_{ii}(\sigma^1_i-\sigma^2_i)(\sigma^3_i-\sigma^4_i)f}\right] = C(\beta) n\nu(f\hat{f}) \le_{\text{CS}} O(1)\,, 
    \]
    where we use~\pref{lem:replicon-moments-from-mgf} and~\pref{lem:hatf-moments}. To bound the final term, note the following simplifications based on replica rewrites,
    \allowdisplaybreaks
    \begin{align*}
        &\frac{2\beta^2}{n}\E\left[\sum_{i<j}(D_{ii}+D_{jj})\left(\sum_k (T_{ikj}m_iP_{kj} + T_{jki}m_jP_{ik})\right)\right] =_{T_{ikj}=\an{(\sigma^1_i-\sigma^2_i)(\sigma^1_k-\sigma^3_k)(\sigma^1_j-\sigma^4_j)}} \\
        &\qquad\qquad\qquad\frac{\beta^2}{n}\E\left[\sum_{i<j}(D_{ii}+D_{jj})\sum_k \an{\sigma^7_i(\sigma^1_i-\sigma^2_i)(\sigma^1_k-\sigma^3_k)(\sigma^1_j-\sigma^4_j)(\sigma^5_k-\sigma^6_k)(\sigma^5_j-\sigma^6_j)}\right] \\
        &\qquad\qquad\qquad + \frac{\beta^2}{n}\E\left[\sum_{i<j}(D_{ii} + D_{jj})\sum_k\an{\sigma^7_j(\sigma^1_k-\sigma^2_k)(\sigma^1_j-\sigma^3_j)(\sigma^1_i-\sigma^4_i)(\sigma^5_i-\sigma^6_i)(\sigma^5_k-\sigma^6_k)}\right] \\
        &\qquad\qquad=\beta^2\E\left[\sum_{i<j}(D_{ii}+D_{jj})\left(\an{\sigma^7_i(\sigma^1_i-\sigma^2_i)(\sigma^1_j-\sigma^4_j)(\sigma^5_j-\sigma^6_j)f_{1356}}\right)\right] \\
        &\qquad\qquad + \beta^2\E\left[\sum_{i<j}(D_{ii}+D_{jj})\left(\an{\sigma^7_j(\sigma^1_i-\sigma^4_i)(\sigma^5_i-\sigma^6_i)(\sigma^1_j-\sigma^3_j)f_{1256}}\right)\right] \\
        &= \frac{\beta^2}{2}\E\left[\sum_j D_{jj}\an{(\sigma^1_j-\sigma^4_j)(\sigma^5_j-\sigma^6_j)f_{1356}\sum_{i \ne j} \sigma^7_i(\sigma^1_i-\sigma^2_i)}\right] + \frac{\beta^2}{2}\E\left[\sum_i D_{ii}\an{(\sigma^7_if_{1356} (\sigma^1_i-\sigma^2_i)}\sum_{j \ne i}(\sigma^1_j-\sigma^4_j)(\sigma^5_j-\sigma^6_j)\right] \\
        &\qquad\qquad\qquad = \frac{\beta^2 n}{2}\E\left[\sum_j D_{jj}\an{(\sigma^1_j-\sigma^4_j)(\sigma^5_j-\sigma^6_j)f_{1256}\left((R^{-j}_{17}-q^*) - (R^{-j}_{27}-q^*)\right)}\right] \\
        &\qquad \qquad \qquad + \frac{\beta^2 n}{2}\E\left[\sum_j D_{jj}\an{\sigma^7_i(\sigma^1_i-\sigma^2_i)f_{1256}\left((R^{-i}_{15}-q^*) + (R^{-i}_{46}-q^*) - (R^{-i}_{16}-q^*) - (R^{-i}_{45}-q^*)\right)}\right]\,.
    \end{align*}
    The two terms above are both finite linear combinations of terms of the form $\E\left[\sum_i D_{ii}\an{(\sigma^a_i-\sigma^b_i)B_sf_{abcd}\Xi^-}\right]$ which are easily shown to be $O(1/n)$ by~\pref{lem:D-oneleg-O(n)} and~\pref{rem:boundedness-handling}. This immediately implies, in conjunction with the prior bounds, that
    \[
        R_2 = O(1)\,.
    \]
    Using the same simplifications for the repeated-indices centered third-moments and fourth cumulants mentioned above, it is not hard to show that
    \allowdisplaybreaks
    \begin{align*}
        R^{(3)}_1 &= \frac{8\beta^2}{n}\E\left[c\sum_{i<j}m_im_j(P^2)_{ij}\right] + \frac{8\beta^2}{n}\E\left[c\sum_{i<j}P_{ij}(P^2)_{ij}\right] + \frac{4\beta^2}{n}\E\left[c\sum_{i<j} m_i\sum_k P_{kj}T_{ikj}\right] \\
        &\le_{|c| \le \beta^2} \frac{C(\beta)}{n}\E\left[c\,\langle m, P^2 m\rangle\right] + \frac{8\beta^2}{n}\E\left[c\Tr[P^3] - c\sum_i P_{ii}(P^2)_{ii}\right] + \frac{4\beta^2}{n}\E\left[c\sum_{i<j} m_i\sum_k P_{kj}T_{ikj}\right]\,.
    \end{align*}
    For the first term, note that the fact that $P \succeq 0$ and an application of $(1,\infty)$-H\"older's in conjunction with the fact that $|c| \le \beta^2$ immediately implies that
    \[
        \frac{C(\beta)}{n}\E\left[c\,\langle m, P^2 m\rangle\right] \le \frac{C'(\beta)}{n}\E\left[\langle m, P^2 m\right] \le O(1)\, ,
    \]
    where the final inequality uses the bound for $\E[\Tr[PVP]]$ in the proof of~\pref{lem:bound-d-remainder}. The second term is already bounded by $O(1)$ by the bounds developed for the tracial terms earlier in the proof. For the final term, note using replica identities that
    \allowdisplaybreaks
    \begin{align*}
        &\frac{4\beta^2}{n}\E\left[c\sum_{i<j} m_i\sum_k P_{kj}T_{ikj}\right] = \frac{2\beta^2}{n}\E\left[c \sum_{i < j}m_i \an{\sum_k (\sigma^1_i-\sigma^2_i)(\sigma^1_k-\sigma^3_k)(\sigma^1_j-\sigma^4_j)(\sigma^5_k-\sigma^6_k)(\sigma^5_j-\sigma^6_j)}\right] \\
        &= 2\beta^2\E\left[c\sum_{i<j}\an{\sigma^7_i(\sigma^1_i-\sigma^2_i)(\sigma^1_j-\sigma^4_j)(\sigma^5_j-\sigma^6_j)f_{1356}}\right] \\
        &= \beta^2\E\left[\an{((1-q^*)-(R_{8,9}-q^*))\left(\sum_i(\sigma^7_i\sigma^1_i - \sigma^7_i\sigma^2_i)\right)\left(\sum_j (\sigma^1_j-\sigma^4_j)(\sigma^5_j-\sigma^6_j)\right)f_{1356}}\right] \\
        &= C(\beta,q^*)n^2\mathsf{span}\left\{\E\left[\an{(R_{8,9}-q^*)\Xi_{ab}\Xi_{cd}f_{1356}}\right],\E\left[\an{\Xi_{ab}\Xi_{cd}f_{1356}}\right]\right\}\, ,
    \end{align*}
    where the first set of terms are directly controlled to be $O(n^{-2})$ by a $(4,4,4,4)$-H\"older's inequality, followed by an invocation of~\pref{lem:replicon-moments-from-mgf} and~\pref{thm:overlap-moment-concentration}. The second set of terms are $O(n^{-2})$ by a direct invocation of~\pref{prop:D15}. \\
    Due to the symmetry and exact same structure of the terms, an exactly analogous argument shows that
    \[
        R^{(3)}_2 = O(1)\,.
    \]
    We now bound the $R^{(3)}_3$ term by first simplifying using the same replica identities and repeated-index centered third moment and fourth cumulant simplifications to obtain that
    \[
        m_iT_{kkj}m_jT_{kki} = m_i(-2m_kP_{kj})m_j(-2m_kP_{ki}) = 4m_im_jm^2_kP_{jk}P_{ki}\,,
    \]
    and
    \[
        \Gamma_{kkij} = -2m_kT_{kij} - 2P_{ik}P_{kj}\,.
    \]
    This gives,
    \allowdisplaybreaks
    \begin{align*}
        &R^{(3)}_3 = -\frac{4\beta^4}{n^2}\sum_{i<j}\E\left[(P^2)_{ij}\sum_k\left(m_iT_{kkj}  m_jT_{kki} + \Gamma_{kkij}\right)\right] \\
        &= \frac{16\beta^4}{n^2}\sum_{i<j}\E\left[(P^2)_{ij}m_im_j(P\mathsf{diag}(m)P)_{ij}\right] - \frac{8\beta^4}{n^2}\sum_{i<j}\E\left[(P^2)_{ij}\sum_k m_kT_{kij}\right] - \frac{8\beta^4}{n^2}\sum_{i<j}\E\left[(P^2)_{ij}(P^2)_{ji}\right] \\
        &\le_{\text{CS}} \frac{16\beta^4}{n^2}\E\left[\sqrt{\sum_{i<j}(P^2)^2_{ij}}\sqrt{\sum_{i<j}(P^2)^2_{ij}}\opnorm{\diag(m^2)}\right] - \frac{4\beta^4}{n^2}\E\left[\Tr[P^4]\right] + \frac{4\beta^4}{n^2}\E\left[\sum_i (P^2)^2_{ii}\right] - \frac{8\beta^4}{n^2}\sum_{i<j}\E\left[(P^2)_{ij}\sum_k m_kT_{kij}\right]\,. 
    \end{align*}
    For the first three terms, note that $\opnorm{\diag(m^2)} \le 1$ and that $\sum_{i<j}(P^2)^2_{ij} = \frac{1}{2}\Tr[P^4] - \sum_i (P^2)^2_{ii}$. Since $\sum_i (P^2)^2_{ii} \le \norm{P^2}^2_F = \Tr[P^4]$, this implies that 
    \[
        \frac{4\beta^4}{n^2}\E\left[\sum_i (P^2)^2_{ii}\right] \le \frac{C(\beta)}{n^2}\E\left[\Tr[P^4]\right]] = C(\beta)n^2 \nu(f_{1234}f_{34556}f_{5678}f_{7812}) = O(1)\,.
    \]
    where we use~\pref{lem:replicon-moments-from-mgf} along with $(4,4,4,4)$-H\"older's. This immediately bounds the first three terms by $O(1)$. For the final term, observe that
    \allowdisplaybreaks
    \begin{align*}
        \frac{8\beta^4}{n^2}\sum_{i<j}\E\left[(P^2)_{ij}\sum_k m_kT_{kij}\right] &=  \frac{C(\beta)}{n}\sum_{i<j}\E\left[\an{(\sigma^1_i-\sigma^2_i)(\sigma^1_j-\sigma^2_j)f_{1234}}\sum_k \an{\sigma^9_k(\sigma^5_k-\sigma^6_k)(\sigma^5_i-\sigma^7_i)(\sigma^5_j-\sigma^8_j)}\right] \\
        &= \mathsf{span}\left\{C(\beta)n^2\nu\left(\Xi_{ab}\,\Xi_{cd}\,\Xi_{ef}f_{1234}\right),C'(\beta)n\nu\left(f_{1234}\Xi_{ab}F\right)\right\}\,,
    \end{align*}
    where the last equality follows after symmetrizing the summations over $i$ and $j$ and then summing then $a^{abcd}$ factors separately with their index, and $F$ is some absolutely bounded function of all replicas. 
    At this point, note that $(4,4,4,4)$-H\"older's along with~\pref{thm:overlap-moment-concentration} and~\pref{lem:replicon-moments-from-mgf} yields the first type of term is $O(1)$, and Cauchy-Schwarz gives $O(1)$ bounds for the second type of term. \\
    We now bound the four type of higher-order remainder terms from the $\beta^2\E[\Tr[PA^2P]]$ term. Note by the repeated index identities, that
    \allowdisplaybreaks
    \begin{align*}
        \sum_b P_{bk}T_{jbj} &= \frac{n}{2}\an{(\sigma^1_k-\sigma^2_k)(\sigma^3_j-\sigma^4_j)(\sigma^3_j-\sigma^6_j)f_{1235}}\,, \\
        \sum_a P_{ja}\Gamma_{akij} &= \frac{n}{2}\an{(\sigma^1_j-\sigma^2_j)(\sigma^3_j-\sigma^7_j)(\sigma^3_k-\sigma^5_k)(\sigma^3_i-\sigma^6_i)f_{1234}} - (P^2)_{jk}P_{ij} - (P^2)_{ji}P_{kj} - (P^2)_{jj}P_{ki}\,.
    \end{align*}
    For the non-repeated third moment, the following identity holds,
    \allowdisplaybreaks
    \begin{align*}
        \sum_a P_{ja}T_{akj} &= \frac{n}{2}\an{(\sigma^1_j-\sigma^2_j)(\sigma^3_j-\sigma^6_j)(\sigma^3_k-\sigma^5_k)f_{1234}}\,.
    \end{align*}
    We will define a ``triplet'' product of individual $f$ site factors as $U_{jk}$ (which will show up repeatedly) as
    \[
        U_{ab} := (\sigma^{r_1}_a-\sigma^{r_2}_a)(\sigma^{r_3}_a - \sigma^{r_4}_a)(\sigma^{r_3}_b-\sigma^{r_5}_b)\,.
    \]
    We handle each of the four terms chronologically. For the first term, using the identity above,
    \allowdisplaybreaks
    \begin{align*}
        &\frac{\beta^4}{n}\sum_{i<j}\sum_k\E\left[A_{ki}\left(m_i\sum_a P_{ja}T_{akj}\right)\right] = \frac{\beta^4}{2}\sum_{i<j}\sum_k\E\left[A_{ik}m_i\an{U_{jk}f_{1234}}\right] \\
        &=_{\eqref{eq:GOE-ibp}} \frac{\beta^4}{2n}\sum_{i<j}\sum_k\E\left[m_i(\partial_{ki}\an{U_{jk}f_{1234}}) + \an{U_{jk}f_{1234}\partial_{ki}m_i}\right] \\
        &=\frac{\beta^4}{4n}\left(\sum_{\ell=1}^6\E\left[\an{\left(\sum_i\sigma^7_i\sigma^{\ell}_i\right)\left(\sum_{j,k}U_{jk}\sigma^\ell_k\right)f_{1234}}\right]-\sum_{\ell=1}^6\an{\sum_i\left(\sigma^7_i\sigma^\ell_i\sum_k U_{ik}\sigma^\ell_k\right)f_{1234}}\right) + \frac{\beta^4}{2n}\sum_{i<j}\sum_k \E\left[\an{U_{jk}f_{1234}\partial_{ki}m_i}\right]\,.
    \end{align*}
    Each term in the first two sums (summed over $\ell$) belongs to (where $B^n$ is a term bounded by $Cn$ and dependent only on the last spins of the replicas)
    \[
        \mathsf{span}\left\{ C(\beta)n^2\,\nu\left(f_{1234}\Xi_{cd}\Xi_{ab}\right), C(\beta)n\,\nu\left(B^n\Xi_{ab}f_{1234}\right)\right\}\,.
    \]
    The first term is $O(1)$ upon an application of~\pref{prop:D15}, and the second sum is $O(1)$ by the fact that $|B^n| \le Cn$ and an application of Cauchy--Schwarz followed by~\pref{thm:overlap-moment-concentration} and~\pref{lem:replicon-moments-from-mgf}. Similar simplifications symmetrizing the sum over $i$ and $j$, and using the facts that $\partial_{ki}m_i = \beta(P_{ii}m_k - m_iP_{ik})$ and $\sum_k (\sigma^3_k-\sigma^5_k)(Pm)_k = \frac{n^2}{2}\an{f_{3578}((R_{79}-q^*)-(R_{89}-q^*))}$, give that
    \allowdisplaybreaks
    \begin{align*}
        &\frac{\beta^4}{2n}\sum_{i<j}\sum_k \E\left[\an{U_{jk}f_{1234}\partial_{ki}m_i}\right] = \\
        & \mathsf{span}\left\{C(\beta)n^2 \nu(f_{1234}\Xi_{ab}\Xi_{cd}),\,C(\beta)n^2\nu(f_{1234}f_{abcd}f_{abef}\Xi_{gh}), C(\beta)n\nu(B^n\Xi_{ab}f_{1234})\right\}\,.
    \end{align*}
    Each of these terms can be bounded by~\pref{prop:D15}, or by invoking $(4,4,4,4)$-H\"older's or Cauchy--Schwarz with~\pref{thm:overlap-moment-concentration} and~\pref{lem:replicon-moments-from-mgf}. \\\
    For the second term, we use the repeated replica identity and~\eqref{eq:GOE-ibp} to simplify as
    \allowdisplaybreaks
    \begin{align*}
        \frac{\beta^4}{n}\sum_{i<j}\E\left[\sum_k A_{ki}m_i\sum_k P_{bk}T_{jbj}\right] = \frac{\beta^4}{2n}\sum_{i<j}\sum_k \E\left[\an{(\partial_{ki}m_i)U_{jk}f_{1235}} + m_i(\partial_{ki}\an{f_{1235}U_{jk})} \right]\,.
    \end{align*}
    This term has, by symmetry, exactly the same structure as the first term and can be bounded as $O(1)$ by the exact same argument. \\
    For the third term, we simplify using the second repeated index identity and~\eqref{eq:GOE-ibp} to obtain
    \allowdisplaybreaks
    \begin{align*}
        &\frac{\beta^4}{n}\sum_{i<j}\sum_k \E\left[A_{ki}\sum_b P_{bk}\Gamma_{jbij}\right] = -\frac{2\beta^4}{n}\sum_{i<j}\sum_k\E\left[A_{ki}P_{ji}(P^2)_{kj}\right] - \frac{2\beta^4}{n}\sum_{i<j}\sum_k \E\left[A_{ki}(m_j\sum_b P_{bk}T_{bij})\right] \\
        &= -\frac{2\beta^4}{n^2}\sum_{i<j}\sum_k \E\left[(P^2)_{kj}\partial_{ki}P_{ji} + P_{ji}\partial_{ki}(P^2)_{kj}\right] - \frac{2\beta^4}{n^2}\sum_{i<j}\sum_k \E\left[\partial_{ki}\left(m_j \sum_b P_{bk}T_{bij}\right)\right]\,.
    \end{align*}
    We first simplify the first term using~\pref{lem:dP-pairing} and the product rule along with some algebra to obtain that
    \allowdisplaybreaks
    \begin{align*}
        &-\frac{2\beta^4}{n^2}\sum_{i<j}\sum_k \E\left[(P^2)_{kj}\partial_{ki}P_{ji} + P_{ji}\partial_{ki}(P^2)_{kj}\right] = -\frac{2\beta^5}{n^2}\sum_{i<j}\sum_k \E\left[(P^2)_{kj}P_{ji}P_{ik}\right] - \frac{2\beta^5}{n^2}\sum_{i<j}\sum_k \E\left[(P^2)_{kj}P_{jk}P_{ii}\right] \\
        & - \frac{2\beta^5}{n^2}\sum_{i<j}\sum_k \E\left[(P_{ji})^2(P^2)_{kk} + (P^2)_{ki}P_{kj}P_{ji} + (P^2)_{ji}P_{kk}P_{ji} + (P^2)_{jk}P_{ki}P_{ij}\right] + \mathrm{HO}_{12} + \mathrm{HO}_{13} + \mathrm{HO}_{14} + R^{(3)}_{\text{type}}\,, 
    \end{align*}
    where
    \allowdisplaybreaks
    \begin{align*}
        \mathrm{HO}_{12} &= -\frac{2\beta^5}{n^2}\sum_{i<j}\sum_k \E\left[(P^2)_{kj}m_kT_{jii}\right]\,, \\
        \mathrm{HO}_{13} &= -\frac{2\beta^5}{n^2}\sum_{i<j}\sum_k \E\left[(P^2)_{kj}m_iT_{jik}\right]\,, \\
        \mathrm{HO}_{14} &= - \frac{2\beta^5}{n^2}\sum_{i<j}\sum_k \E\left[(P^2)_{jk}\Gamma_{jiki}\right]\,, \\
        R^{(3)}_{\text{type}} &= -\frac{2\beta^5}{n^2}\sum_{i<j}\sum_k \E\left[P_{ji}\left(m_k\left(\sum_a P_{ka}T_{aji} + \sum_a P_{ja}T_{kai}\right) + m_i\left(\sum_a P_{ka}T_{ajk} + \sum_a P_{ja}T_{kak}\right) + \sum_a P_{ka}\Gamma_{ajik} + \sum_a \Gamma_{kaki}\right)\right]\,.
    \end{align*}
    The $R^{(3)}_{\text{type}}$ term can be shown to be $O(1)$ by arguments exactly analogous to those used to bound $R^{(3)}_1, R^{(3)}_2$ and $R^{(3)}_3$. We will bound the $\mathrm{HO}_i$ terms later to be  $O(1)$. We show that the remaining ``powers-of-$P$'' terms are $O(1)$ first. Note that
    \[
        |\sum_{i<j}\sum_k (P^2)_{ji}P_{ik}| = |\sum_{i<j}P_{ji}(P^3)_{ji}| \le \frac{1}{2}\Tr[P^4]\, ,
    \]
    which implies that $\frac{2\beta^5}{n^2}\sum_{i<j}\sum_k \E\left[(P^2)_{ji}P_{ik}\right] = O(1) $ by \eqref{eq:c-p4-bound}. Similarly,
    \[
        |\sum_{i<j}\sum_k (P^2)_{kj}P_{jk}P_{ii}| = |\sum_{i<j }P_{ii}(P^3)_{jj}| \le |\sum_i P_{ii}\sum_j (P^3)_{jj}| = \Tr[P]\Tr[P^3]\,, 
    \]
    and so $\frac{C(\beta)}{n^2}\E\left[\Tr[P]\Tr[P^3]\right] = C(\beta)n\E[(n - \an{R_{78}})\an{f_{1234}f_{1256}f_{3456}]} = O(1)$ by \eqref{eq:c-p3-bound} and an invocation of a $(4,4,4,4)$-H\"older's along with~\pref{lem:replicon-moments-from-mgf}. The remaining ``powers-of-P'' terms simplify by symmetrization and trace rewrites as
    \allowdisplaybreaks
    \begin{align*}
       &- \frac{2\beta^5}{n^2}\sum_{i<j}\sum_k \E\left[(P_{ji})^2(P^2)_{kk} + (P^2)_{ki}P_{kj}P_{ji} + (P^2)_{ji}P_{kk}P_{ji} + (P^2)_{jk}P_{ki}P_{ij}\right] = \\
       &\qquad -\frac{\beta^5}{n^2}\E\left[\Tr[P^2](\Tr[P^2]-\sum_i (P_{ii})^2) + \Tr[P](\Tr[P^3]-\sum_i (P^2)_{ii}P_{ii}) + \Tr[P^4]-\sum_iP_{ii}(P^3)_{ii}\right]\,.
    \end{align*}
    Since $0 \le P_{ii} \le 1$ and $(P^k)_{ii} \ge 0$, we obtain the facts that
   \allowdisplaybreaks
    \begin{align*}
        \sum_{i}P^2_{ii} &\le \Tr[P] \le n\,,  \\
        \sum_{i} P_{ii}(P^2)_{ii} &\le \sum_i (P^2)_{ii} = \Tr[P^2] =_{\text{\eqref{eq:c-p2-bound}}} \frac{n^2}{4}\an{f^2}\,, \\
        \sum_i P_{ii}(P^3)_{ii} &\le \Tr[P^3] =_{\text{\eqref{eq:c-p3-bound}}} \frac{n^3}{8}\an{f_{1234}f_{1256}f_{3456}}\,.
    \end{align*}
    Now, note \eqref{eq:c-p2-bound} implies that $\Tr[P^2]^2 = \frac{n^4}{16}\an{f^2_{1234}f^2_{5678}}$. Putting these bounds together gives that
    \allowdisplaybreaks
    \begin{align}\label{eq:p-trace-power-bounds}
        \frac{C(\beta)}{n^2}\left|\Tr[P^2]\left(\Tr[P^2]-\sum_i (P_{ii})^2\right)\right| &\le C(\beta)n^2\an{f^2_{1234}f^2_{5678}} + C(\beta)n\an{f^2}\,, \nonumber\\
        \frac{C(\beta)}{n^2}\left|\Tr[P]\left(\Tr[P^3]-\sum_i P_{ii}(P^2)_{ii}\right)\right| &\le C(\beta)n^2\an{f_{1234}f_{1256}f_{3456}} + C(\beta)n\nu(f^2)\, \nonumber \\
        \frac{C(\beta)}{n^2}\left|\Tr[P^4]-\sum_i P_{ii}(P^3)_{ii}\right| &\le_{\text{\eqref{eq:c-p4-bound}}} C(\beta)n^2\an{f_{1234}f_{5678}f_{1278}f_{3456}} + C(\beta)n\an{f_{1234}f_{1256}f_{3456}}\,. 
    \end{align}
    Taking expectation $\E[\cdot]$ and applying~\pref{prop:D15} for the cubic terms, along with Cauchy--Schwarz and~\pref{lem:replicon-moments-from-mgf} for the other terms yields the $O(1)$ bound for the ``powers-of-P'' terms. \\
    We now bound the final contribution to this term which comes from $\frac{2\beta^4}{n^2}\sum_{i<j}\sum_k \E\left[\partial_{ki}\left(m_j\sum_b P_{kb}T_{bij}\right)\right]$ by first using \eqref{eq:GOE-ibp} and the identity that $\sum_b P_{kb}T_{bik} = \frac{n}{2}\an{\Delta^{12}_k\Delta^{35}_i\Delta^{36}_jf_{1234}}$ to obtain
    \allowdisplaybreaks
    \begin{align*}
        &\frac{2\beta^4}{n^2}\sum_{i<j}\sum_k \E\left[\partial_{ki}\left(m_j\sum_b P_{kb}T_{bij}\right)\right] =  \frac{\beta^4}{n^2}\sum_{i,j}\sum_k \E\left[\partial_{k,i}\left(m_j\sum_b P_{kb}T_{bij}\right)\right] - \frac{\beta^4}{n^2}\sum_{i,k}\E\left[\partial_{ki}\left(m_i\sum_b P_{kb}T_{bii}\right)\right]\\
        &\qquad\qquad\qquad\qquad= \frac{\beta^4}{n}\sum_{i,j}\sum_k \E\left[\partial_{ki}\an{\sigma^7_j\Delta^{12}_k\Delta^{35}_i\Delta^{36}_j f_{1234}}\right] - \frac{C(\beta)}{n^2}\sum_{i,k}\E\left[m_i\left(m_i\partial_{ki}(P^2)_{ki} + (P^2)_{ki}\partial_{ki}m_i\right)\right]\,.
    \end{align*}
    We bound the first term by grouping sums and some algebra multiplying replica indices, which yields
    \allowdisplaybreaks
    \begin{align*}
        &\frac{\beta^4}{n}\sum_{i,j}\sum_k \E\left[\partial_{ki}\an{\sigma^7_j\Delta^{12}_k\Delta^{35}_i\Delta^{36}_j f_{1234}}\right] = \frac{\beta^5}{n}\sum_{\ell=1}^7\sum_{i,j}\sum_k\E\left[\an{\sigma^7_j\Delta^{12}_k\Delta^{35}_i\Delta^{36}_jf_{1234}(\sigma^\ell_k\sigma^\ell_i-\an{\sigma_k\sigma_i})}\right] \\
        &= C(\beta) \sum_{\ell=1}^7 \sum_{i,k} \E\left[\an{(R_{37}-R_{67})\left(R_{1\ell}-R_{2\ell}\right)\left(R_{3\ell}-R_{5\ell}\right)f_{1234}}\right] \\
        &\qquad\qquad\qquad\qquad- C(\beta)\sum_{\ell}\sum_{i,k}\E\left[\an{(R_{37}-R_{67})(R_{18}-R_{28})(R_{39}-R_{59})f_{1234}}\right]\,.
    \end{align*}
    Adding and subtracting $q^*$ in all the overlap difference terms, one observes that each of these terms can be bounded by a $(4,4,4,4)$-H\"older's in conjunction with \pref{thm:overlap-moment-concentration} and \pref{lem:replicon-moments-from-mgf} or by~\pref{prop:D15} when $\ell \in \{1,2,3,5\}$.
    We bound the second term in two parts. For the $(P^2)_{ki}\partial_{ki}m_i$ contribution, note that $\partial_{ki}m_i = \beta(m_kP_{ii}-m_iP_{ki})$ and so $|\partial_{ki}m_i| \le 2\beta$ to obtain
    \allowdisplaybreaks
    \begin{align*}
        \left|\frac{C(\beta)}{n^2}\sum_{i,k}m_i(P^2)_{ki}\partial_{ki}m_i\right| \le \frac{C'(\beta)}{n^2}\sum_{i,k}|(P^2)_{ki}| \le_{\text{CS}} n\left(\sum_{i,k}(P^2)^2_{ki}\right)^{1/2} = n\sqrt{\Tr[P^4]}\,.
    \end{align*}
    Taking $\E\left[\cdot\right]$ and using Jensen's inequality in conjunction with \eqref{eq:c-p4-bound} yields $\frac{C(\beta)}{n^2}\sum_{i,k}\E\left[(P^2)_{ki}\partial_{ki}m_i\right] = O(1)$. Now for the second part, using the product rule and \pref{lem:dP-pairing} gives
    \allowdisplaybreaks
    \begin{align*}
        &\frac{C(\beta)}{n^2}\sum_{i,k}\E\left[m^2_i\partial_{ki}(P^2)_{ki}\right] = \frac{C'(\beta)}{n^2}\sum_{i,k}\E\left[m^2_i(P^2)_{ii}P_{kk} + m^2_iP_{kk}(P^2)_{ii}+2m^2_iP_{ik}(P^2)_{ki}\right] \\
        &\qquad\qquad + \underbrace{\frac{C'(\beta)}{n^2}\sum_{i,k}\E\left[m_i^3\left(\sum_a T_{kak}P_{ai} + \sum_a P_{ka}T_{aik}\right) + m_i^2m_k\left(\sum_a T_{kai}P_{ai} + \sum_a P_{ka}T_{aii}\right)\right]}_{:= \mathrm{HO}_{15}} \\
        &\qquad\qquad+ \underbrace{\frac{C(\beta)}{n^2}\sum_{i,k}\E\left[m_i^2\sum_a \Gamma_{kaki}P_{ai} + m_i^2\sum_aP_{ka}\Gamma_{aiki}\right]}_{:= \mathrm{HO}_{16}}\,. 
    \end{align*}
    Once again, we show that $\mathrm{HO}_{15}$ and $\mathrm{HO}_{16}$ are $O(1)$ later, and focus on bounding the ``product-of-P'' terms first. Note that since $0 \le P_{ii} \le 1$ and $m^2_i \le 1$, we have
    \allowdisplaybreaks
    \begin{align*}
        \sum_{i,k}\left((m^2_i\left(P^2)_{ii}P_{kk} + (P^2)_{kk}P_{ii}\right)\right) \le \sum_{i,k}\left(P^2)_{ii}P_{kk} + (P^2)_{kk}P_{ii}\right) = 2\Tr[P]\Tr[P^2]\,. 
    \end{align*}
    Since $\Tr[P] \le n$, this immediately implies that $\frac{C(\beta)}{n^2}\sum_{i,k}\E\left[(P^2)_{ii}P_{kk} + (P^2)_{kk}P_{ii}\right] \le O(1)$ by \eqref{eq:c-p2-bound}. For the remaining term, note that
    \allowdisplaybreaks
    \begin{align*}
    \sum_{i,k}m^2_i(P^2)_{ki}P_{ik} &= \sum_i m^2_i \sum_{k}(P^2)_{ik}P_{ki} = \sum_i m^2_i (P^3)_{ii} = \Tr\left[\diag(m^2)P^3\right]\,.
    \end{align*}
    Clearly $0 \le \Tr[\diag(m^2)P^3]$, and since $\diag(M) \preceq \Id_n$, $\Tr[\diag(m^2)P^3]\le \Tr[P^3]$. Combining the upper bound with \eqref{eq:c-p3-bound} immediately implies that $\frac{C(\beta)}{n^2}\sum_{i,k}\E\left[m^2_i(P^2)_{ki}P_{ik}\right] = O(1)$. \\ 
    Finally, just as the first and second third-centered cumulant terms were similar by symmetry, the final fourth-cumulant remainder term from $\beta^2\E[\Tr[PA^2P]]$ has exactly the same structure as the one bounded above (by symmetry) and can be shown to be $O(1)$ by an exactly analogous argument involving bounding similar terms.
    
    \ppart{The $\mathrm{HO}_i$ terms} The sixteen $\mathrm{HO}$ terms can be shown to be $O(1)$ using replica identities for the centered third moments and fourth moments followed by tedious symmetrization algebra, ultimately reducing to applications of~\pref{prop:D15} or H\"older's inequality depending on the type of term arising from the simplification. Since the techniques are the same as the ones used for the four terms above, but every term is only slightly different, the details are delegated to~\pref{lem:bounds-for-ho-terms}. \qedhere
\end{proof}
        
\begin{lemma}[Bounding the $\calE_R$ term as $O(1)$]\label{lem:rank-1-O1}
    Define
    \[
        \calE_R := \frac{\beta^4}{n^2}\E[\Tr[P(x_0x_0^\sT)^2P]] - \frac{2\beta^2}{n}\E[PS(x_0x_0^\sT)P] + \frac{2\beta^3}{n}\E\left[PA(x_0x_0^\sT)P\right] +  \frac{2\beta^2}{n}\E\Tr[x_0x_0^\sT P]\,.
    \]
    Then,
    \[
        \calE_R = O_{\beta,t}(1)\,.
    \]
\end{lemma}
\begin{prf}
    For the entire proof, we will set (WLOG) the plant $x_0 = \mathbf{1}_n$. Since the function $\|\hat{Q}^{-1}P - \Id\|_F^2$ is a continuous and bounded function of the replicas, the proof of~\pref{thm:mag-conc-planted-sk} implies that the gauge symmetry and Bayes rule leave it invariant under any fixed plant. Furthermore, observing that $\hat{Q}^{-1}P - \Id = \left((S-\beta A )P - \Id\right) - \frac{\beta^2}{n}x_0x_0^\sT P$, simplifies $\calE_R$ to the following
    \begin{align*}
        \calE_R &= \E\norm{\frac{\beta^2}{n}\mathbf{1}_n\mathbf{1}_n^\sT P }_F^2 - 2\E\left\langle (S - \beta A)P - \Id, \frac{\beta^2}{n}\mathbf{1}_n\mathbf{1}_n^\sT P\right\rangle \\
        &\le_{\text{Cauchy-Schwarz}} \E\norm{\frac{\beta^2}{n}\mathbf{1}_n\mathbf{1}_n^\sT P }_F^2 + 2\left(\E\norm{\frac{\beta^2}{n}\mathbf{1}_n\mathbf{1}_n^\sT P }^2_F\right)^{1/2}\left(\E\norm{\left((S-\beta A )P - \Id\right)}^2_F\right)^{1/2}\,.
    \end{align*}
    Note that~\pref{lem:diag-frob-term},~\pref{lem:bound-d-remainder},~\pref{lem:EA-term} and~\pref{lem:remainder-bound-ea} imply that 
    \[
        \E\norm{\left((S-\beta A )P - \Id\right)}^2_F = \calE_{A} + \calE_D = O_{\beta,t}(1)\,.
    \]
    Therefore, it suffices to prove that $\E\norm{\frac{\beta^2}{n}\mathbf{1}_n\mathbf{1}_n^\sT P}_F^2 = O_{\beta,t}(1)$ to conclude the proof. This will follow by applying replica identities for $P_{ij}$ and $(P^2)_{ij}$ using the proof of~\pref{lem:remainder-bound-ea} with the cubic bounds obtained in~\pref{prop:D15}. Specifically,
    \begin{align*}
        \E\norm{\frac{\beta^2}{n}\mathbf{1}_n\mathbf{1}_n^\sT P}_F^2 &= \E\left[\Tr\left[\frac{\beta^2}{n}\mathbf{1}_n\mathbf{1}_n^\sT P^2 \frac{\beta^2}{n}\mathbf{1}_n\mathbf{1}_n^\sT\right]\right] = \frac{\beta^4}{n}\E\left[\Tr\left[\mathbf{1}_n\mathbf{1}_n^\sT P^2\right]\right] \\
        &= \frac{\beta^4}{n}\E\left[\sum_{i,j}(P^2)_{ij}\right] = \frac{\beta^4}{4}\E\left[\an{\sum_{ij}\left(\sigma^1_i-\sigma^2_i\right)\left(\sigma^3_i-\sigma^4_i\right)f_{1234}}\right] \\
        &= \sum_{(a,b)\in \{(1,3),(1,4),(2,3),(2,4)\}}\frac{\beta^4 n^2}{4}\E\left[\an{\left(M_a-m^*\right)\left(M_b-m^*\right)f_{1234}}\right] =_{\text{\pref{prop:D15}}} O(1)\, ,
    \end{align*}
    where we used the fact that 
    \begin{align*}
        \nu((M_a - m^*)(M_b - m^*)f_{1234}) &= \nu((M^-_a-m^*)(M^-_b-m^*)f_{1234}) + \frac{1}{n}\nu((M^-_a-m^*)\eps^b f_{1234}) \\
        &+ \frac{1}{n}\nu((M^-_b-m^*)\eps^a f_{1234}) + \frac{1}{n^2}\nu(\eps^a\eps^b f_{1234})\,, 
    \end{align*}
    and the first term can be bounded directly by~\pref{prop:D15}, the second and third term using Cauchy-Schwarz along with the facts that $|\eps^a| \le 1$, $|\eps^b| \le 1$,~\pref{lem:replicon-moments-from-mgf} and overlap/magnetization concentration, and the last term is calculated by trivially using the boundedness of $|\eps^a\eps^b| \le 1$. \qedhere
\end{prf}

\begin{corollary}[Final bound on the Frobenius norm]\label{cor:final-frobenius-bound}
    The following bounds hold for the three terms constituting the Frobenius norm,
    \[
        \calE_A + \calE_D + \calE_R = O(1)\,.
    \]
\end{corollary}
\begin{prf}
    The proof follows from \pref{lem:diag-frob-term}, \pref{lem:bound-d-remainder}, \pref{lem:EA-term}, \pref{lem:remainder-bound-ea} and \pref{lem:rank-1-O1}.
\end{prf}

\section{Exact moments estimates via planted cavity interpolation on the Nishimori line}\label{sec:cavity-interpolations-trace-identities}
In this section, we first state the main concentration estimates of the magnetization and overlap (\pref{thm:overlap-moment-concentration}, \pref{thm:mag-conc-planted-sk}) that are known for the planted SK model under SL tilt up to $\beta < 1$. These mainly reduce to showing that a mapping from the planted SK model under SL tilt to a Guerra interpolated model in \cite[\S 4]{alaoui2017finite} preserves certain key properties. We then establish moment bounds (under the disorder) on every entry of the diagonal matrix $D$ that appears in the Hessian of the TAP free energy (\pref{lem:D-moments}). With access to these lemmata and a few others regarding replica algebra (\pref{sec:app-a}), we introduce the cavity interpolation for the planted SK model with SL tilt (\pref{sec:cavity-interpolation}). Under this interpolation, we systematically develop various exact moment estimates and strengthened concentration bounds, ultimately culminating in the four \emph{critical} results of this section. These are:
\begin{itemize}[itemsep=0.3em]
    \item \pref{prop:D17} and \pref{lem:cTrPDP-identity}, which provide \emph{exact} estimates for $O(1/n)$ terms that appear in $\E\norm{\hat{Q}^{-1}P - \Id_n}_F^2$ to account for the ``heavy'' cancellations, and
    \item \pref{prop:D15} and \pref{lem:D-oneleg-O(n)}, which are concentration estimates that beat ``naive'' H\"older bounds by using extra cancellations. These are critically used to bound various remainder terms that come from an interaction with the diagonal matrix $D$ or mixed cubic moments of ``bulk'' deviation terms.
\end{itemize}

\subsection{Overlap concentration in the planted SK Gibbs measure}\label{subsec:multioverlaps-planted-sk} We state a critical theorem in this section, which implies that the moments of the overlap variance for the Gibbs measure of the planted SK model can be bounded by a constant that depends only on $\beta$ and the moment parameter $k$. 

\begin{theorem}[Overlap concentration for the planted SK model,~{\cite[Mild extension of Theorem 7]{alaoui2017finite}}]\label{thm:overlap-moment-concentration}
    For $G_{A,\beta,y_t}$, $0 < \beta < 1$, and every positive integer $k \in \Z_+$,
    \[
        \E\an{(R_{12}-q^*)^{2k}} \le \frac{C(k)}{n^k} + C(k)e^{-c(k,t)N}\, , 
    \]
    for positive constants $C(k)$ and $c(k,t)$, where
    \[
        q^* = \E_{h \sim \calN(0,1)}\left[\tanh^2\left(\beta^2 m^* + t + \sqrt{\beta^2q^* + t}h\right)\right]\, ,
    \]
    and
    \[
        m^* = \E_{h \sim \calN(0,1)}\left[\tanh\left(\beta^2 m^* + t + \sqrt{\beta^2 q^* + t}h\right)\right]\,.
    \]
\end{theorem}
\begin{prf}
    This is a mild generalization of~\cite[Theorem 7]{alaoui2017finite} and the remark on~\cite[Pg. 13]{alaoui2017finite} after the statement of~\cite[Theorem 8]{alaoui2017finite}. El Alaoui, Krzakala and Jordan consider a general setup that handles the planted SK model without SL tilt and observe that their Guerra interpolation over $s \in [0,1]$ preserves a \emph{Nishimori property} such that the time-$s$ interpolated Gibbs measure corresponds to the posterior distribution of the plant vector given \emph{additional} scalar side channel information 
    \[ y_i \sim \calN\left((1-s)r, (1-s)r\right)\,, \]
    where $r= \lambda q^*(\lambda)$, with $q^*(\lambda)$ being the unique solution achieving $\sup_{q \ge 0} F(\lambda,q)$ \cite[(3) \& \S 2.3]{LM19}. \\   
    We quickly show that for every $0 < \beta < 1$ and $t \in [0,T(\beta)]$ there exists a unique pair $(\lambda, s)$ that satisfies the the following requirements:
    \begin{enumerate}[itemsep=0.3em]
        \item $s\in [0,1]$\,,
        \item $\lambda s = \beta^2$\,,
        \item $q^*(\lambda) = q^*$\,,
        \item $0 < \beta^2 \le \lambda < \infty$ (and so $\lambda \in \calD$)\,,
        \item $r=\lambda q^*(\lambda)$\,,
        \item $t = (1-s)r$\,, and
        \item $H_s \overset{d}{=} H_{\text{planted}}$\,.
    \end{enumerate}
    \ppart{Parameter map} The map sets $(\beta, t) \to \left(\lambda = \beta^2 + \frac{t}{q^*}, s = \frac{\beta^2 q^*}{\beta^2 q^* + t}\right)$ when $t > 0$. For the case of no tilt at $t= 0$, we map $(\beta, t) \to (\lambda = \beta^2, s=1)$. In both cases, $r := \beta^2 q^* + t$. Since $0 \le t \le T(\beta) < \infty$, it is immediate that $s \in [0,1)$. Furthermore, $\lambda s = \beta^2$ for every $t \in [0,T(\beta)]$.

    \ppart{Fixed-point equivalence} The case of $t = 0$ and $t > 0$ are handled separately. 
    \begin{enumerate}[itemsep=0.3em]
        \item When $t=0$, note that $\lambda = \beta^2 < 1$ and so the $q^*(\lambda)$ that extermizes $F_{\lambda,q}$ is at $q^*(\lambda)=0$ \cite[\S 2.3]{LM19}. For $t=0$, it is also easy to see that $q^* = \E_h\left[\tanh^2\left(\beta^2m^* + \sqrt{\beta^2q^*}h\right)\right]$. Since $m^*=q^*$ by \pref{lem:nishimori-condition}, it is immediate that $q^* = 0$. This yields $q^*(\lambda) = q^*$. 
        \item At $t > 0$, we use \pref{lem:nishimori-condition} and the fact that $q^* = f(\beta^2q^* + t) < \beta^2q^* + t = r = \lambda q^*$ to conclude that $\lambda > 1$ and so \cite[(9) \& \S 2.3]{LM19} implies that $q^*(\lambda) = \E[X^2_0] - \left(1-\E\left[\tanh^2(\langle \lambda q^*(\lambda) + \sqrt{\lambda q^*(\lambda)}h)\right]\right)$. Note that $x_0 \sim \mathsf{Unif}\left(\{-1,1\}^n\right)$ corresponds to the case that $\E[X_0^2]=1$ which immediately implies that $q^*(\lambda) = q^*$. 
    \end{enumerate} 
    We now show that $f(r) = \E_h[\tanh^2(r + \sqrt{r}h)] < r$. Note that, for $g(x) = x - \tanh(x)$ (which is odd) we have
    \allowdisplaybreaks
    \begin{align*}
        g'(x) = 1 - \sech^2(x) = \tanh^2(x) \ge 0\, ,
    \end{align*}
    which implies that $g(x) > 0$ when $x > 0$, and $g(x) < 0$ for $x < 0$. Then, for a Gaussian $Y \sim \calN(r,r)$ with density $\rho_r(x)dx$, we obtain that
    \[
        \E_Y\left[g(Y)\right] = \int_0^\infty g(x)\left(\rho_r(x) - \rho_r(-x)\right)\,dx\,,
    \]
    and since 
    \[
        \frac{\rho_r(x)}{\rho_r(-x)} = e^{2x} > 1\, ,
    \]
    and $g(x) > 0$ for $x > 0$, we immediately obtain $\E_Y[g(Y)] = \E[Y - \tanh(Y)] = r - \E_Y[g_Y] > 0$. By the proof of \pref{lem:nishimori-condition}, note that $\E_Y[g(Y)] = \E_Y[g^2(Y)] = f(r)$, which immediately gives $f(r) < r$ for $r > 0$.
    
    \ppart{Invariance of constraints on $r$ and $t$} Once again, we verify the constraints separately for $t = 0$ and $t > 0$.
    \begin{enumerate}[itemsep=0.3em]
        \item At $t=0$, note that $\lambda = \beta^2$ and the mapping sets $r = \beta^2 q^* + 0 = 0$ (since $q^* = 0$ at $t=0$). Since $\lambda = \beta^2$ and $q^*(\lambda) = q^*$, we also obtain that $r = \lambda q^*(\lambda) = \beta^2\cdot 0 = 0$, which preserves the desired invariance. Similarly, with $r = 0$, we obtain that $t = (1-s)r = 0$, which is also consistent.
        \item At $t > 0$, we have that the mapping sets $\lambda = \beta^2 + t/q^*$ and $r = \beta^2 q^* + t = (\beta^2 + t/q^*)q^* = \lambda q^* = \lambda q^*(\lambda)$, where the last equality follows because $q^* = q^*(\lambda)$. Similarly, $(1-s)r = (1-s)\lambda q^*(\lambda) = \left(1-\frac{\beta^2q^*}{\beta^2q^* + t}\right)(\beta^2 + t/q^*)q^* = t$, as required.
    \end{enumerate}

    \ppart{Distributional equivalence of the Hamiltonians} Given a fixed plant $x_0 \in \{-1,1\}^n$, note that by the distributional equivalences used for $(\mathsf{diag}(x_0)J\mathsf{diag}(x_0),\mathsf{diag}(x_0)y_t,\mathsf{diag}(x_0)\,\sigma)$ in the proof of \pref{thm:mag-conc-planted-sk}, we have that 
    \[
        H_{\text{planted}}(\tau) \overset{d}{=} \frac{\beta}{\sqrt{n}}\sum_{i<j}A_{ij}\tau_i\tau_j + \frac{\beta^2}{2n}\left(\sum_i u_i\right)^2 + t\sum_i u_i + \sqrt{t}\sum_i z_iu_i\,,
    \]
    where $\tau = \diag(x_0)\sigma$, $z = (z_1,\dots,z_n) \sim \calN(0,\Id_n)$, and we remove the diagonal entries because of \pref{lem:diag-does-not-enter-gibbs}. Furthermore, by \cite[(22)]{alaoui2017finite}, we have
    \[
        H_s(\tau) \overset{d}{=} \sqrt{\frac{s\lambda}{n}} \sum_{i<j}A_{ij}\tau_i\tau_j + \frac{s\lambda}{n}\sum_{i<j}u_iu_j + \sqrt{(1-s)r}\sum_i z_i u_i + (1-s)r\sum_i u_i\, , 
    \]
    after the dropping the $u_i$-independent constant terms (using the fact that $u_i^2 = 1$). After symmetrizing the sum over $i$ and $j$, ignoring constants, and using $s\lambda = \beta^2$ and $(1-s)r = t$, we immediately obtain that $H_s(\tau) \overset{d}{=} H_{\text{planted}}(\tau)$. \\
    Given the proof of the requirements above, it immediately follows that \cite[Theorem 7]{alaoui2017finite} and the remark on \cite[Pg. 13]{alaoui2017finite} apply verbatim to $H_{\text{planted}}$.
\end{prf}

\subsection{Magnetization concentration in the planted SK Gibbs measure} This is obtained as a consequence of the overlap concentration achieved by the Guerra interpolation in~\cite[\S 4]{alaoui2017finite}. The key insight is that the symmetry implied by the map $\sigma \to x_0 \circ \sigma$ being a bijection on $\{-1,1\}^n$ combines with the Nishimori identity ($m^* = q^*$) and the Bayes rule to give magnetization concentration. This is, in fact, a remark made in~\cite[\S 3]{alaoui2017finite} and we prove it for completeness.

\begin{theorem}[Magnetization concentration for the planted SK model]\label{thm:mag-conc-planted-sk}
    Fix the plant $x_0 \in \{-1,1\}^n$ and let $\beta \in (0,1/2)$. Then, setting  $M := \frac{1}{n}\sum_{i=1}^n (x_0)_i\sigma_i$, under the Nishimori condition that $m^* = q^*$ where $(m^*,q^*)$ solve the fixed-point equations in~\pref{thm:overlap-moment-concentration},  for every positive integer $k \in \Z_+$, there is a constant $C(\beta, k)$ so that
    \[
        \E\an{(M - m^*)^{2k}} \le \frac{C(k)}{n^k} + C(\beta,k)e^{-c(k,t)N}\,.
    \]
\end{theorem}
\begin{prf}
    Note that $M = R_{1,x_0}$, which is the overlap between a single replica and the planted vector $x_0$. Denote the contribution to the quadratic term for every $ i < j$ as
    \[
        J_{ij} := \beta A_{ij} + \frac{\beta^2}{n}(x_0)_i(x_0)_j\,.
    \]
    Since, for $i = j$, $\sigma_i\sigma_j=\sigma^2_i=1$, by~\pref{lem:diag-does-not-enter-gibbs} we have that $\E\an{f}$ is unchanged under the diagonal entries of $A$ for any measurable function $f$. Given the fixed plant $x_0$, the following distributional equivalences hold
    \begin{align*}
        J_{ij} &\sim \calN\left(\frac{\beta^2}{n}(x_0)_i(x_0)_j, \frac{\beta^2}{n}\right)\,, \\
        (y_t)_i &\sim \calN\left(t(x_0)_i, t\right)\,.
    \end{align*}

    \ppart{Equivalence for a uniform plant via Bayes' rule}
    At this point, choose $x \sim \mathsf{Unif}\left(\{-1,1\}^n\right)$, and conditioned on the disorder $\{J_{ij}\}_{i<j}$ and the field $y_t$, observe the following via Bayes' rule
    \begin{align*}
        \P\left\{x = \sigma \mid \{J_{ij}\}_{i<j}, y_t\right\} &= \frac{2^{-n}\P\left\{\{J_{ij}\}_{i<j}, y_t \mid x = \sigma\right\}}{2^{-n}\sum_x \P\left\{\{J_{ij}\}_{i<j}, y_t \mid x = \sigma\right\}} = \frac{2^{-n} e^{\iprod{\sigma, J\sigma} + \iprod{y_t,\sigma}}}{Z}\,.
    \end{align*}
    Consequently, given iid replicas $(\sigma^1,\sigma^2)$ drawn from  $\P\left\{x \mid \{J_{ij}\}_{i<j}, y_t\right\}$,
    \[
        \left(J,y_t,\sigma^1,\sigma^2\right) \overset{d}{=} \left(J,y_t,\sigma^1,x\right)\,,
    \]
    which follows immediately by~\cite[Lemma 1]{barbier2021overlap}. Since $(M - m^*)^{2k}$ is a continuous and bounded function, it immediately implies that
    \[
        \E_{x}\E_A\an{(R_{12}-q^*)^{2k}} = \E_x\E_{A}\an{(M - m^*)^{2k}} \le C(\beta,k)\,,
    \]
    since $m^* = q^*$ (\pref{lem:nishimori-condition}).

    \ppart{Equivalence for a fixed plant via gauge symmetry} Using the closure of the diagonal algebra and the invariance of Gaussians to signs, we can extend the result above to \emph{every} fixed plant $x_0$. Set
    \allowdisplaybreaks
    \begin{align*}
        J' &:= \mathsf{diag}\left(x_0\right)J\mathsf{diag}\left(x_0\right)\,, \\
        y'_t &:= \mathsf{diag}\left(x_0\right)y_t\,, \\
        \tau &:= \mathsf{diag}\left(x_0\right)\sigma\,.
    \end{align*}
    The matrix $\mathsf{diag}(x_0)$ is clearly its own inverse, and consequently
    \[
        \iprod{\sigma, J\sigma} + \iprod{y_t,\sigma} = \iprod{\tau,J'\tau} + \iprod{h',\tau}\,.
    \]
    Conditioning on $x = x_0$ gives the following law for $J' $ and $y'_t$
    \begin{align*}
        J' &= \mathsf{diag}(x_0)\frac{\beta^2}{n}xx^\sT\diag(x_0) + \frac{\beta}{\sqrt{n}}\diag(x_0)A\diag(x_0) \\
        &\overset{d}{=} \frac{\beta^2}{n}\mathbf{1}_n\mathbf{1}_n^\sT + \frac{\beta}{\sqrt{n}}A'\, , 
    \end{align*}
    where $A' \overset{d}{=} A \sim \mathsf{GOE}(n)$, and
    \[
        y'_t = \diag(x_0)y_t \overset{d}{=} t\mathbf{1}_n + B'_t\, ,
    \]
    where $B'_t \overset{d}{=} B_t$. This immediately implies that the following conditional law is the same,
    \[
        (J', y'_t) \mid x = x_0 \overset{d}{=} (J,y_t) \mid x = \mathbf{1}_n\,.
    \]
    Since $x \sim \mathsf{Unif}\left(\{-1,1\}^n\right)$, one can substitute the distributional equivalence above into the Nishimori expectation equivalence over $x$ to conclude that
    \[
        \E_{A}\an{(M - m^*)^{2k}} = \E_{A}\an{ (R_{1,2}-m^*)^{2k}} \le C(\beta,k)\, ,
    \]
    for \emph{every} $x_0 \in \{-1,1\}^n$. \qedhere
\end{prf}

\subsection{Uniform bounds on bulk statistics of $D(\mg)$} It will be critical to have cavity interpolations for bulk observables that are weighted by terms of $D(\mg)$, which will require uniform control over the moments of the entries of the diagonal matrix.  This will be useful in various applications of H\"older's inequality in statements that control odd moments of the form $\nu_s(\hat{f}^- G)$, where $G$ is a quadratic term involving centered overlaps or magnetizations.
\begin{lemma}[Uniform-in-$s$ moments of $D_{ii}$]\label{lem:D-moments}
Fix $t\in[0,T]$ and $0<\beta<1/2$. For each $s\in[0,1]$ let $\langle\cdot\rangle_s$
denote the Gibbs average along the cavity interpolation Hamiltonian defined in~\pref{sec:cavity-interpolation}. Set
\[
m_i(s):=\langle\sigma_i\rangle_s,\qquad D_{ii}(s):=\frac{1}{1-m_i(s)^2} = \cosh^2(\operatorname{arctanh}(m_i(s)))\,.
\]
Then, for every $p\in \Z_+$, there exists a constant $C_{p,\beta,t}$ such that
\begin{equation}\label{eq:D-moments-uniform-s}
    \sup_{s\in[0,1]}\ \sup_{i\in[n]}\ \E\big[D_{ii}(s)^p\big]\ \le\ C_{p,\beta,t}.
\end{equation}
Consequently,
\[
\sup_{s\in[0,1]}\ \E\Big[\Tr\big(D(m(s))^p\big)\Big]
=
\sup_{s\in[0,1]}\ \sum_{i=1}^n \E\big[D_{ii}(s)^p\big]
\le C_{p,\beta,t}\,n.
\]
\end{lemma}

\begin{proof}
Let the time of the interpolation $s\in[0,1]$, and fix a spin $i\in[n]$. Split the Hamiltonian to isolate the cavity field for $i$ as
\begin{equation}\label{eq:cavity-decomp-general}
\beta H_s(\sigma)
=
\beta H^{(i)}_s(\sigma_{-i})
+\sigma_i\Big(y_{t,i}+H^{\mathrm{row}}_{i,s}(\sigma_{-i})+H^{\mathrm{bd}}_{i}(\sigma_{-i})\Big),
\end{equation}
where, conditional on all disorder variables not in the $i$-th row/column,
\[
H^{\mathrm{row}}_{i,s}(\sigma_{-i})=2\sqrt{s}\,\beta\sum_{k\neq i}(A_{ik})\sigma_k
\quad\text{with}\quad (A_{ik})_{k\neq i}\overset{iid}{\sim} \mathcal N(0,1/n),
\]
and the remaining $\sigma_i$--linear contribution is uniformly bounded:
\begin{equation}\label{eq:Hi-bd}
\sup_{\sigma_{-i}\in\{\pm1\}^{n-1}}\big|H^{\mathrm{bd}}_{i}(\sigma_{-i})\big|
\le C_{\beta,t}<\infty\,.
\end{equation}
Note that with a fixed plant $x_0\equiv \mathbf 1$, one can take $H^{\mathrm{bd}}_{i}(\sigma_{-i})=(\beta^2/n)\sum_{k\neq i}\sigma_k$. This immediately implies that $C_{\beta,t}=\beta^2$.

\ppart{Single-spin marginal and the $\cosh^2$ representation}
Define the constrained (bulk) partition functions
\[
Z_{i,s}(\pm)
:=
\sum_{\sigma_{-i}\in\{\pm1\}^{n-1}}
\exp\Big(\beta H^{(i)}_s(\sigma_{-i}) \pm \big(H^{\mathrm{row}}_{i,s}(\sigma_{-i})+H^{\mathrm{bd}}_{i}(\sigma_{-i})\big)\Big).
\]
Then, the computation of the single-spin marginal using the identity $\tan(x) = (e^{2x}-1)/(e^{2x}+1)$ gives
\[
    m_i(s)=\langle \sigma_i\rangle_s=\tanh\big(y_{t,i}+U_{i,s}\big)\,
\]
where
\[
    U_{i,s}:=\frac12\log\frac{Z_{i,s}(+)}{Z_{i,s}(-)}\,.
\]
This immediately gives
\[
D_{ii}(s)=\sech^{-2}\big(y_{t,i}+U_{i,s}\big)=\cosh^2\big(y_{t,i}+U_{i,s}\big)\, ,
\]
and the following upper bound for any $p$-th power,
\begin{equation}\label{eq:D-moment-reduce-exp}
D_{ii}(s)^p=\cosh^{2p}\big(y_{t,i}+U_{i,s}\big)
\le \exp\big(2p\,|y_{t,i}+U_{i,s}|\big)
\le \exp\big(2p|y_{t,i}|\big)\exp\big(2p|U_{i,s}|\big).
\end{equation}
Note that the Gaussians in $y_{t,i}$ and $U_{i,s}$ are independent of each other, and so it suffices to show that $y_{t,i}$ and $U_{i,s}$ are sub-Gaussian random variables to ensure that the MGFs in the RHS of the inequality above are bounded.

\ppart{Bounding $U_{i,s}$ via conditional Gaussian concentration}
Condition on all the disorder \emph{except} the Gaussian row $(A_{ik})_{k\neq i}$, and exclude $y_{t,i}$.
Under this conditioning, the vector
\[
g:=(g_k)_{k\neq i},\qquad g_k:=\sqrt{n}\,A_{ik},
\]
is standard Gaussian in $\R^{n-1}$, and $U_{i,s}$ is a smooth function of $g$. 

Let $\langle\cdot\rangle_{i,s}^{\pm}$ denote the Gibbs average on $\sigma_{-i}$ with weight proportional to
$\exp(\beta H^{(i)}_s(\sigma_{-i}) \pm (H^{\mathrm{row}}_{i,s}+H^{\mathrm{bd}}_i))$.
Differentiating $\log Z_{i,s}(\pm)$ w.r.t.\ $A_{ik}$ gives
\[
\partial_{A_{ik}} \log Z_{i,s}(+) = 2\sqrt{s}\,\beta\,\langle \sigma_k\rangle_{i,s}^{+},
\qquad
\partial_{A_{ik}} \log Z_{i,s}(-) = -2\sqrt{s}\,\beta\,\langle \sigma_k\rangle_{i,s}^{-}.
\]
Therefore
\[
\partial_{A_{ik}} U_{i,s}
=
\frac12\Big(\partial_{A_{ik}} \log Z_{i,s}(+) - \partial_{A_{ik}} \log Z_{i,s}(-)\Big)
= 2\sqrt{s}\,\beta\Big(\langle \sigma_k\rangle_{i,s}^{+}+\langle \sigma_k\rangle_{i,s}^{-}\Big),
\]
and since $|\langle\sigma_k\rangle_{i,s}^{\pm}|\le1$ we have the uniform bound
\[
\big|\partial_{A_{ik}} U_{i,s}\big|\le 2\sqrt{s}\,\beta\le 2\beta.
\]
Passing to the standard Gaussian coordinates $g_k=\sqrt{n}A_{ik}$ yields the following per-coordinate gradient bound,
\[
\big|\partial_{g_k} U_{i,s}\big|
=
\frac{1}{\sqrt{n}}\big|\partial_{A_{ik}} U_{i,s}\big|
\le \frac{2\beta}{\sqrt{n}}.
\]
This gives the following bound for the Lipschitz-constant of the gradient,
\[
\|\nabla_g U_{i,s}\|_2^2
=
\sum_{k\neq i} \big(\partial_{g_k} U_{i,s}\big)^2
\le
\sum_{k\neq i}\frac{4\beta^2}{n}
\le 4\beta^2.
\]
By concentration for Lipschitz functions of standard Gaussians~\cite[Theorem 5.5]{boucheron2012},
this implies that conditioning on all $\{A_{kj}\}_{k\ne i, j \ne i}$,
\begin{equation}\label{eq:Ui-subg}
\E\Big[\exp\big(\lambda \left(U_{i,s}-\E[U_{i,s}\mid \{A_{kj}\}_{k\ne i, j \ne i}]\right)\big)\,\Big|\,\{A_{kj}\}_{k\ne i, j \ne i}\Big]
\le \exp\Big(2\lambda^2\beta^2\Big)\,,
\end{equation}
for every $\lambda \in \R_{\ge 0}$.
In particular, for any $\alpha>0$, using $e^{\alpha|x|}\le e^{\alpha x}+e^{-\alpha x}$ and \eqref{eq:Ui-subg},
\begin{equation}\label{eq:Ui-exp-moment-around-mean}
\E\Big[\exp\big(\alpha \left|U_{i,s}-\E[U_{i,s}\mid\{A_{kj}\}_{k\ne i, j \ne i}]\right|\big)\,\Big|\,\{A_{kj}\}_{k\ne i, j \ne i}\Big]
\le 2\exp\Big(2\alpha^2\beta^2\Big).
\end{equation}

\ppart{Bounding the conditional mean via the bounded $\sigma_i$--linear remainder}
Let $U^{(0)}_{i,s}$ denote the same constrained free energy ratio as $U_{i,s}$ but with $H^{\mathrm{bd}}_i\equiv 0$, which comes by removing the ``planted'' cavity term from \eqref{eq:cavity-decomp-general}. Then $Z^{(0)}_{i,s}(\pm)$ satisfy
\[
Z^{(0)}_{i,s}(+;A_{i\cdot})=Z^{(0)}_{i,s}(-;-A_{i\cdot}),
\]
and therefore $U^{(0)}_{i,s}(-A_{i\cdot})=-U^{(0)}_{i,s}(A_{i\cdot})$.
Since $(A_{ij})_{i=k\text{ or }j=k}$ are centered Gaussians, it follows that, conditional on $\{A_{kj}\}_{k\ne i, j \ne i}$,
\begin{equation}\label{eq:Ui0-mean0}
\E\big[U^{(0)}_{i,s}\mid \{A_{kj}\}_{k\ne i, j \ne i}\big]=0\, ,
\end{equation}
since it is an odd function.

Now, we compare $U_{i,s}$ and $U^{(0)}_{i,s}$.
Using \eqref{eq:Hi-bd}, for every configuration $\sigma_{-i}$, one obtains that
$e^{-C_{\beta,t}}\le e^{\pm H^{\mathrm{bd}}_i(\sigma_{-i})}\le e^{C_{\beta,t}}$.
Therefore,
\[
e^{-C_{\beta,t}} Z^{(0)}_{i,s}(+)\le Z_{i,s}(+)\le e^{C_{\beta,t}} Z^{(0)}_{i,s}(+),
\qquad
e^{-C_{\beta,t}} Z^{(0)}_{i,s}(-)\le Z_{i,s}(-)\le e^{C_{\beta,t}} Z^{(0)}_{i,s}(-).
\]
Taking ratios and logs gives the deterministic bound
\begin{equation}\label{eq:Ui-vs-Ui0}
\big|U_{i,s}-U^{(0)}_{i,s}\big|
=
\frac12\left|\log\frac{Z_{i,s}(+)/Z^{(0)}_{i,s}(+)}{Z_{i,s}(-)/Z^{(0)}_{i,s}(-)}\right|
\le C_{\beta,t}.
\end{equation}
Combining \eqref{eq:Ui0-mean0} and \eqref{eq:Ui-vs-Ui0} via a triangle inequality gives the following bound on the conditional mean,
\begin{equation}\label{eq:Ui-mean-bd}
\big|\E[U_{i,s}\mid \{A_{kj}\}_{k\ne i, j \ne i}]\big|
\le
\E\big[|U_{i,s}-U^{(0)}_{i,s}|\mid \{A_{kj}\}_{k\ne i, j \ne i}\big]
\le C_{\beta,t}.
\end{equation}

\ppart{Exponential moment bound for $|U_{i,s}|$} By the triangle inequality,
\[
|U_{i,s}|
\le
|U_{i,s}-\E[U_{i,s}\mid\{A_{kj}\}_{k\ne i, j \ne i}]| + |\E[U_{i,s}\mid\{A_{kj}\}_{k\ne i, j \ne i}]|.
\]
Hence, for any $\alpha>0$, using \eqref{eq:Ui-mean-bd} and then \eqref{eq:Ui-exp-moment-around-mean},
\[
\E\Big[e^{\alpha |U_{i,s}|}\,\Big|\,\{A_{kj}\}_{k\ne i, j \ne i}\Big]
\le
e^{\alpha C_{\beta,t}}\,
\E\Big[e^{\alpha |U_{i,s}-\E[U_{i,s}\mid\{A_{kj}\}_{k\ne i, j \ne i}]|}\,\Big|\,\{A_{kj}\}_{k\ne i, j \ne i}\Big]
\le
2\,\exp\Big(\alpha C_{\beta,t}+\frac{\alpha^2\beta^2}{2}\Big).
\]
Taking expectation over $\{A_{kj}\}_{k\ne i, j \ne i}$ gives the unconditional bound
\begin{equation}\label{eq:Ui-exp-moment}
\E\big[e^{\alpha |U_{i,s}|}\big]
\le
2\,\exp\Big(\alpha C_{\beta,t}+\frac{\alpha^2\beta^2}{2}\Big),
\qquad \forall \alpha>0\,.
\end{equation}

Using \eqref{eq:D-moment-reduce-exp}, the independence of $y_{t,i}$ from $\{A_{ij}\}$,
and \eqref{eq:Ui-exp-moment} with $\alpha=2p$,
\[
    \E\big[D_{ii}(s)^p\big]
    \le
    \E\big[e^{2p|y_{t,i}|}\big]\;\E\big[e^{2p|U_{i,s}|}\big]
    \le
    \E\big[e^{2p|y_{t,i}|}\big] \left(2\exp\Big(2p\,C_{\beta,t}+2p^2\beta^2\Big)\right)
    =: C_{p,\beta,t}\,.\qedhere
\]
\end{proof}

\subsection{Exact overlap moment estimates via a modified cavity interpolation}\label{sec:cavity-interpolation}

Since $\hat R_{\ell,\ell'}$ is centered replica-by-replica (sitewise, for each disorder), all mixed covariances that force Talagrand to introduce the auxiliary variables $T_\ell,T$ in~\cite[\S 1.8]{talagrand2010mean} vanish \emph{exactly} as in~\pref{lem:orthog-hatR}. This substantially simplifies the linear system for joint fluctuations: the only nontrivial constant that survives in the cavity recursions is the replicon coefficient $\rho=\E[\sech^4(Y^*)]$. Thus, to generalize~\cite[\S 1.8]{talagrand2010mean} to the setting of the planted SK model under the SL tilt, we work directly with the family $\left\{\sqrt{N}\hat R_{\ell,\ell'}\right\}_{\ell<\ell'}$.

\paragraph{Rectangular sums for recentered and  rescaled overlaps} Let
\[
\beta H_N(\sigma)
:= \frac{\beta}{2}\langle \sigma, A\sigma\rangle
   + \frac{\beta^2}{2N}\Big(\sum_{i=1}^N \sigma_i\Big)^2
   + \langle y_t,\sigma\rangle,
\qquad \sigma\in\{-1,+1\}^N,
\]
where $A \sim\mathsf{GOE}(n)$, $y_t \sim \calN(t,t)$ 
and denote the (expected) Gibbs average under $G_{A,y_t,\beta}$ as $\nu(\cdot):=\E\langle\cdot\rangle$. Set $0<\beta<1/2$ and assume:
\begin{enumerate}[itemsep=0.3em]
\item[(H1)] \textbf{Overlap concentration:} For every $k \in \Z_+$,
$\E\an{\left(R_{1,2}-q^*\right)^{2k}}\le C(\beta,k)/n^k $.
\item[(H2)] \textbf{Magnetization concentration:} An analogous bound holds for
the global magnetization $M:=\frac1N\sum_i\sigma_i$ around $m^*$.
\item[(H3)] \textbf{Nishimori condition:} $m^*=q^*$ and $(m^*,q^*)$ solve the fixed-point equations in~\pref{thm:overlap-moment-concentration}.
\end{enumerate}

For each site,
\[
m_i := \langle \sigma_i\rangle,\qquad s_i^2:=1-m_i^2,\qquad
\hat\sigma_i := \frac{\sigma_i-m_i}{s_i}.
\]
For replicas $\sigma^1,\sigma^2,\dots$ sampled i.i.d.\ from $G_{A,y_t,\beta}$, set
\[
\hat R_{\ell,\ell'} := \frac1N\sum_{i=1}^N \hat\sigma_i^\ell\hat\sigma_i^{\ell'}.
\]
Define the two ``rectangular sums'' as in~\cite[\S 1.8]{talagrand2010mean}:
\begin{align}
f &:= R_{1,3}-R_{1,4}-R_{2,3}+R_{2,4}
= \frac1N\sum_{i=1}^N a_i,
\qquad a_i := (\sigma_i^1-\sigma_i^2)(\sigma_i^3-\sigma_i^4), \label{eq:def-f}\\
\hat f &:= \hat R_{1,3}-\hat R_{1,4}-\hat R_{2,3}+\hat R_{2,4}
= \frac1N\sum_{i=1}^N \hat a_i,
\qquad \hat a_i := (\hat\sigma_i^1-\hat\sigma_i^2)(\hat\sigma_i^3-\hat\sigma_i^4)\,, \label{eq:def-fhat}
\end{align}
and similarly for the ``cavity'' overlaps,
\[
 f^-:= R^-_{1,3}- R^-_{1,4}- R^-_{2,3}+ R^-_{2,4},
\qquad
\hat f^-:=\hat R^-_{1,3}-\hat R^-_{1,4}-\hat R^-_{2,3}+\hat R^-_{2,4},
\]
where 
$\hat R^-_{\ell,\ell'}=\frac{1}{n}\sum_{i\le n-1}\hat\rho_i^\ell\hat\rho_i^{\ell'}$ and $R^-_{\ell,\ell'}:=\frac1n\sum_{i\le n-1}\rho_i^\ell\rho_i^{\ell'}$ with
$\hat\rho_i^\ell=(\rho_i^\ell-\langle\rho_i\rangle)/\sqrt{1-\langle\rho_i\rangle^2}$.

We also use the site-decomposition
\[
f=\frac1n\sum_{i=1}^n a_i,\qquad a_i:=(\sigma_i^1-\sigma_i^2)(\sigma_i^3-\sigma_i^4),
\qquad
\hat f=\frac1n\sum_{i=1}^n \hat a_i,\qquad
\hat a_i:=(\hat\sigma_i^1-\hat\sigma_i^2)(\hat\sigma_i^3-\hat\sigma_i^4).
\]

\paragraph{O-notation}
As in Talagrand~\cite[\S 1.6 \& \S 1.8]{talagrand2010mean} we write $O(k)$ for a quantity bounded by $K n^{-k/2}$,
where $K$ may depend on $(\beta,t)$ but not on $n$.

\begin{lemma}[Exact orthogonality for $\hat R$]\label{lem:orthog-hatR}
Fix disorder $(A,y_t)$. For i.i.d.\ replicas $\sigma^1,\sigma^2,\sigma^3,\sigma^4$,
\[
\langle \hat R_{1,2}\hat R_{1,3}\rangle=0,
\qquad
\langle \hat R_{1,2}\hat R_{3,4}\rangle=0.
\]
Consequently,
\begin{equation}\label{eq:fhat-square-4}
\langle \hat f^2\rangle = 4\langle \hat R_{1,2}^2\rangle.
\end{equation}
\end{lemma}

\begin{proof}
Condition on $\sigma^1$ and use that replicas are independent given disorder and
$\langle \hat\sigma_i\rangle=0$ by definition; thus
$\E[\hat R_{1,2}\mid \sigma^1]=0$, giving $\langle \hat R_{1,2}\hat R_{1,3}\rangle=0$.
The disjoint-pair identity follows from $\langle \hat R_{1,2}\rangle=0$.
Expanding $\hat f^2$ and using these vanishing cross-covariances leaves
four identical square terms, giving \eqref{eq:fhat-square-4}.
\end{proof}

\paragraph{Modified cavity interpolation for the planted SK model} Write $\sigma=(\rho,\varepsilon)$ with $\rho\in\{-1,1\}^{N-1}$ and $\varepsilon=\sigma_N$.
Decompose $A$ into blocks
\[
A=\begin{pmatrix} A^{--}&a\\ a^{\sT}&A_{nn}\end{pmatrix},
\qquad a_i:=A_{in}.
\]
Then (up to additive terms of scale $O(1/n)$ independent of $\varepsilon$)
\begin{equation}\label{eq:cavity-decomp}
\beta H_n(\rho,\varepsilon)
=\beta H_{n-1}^{(n)}(\rho)
+\varepsilon\left[\beta \sum_{i\le n-1} a_i\rho_i + \frac{\beta^2}{n}\sum_{i\le n-1}\rho_i + y_{t,n}\right],
\end{equation}
where $\beta H_{N-1}^{(N)}(\rho)$ contains the $(N-1)$-spin SK Hamiltonian,
the CW term on $\rho$, and the external field on $\rho$.

The modified cavity interpolation interpolates the SK field
$\sum a_i\rho_i$ with an independent Gaussian field matched at $q^*$, and linearly couples the magnetization term $\sum_i \rho^i$ with $m^*$. \\
Let $g_i\sim\calN(0,1)$ i.i.d.\ and $Z\sim\calN(0,1)$ independent.  For $s\in[0,1]$, denote the interpolated Hamiltonian as
\begin{equation}\label{eq:smartpath-corrected}
\beta H_{n,s}(\rho,\varepsilon)
:=\beta H_{n-1}^{(n)}(\rho)
+\varepsilon\left[
\sqrt{s}\,\frac{\beta}{\sqrt n}\sum_{i\le n-1} g_i\rho_i
+\sqrt{1-s}\,\beta\sqrt{q^*}\,Z
+\beta^2\left(\frac{s}{n}\sum_{i\le n-1}\rho_i + (1-s)m^*\right)
+y_{t,n}
\right].
\end{equation}
Let $\langle\cdot\rangle_s$ denote the Gibbs average and $\nu_s(\cdot)=\E\langle\cdot\rangle_s$.

\begin{lemma}[Exact factorization at $s=0$~{\cite[Lemma 1.6.2 for modified cavity interpolation]{talagrand2010mean}}]\label{lem:fact-s0}
Let $F^-$ be any bounded measurable function of the first $n-1$ spins in $k$ replicas
$\rho^1,\dots,\rho^k$, and let $\Phi$ be any bounded measurable function of the last spins
$(\varepsilon^1,\dots,\varepsilon^k)$. Then under the modified cavity interpolation ,
\[
\nu_0\left(F^-\,\Phi(\varepsilon^1,\dots,\varepsilon^k)\right)
=
\nu_0(F^-)\E\Big[\big\langle \Phi(\tilde\varepsilon^1,\dots,\tilde\varepsilon^k)\big\rangle_{Y^*}\Big],
\]
where conditionally on $Y^*$, $\tilde\varepsilon^1,\dots,\tilde\varepsilon^k$ are i.i.d.\ $\pm1$
with mean $\tanh(Y^*)$.
In particular, for any index set $I\subset\{1,\dots,k\}$,
\[
\nu_0\left(F^-\prod_{\ell\in I}\varepsilon^\ell\right)
=
\nu_0(F^-)\E\big[\tanh^{|I|}(Y^*)\big].
\]
\end{lemma}

\begin{proof}
At $s=0$ the interpolated Hamiltonian becomes
\[
\beta H_{n,0}(\rho,\varepsilon)=\beta H_{n-1}(\rho)+\varepsilon\,Y^*,
\qquad Y^*=y_n+\beta^2 m^*+\beta\sqrt{q^*}\,Z,
\]
and $Y^*$ is independent of $\rho$. Applying the proof of~\cite[Lemma 1.6.2]{talagrand2010mean} directly yields the result.
\end{proof}

There is a similar analogue for~\cite[Lemma 1.6.3]{talagrand2010mean} with the exception that the linear part of the interpolation for the magnetization term pulls out a ``drift'' term. Denote $M^{-}_\ell := \frac{1}{n}\sum_{i < n}\rho^\ell_i$.

\begin{lemma}[Derivative identity with CW drift~{\cite[Lemma 1.6.3 for modified cavity interpolation]{talagrand2010mean}}]\label{lem:deriv}
Let $F$ be bounded and depend on $k$ replicas $(\rho^1,\varepsilon^1),\dots,(\rho^k,\varepsilon^k)$. Then, for $s\in(0,1)$,
\begin{align}
&\nu'_s(F)
= \\
&\underbrace{
\beta^2\sum_{1\le \ell<\ell'\le k}\nu_s\left(F\,\varepsilon^\ell\varepsilon^{\ell'}(R^-_{\ell,\ell'}-q^*)\right)
-\beta^2 k\sum_{\ell=1}^k\nu_s\left(F\,\varepsilon^\ell\varepsilon^{k+1}(R^-_{\ell,k+1}-q^*)\right)
+\frac{\beta^2 k(k+1)}{2}\nu_s\left(F\,\varepsilon^{k+1}\varepsilon^{k+2}(R^-_{k+1,k+2}-q^*)\right)
}_{\textnormal{(SK overlap)}} \nonumber
\\
&\quad
+\underbrace{
\beta^2\left[
\sum_{\ell=1}^k \nu_s\left(F\,\varepsilon^\ell( M_\ell^- - m^*)\right)
-k\,\nu_s\left(F\,\varepsilon^{k+1}( M_{k+1}^- - m^*)\right)
\right]
}_{\textnormal{(CW drift)}}.
\label{eq:deriv}
\end{align}
\end{lemma}

\begin{proof}
Differentiate $\nu_s(F)=\E\langle F\rangle_s$ with respect to $s$.
Using $\mathcal G_s$ for the $s$-dependent part of the Hamiltonian in the cavity interpolation gives
\[
\mathcal G_s(\rho,\varepsilon)
:=\varepsilon\Big[
\sqrt{s}\,\frac{\beta}{\sqrt n}\sum_{i\le n-1} g_i\rho_i
+\sqrt{1-s}\,\beta\sqrt{q^*}\,Z
+\beta^2\big((1-s)m^*+s M^-(\rho)\big)
\Big],
\]
so that $\beta H_{n,s}=\beta H_{n-1}+\varepsilon y_n+\mathcal G_s$.
The standard Gibbs differentiation identity (using the division rule) gives
\[
\frac{d}{ds}\langle F\rangle_s
=
\sum_{\ell=1}^k \Big\langle F\,\partial_s\mathcal G_s(\rho^\ell,\varepsilon^\ell)\Big\rangle_s
-\Big\langle F\Big\rangle_s\sum_{\ell=1}^k\Big\langle \partial_s\mathcal G_s(\rho^\ell,\varepsilon^\ell)\Big\rangle_s.
\]
Now split $\partial_s\mathcal G_s$ into the Gaussian SK part and the CW drift:
\[
\partial_s\mathcal G_s(\rho,\varepsilon)
=
\varepsilon\left[
\frac{\beta}{2\sqrt{s\,n}}\sum_{i\le n-1} g_i\rho_i
-\frac{\beta\sqrt{q^*}}{2\sqrt{1-s}}\,Z
+\beta^2(M^-(\rho)-m^*)
\right].
\]
\ppart{The SK contribution} For the Gaussian part, apply~\cite[Lemma 1.1]{panchenko2013sherrington} exactly as in~\cite[Proof of Lemma 1.6.3]{talagrand2010mean}:
this produces the ``SK overlap part'' contribution in \eqref{eq:deriv}, written
with replicas $k+1, k+2$ using the identity
$\langle F\rangle_s\langle G\rangle_s=\langle F\,G(\sigma^{k+1})\rangle_s$.
\ppart{The CW drift contribution} For the CW drift term $\varepsilon\beta^2(M^-(\rho)-m^*)$, no Gaussian analysis techniques are utilized:
the magnetization terms contribute the covariance of $F$ with $\sum_{\ell\le k}\varepsilon^\ell(M^-_\ell - m^*)$.
Using the same replica-augmentation trick to express $\langle F\rangle_s\langle \varepsilon(M-m^*)\rangle_s = \an{F G(\sigma^{k+1})}_s$
yields the ``CW drift part'' in \eqref{eq:deriv}.
\end{proof}

\paragraph{Stability of bounded functions along the cavity interpolation} Pinning down the precise limiting estimates for various tracial statistics reduced to replica form will require invoking magnetization and overlap concentration along the \emph{entire} interpolation $s \in [0,1]$. Furthermore, various other ``polynomial-in-finite-number-of-replicas'' type statistics will be ``stable'' along the cavity interpolation. Showing this will follow a standard template:
\begin{enumerate}
    \item Prove the required estimate at $s=1$ via an argument that uses site-exchangeability and explicit decoupling arguments, and
    \item Invoke ``stability'' of $\nu_1(F)$ by using the fact $|\nu'_s(F)| \le C(\beta,s) \nu_1(F)$.
\end{enumerate}

We begin by showing that overlaps and magnetization are concentrated along the entire interpolation.

\begin{lemma}[Stability of overlap/magnetization concentration along the cavity interpolation]\label{lem:stab-dev-along-s}
Fix $t\in[0,T]$ and $\beta>0$.  Let
\[
\bar R^-_{ab}:=R^-_{ab}-q^*,\qquad \bar M^-_a:=M^-_a-m^*,
\]
for any replica indices $a,b \in \{1,\dots,K\}$.
Assume, for every $p \in \Z_+$ there exist constants $C^{(R)}_p,C^{(M)}_p<\infty$ such that
\begin{equation}\label{eq:stab-assume-endpoint}
\sup_{a<b\le K}\ \nu_1\big(|\bar R^-_{ab}|^p\big)\ \le\ C^{(R)}_p\,n^{-p/2},
\qquad
\sup_{a\le K}\ \nu_1\big(|\bar M^-_{a}|^p\big)\ \le\ C^{(M)}_p\,n^{-p/2}.
\end{equation}
Then, for every $s\in[0,1]$,
\begin{equation}\label{eq:stab-conclude-uniform}
\sup_{a<b\le K}\ \nu_s\big(|\bar R^-_{ab}|^p\big)\ \le\ e^{24\beta^2}\,C^{(R)}_p\,n^{-p/2},
\qquad
\sup_{a\le K}\ \nu_s\big(|\bar M^-_{a}|^p\big)\ \le\ e^{8\beta^2}\,C^{(M)}_p\,n^{-p/2}.
\end{equation}
\end{lemma}
\begin{proof}
We use~\pref{lem:deriv} and the elementary bounds
\begin{equation}\label{eq:bounded-deviations}
|\bar R^-_{ab}|\le 2,\qquad |\bar M^-_a|\le 2,\qquad |\varepsilon^\ell|\le 1\,.
\end{equation}
We treat overlaps and magnetizations separately.
\ppart{Stability of the overlaps along the interpolation}
Fix $a<b$ and define for $s\in[0,1]$
\[
u_{p}^{(a,b)}(s):=\nu_s\big(|\bar R^-_{ab}|^p\big).
\]
Note that the replica indices are \emph{always} exchangeable, and so it suffices to
bound a single pair $u_p(s):=\nu_s(|\bar R^-_{12}|^p)$. Applying~\pref{lem:deriv} with $k=2$ to the function $F:=|\bar R^-_{12}|^p$ gives
\begin{align*}
u_p'(s)
&=
\beta^2\,\nu_s\!\left(F\,\varepsilon^1\varepsilon^2\,\bar R^-_{12}\right)
-2\beta^2\sum_{\ell=1}^2\nu_s\!\left(F\,\varepsilon^\ell\varepsilon^{3}\,\bar R^-_{\ell,3}\right)
+3\beta^2\,\nu_s\!\left(F\,\varepsilon^{3}\varepsilon^{4}\,\bar R^-_{3,4}\right)\\
&\quad
+\beta^2\left[
\sum_{\ell=1}^2\nu_s\!\left(F\,\varepsilon^\ell\,\bar M^-_{\ell}\right)
-2\nu_s\!\left(F\,\varepsilon^{3}\,\bar M^-_{3}\right)
\right].
\end{align*}
Taking absolute values term-by-term and using \eqref{eq:bounded-deviations} is sufficient to get the Lipschitz bound on the derivative. For the first two overlap terms, note that
\begin{align*}
    \Big|\nu_s(F\,\varepsilon^1\varepsilon^2\,\bar R^-_{12})\Big|
&=\nu_s\big(|\bar R^-_{12}|^{p+1}\big)
\le 2\,\nu_s\big(|\bar R^-_{12}|^{p}\big)=2u_p(s)\,, \\
\Big|\nu_s(F\,\varepsilon^\ell\varepsilon^{3}\,\bar R^-_{\ell,3})\Big|
&\le \|\bar R^-_{\ell,3}\|_{\infty}\,\nu_s(F)\le 2u_p(s)\, ,
\end{align*}
when $\ell\in\{1,2\}$. Since there are two such terms and each has coefficient $2\beta^2$, these contribute at most $8\beta^2u_p(s)$.
The final overlap term satisfies the bound,
\[
    \Big|\nu_s(F\,\varepsilon^{3}\varepsilon^{4}\,\bar R^-_{3,4})\Big| \le 2u_p(s)\,,
\]
so the contribution is at most $6\beta^2u_p(s)$. Similarly, the magnetization terms can be bounded as
\[
\Big|\nu_s(F\,\varepsilon^\ell\,\bar M^-_{\ell})\Big|\le 2u_p(s),
\qquad
\Big|\nu_s(F\,\varepsilon^{3}\,\bar M^-_{3})\Big|\le 2u_p(s),
\]
so the magnetization bracket contributes at most
$2\beta^2\big(2u_p(s)+ 2u_p(s)\big)=8\beta^2u_p(s)$.

Putting all the bounds together yields that,
\[
|u_p'(s)|
\le (2+8+6+8)\beta^2\,u_p(s)=24\beta^2\,u_p(s).
\]
By Gr\"onwall's inequality,
\[
    u_p(s)\le e^{24\beta^2(1-s)}u_p(1)\le e^{24\beta^2}u_p(1) \le_{\text{\eqref{eq:stab-assume-endpoint}}} C^{(R)}(\beta,p) O(n^{-p/2})\,.
\]

\ppart{Stability of the magnetization along the interpolation} Fix $a$ and set $v_p(s):=\nu_s(|\bar M^-_{1}|^p)$ (using replica symmetry again). Applying~\pref{lem:deriv} (again) with $k=1$ and $F:=|\bar M^-_1|a ^p$ gives
\[
v_p'(s)
=
-\beta^2\,\nu_s(F\,\varepsilon^1\varepsilon^2\,\bar R^-_{1,2})
+\beta^2\,\nu_s(F\,\varepsilon^2\varepsilon^3\,\bar R^-_{2,3})
+\beta^2\Big(\nu_s(F\,\varepsilon^1\,\bar M^-_1)-\nu_s(F\,\varepsilon^2\,\bar M^-_2)\Big).
\]
Using \eqref{eq:bounded-deviations} gives the following bounds for the overlap terms,
\[
\big|\nu_s(F\,\varepsilon^1\varepsilon^2\,\bar R^-_{1,2})\big|\le 2v_p(s),\qquad
\big|\nu_s(F\,\varepsilon^2\varepsilon^3\,\bar R^-_{2,3})\big|\le 2v_p(s),
\]
and the following bounds for the magnetization terms
\[
\big|\nu_s(F\,\varepsilon^1\,\bar M^-_1)\big|
=\nu_s(|\bar M^-_1|^{p+1})
\le 2v_p(s),
\qquad
\big|\nu_s(F\,\varepsilon^2\,\bar M^-_2)\big|\le 2v_p(s).
\]
Therefore,
\[
|v_p'(s)|\le (2+2+2+2)\beta^2 v_p(s)=8\beta^2 v_p(s),
\]
and Gr\"onwall yields $v_p(s)\le e^{8\beta^2}v_p(1) \le_{\text{\eqref{eq:stab-assume-endpoint}}} C^{(M)}(\beta,p)O(n^{-p/2})$. \qedhere
\end{proof}

We now show that, akin to the statement of~\pref{prop:D15} showing that odd cubic moments of $f$ and $\hat{f}$ are $O(4)$ which defeats naive concentration by $O(1)$, a similar thing happens for the magnetization and overlap fluctuations at $s=1$. This is critical in the evaluation of the ``self-similar'' terms in the proof of~\pref{prop:D17}.

\begin{lemma}[Overlap and magnetization biases are $O(n^{-1})$]\label{lem:mean-bias-On-1}
Fix $t\in[0,T]$ and $0 < \beta < 1/2$. Let
\[
\delta q_n:=\nu(R_{12})-q^*,
\qquad
\delta m_n:=\nu(M_1)-m^*.
\]
Assume the hypotheses \textbf{(H1)--(H3)}, and let $Y^*$ be the cavity field and set
\[
\mu_r:=\E[\tanh^r(Y^*)],\qquad r\in\{1,2,3,4\}\,,
\]
and
\[
\mathbf M:=
\begin{pmatrix}
1-4\mu_2+3\mu_4 & 2(\mu_1-\mu_3)\\
\mu_3-\mu_1 & 1-\mu_2
\end{pmatrix}.
\]
Then, the following hold:
\begin{enumerate}[itemsep=0.3em]
    \item The deviations satisfy the stability equation 
        \begin{equation}\label{eq:mean-bias-stability}
            \begin{pmatrix}\delta q_n\\ \delta m_n\end{pmatrix} =
                \beta^2\,\mathbf M
            \begin{pmatrix}\delta q_n\\ \delta m_n\end{pmatrix}
                \;+\;O(n^{-1}),
         \end{equation} 
        where the $O(n^{-1})$ is uniform in $n$. 
    \item Furthermore, if  $\rho(\beta^2\mathbf M)<1$, then
        \[
            \delta q_n=O(n^{-1}),
            \qquad
            \delta m_n=O(n^{-1}).
        \]
\end{enumerate}
\end{lemma}

\begin{proof} We will use exchangebility and last-site decoupling to set up a ``self-consistent'' linear system described by a specific operator $M$. Provided that $\beta^2\rho(M) < 1$, it will be the case that the fluctuations of the mean and overlap are $O(n^{-1/2})$ stronger than predicted by naive concentration. Invoking~\pref{lem:stability-s} (via Gr\"onwall) will then immediately imply this is true at any $s \in [0,1]$.

\ppart{Reduce $\nu(R_{12})$ and $\nu(M_1)$ to last-spin observables via exchangeability}
By site-exchangeability of the disorder and the Gibbs measure,
\[
\nu(R_{12})=\nu(\varepsilon^1\varepsilon^2),
\qquad
\nu(M_1)=\nu(\varepsilon^1).
\]
This yields that
\[
\delta q_n=\nu(\varepsilon^1\varepsilon^2)-q^*,
\qquad
\delta m_n=\nu(\varepsilon^1)-m^*.
\]
\ppart{Two-term Taylor expansions around $s=0$}
Applying the exact two-term Taylor identity for $F \in \{\eps^1,\eps^1\eps^2\}$ yields that
\[
\nu(F)=\nu_0(F)+\nu'_0(F)+\int_0^1 (1-s)\,\nu_s''(F)\,ds.
\]
\ppart{Second-order remainder bounds of $O(n^{-1})$ }
Fix $F\in\{\varepsilon^1\varepsilon^2,\varepsilon^1\}$ and note $|F|\le 1$. Differentiating twice using~\pref{lem:deriv} yields a finite linear combination of terms where every summand in $\nu_s''(F)$ is of the form
\[
\nu_s\left(F\;\mathcal P(\varepsilon)\;\Xi_1\Xi_2\right),
\qquad  
\Xi_1,\Xi_2\in\{\,R^-_{ab}-q^*,\,M^-_a-m^*\,\},
\]
where $\mathcal P(\varepsilon)$ is a bounded polynomial in finitely many last spins (hence $|\mathcal P|\le C$).
Thus, by Cauchy--Schwarz and boundedness of $F\mathcal P$,
\[
\Big|\nu_s\!\big(F\mathcal P(\varepsilon)\Xi_1\Xi_2\big)\Big|
\le C\;\|\Xi_1\|_{L^2(\nu_s)}\;\|\Xi_2\|_{L^2(\nu_s)}.
\]
By \textbf{(H1)--(H2)} and~\pref{lem:stab-dev-along-s}, we have $\|\Xi_i\|_{L^2(\nu_s)}=O(n^{-1/2})$.
Therefore, we obtain that
\[
\int_0^1(1-s)\,\nu_s''(F)\,ds=O(n^{-1})\, ,
\]
for $F=\varepsilon^1\varepsilon^2$ and $F=\varepsilon^1$.

\ppart{Evaluating $\nu_0(F)$ for $\varepsilon^1\varepsilon^2$ and $\varepsilon^1$} At $s=0$,~\pref{lem:fact-s0} implies that the last-spin replicas are conditionally i.i.d.\ given the cavity field $Y^*$. This immediately yields,
\begin{align*}
   \nu_0(\varepsilon^1)=\E[\tanh(Y^*)]=\mu_1=m^*\,,   
\end{align*}
and
\begin{align*}
    \nu_0(\varepsilon^1\varepsilon^2)=\E[\tanh^2(Y^*)]=\mu_2=q^*\,.
\end{align*}

\ppart{Evaluating $\nu^0_0(\varepsilon^1\varepsilon^2)$ via exact symmetries}
Applying~\pref{lem:deriv} with $k=2$ and $F=\varepsilon^1\varepsilon^2$ gives the following linear combination of terms, evaluated at $s=0$ by the factorization implied by~\pref{lem:fact-s0}.
With $k=2$, the ``SK overlap'' part is given by
\begin{align*}
\text{(SK)}\ &=\
\beta^2\,\nu_s\left(F\,\varepsilon^1\varepsilon^2\,(R^-_{12}-q^*)\right)
-2\beta^2\sum_{\ell=1}^2 \nu_s\left(F\,\varepsilon^\ell\varepsilon^{3}\,(R^-_{\ell,3}-q^*)\right)
+3\beta^2\,\nu_s\left(F\,\varepsilon^{3}\varepsilon^{4}\,(R^-_{34}-q^*)\right).
\end{align*}
With $F=\varepsilon^1\varepsilon^2$, the last-spin products simplify:
$F\varepsilon^1\varepsilon^2=1$,
$F\varepsilon^1\varepsilon^3=\varepsilon^2\varepsilon^3$,
$F\varepsilon^2\varepsilon^3=\varepsilon^1\varepsilon^3$,
$F\varepsilon^3\varepsilon^4=\varepsilon^1\varepsilon^2\varepsilon^3\varepsilon^4$.
Using these explicit simplifications for the last spin terms along with~\pref{lem:fact-s0} at $s=0$, gives
\[
\nu_0(\varepsilon^{i_1}\cdots \varepsilon^{i_r}X)=\E[\tanh^r(Y^*)]\;\nu_0(X)=\mu_r\,\nu_0(X).
\]
Hence, using replica symmetry for the bulk overlaps under $\nu_0$,
\[
\nu_0(R^-_{12}-q^*)=\nu_0(R^-_{13}-q^*)=\nu_0(R^-_{34}-q^*) =:\delta q_n^-,
\]
we obtain
\[
\text{(SK)}\big|_{s=0}
=
\beta^2\Big(1-4\mu_2+3\mu_4\Big)\,\delta q_n^-.
\]
A similar calculation for the ``CW drift'' part gives that, with $k=2$,
\[
    \text{(CW)}\ =\
    \beta^2\Big(
    \nu_s(F\varepsilon^1(M^-_1-m^*))
    +\nu_s(F\varepsilon^2(M^-_2-m^*))
    -2\nu_s(F\varepsilon^3(M^-_3-m^*))
    \Big).
\]
With $F=\varepsilon^1\varepsilon^2$, we have
$F\varepsilon^1=\varepsilon^2$, $F\varepsilon^2=\varepsilon^1$, $F\varepsilon^3=\varepsilon^1\varepsilon^2\varepsilon^3$.
At $s=0$, factorization yields
\[
    \nu_0(\varepsilon^2(M^-_1-m^*))=\mu_1\,\nu_0(M^-_1-m^*),
\quad
    \nu_0(\varepsilon^1(M^-_2-m^*))=\mu_1\,\nu_0(M^-_2-m^*),
\]
\[
\nu_0(\varepsilon^1\varepsilon^2\varepsilon^3(M^-_3-m^*))=\mu_3\,\nu_0(M^-_3-m^*).
\]
By replica symmetry of the bulk magnetizations under $\nu_0$, one obtains
\[
    \delta m_n^-:=\nu_0(M^-_1-m^*)=\nu_0(M^-_2-m^*)=\nu_0(M^-_3-m^*)\,,
\]
which gives the simplifications at $s=0$ as
\[
\text{(CW)}\big|_{s=0}
=
2\beta^2(\mu_1-\mu_3)\,\delta m_n^-.
\]
Combining these two gives the final bound
\begin{equation}\label{eq:d0-eps1eps2}
\nu^0_0(\varepsilon^1\varepsilon^2)
=
\beta^2\Big(1-4\mu_2+3\mu_4\Big)\,\delta q_n^-
\;+\;
2\beta^2(\mu_1-\mu_3)\,\delta m_n^-.
\end{equation}

\ppart{Evaluating $\nu^0_0(\varepsilon^1)$ via exact symmetries} Similar to the overlap case, applying~\pref{lem:deriv} with $k=1$ and $F=\varepsilon^1$ yields
\begin{align*}
\nu'_s(\varepsilon^1)
&=
-\beta^2\,\nu_s\left(\varepsilon^1\varepsilon^1\varepsilon^2(R^-_{12}-q^*)\right)
+\beta^2\,\nu_s\left(\varepsilon^1\varepsilon^2\varepsilon^3(R^-_{23}-q^*)\right)\\
&\quad +\beta^2\Big(\nu_s(\varepsilon^1\varepsilon^1(M^-_1-m^*))-\nu_s(\varepsilon^1\varepsilon^2(M^-_2-m^*))\Big).
\end{align*}
Simplify $\varepsilon^1\varepsilon^1\varepsilon^2=\varepsilon^2$ and $\varepsilon^1\varepsilon^1=1$.
At $s=0$, factorization via~\pref{lem:fact-s0} yields
\[
\nu_0(\varepsilon^2(R^-_{12}-q^*))=\mu_1\,\delta q_n^-,
\qquad
\nu_0(\varepsilon^1\varepsilon^2\varepsilon^3(R^-_{23}-q^*))=\mu_3\,\delta q_n^-,
\]
\[
\nu_0(M^-_1-m^*)=\delta m_n^-,
\qquad
\nu_0(\varepsilon^1\varepsilon^2(M^-_2-m^*))=\mu_2\,\delta m_n^-.
\]
Thus
\begin{equation}\label{eq:d0-eps1}
\nu'_0(\varepsilon^1)
=
\beta^2(\mu_3-\mu_1)\,\delta q_n^-
\;+\;
\beta^2(1-\mu_2)\,\delta m_n^-.
\end{equation}

\ppart{Proving $\delta q_n^-,\delta m_n^-$ equivalent to $\delta q_n,\delta m_n$ up to $O(n^{-1})$} The goal is to show that
\begin{equation}\label{eq:bulk-to-full-bias}
\delta q_n^-=\delta q_n+O(n^{-1}),
\qquad
\delta m_n^-=\delta m_n+O(n^{-1}).
\end{equation}
First, note that
\[
    R_{12}=R^-_{12}+\frac1n\,\varepsilon^1\varepsilon^2,
    \qquad
    M_1=M^-_1+\frac1n\,\varepsilon^1.
\]
Taking $\nu$ and using $|\varepsilon^1\varepsilon^2|\le 1$ gives
\[
\nu(R^-_{12})=\nu(R_{12})+O(n^{-1}) = \nu(\eps^1\eps^2) + O(n^{-1}),\qquad \nu(M^-_1)=\nu(M_1) + O(n^{-1})= \nu(\eps^1)+O(n^{-1}).
\]
It remains to compare $\nu_0(\cdot)$ with $\nu(\cdot)$ on centered bulk observables by showing stability along the cavity interpolation. Let $H \in \{R^-_{12}-q^*, M^-_1-m^*\}$ and write
\[
    \nu(H)-\nu_0(H)=\int_0^1 \nu'_s(H)\,ds.
\]
Applying~\pref{lem:deriv} to $H$ shows that every term in $\nu'_s(H)$ contains at least one additional
centered bulk deviation factor $\Xi \in \{R^-_{ab}-q^*, M^-_a - m^*\}$. An application of $(2,2)$-H\"older's then yields
\[
    \nu_s(H\Xi) \le \sqrt{\nu_s(H^2)}\sqrt{\nu_s(\Xi^2)} \le_{\text{\pref{lem:stab-dev-along-s}}} O(n^{-1})\,.
\]
This implies that
\[
    \
    \delta q_n = \nu(R_{12}) - q^* = \nu(R^-_{12} - q^*)  + O(n^{-1}) = \nu_0(R^-_{12} - q^*)  + O(n^{-1}) = \delta q^-_n + O(n^{-1})\,, 
\]
and similarly that
\[
    \delta m_n = \nu(M_1) - m^* = \delta m^-_n + O(n^{-1})\,.
\]
yielding~\eqref{eq:bulk-to-full-bias}.

\ppart{Finalizing the self-consistent linear-system}
Combining the second-order Taylor expansion with the time $s=0$ estimates yields
\[
\nu(R_{12}) = \nu(\varepsilon^1\varepsilon^2)
=
\nu_0(\varepsilon^1\varepsilon^2)+\nu^0_0(\varepsilon^1\varepsilon^2)+O(n^{-1})
=
\mu_2+\nu^0_0(\varepsilon^1\varepsilon^2)+O(n^{-1}),
\]
\[
\nu(M_1) = \nu(\varepsilon^1)
=
\nu_0(\varepsilon^1)+\nu^0_0(\varepsilon^1)+O(n^{-1})
=
\mu_1+\nu^0_0(\varepsilon^1)+O(n^{-1}).
\]
Insert \eqref{eq:d0-eps1eps2} and \eqref{eq:d0-eps1} into the above, and use \eqref{eq:bulk-to-full-bias} to replace $\delta q_n^-,\delta m_n^-$
by $\delta q_n,\delta m_n$ up to $O(n^{-1})$ factors. Since $\mu_2=q^*$ and $\mu_1=m^*$, we obtain the following self-consistent equations
\[
\begin{pmatrix}\delta q_n\\ \delta m_n\end{pmatrix}
=
\beta^2
\begin{pmatrix}
1-4\mu_2+3\mu_4 & 2(\mu_1-\mu_3)\\
\mu_3-\mu_1 & 1-\mu_2
\end{pmatrix}
\begin{pmatrix}\delta q_n\\ \delta m_n\end{pmatrix}
+O(n^{-1})\,,
\]
which prove \eqref{eq:mean-bias-stability}. Furthermore, assuming $\rho(\beta^2 \mathbf M) < 1$, multiply by $(I-\beta^2\mathbf M)^{-1}$ to conclude
\[
\begin{pmatrix}\delta q_n\\ \delta m_n\end{pmatrix}
=O(n^{-1})\,.\qedhere
\]
\end{proof}

\begin{remark}
    Note that the linear deviation bound for the magnetization and overlaps above relies on the fact that $\rho(\beta^2 \mathbf{M}) < 1$ to allow for the invertibility of the underlying linear system. For $\beta<1/2$, the bound is proved in~\pref{lem:beta2M-spr-bound}. Therefore, for the remainder of the section (and paper), the bound in~\pref{lem:mean-bias-On-1} is used without further qualification.  
\end{remark}

Below we prove a stability result for bounded functions $F$ of finite numbers of replicas along the cavity interpolation. The results below essentially use Gr\"onwall's lemma for propagating stability across the interpolation, invoking~\pref{lem:stab-dev-along-s} to use uniform concentration for the bulk deviation parameters.

\begin{lemma}[Stability of bounded continuous functions along the cavity interpolation]\label{lem:stability-s}
Fix the number of replicas $k\in\mathbb N$, and let $F: \Sigma^k \to \R$ be a bounded function of
$(\rho^1,\varepsilon^1),\dots,(\rho^k,\varepsilon^k)$. Assume that
\begin{equation}\label{eq:stability-conc-assump}
    \sup_{s\in[0,1]}\ \max_{a<b\le k+2}\ \big\|R^-_{ab}-q^*\big\|_{L^p(\nu_s)}
    \;+\;
    \sup_{s\in[0,1]}\ \max_{a\le k+1}\ \big\|M^-_{a}-m^*\big\|_{L^p(\nu_s)}
    \;\le\; C_{p,k}\,n^{-1/2}.
\end{equation}
Then, the following two stability statements hold for $F$ along the cavity interpolation,
\begin{enumerate}
    \item \textbf{Small perturbation for bulk functions}: For bulk functions the following holds
    \begin{equation}\label{eq:stability-transfer} 
        \big|\nu_s(F)-\nu_1(F)\big|
        \;\le\;
        C_{\beta,k}\,n^{-1/2}\int_s^1 \|F\|_{L^2(\nu_u)}\,du
        \;\le\;
        C_{\beta,k}\,n^{-1/2}\,\sup_{u\in[s,1]}\|F\|_{L^2(\nu_u)}.
        \end{equation}
        In particular, if $\sup_{u\in[0,1]}\|F\|_{L^2(\nu_u)} = O(n^{-\alpha})$, then
        \[
        \sup_{s\in[0,1]}\big|\nu_s(F)-\nu_1(F)\big| = O\big(n^{-(\alpha+1/2)}\big).
        \]   
    \item \textbf{Gr\"onwall stability}: Via Gr\"onwall's inequality,
    \begin{equation}\label{eq:stability-diff}
        \big|\nu_t(F)-\nu_1(F)\big|
        \;\le\;
        \Big(e^{4\beta^2 r(r+1)(1-t)}-1\Big)\,\nu_1(F).
    \end{equation}
\end{enumerate}
\end{lemma}

\begin{proof}
We prove both statements in chronological order.
\ppart{Small perturbation for bulk functions}
Applying~\pref{lem:deriv} to $F$ and noting that $|\varepsilon^\ell|=1$, implies that each term in $\nu_s'(F)$ can be bounded (up to absolute constants) by
$\varepsilon$ factors in absolute value. This gives, via a $(2,2)$-H\"older's inequality,
\[
    \big|\nu_s\big(F\,\varepsilon^\ell\varepsilon^{\ell'}(R^-_{\ell,\ell'}-q^*)\big)\big|
    \le \nu_s\big(|F|\,|R^-_{\ell,\ell'}-q^*|\big)
    \le \|F\|_{L^2(\nu_s)}\,\|R^-_{\ell,\ell'}-q^*\|_{L^2(\nu_s)}\,.
\]
The same bound holds for the terms involving $M^-_a-m^*$. Since the ultimate bound on the derivative is a finite linear combination of these bulk deviation terms,
\[
|\nu_s'(F)|
\le
\beta^2\,C_k\,\|F\|_{L^2(\nu_s)}
\Big(
\max_{a<b\le k+2}\|R^-_{ab}-q^*\|_{L^2(\nu_s)}
+\max_{a\le k+1}\|M^-_a-m^*\|_{L^2(\nu_s)}
\Big) \le C(\beta,k)n^{-1/2}\,.
\] 
\eqref{eq:stability-transfer} then immediately follows from FTOC, yielding that $\nu_s(F)-\nu_1(F)=\int_s^1 \nu_u'(F)\,du$.
\ppart{Gronwall stability}  The proof here is the same as above, except the use of $(2,2)$-H\"older's inequality, which is substituted by the fact that $|\Xi\eps^\ell\eps^{\ell'}| \le C$. Therefore, summing over the terms after the application of~\pref{lem:deriv} gives that
\[
    |\nu'_s(F)| \le C(\beta,k)\nu_s(F)\,.
\]
An application of Gr\"onwall's inequality immediately gives~\pref{eq:stability-diff}.
\end{proof}

A corollary of this is that the deviation of the mean and the overlap is $O(n^{-1})$ for any $s \in [0,1]$.
\begin{corollary}[Mean and overlap deviation along the cavity interpolation]
    The following deviation scales hold at any $s \in [0,1]$,
    \[
        \nu_s(R_{12}-q^*) = O(n^{-1})\, ,
    \]
    and
    \[
        \nu_s(M_1 -m^*) = O(n^{-1})\,.
    \]
\end{corollary}
\begin{prf}
    The proof follows by combining~\pref{lem:mean-bias-On-1} with~\pref{lem:stability-s}, specifically \eqref{eq:stability-diff}.
\end{prf}

\paragraph{Analogues to the ``rectangular sum'' lemma in~\cite[\S 1.8]{talagrand2010mean}} At $s=0$, there is an exact factorization
(\pref{lem:fact-s0}) with the RS field
\[
Y^*:=y_{t,n}+\beta^2m^*+\beta\sqrt{q^*}\,Z,
\]
and the following constants show up in the scaling
\begin{equation}\label{eq:rho2rho4}
    \rho_2:=\E[\sech^2(Y^*)],\qquad \rho_4:=\E[\sech^4(Y^*)].    
\end{equation}
Finally, let $\nu_1=\nu$ (the true model), and define the following ``centered'' second moments:
\[
    U_n:=\nu(f^2),\qquad V_n:=\nu(f\hat f),\qquad W_n:=\nu(\hat f^2)\,.
\]
A trivial consequence of \textbf{(H1)} and \textbf{(H2)} for the planted SK model under the SL external field is:
\begin{equation}\label{eq:conc-assumptions}
\nu\big(|R_{1,2}-q^*|^{2p}\big)\le C_{p}\,n^{-p},
\qquad
\nu\big(|M-m^*|^{2p}\big)\le C_{p}\,n^{-p},
\qquad
M:=\frac1n\sum_{i=1}^n\sigma_i.
\end{equation}


Define the bulk rectangular sum
\[
f^- := R^-_{1,3}-R^-_{1,4}-R^-_{2,3}+R^-_{2,4}.
\]

Lastly, the derivative of the (expected) Gibbs average at time $s = 0$ is given as
\[
    \nu^0_0(f) := \frac{d}{ds}\nu_s(f)\vert_{s=0} = \lim_{s \to 0^+}\frac{\nu_s(f) - \nu_0(f)}{s}\,.
\]

\begin{lemma}[Unscaled rectangular cancellation~{\cite[Extension of Lemma 1.8.4]{talagrand2010mean}}]\label{lem:rect-unscaled}
Let $F^-$ be any bounded function of the bulk spins in replicas $1,2,3,4$ and independent of
$(\varepsilon^1,\dots,\varepsilon^4)$. Then
\[
\nu_0^0\Big((\varepsilon^1-\varepsilon^2)(\varepsilon^3-\varepsilon^4)\,F^-\Big)
=
\beta^2\rho_4\;\nu_0(F^-\,f^-).
\]
\end{lemma}

\begin{proof}
Apply~\pref{lem:deriv} at $s=0$ with $k=4$ and
$F=(\varepsilon^1-\varepsilon^2)(\varepsilon^3-\varepsilon^4)F^-$.

\ppart{SK part}
Exactly as in the proof of~\cite[Lemma 1.8.4]{talagrand2010mean}, the only surviving replica pairs are $(1,3),(1,4),(2,3),(2,4)$,
and the last-spin factor produces the rectangular combination $f^-$.
The coefficient is
\[
\nu_0\Big((\varepsilon^1-\varepsilon^2)(\varepsilon^3-\varepsilon^4)\varepsilon^1\varepsilon^3\Big)
=\E\big[\an{(1-\varepsilon^1\varepsilon^2)(1-\varepsilon^3\varepsilon^4)}_{Y^*}\big]
=\E[\sech^4(Y^*)]=\rho_4,
\]
where we use~\pref{lem:fact-s0} to factorize at $s=0$ and conditional independence given $Y^*$.

\ppart{CW drift part}
Each drift term in \eqref{eq:deriv} contains a factor
$\nu_0\big((\varepsilon^1-\varepsilon^2)(\varepsilon^3-\varepsilon^4)\varepsilon^\ell\cdot(\text{bulk})\big)$.
By~\pref{lem:fact-s0} this factors into a bulk expectation times
$\E\left[\an{(\varepsilon^1-\varepsilon^2)(\varepsilon^3-\varepsilon^4)\varepsilon^\ell}_{Y*}\right]$. This term is $0$
for every $\ell$, since for \emph{any} choice of $\ell \in \{1,2,3,4,5\}$, by the independence of the replicas under the Gibbs (conditioned on $Y^*$) at least one of the centered terms factorizes.
Therefore, both terms in (CW drift) from~\pref{lem:deriv} contribute identically $0$ at $s=0$. \qedhere
\end{proof}

\begin{lemma}[Scaled rectangular cancellation~{\cite[Recentered and rescaled version of Lemma 1.8.4]{talagrand2010mean}}]\label{lem:rect-scaled}
Let $F^-$ be as in~\pref{lem:rect-unscaled}. Define the recentered and rescaled last-spin variables at $s=0$ as%
\[
\hat\varepsilon^\ell:=\frac{\varepsilon^\ell-\tanh(Y^*)}{\sech(Y^*)}\,. 
\]
Then
\[
\nu_0^0\Big((\hat\varepsilon^1-\hat\varepsilon^2)(\hat\varepsilon^3-\hat\varepsilon^4)\,F^-\Big)
=
\beta^2\rho_2\;\nu_0(F^-\,f^-).
\]
\end{lemma}

\begin{proof}
At $s=0$, $\hat\varepsilon^1-\hat\varepsilon^2 = (\varepsilon^1-\varepsilon^2)/\sech(Y^*)$,
hence
\[
(\hat\varepsilon^1-\hat\varepsilon^2)(\hat\varepsilon^3-\hat\varepsilon^4)
=\sech^{-2}(Y^*)\,(\varepsilon^1-\varepsilon^2)(\varepsilon^3-\varepsilon^4).
\]
Apply~\pref{lem:deriv} at $s=0$ with
$F=\sech^{-2}(Y^*)\,(\varepsilon^1-\varepsilon^2)(\varepsilon^3-\varepsilon^4)F^-$.
The same cancellation as before yields $f^-$ on the bulk. The last-spin coefficient becomes
\[
\E\big[\sech^{-2}(Y^*)\an{(1-\varepsilon^1\varepsilon^2)(1-\varepsilon^3\varepsilon^4)}_{Y^*}\big]
=\E[\sech^2(Y^*)]=\rho_2.
\]
The CW drift part vanishes for the same reason as in the proof of~\pref{lem:rect-unscaled}.
\end{proof}

We will also need the following ``twice-weighted'' scaled cancellation lemma to eventually simplify the $\E[\Tr[PD^2P]]$ term.

\begin{lemma}[Doubly--scaled rectangular cancellation]\label{lem:rect-doublyscaled}
Let $F^-$ be as in~\pref{lem:rect-unscaled}.  Define
\[
\kappa(Y^*):=\sech^{-4}(Y^*),
\qquad
a_n:=(\varepsilon^1-\varepsilon^2)(\varepsilon^3-\varepsilon^4).
\]
Then
\[
\nu_0^0\Big(\kappa(Y^*)\,a_n\,F^-\Big)
=
\beta^2\;\nu_0(F^-\,f^-).
\]
\end{lemma}

\begin{proof}
Apply~\pref{lem:deriv} at $s=0$ with $k=4$ and
\[
F=\kappa(Y^*)\,(\varepsilon^1-\varepsilon^2)(\varepsilon^3-\varepsilon^4)\,F^-.
\]

\ppart{SK overlap part}
Exactly as in~\pref{lem:rect-unscaled}, the only surviving replica pairs
in the SK part are $(1,3),(1,4),(2,3),(2,4)$, and the resulting bulk overlap combination is $f^-$.
It remains to compute the last--spin coefficient; taking the pair $(1,3)$ as representative, this coefficient equals
\begin{align*}
\nu_0\Big(\kappa(Y^*)\,(\varepsilon^1-\varepsilon^2)(\varepsilon^3-\varepsilon^4)\,\varepsilon^1\varepsilon^3\Big)
&=
\E\Big[\kappa(Y^*)\;\big\langle (1-\varepsilon^1\varepsilon^2)(1-\varepsilon^3\varepsilon^4)\big\rangle_{Y^*}\Big].
\end{align*}
By~\pref{lem:fact-s0}, conditionally on $Y^*$ the last spins are i.i.d.\ $\pm1$ with mean $\tanh(Y^*)$,
so $\langle \varepsilon^1\varepsilon^2\rangle_{Y^*}=\tanh^2(Y^*)$ and hence
\[
\big\langle 1-\varepsilon^1\varepsilon^2\big\rangle_{Y^*}=1-\tanh^2(Y^*)=\sech^2(Y^*),
\]
and similarly for $(3,4)$. Therefore
\[
\big\langle (1-\varepsilon^1\varepsilon^2)(1-\varepsilon^3\varepsilon^4)\big\rangle_{Y^*}
=
\sech^4(Y^*).
\]
Multiplying by $\kappa(Y^*)=\sech^{-4}(Y^*)$ gives coefficient $1$, which gives the final contribution
\[
    \beta^2 \nu_0(F^-\,f^-).
\]

\ppart{CW drift part}
Each drift term in \eqref{eq:deriv} contains a factor of the form
\[
\nu_0\Big(\kappa(Y^*)\,(\varepsilon^1-\varepsilon^2)(\varepsilon^3-\varepsilon^4)\,\varepsilon^\ell\cdot(\text{bulk})\Big).
\]
By~\pref{lem:fact-s0} this factors into a bulk expectation times
\[
\E\Big[\big\langle \kappa(Y^*)\,(\varepsilon^1-\varepsilon^2)(\varepsilon^3-\varepsilon^4)\,\varepsilon^\ell\big\rangle_{Y^*}\Big] = 0,
\]
for every $\ell$ because, conditional on $Y^*$, the replicas are independent and
at least one centered difference factor has mean $0$.
\end{proof}

\paragraph{Invariant symmetries and concentration of the ``rectangular sum''} At this point, it is critical to observe that $\nu(f(\dots)G)$ is left invariant (or only flips signs) under symmetries that permute the indices of replicas provided that the function $G$ is invariant itself.

\begin{lemma}[Replica antisymmetry ``kills'' many terms]\label{lem:antisym}
Let $G$ be any bounded function of replicas $1,2,3,4$, invariant under swapping indices $3\leftrightarrow 4$.
Then
\[
\nu(f\,G)=0,\qquad \nu(f^-\,G)=0,\qquad \nu(a_n^3\,G)=0,\qquad \nu(a_n^2\,f\,G)=0,
\]
and the same holds with $f,\ f^-,\ a_n$ replaced by $\hat f,\ \hat f^-,\ \hat a_n$.
\end{lemma}

\begin{proof}
Under $\nu_s$ the joint law is invariant under permutations of replica labels by exchangeability, since $\sigma^1,\sigma^2,\sigma^3,\sigma^4 \sim G^{\ot 4}$. Under the index swap $3\leftrightarrow 4$,
\[
f\mapsto -f,\qquad f^-\mapsto -f^-,
\]
while $G$ is invariant by assumption and $a_n^2$ is invariant since it is jointly even in $\eps^3$ and $\eps^4$. 
Therefore, $\nu_s(fG)=-\nu_s(fG) = 0$ and $\nu(a^2_n f G) = -\nu(a^2_n f G) = 0$. By the fact that $a_n^3$ is odd, one immediately also obtain $\nu(a_n^3 G) = 0$.

The same argument applies to the situation where the spins are recentered and rescaled, as the rescaling and recentering is via constants with respect to the Gibbs measure.
\end{proof}

In fact, the marvelous symmetry of the ``rectangular sum'' is revealed in the the fact it is \emph{implicitly} centered around the concentrating value as
\[
    f = R_{1,3} - R_{1,4} - R_{2,3}  + R_{2,4} = \left(R_{1,3} - q^*\right) - \left(R_{1,4} - q^*\right) - \left(R_{2,3} - q^*\right) + \left(R_{2,4} - q^*\right)\,.
\]
This allows (after some effort using~\pref{lem:deriv} and the ``exact cancellations'' at $s = 0$) us to recover an extra factor of $O(n^{-1/2})$ in the cubical moments for various bulk statistics. Doing so requires obtaining concentration bounds for $f$ and $\hat{f}$, which is accomplished in the lemmata below.

\begin{lemma}[Even moments of ``rectangular sum'' via concentration]\label{lem:replicon-moments-from-mgf}
Fix $s\in[0,1]$. Assume the concentration implied by~\eqref{eq:conc-assumptions}.
Then for every integer $k\ge 1$,
\begin{equation}\label{eq:replicon-2k-moment}
\nu_s\!\left(|f|^{2k}\right)\le 4^{2k}\,\frac{k!}{\lambda^k}\,C_\lambda\,n^{-k}.
\end{equation}
Equivalently, 
\[
\nu_s(|f|^{2k}) = O(2k).
\]
\end{lemma}
\begin{prf}
Since the coefficients sum to zero, we have the exact identity
\[
    f = (R_{1,3}-q^*)-(R_{1,4}-q^*)-(R_{2,3}-q^*)+(R_{2,4}-q^*).
\]
Furthermore, for any $p\ge 1$ and any real numbers $x_1,\dots,x_m$,
\begin{equation}\label{eq:m-sum-ineq}
    \Big|\sum_{j=1}^m x_j\Big|^{p}\le m^{p-1}\sum_{j=1}^m |x_j|^{p}\,,
\end{equation}
which follows by a triangle and a $(1-1/p,p)$-H\"older inequality.
Apply \eqref{eq:m-sum-ineq} with $m=4$, $p=2k$, and
\[
x_1=R_{1,3}-q^*,\quad x_2=-(R_{1,4}-q^*),\quad x_3=-(R_{2,3}-q^*),\quad x_4=R_{2,4}-q^*.
\]
Then
\[
|f|^{2k}
\le
4^{2k-1}\Big(
|R_{1,3}-q^*|^{2k}+|R_{1,4}-q^*|^{2k}+|R_{2,3}-q^*|^{2k}+|R_{2,4}-q^*|^{2k}
\Big).
\]
Taking $\nu_s$ and using replica exchangeability (each pair overlap has the same law),
\[
\nu_s(|f|^{2k})
\le
4^{2k}\,\nu_s\big(|R_{1,2}-q^*|^{2k}\big).
\]
Using~\eqref{eq:conc-assumptions} on the above implies \eqref{eq:replicon-2k-moment}.
\end{prf}

A similar lemma can be shown to hold for the \emph{recentered} and \emph{rescaled} rectangular sums $\hat{f}$. This will be critical in the proof of~\pref{prop:D15} which will, after invocations of H\"older's inequality, rely on $L^{2p}$ bounds on $\hat{f}$ under $\nu_s(\cdot)$. The proof for this lemma, however, is far more involved and requires using the cavity interpolation itself to isolate the diagonal weighting in the scaled overlaps, and then using decoupling lemmata along with stability arguments to get the desired concentration.

\begin{lemma}[Even moments of the singly--scaled rectangular sum]\label{lem:hatf-moments}
Assume that for every $p\ge1$,
\begin{equation}\label{eq:Dmom-assump}
\sup_{1\le i\le n}\sup_{s\in[0,1]}\ \nu_s\!\big(D_{ii}^{\,p}\big)\le C(k,\beta,t) < \infty\,.
\end{equation}
Then, for every $k\in \Z_+$ there exists $C_k<\infty$, such that 
\begin{equation}\label{eq:hatf-even-moments}
\nu_s\left(|\hat{f}|^{2k}\right) \le \sup_{t\in[0,1]}\nu_t\big(|\hat f|^{2k}\big)\ \le\ C_k\,n^{-k}.
\end{equation}
\end{lemma}

\begin{proof}
This proof is by induction and relies on the cavity interpolation. To this end, note that the base case is considered as $k=0$ and the bound that $\nu_1(|\hat{f}|^{2\cdot 0}) = \nu_1(1) = 1 = O(n^{-0}) = O(1)$ trivially holds true. Therefore, the entire effort of the proof focuses on establishing the inductive step via a combination of H\"older-type bounds along with~\pref{lem:deriv} and~\pref{lem:rect-scaled} to handle the time $0$ terms in the Taylor expansion, and some uses of Young's inequality to handle the second-derivative term (and a few others).

Assume that, for every even integer less than $r < 2k$, $\sup_{s\in[0,1]}\nu_s(\hat{f}^{r}) \le O(n^{-r/2})$. Note the following simplifications using equality of distribution under index exchangeability and last-spin decoupling,
\begin{align*}
    \nu_1(\hat{f}^{2k}) &= \nu_1(\hat{f}^{2k-1}\hat{f}) = \nu_1\left(\hat{f}^{2k-1}\frac{1}{n}\sum_i D_{ii}a_i\right) = \nu_1\left(\hat{f}^{2k-1}D_{nn}a_n\right) \\
    &= \nu_1\left(\left(\hat{f}^- + \frac{1}{n}D_{nn}a_n\right)^{2k-1}D_{nn}a_n\right) \\
    &= \underbrace{\nu_1\left((\hat{f}^-)^{2k-1}D_{nn}a_n\right)}_{:= \mathrm{A}} + \underbrace{\frac{2k-1}{n}\nu_1\left(D^2_{nn}a^2_{n}(\hat{f}^-)^{2k-2}\right)}_{:= \mathrm{B}} + \underbrace{\sum_{r = 2}^{2k-1} \binom{2k-1}{r}\frac{1}{n^r} \nu_1\left((\hat{f}^-)^{2k-1-r}(D_{nn}a_n)^{r+1}\right)}_{:= \mathrm{C}}\,.  
\end{align*}
We now bound all three terms separately. 
\ppart{Bounding $\mathrm{C}$ via inductive hypothesis and H\"older's inequality} Applying $(p,q)$-H\"older's inequality to every term in the summand in $\mathrm{C}$ and using the fact that $|a_n| \le 4$ gives
\allowdisplaybreaks
\begin{align*}
    \sum_{r=2}^k \frac{1}{n^r}\nu_1\left((\hat{f}^-)^{2k-1-r}(D_{nn}a_n)^{r+1}\right) &\le \sum_{r=2}^k \frac{1}{n^r} \nu_1((\hat{f}^-)^{(2k-1-r)q})^{1/q}\nu_1(D^{p(r+1)}_{nn})^{1/p}\sup |a^{p(r+1)}_{n}| \\
    &\le_{q = \frac{2k-2}{2k-r-1},\,p=q/(q-1)}\sum_{r=2}^k \frac{1}{n^r} \nu_1(\hat{f}^{2k-2})^{\frac{2k-1-r}{2k-2}}\nu_1\left(D_{nn}^{q(r+1)/(q-1)}\right)^{(q-1)/q}C(k) \\
    &\le_{\text{\pref{lem:D-moments}\,+\,I.H.}} C'(k,\beta,t) \sum_{r=2}^k\frac{1}{n^r} O\left(\frac{1}{n^{(k-1)\frac{(2k-1-r)}{2k-2}}}\right) \\
    &= C'(k,\beta,t)\sum_{r=2}^k O\left(\frac{1}{n^{(2k-1+r)/2}}\right) \le_{r \ge 1} \frac{C'(k,\beta,t)}{n^k}\,.  
\end{align*}
\ppart{Bounding $\mathrm{B}$ via H\"older's and Young's inequality} Applying a $(p,q)$-H\"older's inequality again with $p=k$ and $q=k/(k-1)$, we obtain the following bound using~\pref{lem:D-moments} and Young's inequality (with $\eta > 0$ fixed later)
\allowdisplaybreaks
\begin{align*}
    \frac{1}{n}\nu_1\left(D^2_{nn}a^2_n (\hat{f}^-)^{2k-2}\right) &\le \sup |a^2_n|\frac{1}{n}\nu_1(D^{2k}_{nn})^{1/k}\nu_1\left((\hat{f}^-)^{2(k-1)\frac{k}{k-1}}\right)^{(k-1)/k} \\
    &\le_{\text{\pref{lem:D-moments}}\,+\,|a^2_n| \le 16} \frac{C(k,\beta,t)}{n}\nu_1\left((\hat{f}^-)^{2k}\right)^{(k-1)/k} \\
    &\le_{ab \le a^p/p + b^q/q} \frac{C^k(k,\beta,t)}{\eta^{k-1} n^k} + \eta\left(1-\frac{1}{k}\right)\nu_1((\hat{f}^-)^{2k}) \le \eta\nu_1(|\hat{f}^-|^{2k}) + O_\eta(n^{-k})\,,
\end{align*}
where Young's inequality is applied with $a=\frac{1}{\eta^{(k-1)/k}}\frac{C(k,\beta,t)}{n}$ and $b=\eta^{(k-1)/k}\nu_1\left((\hat{f}^-)^{2k}\right)^{(k-1)/k}$ and $p=k$ along with $q=k/(k-1)$.
\ppart{Bounding $\mathrm{A}$ via Gibbs-derivative stability along the cavity interpolation} We take the second-order Taylor expansion along the cavity interpolation, and invoke~\pref{lem:deriv} and~\pref{lem:rect-scaled} to evaluate the first two terms, and H\"older's inequality to evaluate the second derivative.
\allowdisplaybreaks
\begin{align*}
    \nu_1\left((\hat{f}^-)^{2k-1}D_{nn}a_n\right) &= \nu_0\left((\hat{f}^-)^{2k-1}D_{nn}a_n\right) + \nu^0_0\left((\hat{f}^-)^{2k-1}D_{nn}a_n\right) + \int_0^s(1-t)\nu''_t\left((\hat{f}^-)^{2k-1}D_{nn}a_n\right)\,dt\,.
\end{align*}
Using~\pref{lem:fact-s0} with $F^- = (\hat{f}^-)^{2k-1}$ on the first term gives
\[
    \nu_0\left((\hat{f}^-)^{2k-1}D_{nn}a_n\right) = \nu_0\left((\hat{f}^-)^{2k-1}\right)\E\left[\an{D_{nn}}_{Y^*}\an{(\eps^1-\eps^2)(\eps^3-\eps^4)}_{Y^*}\right] = 0\,.
\]
For the second term, we apply~\pref{lem:rect-scaled} with $F^- = (\hat{f}^-)^{2k-1}$ to obtain
\[
    \nu^0_0\left(\hat{f}^-)^{2k-1}D_{nn}a_n\right) = \beta^2\rho_2\nu_0\left( (\hat{f}^-)^{2k-1}f^-\right)\,.
\]
At this point, we note that
\[
    -\beta^2\rho_2\nu_0\left(|(\hat{f}^-)^{2k-1}f^-|\right) \le \beta^2\rho_2 \nu_0\left((\hat{f}^-)^{2k-1}f^-\right) \le \beta^2\rho_2\nu_0\left(|(\hat{f}^-)^{2k-1}f^-|\right)\,,
\]
and upper bound $\nu_0\left(|(\hat{f}^-)^{2k-1}f^-|\right)$ via an application of Young's inequality on $|(\hat{f}^-)^{2k-1}f^-|$ with $p=2k/(2k-1)$ and $q=2k$, along with linearity of expectation, as
\allowdisplaybreaks
\begin{align*}
    \beta^2\rho_2\nu_0\left(|(\hat{f}^-)^{2k-1}f^-|\right) &\le \beta^2\rho_2\eta_1\left(1-\frac{1}{2k}\right)\nu_0\left((\hat{f}^-)^{2k}\right)  + \frac{\beta^2\rho_2}{2k \eta^{2k-1}}\nu_0\left((f^-)^{2k}\right) \\
    &\le_{\text{\pref{lem:replicon-moments-from-mgf}}} \beta^2\rho_2\eta_1\nu_0\left((\hat{f}^-)^{2k}\right) + O_{\eta_1}(n^{-k}) \\
    &=_{\text{\pref{lem:stability-s}}} \beta^2\rho_2e^{24\beta^2  k}\eta_1\nu_1\left(|\hat{f}^-|^{2k}\right) + O(n^{-k})\,. 
\end{align*}
For the final term, we invoke~\pref{lem:deriv} twice and note that every term is a finite linear combination of terms of the form
\[
    \nu_t\left((\hat{f}^-)^{2k-1}D_{nn}B(\eps^\ell)\Xi_1\Xi_2\right)\, ,
\]
where $B(\eps^\ell)$ is an absolutely bounded function of the last spin in a finite number of replicas, and $\Xi_1,\Xi_2 \in \{ R^-_{ab}-q^*, M^-_a - m^*\}$. Overlap and magnetization concentration are stable along the cavity interpolation, and so by a $(2k/(2k-1),6k,6k,6k)$-H\"older's inequality, we obtain the following
\begin{align*}
    \nu_t\left((\hat{f}^-)^{2k-1}D_{nn}B(\eps^\ell)\Xi_1\Xi_2\right) &\le \nu_t\left(D^{6k}_{nn}B^6(\eps_\ell)\right)^{1/6k}\nu_t\left(\Xi_1^{6k}\right)^{1/6k}\nu_t\left(\Xi_2^{6k}\right)^{1/6k}\nu_t\left((\hat{f}^-)^{2k}\right)^{(2k-1)/2k} \\
    &\le_{\text{\pref{lem:D-moments}}\,+\,B(\eps^\ell) < \infty} \frac{C(k,\beta,t)}{n} \nu_t\left((\hat{f}^-)^{2k}\right)^{(2k-1)/2k}
\end{align*}
At this point, the exact same invocation of a Young's inequality (with $\eta_2 > 0$ fixed later) gives that
\allowdisplaybreaks
\begin{align*}
    \nu_t\left((\hat{f}^-)^{2k-1}D_{nn}B(\eps^\ell)\Xi_1\Xi_2\right) &\le \frac{C^k(k,\beta,t)}{\eta^{k-1}_2n^k} + \eta_2\left(1-\frac{1}{k}\right)\nu_t\left(|\hat{f}^-|^{2k}\right) \\
    &=_{\text{\pref{lem:stability-s}}} O_{\eta_2}(n^{-k}) + \eta_2\nu_1e^{24\beta^2 k}\left(|\hat{f}^-|^{2k}\right)\,.
\end{align*}
This immediately gives the final bound
\[
    \int_0^s(1-t)\nu''_t\left((\hat{f}^-)^{2k-1}D_{nn}a_n\right)\,dt \le  O_{\eta_2}(n^{-k}) + C'\eta_2\nu_1\left(|\hat{f}^-|^{2k}\right)\,,
\]
for some absolute constant $C' > 1$. Putting it together gives
\[
    \nu_s\left((\hat{f}^-)^{2k-1}D_{nn}a_n\right) \le \left(\beta^2\rho_2 + C'\eta_2\right)\nu_1\left(|\hat{f}^-|^{2k}\right) + O_{\eta_2}(n^{-k})\,. 
\]

\ppart{Combining the bounds} At this point, it merely remains to combine all the bounds and do a little bit of algebra\footnote{Note that $\nu_1(|\hat{f}^-|^{2k}) = \nu_1(|\hat{f}|^{2k}) + O(k)$. This is easy to see via a last-spin decoupling $\hat{f} = \hat{f}^- + \frac{1}{n}D_{nn}a_n$, followed by the binomial identity expansion. All higher-order terms are $O(n^{-k})$ via an argument exactly analogous to that used to bound $\mathrm{C}$, and the remaining terms are $\nu_1(|\hat{f}^-|^{2k})$ and $\frac{1}{n}\nu_1(|\hat{f}^-|^{2k-1}D_{nn}a_n) = \pm\frac{1}{n}\left(\beta^2\rho_2 + C'\eta_2\right)\nu_1\left(|\hat{f}^-|^{2k}\right) + O_{\eta_2}(n^{-k})$ by the argument to bound $\mathrm{A}$.}.  Note by all the bounds above (and the footnote) that
\allowdisplaybreaks
\begin{align*}
    \nu_1\left(\hat{f}^{2k}\right) &\le \left(\eta_1C_1(k,\beta)+ C_2(k,\beta)\eta_2 + \eta\right)\nu_1\left((\hat{f}^-)^{2k}\right) + O(n^{-k}) \\
    &=  \left(\eta_1C_1(k,\beta)+ C_2(k,\beta)\eta_2 + \eta\right)\left(1 + O\left(\frac{1}{n}\right)\right) \nu_1\left(\hat{f}^{2k}\right) + O(n^{-k})\, ,
\end{align*}
which immediately gives
\[
    \nu_1\left(|\hat{f}|^{2k}\right) \le \frac{1}{1- \eta_1C_1(k,\beta)+ C_2(k,\beta)\eta_2 + \eta - O(1/n)} O_{\eta,\eta_2}(n^{-k})\, ,
\]
which is well-defined with a choice of $\eta_2 = \eta = \eta_1 = 1/(10C(k,\beta))$ since $\beta^2\rho_2 < 1$ for $\beta < 1$. \qedhere
\end{proof}

\begin{remark}
    Note that one can transfer the bound $\nu_1(\hat{f}^{2k}) = O(k)$ to any point $s \in [0,1)$ using~\pref{lem:stability-s}.
\end{remark}

\paragraph{$O(4)$ and $O(5)$ bounds for bulk observables via H\"older estimates} In~\pref{lem:D16-corrected}, we prove that $\sup_{s \in [0,1]} |\nu''_s(F)| = O(4)$ whenever $F$ is a (scaled or unscaled) ``rectangular sum''. Before that, it will be critical to use various applications of H\"older's inequality to reason about the fluctuations of certain bulk observables.

\begin{lemma}[Hölder bounds for fluctuation monomials]\label{lem:holder-bookkeeping}
Fix a finite set of replica indices.
For $s\in[0,1]$, let
\[
\bar R^-_{ab}:=R^-_{ab}-q^*,
\qquad
\bar M^-_a:=M^-_a-m^*,
\qquad
f^-:=\bar R^-_{13}-\bar R^-_{14}-\bar R^-_{23}+\bar R^-_{24},
\]
and let $\hat f^-$ be the corresponding standardized rectangular sum.

Assume that for every $p\ge 1$ there exists $C_p<\infty$ such that uniformly in $s\in[0,1]$,
\begin{equation}\label{eq:holder-hyp-1}
\|\bar R^-_{ab}\|_{L^p(\nu_s)}+\|\bar M^-_a\|_{L^p(\nu_s)} \le C_p\,n^{-1/2},
\end{equation}
and that
\begin{equation}\label{eq:holder-hyp-2}
\|f^-\|_{L^p(\nu_s)} \le C_p\,n^{-1/2},
\qquad
\|\hat f^-\|_{L^2(\nu_s)} \le C\,n^{-1/2}.
\end{equation}

Let $\mathcal Q^-$ denote the linear span of quadratic bulk observables
\[
\bar R^-_{ab}\bar R^-_{cd},\qquad \bar R^-_{ab}\bar M^-_c,\qquad \bar M^-_a\bar M^-_b.
\]
Then:

\begin{enumerate}
\item[(i)] For every $p\ge1$ and every $G^-\in \mathcal Q^-$,
\begin{equation}\label{eq:G-Lp}
\|G^-\|_{L^p(\nu_s)} \le C_p\,n^{-1}\,.
\end{equation}

\item[(ii)] Let $B$ be any bounded random variable 
with $\|B\|_\infty\le C_B$.
Then for every $G^-\in\mathcal Q^-$ and every centered bulk fluctuation
$\Xi,\Xi_1,\Xi_2\in\{\bar R^-_{ab},\bar M^-_a\}$,
\begin{align}
\big|\nu_s(B\,G^-\,\Xi_1\,\Xi_2)\big| &\le C\,n^{-2} = O(4), \label{eq:holder-1}\\
\big|\nu_s(B\,f^-\,G^-\,\Xi)\big| &\le C\,n^{-2} = O(4), \label{eq:holder-2}\\
\big|\nu_s(B\,f^-\,G^-\,\Xi_1\,\Xi_2)\big| &\le C\,n^{-5/2} = O(5), \label{eq:holder-3} \\
\big|\nu_s(B\,\hat f^-\,G^-\,\Xi)\big| &\le C\,n^{-2} = O(4), \label{eq:holder-4}\\
\big|\nu_s(B\,\hat f^-\,G^-\,\Xi_1\,\Xi_2)\big| &\le C\,n^{-5/2} = O(5). \label{eq:holder-5}
\end{align}
\end{enumerate}
\end{lemma}

\begin{proof} All of the above bounds follow by applying H\"older's inequality with the appropriate set of H\"older conjugates.
\ppart{H\"older bound for \eqref{eq:G-Lp}} It is enough to prove the bound for a monomial $G^-=XY$ with $X,Y\in\{\bar R^-_{ab},\bar M^-_a\}$,
since $\mathcal Q^-$ is finite-dimensional and linear combinations scale the same via a triangle inequality.
For any $p\ge1$, H\"older gives
\[
\|XY\|_{L^p(\nu_s)}
\le \|X\|_{L^{2p}(\nu_s)}\|Y\|_{L^{2p}(\nu_s)}.
\]
Using \eqref{eq:holder-hyp-1},
\[
\|X\|_{L^{2p}(\nu_s)}\|Y\|_{L^{2p}(\nu_s)}
\le C_{2p}^2\,n^{-1}.
\]
This proves \eqref{eq:G-Lp} for monomials, and hence for every $G^-\in\mathcal Q^-$.

\ppart{H\"older bound for \eqref{eq:holder-1}}
As $B$ is uniformly bounded,
\[
\big|\nu_s(B\,G^-\,\Xi_1\,\Xi_2)\big|
\le \|B\|_\infty\,\nu_s(|G^-\,\Xi_1\,\Xi_2|).
\]
Apply Hölder with exponents $(2,4,4)$:
\[
\nu_s(|G^-\,\Xi_1\,\Xi_2|)
\le \|G^-\|_{L^2(\nu_s)}\|\Xi_1\|_{L^4(\nu_s)}\|\Xi_2\|_{L^4(\nu_s)}.
\]
By \eqref{eq:holder-hyp-1} and \eqref{eq:G-Lp},
\[
\|G^-\|_{L^2(\nu_s)}=O(n^{-1}),
\qquad
\|\Xi_1\|_{L^4(\nu_s)}=O(n^{-1/2}),
\qquad
\|\Xi_2\|_{L^4(\nu_s)}=O(n^{-1/2}),
\]
and so,
\[
\big|\nu_s(B\,G^-\,\Xi_1\,\Xi_2)\big| \le C\,n^{-1}\,n^{-1/2}\,n^{-1/2}=C\,n^{-2}\,.
\]

\ppart{H\"older bound for \eqref{eq:holder-2}}
Again using boundedness of $B$,
\[
\big|\nu_s(B\,f^-\,G^-\,\Xi)\big|
\le \|B\|_\infty\,\nu_s(|f^-\,G^-\,\Xi|).
\]
Applying H\"older with conjugates $(4,2,4)$:
\[
\nu_s(|f^-\,G^-\,\Xi|)
\le \|f^-\|_{L^4(\nu_s)}\|G^-\|_{L^2(\nu_s)}\|\Xi\|_{L^4(\nu_s)}.
\]
By \eqref{eq:holder-hyp-1}, \eqref{eq:holder-hyp-2} and \eqref{eq:G-Lp},
\[
\|f^-\|_{L^4(\nu_s)}=O(n^{-1/2}),
\qquad
\|G^-\|_{L^4(\nu_s)}=O(n^{-1}),
\qquad
\|\Xi\|_{L^4(\nu_s)}=O(n^{-1/2}),
\]
so
\[
\big|\nu_s(B\,f^-\,G^-\,\Xi)\big|\le C\,n^{-1/2}n^{-1}n^{-1/2}=C\,n^{-2}.
\]

\ppart{H\"older bound for \eqref{eq:holder-3}}
Similarly,
\[
\big|\nu_s(B\,f^-\,G^-\,\Xi_1\,\Xi_2)\big|
\le \|B\|_\infty\,\nu_s(|f^-\,G^-\,\Xi_1\,\Xi_2|).
\]
Apply Hölder with exponents $(4,4,4,4)$:
\[
\nu_s(|f^-\,G^-\,\Xi_1\,\Xi_2|)
\le \|f^-\|_{L^4(\nu_s)}\|G^-\|_{L^4(\nu_s)}\|\Xi_1\|_{L^4(\nu_s)}\|\Xi_2\|_{L^4(\nu_s)}.
\]
Using the same bounds as before,
\[
\big|\nu_s(B\,f^-\,G^-\,\Xi_1\,\Xi_2)\big|
\le C\,n^{-1/2}\,n^{-1}\,n^{-1/2}\,n^{-1/2}=C\,n^{-5/2}.
\]

\ppart{H\"older bound for \eqref{eq:holder-4}}
Now we use only the $L^2$ control on $\hat f^-$.
Again,
\[
\big|\nu_s(B\,\hat f^-\,G^-\,\Xi)\big|
\le \|B\|_\infty\,\nu_s(|\hat f^-\,G^-\,\Xi|).
\]
Apply H\"older with conjugates $(2,4,4)$:
\[
\nu_s(|\hat f^-\,G^-\,\Xi|)
\le \|\hat f^-\|_{L^2(\nu_s)}\|G^-\|_{L^4(\nu_s)}\|\Xi\|_{L^4(\nu_s)}.
\]
By \eqref{eq:holder-hyp-2}, \eqref{eq:G-Lp}, and \eqref{eq:holder-hyp-1},
\[
\|\hat f^-\|_{L^2(\nu_s)}=O(n^{-1/2}),
\qquad
\|G^-\|_{L^4(\nu_s)}=O(n^{-1}),
\qquad
\|\Xi\|_{L^4(\nu_s)}=O(n^{-1/2}),
\]
which immediately yields
\[
\big|\nu_s(B\,\hat f^-\,G^-\,\Xi)\big|\le C\,n^{-1/2}n^{-1}n^{-1/2}=C\,n^{-2}.
\]

\ppart{H\"older bound for \eqref{eq:holder-5}}
Finally, note that
\[
\big|\nu_s(B\,\hat f^-\,G^-\,\Xi_1\,\Xi_2)\big|
\le \|B\|_\infty\,\nu_s(|\hat f^-\,G^-\,\Xi_1\,\Xi_2|).
\]
Apply H\"older with conjugates $(2,6,6,6)$:
\[
\nu_s(|\hat f^-\,G^-\,\Xi_1\,\Xi_2|)
\le \|\hat f^-\|_{L^2(\nu_s)}\|G^-\|_{L^6(\nu_s)}\|\Xi_1\|_{L^6(\nu_s)}\|\Xi_2\|_{L^6(\nu_s)}.
\]
Using \eqref{eq:holder-hyp-2}, \eqref{eq:G-Lp}, and \eqref{eq:holder-hyp-1},
\[
\|\hat f^-\|_{L^2(\nu_s)}=O(n^{-1/2}),
\qquad
\|G^-\|_{L^6(\nu_s)}=O(n^{-1}),
\qquad
\|\Xi_j\|_{L^6(\nu_s)}=O(n^{-1/2}),
\]
which gives the final inequality,
\[
\big|\nu_s(B\,\hat f^-\,G^-\,\Xi_1\,\Xi_2)\big|
\le C\,n^{-1/2}\,n^{-1}\,n^{-1/2}\,n^{-1/2}
= C\,n^{-5/2}. \qedhere
\]
\end{proof}


\paragraph{Cubic moments are $O(4)$}
This is the central input akin to the one proved in~\cite[\S 1.10]{talagrand2010mean}. Namely, odd moments of the ``rectangular sums'' are one order smaller than naive concentration predicts, because the cavity equation forces them to be nearly multiplied by
$\beta^2\E[\sech^4(Y^*)]<1$ (which is a memory of the so-called dAT line).

For any bounded observable $F$ depending on finitely many replicas, the fundamental theorem of calculus implies
\[
    \nu(F)-\nu_0(F)-\nu_0^0(F)=\int_0^1(1-s)\,\nu_s''(F)\,ds,
\qquad
    \nu_0^0(F):=\left.\frac{d}{ds}\nu_s(F)\right|_{s=0+}.
\]
Therefore it suffices to prove
\[
    \sup_{s\in[0,1]}|\nu_s''(F)|\le C\,n^{-2}\,.
\]

Unfortunately,~\pref{lem:holder-bookkeeping} does \textbf{not} imply
\[
    \nu_s(f^-G^-)=O(4)
    \quad\text{or}\quad
    \nu_s(\hat f^-G^-)=O(4),
\]
because H\"older's inequality alone only gives $O(3)$ for those cubic expressions. The stronger $O(4)$ estimate for these \emph{replicon-weighted cubic terms} is exactly the content of~\pref{prop:D15}. We first prove a slightly easier lemma that has the same structure (and similar content) as the proof of~\pref{prop:D15}. The proof for~\pref{prop:D15} is then transparent in structure, and follows from this proof with exactly the same structure and similar invocations of H\"older's inequality.  

\begin{proposition}[Cubic moments of ``rectangular sums'' are $O(4)$]\label{prop:cubic-O4}
The ``rectangular sums'' satisfy the following scaling of (mixed and pure) third moments:  
\[
\nu(f^3)=O(4),\qquad \nu(f^2\hat f)=O(4),\qquad \nu(f\hat f^2)=O(4),\qquad \nu(\hat f^3)=O(4).
\]
\end{proposition}

\begin{proof}
We prove $\nu(f^3)=O(4)$; the other three bounds follow by the same argument, using invocations of~\pref{lem:D-moments},~\pref{lem:rect-scaled},~\pref{lem:rect-doublyscaled} and~\pref{lem:hatf-moments} every time $\hat f$ appears.

\ppart{Reducing to an off-diagonal term}
By exchangeability~\cite[\S 1.8, (1.240)]{talagrand2010mean},
\[
\nu(f^3)=\nu\left(\frac1n\sum_{i=1}^n a_i f^2\right)=\nu(a_n f^2).
\]
Write $f=f^-+\frac1n a_n$, hence
\[
f^2=(f^-)^2+\frac{2}{n}a_n f^-+\frac{1}{n^2}a_n^2.
\]
Thus
\begin{equation}\label{eq:f3-split}
\nu(f^3)=\nu\big(a_n(f^-)^2\big)+\frac{2}{n}\nu(a_n^2 f^-)+\frac{1}{n^2}\nu(a_n^3).
\end{equation}
By ~\pref{lem:antisym}, $\nu(a_n^2 f^-)=0$ and $\nu(a_n^3)=0$. 
So,
\begin{equation}\label{eq:f3-main}
\nu(f^3)=\nu\big(a_n(f^-)^2\big).
\end{equation}

\ppart{Applying the second-order Taylor expansion to $\nu(a_n(f^-)^2)$}
Set
\[
F:= (\varepsilon^1-\varepsilon^2)(\varepsilon^3-\varepsilon^4)\,(f^-)^2,
\quad\text{so that}\quad
\nu(F)=\nu(a_n(f^-)^2).
\]
The main claim is that
\begin{equation}\label{eq:two-term-O4}
\nu(F)=\nu_0(F)+\nu_0^0(F)+O(4).
\end{equation}

\ppart{Bounding $\nu''_s(F)$ uniformly} We prove that $\sup_{s \in [0,1]} |\nu''_s(F)| = O(4)$.
Differentiating $\nu(s)$ twice and applying~\pref{lem:deriv} twice gives the following
\begin{align*}
    \nu''_s(f^3) &= \beta^4\Bigg[\sum_{1 \le \ell_3 < \ell_4 < 6}\nu_s\left(H_6\eps^{\ell_3}\eps^{\ell_4}(R^-_{\ell_3,\ell_4}-q^*)\right) - 6\sum_{\ell_3=1}^6 \nu_s\left(H_6\eps^{\ell_3}\eps^7(R^-_{\ell_3,7}-q^*)\right) + 21\nu_s\left(H_6\eps^7\eps^8(R^-_{7,8}-q^*)\right) \\
    &\qquad\quad + \sum_{\ell_3=1}^6 \nu_s\left(H_6\eps^{\ell_3}(M^-_{\ell_3}-m^*)\right) - 6\sum_{\ell_3 = 1}^6 \nu_s\left(H_6\eps^7(M^-_7-m^*)\right)\Bigg]\, , 
\end{align*}
where $H_6$ is such that $\nu'_s(f^3) = \nu_s\left(H_6\right)$ (inferred from~\eqref{lem:deriv}) as
\begin{align*}
    H_6 &= \Bigg(\sum_{1 \le \ell < \ell' \le 4} f^3 \eps^\ell\eps^{\ell'}(R^-_{\ell,\ell'}-q^*) + 4\sum_{\ell = 1}^4 f^3\eps^{\ell}\eps^5(R_{\ell,5}^- - q^*) + 10f^3\eps^5\eps^6(R^-_{5,6}-q^*) \\
    &\qquad\qquad\qquad\qquad\quad + \sum_{\ell=1}^4 f^3\eps^\ell (M^-_\ell - m^*) - 4\left(f^3\eps^5(M^-_5 - m^*)\right)\Bigg)\,.
\end{align*}

Each term in $\nu_u''(F)$ is a finite linear combination
of expectations of $f^3$ multiplied by a product of two deviations, each deviation being either
$(R^-_{ab}-q^*)$ or $(M^-_a-m^*)$ and (up to) $4$ replicas. So, note that
\[
    \nu''_s(f^3) = \sum_{k=1}^{A} C_k\nu_s\left(f^3 p_k\left(\eps^1,\dots,\eps^8\right) a\,b\right)\, ,
\]
where $A, C_k > 0$, $a,b \in \{(R^-_{\ell,\ell'}-q^*), (M^-_{\ell}-m^*)\}_{(\ell, \ell') \in [8]^2}$ and $p_k(\cdot)$ are finite-degree polynomials in no more than $8$ variables. Since $\|p_k\| \le C$ (for some $C > 0$), applying a $(2,4,4)$-H\"older inequality yields
\[
    \left|\nu_s\left(f^3 p\left(\eps^1,\dots,\eps^8\right) a\,b\right)\right| \le C\left|\nu_s\left(f^3\,a\,b\right)\right| \le C\left(\nu_s(f^6)\right)^{1/2}\left(\nu_s(a^4)\right)^{1/4}\left(\nu_s(b^4)\right)^{1/4}\le_{\text{\eqref{eq:conc-assumptions}}} Cn^{-1}\left(\nu_s(f^6)\right)^{1/2}\,.
\]

Invoking~\pref{lem:replicon-moments-from-mgf} implies that $\nu_s(f^6)^{1/2} = O(n^{-3/2})$, and combining the bounds proves \eqref{eq:two-term-O4}.

\ppart{Using exact factorization to show $\nu_0(a_n (f^-)^2) = 0$} Note that $\nu_0(F)=0$ by~\pref{lem:fact-s0} because, conditional on $Y^*$,
$\E[(\varepsilon^1-\varepsilon^2)(\varepsilon^3-\varepsilon^4)\mid Y^*]=0$. Therefore, $\nu_0(f^3) = \nu_0(a_n (f^-)^2) = 0$.

\ppart{Using exact factorization to bound $\nu^0_0(a_n(f^-)^2) = O(4)$} For the derivative term, apply~\pref{lem:rect-unscaled} with $F^-=(f^-)^2$:
\[
\nu_0^0(F)=\beta^2\rho_4\;\nu_0\big((f^-)^3\big).
\]
\ppart{Comparing $\nu_0\left((f^-)^3\right)$ and $\nu_1\left(f^^3\right)$} We now show that
\[
    \nu_0((f^-)^3) - \nu(f^3) = \left(\nu_0((f^-)^3) - \nu_1((f^-)^3)\right) + \left(\nu_1((f^-)^3) - \nu_1(f^3)\right) = O(4)\,.
\]
Bounding the second term relies on simple algebra and applications of \textbf{(H1)} and~\pref{lem:antisym}. Specifically,
\begin{align*}
    \nu_1\left((f)^3 - (f^-)^3\right) &= \nu_1\left(\left(f^- + \frac{1}{n}a_n\right)^3 - (f^-)^3\right) = \frac{1}{n^3}\nu_1\left(a_n^3\right) + \frac{3}{n}\nu_1\left(a_n (f^-)^2\right) + \frac{3}{n^2}\nu_1\left(a^2_n f^-\right) \\
    &\le_{|a_n| \le 4,~\text{\pref{lem:antisym}}} \frac{C}{n^3} + \frac{C'}{n}\nu_1\left((f^-)^2\right) + 0 \\
    &\le_{\text{\textbf{(H1)}}} \frac{C}{n^3} + \frac{C'}{n^2} = O(4)\,.
\end{align*}
For the first term one applies the fundamental theorem of calculus along with~\pref{lem:deriv} and invokes Cauchy-Schwarz on every summand. Expanding,
\begin{align*}
    &\nu_1\left((f^-)^3\right) - \nu_0\left((f^-)^3\right) = \int_0^1 \frac{d}{ds}\nu_s\left((f^-)^3\right)ds \le \sup_{s \in [0,1]} \nu'_s\left((f^-)^3\right) \\
    &=_{\text{\pref{lem:deriv}}} \sup_{s \in [0,1]} \beta^2 \Bigg(\sum_{1 \le \ell < \ell' \le 4} \nu_s\left((f^-)^3 \eps^\ell\eps^{\ell'}(R^-_{\ell,\ell'}-q^*)\right) + 4\sum_{\ell = 1}^4 \nu_s\left((f^-)^3\eps^{\ell}\eps^5(R_{\ell,5}^- - q^*)\right) + 10\nu_s\left((f^-)^3\eps^5\eps^6(R^-_{5,6}-q^*)\right) \\
    &\qquad\qquad\qquad\qquad\quad + \sum_{\ell=1}^4 \nu_s\left((f^-)^3\eps^\ell (M^-_\ell - m^*)\right) - 4\nu_s\left((f^-)^3\eps^5(M^-_5 - m^*)\right)\Bigg)\,.
\end{align*}
Using \textbf{(H1)} and \textbf{(H2)} in conjunction with the fact that $\eps^\ell \in \{-1,1\}$ for every $\ell \ne \ell' \in \{1,2,3,4\}$, invocations of Cauchy-Schwarz give
\begin{align*}
    |\nu_s\left((f^-)^3 \eps^\ell\eps^{\ell'}(R^-_{\ell,\ell'}-q^*)\right)| &\le |\nu_s\left((f^-)^3(R^-_{\ell,\ell'}-q^*)\right)| \le \sqrt{\nu_s\left((f^-)^6\right)}\sqrt{\nu_s\left((R^-_{\ell,\ell'}-q^*)^2\right)} \le C n^{-3/2}n^{-1/2} = O(4),  \\
    |\nu_s\left((f^-)^3\eps^{\ell}\eps^5(R_{\ell,5}^- - q^*)\right)| &\le |\nu_s\left((f^-)^3(R_{\ell,5}^- - q^*)\right)| \le \sqrt{\nu_s\left((f^-)^6\right)}\sqrt{\nu_s\left(R^-_{\ell,5}-q^*)^2\right)} \le C n^{-3/2}n^{-1/2} = O(4), \\
    |\nu_s\left((f^-)^3\eps^5\eps^6(R^-_{5,6}-q^*)\right)| &\le |\nu_s\left((f^-)^3(R^-_{5,6}-q^*)\right)| \le \sqrt{\nu_s\left((f^-)^6\right)}\sqrt{\nu_s\left(R^-_{5,6}-q^*)^2\right)} \le Cn^{-3/2}n^{-1/2} = O(4),\\
    |\nu_s\left((f^-)^3\eps^\ell (M^-_\ell - m^*)\right)| &\le |\nu_s\left((f^-)^3 (M^-_\ell - m^*)\right) | \le \sqrt{\nu_s\left((f^-)^6\right)}\sqrt{\nu_s\left((M^-_\ell - m^*)^2\right)} \le Cn^{-3/2}n^{-1/2} = O(4), \\
    |\nu_s\left((f^-)^3\eps^5(M^-_5 - m^*)\right)| &\le |\nu\left((f^-)^3 (M^-_5 - m^*)\right)| \le \sqrt{\nu_s\left((f^-)^6\right)}\sqrt{\nu_s\left((M^-_5 - m^*)^2\right)} \le Cn^{-3/2}n^{-1/2} = O(4).
\end{align*}
Adding the bounds on the first and second term gives
\[
    \nu_1(f^3) = \nu_0\left((f^-)^3\right) + O(4)\,. 
\]
 
\ppart{Final bound on $\nu_1(f^3)$} Putting all the bounds together yields
\[
    \nu_1(f^3) = \nu_1(a_n(f^-)^2) = \nu_0(f^3) + \nu^0_0(f^3) + \beta^2\rho_4\,\nu(f^3)+ O(4).
\]
Since $\beta^2\rho_4<1$ and $\nu_0(f^3) = 0$, we conclude
\begin{align*}
    \nu_1(f^3) &=  0 + \beta^2 \rho_4\nu_1(f^3) + C\sup_{s \in [0,1]}|\nu''_s(f^3)| + O(4) \\
    &= \frac{C}{1-\beta^2\rho_4} O(4)\, ,
\end{align*}
which is well-defined below the dAT line ($\beta^2 \rho_4 < 1$).

\ppart{Mixed cubic moments}
The proofs of $\nu(f^2\hat f)=O(4)$ and $\nu(f\hat f^2)=O(4)$ and $\nu(\hat f^3)=O(4)$ are identical,
replacing one or more factors $f$ by $\hat f$ and using~\pref{lem:rect-scaled} (scaled cancellation)
whenever the last-spin factor is standardized. The antisymmetry cancellations provided by~\pref{lem:antisym}
continue to remove all ``diagonal'' contributions, and the Taylor remainder is $O(4)$ because the bulk factor
is always at least quadratic in rectangular sums. Concretely, note that
\[
    \nu(f^2\hat{f}) = \nu(D_{nn}a_n(f^-)^2),\quad \nu(f\hat{f}^2) = \nu(D_{nn}a_n f^-\hat{f}^-),\quad \nu(\hat{f}^3) = \nu(D_{nn}a_n (\hat{f}^-)^2)\,,
\]
where one uses exchangeability and the facts that $\nu(D_{nn}a^2_nf^-) = 0 = \nu(D_{nn}a^2_n\hat{f})$ and $\nu(D_{nn}a_n^3) = 0$.
Let $F$ denote the symmetrized function. Two applications of~\pref{lem:deriv} will yield, for any of the above terms, that $\nu''_s(F)$ is a linear combination of terms of the form
\[
    \nu_s(B_s D_{nn} G^- \Xi_1\Xi_2)\, ,
\]
where $B_s$ is an absolutely bounded function of the last spins, $G^- \in \{ (f^-)^2, (\hat{f}^-)^2, f^-\hat{f}^-\}$ and $\Xi_1, \Xi_2$ are centered bulk deviation terms. By applying a $(2,8,8,8,8)$-H\"older's inequality in conjunction with~\pref{lem:D-moments},~\pref{lem:replicon-moments-from-mgf}, \pref{lem:hatf-moments} and overlap/magnetization concentration, one obtains that $|\nu''_s(F)| = O(4)$. At $s=0$, $\nu_0(F) = 0$ since $\an{a_n \mid Y^*}_0 = \an{(\eps^1 - \eps^2)(\eps^3-\eps^4)\mid Y^*} = 0$, for any $F \in G^-$.
Note from the proof for $\nu(f^3) = O(4)$ that
\[
    \nu(f^3) = \nu_0((f^-)^3) + O(4)\, .
\]
The final step is to invoke~\pref{lem:rect-scaled} and observe the following ``towering'' chain of implications using the equality above for every $F \in G^-$,
\begin{align*}
    \nu^0_0(f^2\hat{f}) &= \nu^0_0(D_{nn}a_n (f^-)^2) =_{\text{\pref{lem:rect-scaled}}} \beta^2\rho_2\nu_0((f^-)^3) = \beta^2\rho_2\nu((f^-)^3) + O(4) = O(4)\,, \\
    \nu^0_0(f\hat{f}^2) &= \nu^0_0(D_{nn}a_nf^-\hat{f}^-) =_{\text{\pref{lem:rect-scaled}}} \beta^2\rho\nu_0(\hat{f}^-(f^-)^2) = \beta^2\rho_2\nu(f^2\hat{f}) + O(4) = O(4)\, , \\
    \nu^0_0(\hat{f}^3) &= \nu^0_0(D_{nn}a_n(\hat{f}^-)^2) =_{\text{\pref{lem:rect-scaled}}} \beta^2\rho_2\nu_0(f^-(\hat{f}^-)^2) = \beta^2\rho_2\nu(f\hat{f}^2) + O(4) = O(4)\, , 
\end{align*}
where the last equalities in the second and third equations follow by inserting the $O(4)$ bound implied by the second-order Taylor expansion for the previous $F$ term. The second-to-last equalities in the last two equations follow by expressing
\begin{align*}
    \nu(f^2\hat{f}) - \nu_0((f^-)^2\hat{f}^-) &= \left(\nu(f^2\hat{f}) - \nu((f^-)^2\hat{f}^-)\right) + \underbrace{(\nu((f^-)^2\hat{f}^-) - \nu_0((f^-)^2\hat{f}^-))}_{:=\int_0^1 \nu'_s(F) ds}\, , \\
    \nu(f\hat{f}^2) - \nu_0(f^-(\hat{f}^-)^2) &= \left(\nu(f\hat{f}^2) - \nu(f^-(\hat{f}^-)^2)\right) + \underbrace{\nu(f^-(\hat{f}^-)^2) - \nu_0(f^-(\hat{f}^-)^2)}_{:= \int_0^1 \nu'_s(F)ds}\, ,
\end{align*}
where the first terms in the equalities above can be shown to be small by decoupling the last spin and using symmetry under replica-index exchange to show certain terms are exactly $0$ while the remaining terms are the first derivatives, which can be controlled by applying~\pref{lem:deriv} and H\"older's inequality with~\pref{lem:D-moments},~\pref{lem:replicon-moments-from-mgf} and~\pref{lem:hatf-moments} to control the remainder term to be $O(4)$. 
\end{proof}

We now ``boost'' the lemma above to deal with the situation where we have a \emph{single} copy of the rectangular sum weighted by quadratic powers of linear combinations of bulk deviation terms such as $(R_{ab}-q^*)$ or $(M_{ab}-m^*)$. This lemma will be crucial to bound remainders in the \emph{exact} limiting estimate developed for cavity terms like $W^{(2)}_n$. The structure of the proof will be \emph{exactly} the same as that of the proof of~\pref{prop:cubic-O4}. 

\begin{proposition}[Cubic moments of ``rectangular sums'' with quadratic weights are $O(4)$]\label{prop:D15}
Fix a finite number of replicas and let $G^-\in \calQ$ be any quadratic bulk observable, where
\[
    \calQ := \mathsf{span}\left\{(R^-_{\ell,\ell'}-q^*)^2, (M^-_\ell - m^*)^2, (R^-_{\ell,\ell'}-q^*)(M^-_\ell - m^*)\right\}\,.
\]

Then, for every $s\in[0,1]$,
\begin{align*}
    \nu_s\big(f^-\,G^-\big)&=O(4)\,, \\
    \nu_s\left(\hat{f}^-\,G^-\right)&= O(4)\,.    
\end{align*}
The same statement holds with $f^-$ replaced by $\hat f^-$ (and with $G^-$ allowed to include standardized overlaps).
\end{proposition}
\begin{proof}
We prove the claim for $s=1$; the same argument holds for any $s\in[0,1]$
with constants uniform in $s$, since~\pref{lem:deriv} and the $s=0$ factorization/cancellation (\pref{lem:fact-s0}) are unchanged.

Fix $G^-\in\mathcal Q$. Since each centered bulk deviation has size $O(n^{-1/2})$ in all $L^p$ by \eqref{eq:conc-assumptions},
we have
\begin{equation}\label{eq:Gsize}
\|G^-\|_{L^p(\nu_s)} = O(n^{-1})\,.
\end{equation}

\ppart{Exchangeability converts $\nu(f^-G^-)$ into a cavity-site expectation}
Note that $f$ is permutation-invariant in the site index.
Also $G^-$ is a function of bulk overlaps/magnetizations, hence, also permutation-invariant in the site index -- this is because its definition takes an average of all terms $i \in [n-1]$, which is invariant under permutations of $\{1,\dots,n-1\}$,
and the difference to using $\{1,\dots,n\}$ is $O(n^{-1})$.

Therefore, by site exchangeability,
\begin{equation}\label{eq:exchange}
\nu(f\,G^-)=\nu\left(\frac1n\sum_{i=1}^n a_i\,G^-\right)=\nu(a_n\,G^-).
\end{equation}
Moreover, $f=f^-+\frac1n a_n$, hence using \eqref{eq:Gsize} and the fact that $|a_n| \le 4$ and an application of H\"older's inequality,
\begin{equation}\label{eq:fminus-replace}
\nu(f\,G^-)=\nu(f^-\,G^-)+\frac1n\nu(a_nG^-) = \nu(f^-\,G^-)+O(4). 
\end{equation}
Combining \eqref{eq:exchange} and \eqref{eq:fminus-replace} gives
\begin{equation}\label{eq:key-reduction}
\nu(f^-\,G^-)=\nu(a_n\,G^-)+O(4)\,.
\end{equation}

\ppart{Two-term Taylor expansion for $\nu(a_nG^-)$ with $O(4)$ remainder}
Set
\[
F:=(\varepsilon^1-\varepsilon^2)(\varepsilon^3-\varepsilon^4)\,G^- = a_nG^-.
\]

Write the exact Taylor identity
\begin{equation}\label{eq:Taylor}
\nu(F)=\nu_0(F)+\nu_0^0(F)+\int_0^1 (1-u)\,\nu_u''(F)\,ds.
\end{equation}

The goal is to prove that the remainder is
\begin{equation}\label{eq:remainder-O4}
\int_0^1 (1-s)\,\nu_s''(F)\,du = O(4).
\end{equation}
First, obtain $\nu_u'(F)$ via~\pref{lem:deriv}. Every term in $\nu_u'(F)$ is a finite linear combination of expectations of the form
\[
\nu_s\big(G^- \mathcal P(\varepsilon) \Xi\big),
\qquad \Xi\in\{\bar R^-_{ab},\bar M^-_a\},
\]
with $\mathcal P(\varepsilon)$ a bounded polynomial in finitely many replicas of last spins, containing $(\eps^1-\eps^2)(\eps^3-\eps^4)$.
Differentiating once more, every term in $\nu_u''(F)$ is a finite linear combination of expectations of the form
\begin{equation}\label{eq:structure-second-derivative}
\nu_s\big(G^- \mathcal P(\varepsilon) \Xi_1\Xi_2\big),
\qquad \Xi_1,\Xi_2\in\{\bar R^-_{ab},\bar M^-_a\}.
\end{equation}
Since all last-spin factors are bounded by an absolute constant $O(1)$ in absolute value, we bound the magnitude of each term in
\eqref{eq:structure-second-derivative} by
\[
C\,\nu_u(|G^-|\,|\Xi_1|\,|\Xi_2|).
\]
Using \eqref{eq:Gsize} and $\|\Xi_j\|_{L^p(\nu_u)}=O(n^{-1/2})$ from \eqref{eq:conc-assumptions}, an application of $(2,4,4)$-H\"older yields
\[
\nu_u(|G^-|\,|\Xi_1|\,|\Xi_2|)
\le \|G^-\|_{L^2(\nu_u)}\|\Xi_1\|_{L^4(\nu_u)}\|\Xi_2\|_{L^4(\nu_u)}
= O(n^{-1})\cdot O(n^{-1/2})\cdot O(n^{-1/2})=O(n^{-2})=O(4)\,.
\]
Since $\nu_u''(F)$ is a finite linear combination of such terms, $\nu_u''(F)=O(4)$,
and \eqref{eq:remainder-O4} follows.

\ppart{Evaluate $\nu_0(F)$ and $\nu_0^0(F)$ exactly via factorization}
~\pref{lem:fact-s0} at $s=0$ gives, conditionally on $Y^*$, the last spins are i.i.d.\ with mean $\tanh(Y^*)$. So $\E[\varepsilon^1-\varepsilon^2\mid Y^*]=0$ and, therefore,
\[
    \nu_0(F)=0.
\]
For $\nu_0^0(F)$, applying the unscaled rectangular cancellation lemma in~\pref{lem:rect-unscaled} with bulk factor $F^-=G^-$ gives
\[
    \nu_0^0(F)=\beta^2\rho_4\,\nu_0(G^-\,f^-)\,.
\]
Thus \eqref{eq:Taylor} and \eqref{eq:remainder-O4} give
\begin{equation}\label{eq:nu-anG}
\nu(a_nG^-)=\beta^2\rho_4\,\nu_0(G^-f^-)+O(4).
\end{equation}

\ppart{Replacing $\nu_0(G^-f^-)$ by $\nu(G^-f^-)$ up to $O(4)$ fluctuations}
We show
\begin{equation}\label{eq:replace-nu0-by-nu}
\nu(G^-f^-)-\nu_0(G^-f^-)=O(4).
\end{equation}
Indeed,
\[
\nu(G^-f^-)-\nu_0(G^-f^-)=\int_0^1 \nu'_s(G^-f^-)\,ds.
\]
Applying~\pref{lem:deriv} to $H:=G^-f^-$ yields that every term in $\nu'_s(H)$ is a finite linear combination of expectations of the form
\[
\nu_s\big(H\,\mathcal{P}(\eps)\,\Xi\big),
\qquad \Xi\in\{\bar R^-_{ab},\bar M^-_a\}.
\]
Now $H=G^-f^-$ is \emph{cubic} in centered deviations.
Hence, using $(2,4,4)$-H\"older along with \eqref{eq:Gsize} and~\pref{lem:replicon-moments-from-mgf},
\[
\big|\nu_s(H\, \mathcal{P}(\eps)\, \Xi)\big|
\le C\,\|G^-\|_{L^2(\nu_s)}\|f^-\|_{L^4(\nu_s)}\|\Xi\|_{L^4(\nu_s)}
= O(n^{-1}) O(n^{-1/2}) O(n^{-1/2})=O(4)\,.
\]
Integrating over $s$ yields \eqref{eq:replace-nu0-by-nu}.

Therefore, \eqref{eq:nu-anG} becomes
\begin{equation}\label{eq:nu-anG-2}
\nu(a_nG^-)=\beta^2\rho_4\,\nu(G^-f^-)+O(4).
\end{equation}

\ppart{Final stability equation}
Combining \eqref{eq:key-reduction} with \eqref{eq:nu-anG-2} gives
\[
\nu(f^-G^-)=\nu(a_nG^-)+O(4)=\beta^2\rho_4\,\nu(G^-f^-)+O(4)=\beta^2\rho_4\,\nu(f^-G^-)+O(4).
\]
Thus,
\[
(1-\beta^2\rho_4)\,\nu(f^-G^-)=O(4).
\]
Since $0<\beta<1/2$ and $\rho_4\le 1$, we have $1-\beta^2\rho_4 > 0$ and, so,
\[
\nu(f^-G^-)=O(4)\, ,
\]
yielding the final result.

\paragraph{Uniformity in $s\in[0,1]$}
For a general endpoint $s\in[0,1]$, repeat the same argument with $\nu$ replaced by $\nu_s$ throughout;
all estimates remain uniform because:
(i) the derivative identity has the same form,
(ii) $s=0$ factorization and the cancellation lemma are unchanged,
(iii) the moment bounds \eqref{eq:conc-assumptions} are assumed uniform along the interpolation.
\end{proof}

\begin{lemma}[$O(4)$ Taylor remainder for ``rectangular sums'']\label{lem:D16-corrected}
Let
\begin{align*}
F_1&:=(\varepsilon^1-\varepsilon^2)(\varepsilon^3-\varepsilon^4)\,f^-,
\\
F_2&:=(\varepsilon^1-\varepsilon^2)(\varepsilon^3-\varepsilon^4)\,\hat f^-,
\\
F_3&:=(\hat\varepsilon^1-\hat\varepsilon^2)(\hat\varepsilon^3-\hat\varepsilon^4)\,\hat f^-,
\qquad
\hat\varepsilon^\ell:=\frac{\varepsilon^\ell-\tanh(Y^*)}{\sech(Y^*)}\quad(s=0).
\end{align*}
Then
\[
\nu(F_i)=\nu_0(F_i)+\nu_0^0(F_i)+O(4)
\qquad\text{for }i\in\{1,2,3\}.
\]
\end{lemma}

\begin{proof}
We prove the three cases separately, using~\pref{lem:holder-bookkeeping} and~\pref{prop:D15}. The goal will be to, as is the approach in all prior lemmata, evaluate the Gibbs measure (and its derivative) at time $s=0$ and apply uniform bounds on the second derivative, concluding the bound via the Taylor expansion. So, for every type of rectangular sum $F$, we have
\[
\nu(F_1)-\nu_0(F_1)-\nu_0^0(F_1)=\int_0^1(1-s)\,\nu_s''(F_1)\,ds\,,
\]
and it suffices to show
\[
\nu_s''(F_1)=O(4)\,.
\]

\paragraph{Case 1: $F_1=(\varepsilon^1-\varepsilon^2)(\varepsilon^3-\varepsilon^4)\,f^-$} We apply~\pref{lem:deriv} twice to $\nu_s(F_1) = \nu_s(a_nf^-)$ and begin by looking at the structure of the terms in the linear combination.
Each application of~\pref{lem:deriv} inserts one centered bulk fluctuation
$\Xi\in\{\bar R^-_{ab},\bar M^-_a\}$ together with a bounded polynomial $\calP(\eps^1,\dots,\eps^k)$ in the last-spin variables of the replicas.
Hence, every term in $\nu_s''(F_1)$ is a finite linear combination of quantities of the form
\begin{equation}\label{eq:D16-F1-structure}
T_s(B,G^-):=\nu_s\big(B\,f^-\,G^-\big),
\end{equation}
where
\begin{itemize}
\item $B$ is a bounded random variable depending only on finitely many last spins, with $\|B\|_\infty\le C$;
\item $G^-\in\mathcal Q^-$ is quadratic in centered bulk fluctuations.
\end{itemize}
To conclude $\nu''_s(F_1) = O(4)$ it, therefore, suffices to show that every $T_s(B,G^-)$ in \eqref{eq:D16-F1-structure} is $O(4)$.
\ppart{Taylor-expanding the $T_s(B,G^-)$ terms}
Fix $B,G^-$ as above and denote $H:=B\,f^-\,G^-$. Then, $T_s(B,G^-)=\nu_s(H)$. Taylor expanding around time $0$ gives,
\begin{equation}\label{eq:D16-Taylor-H}
\nu_s(H)=\nu_0(H)+\nu_0^0(H)+\int_0^s (s-u)\,\nu_u''(H)\,du.
\end{equation}

\ppart{Bounding the $s=0$ term $\nu_0(H)$ as $O(4)$}
By exact factorization at $s=0$ via~\pref{lem:fact-s0},
\[
\nu_0(H)
=
\underbrace{\E\Big[\big\langle B\big\rangle_{Y^*}\Big]}_{=:c_0(B)}\, \nu_0(f^-\,G^-)\,.
\]
Boundedness of $B$ implies that $|c_0(B)|\le C$ and~\pref{prop:D15} implies $\nu_0(f^-\,G^-)=O(4)$. Therefore,
\begin{equation}\label{eq:D16-F1-s0}
\nu_0(H)=O(4).
\end{equation}

\ppart{Bounding the first derivative term $\nu_0^0(H)$ as $O(4)$}
Applying~\pref{lem:deriv} at $s=0$ to the observable $H=B f^- G^-$ yields a linear combination of terms of the form
\[
\nu_0\big(B'\,f^-\,G^-\,\Xi\big),
\qquad
\Xi\in\{\bar R^-_{ab},\bar M^-_a\},
\]
where $B'$ is again a bounded polynomial in finitely many last-spin variables.
Again, by~\pref{lem:fact-s0} at $s=0$,
\[
\nu_0\big(B'\,f^-\,G^-\,\Xi\big)
=
\underbrace{\E\Big[\big\langle B'\big\rangle_{Y^*}\Big]}_{=:c_1(B')}
\,
\nu_0\big(f^-\,G^-\,\Xi\big)\,.
\]
Using~\pref{lem:holder-bookkeeping}, specifically \eqref{eq:holder-2} with $B\equiv 1$, gives
\[
\nu_0\big(f^-\,G^-\,\Xi\big)=O(4).
\]
Since $\nu_0^0(H)$ is a finite linear combination of $O(4)$ terms,
\begin{equation}\label{eq:D16-F1-deriv0}
\nu_0^0(H)=O(4)\,.
\end{equation}

\ppart{The Taylor remainder is $O(5)$}
Differentiate $\nu_s(H)$ once again via~\pref{lem:deriv} for the observable $ B f^- G^- \Xi$.
Each term in $\nu_u''(H)$ is then a finite linear combination of expectations of the form
\[
\nu_u\big(B''\,f^-\,G^-\,\Xi_1\,\Xi_2\big),
\qquad
\Xi_1,\Xi_2\in\{\bar R^-_{ab},\bar M^-_a\},
\]
with $B''$ bounded. Invoking~\pref{lem:holder-bookkeeping}, specifically \eqref{eq:holder-3},
\[
    \nu_u\big(B''\,f^-\,G^-\,\Xi_1\,\Xi_2\big)=O(5)\,.
\]
Therefore,
\begin{equation}\label{eq:D16-F1-remainder}
\int_0^s (s-u)\,\nu_u''(H)\,du = O(5)\,.
\end{equation}

\ppart{Concluding the final bound for $\nu(F_1)$}
From \eqref{eq:D16-Taylor-H}, \eqref{eq:D16-F1-s0}, \eqref{eq:D16-F1-deriv0}, and
\eqref{eq:D16-F1-remainder}, we get
\[
T_s(B,G^-)=\nu_s(H)=O(4)\,.
\]
Since $\nu_s''(F_1)$ is a finite linear combination of terms $T_s(B,G^-)$, we have $\nu_s''(F_1)=O(4)$.
Therefore,
\[
\nu(F_1)=\nu_0(F_1)+\nu_0^0(F_1)+O(4).
\]

\paragraph{Case 2: $F_2=(\varepsilon^1-\varepsilon^2)(\varepsilon^3-\varepsilon^4)\,\hat f^-$} The proof is (basically) identical, with every occurrence of $f^-$ replaced by $\hat f^-$.
Indeed, after differentiating twice using~\pref{lem:deriv}, every term in $\nu_s''(F_2)$ is a finite linear combination of
\[
\tilde T_s(B,G^-):=\nu_s\big( B\,\hat f^-\,G^-\big),
\]
where $B$ is bounded and $G^-\in\mathcal Q^-$ is quadratic in bulk deviation terms. Taylor expanding $\tilde T_s(B,G^-)$ around time $0$ gives
\[
\tilde T_s(B,G^-)=\nu_0(\tilde H)+\nu_0^0(\tilde H)+\int_0^s (s-u)\,\nu_u''(\tilde H)\,du,
\qquad \tilde H:=B\,\hat f^-\,G^-.
\]

At $s=0$, factorization gives
\[
\nu_0(\tilde H)=\E[\langle B\rangle_{Y^*}]\,\nu_0(\hat f^-\,G^-),
\]
and~\pref{prop:D15}, with $f$ replaced by $\hat{f}$ and $\calQ$ including the rescaled-and-centered bulk deviation terms, yields
\[
    \nu_0(\hat f^-\,G^-)=O(4),
\]
which immediately implies that \(\nu_0(\tilde H)=O(4)\). For the derivative term, every contribution to $\nu_0^0(\tilde H)$ is a bounded last-spin coefficient
times $\nu_0(\hat f^-\,G^-\,\Xi)$, and~\pref{lem:holder-bookkeeping}, specifically \eqref{eq:holder-4}, gives
\[
\nu_0(\hat f^-\,G^-\,\Xi)=O(4)\,,\] 
which implies that \(\nu_0^0(\tilde H)=O(4)\) and it remains to uniformly bound the second derivative. For the second-derivative remainder, every term is of the form
\[
\nu_u(B''\,\hat f^-\,G^-\,\Xi_1\,\Xi_2),
\]
with $B''$ bounded, and~\pref{lem:holder-bookkeeping} (via \eqref{eq:holder-5}) gives
\[
\nu_u(B''\,\hat f^-\,G^-\,\Xi_1\,\Xi_2)=O(5).
\]
Putting this altogether gives that
\[
\nu(F_2)=\nu_0(F_2)+\nu_0^0(F_2)+O(4).
\]

\paragraph{Case 3: $F_3=(\hat\varepsilon^1-\hat\varepsilon^2)(\hat\varepsilon^3-\hat\varepsilon^4)\,\hat f^-$} For this case, the terms that depends on the last-spin is \emph{not} bounded, but it has moments of all orders.
More specifically, at $s=0$,
\[
\hat\varepsilon^\ell=\frac{\varepsilon^\ell-\tanh(Y^*)}{\sech(Y^*)},
\]
which implies that
\[
|\hat\varepsilon^\ell|
\le \frac{|\varepsilon^\ell|+|\tanh(Y^*)|}{\sech(Y^*)}
\le 2\cosh(Y^*).
\]
Therefore,
\[
|\hat{a}_n| = \Big|(\hat\varepsilon^1-\hat\varepsilon^2)(\hat\varepsilon^3-\hat\varepsilon^4)\Big|
\le C\,\cosh^4(Y^*),
\]
and since $Y^* \sim \calN(\beta^2 m^* + t, \beta^2q^* + t)$, the term \(\cosh^4(Y^*)\) has finite moments of all orders, i.e., $\E[\cosh(Y^*)^{4k}] \le C(\beta,t,k) < \infty$ for every $k \in \Z_+$.
Thus, every last-spin factor appearing in derivatives of $F_3$ satisfies a bound of $O(1)$ with an adequate application of H\"older's inequality.

Differentiating $\nu_s(F_3)$ twice using~\pref{lem:deriv} yields that every term in $\nu_s''(F_3)$ is a finite linear combination of
\[
\hat T_s(B_s,G^-):=\nu_s(\hat{a}_n B_s\,\hat f^-\,G^-),
\]
where \(G^-\in\mathcal Q^-\) is quadratic and $B_s$ is a bounded function in a finite number of replicas of the last spin. Taylor expanding \(\hat T_s(B_s,G^-)\) around \(s=0\) gives
\[
\hat T_s(B_s,G^-)=\nu_0(\hat H)+\nu_0^0(\hat H)+\int_0^s (s-u)\,\nu_u''(\hat H)\,du,
\qquad \hat H:=B_s\,\hat f^-\,G^-.
\]

At $s=0$,~\pref{lem:fact-s0}, the boundedness of $B_s$, and the fact that $\E_{Y^*}[|\hat{a}_n|] \le C\E_{Y^*}[\cosh^4(Y^*)] \le C(\beta,t)$ gives
\begin{align*}
    \nu_0\left(B_s\hat{a}_n\hat{f}^-G^-\right) &= \nu_0(\hat{f}^-G^-)\E_{Y^*}\left[\an{\hat{a}_nB_s}_{Y^*}\right] \\
    &\le_{\text{\pref{prop:D15}}} O(4)\sup |B_s| \E_{Y^*}[\cosh^4(Y^*)] \le C(\beta,t) O(4) = O(4)\,.
\end{align*}

For the derivative term, applying~\pref{lem:deriv} implies that \(\nu_0^0(\hat H)\) is a finite linear combination of finite coefficients times terms of the form
\(\nu_0((B_s\hat{a}_n)\hat f^-\,G^-\,\Xi)\). Then, by a $(2,6,6,6)$-H\"older's inequality,
\begin{align*}
    \nu_0((B_s\hat{a}_n)\hat{f}^-G^-\Xi) &\le \nu_0(|\hat{a}_n|^2 B^2_s)\nu_0\left((\hat{f}^-)^6\right)^{1/6}\nu_0\left((G^-)^6\right)^{1/6}\nu_0\left(\Xi^6\right)^{1/6} \\
    &\le_{\text{\pref{lem:hatf-moments},~\pref{lem:D-moments},~\pref{lem:holder-bookkeeping}}}\sup|B^2_s|\nu_0(\hat{a}^2_n)^{1/2}O(n^{-1/2})O(n^{-1})O(n^{-1/2}) \\
    &\le_{\nu_0(|\hat{a}_n|^2)\le \E_{Y^*}(\cosh^8(Y^*))} C(\beta, t) O(n^{-2}) = O(4)\,.
\end{align*}
This implies that \(\nu_0^0(\hat H)=O(4)\). For the second-derivative remainder, two applications of~\pref{lem:deriv} implies that every term is of the form
\[
    \nu_u(\hat{a}_n\,B_u\,\hat f^-\,G^-\,\Xi_1\,\Xi_2),
\]
where \(B_u'\) is absolutely bounded. An application of $(2,8,8,8,8)$-H\"older's inequality with invocations of~\pref{lem:D-moments},~\pref{lem:replicon-moments-from-mgf},~\pref{lem:hatf-moments} and~\pref{lem:holder-bookkeeping} yields that each such term is $O(4)$. Therefore,
\[
    \nu(F_3)=\nu_0(F_3)+\nu_0^0(F_3)+O(4). \qedhere
\]
\end{proof}

\begin{proposition}[CLT bounds for the $(U_n,V_n,W_n)$ system~{\cite[Mild refinement of Theorem 1.10.1]{talagrand2010mean}}]\label{prop:D17}
 Define
\[
U_n:=\nu(f^2),\qquad V_n:=\nu(f\hat f),\qquad W_n:=\nu(\hat f^2).
\]
Then, 
\begin{align*}
(1-\beta^2\rho_4)\,U_n &= \frac{4\rho_4}{n}+O(4),\\
(1-\beta^2\rho_4)\,V_n &= \frac{4\rho_2}{n}+O(4),\\
W_n &= \frac{4}{n}+\beta^2\rho_2\,V_n+O(4).
\end{align*}
Equivalently,
\[
U_n=\frac{4\rho_4}{n(1-\beta^2\rho_4)}+O(n^{-2}),
\qquad
V_n=\frac{4\rho_2}{n(1-\beta^2\rho_4)}+O(n^{-2}),
\qquad
W_n=\frac{4}{n}\Big(1+\frac{\beta^2\rho_2^2}{1-\beta^2\rho_4}\Big)+O(n^{-2}).
\]
\end{proposition}

\begin{proof}
This proof uses site-exchangeability to along with a decoupling of the last spin to obtain \emph{exact} estimates for the time $0$ terms. The remainder term is shown to be $O(4)$ using~\pref{lem:D16-corrected}.

\ppart{Exchangeability and diagonal/off-diagonal decompositions}
Define
\[
a_i:=(\sigma_i^1-\sigma_i^2)(\sigma_i^3-\sigma_i^4),
\qquad
\hat a_i:=(\hat\sigma_i^1-\hat\sigma_i^2)(\hat\sigma_i^3-\hat\sigma_i^4),
\]
so that $f=\frac1n\sum_i a_i$ and $\hat f=\frac1n\sum_i\hat a_i$.
By exchangeability of sites,
\[
U_n=\nu(f^2)=\nu\Big(\frac1n\sum_i a_i f\Big)=\nu(a_n f),
\qquad
V_n=\nu(f\hat f)=\nu(a_n\hat f),
\qquad
W_n=\nu(\hat f^2)=\nu(\hat a_n\hat f).
\]
Write $f=f^-+\frac1n a_n$ and $\hat f=\hat f^-+\frac1n\hat a_n$ where $f^-,\hat f^-$ depend only on the bulk spins.
Then
\begin{align}
U_n&=\nu(a_n f^-)+\frac1n\nu(a_n^2), \label{eq:U-split-D17}\\
V_n&=\nu(a_n\hat f^-)+\frac1n\nu(a_n\hat a_n), \label{eq:V-split-D17}\\
W_n&=\nu(\hat a_n\hat f^-)+\frac1n\nu(\hat a_n^2). \label{eq:W-split-D17}
\end{align}

\ppart{Computing the self-similar terms via exact factorization} We compute the self-similar terms in the expansion above at time $s=0$ via application of~\pref{lem:fact-s0} in conjunction with the observation that $\an{\sigma_n}_{Y^*} = \tanh(Y^*)$ to conclude that
\[
    \nu_0(a_n^2)=4\rho_4,\qquad \nu_0(a_n\hat a_n)=4\rho_2,\qquad \nu_0(\hat a_n^2)=4\,.
\]
We now invoke~\pref{lem:deriv} to write $\nu^0_0(F)$ exactly as a term controlled via~\pref{lem:mean-bias-On-1} and \eqref{eq:stability-diff}, and invoke~\pref{lem:deriv} twice to control $\nu''_s(F)$ via~\pref{lem:stab-dev-along-s} and concentration, for $F \in \{a^2_n, a_n\hat{a}_n, \hat{a}^2_n\}$. Specifically,
\[
    \frac{1}{n}\nu(F) = \frac{1}{n}\nu_0(F) + \frac{1}{n}\nu^0_0(F) + \frac{1}{n}\int_0^1(1-s)\nu''_s(F)ds\,.
\]
\begin{align*}
    \nu^0_0(F) &= \beta^2\sum_{1\le \ell < \ell' \le 4}\nu_0(F\eps^\ell\eps^{\ell'}(R^-_{\ell,\ell'}-q^*)) -4\beta^2\sum_{\ell=1}^4\nu_0(F\eps^\ell (R^{-}_{\ell,5}-q^*)) + 10\beta^2\nu_0(F\eps^5\eps^6(R^-_{5,6}-q^*)) \\
    &\quad\quad + \beta^2\sum_{\ell=1}^4 \nu_0(F\eps^\ell (M^-_{\ell}-m^*)) - 4\beta^2\nu_0(F\eps^5(M^-_5-q^*))\,.
\end{align*}
Now, note that, for any choice of $F \in \{a^2_n,a_n\hat{a}_n,\hat{a}^2_n\}$, the terms $F\eps^\ell$ or $F\eps^\ell\eps^{\ell'}$ are, conditioned on $Y^*$, functions that are independent of the first $(n-1)$ spins at $s=0$. Then, by~\pref{lem:fact-s0},
\[
    \nu_0(F\eps^\ell\eps^{\ell'}\Xi^-) = \nu_0(\Xi^-)\E_{Y^*}\left[F\eps^\ell\eps^{\ell'}\right]\,,
\]
and
\[
    \nu_0(F\eps^\ell\Xi^-) = \nu_0(\Xi^-)\E_{Y^*}\left[F\eps^\ell\right]\,,
\]
where $\Xi^- \in \{R_{ab}-q^*, M^-_a - m^*\}$. Now, in the event that $F = a^2_{nn}$, one has that $|F\eps^{\ell}\eps^{\ell'}| \le C$, and $|F\eps^\ell| \le C'$. In the event that $F$ has a power of $\hat{a}_{n} = D_{nn}a_n$, then one can invoke a $(2,2)$-H\"older's inequality in conjunction with~\pref{lem:D-moments} to conclude that $\E_{Y^*}\left[D^r_{nn}a^r_n\eps^\ell\eps^{\ell'}\right] \le \E[D^{2r}_{nn}]^{1/2}\E[a^{2r}_{nn}]^{1/2} \le  C(\beta,r,t) < \infty$ for $r \in \{1,2\}$. Then, using~\pref{lem:mean-bias-On-1} in conjunction with \eqref{eq:stability-diff} immediately yields that
\[
    \nu_0(\calP \Xi^-) = \nu_0(\Xi^-)\E_{Y^*}[\calP] \le C(\beta,r,t) e^{C'(\beta)}O(1/n) = O(1/n)\,,
\]
where $\calP \in \{F\eps^\ell\eps^{\ell'}, F\eps^\ell\}$. To conclude, note that two applications of~\pref{lem:deriv} yield that $\nu''_s(F)$ is a finite linear sum of terms of the form
\[
    \nu_s( F\, B\, \Xi^-_1\,\Xi^2_- ) \le_{\text{(2,4,4)-H\"older's}} \nu_s((F\,B)^{2})^{1/2}\norm{\Xi^-_1}_{L^4(\nu_s)}\norm{\Xi^-_2}_{L^4(\nu_s)} \le O(n^{-1})\, ,
\]
where we use~\pref{lem:stab-dev-along-s} and~\pref{lem:D-moments} in conjunction with the fact that $B$ is an absolutely bounded function of the last spin of a finite number of replicas. Putting the bounds on the first derivative and second derivative above together with the explicit evaluations of $\nu_0(F)$ immediately yields 
\[
\frac1n\nu(a_n^2)=\frac{4\rho_4}{n}+O(n^{-2}),
\qquad
\frac1n\nu(a_n\hat a_n)=\frac{4\rho_2}{n}+O(n^{-2}),
\qquad
\frac1n\nu(\hat a_n^2)=\frac{4}{n}+O(n^{-2}).
\]

\ppart{Bounding non-similar terms in $U_n$ via~\pref{lem:D16-corrected}} Note that
\[
\nu(a_n f^-)=\nu\big((\varepsilon^1-\varepsilon^2)(\varepsilon^3-\varepsilon^4)\,f^-\big)=\nu(F_1)\,,
\]
and so~\pref{lem:D16-corrected} immediately implies that
\[
    \nu(F_1)=\nu_0(F_1)+\nu_0^0(F_1)+O(4).
\]
By~\pref{lem:fact-s0}, $\nu_0(F_1)=\nu_0(f^-)\E_{Y^*}\left[\an{a_n}_{Y^*}\right] = 0$, and~\pref{lem:rect-scaled} implies that
\[
\nu_0^0(F_1)=\beta^2\rho_4\,\nu_0\big((f^-)^2\big).
\]
At this point, we let $B$ be an absolutely bounded function in the last spin of $4$ replicas, and observe that
\begin{align*}
    &\nu(f^2) - \nu_0((f^-)^2) = \left(\nu(f^2) - \nu((f^-)^2)\right) + \left(\nu((f^-)^2) - \nu_0((f^-)^2)\right) \\\
    &= \left(\frac{2}{n}\nu(a_nf^-) + \frac{1}{n^2}\nu(a^2_n)\right) + \int_0^1 \nu'_s((f^-)^2)\,ds \\
    &=_{\text{\pref{lem:deriv},\,}|a_n| \le 4,\text{\pref{lem:D16-corrected}}} \frac{2\beta^2\rho_4}{n}\nu_0((f^-)^2) + \int_0^1\left(\mathsf{span}\{\nu_s(B f^- (f^-\Xi^-)) \}\right)\,ds + O(4) \\
    &=_{\text{\pref{lem:replicon-moments-from-mgf},~\eqref{eq:stability-diff}}} O(4) + O\left(\nu_s(B\left( f^- \Xi_1^- \Xi^-_2\right))\right)\,,
\end{align*}
where we used the fact that $f^-$ is a sum of 4 centered overlap terms, and the $\Xi^-_a$ represented centered bulk deviation terms. Now, if $\sup_{s \in [0,1]}\left(\nu_s(B(f^-\Xi^-_1\Xi^-_2)) - \nu_s(B)\nu_s(f^-\Xi^-_1\Xi^-_2)\right) = O(4)$, then by~\pref{prop:D15} we can conclude that $\nu_s(Bf^-\Xi^-_1\Xi^-_2) = O(4)$ and so $\nu_s(f^2) = \nu_s((f^-)^2) + O(4)$. So, to obtain the bound note for $s=0$ by~\pref{lem:fact-s0},
\[
   | \nu_0(Bf^-\Xi^-_1\Xi^-_2) | = \E[\an{B}_{Y^*}]\nu_0(f^-\Xi^-_1\Xi^-_2) =_{|B|\le C,\,\text{\pref{prop:D15}}} O(4)\,.
\]
It now suffices to bound $|\left(\nu_s(Bf^-\Xi^-_1\Xi^-_2) - \nu_s(B)\nu_s(f^-\Xi^-_1\Xi^-_2)\right)'| = O(4)$ and apply FTOC. Note that, by~\pref{lem:deriv}, $ \nu'_s(Bf^-\Xi^-_1\Xi^-_2)$ is a finite linear combination of terms of the form
\[
   |\nu_s(B' f^-\Xi^-_1\Xi^-_2\Xi)| \le_{(2,8,8,8,8)\text{-H\"older's}} \nu_s(B'^2)(\nu_s(\Xi^8))^{4/8} \le O(4)\,.
\]
Similarly, applying~\pref{lem:deriv}  gives that $\nu'_s(B)$ is a finite linear combination of terms $\nu_s(B\Xi)$ bounded by $O(n^{-1/2})$, and that $\nu'_s(f^-\Xi^1_i\Xi^2_i)$ is a finite linear combination of terms $\nu_s(f^-\Xi^-_1\Xi^-_2\Xi)$ bounded by $O(n^{-2})$. By the product rule, this implies that
\begin{align*}
    \left(\nu_s(Bf^-\Xi^-_1\Xi^-_2) - \nu_s(B)\nu_s(f^-\Xi^-_1\Xi^-_2)\right)' &= \nu'_s(Bf^-\Xi^-_1\Xi^-_2) - \nu'_s(B)\nu_s(f^-\Xi^-_1\Xi^-_2) - \nu_s(B)\nu'_s(f^-\Xi^1_i\Xi^2_i) \\
    &= O(4) + O(4) + O(4) = O(4)\,,
\end{align*}
which yields the desired bound. Putting it altogether gives that
\[
    \nu(a_n f^-)=\beta^2\rho_4\,U_n+O(4)\,
\]
and inserting into \eqref{eq:U-split-D17} along with the diagonal term immediately gives
\[
    \,U_n=\frac{4\rho_4}{(1-\beta^2\rho_4)n}+O(4).
\]

\ppart{Bounding non-similar terms in $V_n$ via~\pref{lem:D16-corrected}}
Similarly,
\[
\nu(a_n\hat f^-)=\nu\big((\varepsilon^1-\varepsilon^2)(\varepsilon^3-\varepsilon^4)\,\hat f^-\big)=\nu(F_2).
\]
By~\pref{lem:D16-corrected},
\[
\nu(F_2)=\nu_0(F_2)+\nu_0^0(F_2)+O(4).
\]
Again,~\pref{lem:fact-s0} gives that $\nu_0(F_2)=\nu_0(\hat{f}^-)\E_{Y^*}[\an{a_n}_{Y^*}] = 0$ by factorization at $s=0$. Furthermore, using~\pref{lem:rect-unscaled} gives
\[
\nu_0^0(F_2)=\beta^2\rho_4\,\nu_0(\hat f^-\,f^-).
\]
And $\nu_0(\hat f^-\,f^-)=V_n+O(4)$. This follows by a strategy similar to the one used to show that $\nu_0(f^2) = \nu_1(f^2) + O(4)$, and we give a terse outline below.
\begin{align*}
    &\nu(f\hat{f}) - \nu_0(f\hat{f}) = \left(\nu(f\hat{f}) - \nu(f^-\hat{f}^-)\right) + \left(\nu(f^-\hat{f}^-) - \nu_0(f\hat{f}^-)\right) \\
    &= \left(\frac{1}{n}\nu(\hat{a}_nf^-) + \frac{1}{n}\nu(a_n\hat{f}^-) + \frac{1}{n^2}\nu_0(D_{nn}a^2_n)\right) + \int_0^1 \nu'_s(f\hat{f}^-)\,ds \\
    &=_{\text{\pref{lem:deriv},\,\pref{lem:D16-corrected},\,\pref{lem:D-moments}}} \frac{1}{n}\nu_0(\hat{a}_nf^-) + \frac{1}{n}\nu^0_0(\hat{a}_nf^-) + \frac{1}{n}\nu_0(a_n\hat{f}^-) O(4) + \frac{1}{n}\nu^0_0(a_n\hat{f}^-)\int_0^1\mathsf{span}\{\nu_s(B \hat{f}^- G^-)\}\,ds \\
    &=_{\text{\pref{lem:fact-s0},\,\pref{lem:rect-scaled}}} 0 + \frac{1}{n}\nu_0((f^-)^2) + \frac{1}{n}\nu_0(f^-\hat{f}^-) + O(4) =_{\text{\pref{lem:replicon-moments-from-mgf},~\pref{lem:hatf-moments}}} O(4)\, ,
\end{align*}
where one can prove $\nu_s(B\hat{f}G) = O(4)$ by bounding $\sup_{s\in[0,1]}\left(\nu_s(B\hat{f}G)-\nu_s(B)\nu_s(\hat{f}G)\right)' = O(4)$ by a strategy similar to the one invoking the absolute boundedness of $B$ and the moment bounds for $\hat{f}$ proved in~\pref{lem:hatf-moments}.
This finally gives,
\[
    \nu(a_n\hat f^-)=\beta^2\rho_4\,V_n+O(4)\,,
\]
and inserting into \eqref{eq:V-split-D17} allows us to obtain
\[
    V_n=\frac{4\rho_2}{(1-\beta^2\rho_4)n}+O(4).
\]
\ppart{Bounding non-similar terms in $W_n$}
Write
\[
\nu(\hat a_n\hat f^-)=\nu\big((\hat\varepsilon^1-\hat\varepsilon^2)(\hat\varepsilon^3-\hat\varepsilon^4)\,\hat f^-\big)=\nu(F_3).
\]
By~\pref{lem:D16-corrected},
\[
\nu(F_3)=\nu_0(F_3)+\nu_0^0(F_3)+O(4).
\]
As before, by~\pref{lem:fact-s0}, $\nu_0(F_3)=\nu_0(\hat{f}^-)\E_{Y^*}\left[D_{nn}(\eps^1-\eps^2)(\eps^3-\eps^4) \mid Y^*\right] = 0$. Furthermore, by~\pref{lem:rect-scaled},
\[
\nu_0^0(F_3)=\beta^2\rho_2\,\nu_0(\hat f^-\,f^-)=\beta^2\rho_2\,V_n+O(4)\, ,
\]
where we invoke the fact that $\nu(f^-\hat{f}^-) = \nu_0(f^-\hat{f}^-) + O(4)$ as shown previously in the proof. This gives
\[
\nu(\hat a_n\hat f^-)=\beta^2\rho_2\,V_n+O(4).
\]
Inserting the above equation into \eqref{eq:W-split-D17} along with the self-similarity term gives
\[
    W_n=\frac{4}{n}+\beta^2\rho_2\,V_n+O(4).
\]
Solving the first two relations yields
\[
U_n=\frac{4\rho_4}{n(1-\beta^2\rho_4)}+O(n^{-2}),
\qquad
V_n=\frac{4\rho_2}{n(1-\beta^2\rho_4)}+O(n^{-2}),
\]
and substituting it into the third finally gives
\[
W_n=\frac{4}{n}\Big(1+\frac{\beta^2\rho_2^2}{1-\beta^2\rho_4}\Big)+O(n^{-2}). \qedhere
\]
\end{proof}

We use the following slight modification of~\pref{prop:D17} to deal with ``doubly scaled'' rectangular sum limits, which appear in the replica-identity based simplification of $\mathsf{Tr}[PD^2P]$.

\begin{proposition}[Extension of~\pref{prop:D17} to $\nu(f\,f^{(2)})$]\label{prop:D18}
Recall $m_i=\langle\sigma_i\rangle$ and $D_{ii}=(1-m_i^2)^{-1}$.
Define the $D^2$--weighted rectangular sum
\[
a_i^{(2)}:=D_{ii}^2\,a_i,
\qquad
f^{(2)}:=\frac1n\sum_{i=1}^n a_i^{(2)},
\qquad
f^{(2),-}:=\frac1n\sum_{i\le n-1} a_i^{(2)}.
\]
Then
\begin{equation}\label{eq:D18-main}
W_n^{(2)} :=\nu\big(f\,f^{(2)}\big) =\frac{4}{n}+\beta^2\,U_n+O(4)\,.
\end{equation}
Consequently, using \pref{prop:D17},
\begin{equation}\label{eq:D18-solved}
W_n^{(2)}=\frac{4}{n(1-\beta^2\rho_4)}+O(n^{-2})\,.
\end{equation}
\end{proposition}

\begin{proof} The proof is completely similar in structure to that of~\pref{prop:D17} with the only differences being in the evaluation of $\nu_0(a^{(2)}_na_n)$, the explicit evaluation of $\nu^0_0(a^{(2)}_nf^-)$ via~\pref{lem:rect-doublyscaled}, and an invocation of the fact that terms of the form $\nu_s(D^2_{nn}a_n B_s f\Xi^-_1\Xi^-_2)$ are $O(4)$ (which follows directly by observing that the powers of $D_{nn}$ do not change the proof strategy of~\pref{lem:D16-corrected}, case-3).

Note that,
\begin{align*}
    &\nu(ff^{(2}) = \nu(a^{(2}_nf) = \nu(a^{(2}_nf^-) + \frac{1}{n}\nu(a^{(2)}_na_n) \\
    &= \left(\nu_0(a^{(2)}_nf^-) + \nu^0_0(a^{(2)}_nf^-) + \int_0^1 (1-s)\nu''_s(a^{(2)}_nf^-)\,ds\right) + \frac{1}{n}\left(\nu_0(a^{(2)}_na_n) + \nu^0_0(a^{(2)}_na_n) + \int_0^1(1-s)\nu''_s(a^{(2)}_na_n)\,ds\right)\,.
\end{align*}

Since $a_n^{(2)}a_n=D_{nn}^2 a_n^2$, we compute it at $s=0$. Under $\nu_0$, $m_n=\langle \varepsilon\rangle_{0}=\tanh(Y^*)$, and so
\[
    D^2_{nn}=\left(\frac1{1-m_n^2}\right)^2=\sech^{-4}(Y^*).
\]
Moreover,
\[
    a_n^2
    =(\varepsilon^1-\varepsilon^2)^2(\varepsilon^3-\varepsilon^4)^2
    =4(1-\varepsilon^1\varepsilon^2)(1-\varepsilon^3\varepsilon^4),
\]
so conditioning on $Y^*$ and using independence of replicas at $s=0$,
\[
    \big\langle a_n^2\big\rangle_{Y^*}=4\sech^4(Y^*).
\]
Therefore
\begin{equation}\label{eq:D18-diag-s0}
\nu_0(D_{nn}^2a_n^2)=4.
\end{equation}
Applying~\pref{lem:deriv} followed by~\pref{lem:fact-s0} immediately yields that $\nu^0_0(a^{(2)}_na_n) = \mathsf{span}\left\{\E_{Y^*}[D^{2}_nn a^2_n B] \nu_0(\Xi^-) \right\}$ which is immediately $O(1/n)$ by a combination of~\pref{lem:D-moments} and~\pref{lem:mean-bias-On-1}, along with the absolute boundedness of $B$.

Similarly, an application of~\pref{lem:deriv} twice yields that $\nu''_s(D^2_{nn}a^2_n) = \mathsf{span}\left\{\nu_s(D^{2}_{nn}a^2_n\Xi^-_1\Xi^-_2)\right\}$, which are $O(1/n)$ terms by a combination of~\pref{lem:D-moments},~\textbf{(H1)} and~\textbf{(H2)} in conjunction with a $(2,4,4)$-H\"older's inequality. 

For the terms in the first bracket, observe that, by~\pref{lem:fact-s0},
\[
    \nu_0(a^{(2)}_nf^-) = \nu_0(f^-)\E_{Y^*}\left[D^2_{nn}a_{nn}\right] = 0\,, 
\]
and by~\pref{lem:rect-doublyscaled},
\[
    \nu^0_0(a^{(2)}_nf^-) = \beta^2\nu_0((f^-)^2) = \beta^2\nu(f^2) + O(4)\, ,
\]
where the last equality by the same strategy in the proof of~\pref{prop:D17}. The last remaining thing to observe is that $\nu''_s(a^{(2)}_nf^-) = \mathsf{span}\left\{\nu_s(D^2_{nn}a_n B_s f\Xi^-_1\Xi^-_2)\right\}$, and each such term is $O(4)$ by exactly the same strategy as in the proof of~\pref{lem:D16-corrected}, case-3.

Putting the above together immediately gives,
\[
    W^{(2)}_n = \nu(ff^{(2)}) = \beta^2\nu(f^2) + \frac{4}{n} + O(4)\, ,
\]
and substituting $U_n = \nu(f^2)$ closes the proof of the main claim. Observe that substituting the expression for $U_n$ in~\pref{prop:D17} further gives
\[
W_n^{(2)}=\frac4n+\beta^2\frac{4\rho_4}{n(1-\beta^2\rho_4)}+O(n^{-2})
=\frac{4}{n(1-\beta^2\rho_4)}+O(n^{-2})\, . \qedhere
\]
\end{proof}

\paragraph{Bounds on ``diagonally-weighted'' tracial terms} Controlling certain remainder terms in $\calE_D$ and $\calE_A$ will requiring bounding terms where there are ``diagonally-weighted'' (by $D_{ii}$ or $D^2_{ii}$) Gibbs averages (in expectation over the disorder) of a rectangular sum $f$ multiplied by some bounded polynomial in a spin and some power of a bulk deviation (such as centered overlaps and magnetizations).

This is accomplished by extending the cavity estimates of the type in~\pref{prop:D15} and~\pref{prop:D15-full} with diagonal tilts, by explicitly evaluating the Gibbs derivative at $s = 0$ by using the fact that $D_{nn}$ can be exactly estimated in the cavity field, and using H\"older-type estimates on the ``bulk term'' to bound its second derivative.

\begin{lemma}[The $a_n$--term is $O(n^{-2})$ for linear bulk deviations]\label{lem:an-term-On-2}
Assume \textbf{(H1)}--\textbf{(H3)} and~\pref{lem:D-moments}.
Let $\Xi$ be any centered \emph{linear} bulk deviation of the form
\[
\Xi\in\mathsf{span}\{\,M_a-m^*,\;R_{ab}-q^*\,\}\,,
\]
where $a,b$ are among a fixed finite set of replica indices.
Define
\[
B:=D_{nn}\,(\varepsilon^3-\varepsilon^4)\,a_n = D_{nn}(\eps^1-\eps^2)(\eps^3-\eps^4)^2\,.
\]
Then
\begin{equation}\label{eq:an-term-bound}
\frac1n\,\nu(B\,\Xi)=O(n^{-2})\,.
\end{equation}
\end{lemma}

\begin{proof}
The proof will follow by the standard Taylor expansion, and control of the second derivative. We will isolate the ``cavity'' by removing the last spin, using exchangeability, and the exact factorization at $s = 0$ given by~\pref{lem:fact-s0}.

\ppart{Split $\Xi$ into a bulk part plus an explicit $O(1/n)$ last-spin correction}
If $\Xi=M_a-m^*$ then
\[
M_a=\frac1n\sum_{i<n}\rho_i^a+\frac1n\varepsilon^a=:M_a^-+\frac1n\varepsilon^a\,,
\]
which yields $\Xi= (M_a^- - m^*) + \frac{1}{n}\varepsilon^a$.
Similarly, for $\Xi=R_{ab}-q^*$,
\[
R_{ab}=\frac1n\sum_{i<n}\rho_i^a\rho_i^b+\frac1n\varepsilon^a\varepsilon^b=:R^-_{ab}+\frac1n\varepsilon^a\varepsilon^b\,,
\]
which gives $\Xi=(R^-_{ab} - q^*) + \frac{1}{n}\varepsilon^a\varepsilon^b$.
In either case we can write
\[
\Xi=\Xi^-+\frac1n\,\Pi(\varepsilon),
\]
where $\Xi^-$ depends only on the bulk spins $\rho$ and 
$\Pi(\varepsilon)$ is a bounded polynomial in finitely many last spins with $|\Pi|\le 1$.

Therefore
\[
\nu(B\Xi)=\nu(B\Xi^-)+\frac1n\nu(B\,\Pi(\varepsilon)).
\]
Since $|(\varepsilon^3-\varepsilon^4)a_n|\le 8$, $|\Pi|\le 1$ and~\pref{lem:D-moments} implies bounded moments for $D_{nn}$,
\[
|\nu(B\,\Pi)|\le 8\,|\nu(D_{nn})|\le 8\,\|D_{nn}\|_{L^1(\nu)}=O(1)
\]
Hence
\[
\frac1n\nu(B\,\Pi)=O(n^{-1}).
\]
Thus it remains to prove
\begin{equation}\label{eq:goal-bulk}
\nu(B\,\Xi^-)=O(n^{-1}).
\end{equation}

\ppart{Apply the Taylor expansion to $\nu(B\Xi^-)$ around $s=0$}
Let $F:=B\,\Xi^-$, viewed as a function of $k$ replicas. 
Using the two-term Taylor formula along the interpolation $(\nu_s)_{s\in[0,1]}$,
\[
    \nu(F)=\nu_0(F)+\nu_0^0(F)+\int_0^1(1-s)\,\nu_s''(F)\,ds.
\]

\ppart{Prove $\nu_0(F)=0$}
At $s=0$, by~\pref{lem:fact-s0}, the last spins $(\varepsilon^1,\dots,\varepsilon^k)$ are conditionally i.i.d.\ given $Y^*$.
Moreover, $D_{nn}$ is a function of $Y^*$ only at $s=0$, since $m_n=\tanh(Y^*)$ and $D_{nn}=\sech^{-2}(Y^*)$. At the same time, $\Xi^-$ depends only on the bulk spins $\rho^j$, which makes it independent of the last spins given $Y^*$.
Therefore,
\[
\nu_0(F)=\nu_0(\Xi^-)\,\E\Big[ D_{nn}\,\big\langle (\varepsilon^3-\varepsilon^4)a_n\big\rangle_{Y^*}\Big].
\]
But $a_n$ contains the factor $(\varepsilon^1-\varepsilon^2)$, and conditional on $Y^*$,
\[
\big\langle \varepsilon^1-\varepsilon^2\big\rangle_{Y^*}=\tanh(Y^*)-\tanh(Y^*)=0,
\]
hence $\langle(\varepsilon^3-\varepsilon^4)a_n\rangle_{Y^*}=0$ and thus $\nu_0(F)=0$.

\ppart{Bound $\nu_0^0(F)$ by $O(n^{-1})$ using $(2,2,4)$-H\"older}
Apply~\pref{lem:deriv} at $s=0$ to the observable $F$.
Every term in $\nu_0^0(F)$ is a finite linear combination of expectations of the form
\[
    \beta^2\,\nu_0\left(B'\,\Xi^-\,\Delta\right),
\]
where:
\begin{itemize}[itemsep=0.2em]
\item $B'$ is a bounded polynomial in finitely many last spins multiplied by $D_{nn}$, hence
$\|B'\|_{L^p(\nu_0)}=O(1)$ for all fixed $p$ by~\pref{lem:D-moments}, which comes from reasoning similar to that used to bound $\nu(B\Xi)$ above.
\item $\Delta$ is a centered bulk fluctuation of the form
$\Delta\in\{R^-_{uv}-q^*,\;M^-_u-m^*\}$ for some replica indices $u,v$.
\end{itemize}
By H\"older's inequality with exponents $(2,4,4)$,
\[
\big|\nu_0(B'\,\Xi^-\,\Delta)\big|
\le \|B'\|_{L^2(\nu_0)}\,\|\Xi^-\|_{L^4(\nu_0)}\,\|\Delta\|_{L^4(\nu_0)}.
\]
Using the concentration assumptions \eqref{eq:conc-assumptions} (and replica exchangeability),
\[
\|\Xi^-\|_{L^4(\nu_0)}=O(n^{-1/2}),\qquad \|\Delta\|_{L^4(\nu_0)}=O(n^{-1/2}),
\]
and $\|B'\|_{L^2(\nu_0)}=O(1)$, hence each summand is $O(n^{-1})$.
Since there are only finitely many summands,
\[
\nu_0^0(F)=O(n^{-1}).
\]

\ppart{Bounding $|\nu''_s(F)|$ in the Taylor remainder} Differentiating the terms in $\nu'_s(F)$ once more by applying~\pref{lem:deriv} yields that each term in $\nu_s''(F)$ is a finite linear combination of terms of the form
\[
\nu_s\left(B_s\,\Xi^-\,\Delta_1\,\Delta_2\right)\,,
\]
where $B_s$ is again a bounded last-spin polynomial times $D_{nn}$, and $\Delta_1,\Delta_2$ are centered bulk fluctuations.
Applying H\"older with exponents $(4,4,4,4)$ and invoking~\pref{lem:D-moments} gives
\[
\big|\nu_s(B_s\,\Xi^-\,\Delta_1\,\Delta_2)\big|
\le \|B_s\|_{L^4}\,\|\Xi^-\|_{L^4}\,\|\Delta_1\|_{L^4}\,\|\Delta_2\|_{L^4}
=O(1)\cdot O(n^{-1/2})\cdot O(n^{-1/2})\cdot O(n^{-1/2})
=O(n^{-3/2}),
\]
uniformly in $s\in[0,1]$. Integrating against $(1-s)$ preserves the order, so the remainder is $O(n^{-3/2})$.

Combining the previous bounds immediately yields
\[
\nu(B\Xi^-)=\nu_0(F)+\nu_0^0(F)+O(n^{-3/2})=0+O(n^{-1})+O(n^{-3/2})=O(n^{-1}),
\]
which is \eqref{eq:goal-bulk}. This immediately yields $\nu(B\Xi)=O(n^{-1})$, hence
\[
\frac1n\,\nu(B\Xi)=O(n^{-2})\,. \qedhere
\]
\end{proof}

\begin{remark}\label{rem:boundedness-handling}
    Note the following important property about the proof of~\pref{lem:an-term-On-2} -- The \emph{specific} power of $D_{nn}$ and the choice of $(\eps^3-\eps^4)$ is not necessary for the proof to work. The proof would go through for any $B$ of the form
    \[
        B := D^k_{nn}p_d(\eps^3,\eps^4,\dots,\eps^\ell)a_n\,,
    \]
    for $p_d$ being a degree-$d$ polynomial in replicas not including $1$ and $2$, and $k,\ell \in \N$. This is so because $\nu_0(B\Xi^-) = 0$ for such a choice (as the $(\eps^1-\eps^2)$ provides the ``centering'' term in $a_n$, conditional on $Y^*$) and the scaling of the remaining bounds on $\nu^0_0(B\Xi^-)$ and $\nu''_s(B\Xi^-)$ are also unchanged because of the application of H\"older's inequality. This also applies to the proof of~\pref{lem:D-oneleg-O(n)}.
\end{remark}

\begin{lemma}[Diagonal $D$--weighted one-leg term]\label{lem:D-oneleg-O(n)}
Assume \textbf{(H1)}--\textbf{(H3)} and~\pref{prop:D15}.  
Let $f_{1234}$ be the rectangular sum and let $\Xi$ be any centered bulk deviation
depending on finitely many replicas (e.g.\ $\Xi=M_a-m^*$ or $\Xi=R_{ab}-q^*$).
Then
\begin{equation}\label{eq:D-oneleg-target}
\frac{n^2}{4}\,
\Bigg|
\mathbb E\Bigg[\sum_{j=1}^n D_{jj}\,\Big\langle (\sigma_j^3-\sigma_j^4)\,\Xi\,f_{1234}\Big\rangle\Bigg]
\Bigg|
\;=\;O(n).
\end{equation}
More precisely,
\[
\mathbb E\Bigg[\sum_{j=1}^n D_{jj}\,\Big\langle (\sigma_j^3-\sigma_j^4)\,\Xi\,f_{1234}\Big\rangle\Bigg]
\;=\;O(n^{-1}).
\]
\end{lemma}

\begin{proof}
The proof proceeds by ``decoupling'' the final site after invoking exchangeability, writing down the second-order Taylor expansion, applying exact estimates at $s=0$, and bounding the derivative uniformly thereafter by appealing to~\pref{lem:fact-s0}.

\ppart{Reducing the sum over $j$ to a single cavity via exchangeability}
By exchangeability of sites under $(A,y_t)$,
\[
\mathbb E\Bigg[\sum_{j=1}^n D_{jj}\,\big\langle (\sigma_j^3-\sigma_j^4)\,\Xi\,f_{1234}\big\rangle\Bigg]
=
n\;\nu\Big(D_{nn}\,(\varepsilon^3-\varepsilon^4)\,\Xi\,f_{1234}\Big),
\]
where $\sigma=(\rho,\varepsilon)$ with $\varepsilon=\sigma_n$ and $\nu(\cdot)=\mathbb E\langle\cdot\rangle$.

Thus, it suffices to show that
\begin{equation}\label{eq:key-nu-bound}
\nu\Big(D_{nn}\,(\varepsilon^3-\varepsilon^4)\,\Xi\,f_{1234}\Big)=O(n^{-2}).
\end{equation}

\ppart{Separating the bulk from $f_{1234}$}
Write
\[
f_{1234}=f^-_{1234}+\frac1n a_n,
\]
where $f^-_{1234}$ depends only on the bulk spins $\rho$.
Then
\[
\nu\Big(D_{nn}(\varepsilon^3-\varepsilon^4)\Xi f_{1234}\Big)
=
\nu\Big(D_{nn}(\varepsilon^3-\varepsilon^4)\Xi f^-_{1234}\Big)
+\frac1n\,\nu\Big(D_{nn}(\varepsilon^3-\varepsilon^4)\Xi a_n\Big).
\]
The second term is $O(n^{-2})$ by~\pref{lem:an-term-On-2}.
So, the goal is to now show that
\begin{equation}\label{eq:main-term}
\nu\Big(D_{nn}(\varepsilon^3-\varepsilon^4)\Xi f^-_{1234}\Big)=O(n^{-2}).
\end{equation}

\ppart{Apply the Taylor expansion to $\nu(D_{nn}(\eps^3-\eps^4)\Xi f^-_{1234})$}
Let
\[
F:=D_{nn}(\varepsilon^3-\varepsilon^4)\,\Xi\,f^-_{1234}.
\]
Apply the two-term Taylor formula along the interpolation $(\nu_s)_{s\in[0,1]}$ to obtain
\[
\nu(F)=\nu_0(F)+\nu_0^0(F)+\int_0^1(1-s)\nu_s''(F)\,ds.
\]
We prove that each term is $O(n^{-2})$.

\ppart{Prove that $\nu_0(F) = 0$}
At $s=0$, by~\pref{lem:fact-s0}, conditional on $Y^*$ the last spins
$\varepsilon^1,\dots,\varepsilon^4$ are i.i.d.\ with mean $\tanh(Y^*)$.
Hence $\langle \varepsilon^3-\varepsilon^4\rangle_{Y^*}=0$, and since $D_{nn}$ has finite moments by~\pref{lem:D-moments},
\[
\nu_0(F)=\nu_0\Big(D_{nn}\,\Xi\,f^-_{1234}\,\big\langle \varepsilon^3-\varepsilon^4\big\rangle_{Y^*}\Big)=0.
\]

\ppart{Bound the first derivative term $\nu_0^0(F)$ as $O(n^{-2})$}
Apply~\pref{lem:deriv} at $s=0$ to the observable $F$ with $k=4$.
Because the last-spin polynomial contains $(\varepsilon^3-\varepsilon^4)$, the only nonzero contributions
in the SK overlap part are those where the pair $(\ell,\ell')$ involves replica $3$ or $4$.
Similarly, in the CW drift part, only terms with $\ell=3,4$ survive.
Consequently one obtains an identity of the form
\begin{equation}\label{eq:nu0prime-structure}
\nu_0^0(F)
=
\beta^2\sum_{u} c_u\;\nu_0\Big(f^-_{1234}\,\Xi\,\Delta_u\Big),
\end{equation}
where 
\begin{itemize}[itemsep=0.2em]
\item each coefficient $c_u$ is a last-spin expectation at $s=0$ of a bounded polynomial in $(\varepsilon^1,\dots,\varepsilon^6)$
multiplied by $D_{nn}$; hence $|c_u|\le C$ by the fact $c_u = \nu_0(B_u D_{nn})\le \left(\sup B_s\right) |\nu_0(D_{nn})|$ and~\pref{lem:D-moments}, and
\item each $\Delta_u$ is a \emph{centered linear bulk fluctuation}, namely one of
\[
R^-_{3,a}-R^-_{4,a},\qquad M^-_3-M^-_4,
\qquad\text{for some }a\in\{1,2,5\}.
\]
\end{itemize}
In particular, $\Xi \Delta_u$ is a quadratic bulk observable in centered deviations.
Therefore each term in \eqref{eq:nu0prime-structure} is exactly of the type covered by~\pref{prop:D15}, which immeidiately implies that
\[
\nu_0\Big(f^-_{1234}\,\Xi\,\Delta_u\Big)=O(n^{-2}).
\]
Adding the finitely many terms of this form gives
\[
\nu_0^0(F)=O(n^{-2}).
\]

\ppart{Apply H\"older to the Taylor remainder}
We now apply~\pref{lem:deriv} once more to $\nu_s'(F)$ to obtain that each term in $\nu_s''(F)$ is a finite linear combination of expectations of the form
\[
\nu_s\Big(B_s\;D_{nn}\;f^-_{1234}\;\Xi\;\Xi_1\;\Xi_2\Big),
\]
where $B_s$ is a bounded polynomial in finitely many last spins and $\Xi_1,\Xi_2$ are centered bulk fluctuations
of the form $(R^-_{ab}-q^*)$ or $(M^-_a-m^*)$.
By applying $(2,8,8,8,8)$-H\"older's, and using~\pref{lem:D-moments} to control the $D_{nn}$ moments, along with \eqref{eq:conc-assumptions} together with~\pref{lem:replicon-moments-from-mgf} to control $f^-_{1234}$ and the centered deviations, one obtains the following bound via a $(2,8,8,8,8)$-H\"older''s inequality
\[
\big|\nu_s''(F)\big|
\;\le\;
C\;\|D_{nn}\|_{L^2(\nu_s)}\;\|f^-_{1234}\|_{L^8(\nu_s)}\;\|\Xi\|_{L^8(\nu_s)}\;\|\Xi_1\|_{L^8(\nu_s)}\;\|\Xi_2\|_{L^8(\nu_s)}
\;=\;O(n^{-2})\,,
\]
and so
\[
    \int_0^1 (1-s)\nu''(s)ds \le |\nu''(s)| = O(n^{-2})\,.
\]
Combining the bounds above gives \eqref{eq:main-term}. This yields the final result of \eqref{eq:key-nu-bound} as
\[
\mathbb E\Bigg[\sum_{j=1}^n D_{jj}\,\big\langle (\sigma_j^3-\sigma_j^4)\,\Xi\,f_{1234}\big\rangle\Bigg]
=
n\cdot O(n^{-2})=O(n^{-1})\,. \qedhere
\]
\end{proof}

The final lemma we need to assert that the remainder terms in $\calE_A$ are $O(1)$ is a mild generalization of~\pref{prop:D15} to handle cubic moments across \emph{6} replicas (as opposed to $4$).

\begin{proposition}[Replicon-weighted cubic terms with \emph{general} quadratic weights]\label{prop:D15-full}
Fix a finite set of replica indices and $s\in[0,1]$.
Let
\[
\bar R^-_{ab}:=R^-_{ab}-q^*,\qquad \bar M^-_a:=M^-_a-m^*,
\qquad
f^-:=\bar R^-_{13}-\bar R^-_{14}-\bar R^-_{23}+\bar R^-_{24}.
\]
Let $\mathcal Q^-_{\mathrm{full}}$ be the linear span of \emph{all} quadratic bulk monomials built from these centered
deviations:
\[
\mathcal Q^-_{\mathrm{full}}
:=
\mathrm{span}\Big\{
\bar R^-_{ab}\bar R^-_{cd},\ \bar R^-_{ab}\bar M^-_c,\ \bar M^-_a\bar M^-_b
\Big\},
\]
where $a,b,c,d$ range over the fixed finite replica set. Then, for every $G^-\in \mathcal Q^-_{\mathrm{full}}$,
\begin{align*}
    \nu_s\big(f^-\,G^-\big)&= O(4) \\
    \nu_s\left(\hat{f}^-\,\hat{G}^-\right) &= O(4)\,.
\end{align*}
\end{proposition}
\begin{proof}
We follow the proof of~\pref{prop:D15}, and indicate the only point where the enlarged class
$\mathcal Q^-_{\mathrm{full}}$ matters.

\ppart{Bounding  $G^-$ in $L^p$}
Each $\bar R^-_{ab}$ and $\bar M^-_a$ is $O(n^{-1/2})$ in every $L^p(\nu_s)$ uniformly in $s$ by the moment
bounds implied by \textbf{(H1)}--\textbf{(H3)} and~\pref{lem:stability-s}.
This means that any quadratic monomial in these variables is $O(n^{-1})$ in every $L^p$,
and so for every fixed $p$,
\[
\|G^-\|_{L^p(\nu_s)} = O(n^{-1})
\qquad\text{uniformly in }s\in[0,1].
\]

\ppart{Reducing $\nu_s(f^-G^-)$ to a last-spin observable via exchangeability}
Using exchangeability, write $f=\frac1n\sum_{i=1}^n a_i$ with
$a_i=(\sigma_i^1-\sigma_i^2)(\sigma_i^3-\sigma_i^4)$ and $f=f^-+\frac1n a_n$.
Then site exchangeability gives
\[
\nu_s(f\,G^-)=\nu_s(a_n\,G^-),
\]
and since $\|G^-\|_{L^1}=O(n^{-1})$ and $|a_n|\le 4$,
\[
\nu_s(f\,G^-)-\nu_s(f^-\,G^-)=\frac1n\nu_s(a_nG^-)=O(n^{-2})=O(4).
\]
Thus
\[
\nu_s(f^-\,G^-)=\nu_s(a_nG^-)+O(4).
\]

\ppart{Taylor expanding $\nu_s(a_nG^-)$ around $s=0$}
Set $F:=(\varepsilon^1-\varepsilon^2)(\varepsilon^3-\varepsilon^4)\,G^-$.
Then $\nu_s(a_nG^-)=\nu_s(F)$ and
\[
\nu_s(F)=\nu_0(F)+\nu_0^0(F)+\int_0^s (s-u)\,\nu_u''(F)\,du.
\]
At $s=0$, factorization (\pref{lem:fact-s0}) gives $\nu_0(F)=0$.
By the unscaled rectangular cancellation lemma (\pref{lem:rect-unscaled}),
\[
\nu_0^0(F)=\beta^2\rho_4\,\nu_0(f^-G^-).
\]

\ppart{$O(4)$ bound for the Taylor remainder}
Differentiate twice using the derivative identity (\pref{lem:deriv}).
Exactly as in the proof of~\pref{prop:D15},
every term in $\nu_u''(F)$ is a finite linear combination of expectations of the form
\[
\nu_u\Big(B_u\,G^-\,\Xi_1\,\Xi_2\Big),
\]
where $B_u$ is a bounded polynomial in finitely many last spins, and
$\Xi_1,\Xi_2\in\{\bar R^-_{ab},\bar M^-_a\}$.
By Hölder and the $L^p$ size bounds ($\|G^-\|_{L^p}=O(n^{-1})$ and $\|\Xi_j\|_{L^p}=O(n^{-1/2})$),
each such term is $O(n^{-2})=O(4)$.
Thus, the second-order Taylor remainder is $O(4)$.

\ppart{Using linear stability at $\beta < 1$}
We have shown
\[
\nu_s(a_nG^-)=\beta^2\rho_4\,\nu_0(f^-G^-)+O(4).
\]
Replacing $\nu_0(f^-G^-)$ by $\nu_s(f^-G^-)$ costs only $O(4)$ by the same argument in~\pref{prop:D15}.
Therefore
\[
\nu_s(a_nG^-)=\beta^2\rho_4\,\nu_s(f^-G^-)+O(4).
\]
Combining the above with $\nu_s(f^-G^-)=\nu_s(a_nG^-)+O(4)$ immediately gives
\[
(1-\beta^2\rho_4)\,\nu_s(f^-G^-)=O(4),
\]
and since $\beta^2\rho_4<1$ for $\beta<1/2$, we conclude $\nu_s(f^-G^-)=O(4)$.
\end{proof}

\paragraph{Main cavity estimate via single-diagonal tilts} The final estimate that is key to obtaining the bound $\calE_A + \calE_D = O(1)$ is the observation that, up to $O(1)$ factors,
\[
    \E[c\Tr[PDP]] = \E[c\Tr[P]] + \E[c^2\Tr[P^2]]\,.
\]
This requires using the limits for $U_n$ and $V_n$ established in~\pref{prop:D17} in conjunction with a slight modification of the limit for $W^{(2)}_n$ that happens when there is conjugation by the elements of the diagonal matrix $D$.

\begin{lemma}[Singly scaled cancellation at $s=0$]\label{lem:rect-singlyscaled-cancel}
Let the last-spin variable $\varepsilon:=\sigma_n$ and bulk spins
$\sigma^-:=(\sigma_1,\dots,\sigma_{n-1})$.
Let $D_{nn}=(1-m_n^2)^{-1}$ and $P_{nn}=1-m_n^2$ computed under $\nu_s$.
Fix any bounded bulk observable $F^-=F^-(\sigma^{1,-},\dots,\sigma^{4,-})$ depending on finitely
many replicas and only on the first $n-1$ coordinates.
Then,
\begin{equation}\label{eq:singlyscaled-cancel}
\nu_0^0\!\Big(D_{nn}\,a_n\,f^-\,F^-\Big)
=
\beta^2\,\nu_0\!\Big(P_{nn}\,(f^-)^2\,F^-\Big)\,.
\end{equation}
\end{lemma}
\begin{proof}
We use the same differentiation formula for $\nu_s'(\cdot)$ as the one used throughout the cavity CLT proofs, and evaluate it at $s=0$. Note that
\[
    \nu^0_0(D_{nn}a_n f^-F^-) = \left(D'_{nn}\nu_0(a_n f^-F^-)\right) + \left(D_{nn}\nu'_0(a_nf^-F^-)\right) = \nu'_0(D_{nn}a_nf^-F^-)\,.
\]

At $s=0$, apply~\pref{lem:deriv} with the function of evaluation $D_{nn}a_nf^-F^-$. Using~\pref{lem:fact-s0}, since $\an{a_n\eps^\ell}_{Y^*} = 0$ for any $\ell \in \{1,2,3,4,5\}$ by the centering of $a_n$, all the (CW) terms in the derivative vanish. Similarly, since $\an{\eps^\ell\eps^5a_n}_{Y^*} = 0$, all but the first set of the (SK) terms are also $0$. So, we evaluate the remaining list of terms, that is
\[
    \nu^0_0\left(D_{nn}a_n f^-F^-\right) = \beta^2\sum_{1\le\ell<\ell'\le4} \nu_0\left(f^-F^-(R_{\ell,\ell'}-q^*)\right)\E_{Y^*}\left[D_{nn}(Y^*)\an{a_n\eps^\ell\eps^{\ell'}}_{Y^*}\right]\,,
\]
and independence implies the only contributing pairs are $(1,3),(1,4),(2,3),(2,4)$. It's easy to see that the following simplifications ensue
\[
    a_n\eps^\ell\eps^{\ell'} = (-1)^{\ell + \ell'}(1-\eps^1\eps^2)(1-\eps^3\eps^4)\, ,
\]
for $\ell \in \{1,2\}$ and $\ell' \in \{3,4\}$. Now,
\[
    (-1)^{\ell + \ell'}\E_{Y^*}\left[D_{nn}(1-\eps^1\eps^2)(1-\eps^3\eps^4)\right] = (-1)^{\ell + \ell'}\E_{Y^*}\left[P^{-1}_{nn}P^2_{nn}\right] = (-1)^{\ell + \ell'}\E_{Y^*}[P_{nn}]\,.
\]
Putting the above together with the definition of $f^- = (R^-_{1,3}-q^*) + (R^-_{2,4}-q^*) - (R^-_{2,3}-q^*) - (R^-_{1,4}-q^*)$ gives that
\[
    \nu^0_0(D_{nn}a_nf^-F^-) = \beta^2\nu_0\left(P_{nn}(f^-)^2 F^-\right)\,. \qedhere
\]
\end{proof}


\begin{lemma}[Cavity identity for the diagonally-titled trace of the squared-covariance]\label{lem:cTrPDP-identity}
Assume \textbf{(H1)--(H3)} and~\pref{prop:D15}.
Then
\begin{equation}\label{eq:cTrPDP-identity}
\E\big[c\,\Tr(PDP)\big]
=
\E\big[c\,\Tr(P)\big]
+
\E\big[c^2\,\Tr(P^2)\big]
+
O(1).
\end{equation}
\end{lemma}

\begin{proof}
By~\pref{lem:trace-to-replica} and the independence of replicas,
\[
    \E[c\,\Tr(PDP)] = \frac{n^2}{4}\,\nu\left(c\,f\,\hat f\right),
\qquad
\E[c^2\,\Tr(P^2)]
=
\frac{n^2}{4}\,\nu\big(c^2\,f^2\big).
\]
So, it suffices to prove
\begin{equation}\label{eq:target-nu}
\nu\big(c\,f\,\hat f\big)
=
\frac{4}{n^2}\,\E[c\,\Tr(P)]
+
\nu\big(c^2\,f^2\big)
+
O(4)\,.
\end{equation}

\ppart{Exchangeability reduction to the last site using the mean of $\hat f$}
Write $\hat f = \frac1n\sum_{i=1}^n \hat a_i$ with $\hat a_i:=D_{ii}a_i$.
Then, by site exchangeability under $\nu$,
\[
\nu(c f \hat f)
=
\nu\Big(c f\, \frac1n\sum_{i=1}^n \hat a_i\Big)
=
\nu\big(c f \hat a_n\big)
=
\nu\big(c\,D_{nn}\,a_n\,f\big).
\]
Decomposing $f=f^-+\frac1n a_n$ gives
\begin{equation}\label{eq:nu-split}
\nu(c f\hat f)
=
\nu\big(c\,D_{nn}\,a_n\,f^-\big)
+
\frac1n\,\nu\big(c\,D_{nn}\,a_n^2\big).
\end{equation}

\ppart{Identifying the diagonal term as \emph{exactly} $\frac{4}{n^2}\E[c\Tr(P)]$}
As in the proof of~\pref{lem:rect-singlyscaled-cancel}, because replicas are i.i.d.\ under the Gibbs measure, the following holds for any fixed instance of the disorder,
\[
\langle D_{nn}a_n^2\rangle = 4P_{nn}.
\]
Therefore, using the simplification above and exchangeability once more yields
\begin{equation}\label{eq:diag-exact}
\frac1n\,\nu(cD_{nn}a_n^2) = \frac{4}{n}\,\nu(cP_{nn}) = \frac{4}{n^2}\left(n\nu(cP_{nn})\right) = \frac{4}{n^2}\E\left[\sum_{i=1}^ncP_{ii}\right] = \frac{4}{n^2}\E\left[c\Tr[P]\right].
\end{equation}

\ppart{Cavity expansion for the off-diagonal term}
Setting $F:=c\,D_{nn}a_n f^-$ and applying the two-term Taylor formula in the cavity interpolation gives
\[
    \nu(F)=\nu_0(F)+\nu_0^0(F)+ \int_0^1(1-s)\nu''_s(F)\,ds.
\]
At $s=0$, conditional on the last-spin field $Y^*$ one has $\langle a_n\mid Y^*\rangle=0$, while $cD_{nn}$ depends only on $(Y^*,\sigma^-)$ but not on the replica labels,
and so $\nu_0(F)=0$. Furthermore, by~\pref{lem:rect-singlyscaled-cancel} with $F^- = c$\footnote{Note that $c = c^- + \frac{1}{n}\an{\eps^5}^2$, and one can observe that $\nu^0_0(\an{\eps^5}^2D_{nn}a_n f^-) = \beta^2\nu_0\left(\an{\eps^5}^2 P_{nn}(f^-)^2\right)$ by the same calculations as those in the proof of~\pref{lem:rect-singlyscaled-cancel} with the observation that $\an{\eps^5}^2$ is a constant with respect to the symmetries that result in cancellations after the application of~\pref{lem:fact-s0}.},
\begin{equation}\label{eq:off-cancel}
\nu^0_0(cD_{nn}a_nf^-) = \beta^2\,\nu_0\left(c\,P_{nn}\,(f^-)^2\right).
\end{equation}
To bound the second derivative, note that two applications of~\pref{lem:deriv} imply that $\nu''_s(F)$ is given by $\mathsf{span}\left\{\nu_s(B f^- \Xi^-_1\Xi^-_2)\right\}$ with $B_s = D_{nn}a_n\calP(\eps)$. Then, by the same strategy used to bound the deviation $\sup_{s \in [0,1]}\left(\nu_s(B_s f^-\Xi^-_1\Xi^-_2)-\nu_s(B)\nu_s(f^-\Xi^-_1\Xi^-_2)\right) = O(4)$ in the proof of~\pref{prop:D17}, we obtain the fact that $\nu''_s(F) = O(4)$.

\ppart{Replacing $\nu_0(cP_{nn}(f^-)^2)$ by $\nu(cP_{nn}f^2)$ up to $O(4)$}
Simplifying as in the proof of~\pref{prop:D17} gives
\begin{align*}
    &\nu(c P_{nn}f^2) -\nu_0(c P_{nn} (f^-)^2)= \left(\nu(c P_{nn}f^2) - \nu(c P_{nn} (f^-)^2)\right) + \left(\nu(c P_{nn} (f^-)^2) - \nu_0(c P_{nn}(f^-)^2)\right) \\
    &= \left(\frac{2}{n}\nu(c P_{nn}a_nf^-) + \frac{1}{n^2}\nu(c P_{nn} a^2_n)\right) + \int_0^1 \nu'_s(c P_{nn} (f^-)^2)\,ds \\
    &=_{\text{\pref{lem:deriv},}\,|cP_{nn}a^2)n| \le C} \frac{2}{n}\left(\nu_0(cP_{nn}a_nf^-) + \int_0^1 \mathsf{span}\left\{\nu_s(cP_{nn}a_nf^-\Xi^-)\right\}\,ds\right) + \int_0^1 \mathsf{span}\left\{\nu_s\left(cP_{nn}f^-\Xi^-_1\Xi^-_2\right)\right\}\,ds + O(4) \\
    &=_{\text{\pref{lem:fact-s0}\,+\,}(2,2)-\text{H\"older's}} \int_0^1 \mathsf{span}\left\{\nu_s\left(cP_{nn}f^-\Xi^-_1\Xi^-_2\right)\right\}\,ds + O(4)\,.
\end{align*}
Now, observe the following fact using a re-centering argument in conjunction with $(4,4,4,4)$-H\"older's and a final two-term Taylor expansion,
\begin{align*}
    \nu_s\left(cP_{nn}f^-\Xi^-_1\Xi^-_2\right) &=_{c = \beta^2\an{(1-R_{56})}} \beta^2(1+q^*)\nu_s(P_{nn}f^-\Xi^-_1\Xi^-_2) - \beta^2\nu_s\left(P_{nn}f^-\Xi^-_1\Xi^-_2(R_{56}-q^*)\right) \\
    &\le_{|P_{nn}| \le 1, \text{\pref{prop:D17}}} \beta^2(1+q^*)\nu_s(P_{nn}f^-\Xi^-_1\Xi^-_2) + O(4) \\
    &= O(4)\, ,
\end{align*}
where the last equality uses the fact that $\sup_{s \in [0,1]}\left(\nu_S(P_{nn}f^-\Xi^-_1\Xi^-_2)-\nu_s(P_{nn})\nu_s(f^-\Xi^-_1\Xi^-_2)\right) = O(4)$, which follows via uses of H\"older's inequality and~\pref{prop:D15} to control the derivative in a similar way to the way it was used in the proof of~\pref{prop:D17}. 
Inserting this into \eqref{eq:off-cancel} yields
\begin{equation}\label{eq:off-to-nu}
\nu(cD_{nn}a_nf^-)
=
\beta^2\,\nu\big(c\,P_{nn}\,f^2\big)
+O(4).
\end{equation}

\ppart{Invoking exact exchangeability via $\beta^2\nu(cP_{nn}f^2)=\nu(c^2 f^2)$}
Recall $c=\frac{\beta^2}{n}\Tr(P)$, hence
\[
    \beta^2\,\nu(cP_{nn}f^2) = \frac{\beta^4}{n}\,\nu\big(\Tr(P)\,P_{nn}\,f^2\big).
\]
We now expand $\Tr(P)=\sum_{i=1}^n P_{ii}$ and invoke exchangeability of the joint law of $\{P_{ii}\}_{i\le n}$ together with the site-symmetry of $f^2$ to obtain
\[
\nu\big(\Tr(P)\,P_{nn}\,f^2\big)
=
\sum_{i=1}^n \nu(P_{ii}P_{nn}f^2)
=
\nu(P_{nn}^2 f^2) + (n-1)\nu(P_{11}P_{22}f^2).
\]
On the other hand,
\[
\frac1n\,\nu\big(\Tr(P)^2 f^2\big)
=
\frac1n\sum_{i,j}\nu(P_{ii}P_{jj}f^2)
=
\nu(P_{nn}^2 f^2) + (n-1)\nu(P_{11}P_{22}f^2),
\]
so
\[
\nu\big(\Tr(P)\,P_{nn}\,f^2\big)=\frac1n\,\nu\big(\Tr(P)^2 f^2\big).
\]
Therefore
\[
\beta^2\,\nu(cP_{nn}f^2)
=
\frac{\beta^4}{n}\cdot \frac1n\,\nu(\Tr(P)^2 f^2)
=
\nu\Big(\frac{\beta^4}{n^2}\Tr(P)^2\,f^2\Big)
=
\nu(c^2 f^2).
\]
Combining with \eqref{eq:off-to-nu} gives
\begin{equation}\label{eq:off-final}
\nu(cD_{nn}a_nf^-)=\nu(c^2 f^2)+O(4).
\end{equation}

Inserting \eqref{eq:diag-exact} and \eqref{eq:off-final} into \eqref{eq:nu-split} gives
\[
\nu(c f\hat f)
= \nu(c^2 f^2) + \frac{4}{n^2}\E[c\Tr(P)] + O(4).
\]
This is exactly \eqref{eq:target-nu}. Multiplying by $n^2/4$ and using \pref{lem:trace-to-replica} yields \eqref{eq:cTrPDP-identity}.
\end{proof}

\subsection{Uniform spectral bound for $\mathbf{M}$} In~\pref{lem:mean-bias-On-1} it was critical to argue that $\beta^2\rho(M) < 1$ for every $\beta < 1/2$ to allow for the inversion of the operator $\Id - \beta^2 M$ to not cause a divergence when obtaining $O(1)$ fluctuations for the centered bulk observables. In the claim below, we provide a small lemma proving that the spectral radius of this operator arising from the cavity interpolation is indeed bounded as desired at high-temperature.

\begin{lemma}[Spectral-radius bound for $\beta^2\mathbf M$ when $\beta < 1/2$]\label{lem:beta2M-spr-bound}
Let $Y^*$ be a random variable that is almost-surely finite. Set $T:=\tanh(Y^*)\in[-1,1]$.  Define
\[
\mu_r:=\E[T^r],\qquad r\in\{1,2,3,4\},
\]
and
\[
\mathbf M:=
\begin{pmatrix}
1-4\mu_2+3\mu_4 & 2(\mu_1-\mu_3)\\
\mu_3-\mu_1 & 1-\mu_2
\end{pmatrix}.
\]
Then for every $\beta\in(0,1/2)$,
\[
\rho(\beta^2\mathbf M)\;<\;1.
\]

\end{lemma}

\begin{proof}
Factorization yields
\[
    1-4T^2+3T^4=(1-T^2)(1-3T^2),\qquad T-T^3=T(1-T^2)\,,
\]
for $T \in [-1,1]$. Hence, with $S:=1-T^2=\sech^2(Y^*)\in[0,1]$, we can rewrite the entries of $\mathbf M$ as
\[
    1-4\mu_2+3\mu_4=\E[(1-4T^2+3T^4)]=\E[S(1-3T^2)],
\]
\[
    \mu_1-\mu_3=\E[T-T^3]=\E[T(1-T^2)]=\E[TS],
\]
\[
    1-\mu_2=\E[1-T^2]=\E[S],\qquad \mu_3-\mu_1=-\E[TS].
\]
This allows $\mathbf M$ to be written as the expectation of a random matrix depending only on $T$ as
\begin{equation}\label{eq:M-as-expectation}
    \mathbf{M} := \E\left[\mathbf X(T)\right] := \E(1-T^2)\begin{pmatrix}1-3T^2 & 2T\\ -T & 1\end{pmatrix}.
\end{equation}
For any square matrix $\mathbf{M}$, the spectral radius is bounded by the operator norm as
\[
\rho(\mathbf{M})\le \opnorm{\mathbf{M}}.
\]
By the convexity of the operator norm,
\[
\opnorm{\mathbf M} =\opnorm{\E[\mathbf X(T)]}\le_{\text{Jensen's}} \E\big[\opnorm{\mathbf X(T)}\big].
\]
Upper bounding the operator norm by the Frobenius norm further gives
\[
\rho(\mathbf{M}) \le \opnorm{\mathbf M}\le \E\big[\|\mathbf X(T)\|_F\big]\le \sup_{|T|\le 1}\|\mathbf X(T)\|_F.
\]
To bound the supremum over $T$, note that for any $T\in[-1,1]$ and $u:= T^2$, a direct computation yields 
\begin{align*}
\|\mathbf X(T)\|_F^2
&=(1-T^2)^2\Big((1-3T^2)^2+5T^2+1\Big)\\
&=(1-u)^2\Big((1-3u)^2+5u+1\Big)\\
&=(1-u)^2(2-u+9u^2)\\
&=9u^4-19u^3+13u^2-5u+2\,.
\end{align*}
Differentiating $\norm{\mathbf X(T)}^2_F$ with respect to $u$ gives
\[
\frac{d}{du}\norm{\mathbf X(T)}^2_F=(u-1)\,(36u^2-21u+5).
\]
The equation $36u^2-21u+5$ has discriminant $(-21)^2-4\cdot 36\cdot 5=-279<0$, so it is strictly positive for all $u\in\R$. Therefore, using the fact that $(u-1) < 0$ implies that $\frac{d}{du}\norm{\mathbf X(T)}^2_F < 0$, and $\norm{\mathbf X(T)}^2_F$ is strictly decreasing on $[0,1]$. This implies that
\[
\sup_{T \in [-1,1]}\|\mathbf X(T)\|_F=\sqrt2.
\]
This immediately implies that
\[
\rho(\beta^2\mathbf M)\le \beta^2\|\mathbf M\|_2\le_{\norm{M}_2 \le \sqrt{2}} \beta^2\sqrt2 <_{\beta < 1/2} < 1. \qedhere
\]
\end{proof}

\subsection{The Nishimori line for the planted SK model}\label{app:nishimori} A critical fact used to ``convert'' overlap-concentration in the planted SK model to magnetization concentration is the fact that the external fields introduced by the SL process force $m^* = q^*$. This is a specific consequence of the fact that the SL process induces Gaussians with the same mean and variance.

\begin{lemma}[Nishimori condition for $\beta < 1/\sqrt{3}$]\label{lem:nishimori-condition}
Assume $0<\beta<1/\sqrt{3}$. Let $(m^*,q^*)$ solve the fixed-point equations
\begin{align*}
m^* &= \E_{G\sim\calN(0,1)}\Big[\tanh\big(t+\beta^2 m^*+\sqrt{\beta^2 q^*+t}\,G\big)\Big],\\
q^* &= \E_{G\sim\calN(0,1)}\Big[\tanh^2\big(t+\beta^2 m^*+\sqrt{\beta^2 q^*+t}\,G\big)\Big].
\end{align*}
Then $m^*=q^*$.
\end{lemma}

\begin{proof}
Set
\[
    \alpha:=t+\beta^2 q^* \ge 0,\qquad \delta:=\beta^2(m^*-q^*),
\]
so that $t+\beta^2 m^*=\alpha+\delta$ and $\sqrt{\beta^2 q^*+t}=\sqrt{\alpha}$.
Define
\[
h(x):=\tanh(x)-\tanh^2(x).
\]
Then, the fixed-point equations can be rewritten as
\[
m^*=\E\big[\tanh(\alpha+\delta+\sqrt{\alpha}\,G)\big],
\qquad
q^*=\E\big[\tanh^2(\alpha+\delta+\sqrt{\alpha}\,G)\big],
\]
and so
\begin{equation}
m^*-q^*=\E\Big[h\big(\alpha+\delta+\sqrt{\alpha}\,G\big)\Big].
\label{eq:diff-h}
\end{equation}

\ppart{An identity at $\delta=0$}
We claim that for every $\alpha\ge 0$,
\begin{equation}
\E\big[h(\alpha+\sqrt{\alpha}\,G)\big]=0,
\qquad\text{i.e.}\qquad
\E\big[\tanh(\alpha+\sqrt{\alpha}\,G)\big]=\E\big[\tanh^2(\alpha+\sqrt{\alpha}\,G)\big].
\label{eq:h-zero}
\end{equation}
If $\alpha=0$, the identity is trivial. Assume $\alpha>0$.

Let $\varphi_{m,\alpha}$ denote the density of $\calN(m,\alpha)$:
\[
\varphi_{m,\alpha}(x)=\frac{1}{\sqrt{2\pi\alpha}}\exp\!\Big(-\frac{(x-m)^2}{2\alpha}\Big).
\]
Then
\[
\E\big[\tanh(\alpha+\sqrt{\alpha}\,G)\big]=\int_{\R}\tanh(x)\,\varphi_{\alpha,\alpha}(x)\,dx,
\qquad
\E\big[\tanh^2(\alpha+\sqrt{\alpha}\,G)\big]=\int_{\R}\tanh^2(x)\,\varphi_{\alpha,\alpha}(x)\,dx.
\]
Let $\varphi_{-\alpha,\alpha}$ be the density of $\calN(-\alpha,\alpha)$. Note $\varphi_{-\alpha,\alpha}(x)=\varphi_{\alpha,\alpha}(-x)$, so
the sum $\varphi_{\alpha,\alpha}+\varphi_{-\alpha,\alpha}$ is even, while $\tanh$ is odd; hence
\[
\int_{\R}\tanh(x)\big(\varphi_{\alpha,\alpha}(x)+\varphi_{-\alpha,\alpha}(x)\big)\,dx=0,
\]
and therefore
\begin{equation}
\int_{\R}\tanh(x)\,\varphi_{\alpha,\alpha}(x)\,dx
=\frac12\int_{\R}\tanh(x)\big(\varphi_{\alpha,\alpha}(x)-\varphi_{-\alpha,\alpha}(x)\big)\,dx.
\label{eq:tanh-diff}
\end{equation}
Similarly, $\tanh^2$ is even and $\int \tanh^2 \varphi_{\alpha,\alpha}=\int \tanh^2 \varphi_{-\alpha,\alpha}$ (change $x\mapsto -x$), so
\begin{equation}
\int_{\R}\tanh^2(x)\,\varphi_{\alpha,\alpha}(x)\,dx
=\frac12\int_{\R}\tanh^2(x)\big(\varphi_{\alpha,\alpha}(x)+\varphi_{-\alpha,\alpha}(x)\big)\,dx.
\label{eq:tanh2-sum}
\end{equation}

A direct computation gives the likelihood ratio
\[
\frac{\varphi_{\alpha,\alpha}(x)}{\varphi_{-\alpha,\alpha}(x)}=\exp(2x),
\qquad\text{i.e.}\qquad
\varphi_{\alpha,\alpha}(x)=e^{2x}\varphi_{-\alpha,\alpha}(x).
\]
Using \eqref{eq:tanh-diff} and \eqref{eq:tanh2-sum} and rewriting everything in terms of $\varphi_{-\alpha,\alpha}$ yields
\begin{align*}
\E\big[\tanh(\alpha+\sqrt{\alpha}\,G)\big]
&=\frac12\int_{\R}\tanh(x)\,(e^{2x}-1)\,\varphi_{-\alpha,\alpha}(x)\,dx,\\
\E\big[\tanh^2(\alpha+\sqrt{\alpha}\,G)\big]
&=\frac12\int_{\R}\tanh^2(x)\,(e^{2x}+1)\,\varphi_{-\alpha,\alpha}(x)\,dx.
\end{align*}
But the integrands coincide pointwise because $
\tanh(x)=(e^{2x}-1)/(e^{2x}+1)$, implying
\[
\tanh(x)(e^{2x}-1)=\tanh^2(x)(e^{2x}+1).
\]
Hence the two expectations are equal, proving \eqref{eq:h-zero}.

\ppart{Lipschitz bound forces $\delta=0$ for $\beta<\frac{1}{\sqrt{3}}$}
Differentiate
\[
h'(x)=\sech^2(x)\bigl(1-2\tanh(x)\bigr),
\]
so $|h'(x)|\le 1\cdot 3=3$ for all $x$, and thus $h$ is globally $3$-Lipschitz:
\[
|h(u)-h(v)|\le 3|u-v|\quad\forall u,v\in\R.
\]
Using \eqref{eq:diff-h} and \eqref{eq:h-zero},
\begin{align*}
|m^*-q^*|
&=\Big|\E\big[h(\alpha+\delta+\sqrt{\alpha}G)\big]-\E\big[h(\alpha+\sqrt{\alpha}G)\big]\Big|\\
&\le \E\Big|h(\alpha+\delta+\sqrt{\alpha}G)-h(\alpha+\sqrt{\alpha}G)\Big|\\
&\le 3|\delta|
=3\beta^2|m^*-q^*|.
\end{align*}
Since $0<\beta<\tfrac{1}{\sqrt{3}}$ implies $3\beta^2<1$, the inequality forces $|m^*-q^*|=0$, i.e.\ $m^*=q^*$.
\end{proof}

\section{Properties of the algorithmic stochastic process}
\label{sec:alg-properties}

Recall the algorithmic SDE
\begin{align*}\repeatequation{e:ASL-TAP}\end{align*}
simulated by the main sampler (\pref{alg:main}).
Its coefficients are given by the functions
\begin{equation}
\label{eq:DahQ}
\begin{aligned}
    D(m) &:= (\Id_n - \diag(m)^2)^{-1},
    \\
    a(m) &:= \beta^2\tr_n(D(m)^{-1}),
    \\
    \hat{Q}(m) &:= \left(a(m)\Id_n - \beta A + D(m) - \frac{2\beta^2}{n} mm^{\sT}\right)^{-1}.
\end{aligned}
\end{equation}

These coefficients of the SDE for $\hat{m}_t$ (and consequently the SDE itself) turn out to be remarkably well-behaved under certain high-probability events on the GOE matrix $A$.
This is mainly due to the boundedness and regularity of the resolvent $(z,M) \to (z - M)^{-1}$, which \pref{e:ASL-TAP} inherits from its similarity with \pref{e:PHD}, the resolvent-based SDE whose diffusion coefficient is simply the Hessian of the TAP free energy $d\hat{m}_t = \hat{Q}(\hat{m}_t)dB_t$.

Because ASL-TAP is closely related to the Gibbs measure and PHD is closely related to the TAP free energy, the closeness of these two processes is central to the functioning of \pref{alg:main}.

We introduce the high-probability events in \pref{sec:alg-events}, including one simple spectral norm bound and another more sophisticated ``control of the diagonal'' event from free probability theory.
In \pref{sec:alg-basic-bounds-lipschitz}, we collect boundedness and Lipschitzness properties for the resolvent $\hat{Q}(m)$ implied by the simple spectral norm event.
The event controlling the diagonal is stated for the resolvent $\hat{Q}(m)$ without the rank-1 term, so \pref{sec:alg-diagonal-error-rank-1} shows that the addition of the rank-1 term doesn't change the diagonal too much.
\pref{sec:alg-diagonal-control} then gives statements for control of the diagonal that directly control quantities appearing in \pref{e:ASL-TAP}.
\pref{sec:alg-phd-asl-tap-closeness} shows that the SDEs \pref{e:ASL-TAP} and \pref{e:PHD} are close to each other.
Some consequences for the \pref{e:ASL-TAP} trajectories are detailed in \pref{sec:alg-trajectories}, including that the trajectory almost surely stays within the open solid hypercube $(-1,1)^n$ and that Lipschitz functionals of the trajectories obey sub-Gaussian concentration.
Finally, \pref{sec:alg-desiderata} details how these properties of \pref{e:ASL-TAP} contribute to the success of the sampling algorithm \pref{alg:main}.

\subsection{High-probability events}
\label{sec:alg-events}

There are two high-probability events on the GOE matrix $A$ that underpin most of the well-behaved properties of the algorithmic processes and their SDE coefficients.

The first is a simple bound on the spectral norm of $A$. Let $\gamma := (1-2\beta)/2$ so that $\beta = 1/2 - \gamma$, and define the event
\begin{equation}
    \label{eq:event-opnorm}
\Omega^{(A)}_{\gamma} = \left\{\vphantom{\big|}\beta\opnorm{A} \le 1-\gamma\right\}.
\end{equation}
By \cite[Proposition 1]{aubrun2004sharp}, when $\beta < 1/2$, we have $\beta\opnorm{A} \le \beta(2 + \frac{1-2\beta}{2\beta}) = 1 - \gamma$ with probability at least $1 - e^{-C_{\beta}n}$ for some constant $C_{\beta}$.

The second event asserts that the sum $\beta A + D$ for diagonal $D$ hews very closely to the (uniformly) freely independent case, which includes having every diagonal entry of the squared resolvents $(z - \beta A + D)^{-2}$ essentially determined by $z$ and $D$.
For some constant $C_{\mathrm{free}}$,
\begin{equation}
    \label{eq:event-free}
     \Omega^{(\mathrm{free})} = \left\{\forall D\in \calD_n([1,\infty)) \mathrel{.} \schnorm{E_{\calD_n}\left[\left(\beta^2\tr_n(D^{-1}) - \beta A + D\right)^{-2}\right]D^2 - \frac{\Id_n}{1-\beta^2\tr_n(D^{-2})}} \le \frac{C_{\mathrm{free}}\beta}{\gamma^4\sqrt{n}} \right\}.
\end{equation}
As shown in \pref{cor:diagonal-Q2-controlled-in-l2}, this occurs with probability at least $1 - e^{-O((\log n)^{\xi})}$ with some $\xi > 1$.

\subsection{Basic resolvent and Lipschitzness bounds}
\label{sec:alg-basic-bounds-lipschitz}

We start with bounds implied by the simple spectral norm condition $\Omega^{(A)}_{\gamma}$.

\begin{lemma}[Basic resolvent bounds assuming $\Omega^{(A)}_{\gamma}$]
    \label{lem:alg-resolvent-bounds}
    Suppose that $\Omega^{(A)}_{\gamma}$ holds.
    Let $D(m), a(m), \hat{Q}(m)$ be as in \pref{eq:DahQ}.
    Let also
    \[
    \bar{Q}(m) := \left(a(m)\Id_n - \beta A + D(m)\right)^{-1}
    \]
    be the version of $\hat{Q}$ without the rank-1 term.
    Then for all $m \in (-1,1)^n$, we have the following bounds:
    \begin{enumerate}[label=(\alph*), ref=\thetheorem(\alph*)]
    \item \label{lem:alg-resolvent-bounds-Q}
    \[ \opnorm{\bar{Q}(m)} \le \gamma^{-1}, \qquad \opnorm{\hat{Q}(m)} \le \gamma^{-1}. \]
    \item \label{lem:alg-resolvent-bounds-DQ}
    \[ \opnorm{D(m)\bar{Q}(m)} \le 2\gamma^{-1}, \qquad \opnorm{D(m)\hat{Q}(m)} \le 2\gamma^{-1}. \]
    \item \label{lem:alg-resolvent-bounds-DQ2}
    \[ \opnorm{D(m)\bar{Q}^2(m)} \le 2\gamma^{-2}, \qquad \opnorm{D(m)\hat{Q}^2(m)} \le 2\gamma^{-2}. \]
    \item \label{lem:alg-resolvent-bounds-DQ2D}
    \[ \opnorm{D(m)\hat{Q}^2(m)D(m)} \le 4\gamma^{-2}. \]
    \end{enumerate}
\end{lemma}
\begin{proof}
    The $\bar{Q}$ parts of \pref{lem:alg-resolvent-bounds-Q}, \pref{lem:alg-resolvent-bounds-DQ}, and \pref{lem:alg-resolvent-bounds-DQ2} follow directly from \pref{lem:resolvent-bounds} with 
    \[W \;\gets\; \Id_n - \diag(m)^2\qquad\text{and}\qquad M \;\gets\; a(m)\Id_n - \beta A,\]
    so that
    \[\min\Spec(|X|) \;\ge\; a(m) - (1-\gamma) + 1 \;\ge\; \gamma \qquad\text{and}\qquad \norm{M} \;\le\; 2-\gamma.\qquad\]

    For the $\hat{Q}$ parts and for \pref{lem:alg-resolvent-bounds-DQ2D}, first note that $mm^{\sT}/n \preceq \diag(m)^2$, which follows by taking the inequality $\mathbf{1}\mathbf{1}^{\sT}/n \preceq \Id_n$ and multiplying by $\diag(m)$ on both the left and right.
    Also we have $D(m) \succeq \Id_n + \diag(m)^2 $ by convexity of $x \mapsto (1 - x)^{-1}$.
    We can therefore apply \pref{lem:resolvent-bounds} with
    \[W \;\gets\; \Id_n - \diag(m)^2 \qquad\text{and}\qquad M \;\gets\; a(m)\Id_n - \beta A - \frac{2\beta^2}{n} mm^{\sT},\]
    so that \[X \;\;\gets\;\; a(m)\Id_n - \beta A + D(m) - \frac{2\beta^2}{n}mm^{\sT}  \;\;\succeq\;\; 0 - (1-\gamma)\Id_n +  (\Id_n + \diag(m)^2) - 2\beta^2\diag(m)^2 \;\;\succeq\;\; \gamma \Id_n.\]
    So $\min\Spec(|X|) \ge \gamma$ and $\opnorm{M} \le 2-\gamma$.
\end{proof}

\begin{lemma}[Basic resolvent Lipschitzness assuming $\Omega^{(A)}_{\gamma}$]
    \label{lem:alg-resolvent-lipschitz}
    In the same setting as \pref{lem:alg-resolvent-bounds}, over the domain $m \in (-1,1)^m$, we have the following Lipschitz constants with respect to $\schnorm{\cdot}$:
    \begin{enumerate}[label=(\alph*), ref=\thetheorem(\alph*), itemsep=0.1em]
    \item \label{lem:alg-resolvent-lipschitz-Q}
    \[ \norm{\hat{Q}}_{\Lip} \le \frac{\gamma^{-2}}{\sqrt{n}}\left(2 + 8\beta^2 + \frac{16\beta^2}{\sqrt{n}}\right).\]
    \item \label{lem:alg-resolvent-lipschitz-Q2}
    \[ \norm{\hat{Q}^2}_{\Lip} \le \frac{\gamma^{-3}}{\sqrt{n}}\left(4 + 16\beta^2 + \frac{32\beta^2}{\sqrt{n}}\right).\]
    \item \label{lem:alg-resolvent-lipschitz-DQ}
    \[ \norm{D\hat{Q}}_{\Lip} \le \frac{\gamma^{-2}}{\sqrt{n}}\left(4 + 16\beta^2 + \frac{32\beta^2}{\sqrt{n}}\right). \]
    \item \label{lem:alg-resolvent-lipschitz-DQQD}
    \[ \norm{D\hat{Q}^2D}_{\Lip} \le \frac{\gamma^{-3}}{\sqrt{n}}\left(16 + 64\beta^2 + \frac{128\beta^2}{\sqrt{n}}\right). \]
    \end{enumerate}
\end{lemma}
\begin{proof}
    We apply \pref{lem:resolvent-frechet} with
    \[W \;\gets\; \Id_n - \diag(m)^2 \qquad\text{and}\qquad M \;\gets\; a\Id_n - \beta A - \frac{2\beta^2}{n} mm^{\sT},\]
    following the proof of \pref{lem:alg-resolvent-bounds}.
    Then $da = -2\beta^2 m^{\sT} dm\Id_n/n$, so
    \[dW \,=\, -2\diag(m\odot dm) \qquad\text{and}\qquad dM \,=\, -\frac{2\beta^2}{n} \left(m^{\sT} dm\Id_n + (dm)m^{\sT}  + m(dm)^{\sT} \right),\]
    and we have
    \[\frenorm{dW} \;\le\; \frac{2}{\sqrt{n}}\norm{dm} \qquad\text{and}\qquad \frenorm{dM} \;\le\; \frac{2\beta^2}{n}\left(\sqrt{n} + 2\right)\frenorm{dm}.\]
    Then by \pref{lem:resolvent-frechet-dR}, \pref{lem:resolvent-frechet-dRk}, \pref{lem:resolvent-frechet-dDR}, and \pref{lem:resolvent-frechet-dDRkD},
    \begin{align*}
        \frenorm{d\hat{Q}} &{}\le \left(\frac{2\gamma^{-2}}{\sqrt{n}} + \frac{8\beta^2\gamma^{-2}}{\sqrt{n}}\left(1 + \frac{2}{\sqrt{n}}\right)\right)\frenorm{dm},
        \\
        \frenorm{d\hat{Q}^2} &{}\le \left(\frac{4\gamma^{-3}}{\sqrt{n}} + \frac{16\beta^2\gamma^{-3}}{\sqrt{n}}\left(1 + \frac{2}{\sqrt{n}}\right)\right)\frenorm{dm},
        \\
        \frenorm{d(D\hat{Q})} &{}\le \left(\frac{4\gamma^{-2}}{\sqrt{n}} + \frac{16\beta^2\gamma^{-2}}{\sqrt{n}}\left(1 + \frac{2}{\sqrt{n}}\right)\right)\frenorm{dm},
        \\
        \frenorm{d(D\hat{Q}^2D)} &{}\le \left(\frac{16\gamma^{-3}}{\sqrt{n}} + \frac{16(4\gamma^{-1} - 2)\beta^2\gamma^{-2}}{\sqrt{n}}\left(1 + \frac{2}{\sqrt{n}}\right)\right)\frenorm{dm}.
    \end{align*}
     We conclude Lipschitzness by \pref{fact:frechet-lipschitz} since $\frenorm{dm}=1$ when $m$ is the independent variable.
\end{proof}

Finally, we record here the strict convexity of $\mathcal{F}_{\mathrm{TAP}}$, equivalent to the positive-definiteness of $\hat{Q}$---which is a simple spectral bound assuming $\Omega^{(A)}_{\gamma}$---although it will not be used until \pref{cor:alg-desiderata-d2}.
\begin{lemma}[Strict convexity of $\mathcal{F}_{\mathrm{TAP}}$ assuming $\Omega^{(A)}_{\gamma}$]
    \label{lem:tap-convexity}
    Suppose that $\Omega^{(A)}_{\gamma}$ holds.
    Let $D(m), a(m), \hat{Q}(m)$ be as in \pref{eq:DahQ}.
    Then
    \[  \hat{Q}(m) \succ \frac{1}{4} D(m)^{-1}. \]
\end{lemma}
\begin{proof}
    Since $D(m) \succ 0$, the lemma statement is equivalent to 
    \[ Y := D(m)^{1/2}\hat{Q}(m)D(m)^{1/2} \succ \frac{1}{4}\Id_n.\]
    The final line of \pref{lem:alg-resolvent-bounds} implies $\hat{Q}(m) \succ 0$ since $X = \hat{Q}(m)^{-1}$.
    Therefore, $Y$ is positive-definite and 
    \[Y \succeq \opnorm{Y^{-1}}^{-1}\Id_n,\]
    so it is sufficient to show that $\opnorm{Y^{-1}} < 4$. We expand
    \[ Y^{-1} = a(m)D(m)^{-1} - \beta D(m)^{-1/2}AD(m)^{-1/2} + \Id_n - \frac{2\beta}{n} D(m)^{-1/2}mm^{\sT} D(m)^{-1/2}.\]
    By triangle inequality and the fact that $\opnorm{D(m)^{-1/2}} \le 1$,
    \[ \opnorm{Y^{-1}} \le a(m) + \opnorm{\beta A} + 1 + \frac{2\beta}{n} \lpnorm{m}^2 \le 1 + (1 - \gamma ) + 1 + 2\beta < 4. \qedhere\]
\end{proof}

\subsection{Diagonal error from rank-1 term of TAP Hessian}
\label{sec:alg-diagonal-error-rank-1}
The event $\Omega^{(\mathrm{free})}$ in \pref{eq:event-free} captures the control over the diagonal that we get from the free probability theory argument in \pref{sec:free-interpolation}.
But this control is on the inverse $\bar{Q}$ of the TAP Hessian \emph{without the rank-1 term}, whereas the rank-1 term is present throughout the actual algorithmic process.

We need just one fact to capture the discrepancy between these two resolvents.

\begin{lemma}[Small diagonal error from rank-1 terms]
    \label{lem:alg-resolvent-bounds-DQ2-hQ2D}
    In the same setting as \pref{lem:alg-resolvent-bounds}, for all $m \in (-1,1)^n$,
    the diagonal error between $\hat{Q}^2$ and $\bar{Q}^2$ satisfies
    \[ \schnorm{E_{\calD_n}\left[\hat{Q}^2(m) - \bar{Q}^2(m)\right]D^2(m)} \le \frac{16\beta^2\gamma^{-3}}{\sqrt{n}}. \]
\end{lemma}
\begin{proof}
    We factor with the second-order resolvent identity $X^{-2} - Y^{-2} \equiv X^{-2}(Y-X)Y^{-1} + X^{-1}(Y-X)Y^{-2}$ for all invertible matrices $X, Y$. 
    Since $\bar{Q}{}^{-1} - \hat{Q}^{-1} =  2\beta^2 mm^{\sT} /n$, we have
    \begin{align*}
    \hat{Q}^2 - \bar{Q}^2
    {}&\;=\; \hat{Q}^2(\bar{Q}{}^{-1}-\hat{Q}^{-1})\bar{Q} + \hat{Q}(\bar{Q}{}^{-1}-\hat{Q}^{-1})\bar{Q}{}^2
    \\{}&\;=\; \frac{2\beta^2}{n}\left[\hat{Q}^2mm^{\sT} \bar{Q} + \hat{Q}mm^{\sT} \bar{Q}{}^2\right]
    .
    \end{align*}
    Then, applying the identity $E_{\calD_n}[uv^{\sT} ] \equiv \diag(u)\diag(v)$,
    \begin{equation}
    \label{eq:alg-EDbQbQD-EDhQhQD-factorization}
    E_{\calD_n}\left[\hat{Q}^2 - \bar{Q}{}^2\right]D^2 = \frac{2\beta^2}{n}\left[\diag(D\hat{Q}^2m)\diag(D\bar{Q}m) + \diag(D\hat{Q}m)\diag(D\bar{Q}{}^2m)\right].
    \end{equation}
    We now give a quick inequality for generic $x,y \in \R^n$.
    By AM-GM and the $\ell_4$-to-$\ell_2$ norm inequality $\norm{\cdot}_4 \le \norm{\cdot}_2$,
    \begin{align*}
        \schnorm{\diag(x)\diag(y)}^2 = \frac{1}{n}\sum_{i=1}^n (x_iy_i)^2 \le \frac{1}{2n}\sum_j (x_i^4 + y_i^4) = \frac{1}{2n}\left(\lpnorm[4]{x}^4 + \lpnorm[4]{y}^4\right) \le \frac{1}{2n}\left(\lpnorm{x}^4 + \lpnorm{y}^4\right) \,.
    \end{align*}
    Applying this inequality to the two terms in \pref{eq:alg-EDbQbQD-EDhQhQD-factorization} and then using the basic resolvent bounds \pref{lem:alg-resolvent-bounds-DQ} and \pref{lem:alg-resolvent-bounds-DQ2},
    \begin{align*}
    \schnorm{E_{\calD_n}\left[\hat{Q}^2 - \bar{Q}{}^2\right]D^2}^2
    &\le \frac{4\beta^4}{n^2}\schnorm{\diag(D\hat{Q}^2m)\diag(D\bar{Q}m) + \diag(D\hat{Q}m)\diag(D\bar{Q}{}^2m)}^2
    \\&\le \frac{4\beta^4}{n^3}\left(\lpnorm{\gamma^{1/2}D\hat{Q}^2m}^4 + \lpnorm{\gamma^{-1/2}D\bar{Q}m}^4 + \lpnorm{\gamma^{-1/2}D\hat{Q}m}^4 + \lpnorm{\gamma^{1/2}D\bar{Q}{}^2m}^4\right)
    \\&\le \frac{4\beta^4}{n^3}\left(16\gamma^{-6} + 16\gamma^{-6} + 16\gamma^{-6} + 16\gamma^{-6} \right)\lpnorm{m}^4.
    \qedhere
    \end{align*}
\end{proof}

\subsection{Control of the diagonal}
\label{sec:alg-diagonal-control}
The diagonal entries of the instantaneous covariance matrix $(dm_t)(dm_t)^{\sT}  = \hat{Q}^2dt$ hold a special position in the analysis because they determine the instantaneous variance of each coordinate: $\E[dm_{t,i}^2] = (\hat Q^2)_{i,i}dt$.

In particular, the condition $\Omega^{(\mathrm{free})}$ implies that the diagonal entries of $\hat{Q}^2$ are close to a scalar multiple of the diagonal entries of $D^{-2}$, a fact that will be succinctly captured by bounds on the following quantity:
\begin{equation}
\label{eq:delta-diag}
    \delta_{\diag}(m) := E_{\calD_n}[\hat{Q}(m)^2]\,D(m)^2 - (1 + \beta^2\tr_n(\hat{Q}(m)^2))\Id_n
\end{equation}

First, we record operator norm and Lipschitzness bounds on $\delta_{\diag}$ that hold even under $\Omega^{(A)}_{\gamma} \setminus \Omega^{(\mathrm{free})}$.
\begin{lemma}[Diagonal properties assuming only $\Omega^{(A)}_{\gamma}$]
\label{lem:alg-diag-basic}
Suppose that $\Omega^{(A)}_{\gamma}$ holds.
    Let $D(m), a(m), \hat{Q}(m)$ be as in \pref{eq:DahQ} and let $\delta_{\diag}(m)$ be as in \pref{eq:delta-diag}. Then with $m \in (-1,1)^n$,
    \begin{enumerate}[label=(\alph*), ref=\thetheorem(\alph*)]
    \item \label{lem:alg-resolvent-bounds-delta}
    \[ \opnorm{\delta_{\diag}(m)} \le 4\gamma^{-2}. \]
    \item \label{lem:alg-resolvent-lipschitz-delta}
    \[ \norm{\delta_{\diag}}_{\Lip} \le \frac{\gamma^{-3}}{\sqrt{n}}(4 + \beta^2)\left(4 + 16\beta^2 + \frac{32\beta^2}{\sqrt{n}}\right),\]
    where the Lipschitz constant is with respect to $\schnorm{\cdot}$.
    \end{enumerate}
\end{lemma}
\begin{proof}
    \ppart{Part \ref{lem:alg-resolvent-bounds-delta}} $E_{\calD_n}[\cdot]$ and $\tr_n(\cdot)\Id_n$ are both contractions for $\opnorm{\cdot}$, and $\delta_{\diag}$ is a difference of two PSD matrices, so its operator norm is bounded by $\max(\opnorm{D\hat{Q}^2D},\beta^2\opnorm{\hat{Q}^2})$, which is given by \pref{lem:alg-resolvent-bounds-DQ2D} and \pref{lem:alg-resolvent-bounds-Q}.

    \ppart{Part \ref{lem:alg-resolvent-lipschitz-delta}} Again $E_{\calD_n}[\cdot]$ and $\tr_n(\cdot)\Id_n$ are both contractions for $\schnorm{\cdot}$, so it suffices by triangle inequality to add the basic Lipschitz bounds $\lipnorm{D\hat{Q}^2D}$ and $\lipnorm{\hat{Q}^2}$ established in \pref{lem:alg-resolvent-lipschitz-DQQD} and \pref{lem:alg-resolvent-lipschitz-Q2}.
\end{proof}

Then a strong uniform bound on the 2-norm of $\diag_{\delta}$ follows straightforwardly from $\Omega^{(\mathrm{free})}$ and the rank-1 diagonal error bound in \pref{sec:alg-diagonal-error-rank-1}.
\begin{lemma}[Diagonal control assuming $\Omega^{(\mathrm{free})}$]
    In the same setting as \pref{lem:alg-diag-basic}, suppose additionally that $\Omega^{(\mathrm{free})}$ holds.
    Then there is a constant $C_{\diag}$ depending on $C_{\mathrm{free}}$ so that for all $m \in (-1,1)^n$,
    \label{lem:alg-resolvent-bounds-delta-frob}
    \[ \schnorm{\delta_{\diag}(m)} \le  \frac{C_{\diag}}{\gamma^4\sqrt{n}}. \]
\end{lemma}
\begin{proof}
The definition \pref{eq:event-free} of $\Omega^{(\mathrm{free})}$ tells us
\[\schnorm{E_{\calD_n}\left[\bar{Q}{}^2\right]D^2 - \frac{\Id_n}{1-\beta^2\tr_n(D^{-2})}} \le \frac{C_{\mathrm{free}}\beta}{\gamma^4\sqrt{n}}.\]
Then triangle inequality with the smallness of the diagonal error from the rank-1 term of $\hat{Q}$ (\pref{lem:alg-resolvent-bounds-DQ2-hQ2D}) allows us to replace the $\bar{Q}$ with $\hat{Q}$ to get
\begin{equation}
\label{eq:alg-hQ2D2-diagonal-bound}
\schnorm{E_{\calD_n}\left[\hat{Q}{}^2\right]D^2 - \frac{\Id_n}{1-\beta^2\tr_n(D^{-2})}} \le \frac{C_{\mathrm{free}}\beta + 16\beta^2\gamma}{\gamma^4\sqrt{n}}.
\end{equation}
Contractivity of $D^{-1}$ finds
\[\schnorm{E_{\calD_n}\left[\hat{Q}^2\right] - \frac{D^{-2}}{1-\beta^2\tr_n(D^{-2})}} \le \frac{C_{\mathrm{free}}\beta + 16\beta^2\gamma}{\gamma^4\sqrt{n}}.\]
Taking the trace of the above and applying $\tr_n(\cdot) \le \schnorm[1]{\cdot} \le \schnorm[2]{\cdot}$ yields the inequality in
    

    \[ \left|\left(1+\beta^2\tr_n(\hat{Q}^2)\right) - \left(\frac{1}{1 - \beta^2\tr_n(D^{-2})}\right)\right| = \left|\left(1+\beta^2\tr_n(\hat{Q}^2)\right) - \left(1 + \frac{\beta^2\tr_n(D^{-2})}{1 - \beta^2\tr_n(D^{-2})}\right)\right| \le \frac{C_{\mathrm{free}}\beta^3 + 16\beta^4\gamma}{\gamma^4\sqrt{n}}. \]
    Then triangle inequality with \pref{eq:alg-hQ2D2-diagonal-bound} obtains the bound.
\end{proof}

\subsection{Error between ASL-TAP and PHD is small}
\label{sec:alg-phd-asl-tap-closeness}
When we compare the two SDEs
\[ d\hat{m}_t = \hat{Q}(\hat{m}_t)dB_t \]
\[
        d\hat \mg_t = \hQ(\hat \mg_t)\pa{\hat \mg_t - f(\hat \mg_t)} dt + \hQ(\hat \mg_t)dB_t
\]
from \pref{e:PHD} and \pref{e:ASL-TAP} respectively, we see that they only differ by the ASL-TAP drift term 
\[\hQ(\hat \mg_t)\pa{\hat \mg_t - f(\hat \mg_t)} dt.\]

With the basic resolvent bounds and the control over the diagonal in hand, we are ready to show that this error between the two SDEs is in fact small and regular.
This closeness between PHD and ASL-TAP then underpins the functioning of the sampling algorithm, as it enjoys both PHD's regularity and connection to the TAP free energy (as leveraged in \pref{sec:pha-to-score}) and ASL's direct connection to the Gibbs measure.

\begin{lemma}[Closeness of ASL-TAP and PHD]
Suppose that $\Omega^{(A)}_{\gamma}$ and $\Omega^{(\mathrm{free})}$ hold.
Let $D(m), a(m), \hat{Q}(m)$ be as in \pref{eq:DahQ} and let $\delta_{\diag}(n)$ be as in \pref{eq:delta-diag}. 
Let $f(m) := (\delta_{\diag}(m) + \Id_n + 2\beta^2\hat{Q}^2/n)m$ as in \pref{e:Itomag}.
Let $C_{\diag}$ be the constant in \pref{lem:alg-resolvent-bounds-delta-frob}.
Then with $m \in (-1,1)^n$,
    \begin{enumerate}[label=(\alph*), ref=\thetheorem(\alph*)]
    \item
\label{lem:alg-resolvent-lipschitz-f}
    \[ \norm{f}_{\Lip} \;\le\; 34\gamma^{-4} + 4\gamma^{-2} + 1 + O\left(n^{-1/2}\right).\]
    \item \label{lem:alg-resolvent-bounds-m-f}
    \[ \lpnorm{m - f(m)} \;\le\; C_{\diag}\gamma^{-4} + O(n^{-1/2}). \]
    \item \label{lem:alg-resolvent-bounds-drift}
    The drift coefficient of \pref{e:ASL-TAP} is bounded by
    \[ \lpnorm{\hat{Q}(m)(m - f(m))} \;\le\; C_{\diag}\gamma^{-5} + O(n^{-1/2}). \]
    \item \label{lem:alg-resolvent-bounds-dual-drift-1}
    The first dual drift coefficient of \pref{eq:d-ASL-TAP} is bounded by
    \[ \lpnorm{D(m)\hat{Q}(m)(m - f(m))} \;\le\; 2C_{\diag}\gamma^{-5} + O(n^{-1/2}). \]
    \item \label{lem:alg-resolvent-lipschitz-drift}
    The drift coefficient of \pref{e:ASL-TAP} is Lipschitz:
    \[ \norm{m \mapsto \hat{Q}(m)(m - f(m))}_{\Lip} \;\le\; 4C_{\diag}\gamma^{-6} + 34\gamma^{-5} + 4\gamma^{-3} + \gamma^{-1} + 1 + O(n^{-1/2}).\]
    \end{enumerate}
    \end{lemma}
    \begin{proof}
    \ppart{Part \ref{lem:alg-resolvent-lipschitz-f}} We have by product rule that 
    \[\norm{f}_{\Lip} \le \norm{\delta_{\diag}}_{\Lip}\sqrt{n} + \frac{2}{\sqrt{n}}\norm{\hat{Q}^2}_{\Lip\to\infty} + C_{\delta} + 1 + \frac{2}{n}C_{\hat{Q}}^2,\]
    where $C_{\delta}$ is a uniform bound on $\opnorm{\delta_{\diag}}$ and $C_{\hat{Q}}$ is a uniform bound on $\opnorm{\hat{Q}}$ and $\norm{\hat{Q}^2}_{\Lip\to\infty}$ is the Lipschitz constant with respect to the operator norm.
    Note that $\norm{\hat{Q}^2}_{\Lip\to\infty} \le \sqrt{n}\norm{\hat{Q}^2}_{\Lip}$, so by the Lipschitzness of $\delta_{\diag}$ in \pref{lem:alg-resolvent-bounds-delta}, the simple resolvent Lipschitzness bound \pref{lem:alg-resolvent-lipschitz-Q2} for $\hat{Q}^2$, control of the diagonal from $\Omega^{\mathrm{free}}$ in \pref{lem:alg-resolvent-lipschitz-delta} for $C_{\delta}$, and the simple resolvent operator norm bound \pref{lem:alg-resolvent-bounds-Q} for $C_{\hat{Q}}$, 
    \[\norm{f}_{\Lip} \le \gamma^{-4}(4+\beta^2)(4 + 16\beta^2) + 4\gamma^{-2} + 1 + O(1/\sqrt{n}).\]

    \ppart{Part \ref{lem:alg-resolvent-bounds-m-f}} We notice that 
    \[m - f(m) = -\delta_{\diag}(m)m - 2\beta^2\hat{Q}^2m/n.\]
    So after applying triangle inequality, the bound for the first term follows by 
    \[\lpnorm{\delta_{\diag}(m)m}^2 = \sum_{i=1}^n \delta_{\diag}(m)_i^2m_i^2 \le \sum_{i=1}^n \delta_{\diag}(m)_i^2 \le n\schnorm{\delta_{\diag}}^2\]
    and control of the diagonal from  $\Omega^{(\mathrm{free})}$ (\pref{lem:alg-resolvent-bounds-delta-frob}).
    For the second term, $\lpnorm{\hat{Q}^2m} \le \opnorm{\hat{Q}}^2\lpnorm{m}$ and then the simple resolvent operator norm bound \pref{lem:alg-resolvent-bounds-Q} concludes.

    \ppart{Parts \ref{lem:alg-resolvent-bounds-drift} and \ref{lem:alg-resolvent-bounds-dual-drift-1}} We combine the just-established bound on $\lpnorm{m - f(m)} ($\pref{lem:alg-resolvent-bounds-m-f}) with the simple resolvent operator norm bounds \pref{lem:alg-resolvent-bounds-Q} and \pref{lem:alg-resolvent-bounds-DQ} respectively.
    
    \ppart{Part \ref{lem:alg-resolvent-lipschitz-drift}} By the product rule,
    \[\norm{m \mapsto \hat{Q}(m)(m - f(m))}_{\Lip} \le C_{\hat{Q}}\left(1 + \norm{f}_{\Lip}\right) + C_{mf}\norm{\hat{Q}}_{\Lip\to\infty},\]
    where $C_{\hat{Q}}$ is a uniform bound on $\opnorm{\hat{Q}}$ and $C_{mf}$ is a uniform bound on $\lpnorm{m - f(m)}$ and $\norm{\cdot}_{\Lip\to\infty}$ is the Lipschitz constant with respect to the operator norm. 
    From \pref{lem:alg-resolvent-bounds-Q}, we can take $C_{\hat{Q}} = \gamma^{-1}$ and from \pref{lem:alg-resolvent-bounds-m-f}, we can take $C_{mf} \le C_{\diag}\gamma^{-4} + O(n^{-1/2})$.
    By \pref{lem:alg-resolvent-lipschitz-f}, we also have $\norm{f}_{\Lip} \le 34\gamma^{-4} + 4\gamma^{-2} + 1 + O\left(n^{-1/2}\right)$ and from \pref{lem:alg-resolvent-lipschitz-Q} we have $\norm{\hat{Q}}_{\Lip\to\infty} \le \sqrt{n}\norm{\hat{Q}}_{\Lip} \le 4\gamma^{-2} + O(1/\sqrt{n})$.
    \end{proof}

\subsection{Properties of the trajectory}
\label{sec:alg-trajectories}

We now concern ourselves with the \pref{e:ASL-TAP} and the \pref{e:PHD} evolutions for $\hat{\mg}_t$, specifically the initial condition $\hat \mg_0 = 0$ and the SDEs
\begin{equation}
\label{eq:asl-tap-m}
d\hat{\mg}_t
=
\hat Q(\hat{\mg}_t)\big(\hat{\mg}_t-f(\hat{\mg}_t)\big)\,dt
+\hat Q(\hat{\mg}_t)\,dB_t.\end{equation}
\begin{equation}
\label{eq:phd-m}
d\hat{\mg}_t
=
\hat Q(\hat{\mg}_t)\,dB_t.\end{equation}
The boundedness and Lipschitzness properties established so far will be enough to show the basic fact that $\hat{\mg}_t$ remains inside $(-1,1)^n$ at all times (\pref{lem:stay-in-cube}).
Moreover, they will imply transportation inequalities that yield a strong sub-Gaussian concentration statement for Lipschitz functionals of $\{\hat{m}_t\}_t$ (\pref{thm:mSDE-path-subgaussian}).

For the proof of existence, uniqueness, and non-escape of the SDE solutions, it is convenient to make a change of coordinates $\tanh^{-1}:(-1,1)^n\to \R^n$ to the dual space of the magnetization.

We apply this change of coordinates to \eqref{e:PHD}, given by $\hat{\mg}_0 = 0$ and $d\hat \mg_t = \hat Q(\hat\mg_t) dB_t$.
Let
$\hat x_t = \tanh^{-1}(\hat \mg_t)$.  
Applying Itô's formula to the change of variables $\tanh^{-1}:(-1,1)^n\to \R^n$ gives
\begin{align}
    d\hat x_t &= 
    \ED\pa{\hat Q(\hat \mg_t)^2}  \pf{\hat \mg_t}{(1-\hat \mg_t^2)^2}
    dt + 
    \diag\prc{1-\hat \mg_t^2} \hat Q(\hat \mg_t)dB_t
    \nonumber
    \\
    \label{eq:d-PHD}
    &= \ub{ \ED\pa{\hat Q(\hat \mg_t)^2} D(\hat \mg_t)^2\hat \mg_t}{\hat f^*(\hat x_t)} dt + \ub{D(\hat \mg_t)\hat Q(\hat \mg_t)}{\hat \Si(\hat x_t)} dB_t,
    \tag{d-PHD}
\end{align}
where operations on vectors are interpreted component-wise. Define $\hat f^*(\hat x_t)$ and $\hat \Si(\hat x_t)$ as above, noting that $\hat \mg_t = \tanh(\hat x_t)$, and that here, $\hat{f}^*$ is not denoting a convex conjugation. 

We can similarly apply this change of coordinates to \eqref{e:ASL-TAP} to obtain
\begin{equation}
\label{eq:d-ASL-TAP}
d\hx_t = 
\ba{D(\hat\mg_t)\hQ(\hat\mg_t) (\hat \mg_t - f(\hat \mg_t)) + \ED\pa{\hQ(\hat\mg_t)^2} D(\hat \mg_t)^2\hat \mg_t} dt + D(\hat \mg_t)\hQ(\hat \mg_t) dB_t.
\tag{d-ASL-TAP}
\end{equation}

We show that the coefficients of \eqref{eq:d-PHD} and \eqref{eq:d-ASL-TAP} are Lipschitz.
\begin{lemma}[Lipschitzness of dual SDE coefficients]
\label{l:dual-phd}
Suppose that $\Omega^{(A)}_{\gamma}$ and $\Omega^{(\mathrm{free})}$ hold.
Then over $\R^n$,
    \begin{enumerate}
        \item 
        \label{i:lip-hS}
        $\hat \Si$ is Lipschitz (with respect to $\schnorm{\cdot}$) with constant $\frac{\gamma^{-2}}{\sqrt{n}}\left(4 + 16\beta^2 + \frac{32\beta^2}{\sqrt{n}}\right)$.
        \item 
        \label{i:lip-hfs}
     $\hat f^*$ is Lipschitz with constant $\gamma^{-3}\left(16 + 64\beta^2 + \frac{128\beta^2}{\sqrt{n}}\right) + 4\gamma^{-2}$. 
        \item With $m := \tanh(x)$, the function $x \mapsto D(m)\hQ(m) (m - f(m))$ is Lipschitz with constant $8C_{\diag}\gamma^{-6} + 68\gamma^{-5} + 8\gamma^{-3} + 2\gamma^{-1} + 2 + O(n^{-1/2})$.
    \end{enumerate}
\end{lemma}

\begin{proof}
\ppart{Lipschitzness of $\hat \Si$}
This follows by combining \pref{lem:alg-resolvent-lipschitz-DQ} with the Lipschitzness of $\tanh$ and the multiplicativity of Lipschitz constants under composition of functions.

\ppart{Lipschitzness of $\hat f^*$}
Since $\tanh$ is Lipschitz, it suffices to bound the Lipschitz constant of $\hat f^*(x)$ as a function of $m := \tanh(x)$.
By product rule, this Lipschitz constant is bounded by
\[ \norm{m \mapsto \ED\pa{\hat Q(m)^2} D(m)^2}_{\Lip}\sqrt{n} + \sup_{m \in (-1,1)} \opnorm{\ED\pa{\hat Q(m)^2} D(m)^2}. \]
The first term is bounded by \pref{lem:alg-resolvent-lipschitz-DQQD} after noting the contractivity of $\ED$.
The second term is bounded by \pref{lem:alg-resolvent-bounds-DQ2D}.

\ppart{Lipschitzness of $x \mapsto D(m)\hQ(m) (m - f(m))$} This follows by the same product-rule argument as \pref{lem:alg-resolvent-lipschitz-drift}, but replacing \pref{lem:alg-resolvent-bounds-Q} and \pref{lem:alg-resolvent-lipschitz-Q} with \pref{lem:alg-resolvent-bounds-DQ} and \pref{lem:alg-resolvent-lipschitz-DQ}
\end{proof}

\begin{lemma}[ASL-TAP and PHD stay in the cube]
\label{lem:stay-in-cube}
Assume that $\Omega^{(A)}_{\gamma}$ and $\Omega^{(\mathrm{free})}$ hold. Then \eqref{eq:d-PHD}, \eqref{eq:d-ASL-TAP}, \eqref{e:PHD}, and \eqref{e:ASL-TAP} all have unique global strong solutions, and in all of them, 
$\hat{\mg}_t\in(-1,1)^n$ for all $t\ge 0$ almost surely.
\end{lemma}

\begin{prf}
By \pref{l:dual-phd}, the coefficients of $dt$ and $dB_t$ for $d\hat{x}_t$ are globally Lipschitz as functions of $\hat{x}_t$ in both \eqref{eq:d-PHD} and \eqref{eq:d-ASL-TAP}. This implies that the coefficients are also
locally bounded. Therefore, the standard Lipschitz theorem for SDEs yields a
unique global adapted solution for $\hat{x}_t$ driven by the Brownian motion $B_t$, and the solution is strong with respect to the Brownian
filtration \cite[Chapter IX, Theorem 2.1]{Revuz2010-nk}. 

Since $\tanh:\R^n\to(-1,1)^n$ coordinate-wise is $C^\infty$, we can apply It\^o's formula
again to $\hat{\mg}_t:=\tanh(\hat{x}_t)$, showing that this reverses the dual transformation and recovers solutions to \eqref{e:PHD} and \eqref{e:ASL-TAP} respectively. In particular, these solutions for $\hat{\mg}$ are
continuous semimartingale satisfying $\hat{\mg}_t\in(-1,1)^n$ for every $t\ge 0$.

It remains to prove pathwise uniqueness for \eqref{e:PHD} and \eqref{e:ASL-TAP}.
Let $\tilde{\mg}$ be any other continuous primal solution, driven by the same
Brownian motion and started from the same initial condition, and define its exit
time
\[
\tilde\tau:=\inf\{t\ge 0:\tilde{\mg}_t\notin(-1,1)^n\}.
\]
For $t<\tilde\tau$, set
\[
\tilde x_t:=\tanh^{-1}(\tilde{\mg}_t).
\]
Again, since $\tanh^{-1}$ is $C^\infty$ on $(-1,1)^n$, It\^o's
formula shows that $\tilde x$ satisfies the dual SDE on $[0,\tilde\tau)$.

Fix $T>0$. On the event $\{T<\tilde\tau\}$, both $\tilde x$ and $\hat{x}$ solve
the same dual SDE on $[0,T]$ with the same Brownian motion and the same initial
condition. By pathwise uniqueness for the dual SDE, they coincide on $[0,T]$ on
that event. Letting $T$ range over the positive rationals and extending to irrational $t$ by continuity, it follows that
\[
\tilde x_t= \hat{x}_t \qquad\text{for all } t<\tilde\tau \quad\text{a.s.}
\]
Therefore
\[
\tilde{\mg}_t
=
\tanh(\tilde x_t)
=
\tanh(\hat{x}_t)
=
\hat{\mg}_t
\qquad\text{for all } t<\tilde\tau.
\]

If $\tilde\tau<\infty$, then by continuity of $\hat{x}$ and $\tanh$,
\[
\tilde{\mg}_{\tilde\tau}
=
\lim_{t\uparrow\tilde\tau}\tilde{\mg}_t
=
\lim_{t\uparrow\tilde\tau}\tanh(\hat{x}_t)
=
\tanh(\hat{x}_{\tilde\tau})
\in(-1,1)^n,
\]
which contradicts the definition of $\tilde\tau$.
Hence $\tilde\tau=\infty$ a.s., and consequently
\[
\tilde{\mg}_t=\hat{\mg}_t \qquad\text{for all } t\ge 0 \quad\text{a.s.}
\]
Thus the primal SDE is pathwise unique.
\end{prf}

To obtain sub-Gaussian concentration for all Lipschitz functionals of the trajectory of $\hat{m}_t$, we first restate a powerful characterization of transportation inequalities in terms of said sub-Gaussian concentration.
Then we adjust and restate a result of \"Utsu\"nel yielding a transportation inequality for any It\^o process that has a Lipschitz drift coefficient and a bounded and Lipschitz diffusion coefficient.

\begin{definition}
    \label{def:Tp-transport-inequality}
    A probability measure $\mu$ on a metric space $(E,d)$ satisfies a $T_p$ transportation inequality with constant $C$ if \[W_p(\mu,\nu) \le \sqrt{2C\,H(\nu \mid \mu)},\]
    where 
    \[W_p(\mu,\nu) := \inf_{\pi \in \Sigma(\mu,\nu)} \left(\int (d(x,y))^pd\pi(x,y)\right)^{1/p}\]
    is the Wasserstein distance, where  $\Sigma(\mu,\nu)$ is the subset of probability measures over $E \times E$ whose first marginal is $\mu$ and whose second marginal is $\nu$, and $H(\nu \mid \mu)$ is the relative entropy.
\end{definition}

\begin{theorem}[$T_1$ inequality \& sub-Gaussian concentration~{\cite[Theorem 1.1, Bobkov \& G\"otze]{djellout2004}}]
    \label{thm:t1-subgaussian}
    Let $\mu$ be a measure on a metric space $(E,d)$ that satisfies a $T_1$ transport inequality with constant $C > 0$. This is equivalent to $\mu$ satisfying the following concentration under push-forward by any $K$-Lipschitz function $F: (E,d) \to \R$ that is $\mu$-integrable:
    first, in moment-generating function form,
    \[ \E_{\mu}\left[\exp\left(\lambda\left(F - \E_{\mu}[F]\right)\right)\right] \le \exp\left(\frac{\lambda^2CK^2}{2}\right)\qquad \forall \lambda \in \R, \]
    and in sub-Gaussian tail bound form,
    \[
        \mu\left(\left|F - \E_{\mu}[F]\right| > r\right) \le 2\exp\left(-\frac{r^2}{2CK^2}\right) \qquad \forall r > 0.
    \]
\end{theorem}

\begin{theorem}[(Modified) $T_2$ transport inequality for regular It\^o processes~{\cite[Theorem 1]{ustunel2012transportation}}]
\label{thm:t2-transport-SDE}
    Let $0 < C, K < \infty$ and $0 < T < \infty$.
    Consider an It\^o process $\{X_t \in \R^n\}_{t\in[0,T]}$ given by $X_0 = 0$ and
    \[
        dX_t = b_t(X_t) dt + \sigma_t(X_t) dB_t\,.    
    \]
    for $t \in (0, T)$.
    Assume that the diffusion coefficients $\{\sigma_t: \R^n \to \R^{n\times n}\}_{0 \le t \le T}$ are matrix-valued functions satisfying $\opnorm{\sigma_t(x)} \le C$ for every $0 \le t \le T$ and $x \in \R^n$.
    Assume also that both the diffusion coefficients and the drift coefficients $\{b_t: \R^n \to \R^n\}_{0 \le t \le T}$ satisfy the following Lipschitz conditions for all $x,y \in \R^n$ and $0 \le t \le T$:
    \[
        \lpnorm{b_t(x) - b_t(y)} \le K\lpnorm{x-y}\,,
    \]
    \[
        \schnorm{\sigma_t(x) - \sigma_t(y)} \le \frac{K}{\sqrt{n}}\lpnorm{x-y}.
    \] 
    Then $\dist(X)$ satisfies a $T_2$ transportation inequality on $C^0([0,T],\R^n)$ with constant $3TC^2e^{3T(T+4)K^2}$, where $X$ refers to the entire path $\{X_t\}_{t \in [0,T]}$, $C^0([0,T],\R^n)$ is the set of continous paths in $\R^n$ on $[0,T]$, and the metric on $C^0([0,T],\R^n)$ is given by the uniform norm $d(X,Y) := \sup_{0 \le t \le T}\lpnorm{X_t - Y_t}$.
\end{theorem}
\begin{proof}
    The proof is identical to \"Utsu\"nel's original argument, but with minute changes to accommodate $t \in [0,T]$ instead of $t \in [0,1]$ and $\opnorm{\sigma_t(x)} \le C$ instead of $\opnorm{\sigma_t(x)} \le 1$.
    Specifically, after establishing the domination $d_W^2(\dist(X),\mu) \le \E_{\mu}[\sup_{0\le t \le T}\lpnorm{x^v_t - x_t}^2]$, we calculate the difference between the SDEs $dx_t = \sigma_t(x_t)d\beta_t + b_t(x_t)dt$ and $dx^v_t = \sigma_t(x^v_t)dz_t + b_t(x^v_t)dt$ with $d\beta_t = dz_t - dv_t$ to see
    \[ x^v_r - x_r = \int_0^r(\sigma_s(x^v_s) - \sigma_s(x_s))dz_s + \int_0^r(b_s(x^v_s) - b_s(x_s))ds + \int_0^r \sigma_s(x_s)dv_s.  \]
    By $\norm{a+b+c}^2 \le 3(\norm{a}^2 + \norm{b}^2 + \norm{c}^2)$,
    \begin{align*}
        \lpnorm{x^v_r - x_r}^2
        &\le 3\lpnorm{\int_0^r(\sigma_s(x^v_s) - \sigma_s(x_s))dz_s}^2 + 3\lpnorm{\int_0^r(b_s(x^v_s) - b_s(x_s))ds}^2 + 3\lpnorm{\int_0^r \sigma_s(x_s)dv_s}^2.
    \end{align*}
    By Doob's $L^2$ maximal inequality, It\^o isometry, and then the Lipschitzness of $\sigma_s$,
    \begin{align*}
    \E_{\mu}\left[\sup_{r \le t}\lpnorm{\int_0^r(\sigma_s(x^v_s) - \sigma_s(x_s))dz_s}^2\right]
    &\le 4\E_{\mu}\left[\lpnorm{\int_0^t(\sigma_s(x^v_s) - \sigma_s(x_s))dz_s}^2\right]
    \\&= 4n\E_{\mu}\left[\int_0^t\norm{\sigma_s(x^v_s) - \sigma_s(x_s)}_F^2ds\right]
    \\&\le 4K^2\E_{\mu}\left[\int_0^t\lpnorm{x^v_s - x_s}^2ds\right].
    \end{align*}
    By the Cauchy-Schwarz inequality and the Lipschitzness of $b_s$,
    \begin{align*}
    \sup_{r \le t}\lpnorm{\int_0^r(b_s(x^v_s) - b_s(x_s))ds}^2
    \le t\int_0^t\lpnorm{b_s(x^v_s) - b_s(x_s)}^2ds
    \le K^2t\int_0^t\lpnorm{x^v_s - x_s}^2ds.
    \end{align*}
    Then by the Cauchy-Schwarz inequality and the boundedness of $\sigma_s$,
    \begin{align*}
    \sup_{r \le t}\lpnorm{\int_0^r\sigma_s(x_s)dv_s}^2
    \le t\int_0^t\lpnorm{\sigma_s(x_s)\dot v_s}^2ds
    \le C^2t\int_0^t\lpnorm{\dot v_s}^2ds.
    \end{align*}
    Combining the last 4 displayed inequalities and letting $u(t) := \E_{\mu}\left[\sup_{r \le t}\lpnorm{x^v_t - x_t}^2\right]$,
    \[ u(t) \le 3K^2(t+4)\int_0^tu(s)ds + 3C^2t\E_{\mu}\left[\int_0^t\lpnorm{\dot v_s}^2ds\right].\]
    The result follows from Gr\"onwall's inequality and the property of the exponential martingale that $\E_{\mu}\left[\int_0^T\lpnorm{\dot v_s}^2ds\right] = 2\, H(\mu \mid \dist(X))$.
\end{proof}

\begin{theorem}[Transport and sub-Gaussian concentration for ASL-TAP]
\label{thm:mSDE-path-subgaussian}
Suppose that $\Omega^{(A)}_{\gamma}$ from \pref{eq:event-opnorm} and $\Omega^{(\mathrm{free})}$ from \pref{eq:event-free} both hold.

Let $0 < T < \infty$.
Let $\calW:=C^0([0,T],\R^n)$, equipped with the uniform norm
$\|x\|_\infty:=\sup_{0\le t\le T}\lpnorm{x_t}$, and let
$P:=\dist(\hat{m}_\cdot)$ be the law of the full path of the ASL-TAP SDE~\pref{eq:asl-tap-m} on $\calW$.
For some constant $C_K$, define
\[
c_{T_2}:=\frac{3T}{\gamma^2}\exp\left(\frac{3T(T+4)C_K^2}{\gamma^{12}}\right).
\]
Then:

\begin{enumerate}
\item[\textup{(i)}]
For every probability
measure $\mu$ on $\big(\calW,\calB(\calW)\big)$, where $\calB$ is the Borel $\sigma$-algebra on $\calW$,
\[
d_{W,2}(P,\mu)^2\ \le\ 2c_{T_2}\,H(\mu\mid P).
\]

\item[\textup{(ii)}]
For every $P$-integrable $K_F$-Lipschitz functional $F:(\calW,\|\cdot\|_\infty)\to\R$,
    \[ \E_{P}\left[\exp\left(\lambda\left(F - \E_{P}[F]\right)\right)\right] \le \exp\left(\frac{\lambda^2c_{T_2}K_F^2}{2}\right),\qquad \forall \lambda \in \R. \]
    In sub-Gaussian tail bound form,
\[
\P_P\Big(\big|F-\E_P[F]\big|>r\Big)
\ \le\
2\exp\!\left(-\frac{r^2}{2c_{T_2}K_F^2}\right),
\qquad r>0.
\]
\end{enumerate}
\end{theorem}

\begin{proof}
To apply \pref{thm:t2-transport-SDE} to the SDE~\eqref{eq:asl-tap-m}, we need to show that the drift coefficient $\hat Q(\hat{\mg}_t)\big(\hat{\mg}_t-f(\hat{\mg}_t)\big)$ is $K$-Lipschitz and the diffusion coefficient $\hat Q(\hat{\mg}_t)$ is $K/\sqrt{n}$-Lipschitz and bounded everywhere by $C$. By \pref{lem:stay-in-cube} and the Kirszbraun theorem, we can use a Lipschitz extension of $\hat{Q}$ outside of $(-1,1)^n$ without affecting the trajectories, so
the basic resolvent operator norm bound \pref{lem:alg-resolvent-bounds-Q} shows that we can take $C \gets \gamma^{-1}$.
The closeness of PHD and ASL-TAP  (\pref{lem:alg-resolvent-lipschitz-drift}) controls the Lipschitz constant of the drift term, which, when combined with a basic resolvent Lipschitzness bound (\pref{lem:alg-resolvent-lipschitz-Q}) on the diffusion term, shows that we can take $K \gets \max(C_K\gamma^{-6}, 4\gamma^{-2}) + O(1/\sqrt{n})$ for some constant $C_K$.
Hence we find the stated $T_2$ inequality with constant $c_{T_2}=3TC^2e^{3T(T+4)K^2}$.

Since $d_{W,1}\le d_{W,2}$, this implies a $T_1$ inequality with constant $c_{T_2}$.
Applying~\pref{thm:t1-subgaussian} to a $K_F$-Lipschitz $F$ yields the conclusion.
\end{proof}

\subsection{Implications for the sampling algorithm}
\label{sec:alg-desiderata}

We can now index various algorithmic desiderata stated in \pref{s:desiderata} into the algorithmic process properties we just proved.

\begin{corollary}[Algorithmic desiderata]
\label{cor:alg-desiderata}
Suppose that $\Omega^{(A)}_{\gamma}$ from \pref{eq:event-opnorm} and $\Omega^{(\mathrm{free})}$ from \pref{eq:event-free} both hold. Then
\begin{enumerate}[itemsep=0.7em,label=\textup{(d\arabic*)},ref={\thetheorem\,(d\arabic*)},start=0]
\item\label{cor:alg-desiderata-d0}
The algorithmic SDEs \eqref{e:PHD} and \eqref{e:asl-tap-je} have unique strong solutions such that $\hat{m}_t \in (-1,1)^n$ for all $t < \infty$ almost-surely.
\setcounter{enumi}{1}
\item\label{cor:alg-desiderata-d2}
There are constants $\lambda,\Lambda>0$ such that 
    for all $\mg\in (-1,1)^n$,
    \[
        \lambda D(\mg)^{-1} \preceq \hQ(\mg) \preceq \Lambda D(\mg)^{-1}.
    \]
\item\label{cor:alg-desiderata-d3}
    For all $t\in [0,T]$, for $\mg_t$ drawn from \pref{e:HD}, 
    \[
    \E\ve{f(\mg_t)-m_t}_2^2 \le \ep_{\mathrm{drift}}(t)^2,\]
    where $\ep_{\mathrm{drift}}(t)$ is a constant depending only on $\beta$ and $t$ and bounded on $[0,T]$.
\item\label{cor:alg-desiderata-d4}
\begin{enumerate}
        \item There exists a constant $L$  such that for all $\mg, w\in (-1,1)^n$,
    \[
    \ve{\hQ(\mg) - \hQ(w)}_F \le L\ve{\mg-w}_2.
    \]
        \item 
        There exists a constant $L_f$  such that for all $\mg, w\in (-1,1)^n$, 
    \[
    \ve{f(\mg) - f(w)}_2 \le L_f\ve{\mg - w}_2.
    \]
    \end{enumerate}
\item\label{cor:alg-desiderata-d5} As in  \pref{e:asl-tap-je} and with $\hat{m}_t$ solving \pref{eq:asl-tap-m}, let 
\[
\hat{w}_T\ :=\ \frac{1}{2}\int_0^T \left(\Tr\big(\hat Q(\hat{\mg}_t)\big)\;+\;\lpnorm{\hat{\mg}_t}^2\right)dt.
\]
Then for all $T \in [0,\infty)$, $\hat{w}_T$ is sub-Gaussian with variance proxy $\lambda^2c_{T_2}C_{\diag}^2T^2/\gamma^8$, i.e. for all $\lambda \in \R$,
    \[ \E_{B}\left[\exp\left(\lambda\hat{w}_T - \E_{B}[\lambda\hat{w}_T]\right)\right] \le \exp\left(\frac{\lambda^2c_{T_2}C_{\diag}^2T^2}{2\gamma^{8}}\right), \]
    where $c_{T_2}$ is the function of $T$ and $\gamma$ given in \pref{thm:mSDE-path-subgaussian} and $C_{\diag}$ is the constant (for $\beta < 1/2$) in \pref{lem:alg-resolvent-bounds-delta-frob}.
\end{enumerate}
\end{corollary}
\begin{proof}
We give the proofs in sequential order. First, \pref{cor:alg-desiderata-d0} is given by \pref{lem:stay-in-cube}.

\ppart{Part \ref{cor:alg-desiderata-d2}} The lower bound is \pref{lem:tap-convexity}. The upper bound follows from the operator bound on $D\hat{Q}$ given $\Omega^{(A)}_{\gamma}$ (\pref{lem:alg-resolvent-bounds-DQ}) with $\Lambda \gets 2\gamma^{-1}$, since $\opnorm{D\hat{Q}} \le 2\gamma^{-1}$ implies $D^{1/2}\hat{Q}D^{1/2} \preceq 2\gamma^{-1}\Id_n$, and then multiply on the left and right by $D^{-1/2}$.

\ppart{Part \ref{cor:alg-desiderata-d3}} The closeness of PHD and ASL-TAP (\pref{lem:alg-resolvent-bounds-m-f}) implies $\lpnorm{f(m) - m}^2 \le C_{\diag}\gamma^{-4} + O(n^{-1/2})$ uniformly for all $m \in (-1,1)^n$ for a constant $C_{\diag}$, therefore the inequality is also true in expectation over $m_t$ drawn from \pref{e:HD}.

\ppart{Part \ref{cor:alg-desiderata-d4}}
The first item is \pref{lem:alg-resolvent-lipschitz-Q}, obtaining $L \gets 4\gamma^{-2} + O(1/\sqrt{n})$.
The second item is \pref{lem:alg-resolvent-lipschitz-f}, with $L_f \gets 34\gamma^{-4} + 4\gamma^{-2} + 1 + O(1/\sqrt{n})$.

\ppart{Part \ref{cor:alg-desiderata-d5}} 
Let $\omega(m):=\frac{1}{2}\left(\Tr(\hat Q(m))+\lpnorm{m}^2\right)$.
To show that $\omega$ is Lipschitz, we differentiate with respect to $m$ using Fr\'echet derivatives (\pref{def:frechet}):
    \begin{align}
    \nonumber
    d\omega &= -\frac{1}{2}\Tr\left[\hat Q\left(-\frac{2\beta^2}{n}m^{\sT} dm\cdot\Id_n + 2D\diag(m\odot dm)D - \frac{2\beta^2}{n}\left(m\,dm^{\sT}  + (dm)m^{\sT} \right)\right)\hat Q\right] + m^{\sT} dm
    \\\nonumber&= -\Tr\left[\hat{Q}^2D^2\diag(m)\diag(dm)\right] - \beta^2m^{\sT} \left(\tr_n[\hat{Q}^2]\Id_n + \frac{2}{n}\hat{Q}^2\right)dm + m^{\sT} dm
    \\\nonumber&= -m^{\sT} E_{\mathcal{D}_n}\left[\hat{Q}^2D^2\right]dm + \beta^2m^{\sT} \left(\tr_n[\hat{Q}^2]\Id_n + \frac{2}{n}\hat{Q}^2\right)dm + m^{\sT} dm
    \\\nonumber&= m^{\sT} \left(\beta^2\left(\tr_n(\hat{Q}^2)\Id_n + \frac{2}{n}\hat{Q}^2\right) + \Id_n - E_{\mathcal{D}_n}\left[\hat{Q}^2\right]D^2\right)dm
    \\&= (m - f(m))^{\sT} dm
    .\label{e:d-om}\end{align}
    Therefore, by bounded Fr\'echet derivatives implying Lipschitzness (\pref{fact:frechet-lipschitz}) and the closeness of PHD and ASL-TAP (\pref{lem:alg-resolvent-bounds-m-f}), $\omega$ is $C_{\diag}\gamma^{-4}$-Lipschitz.
So for $x,y\in C^0([0,T],\R^n)$,
\[
\left|\hat{w}_T(x)-\hat{w}_T(y)\right|
\le \int_0^T |\omega(x_t)-\omega(y_t)|\,dt
\le TC_{\diag}\gamma^{-4}\,\sup_{t\in[0,T]}\lpnorm{x_t-y_t},
\]
so $\lipnorm{\hat{w}_T}\le TC_{\diag}\gamma^{-4}$, with $\hat{w}_T$ as a functional of the path $\{\hat{m}_t\}_{t \in [0,T]}$ metricized by the uniform norm.
Then apply sub-Gaussian concentration for Lipschitz path functionals of ASL-TAP (\pref{thm:mSDE-path-subgaussian}) with $F \gets \hat{w}_T$ and $K_F \gets TC_{\diag}\gamma^{-4}$.
\end{proof}

\section{Sampling from the time-\texorpdfstring{$T$}{T} localized distribution}
\label{s:localized}


Our main theorem is the following: the stochastic localization process for SK will concentrate the measure on a small wedge $B_{\ep n}$ after large constant time, and when the wedge is small enough, 
the measure will allow rapid mixing for an appropriate Markov chain. For technical reasons, we prove fast mixing for the polarized walk (see \Cref{d:pw}) rather than Glauber dynamics, although we conjecture the same result holds for Glauber dynamics. Note that entropy contraction implies a modified log-Sobolev inequality for the corresponding Markov chain, in this case the polarized walk.
\begin{theorem}
\label{t:ls-sk}
Consider the stochastic localization process for the Gibbs distribution of the SK model with inverse temperature $\be<\rc 2$. 
There are constants $T(\beta)$, $\ep(\beta)$, and $c(\beta)$ such that 
    with probability $\ge 1-e^{-c(\be) n}$, for $T\ge T(\beta)$, 
    \begin{enumerate}
        \item $\mu_{\be A,y_T}(\ball{\ep(\be) n}{\sign(y_T)}) \ge 1-e^{-c(\be) n}$.
        \item For any $x_0$, 
        the polarized walk on $\mu_{\be A,y_T}|_{\ball{\ep(\be)n}{x_0}}$ satisfies entropy contraction with constant $\ge\fc{\rh_\be}{n}$.
    \end{enumerate}
\end{theorem}
This will follow from \Cref{l:ls-sk-wedge-ii} (entropy contraction 
on the wedge) and \Cref{l:conc-x0} (concentration of Gibbs measure on the wedge). The main intuition for (2) is that a wedge has smaller diameter, which effectively counteracts a large operator norm of $\be A$. Concretely, this is formalized through the fact that $A_{S\times S}$ has smaller operator norm for smaller sets $S$.

\begin{remark}[Necessity of $\be<\rc2$]
    We conjecture that the theorem holds for all $\be>0$. Our proof relies on $\be<\rc2$, because the needle decomposition of \cite{eldan2022spectral} into rank-1 Ising models (\Cref{t:needle}) only works for PSD matrices. Our $A$ is not PSD, which we fix by adding a multiple of the identity. However, this breaks the property of $A$ having small operator norm on subsets of $O(\ep n)$ indices, which we need for the resulting rank-1 Ising models to have bounded covariance. When $\be<\rc 2$, we can add $\ga I$ with $\ga<1$ to make it PSD, and then choose a small enough $\ep>0$ so the positive part of $A$ has arbitrarily small operator norm on subsets of $\ep n$ indices, for small enough $\ep$.
    Our threshold for $\be$ is 2 times the threshold for \cite{eldan2022spectral} because we only have to ensure the most negative eigenvalue of $A$ is bounded by 1, rather than the difference between the extreme eigenvalues.
\end{remark}
After giving preliminaries in \pref{s:fi}, we first prove in \pref{s:lsi-1} 
entropy contraction
when $\be<\rc{2\sqrt{2}}$. Our proof is significantly simpler in this case, and showcases the main idea, the two-stage decomposition of \pref{lem:LS-2-stage}. This reduces to bounding the variance of any 1-dimensional marginal arising in the decomposition, and central to this is the  bound on the covariance of a product distribution restricted to a wedge by its diagonal (\pref{l:prod-dist-cov}). Unfortunately, this naive bound is off by a factor of 2 compared to the unrestricted setting, and we have to do some significant work in proving and then applying a tighter version of this bound that holds for most tilts. We do this in \pref{s:lsi-2}. In \pref{s:hi-prob-wedge}, we show that for large enough $T$, with high probability that most of the mass is concentrated in the wedge. Finally, in \pref{s:samp-Z}, we show as a corollary that we can efficiently sample and compute the partition function for $\mu_{A,y_T}$.

We use the polarized walk because our proof requires a Markov chain that is fast mixing for all product distributions on wedges. This is false for Glauber dynamics, as the mass of such a distribution may be concentrated on the boundary slice of the wedge, where any two points differ in at least 2 coordinates. The polarized walk, defined in terms of a down-up walk, allows such transitions.

\subsection{Preliminaries on Markov kernels, Markov chain mixing, and measures on the hypercube}
\label{s:fi}

Throughout, we assume that the state space $\Om$ is finite.
We follow the convention that Markov kernels $P$ operate on the left on functions ($Pf$) and on the right for measures ($\nu P$).
\begin{definition}
    Let $P$ be a Markov kernel from $\Om$ to $\Om'$. We say that $P$ satisfies \vocab{$(1-\rh)$-entropy contraction} with respect to $\mu$ if for any if for any probability measure $\nu\ll \mu$, 
    \[
\KL(\nu P \|\mu P) \le (1-\rh) \KL(\nu\|\mu). 
    \]
    We call the largest $\rh$ for which this holds the entropy contraction constant, and denote it by $\rec(\mu)$ when $P$ is understood.
\end{definition}
For Markov chains, iterating this inequality gives exponential convergence. 
We quantify the time it takes the Markov chain to mix by the \vocab{mixing time},
\begin{align*}
    \tmix(P,\nu,\ep) &= \min \set{t}{\TV(\nu P^t,\mu)\le \ep}, &
    \tmix(P,\ep) &= \sup_{\nu} \tmix(P,\nu,\ep).
\end{align*}
\begin{lemma}[Mixing time]
\label{l:tmix}
    Suppose $(\Om, P)$ is an ergodic, reversible Markov chain with stationary distribution $\mu$. If $P$ satisfies entropy contraction with constant $\rec$, then for all $\nu$ and $\ep>0$, 
    \begin{align*}
    \tmix(P,\nu, \ep) \le \ce{\rc{\rec}\pa{\log\log \max_{x\in \Om}\fc{\nu(x)}{\mu(x)} + \log\prc{2\ep^2} }}\\
    \tmix(P, \ep) \le \ce{\rc{\rec}\pa{\log\log \max_{x\in \Om}\fc{1}{\mu(x)} + \log\prc{2\ep^2} }}.
    \end{align*}
\end{lemma}
For instance, if $\rec = \Om\prc{n}$ and $\min_{x\in \Om}\mu(x) = e^{-\poly(n)}$, then we have optimal mixing, $\tmix(P,\ep)= O(n\log n)$. 

\begin{definition}
    Let $\mu$ be a probability measure on $\Om$ and $P$ be a Markov kernel from $\Om$ to $\Om'$. Define the \vocab{adjoint} (or \vocab{time-reversal}) of $P$ with respect to the distribution $\mu$ to be the kernel $P^*_\mu$ from $\Om'$ to $\Om$, where $P^*_\mu(Y,\cdot)$ is the regular conditional distribution of $X|Y$ when $X\sim \mu$ and $Y\sim P(X,\cdot)$. 
\end{definition}
The adjoint can be interpreted as a posterior distribution of $X$ where $\mu$ is the prior and the output $Y$ of the channel $P$ is observed.
\begin{proposition}[Equivalent criterion for entropy contraction, {\cite[Proposition 26--27]{anari2024trickle}}]
\label{p:EC-equiv}
    Suppose $\mu$ is a probability measure on $\Om$ and $P$ is a Markov kernel from $\Om$ to $\Om'$. 
    If $\nu\ll \mu$, then 
    \[
\KL(\nu\|\mu) - \KL(\nu P \|\mu P) 
= \E_{Y \sim \mu P} \ba{\Ent_{X\sim (P_\mu^*)(Y,\cdot)}\ba{\dd{\nu}{\mu}(X)}}.
    \]
    Hence, $P$ satisfies $(1-\rh)$-entropy contraction with respect to $\mu$ iff for all non-negative functions $f$ on $\Om$,
    \[
\E_{Y \sim \mu P} \ba{\Ent_{X\sim (P_\mu^*)(Y,\cdot)}[f(X)]}
\ge \rh \cdot \Ent_\mu[f].
    \]
\end{proposition}
\begin{lemma}[Concavity of entropy decrement, {\cite[Lemma 29]{anari2024trickle}}]
\label{l:concave-entropy-dec}
    Let $P$ be a Markov kernel from $\Om$ to $\Om'$ and $f$ be a non-negative function on $\Om$. If $\mu = \E_{\te\sim \pi}[\mu_\te]$ is a mixture distribution, then 
    \[
\E_{Y\sim \mu P}\ba{\Ent_{X\sim (P_\mu)^*(Y,\cdot)}[f(X)]}
\ge \E_{\te\sim \pi} \ba{\E_{Y\sim \mu_\te P} \ba{\Ent_{X\sim (P_{\mu_\te})^*(Y,\cdot)}[f(X)]}}.
    \]
\end{lemma}
We first introduce the down-up walk for distributions on $\slice nk$. 
\begin{definition}
    Let $\mu$ be a distribution on $\slice nk$. Define the \vocab{down-up walk} as the Markov chain defined by the composition of the \vocab{down} and \vocab{up operators} $P_\mu = D_{k\to k-1}U_{k-1\to k,\mu}$, defined below.
    \begin{enumerate}
        \item The down operator $D_{k\to k-1}$ is defined by
        \[
D_{k\to k-1} (S,T) = \begin{cases}
    \rc k, & \text{if $T= S\bs \{i\}$ for some $i\in S$}\\
    0,&\text{otherwise.}
\end{cases}
        \]
    \item The up operator $U_{k-1\to k, \mu}$ is defined as the adjoint of $D_{k\to k-1}$ with respect to $\mu$. Explicitly, 
    \[
U_{k-1\to k,\mu}(T,S) =  \begin{cases}
    \fc{\mu(S)}{\sum_{j\in [n]\bs T}\mu(T\cup \{j\})}, & \text{if $S = T\cup \{i\}$ for some $i\in [n]\bs T$}\\
    0,&\text{otherwise.}
\end{cases}
    \]
    \end{enumerate}
\end{definition}
We recall the polarized walk for a distribution on $\{\pm 1\}^n$ given in \cite{anari2024trickle}. It can be defined in terms of the down-up walk after homogenizing.
\begin{definition}[{\cite[Definition 47, 84]{anari2024trickle}}]
\label{d:pw}
    Let $\mu$ be a distribution on $\calP([n])$. Define the polar homogenization kernel $\Pi$ from $\calP([n])$ to $\binom{[n]\sqcup [n]}{n}$, by taking $S_1$ to $S_1\sqcup S_2$ where $S_2\in \binom{[n]}{n-|S_1|}$ is chosen uniformly at random. Note that its adjoint $\Pi^*$ takes $S_1\sqcup S_2$ (where $S_1,S_2$ are subsets of the disjoint copies of $[n]$) to $S_1$.

    Define the \vocab{polarized walk} for $\mu$ as the Markov chain defined as the composition of the down and up operators $P_\nu\pol = D\pol U_\nu \pol$, defined below. These correspond to the down and up operators on the polar homogenized distribution.
    \begin{enumerate}
        \item The down operator $D\pol$ is defined as $D\pol = \Pi D_{n\to n-1} \Pi^*$. Explicitly, 
        \begin{align}
            D\pol(S,T) =\begin{cases}
                \rc{n} , & \text{if $T=S\bs \{i\}$ for some $i\in S$,}\\
                \fc{n-|S|}{n}, & \text{if $T=S$,}\\
                0,&\text{otherwise.}
            \end{cases}
        \end{align}
        \item The up operator $U_\mu \pol$ is defined as $U\pol = \Pi U_{n-1\to n,\mu}\Pi^*$. Explicitly,
        \begin{align}
            U_\nu\pol(T,S) = \begin{cases}
                \rc{n}\cdot \fc{\mu(S)}{\mu D\pol(T)}, & \text{if $S=T\cup \{i\}$ for some $i\in [n]\bs T$}\\
                \fc{n-|T|}{n} \cdot \fc{\mu(S)}{\mu D\pol(T)}, & \text{if $S=T$,}\\
                0, &\text{otherwise.}
            \end{cases}
        \end{align}
    \end{enumerate}
    For a distribution $\mu$ on $\{\pm 1\}^n$, define the polarized walk centered at $x_0$ by the above with the bijection $\{\pm 1\}^n\to \calP([n])$ defined by 
    $x\mapsto x_+:= \set{i}{x_i \ne (x_0)_i}$. 
    Let $\rh_{\text{EC}}(\mu)$ denote the maximum $\rh$ for which $D\pol$ satisfies $(1-\rh)$-entropy contraction with respect to $\mu$.\footnote{Note that the dependence on $x_0$ is implicit in the notation. When $\mu$ is supported on the wedge $B_r(x_0)$, we take $x_0$ as the center for the polarized walk.}
\end{definition}
Note that the down operators do not depend on $\mu$, but the up operators do. 
Moreover, $U_\mu\pol$ is the adjoint of $D\pol$ with respect to $\mu$. Also, note that by the data processing inequality, $P_{\mu}\pol$ has at least as much entropy contraction as $D\pol$ with respect to $\mu$.

We also require properties about measures on $\{0,1\}^n$ related to the geometry of their generating polynomials. We refer to \cite{borcea2009negative} and \cite{anari2018log} for background.
\begin{definition}\label{d:gen-poly}
    Let $\mu$ be a distribution on $\{0,1\}^n$. The \vocab{generating polynomial} of $\mu$ is the polynomial in $z=(z_1,\ldots, z_n)$ given by 
    \[
p_\mu(z) = \sum_{\si\in \{0,1\}^n} \mu(z) z^\si,
    \]
    where $z^\si:=\prod_{i=1}^n z_i^{\si_i}$.
\end{definition}
We will also write $z^S = \prod_{i\in S} z_i$ for $S\subeq [n]$. We write $\pl_{i}:=\pl_{z_i}$. 
\begin{definition}\label{d:rs-lc}
    Let $p\in \R[z_1,\ldots, z_n]$. 
    \begin{enumerate}[itemsep=0.1em]
        \item We say that $p$ is \vocab{real stable} if $p$ has no zero in the upper-half plane $\bbH^n$, where $\bbH = \set{z\in \C}{\Im z>0}$. 
        A distribution $\mu$ on $\{0,1\}^n$ is \vocab{strongly Rayleigh} if $p_\mu$ is real stable.
        \item We say that $p$ is \vocab{log-concave} if $\log p$ is concave on $\R_{>0}^n$. We consider 0 to be log-concave.
        \item We say that $p$ is \vocab{strongly log-concave} if for all $k\ge 0$ and $i_1,\ldots, i_k\in [n]$, $\pl_{i_1}\cdots \pl_{i_k}p$ is log-concave.
        \item We say that $p$ is \vocab{completely log-concave} if for all $k\ge 0$ and $V\in \R_{\ge 0}^{n\times k}$,
        \[
D_Vp(z):= \pa{\prodo jk \sumo in V_{ij}\pl_{i}}p(z)
        \]
        is non-negative and log-concave.
    \end{enumerate}
     A distribution $\mu$ on $\{0,1\}^n$ is (strongly/completely) log-concave if $p_\mu$ is (strongly/completely) log-concave.
\end{definition}
\begin{proposition}
\label{p:rs-lc}
The following implications hold for a polynomial $p\in \R[z_1,\ldots, z_n]$:
\begin{align*}
    \textup{Real stable} \implies
    \textup{Completely log-concave} \implies
    \textup{Strongly log-concave} \implies
    \textup{Log-concave}
\end{align*}
If $p$ is homogeneous, $p$ is completely log-concave iff it is strongly log-concave. 
\end{proposition}
\begin{proposition}\label{p:preserve}
    The following operations preserve real stability and (strong/complete) log-concavity.\footnote{(s/c) log-concave is to be read as strongly log-concave or completely log-concave or log-concave, respectively.}
    \begin{enumerate}[itemsep=0.1em]
        \item Product: If $p,q$ are real stable (respectively, log-concave) then $p\cdot q$ is real stable (log-concave).
        \item Symmetrization: If $p(z_1,\ldots, z_n)$ is real stable ((s/c) log-concave), then $p(z_1,z_1,z_3,\ldots, z_n)$ is real-stable ((s/c) log-concave).
        \item 
        Specialization: If $p(z_1,\ldots, z_n)$ is real stable ((s/c) log-concave), then $p(a,z_1,\ldots, z_n)$ is real-stable ((s/c) log-concave) for $a\in \R$ ($a\ge 0$).
        \item External field: If $p(z_1,\ldots, z_n)$ is real stable ((s/c) log-concave), then so is $q(z_1,\ldots, z_n) = p(w_1 z_1,\ldots, w_n z_n)$ for any $w_1,\ldots, w_n\ge 0$.
        \item Differentiation: If $p$ is real stable or strongly or completely log-concave, then so is $\pl_{i} p$. 
        \item Homogenization: If $p$ is real stable, then $p\homg(z_1,\ldots, z_{n+1}) = z_{n+1}^{\deg p}p\pa{\fc{z_1}{z_{n+1}},\cdots, \fc{z_n}{z_{n+1}}}$ is real stable. 
        \item Polarization (\cite[Proposition 43]{anari2024trickle}, \cite{branden2020lorentzian}): Suppose $p(z_1,\ldots, z_n) = \sum_{\al \in \N_0^n} c_\al z^\al$ is homogeneous and strongly log-concave with degree at most $\ka_i$ in $z_i$. Then so is its polarization with degree $\ka$, $\polar_\ka(p)$, defined by 
        \[
        \polar_\ka(p) = \sum_{\al \in \N_0^n} c_\al \prodo in \rc{\binom{\ka_i}{\al_i}} e_{\al_i}(z_{i1},\ldots, z_{i\ka_i})
        \]
        as a polynomial in $z_{ij}$, $i\in [n]$, $j\in [\ka_i]$, where $e_k(y_1,\ldots, y_n) = \sum_{S\in \slice{n}{k}} y^S$ is the $k$th elementary symmetric polynomial.
    \end{enumerate}
\end{proposition}

\begin{lemma}\label{l:esp}
The following hold.
\begin{enumerate}
    \item The elementary symmetric polynomial $e_k(z_1,\ldots, z_n) = \sum_{S\in \binom{[n]}{k}} z^S$ is real stable.
    \item 
    For any matroid $M$, $p_M(z) = \sum_{B\in M}z^B$ is completely log-concave \cite[Theorem 4.2]{anari2018log}, and its homogenization
    \[
    p_M\homg (z_1,\ldots, z_n,y) 
    = \sum_{B\in M}z^B y^{n - |B|}
    \]
    is also completely log-concave \cite[Theorem 4.1]{anari2024log}.
\end{enumerate}
\end{lemma}
This allows us to show the following.
\begin{lemma}
\label{l:prod-sr}
    Consider a restricted product distribution on $\{0,1\}^n$, $\mu = \mu_h|_{\Om}$.
\begin{enumerate}
    \item If $\Om = \slice{n}{k}$, then $\mu$ is strongly Rayleigh.
    \item If $\Om = \wdg{n}{k}$, then $\mu$ is completely log-concave.
\end{enumerate}
\end{lemma}
\begin{prf}
Let $w_i = e^{h_i}$.
\ppart{Slice $\Om = \slice{n}{k}$}
The generating polynomial of $\mu$ is $p(z_1,\ldots, z_n) \propto \sum_{S\in \slice nk} \prod_{i\in S}w_i z_i$. 
This is an elementary symmetric polynomial under the substitution $z_i\mapsfrom w_iz_i$ and rescaling, hence real stable by \Cref{l:esp} and \Cref{p:preserve}.

\ppart{Wedge $\Om = \wdg{n}{k}$}
    The generating polynomial of $\mu$ is $p(z_1,\ldots, z_n) \propto \sum_{S\in \wdg nk} \prod_{i\in S}w_i z_i$. 
    This is a matroid polynomial under the substitution $z_i\mapsfrom w_iz_i$ and rescaling, hence completely log-concave by \Cref{l:esp} and  \Cref{p:preserve}.
\end{prf}

Real stability and log-concavity give some useful facts about the corresponding measures.
\begin{lemma}[Negative correlation for strongly Rayleigh measures, \cite{borcea2009negative}]
\label{l:nd-sr}
    If $\mu$ on $\{0,1\}^n$ is strongly Rayleigh, then it satisfies negative correlation: for any $i\ne j$,
    \[
\mu(\si_i=\si_j=1)\le \mu(\si_i=1)\mu(\si_j=1).
    \]
\end{lemma}
\begin{lemma}[Ultra-log-concavity of rank sequence, \cite{gurvits2009multivariate}]
\label{l:ulc}
If $p(z_1,\ldots, z_n) = \sum_{\si} c_\si x^\si$, $c_\si\ge 0$ is strongly log-concave, then letting $a_k = \sum_{\ve{\si}_1 = k}c_\si$, the sequence $(a_0,\ldots, a_n)$ is ultra-log-concave, that is, 
\[
\pf{a_k}{\binom nk}^2 \ge \fc{a_{k-1}}{\binom n{k-1}}\cdot \fc{a_{k+1}}{\binom n{k+1}},
\]
and hence log-concave, $a_k^2 \ge a_{k-1}a_{k+1}$.
\end{lemma}
\begin{lemma}[Entropy contraction for log-concave distributions, \cite{cryan2019modified}]
\label{l:ec-lc}
    Let $\mu$ be a log-concave distribution on $\slice nk$. Then $D_{n\to n-1}$ satisfies $(1-\rc k)$-entropy contraction with respect to $\mu$: If $\nu\ll \mu$, then 
    \[
\KL(\nu D_{k\to k-1}\| \mu D_{k\to k-1}) 
\le \KL(\nu\|\mu).
    \]
\end{lemma}

We also introduce the 
functional inequalities for Langevin diffusion on $\R^n$. Given a probability density $\mu(x) \propto e^{-U(x)}$ with $U\in C^1$, define the generator by $\sL f = -\an{\gd f, \gd U} + \De f$ and the Dirichlet form by 
\[
\sE(f,g) = 
\an{f, \sL g}_{\mu} = 
\int_{\R^n} \an{\gd f, \gd g}d\mu. 
\]

\begin{theorem}[Bakry--\'Emery, \cite{bakry2006diffusions}]
\label{thm:gaussian-poincare}
    Suppose $\mu$ is $\al$-log-concave, that is, $\mu(x)\propto e^{-U(x)}$ where $U$ is $\al$-strongly convex. Then $\mu$ satisfies a log-Sobolev and hence Poincaré inequality with constant $\al$, that is,
    \begin{align*}
    \forall f&:\R^n\to \R_{\ge0},& \Ent_\mu(f)&\le \rc{2\al} \sE(f,\log f)=\rc{2\al} \int_{\R^n} \fc{\ve{\gd f}^2}{f} \,d\mu\\
        \forall f&:\R^n\to \R,& \Var_\mu(f)&\le \rc{\al} \sE(f,f) =\rc{\al} \int_{\R^n} \ve{\gd f}^2 \,d\mu.
    \end{align*}
    In particular, $\calN(0,\si^2I)$ satisfies a log-Sobolev and hence Poincaré inequality with constant $\rc{\si^2}$. 
\end{theorem}

\subsection{Entropy contraction for Ising models on the slice/wedge} 
\label{s:lsi-1}

We use the needle decomposition into rank-1 Ising models by \cite{eldan2022spectral}, then decompose a rank-1 into product measures using Hubbard--Stratonovich \cite{hubbard1959calculation} (see also \cite{bauerschmidt2019very}, \cite{koehler2022sampling}), and use the localization framework of \cite{CE25} to derive the functional inequality. 
We follow \cite{anari2024trickle} which applies the localization framework to the polarized walk.
We draw inspiration from \cite{huang2024sampling,huang2025weak} who show a weak Poincaré inequality for localized distributions for the spherical $p$-spin model. Whereas they proceed by showing a Poincaré inequality on spherical caps, we show entropy contraction on wedges. Our proof shares some features with \cite{mikulincer2025fastmixingisingmodels}, who also use Hubbard--Stratonovich and Gaussian approximation, and prove a similar bound on the covariance matrix in their setting.

\begin{lemma}[Entropy contraction for SK on wedges]\label{l:ls-sk-wedge}
    Let $A\sim \GOE(n)$ and $0<\be<\rc{2\sqrt 2}$. 
    There is a constant $c_1$ such that for small enough $\ep$, with probability $1-e^{-\Om(n)}$, 
    \[
\rec(\mu_{A,h}|_{\ball{\ep n}{x_0}}) \ge \rc n\exp\pa{-4\pa{1+ \fc1{1-2\pa{\be + \rc{2\sqrt 2} + \be c_1\sqrt{h(4\ep)}}^2}}}.
    \]
\end{lemma}
Here $h(\ep) = -\ep \log \ep - (1-\ep)\log (1-\ep)$ is the binary entropy function. 
We first understand product distributions restricted to wedges, which will be the basic components in our decomposition.

\begin{lemma}[Covariance bound for product distributions on the slice/wedge]
\label{l:prod-dist-cov}
    Consider a product distribution on $\{0,1\}^n$ restricted to the slice or wedge, $\mu = \mu_h|_{\Om}$ where $\Om = \slice{n}{k}$ or $\wdg{n}{k}$. 
    Let $p_i = \an{\si_i}_\mu$. 
    Then 
        \[
    \Cov(\mu) \preceq 2 \cdot \diag(p_i(1-p_i)).
    \]
\end{lemma}
The best bound we can obtain on the operator norm of the covariance is $\opnorm{\Cov(\mu)}\le \rc2$, which by itself is not useful. This lemma gives us the useful fact that the covariance cannot have many large eigenvalues. As a corollary, $\Tr[\Cov(\mu)]\le 2k$. 
\begin{prf}
\ppart{Slice $\Om = \slice{n}{k}$}
Since $\mu$ is strongly Rayleigh by \Cref{l:prod-sr}, we have negative correlation (\Cref{l:nd-sr}): for any $i\ne j$,
    \[
\mu(\si_i=\si_j=1) \le 
\mu(\si_i=1) \mu(\si_j=1).
    \]
    Hence, the off-diagonal entries in $\Cov_{\mu_{h}}(\si)$ are all $\le 0$. Let 
    $D=\diag(\an{\si}_{\mu})\diag(1-\an{\si}_{\mu})$, $C=\Cov(\mu)$,  and $M = 2D-C$. All entries of $M$ are non-negative. 
    By linearity of expectation, $\sumo jn \mu(\si_j=1) =  \sumo jn \E\si_j = k$.
    Hence for all $i$, 
    \begin{align*}
        \sum_{j\ne i}
        M_{ij} &= 
        \sum_{j\ne i}
        [\mu(\si_i=1)\mu(\si_j=1) - \mu(\si_i=1, \si_j=1)]
        = \mu(\si_{i}=1)\sum_{j\ne i}
        [\mu(\si_j=1) - \mu(\si_j=1|\si_i=1)]\\
        &= \mu(\si_{i}=1)[(k-\mu(\si_i=1)) - (k-1)]
        = \mu(\si_{i}=1) (1-\mu(\si_{i}=1)) 
        = M_{ii}.
    \end{align*}
    Therefore, by Gerschgorin's disc theorem, all eigenvalues of $M$ are non-negative, so $C\preceq D$.

\ppart{Wedge $\Om = \wdg{n}{k}$}
    We have
    \begin{align*}
    \sum_{j\ne i} M_{ij} &= 
    \mu(\si_{i}=1)\sum_{j\ne i}
        [\mu(\si_j=1) - \mu(\si_j=1|\si_i=1)]
    = \mu(\si_i=1)\mu(\si_i=0) \sum_{j\ne i} 
        [\mu(\si_j=1|\si_i=0) - \mu(\si_j=1|\si_i=1)]\\
    &= \mu(\si_i=1)\mu(\si_i=0) 
        \ba{ \sum_{j\ne i}\E[\si_j|\si_i=0]
        - \E[\si_j|\si_i=1]}.
    \end{align*}
    Let $a_m = \mu(\sum_{j\ne i}\si_j = m |\si_i=0)$.
    Then 
    \begin{align*}
        \sum_{j\ne i}\E[\si_j|\si_i=0]
        - \E[\si_j|\si_i=1]
        &= \fc{a_1+\cdots +a_kk}{a_0+\cdots +a_k} - \fc{a_1+\cdots +a_{k-1}(k-1)}{a_0+\cdots+a_{k-1}} = \fc{a_k(kS-T)}{S(S+a_k)}
    \end{align*}
    where $S = \sumz m{k-1}a_m$ and $T=\sumo m{k-1}ma_m$. We claim this is $\le 1$. Indeed, this is equivalent to
    \[
a_k (a_{k-2}+\cdots + (k-1)a_0)\le (a_0+\cdots +a_{k-1})^2.
    \]
    Note 
    $\mu_{h}|_{\Om\cap \{\si_i=0\}}$ is completely log-concave by \Cref{l:prod-sr}. 
    Hence, the sequence $(a_0,\ldots, a_k)$ is log-concave (\Cref{l:ulc}), so $a_ia_j\ge a_{i-1}a_{j+1}$ for any $0<i\le j<k$. Using this repeatedly shows the desired inequality, and we also obtain $\sum_{j\ne i}M_{ij}\le M_{ii}$.
\end{prf}
\pref{l:prod-dist-cov} allows us the following calculation. For any $x$ such that $xx^{\sT} \preceq A$, for $\Om = \ball{k}{x_0}$ (noting a factor of 4 from rescaling from $\{0,1\}^n$ to $\{\pm 1\}^n$),
\begin{align}
\nonumber
    \Var_{\mu_{h}|_{\Om}}(\an{x,\si})
&\le 4\an{x,2\diag\pa{\mu(\si_i=1)\mu(\si_i=-1)} x}\\
&\le
\sumo in 
(8\mu(\si_i=1)\wedge 2) x_i^2 
\le 2\max_{i_1<\cdots <i_{4k}}
\sumo j{4k} x_{i_j}^2
\le  2\max_{|S|=4k} \opnorm{A_{S\times S}}^2 .
\label{e:var-by-AS}
\end{align}
We would like to rule out the bad example where $A$ is concentrated on a subset of coordinates $S$, so there is no reduction in the interaction strength of $A$ on $S$.
For this, we show that with high probability, we can bound the operator norms by all submatrices of size $\ep n \times \ep n$ by some function of $\ep$ that is converging to 0 (which is in fact $O\pa{\ep \log \prc{\ep}}$. The idea that the worst subsystem is only a small factor worse is also used by \cite{adhikari2024spectral} for $p$-spin systems\footnote{This will not give sharp bounds, but this is fine because we can always run for longer time to concentrate the localized distribution on smaller slices.}.
\begin{theorem}[{\cite[Corollary 7.3.2]{vershynin2018high}}]
\label{t:GOE-opnorm}
    Let $A\sim \GOE(n)$. For some constant $c$,
    \[
P(\opnorm{A} \ge 2+t)\le 2e^{-cnt^2}.
    \]
\end{theorem}
\begin{lemma}[High-probability bound on operator norm of small submatrices]
\label{l:AS}
There is a constant $c_1$ such that the following holds. 
    Let $A\sim \GOE(n)$. Then 
    \[
P\pa{\forall |S|\le \ep n, \opnorm{A_{S\times S}}\le c_1\sqrt{h(\ep)}} \ge 1-2e^{-h(\ep) n}.
\]
\end{lemma}
\begin{prf}
We can apply \Cref{t:GOE-opnorm} with a union bound over all $e^{(h(\ep)+o(1))n}$ sets of $\le \ep n$ coordinates of $[n]$: 
\[
P\pa{\forall |S|\le \ep n, \opnorm{A_{S\times S}}\le (2+t)\sqrt\ep} \ge 1-2e^{(h(\ep) + o(1) - c\ep t^2)n}.
\]
We can choose $t=\Theta\pa{\sfc{h(\ep)}{\ep}}$ with appropriate constants to get the result.    
\end{prf}
The following gives a decomposition of Ising models into rank-1 Ising models (with interaction matrix bounded by the original) which preserve the expectation of a given test function $\ph$. This means that to prove a functional inequality (Poincaré or log-Sobolev) or the associated contraction of variance or entropy, it suffices to bound a functional inequality for these rank-1 Ising models, since no new variance/entropy is added by the decomposition (i.e., we have conservation of variance/entropy). 
\begin{theorem}[Needle decomposition for Ising models, \cite{eldan2022spectral}]
\label{t:needle}
Let $\Omega$ be any finite subset of $\set{v}{\ve{v}=R}$, and consider Ising models on $\Omega$.
    Given a test function $\ph:\Omega\to \R^n$ and an Ising model $\mu_{A,h}$ where $A$ is PSD, there exists a measure $\nu_{\ph}$ on $\R^{2n}$ such that 
    \[
\mu_{A,h} = \int_{\R^{2n}} \mu_{uu^{\sT} , v}
d\nu_\ph (u,v)
    \]
    where $\nu_\ph$-a.s., 
    $uu^{\sT} \preceq A$, $\E_{\mu_{A,h}}\ph = \E_{\mu_{uu^{\sT} , v}} \ph$. 
\end{theorem}
We note although \cite{eldan2022spectral} only consider $\Omega=\{\pm 1\}^n$, the proof does not rely on the particular structure of $\Omega$, other than the fact that the norm is constant on $\Omega$. 
While \cite{eldan2022spectral} prove the Poincaré inequality for mixture components through an influence argument, we make a further decomposition into a log-concave mixture of product distributions. This is given by the Hubbard--Stratonovich transform (\cite{hubbard1959calculation,koehler2022sampling}) (see \cite{bauerschmidt2019very} for an alternative log-concave decomposition).
(The localization process which gives the rank-1 decomposition can be thought of as an annealing step \cite[\S4]{CE25}.)
\begin{theorem}[Hubbard--Stratonovich transform, \cite{hubbard1959calculation,koehler2022sampling}]
\label{t:hs}
Let $A = XX^{\sT} $ be PSD, and let $\Omega$ be any finite subset of $\set{v}{\ve{v}=R}$. 
Then we have the following decomposition of $\mu_A$ into ``product measures'' restricted to $\Om$, where $\lm_{H}$ denotes Lebesgue measure on the hyperplane $H$:
    \begin{align*}
\mu_{A,h}|_{\Omega}(\si) 
&=\int_{\im(A)} p(u) 
\mu_{X^{\sT} u + h}|_{\Om}(\si) 
d\lambda_{\im(A)}(u) \\
\text{where }p(u) &\propto \exp\pa{-\fc 12 \ve{u}^2}
\sum_{\si\in \Omega} \exp\pa{\an{
Xu+h, \si}
}.
\end{align*}
\end{theorem}

\begin{lemma}[Entropy contraction for polarized walk on product distribution on wedge]
\label{l:ec-pw-prod}
    A product distribution restricted to the wedge $\Om = B_r(x_0)$, $\mu_{v}|_{\Om}(x) \propto \exp(\an{v,x})$ on $\Om$, satisfies $\rec(\mu_v|_{\Om})\ge \rc n$.
\end{lemma}
\begin{prf}
    Associating $\mu:= \mu_v|_{\Om}$ with a measure on $\calP([n])$ via the map $x\mapsto x_+:=\set{i}{x_i\ne (x_0)_i}$
    and letting $w_i = e^{-2v_i(x_0)_i}$, 
    the generating polynomial of $\mu_v$ is $g(z_1,\ldots, z_n) \propto  \sum_{|S|\le r} \prod_{i\in S}w_iz_i$. 
    This is a matroid polynomial under the substitution $z_i\mapsfrom w_iz_i$. 
    By \Cref{l:esp} and preservation of complete log-concavity under external fields (\Cref{p:preserve}), its homogenization
    \[
p\homg(z_1,\ldots, z_n, y) \propto 
\sum_{|S|\le r} \prod_{i\in S}w_iz_i
\cdot y^{n-|S|}.
    \]
    is completely (and hence strongly) log-concave.
    By \Cref{p:preserve}, 
    the polarization (with degree $(1,\ldots, 1,n)$)
    \[
\polar p\homg(z_1,\ldots, z_n, y_1,\ldots, y_n) \propto \sum_{|S|\le r} \prod_{i\in S}w_iz_i
\cdot \rc{\binom n{|S|}} e_{n-|S|}(y_1,\ldots, y_n)
    \]
    is strongly log-concave.
    The corresponding measure on $\binom{[n]\sqcup [n]}{n}$ is exactly the polar homogenization $\mu \Pi$ (\Cref{d:pw}).      
    Therefore, by \Cref{l:ec-lc}, $D_{n\to n-1}$ has $(1-\rc n)$-entropy contraction with respect to $\mu \Pi$, and this implies that $D\pol$ satisfies the same contraction with respect to $\mu$ by noting that $D\pol$ is defined so as to be the ``projection'' of $D_{n\to n-1}$ via $\Pi_\mu^*$: for any $\nu\ll \mu$,
    \begin{align*}
    \KL(\nu D\pol \| \mu D\pol) 
    &\le 
    \KL(\nu \Pi D_{n\to n-1} \| \mu \Pi D_{n\to n-1})\\
    &\le \pa{1-\rc n}\KL(\nu \Pi \|\mu \Pi) = \pa{1-\rc n} \KL(\nu\|\mu). \qedhere
    \end{align*}
    For the last equality, we note that the kernel $\Pi$ applied to distinct points has disjoint supports.
\end{prf}

\begin{theorem}[Entropy contraction from covariance bound under localization, {cf. \cite[Theorem 86]{anari2024trickle}}]
\label{t:sl-lsi}
    Let $\nu$ be a measure on 
    $B_r(x_0)$
    and $J$ be a PSD $n\times n$ matrix. For any $0\le \lm\le 1$ and every $v\in \R^n$, consider the probability measure defined by 
    \[
    \dd{\mu_{\lm,v}}{\nu}(x) \propto \exp\pa{-\fc{\lm}2 \an{x,Jx}+\an{v,x}}.
    \]
    Suppose that for some $\al:[0,1]\to \R_+$,
    \[
\opnorm{\Cov(\mu_{\lm,v})}\le \al(\lm), \quad \forall \lm\in [0,1],\, v\in \R^n
    \]
    and for some $\ep>0$, 
    \[
\rec(\mu_{1,v})\ge \ep, \quad \forall v\in \R^n.
    \]
    Then 
    \[
\rec(\nu)\ge \ep\exp\pa{-\opnorm{J}\int_0^1\al(\lm)\,d\lm}.
    \]
\end{theorem}
Note that \cite[Theorem 86]{anari2024trickle} is written for a specific measure, but the proof works for any measure as stated above. This is the analogue of \cite[Theorem 49]{CE25}, which is for Glauber dynamics.
\begin{prf}
    Stochastic localization with driving matrix $C_t = J^{1/2}$ gives the decomposition of measure
    \[
\nu = \int_{\R^n} \mu_{1,v} \,d\pi(v) 
    \]
    where $\pi$ is the distribution of the resulting $v$ in the stochastic localization process. Using the given operator norm bound on the covariance and \cite[Lemma 40]{CE25}, $\mu_{\lm,v}$ is $\al(\lm)\opnorm{J}$-entropically stable with respect to $\psi(x,y) = \rc 2 \ve{J^{1/2}(x-y)}^2$. Hence by \cite[Proposition 39]{CE25}, for a test function $f$, 
    \begin{align}
    \label{e:EC-decomp-1}
    \exp\pa{-\opnorm{J} \int_0^1 \al(\lm)\,d\lm} \Ent_{\nu}[f]
    \le 
\E_{v\sim \pi}[\Ent_{\mu_{1,v}}[f]].
    \end{align}
    Next, note by the equivalent definition of entropy contraction in \Cref{p:EC-equiv} and 
    \Cref{l:concave-entropy-dec},
    \begin{align}
    \label{e:EC-decomp-2}
\ep 
\E_{v\sim \pi} [\Ent_{\mu_{1,v}}[f]]
\le_{\text{\Cref{p:EC-equiv}}} 
\E_{v\sim \pi} 
\left[\E_{z\sim \mu_{1,v}D\pol} \left[\Ent_{U_{\mu_{1,v}}\pol(z,\cdot)}[f]\right]\right]
\le_{\text{\Cref{l:concave-entropy-dec}}} \E_{z\sim \nu D\pol} \Ent_{U_{\nu}\pol(z,\cdot)}[f].
    \end{align}
    Then combining \eqref{e:EC-decomp-1} and \eqref{e:EC-decomp-2} gives
    \[
\ep \exp\pa{-\opnorm{J} \int_0^1 \al(\lm)\,d\lm} \Ent_{\nu}[f]
    \le \E_{z\sim \nu D\pol} \Ent_{U_{\nu}\pol(z,\cdot)}[f]
    \]
and the conclusion follows using the equivalent definition of entropy contraction in \Cref{p:EC-equiv}.
\end{prf}

\begin{lemma}[Entropy contraction for Ising model on wedges via two-stage decomposition]
\label{lem:LS-2-stage}
    Let $\Om = \ball{r}{x_0}$ for some $x_0\in \{\pm 1\}^n$.
    Suppose $A\in \R^{n\times n}$ is a PSD matrix and $h\in \R^n$ such that for all $xx^{\sT} \preceq A$, $h\in \R^n$, if we define
    \begin{align}\label{e:1d-p}
p(u) &\propto  \exp\pa{-\fc{u^2}2}
\sum_{\si\in \Omega} \exp\pa{\an{
u\sqrt{\lm}x+h, \si}},
    \end{align}
    then $\Var(p)\le V(\lm)$, for some non-decreasing function $V:[0,1]\to \R_{\ge 0}$. 
    Then 
    \[
\rec(\mu_{A,h}|_{\ball{r}{x_0}})
\ge \rc n \exp\pa{-2\opnorm{A}\pa{1 + \int_0^1 V(\lm)\,d\lm}}
\ge \rc n  \exp\pa{-2\opnorm{A}\pa{1+V(1)}}.
    \]
\end{lemma}
\begin{prf}
    Given a test function $f:\Om\to \R^n$, 
    decompose $\mu_{A,h}|_{\Om}$ using \iffocs{Theorem 6.15}{\Cref{t:needle}},
    \begin{align}\label{e:LS-from-AS-1}
\mu_{A,h}|_{\Om} = \int_{\R^{2n}} \mu_{xx^{\sT} , v}|_{\Om}\,
d\nu_f (x,v),
    \end{align}
where $\nu_f$ is supported on $x$'s such that $xx^{\sT} \preceq A$. 
Suppose $x$ is such that $xx^{\sT} \preceq  A$, and consider $0\le \lm\le 1$. The Hubbard--Stratonovich transform (\iffocs{Theorem 6.16}{\Cref{t:hs}}) gives us 
    \begin{align}\label{e:LS-from-AS-2}
\mu_{\lm xx^{\sT} ,v}|_{\Om}(\si) 
&=\int_{\R} p_\lm(u) 
\mu_{u\sqrt{\lm} x + v}|_{\Om}(\si) 
du 
\end{align}
where $p=p_\lm$ is defined by \eqref{e:1d-p}. 
By assumption, $\Var(p)\le V(\lm)$.
Note 
\[\ddd u \E_{\mu_{u\sqrt{\lm}x+v}}\si = \Var_{\mu_{u\sqrt{\lm}x+v}}\pa{\an{\si,\sqrt\lm x}}\le 2\lm \ve{x}^2
\]
by \iffocs{Theorem 6.12}{\Cref{l:prod-dist-cov}}, so $\E_{\mu_{u\sqrt{\lm}x+v}}\si$ is $2\lm \ve{x}^2$-Lipschitz.
Therefore, $\E_{\mu_{u\sqrt \lm x+v}}\si$ has covariance at most $2\lm \ve{x}^2V(\lm)$. 
By covariance decomposition,
\begin{align}
\nonumber
    \Cov_{\mu_{\lm xx^{\sT} ,v}|_{\Om}}(\si) &= 
    \Cov_p(\E_{\mu_{u\sqrt{\lm}x + v}|_{\Om}}\si)
    + 
    \E_p \Cov_{\mu_{u\sqrt{\lm}x + v}|_{\Om}}(\si)\\
    &\le \pa{2\lm \ve{x}^2V(\lm) + 2}I_n.
\label{e:cov-decomp-1d}
\end{align}
Defining $(\cdot)_{\lm,v'}$ as in \Cref{t:sl-lsi} with $J=xx^{\sT} $, we have that for all $0\le \lm \le 1$ that 
\[
\opnorm{\Cov((\mu_{xx^{\sT} ,v}|_{\Om})_{1-\lm,v'})}  \le 2(1+\lm \ve{x}^2V(\lm)) \quad \forall \lm\in [0,1],\, v,v'\in \R^n.
\]
because \eqref{e:cov-decomp-1d} applies to $(\mu_{xx^{\sT} ,v}|_{\Om})_{1-\lm,v'}= \mu_{\lm xx^{\sT} , v+v'}|_{\Om}$. Hence by \Cref{t:sl-lsi}, noting $\ve{x}^2\le \opnorm{A}$, 
\begin{align}
\nonumber
\rec(\mu_{xx^{\sT} ,v}|_{\Om})
&\ge \inf_{v'\in \R^n}\rec(\mu_{v'}|_{\Om}) \cdot \exp\pa{-2\opnorm{A}\pa{1+\opnorm{A}\int_0^1 \lm V(\lm)\,d\lm}}\\
&\ge \rc n \exp\pa{-2\opnorm{A}\pa{1+\opnorm{A}\int_0^1 \lm V(\lm)\,d\lm}},
\label{e:mls-rank-1}
\end{align}
by \Cref{l:ec-pw-prod}. 
Now we proceed as in \Cref{t:sl-lsi}, noting that we have preservation of entropy in the decomposition: for $f\ge 0$ such that $\E_{\mu_{A,h}}f>0$, by \Cref{t:needle},
\begin{align*}
\Ent_{\mu_{A,h}}[f]
&=\E_{\mu_{A,h}}[f\log f] - \E_{\mu_{A,h}}[f] \log \ba{\E_{\mu_{A,h}}[f]}\\
&=\E_{(x,v)\sim \nu_f} \ba{\E_{\mu_{xx^\sT,v}}[f\log f] - \E_{\mu_{xx^\sT,v}}[f] \log \ba{\E_{\mu_{xx^\sT,v}}[f]}}
=\E_{(x,v)\sim \nu_f}\Ent_{\mu_{xx^\sT,v}}[f] .
\end{align*}
Hence, 
\begin{align*}
\Ent_{\mu_{A,h}}[f]
&= \E_{(x,v)\sim \nu_f}\ba{\Ent_{\mu_{xx^\sT,v}}[f]} \le_{\text{\Cref{p:EC-equiv}}} \rc{\inf_{xx^{\sT} \preceq A, v\in \R^n} \rec(\mu_{xx^{\sT} ,v}|_{\Om})}
\E_{(x,v)\sim \nu_f} \ba{\E_{z\sim \mu_{xx^{\sT} ,v}D\pol}\ba{\Ent_{ U_{\mu_{xx^{\sT} ,v}}\pol(z,\cdot)}[f]}}\\
&\le_{\text{\Cref{l:concave-entropy-dec}}} \rc{\inf_{xx^{\sT} \preceq A, v\in \R^n}\rec(\mu_{xx^{\sT} ,v}|_{\Om})}
\E_{z\sim \mu_{A,h}}\ba{\Ent_{ U_{\mu_{A,h}}\pol(z,\cdot)}[f]}.
\end{align*}
By \Cref{p:EC-equiv}, this shows that $\rec(\mu_{A,h}|_{\Om})\ge \inf_{xx^{\sT} \preceq A, v\in \R^n}\rec(\mu_{xx^{\sT} ,v}|_{\Om})$, which, combined with \eqref{e:mls-rank-1}, gives the desired result. 
\end{prf}

\begin{lemma}[Entropy contraction for Ising model on wedges]\label{l:LS-from-AS}
    Suppose $A$ is a PSD matrix such that for all $|S|\le 4\ep n$, $\opnorm{A_{S\times S}}\le C<\rc{\sqrt 2}$. Let $\Om = \wdg{n}{\ep n}$. Then for any $h\in \R^n$, 
    \[
\rec(\mu_{A,h}|_{\Om}) \ge \exp\pa{-2\opnorm{A}\pa{1+ \fc{\opnorm{A}}{1-2C^2}}}.
    \]
\end{lemma}
\begin{proof}
It suffices to compute an upper bound $V(\lm)$ in \pref{lem:LS-2-stage}. 
Let $xx^{\sT} \preceq A$, and define $p$ as in \eqref{e:1d-p}. 
We compute that the mixture distribution is log-concave by \eqref{e:var-by-AS} (note $\lm xx^{\sT} \preceq \lm A$),
\begin{align*}
    -\log p''(u)
    = 1 - \Var_{\mu_{u\sqrt\lm x+v}|_{\Om}}\pa{\an{\sqrt\lm x,\si}}
    \ge 1-2\lm^2 C^2.
\end{align*}
Therefore, $p$ has variance at most $V(\lm):=\rc{1-2\lm^2C^2}$, and plugging this into \pref{lem:LS-2-stage} gives the result.
\end{proof}

\begin{prf}[Proof of \Cref{l:ls-sk-wedge}]
    For $A\sim \GOE(n)$, consider $\be A + \ga I$, where $2\be < \ga < \rc{\sqrt2}$ (e.g. $\ga=\be+\rc{2\sqrt 2}$). With probability $\ge 1-e^{-\Om(n)}$, $\be A + \ga I$ is PSD and has operator norm bounded by $\sqrt 2$. 
Also with probability $1-e^{-\Om(n)}$, by \Cref{l:AS}, for $|S|\le 4\ep n$,
\[
\opnorm{(\be A+\ga I)_{S\times S}} \le \ga + 
\be c_1 \sqrt{h(4\ep)}
< \rc{\sqrt 2}
\]
for small enough $\ep$. 
The result than follows from \Cref{l:LS-from-AS}.
\end{prf}

\subsection{Attaining $\beta<\rc2$ using refined covariance bound}\label{s:lsi-2}
Extending entropy contraction to all $\be<\rc 2$ hinges on \Cref{l:cov-tight-off-diag}, which gives a better covariance bound when the uncertainty isn't limited to a small number of coordinates.

\begin{lemma}[Entropy contraction for SK on wedges, II]\label{l:ls-sk-wedge-ii}
    Fix $x_0\in \{\pm 1\}^n$. 
    Let $A\sim \GOE(n)$ and $0<\be<\rc{2}$. 
    There is a constant $\ep_0$ such that for all $\ep<\ep_0$, with probability $1-e^{-\Om(n)}$ (constants depending on $\be$), $\rec(\mu_{A,h}|_{B_{\ep n}(x_0)}) \ge \fc{\rh_\be}{n}$ where $\rh_\be>0$ is a constant depending only on $\be$. 
\end{lemma}

We have the following elementary fact about product distributions on slices.
\begin{lemma}[Monotonicity]
\label{l:monotonicity}
    Consider a product distribution $\mu_h$ on $\{0,1\}^n$, and let $\mu^{(k)}= \mu_h|_{\binom{[n]}{k}}$. 
    Let $p_j^k = \mu^{(k)}(\si_j = 1)$. Then for $k'\le k$, $p_j^{k'}\le p_j^k$. 
\end{lemma}
\begin{proof}
    It suffices to show that $p_j^k\le p_j^{k+1}$. This is equivalent to 
    \begin{align*}
        \fc{\sum_{\ve{\si_{-j}}_1= k-1} e^{h_j + \an{h_{-j}, \si_{-j}}}}{\sum_{\ve{\si_{-j}}_1= k-1} e^{h_j + \an{h_{-j}, \si_{-j}}}+ \sum_{\ve{\si_{-j}}_1\le k} e^{ \an{h_{-j}, \si_{-j}}}}
        &\le  \fc{\sum_{\ve{\si_{-j}}_1= k} e^{h_j + \an{h_{-j}, \si_{-j}}}}{\sum_{\ve{\si_{-j}}_1= k} e^{h_j + \an{h_{-j}, \si_{-j}}}+ \sum_{\ve{\si_{-j}}_1= k+1} e^{ \an{h_{-j}, \si_{-j}}}}\\
        \iff 
        \fc{\sum_{\ve{\si_{-j}}_1= k} e^{ \an{h_{-j}, \si_{-j}}}}{\sum_{\ve{\si_{-j}}_1= k-1}e^{h_j + \an{h_{-j}, \si_{-j}}}}
        &\ge \fc{\sum_{\ve{\si_{-j}}_1= k+1} e^{ \an{h_{-j}, \si_{-j}}}}{\sum_{\ve{\si_{-j}}_1= k}e^{h_j + \an{h_{-j}, \si_{-j}}}}\\
        \iff a_k^2 &\ge a_{k-1}a_{k+1}
    \end{align*}
    where $a_k=\mu_{h,-j}(\ve{\si}_1=k)$. This holds by log-concavity of $\mu_h$.
\end{proof}

We bound the differences between the ``effective uncertainty'' (sum of covariances of individual coordinates) for the unrestricted distribution (with appropriate tilts) and the distribution restricted to a slice or wedge. We will frequently need to switch between these quantities in the coming proofs.
\begin{lemma}
\label{l:mm'}
    Consider a product distribution $\mu_h$ on $\{0,1\}^n$. Let $f(p) = p(1-p)$.
    For $k\in [0,n]$, let $t(k)$ be such that $\E_{\mu_{h+t(k)\one}}$ be such that $\E_{\mu_{h+t(k)\one}} \ve{\si}_1 = k$, $q_i^{(k)} = \mu_{h+t(k)\one}(\si_i=1)$, and $n^{(k)} = \sumo in f(q_i^{(k)})$. (When $k=0$ or $n$, set $t=-\iy$, $t=\iy$, respectively.) For $k\in \{0,\ldots n\}$, let $p_i^{(k)} = \mu_h (\si_i = 1|\ve{\si}_1=k)$ and $m^{(k)} = \sumo in f(p_i^{(k)})$; 
    let $p_i^{(\le k)} = \mu_h(\si_i = 1 | \ve{\si}_1\le k)$ and $m^{(\le k)} = \sumo in f(p_i^{(\le k)})$. 
    Then the following hold.
\begin{enumerate}
    \item 
    We have for $k<l$ that 
    \begin{align*}
     |m^{(k)}-m^{(l)}| \le \sumo in \ab{f(p_i^{(k)}) - f(p_i^{(l)})} &\le 
     \sumo in \pa{\max_{p\in [p_i^{(k)}, p_i^{(l)}]} f(p) - \min_{p\in [p_i^{(k)}, p_i^{(l)}]}f(p)}\le 
     l-k\\
     |n^{(k)} - n^{(l)}| \le \sumo in \ab{f(q_i^{(k)}) - f(q_i^{(l)})} &\le 
     \sumo in \pa{\max_{p\in [q_i^{(k)}, q_i^{(l)}]} f(p) - \min_{p\in [q_i^{(k)}, q_i^{(l)}]}f(p)}\le 
     l-k.
    \end{align*}
    \item 
    Let $\bar k = \E_{\mu_h} \ve{\si}_1$. 
    We have 
    \begin{align*}
        |m^{(k)}-n^{(k)}| \le \sumo in \ab{f(p_i^{(k)}) - f(q_i^{(k)})} &= O(\sqrt {n^{(k)}}) = O(\sqrt {m^{(k)}})\\
        |m^{(\le k)}-m^{(k\wedge \fl{\bar k})}| \le \sumo in \ab{f(p_i^{(\le k)}) - f(p_i^{(k\wedge \fl{\bar k})})} &= O(\sqrt {n^{(k\wedge \fl{\bar k})}}) = O(\sqrt {m^{(\le k)}})
    \end{align*} 
\end{enumerate}
\end{lemma}
\begin{prf} The proof proceeds in two parts, where we split the second part further into two cases. 
\ppart{Part 1}
    For (1), note first that $t$ exists because $\E_{\mu_{h+t\one}} \ve{\si}_1$ is a continuous and monotonically increasing function tending to 0 and $n$ as $t\to -\iy$ and $t\to \iy$, respectively. (It is continuous and increasing because $\ddd t \E_{\mu_{h+t\one}}\ve{\si}_1 = \Cov_{\mu_{h+t\one}}(\ve{\si}_1)>0$.) 
    Also note that $p^{(k)}_i$ is monotonically increasing in $k$ by \pref{l:monotonicity}. Note that $f(p)$ is 1-Lipschitz on $[0,1]$. Hence 
    \[
\sumo in \ab{f(p_i^{(k)}) - 
    f(p_i^{(l)})}
    \le \sumo in \pa{\max_{p\in [p_i^{(k)}, p_i^{(l)}]} f(p) - \min_{p\in [p_i^{(k)}, p_i^{(l)}]}f(p)}
    \le \sumo in p_i^{(l)}-p_i^{(k)}
    = l-k.
    \]
    Similarly, for $k,l\in [0,n]$, we obtain the second inequality in the same way.

    \ppart{Part 2a}
    Now consider (2). In the following, all expectations are with respect to $\mu_{h+t(k)\one}$.  
    Note that $q_i^{(k)} = \E[\E[\si_i|\ve{\si}_1]]=\E[p_i^{(\ve{\si}_1)}]$. Then noting $f(y) = f(x) + (1-2x) (y-x) - (y-x)^2$, taking $x=\E[p_i^{(\ve{\si}_1)}]$ and $y=f(p_i^{(\ve{\si}_1)})$ and noting the linear term disappears under expectation,
    \begin{align*}
\ab{f(q_i^{(k)}) - f(p_i^{(k)})} 
= \ab{f(\E[p_i^{(\ve{\si}_1)}]) - f(p_i^{(k)})}
&= \ab{\E \ba{f(p_i^{(\ve{\si}_1)}) + \pa{p_i^{(\ve{\si}_1)} - \E[p_i^{(\ve{\si}_1)}]}^2} - f(p_i^{(k)})}\\
&= \E\ab{f(p_i^{(\ve{\si}_1)}) - f(p_i^{(k)})} +  \E \pa{p_i^{(\ve{\si}_1)} - \E[p_i^{(\ve{\si}_1)}]}^2\\
&= \E\ab{f(p_i^{(\ve{\si}_1)}) - f(p_i^{(k)})} + \rc2\E \pa{p_i^{(\ve{\si}_1)} - p_i^{(\ve{\si'}_1)}}^2\\
&\le\E\ab{f(p_i^{(\ve{\si}_1)}) - f(p_i^{(k)})} + \rc2\E \ab{p_i^{(\ve{\si}_1)} - p_i^{(\ve{\si'}_1)}}
    \end{align*}
    where $\si'$ is an independent copy of $\si$. 
    Therefore, by the first part, 
\begin{align*}
    \sumo in \ab{f(q_i^{(k)}) - f(p_i^{(k)})} 
    &\le \E \ba{\sumo in \ab{f(p_i^{(\ve{\si}_1)}) - f(p_i^{(k)})} + \rc2\ab{p_i^{(\ve{\si}_1)} - p_i^{(\ve{\si'}_1)}}}\\
    &\le \E \ab{\ve{\si}_1 - k} + \rc2\ab{\ve{\si}_1 - \ve{\si'}_1}
    \le 2 \pa{\E \pa{\ve{\si}_1 - k}^2}^{1/2}\\
    & = 2\pa{\sumo in \Var_{\mu_{h+t(k)\one}}(\si_i)}^{1/2} = 2 \pa{\sumo in f(q_i)}^{1/2}= 2\sqrt{n^{(k)}}.
\end{align*}
This shows that $|m^{(k)}-n^{(k)}|\le 2\sqrt{n^{(k)}}$. Thus, $\sqrt{n^{(k)}} = O(\sqrt{m^{(k)}})$ as $m^{(k)}\to \iy$, which shows the first equation in the second part.

\ppart{Part 2b}
Now consider $p_i^{(\le k)}$. First suppose $k \le \E_{\mu_h}\ve{\si}_1$. 
We first consider $\ve{\si}_1$ under $\mu_{h+t(k)\one}^{(\le k)}$.
Note by Cauchy-Schwarz
\begin{align*}
    \E_{\mu_{h+t(k)\one}}[k-\ve{\si}_1 | \ve{\si}_1\le k]
    &= \rc{\mu_{h+t(k)\one}(\ve{\si}_1\le k)} \E_{\mu_{h+t(k)\one}}[(k-\ve{\si}_1)\one_{\ve{\si}_1\le k}]\\
    &\le \rc{\mu_{h+t(k)\one}(\ve{\si}_1\le k)^{1/2}}
    \Var_{\mu_{h+t(k)\one}}(\ve{\si}_1)^{1/2} = 
    O(1)\cdot \sqrt{n^{(k)}}.
\end{align*}
where the bound $\mu_{h+t(k)\one}(\ve{\si}_1\le k)=\Om(1)$ holds by Berry-Esseen (or \Cref{t:llt}). 
Therefore, 
\begin{align*}
    k - O(\sqrt{n^{(k)}}) \le \E_{\mu_{h+t(k)\one}}[\ve{\si}_1 | \ve{\si}_1\le k] \le k 
\end{align*}
Now since $k \le \E_{\mu_h}\ve{\si}_1$, $t(k)\le 0$. We have that $\mu_{h+t(k)\one}^{(\le k)}(\si)\propto \mu_h^{(\le k)}(\si)e^{t(k)\ve{\si}_1}$, so $\ve{\si}_1$ under $\mu_h^{(\le k)}$ stochastically dominates $\ve{\si}_1$ under $\mu_{h+t(k)\one}^{(\le k)}$. Hence 
\begin{align*}
    k - O(\sqrt{n^{(k)}}) \le \E_{\mu_{h}}[\ve{\si}_1 | \ve{\si}_1\le k] \le k.
\end{align*}
As in part 1, by 1-Lipschitzness of $f$ and noting $p_i^{(k)}\ge p_i^{(\le k)}$ by monotonicity,
\begin{align}
\label{e:pk-plek}
    \sumo in \ab{f(p_i^{(\le k)}) - f(p_i^{(k)})} 
    &\le \sumo in (p_i^{(k)} - p_i^{(\le k)})
    = k -  \E_{\mu_{h}}[\ve{\si}_1 | \ve{\si}_1\le k] = O(\sqrt{n^{(k)}}). 
\end{align}
This shows $|m^{(\le k)}-m^{(k)}|= O(\sqrt{n^{(k)}})$ as $n^{(k)}\to \infty$. Combining with part 2a, we get $|m^{(\le k)}-n^{(k)}|= O(\sqrt{n^{(k)}})$, and hence this is also $O(\sqrt{m^{(\le k)}})$ as $m^{(\le k)}\to \iy$.

Now suppose $k> \bar k$. By monotonicity, $p_i^{(\le \fl{\bar k})}\le p_i^{(\le k)}\le q_i^{(\bar k)}$. Hence by 1-Lipschitzness of $f$,
\begin{align*}
    \sumo in \ab{f(p_i^{(\le k)}) - f(p_i^{(\fl{\bar k})})}
    &\le
    \sumo in \ab{f(p_i^{(\le k)}) - f(p_i^{(\le \fl{\bar k})})}
    + 
    \sumo in \ab{f(p_i^{(\le \fl{\bar k})}) - f(p_i^{(\fl{\bar k})})}\\
    &\le \sumo in (q_i^{(\bar k)} - p_i^{(\le \fl{\bar k})})
    + \sumo in \ab{f(p_i^{(\le \fl{\bar k})}) - f(p_i^{(\fl{\bar k})})}\\
    &\le \pa{\bar k - \E_{\mu_h}[\ve{\si}_1|\ve{\si}_1\le \bar k]} + \sumo in \ab{f(p_i^{(\le \fl{\bar k})}) - f(p_i^{(\fl{\bar k})})}\\
    &\le (\bar k - \fl{\bar k}) + O(\sqrt{n^{(\fl{\bar{k}})}}) + O(\sqrt{n^{(\fl{\bar{k}})}}) = O(\sqrt{n^{(\fl{\bar{k}})}})
\end{align*}
where 
both terms are bounded by \eqref{e:pk-plek} applied to $\fl{\bar k}$. Since $|m^{(\fl{\bar k})} - n^{(\fl{\bar k})}| = O(\sqrt{n^{(\fl{\bar k})}})$, this is also $O(\sqrt{m^{(\fl{\bar k})}})$, and hence $O(\sqrt{m^{(\le k)}})$.
\end{prf}

\begin{lemma}[Cross-ratio for product distribution]\label{l:cross-ratio}
    Consider a product distribution $\mu_h$ on $\{0,1\}^n$. Let $p_i = \mu_h(\si_i=1)$ and $m = \sumo in p_i (1-p_i)$. 
    Let $a_k = \mu_h (\ve{\si}_i=k)$. 
Then for any $k_1\le k_2$ with $k_i - \E_{\mu_h} \ve{\si}_1 = o(m)$ for $i=1, 2$, we have 
    \[
\fc{a_{k_1-1}a_{k_2+1}}{a_{k_1} a_{k_2}}
= 1 - \fc{k_2-k_1+1}{m} (1+o(1)).
    \]
\end{lemma}
In the proof we also show this is $1-\fc{k_2-k_1+1}{m^{(k)}}(1+o(1))$, where $m^{(k)}$ is defined as in \eqref{e:mk}.
\begin{prf}
Let $f(p)=p(1-p)$. 
Let $\mu_h^{(k)}:= \mu_h|_{\slice{n}{k}}$.
The proof hinges on a combinatorial interpretation. 
    We have
    \begin{align*}
\fc{a_{k_1-1}a_{k_2+1}}{a_{k_1} a_{k_2}}
&= \fc{\pa{\sum_{\ve{\si}_1=k_1-1} e^{\an{h,\si}}}\pa{\sum_{\ve{\si}_1=k_2+1} e^{\an{h,\si}}}}{\pa{\sum_{\ve{\si}_1=k_1} e^{\an{h,\si}}}\pa{\sum_{\ve{\si}_1=k_2}e^{\an{h,\si}}}}\\
    &= \E_{(\si_1,\si_2)\sim \mu_h^{(k_1)}\ot \mu_h^{(k_2)}} \fc{\set{(\si_1',\si_2')}{\si_1'\in \binom{[n]}{k_1-1}, \si_2'\in \binom{[n]}{k_2+1}, \si_1'+\si_2'=\si_1+\si_2}}{\set{(\si_1'',\si_2'')}{\si_1'' \in \binom{[n]}{k_1}, \si_2''\in \binom{[n]}{k_2}, \si_1''+\si_2''=\si_1+\si_2}}.
    \end{align*}
For $\si_1, \si_2\in \{0,1\}^{n}$, let $s_1=s_1(\si_1,\si_2):=\ab{\set{j}{(\si_1)_j=1, \,(\si_2)_j=0}}$ and $s=s(\si_1,\si_2):=\ab{\set{j}{(\si_1+\si_2)_j=1}}$. 
(If we think of $\si_1,\si_2$ as sets, then $s_1 = |\si_1\bs\si_2|$ and $s=|\si_1\triangle\si_2|$.) 
Let $r=k_2-k_1+1$. 
Note that for $(\si_1'',\si_2'')$ as in the denominator, $|\si_1''\triangle \si_2''| = |\si_1\triangle \si_2| = s = s_1+(k_2-k_1+s_1)=2s_1+r-1$; we must choose $\si_1''\bs \si_2''$ with cardinality $s_1$ from $\si_1\triangle \si_2$. In the numerator we choose it of size $s_1-1$. Hence the integrand equals
$\binom{2s_1+r-1}{s_1-1}\big/ \binom{2s_1+r-1}{s_1} = \fc{s_1}{s_1+r} = \fc{s-r+1}{s+r+1}$. 
We can sample from $\mu_h^{(k')}$ as follows. Choose $t(k')$ such that $\E_{\mu_{h+t(k')\one}}\ve{\si}_1=k'$. Now sample from the product distribution $\mu_{h+t(k')\one}$ and accept a sample only if $\ve{\si}_1=k'$. This formulation allows us to use concentration for sums of independent random variables. 
Let $p_{j}^{(k')} = \mu_{h}^{(k')}(\si_j=1)$ and 
\begin{align}\label{e:mk}
m^{(k')} = \sum_{j=1}^n f(p_{j}^{(k')})
\end{align}
be the effective uncertainty for $\mu_{h}^{(k')}$.
Let $q_{j}^{(k')} = \mu_{h+t(k')\one}(\si_j=1)$ and $n^{(k')} = \sum_{j=1}^n f(q_{j}^{(k')})
$ be the effective uncertainty for $\mu_{h+t(k')\one}$. 
By \Cref{t:llt}, with $\ga$ as the pdf of the standard normal, 
\[
\mu_{h+t(k')\one}\pa{\set{\si}{\ve{\si}_1=k'}}
\ge 
\rc{\sqrt{n^{(k')}}} \ga(0) - o\prc{\sqrt{n^{(k')}}}
= \Om\prc{\sqrt{n^{(k')}}}.
\]
Then 
$R(\si_1,\si_2) := \dd{\mu_h^{(k_1)}\ot \mu^{(k_2)}}{\mu_{h+t(k_1)\one}\ot \mu_{h+t(k_2)\one}}(\si_1,\si_2) \le 
C \sqrt{n^{(k_1)}n^{(k_2)}}$ for some $C$ for $n^{(k_1)}, n^{(k_2)}$ large enough.
By \Cref{l:mm'}, 
\allowdisplaybreaks
\begin{align*}
    |m-n^{(k_1)}| &\le \ab{k_1-\E_{\mu_h} \ve{\si}_1} = o(m)\\
    |n^{(k_1)} - n^{(k_2)}| &\le r-1\\
    |n^{(k_1)} - m^{(k_1)}| &= O\pa{\sqrt{m^{(k_1)}}}\\
    \implies\ab{m - m^{(k_1)}}&= O\pa{\sqrt{m^{(k_1)}}}
    + O\pa{|k_1 - \E_{\mu_h}\ve{\si}_1|}.
\end{align*}
Therefore, when $r=o(m^{(k_2)})$, 
we have $n^{(k_1)}, n^{(k_2)} = m^{(k_2)} (1+o(1))$.
Let $2m' = \sum_{j=1}^n q_j^{(k_1)}(1 - q_j^{(k_2)}) + (1-q_j^{(k_1)}) q_j^{(k_2)}$. 
Note $2m' = \E_{(\si_1,\si_2)\sim \mu_h^{(k_1)}\ot \mu_h^{(k_2)}} s$ and 
\allowdisplaybreaks
\begin{align*}
    2m' &= \sumo jn \ba{q_j^{(k_1)}(1-q_j^{(k_1)}) + q_j^{(k_1)}(q_j^{(k_2)}-q_j^{(k_1)})} + 
    \sumo jn \ba{(1-q_j^{(k_2)}) q_j^{(k_2)} + (q_j^{(k_2)}-q_j^{(k_1)}) q_j^{(k_2)}} \\
    &\in 
    \sumo jn f(q_j^{(k_1)}) + f(q_j^{(k_2)}) + \ba{-2|r-1|, 2|r-1|} = 
    m^{(k_1)} + m^{(k_2)} + \ba{-2|r-1|, 2|r-1|}
\end{align*}
by \Cref{l:mm'}.
We moreover have 
\[
\Var_{(\si_1,\si_2)\sim \mu_h^{(k_1)}\ot \mu_h^{(k_2)}}s = \sum_{j=1}^n f\pa{q_j^{(k_1)}(1 - q_j^{(k_2)}) + (1-q_j^{(k_1)}) q_j^{(k_2)}} \le 2m'.
\]
Let $I=\{r-1,r+1,\ldots,2\fl{m'  - \sqrt{m'}\log m'}-r-3\}$. If $s\nin I$, we have $\fc{2r}{s+r+1}\le \fc{r}{\fl{m' - \sqrt{m'}\log m'}}$.
Now by Bernstein's inequality as $m^{(k_2)}\to \iy$,
\allowdisplaybreaks
\begin{align*}
    &\E_{(\si_1,\si_2)\sim \mu_h^{(k_1)}\ot \mu_h^{(k_2)}}  \ba{\fc{s-r+1}{s+r+1}} = 
    \E_{(\si_1,\si_2)\sim \mu_h^{(k_1)}\ot \mu_h^{(k_2)}}  \ba{1-\fc{2r}{s+r+1}}
    \\
    &\ge 
    \pa{1-\fc{r}{\fl{m' - \sqrt{m'}\log m'}}} -4r \sum_{j\in I}\P_{(\si_1,\si_2)\sim \mu_h^{(k_1)}\ot \mu_h^{(k_2)}}(s\le j) \rc{(j+r+1)(j+r+3)}\\
    &\ge 
    \pa{1-\fc{r}{\fl{m' - \sqrt{m'}\log m'}}} -4r \sum_{j\in I}
     C \sqrt{n^{(k_1)}n^{(k_2)}}
    \P_{(\si_1,\si_2)\sim \mu_{h+t(k_1)\one}\ot \mu_{h+t(k_2)\one}}(s\le j) \rc{(j+r+1)(j+r+3)}\\
    &\ge 1 - \fc{r}{m'} (1+o(1)) - O(1)r \sum_{j\in I}
    m'
    \exp\pa{ - \fc{(j-2m')^2/2}{|j-2m'|/3+2m'}} \rc{(j+r+1)(j+r+3)}\\
    & \ge 1 - \fc{r}{m'} (1+o(1)) - O(1)r \sum_{j\in I}
    m'
    \exp\pa{ - \Om\pa{\fc{(j-2m')^2}{m'}}} \rc{(j+r)^2}\\
    &\ge 1 - \fc{r}{m'} (1+o(1)) - O(1)r 
    m^{\prime 2} e^{-\Om(\log^2 m')}\\
    &\ge 1-\fc{r}{m'} (1+o(1))
\end{align*}
by considering the sum as a Gaussian tail. For the upper bound, for $j>2m'$,
\begin{align*}
    \E_{(\si_1,\si_2)\sim \mu_h^{(k_1)}\ot \mu_h^{(k_2)}}  \ba{\fc{s-r+1}{s+r+1}}
    &= 1- \E_{(\si_1,\si_2)\sim \mu_h^{(k_1)}\ot \mu_h^{(k_2)}} \ba{\fc{2r}{s+r+1}}\\
    &\le 1- \P_{(\si_1,\si_2)\sim \mu_h^{(k_1)}\ot \mu_h^{(k_2)}}(s\le j) \fc{2r}{j+r+1}\\
    &\le 1- \pa{1- 4m'(1+o(1))
    \exp\pa{ - \fc{(j-2m')^2/2}{(j-2m')/3+2m'}}} \fc{2r}{j+r+1}.
\end{align*}
Choosing $j=2m'+\sqrt{m'}\log m'$, we get an upper bound of $1-\fc{r}{m'}(1-o(1))$. 
Finally, note $m'=m^{(k_2)}(1+o(1)) = m(1+o(1))$.
\end{prf}
The following computation is useful.
\begin{lemma}\label{l:E-m-evol}
    Let $\si(x) = \fc{e^x}{1+e^x}$. Let $h\in \R^n$, $E(t) = \sumo in \si(h_i+t)$ and $m(t) = \sumo in \si(h_i+t)(1-\si(h_i+t))$. Then 
    \begin{align*}
        E'(t) &= m(t)\\
        m'(t) &= \sumo in \si(h_i+t)(1-\si(h_i+t))(1-2\si(h_i+t)) \in [-m(t),m(t)].
    \end{align*}
\end{lemma}
\begin{prf}
    Immediate.
\end{prf}

We can now obtain tight control of the off-diagonal entries of the covariance matrix when the effective uncertainty is large.
\begin{lemma}[Covariance matrix of product distribution on slice/wedge under high effective uncertainty]\label{l:cov-tight-off-diag}
Consider a product distribution $\mu_h$ on $\{0,1\}^n$. 
\begin{enumerate}
    \item 
    Consider the restriction to a slice, $\mu^{(k)} := \mu_h|_{\slice{n}{k}}$. 
    Let $p_i^{(k)} = \mu^{(k)}(\si_i=1)$ and define the effective uncertainty as $m^{(k)} = \sumo in p_i^{(k)}(1-p_i^{(k)})$. 
    Then as $m^{(k)}\to \iy$, 
\begin{align*}
    \Cov_{\mu^{(k)}}(\si_i,\si_j) = -p_i^{(k)}(1-p_i^{(k)}) p_j^{(k)}(1-p_j^{(k)}) \cdot \rc{m^{(k)}}(1+o(1)).
\end{align*}
    \item 
    Consider the restriction to a wedge, $\mu^{(\le k)} := \mu_h|_{\wdg{n}{k}}$. 
    Let $p_i^{(\le k)} = \mu^{(\le k)}(\si_i=1)$ and define the effective uncertainty as $m^{(\le k)} = \sumo in p_i^{(\le k)}(1-p_i^{(\le k)})$. 
    Then as $m^{(\le k)}\to \iy$, there exists $0<c<1$ depending on $h$, $k$ such that for all $i\ne j$,
\begin{align*}
    \Cov_{\mu^{(\le k)}}(\si_i,\si_j) = -p_i^{(\le k)}(1-p_i^{(\le k)}) p_j^{(\le k)}(1-p_j^{(\le k)}) \cdot \rc{m^{(\le k)}}(c+o(1)).
\end{align*}
\end{enumerate}
\end{lemma}
\begin{proof}
Let $f(p)=p(1-p)$. 
For a distribution $\mu$ on $\{0,1\}^n$, let 
$\mu_{-i}$ denote the distribution of all coordinates except the $i$th one (for $\mu_h$, denote this by $\mu_{h,-i}$).
We use the notation of \Cref{l:mm'}, further using a subscript $-i$ to denote the relevant quantities for the distribution $\mu_{h,-i}$. 
Define the odds ratio by $\mathrm{odds}(p) = \fc{p}{1-p}$ for $0<p<1$. 
Let $a_{k'} = \mu_h (\ve{\si}_1=k')$ and $a_{k',-i} = \mu_{h,-i}(\ve{\si}_1=k')$. (Here $\si\in \{0,1\}^{[n]\bs \{i\}}$.)

\ppart{Part (1)}
Let $1\le r< k$. First note that we can approximate the following quotient of odds ratios with \Cref{l:cross-ratio}:
\begin{align}
\nonumber 
    \fc{\mathrm{odds}(\mu^{(k-r)}(\si_i=1))}{\mathrm{odds}(\mu^{(k)}(\si_i=1))} &= 
    \fc{\mu^{(k-r)}(\si_i=1)}{\mu^{(k-r)}(\si_i=0)} \Bigg/
    \fc{\mu^{(k)}(\si_i=1)}{\mu^{(k)}(\si_i=0)}\\
    &= 
    \fc{a_{k-r-1,-i} a_{k,-i}}{a_{k-r,-i} a_{k-1,-i}} = 1-\fc{r}{m_{-i}^{(k-1)}}(1+o(1))
\end{align}

Let $p_j^{(k')} = \mu^{(k')}(\si_j=1)$ and $m^{(k')} = \sum_{j=1}^r p_j^{(k')}(1-p_j^{(k')})$.
We also have the equation
\begin{align}\label{e:p-cvx-comb}
    p_j^{(k)} &= (1-p_i^{(k)}) p_{j,-i}^{(k)} + p_i^{(k)} p_{j,-i}^{(k-1)}.
\end{align}
Hence writing $ m^{(k')}= p_i^{(k')}(1-p_i^{(k')}) + \sum_{j\ne i} p_j^{(k')}(1-p_j^{(k')})$,
\[
\sum_{j\ne i} \min_{p\in [p_{j,-i}^{(k-1)}, p_{j,-i}^{(k)}]} f(p)
\le m^{(k)}\le 1 + \sum_{j\ne i} \max_{p\in [p_{j,-i}^{(k-1)}, p_{j,-i}^{(k)}]} f(p).
\]
By \Cref{l:mm'}(1), this interval has length at most 1 and contains $m_{-i}^{(k-1)}$ and $m_{-i}^{(k)}$, so that $m^{(k)}=m_{-i}^{(k-1)}+O(1)=m_{-i}^{(k)}+O(1)$. Hence 
\begin{align}\label{e:odds-ratio-r}
\fc{p_i^{(k-r)}/(1-p_i^{(k-r)})}{p_i^{(k)}/(1-p_i^{(k)})} = 
    \fc{\mathrm{odds}(\mu^{(k-r)}(\si_i=1))}{\mathrm{odds}(\mu^{(k)}(\si_i=1))} &= 1-\fc{r}{m^{(k)}}(1+o(1))
\end{align}
Applying this for $r=1$, $i$ replaced by $j$ and  $\mu$ replaced with $\mu|_{-i}$ gives that 
\begin{align}\label{e:odds-ratio-1m}
    \fc{p_{j,-i}^{(k-1)}/(1-p_{j,-i}^{(k-1)})}{p_{j,-i}^{(k)}/(1-p_{j,-i}^{(k)})} = 1-\fc{1}{m_{-i}^{(k)}}(1+o(1)) = 1-\fc{1}{m^{(k)}}(1+o(1)) 
\end{align}
We solve~\eqref{e:p-cvx-comb} and~\eqref{e:odds-ratio-1m}. For ease of notation, let $q=p_j^{(k)}$, $p=p_{j,-i}^{(k-1)}$, $r=p_{j,-i}^{(k)}$, $a=p_i^{(k)}$, $\ep= \rc{m^{(k)}}$. Then~\eqref{e:p-cvx-comb} becomes $q=(1-a)r + qp$ and~\eqref{e:odds-ratio-1m} becomes $\fc{p/(1-p)}{r/(1-r)} = 1-\ep (1+o(1))$ as $\ep\to 0$. This gives
\begin{align*}
\fc{1-r}{r} &= \fc{(1-p)(1-\ep(1+o(1)))}{p}\\
r&= \fc{p}{p+(1-p)(1-\ep + o(\ep))}\\
\implies    q &= \fc{(1-a)p}{p+(1-p) (1-\ep+o(\ep))} + ap\\
    &= (1-a)p(1+\ep(1+o(1))(1-p)+o(\ep(1-p))) + ap\\
    &= p+(1-a)\ep p(1-p)(1+o(1))\\
    \implies 0 &= up^2 - (1+u) p + q
\end{align*}
where $u=(1-a)\ep(1+o(1))$. Then (as $u\to 0$)
\begin{align*}
    \rc p &= \fc{(1+u)+(1+u)\sqrt{1-\fc{4uq}{(1+u)^2}}}{2q} 
    = \fc{1+u}{q} \pa{1-\pf{qu}{(1+u)^2}+O(q^2u^2)}\\
    \implies
    p &= q(1-u+O(u^2))\pa{1+\fc{qu}{(1+u)^2} + O(u^2)} = q - q(1-q)u(1+o(1)).
\end{align*}
Substituting back in, 
\begin{align*}
    p_{j,-i}^{(k-1)} &= p_j^{(k)} - (1-p_i^{(k)})p_j^{(k)}(1-p_j^{(k)})\rc{m^{(k)}} (1+o(1))\\
    \Cov_{\mu^{(k)}}(\si_i,\si_j) &= 
    p_i^{(k)}(p_{j,-i}^{(k-1)} - p_j^{(k)}) = -p_i^{(k)}(1-p_i^{(k)})p_j^{(k)}(1-p_j^{(k)}) \rc {m^{(k)}}(1+o(1)).
\end{align*}
\ppart{Part (2)}
For the second part, we use the superscripts ${}^{(\le k')}$ as the obvious analogue of the notations in the first part. 
To compute expectations, we will condition on which slice the sample is in. Towards this, we let $a_{k'} = \mu_h(\ve{\si}_1=k')$ and $S_{k} = \sumz{k'}{k} a_{k'} = \mu_h(\ve{\si}_1\le k')$. For $k'\le k$, note $\mu^{(\le k)}(\ve{\si}_1=k') = \fc{a_{k'}}{S_k}$. 
The analogue of~\eqref{e:p-cvx-comb} is
\begin{align*}
    p_j^{(\le k)} = (1-p_i^{(\le k)}) p_{j,-i}^{(\le k)} + p_i^{(\le k)} p_{j,-i}^{(\le k-1)}. 
\end{align*}
Then
\begin{align}\label{e:pj-i-pj}
    p_{j,-i}^{(\le k-1)} - p_j^{(\le k)} = (1-p_i^{(\le k)}) (p_{j,-i}^{(\le k-1)} - p_{j,-i}^{(\le k)}).
\end{align}
The rest of the proof will be devoted to computing $p_{j,-i}^{(\le k-1)} - p_{j,-i}^{(\le k)}$.
We have 
\begin{align}
\nonumber
    p_{j}^{(\le k)} &= \rc{S_k}\sumz{k'}k a_{k'} p_j^{(k')}\\
    \nonumber
    p_{j}^{(\le k-1)} &= \rc{S_{k-1}}\sumz{k'}{k-1} a_{k'} p_j^{(k')}\\
\label{e:diff-k-1-k}
    p_{j}^{(\le k-1)} - p_{j}^{(\le k)}
    & = \fc{a_{k}}{S_{k}} \sumz{k'}{k-1} \fc{a_{k'}}{S_{k-1}} \pa{ p_j^{(k')}  - p_j^{(k)}}
\end{align}
Suppose $\fc{\mathrm{odds}(p_j^{(k')})}{\mathrm{odds}(p_j^{(k)})} = R_{k'}$. Then
$p_j^{(k')} = \fc{R_{k'}}{R_{k'}+\rc{p_j^{(k)}}-1}$ and 
\begin{align}
\label{e:pk'-pk}
p_j^{(k')} - p_j^{(k)} = -(1-p_j^{(k)})p_j^{(k)} \fc{1-R_{k'}}{1-p_j^{(k)} (1-R_{k'})}
= 
-(1-p_j^{(k)})p_j^{(k)}
(1-R_{k'})(1 + O(1-R_{k'}))
\end{align}
as $R_{k'}\to 1$. 

\ppart{Computing $p_j^{(\le k-1)}-p_j^{(\le k)}$}
Let $p_i = \mu_h(\si_i=1)$, $m= \sumo in p_i(1-p_i)$, and $\bar k = \E_{\mu_h} \ve{\si}_1$. We show that
\begin{align}
\label{e:pjk1-pjk}
p_j^{(\le k-1)}-p_j^{(\le k)} &= - \fc{c + o(1)}{m^{(\le k)}}p_j^{(\le k)}(1-p_j^{(\le k)}) 
\end{align}
for some $0\le c\le 1$ depending on $h, k, j$, as $m^{(\le k)}\to \iy$, where
\begin{enumerate}
    \item If $k= \bar k + \Om(\sqrt m \log m)$, then $c=0$.
    \item If $k< \bar k + \sqrt m \log m$, then $c=\fc{a_k}{S_k} \sumo r{k - ((k\wedge \bar k) - L(n^{(k \wedge \bar k)}))} \fc{a_{k-r}}{S_{k-1}}r$, for any function $L(x)=\om(\sqrt x)$. 
\end{enumerate}
Consider these 2 cases.
\begin{enumerate}
\item Suppose $k=\bar k  +\Om\pa{\sqrt m \log m}$ (for any choice of constant). We note for any $L$
\begin{align}
\label{e:asp1}
    \fc{a_k}{S_k} p_j^{(k)}= 
\mu^{(\le k)}(\si_j=1 \text{ and } \ve{\si}_1=k)&=
p_j^{(\le k)} \mu_{-j}^{(\le k-1)} (\ve{\si}_1= k-1) \\
\nonumber
\sum_{k'<L} \fc{a_{k'}}{S_{k}} (1-p_j^{(k')})
= \mu^{(\le k)}(\si_j = 0\text{ and }\ve{\si}_1<L)
&= (1-p_j^{(\le k)}) \mu_{-j}^{(\le k)}(\ve{\si}_1<L)\\
\label{e:asp2}
&=(1-p_j^{(\le k)}) \fc{\mu_{h,-j}(\ve{\si}_1<L)}{\mu_{h,-j}(\ve{\si}_1\le k)}.
\end{align}
We note by monotonicity (\Cref{l:monotonicity}) that $p_j^{(k')}\le p_j^{(k)}$ for $k'\le k$. 
Hence, 
\begin{align*}
0\ge_{\eqref{e:diff-k-1-k}} p_j^{(\le k-1)} - p_j^{(\le k)} 
&\ge_{\eqref{e:diff-k-1-k}} - \fc{a_k}{S_k} p_j^{(k)}\\
&=_{\eqref{e:asp1}} -p_j^{(\le k)} \mu_{-j}^{(\le k-1)} (\ve{\si}_1= k-1)
= o\prc{m} p_j^{(\le k)}    
\end{align*}
by noting that $\mu_{-j}^{(\le k-1)} (\ve{\si}_1= k-1)=\fc{\mu_{h,-j}(\ve{\si}_1=k-1)}{\mu_{h,-j}(\ve{\si}_1\le k-1)} = \fc{o(1/m)}{1-o(1)} = o\prc m$ by Bernstein's inequality, noting that the mean and standard deviation for $\ve{\si}_1$ under $\mu_{h}$ differ from $\ve{\si}_1$ under  $\mu_{h,-j}$ by at most 1.
We also have
\begin{align*}
p_j^{(\le k-1)} - p_j^{(\le k)} &= 
\fc{a_k}{S_k}\sum_{k'=0}^{k-1} \fc{a_{k'}}{S_{k-1}} \pa{ p_j^{(k')}  - p_j^{(k)}}
= \fc{a_k}{S_k}\sum_{k'=0}^{k-1} \fc{a_{k'}}{S_{k-1}} \pa{ (1 - p_j^{(k)}) - (1-p_j^{(k')})}\\
&\ge -\fc{a_k}{S_{k-1}}\sum_{k'=0}^{k-1} \fc{a_{k'}}{S_{k}} (1-p_j^{(k')})\\
&=_{\eqref{e:asp2}} -\fc{\mu_h(\ve{\si}_1 = k)}{\mu_h(\ve{\si}_1 \le k-1)} (1-p_j^{(\le k)}) \fc{\mu_{h,-j}(\ve{\si}_1\le k-1)}{\mu_{h,-j}(\ve{\si}_1\le k)}\\
&= -o\prc{m} (1-p_j^{(\le k)})
\end{align*}
by Bernstein's inequality for $\ve{\si}_1$ under $\mu_h$. 
    Putting the last two inequalities together, 
    \begin{align*}
        p_j^{(\le k-1)}-p_j^{(\le k)} = o\prc m p_j^{(\le k)}(1-p_j^{(\le k)}) = o\prc{m^{(\le k)}} p_j^{(\le k)}(1-p_j^{(\le k)})
    \end{align*}
    where the last equality holds because by \Cref{l:mm'}(2), noting $k>\bar k$ and $m=n^{(\bar k)}$, 
    \begin{align*}
        |m^{(\le k)} - m| 
        &\le |m^{(\le k)} - m^{(\fl{\bar k})}| + |m^{(\fl{\bar k})} - n^{(\bar k)}|
        = O(\sqrt{m^{(\le k)}}).
    \end{align*}
    \item Suppose $k< \bar k +\sqrt m \log m$. 
    We would like to use \eqref{e:diff-k-1-k} with \eqref{e:pk'-pk}, but because \eqref{e:pk'-pk} only holds for $R_{k'}\to 1$, we split the sum and first show the tail is negligible.
    Let $L'=\fl{2\sqrt{n^{(k)}}\log n^{(k)}}$. We have using \eqref{e:diff-k-1-k} that
\begin{align}
\nonumber
p_{j}^{(\le k-1)} - p_{j}^{(\le k)}
    & =_{\eqref{e:pk'-pk}} \fc{a_{k}}{S_{k}} \pa{\sum_{k'< k-L'} \fc{a_{k'}}{S_{k-1}} \pa{ p_j^{(k')}  - p_j^{(k)}} - 
 \sum_{k'=k-L'}^{k-1} \fc{a_{k'}}{S_{k-1}}(1-p_j^{(k)})p_j^{(k)} (1-R_{k'})(1 + O(1-R_{k'}))}\\
&=_\eqref{e:odds-ratio-r} 
\fc{a_{k}}{S_{k}} \sum_{k'< k-L'} \fc{a_{k'}}{S_{k-1}} \pa{ p_j^{(k')}  - p_j^{(k)}} -
\fc{a_{k}}{S_{k}} \sum_{r=1}^{L'} \fc{a_{k-r}}{S_{k-1}} (1-p_j^{(k)})p_j^{(k)}\fc{r}{m^{(\le k)}}(1+o(1)),
\label{e:p-diff-break-sum}
\end{align}
where we use $R_{k-r} = 1-\fc r{m^{(k)}} (1+o(1))
= 1-\fc r{m^{(\le k)}} (1+o(1))$ when $r=o(m^{(k)})$ and as $m^{(\le k)}\to \iy$ by \eqref{e:odds-ratio-r} and \Cref{l:mm'}(2). 
For the first sum (tail term) in \eqref{e:p-diff-break-sum}, we note by monotonicity that $p_j^{(k')}\le p_j^{(k)}$ for $k'\le k$, so 
\begin{align}
\nonumber
0\ge \fc{a_k}{S_k}\sum_{k'<k-L'} \fc{a_{k'}}{S_{k-1}} \pa{ p_j^{(k')}  - p_j^{(k)}}
&\ge -\fc{\mu_h (\ve{\si}_1<k-L')}{\mu_h (\ve{\si}_1\le k-1)} \fc{a_k}{S_k} p_j^{(k)}\\
\nonumber
&=_{\eqref{e:asp1}} -\fc{\mu_h (\ve{\si}_1<k-L')}{\mu_h (\ve{\si}_1\le k-1)} p_j^{(\le k)} \mu_{-j}^{(\le k-1)} (\ve{\si}_1= k-1)\\
\nonumber
&\ge -\fc{\mu_{h+t(k-1)\one} (\ve{\si}_1<k-L')}{\mu_{h+t(k-1)\one} (\ve{\si}_1\le k-1)} p_j^{(\le k)} \mu_{-j}^{(\le k-1)} (\ve{\si}_1= k-1)\\
\label{e:1st-term-p}
&=_{\textup{\eqref{e:tail-k-L'}}} o\prc{n^{(k)}} p_j^{(\le k)}
=_{\textup{\Cref{l:mm'}(2)}} o\prc{m^{(\le k)}}p_j^{(\le k)}.
\end{align}
We claim that
\begin{align}
\label{e:tail-k-L'}
    \fc{\mu_h (\ve{\si}_1<k-L')}{\mu_h (\ve{\si}_1\le k-1)} = o\prc{n^{(k)}}.
\end{align}
We prove this by considering two cases.
\begin{enumerate}
    \item $k\ge \bar k$: As $m = n^{(\bar k)}$, we have by \Cref{l:mm'} that
    \begin{align*}
    |m^{(\le k)} - m| 
    &\le |m^{(\le k)} - n^{(k)}| + |n^{(k)} - m|\\
    &\le O(\sqrt{m^{(\le k)}}) + |k-\bar k| = O(\sqrt{m^{(\le k)}}) + O(\sqrt m\log m)
    = o(m^{(\le k)}) = o(n^{(k)})
    \end{align*}
    as $m^{(\le k)}\to \iy$ (equivalently, $m\to \iy$ or $n^{(k)}\to \iy$).
    Hence $2\sqrt{n^{(k)}}\log n^{(k)}\sim 2 \sqrt m\log m$. Then $\mu_h(\ve{\si}_1\le k-1) = \Om(1)$ by Berry-Esseen (or \Cref{t:llt}) and $\mu_h(\ve{\si}_1\le k-L') = o\prc{m} = o\prc{n^{(k)}}$ by Bernstein's inequality, so
    \eqref{e:tail-k-L'} holds.
    \item $k<\bar k$: 
    This means $t(k)<0$, so we have, with Bernstein's inequality for $\mu_{h+t(k)\one}$,
    \begin{align*}
    \fc{\mu_{h}(\ve{\si}_1<k-L')}{\mu_{h}(\ve{\si}_1\le k-1)}
    &\le 
    \fc{\mu_{h+t(k)\one}(\ve{\si}_1<k-L')}{\mu_{h+t(k)\one}(\ve{\si}_1\le k-1)} = o\prc{n^{(k)}}. 
\end{align*}
\end{enumerate}
Now we show \eqref{e:1st-term-p} but with $1-p_j^{(\le k)}$ instead. 
Note that 
\begin{align}
\label{e:bd-w-wo-j}
\mu_h(\ve{\si}_1\le k') \le \mu_{h,-j}(\ve{\si}_1\le k')\le \mu_h(\ve{\si}_1\le k'+1). 
\end{align}
We also have
\begin{align}
\nonumber
\fc{a_k}{S_k}\sum_{k'<k-L'} \fc{a_{k'}}{S_{k-1}} \pa{ p_j^{(k')}  - p_j^{(k)}}
&= \fc{a_k}{S_k}\sum_{k'<k-L'} \fc{a_{k'}}{S_{k-1}} \pa{ (1 - p_j^{(k)}) - (1-p_j^{(k')})}\\
\nonumber
&\ge -\fc{a_k}{S_{k-1}}\sum_{k'<k-L'} \fc{a_{k'}}{S_{k}} (1-p_j^{(k')})\\
\nonumber
&=_{\eqref{e:asp2}} -\fc{\mu_h(\ve{\si}_1 = k)}{\mu_h(\ve{\si}_1 \le k-1)} (1-p_j^{(\le k)}) \fc{\mu_{h,-j}(\ve{\si}_1<k-L')}{\mu_{h,-j}(\ve{\si}_1\le k)}\\
\nonumber
&\ge_{\eqref{e:bd-w-wo-j}} -\fc{\mu_h(\ve{\si}_1 = k)}{\mu_h(\ve{\si}_1 \le k)} (1-p_j^{(\le k)}) \fc{\mu_{h,-j}(\ve{\si}_1<k-L')}{\mu_{h,-j}(\ve{\si}_1\le k-2)}\\
\label{e:1st-term-1-p}
&= o\prc{m^{(\le k)}} (1-p_j^{(\le k)})
\end{align}
where the final line uses $\fc{\mu_{h,-j}(\ve{\si}_1<k-L')}{\mu_{h,-j}(\ve{\si}_1\le k-2)} = o\prc{m^{(\le k)}}$,
established in a similar way to before by considering 2 cases, noting that the relevant quantities for $\mu_h$ and $\mu_{h,-j}$ differ by at most 1 ($|m_{-j} - m|\le 1$, etc.).

Combining \eqref{e:1st-term-p} and \eqref{e:1st-term-1-p}, we obtain the bound on the first term of \eqref{e:p-diff-break-sum}:
\[\fc{a_{k}}{S_k} \sum_{k'< k-L'} \fc{a_{k'}}{S_{k-1}} \pa{ p_j^{(k')}  - p_j^{(k)}} = o\prc{m^{(\le k)}} p_j^{(\le k)}(1-p_j^{(\le k)}).\] 

Now we address the second term in \eqref{e:p-diff-break-sum}.
We have by log-concavity of the sequence $a_0,\ldots, a_k$ (\Cref{l:prod-sr} with \Cref{l:ulc}) that 
\begin{align*}
 \fc{a_k}{S_k} \sumo r{L'} \fc{a_{k-r}}{S_{k-1}}r
 \le 
 \fc{a_k}{S_k} \sumo rk \fc{a_{k-r}}{S_{k-1}}r
 \le \rc{S_kS_{k-1}}\sumo rk (a_{k}a_{k-r}+\cdots + a_{k-r+1}a_{k-1}) \le \fc{S_kS_{k-1}}{S_kS_{k-1}} = 1,
\end{align*}
so the second term is bounded by $\fc{1}{m^{(\le k)}}(1-p_j^{(k)})p_j^{(k)}(1+o(1))$. 

Now we show we can further truncate the sum. Note that $\sumo rs \fc{a_{k-r}}{S_{k-1}} r = \E_{\mu^{(\le k-1)}} (k-\ve{\si}_1)\one_{\ve{\si}_1\ge k-s}$. In the following, we consider when $n^{(k\wedge \bar k)}\to \iy$; note by \Cref{l:mm'}(1) that $m^{(\le k)} \sim m^{(k\wedge \fl{\bar k})} \sim n^{(k\wedge \fl{\bar k})} \sim n^{(k\wedge \bar k)}$. We will show that for any $L = \om (\sqrt{n^{(k \wedge \bar k)}})$ that 
\[
\sum_{r=1}^{(0\vee (k-\bar k)) + L} \fc{a_{k-r}}{S_{k-1}}r = (1-o(1)) \sum_{r=1}^{k} \fc{a_{k-r}}{S_{k-1}}r
\]
Using this for both $L=L' - (0\vee (k-\bar k))$ and for $L=L(n^{(k\wedge \bar k)})$ establishes that 
\[
\sum_{r=1}^{L'} \fc{a_{k-r}}{S_{k-1}}r = 
(1\pm o(1)) \sumo r{k - ((k\wedge \bar k) - L(n^{(k \wedge \bar k)}))} \fc{a_{k-r}}{S_{k-1}}r
\]
and proves the claim about the constant $c$.
\begin{enumerate}
    \item Suppose $k\ge \bar k - \sqrt m$. 
    By the local limit \Cref{t:llt}, by comparing $\mu^{(\le k-1)}$ to a truncated Gaussian centered around $\bar k$,
\begin{align*}
    \E_{\mu^{(\le k-1)}} (k-\ve{\si}_1) &= 
    \E_{\mu^{(\le k-1)}} (\bar k - \ve{\si}_1) + (k-\bar k) = 
    \Om(\sqrt m + |k-\bar k|).
\end{align*}
Note that by \Cref{l:mm'}(1), 
    $n^{(k\wedge \bar k)} \sim n^{(\bar k)}= m$, 
    so we have $L = \om(\sqrt m)$. 
By Cauchy-Schwarz and Bernstein's inequality, 
\begin{align*}
    \E_{\mu^{(\le k-1)}} (k-\ve{\si}_1)\one_{\ve{\si}_1 < (k \wedge \bar k) - L}
    &= \rc{\mu_h(\ve{\si}_1\le k-1)}
    \E_{\mu_h} (k-\ve{\si}_1) \one_{\ve{\si}_1 < (k \wedge \bar k) - L}\\
    &= O(1) \sqrt{\Var_{\mu_h}(\ve{\si}_1) + (k-\bar k)^2} \sqrt{\mu_h (\ve{\si}_1 < (k \wedge \bar k) - L)} \\
    &= o(\sqrt m + k-\bar k) = o\pa{\E_{\mu^{(\le k-1)}} (k-\ve{\si}_1)}
\end{align*}
as $n^{(k\wedge \bar k)} \sim m\to \iy$.
Therefore, $\E_{\mu^{(\le k-1)}} (k-\ve{\si}_1) \one_{\ve{\si}_1 \ge (k\wedge\bar k) - L} = (1-o(1))\E_{\mu^{(\le k-1)}} (k-\ve{\si}_1)$.
\item 
For $k<\bar k - \sqrt m$, let $t$ be such that $\E_{\mu_{h-t\one}}\ve{\si}_1=k$. 
\begin{align*}
    \fc{\sumo rk a_{k-r}r}{S_{k-1}}
    &= \fc{\sumo sk a_{k-s}'}{\sumo sk a_{k-s}'e^{-ts}}
    \sumo rk \pf{a_{k-r}'}{\sumo sk a_{k-s}'} e^{-tr}r\\
    &= \fc{\sumo sk a_{k-s}'}{\sumo sk a_{k-s}'e^{-ts}} \E_{\mu_{h-t\one}} [e^{-tR} R | \ve{\si}_1\le k-1 ]
\end{align*}
where $R=k-\ve{\si}_1$. 
Consider two cases. Define $E(t)$, $m(t)$ as in \Cref{l:E-m-evol}, so we have $E(-t)=k$, $E(0)=\bar k$, $m(-t) = n^{(k)}$, $m(0) = m$.
\begin{enumerate}
    \item $\fc{n^{(k)}}{m}\in [\rc e, e]$: Within this case, as $n^{(k)}\to \iy$, $m\to \iy$. 
    Then \Cref{l:E-m-evol} gives
    \[
\sqrt m = E(0)-E(-t) \le \int_{-t}^0 m(s)ds \le \int_0^t m(0) e^s ds = m(e^t-1)
    \]
    so $t=\Om\prc{\sqrt m}=\Om \prc{\sqrt{n^{(k)}}}$. 
    \item $\fc{n^{(k)}}{m}\nin [\rc e, e]$: \Cref{l:E-m-evol} gives
        \[
m = m(0) \in m(-t)\cdot [e^{-t}, e^t] = n^{(k)} \cdot  [e^{-t}, e^t] 
    \]
    so $t>1$. 
\end{enumerate}
    In either case, we have 
    $\min_{|r|\le \sqrt{n^{(k)}}} e^{-tr}r  = O(1)\cdot \max_{|r|\ge \sqrt{n^{(k)}}}e^{-tr}r$ .
When $L\ge \sqrt{n^{(k)}}$, by Berry-Esseen, 
\begin{align*}
    E_1 := \E_{\mu_{h-t\one}}[e^{-tR}R \one_{k-L\le \ve{\si}_1\le k-1}]
    &\ge \min_{|r|\le \sqrt{n^{(k)}}} e^{-tr}r\cdot  \mu_{h-t\one}(k-\sqrt{n^{(k)}} \le \ve{\si}_1\le k-1)\\
    &= \min_{|r|\le \sqrt{n^{(k)}}} e^{-tr}r \cdot \Om(1)\\
    E_2:=\E_{\mu_{h-t\one}}[e^{-tR}R \one_{ \ve{\si}_1< k-L}] & \le \max_{|r|\ge L} e^{-tr}r \cdot \mu_{h-t\one}(\ve{\si}_1 < k-L)\\
    &= \min_{|r|\le \sqrt{n^{(k)}}} e^{-tr}r \cdot o(1).
\end{align*}
Thus $E_2=o(E_1)$ and dividing by $\mu_{h-t\one}(\ve{\si}_1\le k-1)$, we obtain that 
\[
\fc{\sumo r{k-L} a_{k-r}r}{S_{k-1}} = (1-o(1))
\fc{\sumo r{k} a_{k-r}r}{S_{k-1}}.
\]
as $n^{(k)}\to \iy$, as needed.
\end{enumerate}
\end{enumerate}
\ppart{Comparing with $\mu_{-i}^{(\le k)}$}
We now apply the above analysis to the distribution $\mu_{-i}^{(\le k)}$, showing that we get a constant $c$ that is independent of $i$.
Again, we split into 2 cases:
\begin{enumerate}
    \item If $k\ge \bar k + \sqrt{m}\log m$, then we also have 
    $k-1\ge \E_{\mu_{h,-i}}\ve{\si}_1 + (m_{-i}(1+o(1))\sqrt{(m_{-i}(1+o(1))}-O(1)$ for each $i$, so using \eqref{e:pj-i-pj} and case 1 above for \eqref{e:pjk1-pjk}, we obtain
    \[
p_{j,-i}^{(\le k-1)} - p_j^{(\le k)} = (1-p_i^{(\le k)}) \cdot o\prc{m_{-i}^{(\le k)}} p_{j,-i}^{(\le k)} (1-p_{j,-i}^{(\le k)})
= (1-p_i^{(\le k)}) \cdot o\prc{m^{(\le k)}} p_{j}^{(\le k)} (1-p_{j}^{(\le k)})
    \]
    For the last equality, we reason as follows. We 
    note by \eqref{e:pjk1-pjk} that $p_{j,-i}^{(\le k-1)} - p_{j,-i}^{(\le k)} = o \prc{m^{(\le k)}}p_{j,-i}^{(\le k)}(1-p_{j,-i}^{(\le k)})$, and use the following observation twice: 
    If $x-y = o(y(1-y))$, then because this is $o(\min\{y,1-y\})$, we have $x(1-x) = y(1-y)(1+o(1))$. This observation gives us 
    \begin{align}\label{e:p(1-p)-i}
    p_{j,-i}^{(\le k)}(1-p_{j,-i}^{(\le k)}) \sim p_{j,-i}^{(\le k-1)}(1-p_{j,-i}^{(\le k-1)}) \sim p_{j}^{(\le k)}(1-p_{j}^{(\le k)}).
    \end{align}
    \item Suppose $k< \bar k  +\sqrt m \log m$. 
We now show for a certain choice $L = \om(\sqrt{n^{(k\wedge \bar k)}})$ that 
\begin{align}
\label{e:sum-i-sum}
\fc{a_{k,-i}}{S_{k,-i}} \sumo r{k - ((k \wedge \bar k) - L)} \fc{a_{k-r,-i}}{S_{k-1,-i}}r
= \fc{a_k}{S_k} \sumo r{k - ((k \wedge \bar k) - L)} \fc{a_{k-r}}{S_{k-1}}r (1+o(1)).
\end{align}
We have
\begin{align}
\label{e:aS-i-approx}
\fc{a_{k-r}}{S_k} = (1-p_i^{(k-r)}) \fc{a_{k-r,-i}}{S_{k,-i}} + p_i^{(k-r)} \fc{a_{k-r-1,-i}}{S_{k-1,-i}}
\in \ba{\min\bc{\fc{a_{k-r,-i}}{S_{k,-i}}, \fc{a_{k-r-1,-i}}{S_{k-1,-i}}}, 
\max\bc{\fc{a_{k-r,-i}}{S_{k,-i}}, \fc{a_{k-r-1,-i}}{S_{k-1,-i}}}}.
\end{align}
We will use this equation for $r=0$ to get $\fc{a_{k,-i}}{S_{k,-i}}\approx \fc{a_k}{S_k}$ and for $k\mapsfrom k-1$, $r\mapsfrom r-1$ to get $\fc{a_{k-r,-i}}{S_{k-1,-i}} \approx \fc{a_{k-r}}{S_{k-1}}$. 
For this, we bound the length of the interval in \eqref{e:aS-i-approx}, $\ab{\fc{a_{k-r,-i}}{S_{k,-i}} - \fc{a_{k-r-1,-i}}{S_{k-1,-i}}}$. 
For ease of notation, we consider the full distribution $\mu_h$ and replace by $\mu_{h,-i}$ at the end.

We bound $\fc{a_{k-r}}{S_k}-\fc{a_{k-r-1}}{S_{k-1}}$ by considering two cases. 
\begin{enumerate}
    \item 
First, suppose $k\le \bar k$. 
We would like to use \Cref{t:llt}; however it does not give good bounds in the tails. We can use it to get bounds on ratios if we first tilt the measure so the mean is $k$. Towards this, let $t$ be such that $\E_{\mu_{h-t\one}}\ve{\si}_1=k$; note $t\ge 0$. Let $a_{k'}' = \mu_{h-t\one}(\ve{\si}_1 = k')$. Note that $a_{k-r}\propto e^{-tr} a_{k-r}'$. Choose $L$ such that $\ga\pa{-\fc L{\sqrt{n^{(k)}}}} \ge \rc{(n^{(k)})^{1/4}} \ge \ga\pa{-\fc{L+1}{\sqrt{{n^{(k)}}}}}$. Note by Bernstein's inequality that 
\begin{align*}
    a_{k-L-1}+\cdots + a_0
    &\le e^{-t(L+1)} (a_{k-L-1}'+\cdots + a_0') = e^{-t(L+1)} (a_k'+\cdots + a_{k-L}') 
    \le (a_k+\cdots + a_{k-L}) o(1)\\
    \implies
    a_k+\cdots+ a_{k-L} &= (1-o(1))(a_k+\cdots + a_0).
\end{align*}
Then by \Cref{t:llt}, for $r\le L$
\begin{align*}
    \fc{a_{k-r}}{S_k} &= (1-o(1)) \fc{a_{k-r}}{a_k+\cdots + a_{k-L}}\\
    &=(1+o(1))  \fc{\rc{\sqrt{n^{(k)}}} \ga\pa{-\fc{r}{\sqrt{n^{(k)}}}} + \td O\prc{n^{(k)}}}{\sumz sL \pa{\rc{\sqrt{n^{(k)}}} \ga\pa{-\fc s{\sqrt{n^{(k)}}}} + \td O\prc{n^{(k)}}}e^{-tr}} = (1+o(1))\fc{\ga\pa{-\fc{r}{\sqrt{n^{(k)}}}}e^{-tr}}{\sumz sL \ga\pa{-\fc s{\sqrt{n^{(k)}}}}e^{-ts}}
\end{align*}
and similarly
\begin{align*}
    \fc{a_{k-r-1}}{S_{k-1}} &= (1+o(1))\fc{\ga\pa{-\fc{r+1}{\sqrt{n^{(k)}}}}e^{-t(r+1)}}{\sumo r{L+1} \ga\pa{-\fc r{\sqrt{n^{(k)}}}}e^{-tr}}.
\end{align*}
Then
\begin{align*}
    \fc{a_{k-r}}{S_k}\Big / \fc{a_{k-r-1}}{S_{k-1}} &= (1+o(1)) \fc{\ga\pa{-\fc{r}{\sqrt{n^{(k)}}}}}{\ga\pa{-\fc{r+1}{\sqrt{n^{(k)}}}}}  \fc{{\sumo r{L+1} \ga\pa{-\fc r{\sqrt{n^{(k)}}}}e^{-tr}}}{e^{-t}\sumz rL \ga\pa{-\fc r{\sqrt{n^{(k)}}}}e^{-tr}} \\
    &= (1+o(1))
    \fc{{\sumz r{L} \ga\pa{-\fc r{\sqrt{n^{(k)}}}}(1+o(1))e^{-tr}}}{\sumz rL \ga\pa{-\fc r{\sqrt{n^{(k)}}}}e^{-tr}} = 1+o(1).
\end{align*}
Applying this to $\mu_{h,-i}$, in light of \eqref{e:aS-i-approx}, we obtain $\fc{a_{k-r,-i}}{S_{k,-i}} = \fc{a_{k-r}}{S_{k}} (1+o(1))$.
\item 
Next suppose $k>\bar k$. 
Note $\fc{a_{\fl{\bar k}}}{a_{\fl{\bar k}+1}} = 1\pm o(1)$ as $n^{(\bar k)}\to\iy$ by \Cref{t:llt}. 
Because $k<\bar k + O(\sqrt m\log m)$, we have by \Cref{l:cross-ratio} that for $\bar k< k' \le k$, 
\begin{align*}
    \fc{a_{k'}}{a_{k'-1}} \fc{a_{\fl{\bar k}}}{a_{\fl{\bar k}+1}} = 1-o(1) \implies \fc{a_{k'}}{a_{k'-1}} = 1\pm o(1).
\end{align*}
Similarly, or using \Cref{t:llt}, we have for $\bar k-L\le k' < \bar k$ that $\fc{a_{k'}}{a_{k'+1}} = 1\pm o(1)$. 

We note $S_k = \Om(1)$ and $S_{k}-S_{k-1} = a_k = o(1)$. Hence we also have $\fc{a_{k-r}}{S_k}\Big/\fc{a_{k-r-1}}{S_k} = 1+o(1)$. In the same way, we obtain $\fc{a_{k-r,-i}}{S_{k,-i}} = \fc{a_{k-r}}{S_{k}} (1+o(1))$. 
\end{enumerate}

Applying the above analysis for $\mu_h$ replaced by $\mu_{h,-i}$ and noting the relevant quantities $\bar k$, $m$ differ by $O(1)$, we get that the corresponding terms of \eqref{e:sum-i-sum} differ by a factor of $1\pm o(1)$, establishing that equation.
\end{enumerate}
Hence by \eqref{e:pj-i-pj}, for some $c\le 1$ independent of $i$, 
    \begin{align*}
p_{j,-i}^{(\le k-1)} - p_j^{(\le k)} &= 
-(1-p_i^{(\le k)}) (1-p_{j,-i}^{(\le k)}) p_{j,-i}^{(\le k)} 
\pa{\rc{m_{-i}^{(\le k)}} (c+o(1))}\\
&= -(1-p_i^{(\le k)}) (1-p_{j}^{(\le k)}) p_{j}^{(\le k)} 
\pa{\rc{m^{(\le k)}} (c+o(1))}
    \end{align*}
using \eqref{e:p(1-p)-i} and 
\begin{align*}
    \Cov_{\mu^{(\le k)}}(\si_i,\si_j) &= 
    p_i^{(\le k)}(p_{j,-i}^{(\le k-1)} - p_j^{(\le k)}) = 
    -p_i^{(\le k)}(1-p_i^{(\le k)}) p_j^{(\le k)}(1-p_j^{(\le k)}) \cdot \rc{m^{(\le k)}}(c+o(1)). \qedhere
\end{align*}
\end{proof}

\begin{corollary}[Covariance bound under high total uncertainty]
\label{cor:cov-hi-under-uncert}
Consider a product distribution $\mu_h$ on $\{0,1\}^n$.
    Let $\Om = \slice{n}{k}$ or $\wdg{n}{k}$. Let $p_i = 
    \E_{\mu_h|_\Om} \si_i$. 
    Given $\de$, there exists $m(\de)$ such that if $\sumo in p_i(1-p_i)\ge m(\de)$, then 
    \[
\Cov(\mu_h|_{\Om}) \preceq (1+\de)
\diag(p_i(1-p_i)).
    \]
\end{corollary}
This follows from \Cref{l:cov-tight-off-diag} by noting that it says that the deviation from the covariance being a diagonal matrix of the variances is approximately a negative rank-1 matrix, plus a negligible part.
\begin{prf}
$u\in \R^n$ be such that $u_i = p_i(1-p_i)$ and let 
$D = \diag(u)$. 
By \Cref{l:cov-tight-off-diag}, there is a constant $c\in [0,1]$ depending on $h, n, k$ such that 
\[
\Cov(\mu_h|_{\Om}) = D + \fc{c}{m}\diag \pa{p_i^2(1-p_i)^2} - \fc{uu^{\sT} }{m}c + \fc{uu^{\sT} }{m}\circ E
\preceq \pa{1+\fc{c}{m}}D + \fc{uu^{\sT} }{m}\circ E
\]
where all entries of $E$ are $o(1)$. 
Now, all rows of $\fc{uu^{\sT} }{m}\circ E$ sum to $o(1)$, so its maximum eigenvalue goes to 0 as $m\to \iy$. Hence, choosing large enough, this is $\preceq (1+\delta)D$. \qedhere
\end{prf}

The following bounds the marginal probability of a coordinate of $\mu_h|_{\binom{[n]}{k}}$ in terms of $\mu_{h+t\one}$ for appropriate $t$. 
\begin{lemma}\label{l:product-marginals}
Let $\si(x) = \fc{e^x}{1+e^x}$. 
Let $\ep>0$ be given and let $t_0$ solve
$
\sumo in \si(h_i+t_0) = k.
$
There exist $C_1,C_2$ depending on $\ep$ such that 
     \[
     \fc{\si(h_i+t_0-C_1) - \ep}{1-\ep} \le 
\mu_h^{(k)}(\si_i=1) \le \fc{\si(h_i+t_0+C_2)}{1-\ep}.
     \]
\end{lemma}
\begin{prf}
    The conclusion follows readily from the following two claims.
    \begin{enumerate}
        \item There is $C_1$ such that the $(1-\ep)$-quantile of $\ve{\si}_1$ under $\mu_{h+(t_0-C_1)\one}$ is $\le k$.
        \item There is $C_2$ such that the $\ep$-quantile of $\ve{\si}_1$ under $\mu_{h+(t_0+C_2)\one}$ is $\ge k$.
    \end{enumerate}
    Indeed, under these claims, by monotonicity of $\E[\si_i|\ve{\si}_1=k]$ in $k$ (\Cref{l:monotonicity}), 
    \begin{align*}
\si(h_i+t_0-C_1)&=
\mu_{h+(t_0-C_1)\one}(\si_i=1) = 
\E_{\mu_{h+(t_0-C_1)\one}}[\E[\si_i |\ve{\si}_1]]\\
&\le_{(1)} (1-\ep)\E[\si_i|\ve{\si_1}=k] + \ep = (1-\ep) \mu_h^{(k)}(\si_i=1) + \ep
    \end{align*}
    and 
    \begin{align*}
\si(h_i+t_0+C_2)&= 
\mu_{h+(t_0+C_2)\one}(\si_i=1) = 
\E_{\mu_{h+(t_0+C_2)\one}}[\E[\si_i |\ve{\si}_1]]\\
&\ge_{(2)} (1-\ep)\E[\si_i|\ve{\si_1}=k] = (1-\ep)\mu_h^{(k)}(\si_i=1).
    \end{align*}
    Solving these inequalities gives the result.
    
    Now we prove the claims. Without loss of generality, suppose $t_0=0$.
    Define $E(t)$ and $m(t)$ as in \Cref{l:E-m-evol}.
First, consider if $m(0)<\ep$. Then by Chebyshev's inequality,
\[
\max\{\mu_{h+t\one}(\ve{\si}_1\ge k+1),
\mu_{h+t\one}(\ve{\si}_1\le k-1)\}
\le m(0)<\ep,
\]
which shows the claims for any $C_1,C_2\ge 0$.

Suppose now that $m(0)\ge \ep$. We first show the existence of $C_2$. 
Let $\bar E(t) = E(t)-E(0)$. Note $m'(t)\le m(t)= \bar E'(t)$, so integrating gives $m(t)-m(0)\le \bar E(t)$. We have $\bar E(1) = \int_0^1 m(t)dt  \ge \int_0^1 e^{-t}m(0)dt \ge (1-e^{-1}) m(0)$. Hence for $t\ge 1$,
\begin{align*}
    \fc{m(t)}{\bar E(t)} = \fc{m(t)}{\bar E(t)+m(0)} \fc{\bar E(t)+m(0)}{\bar E(t)}
    \le 1\cdot \fc{2-e^{-1}}{1-e^{-1}}.
\end{align*}
Consider 2 cases. 
\begin{enumerate}
    \item $m(t) \ge \fc{\bar E(t)}{C}$ for $1\le t\le T:=1+C\log \pf{2-e^{-1}}{(1-e^{-1})^2\ep^2}$. Then $\bar E'(t)\ge \fc{\bar E(t)}{C}$ in this interval so 
    \[
\bar E(T) \ge \bar E(1) e^{\fc{T-1}C} \ge  (1-e^{-1})m(0) e^{\fc{T-1}C}\ge \rc{\ep}\cdot \fc{2-e^{-1}}{1-e^{-1}}
    \]
    and by Chebyshev's inequality,
    \[
\mu_{h+
T\one}(\ve{\si}_1\le k)\le 
\fc{m(T)}{\bar E(T)^2}\le \ep.
    \]
    \item There is $1\le t\le T$ such that $m(t)<\fc{\bar E(t)}{C}$, where we take $C=\rc{1-e^{-1}} \rc{\ep^2}$. 
    Then 
    \[
    \mu_{h+T\one}(\ve{\si}_1\le k) \le 
\mu_{h+
t\one}(\ve{\si}_1\le k)\le 
\fc{m(t)}{\bar E(t)^2}
\le \rc{C} 
\rc{\bar E(1)} 
\le \rc{C}
\rc{(1-e^{-1})\ep}
\le \ep.
    \]
\end{enumerate}
For the existence of $C_1$, we consider $E_-(t):= -E(-t)$ and $m_-(t)=m(-t)$. Note that we also have $E_-'(t) = m_-(t)$ and $m_-'(t) \in [-m_-(t), m_-(t)]$. Applying the same argument to $E_-$ and $m_-$ hence furnishes the lower bound. 
\end{prf}
The purpose of the lemma above is in the following corollary.
\begin{corollary}\label{cor:tilt-far-marginal-far}
    Given $0<q<\rc 2$, if $h_i>h_j+C_q$ for a large enough constant $C_q$, then $\mu_h^{(k)}(\si_i=1)$ and $\mu_h^{(k)}(\si_j=1)$ cannot both be in $[q,1-q]$. The same holds true for $\mu_h^{(\le k)}$. 
\end{corollary}
\begin{prf}
    First consider $\mu_h^{(k)}$. Suppose both $\mu_h^{(k)}(\si_i=1)$ and $\mu_h^{(k)}(\si_j=1)$ are in $[q,1-q]$. 
    By \Cref{l:product-marginals} with $\ep=\fc q2$, there exist $C_1, C_2$ such that letting $t$ be such that $\E_{\mu_{h+t\one}} \ve{\si}_1 = k$, 
    \begin{align*}
        \fc{\si(h_i + t-C_1) - q/2}{1-q/2} &\le \mu_h^{(k)}(\si_i=1) \le 1-q\\
        q &\le \mu_h^{(k)}(\si_j=1) \le 
        \fc{\si(h_j + t+C_2)
        }{1-q/2}.
    \end{align*}
    Hence
    \begin{align*}
        \si(h_i+t-C_1) &\le (1-q)\pa{1-\fc q2}+\fc q2 =: q_1\\
        q_2 := q\pa{1-\fc q2}
        &\le \si(h_j+t+C_2)\\
        \implies 
        h_i-h_j-(C_1+C_2) &= (h_i+t-C_1)-(h_j+t+C_2) \le \si^{-1}(q_1)-\si^{-1}(q_2).
    \end{align*}
    Hence we can take $C_q= \si^{-1}(q_1)-\si^{-1}(q_2)+ (C_1+C_2)$. 
    
    Next consider $\mu_h^{(\le k)}$. Suppose both $\mu_h^{(\le k)}(\si_i=1)$ and $\mu_h^{(\le k)}(\si_j=1)$ are in $[q,1-q]$. We have 
    \begin{align*}
        \mu_h^{(\le k)}(\si_i=1) &=\fc{\mu_h(\si_i=1\text{ and }\ve{\si}_1\le k)}{\mu_h(\ve{\si}_1\le k)}
        = \fc{\mu_h(\si_i=1) \mu_{-i}(\ve{\si}_1\le k-1)}{\mu_h(\ve{\si}_1\le k)}
    \end{align*}
    so
    \begin{align*}
        \fc{\mu_h^{(\le k)}(\si_i=1)}{\mu_h^{(\le k)}(\si_j=1)} &= \fc{\mu_h(\si_i=1)}{\mu_h(\si_j=1)} \fc{\mu_{-i}(\ve{\si}_1\le k-1)}{\mu_{-j}(\ve{\si}_1\le k-1)}.
    \end{align*}
    Note that 
    \[
    \mu_h(\ve{\si}_1\le k-1)\le \mu_{-i}(\ve{\si}_1\le k-1) \le \mu_h(\ve{\si}_1\le k).
    \]
    Consider two cases.
    \begin{enumerate}
        \item $\mu^{(\le k)}(\ve{\si}_1=k)\ge 1-\fc q2$. Then
        \begin{align*}
            \pa{1-\fc q2} \mu_h^{(k)}(\si_i=1) \le \mu_h^{(\le k)} (\si_i=1) \le \fc q2 + \pa{1-\fc q2} \mu_h^{(k)}(\si_i=1)
        \end{align*}
        and similarly for $j$. Hence if $\mu_h^{(\le k)}(\si_i=1),\mu_h^{(\le k)}(\si_j=1)\in [q,1-q]$, then 
        $\mu_h^{(k)}(\si_i=1), \mu_h^{(k)}(\si_j=1)\in \ba{\fc{q/2}{1-q/2}, \fc{1-q}{1-q/2}}$. The result then follows from the first part with $q$ replaced by $\fc{q/2}{1-q/2}$.
        \item $\mu^{(\le k)}(\ve{\si}_1=k)< 1-\fc q2$. Then $\fc{\mu_{-i}(\ve{\si}_1\le k-1)}{\mu_{-j}(\ve{\si}_1\le k-1)} \in
        \ba{\fc{\mu_h(\ve{\si}_1\le k-1)}{\mu_h(\ve{\si}_1\le k)}, \fc{\mu_h(\ve{\si}_1\le k)}{\mu_h(\ve{\si}_1\le k-1)}}\subseteq 
        [\fc q2, \fc 2q]$. Hence if $\mu_h^{(\le k)}(\si_i=1),\mu_h^{(\le k)}(\si_j=1)\in [q,1-q]$, then 
        $\fc{\mu_h(\si_i=1)}{\mu_h(\si_j=1)}\in \ba{\fc{q}{1-q}\cdot \fc{q}2, \fc{1-q}{q}\cdot \fc{2}{q}}$. We can take $C_q = \si^{-1}\pa{\fc{1-q}{q}\cdot \fc{2}{q}} - \si^{-1}\pa{\fc{q}{1-q}\cdot \fc{q}2}$. \qedhere
    \end{enumerate}
\end{prf}
In preparation for applying the results to the Ising model, we now consider product distributions on $\{-1,1\}^n$ rather than on $\{0,1\}^n$. The previous results still hold up to scaling.
The following says that in order for the variance of $\an{\si,x}$ on a product distribution on the slice or wedge to be large, there must be a coordinate on which $x$ is positive and bounded away from 0, and where $x$ is negative and bounded away from 0. 
\begin{lemma}\label{l:exist-small-uncertain}
Let $x_0\in \{\pm 1\}^n$ and $\Om = \ball{k}{x_0}$ or $\bball{k}{x_0}$. 
Let $B<C$ and $\ep,\de>0$ be given. 
There exist $q, c, C_\de>0$ such that the following hold. 
    Suppose that $\ve{x}\le C$, $\max_{|S|\le 4k \vee \fc{48m(\de)}{\de}} \ve{x_S}\le B$ and $\Var_{\mu_{h}|_\Om}(\an{\si,x})> (1+\de) B^2$, where $m(\de)$ is as in \pref{cor:cov-hi-under-uncert}.
    Then there exists $i,j$ with $x_i\le -c$ and $x_j\ge c$, with $p_i=\mu_h|_{\Om}(\si_i=1)\in [q,1-q]$ and $p_j = \mu_h|_{\Om}(\si_j=1)\in [q,1-q]$.
\end{lemma}
The idea is that from \pref{cor:cov-hi-under-uncert}, we must have a small block of coordinates with large variance. Then we use the fact that off-diagonal entries of the covariance are negative, which forces us to have both positive and negative entries in that block of $x$.
\begin{prf}
 Let $\Si = \Cov_{\mu_h|_{\Om}}(\si)$ and $D = 4\diag(\mu_h|_{\Om}(\si_i=1)\mu_h|_{\Om}(\si_i=-1))$. 

First we claim that $\sumo in p_i(1-p_i)< m(\de)$. Indeed, if $\sumo in p_i(1-p_i)\ge m(\de)$, then by \pref{cor:cov-hi-under-uncert}, because $4p_i(1-p_i) \le 1$ and $4\sumo in p_i(1-p_i) \le 4k$, 
\begin{align*}
    \Var_{\mu_h|_\Om}(\an{\si, x}) &= x^{\sT} \Si x \le 
    (1+\de) x^{\sT} D x  = 
    4(1+\de) \sumo in p_i(1-p_i)x_i^2
    \le (1+\de) B^2,
\end{align*}
which is a contradiction.
    
    Thus there are at most $\fc{m(\de)}{q(1-q)}$ coordinates such that $p_i\in [q,1-q]$.
    Let $S= \set{i}{p_i\in [q,1-q]}$, $S_+=\set{i}{i\in S, x_i>0}$, $S_-=\set{i}{i\in S, x_i<0}$. 
    We first bound with a block-diagonal matrix by the inequality $2ab \le \lm a^2 + \rc{\lm}b^2$, for $\lm>0$:
    \[
    \Si \preceq  
    (\Si + \lm D)\circ J_S + \pa{\Si + \rc{\lm}D} \circ J_{S^c},
    \]
    where $J_S$ denotes the matrix that is all-1 on $S\times S$ and 0 elsewhere, and $\circ$ denotes Hadamard product.
    Choose $\lm=\fc{\de}{2}$, $q=\fc{\de}{16}$. 
    Let $\ve{x_{S_+}}=r_+$, $\ve{x_{S_-}}=r_-$, $\ve{x_{S^c}}=r_0$. Note that 
    $x_{S_{\pm}}^{\sT} \Si_{S_{\pm}}x_{S_{\pm}}^{\sT} \le 
    x_{S_{\pm}}^{\sT} D_{S_{\pm}}x_{S_{\pm}}^{\sT} \le r_{\pm}^2$ because off-diagonal entries only contribute negatively. 
    Then
    \begin{align*}
    x^{\sT} \Si x 
    &\le x_S^{\sT} (\Si + \lm D)_{S\times S} x_S + x_{S^c}^{\sT}  \pa{\Si + \rc{\lm}D}_{S^c\times S^c} x_{S^c}\\
    &\le \lm_+(r_+^2+r_-^2) +2\opnorm{\Si+\lm D} r_+r_- + \pa{2+\rc{\lm}}x_{S^c}^{\sT} D_{S^c\times S^c}x_{S^c}.
    \end{align*}
    For $i\in S^c$, $\pa{2+\rc{\lm}} \cdot 4p_i(1-p_i)\le \fc{4}{\de} \cdot 4q= 1$. Hence, noting $\sum_{i\in S^c} \pa{2+\rc{\lm}} \cdot 4p_i(1-p_i) \le
    \sum_{i\in S^c} \pa{2+\rc{\lm}} \cdot 4p_i \le  
    \fc{4}{\de}\cdot 4 m(\de)$, 
    \begin{align*}
        \pa{2+\rc{\lm}}x_{S^c}^{\sT} D_{S^c\times S^c}x_{S^c}
        &= \sum_{i\in S^c}\pa{2+\rc{\lm}} \cdot 4p_i(1-p_i) x_i^2 
        \le \max_{\substack{T\subseteq S^c\\|T|\le \fc{16}{\de}m(\de)}} \ve{\hx_T}^2.
    \end{align*}
    Note $\lm_+(r_+^2+r_-^2)= \lm_+\sum_{i\in S}x_i^2$ with $|S|\le \fc{m(\de)}{q(1-q)}$ and $\opnorm{\Si+\lm D}\le 2+\lm$, so 
    noting that $\fc{16}{\de}m(\de) + \fc{m(\de)}{q(1-q)}\le \fc{48m(\de)}{\de}$,
    \[
x^{\sT} \Si x \le \lm_+ \max_{|T|\le \fc{16}{\de}m(\de)+ \fc{m(\de)}{q(1-q)}} \ve{x_T}^2 + 2(2+\lm)r_+r_-
\le \lm_+ B^2 + 2(2+\lm)r_+r_-.
    \]
    Because $\lm=\fc{\de}2$, in order the inequality $\Var_{\mu_{h}|_\Om}(\an{\si,x})\ge (1+\de) B^2$ to hold, we must have $8r_+r_- \ge 2\pa{2+\fc{\de}{2}} r_+r_- \ge \fc{\de}2 B^2$. Because $r_+,r_-\le C$, we must have $r_+, r_-\ge \fc{\de B^2}{16C}$. Then there must be entries in $x_{S_+}$ and $x_{S_-}$ that are $\ge \rc{\sqrt{|S|}}\fc{\de B^2}{16C}\ge \sfc{q(1-q)}{m(\de)}\fc{\de B^2}{16C}=:c$.
\end{prf}
\begin{lemma}\label{l:var-large-exist-sets}
    Let $\de>0$ be given. Then there is a constant $K(\de, B, C)$ such that for $x$ such that $\ve{x}\le C$ and $\max_{|S|\le 4k \vee \fc{48m(\de)}{\de}} \ve{x_S}\le B$ (with $m(\de)$ as in \Cref{cor:cov-hi-under-uncert}),
    \[
\lm\ba{\set{u}{\Var_{\mu_{h+ux}|_\Om}(\an{\si,x})> (1+\de) B^2}}\le K(\de, B, C),
    \]
    where $\lm$ denotes Lebesgue measure over $\R$.
\end{lemma}
The idea of the proof is that the set of $u$'s for which the condition in \Cref{l:exist-small-uncertain} is satisfied for a given 
$(i,j)$ has at most constant measure, because the probabilities 
$p_i$ and $p_j$ will drift apart. 
Since there are a constant number of possible 
pairs $(i,j)$ for a given $x$, the set of $u$'s for which the variance is large also has at most constant measure.
\begin{prf}
    By \Cref{l:exist-small-uncertain}, there exist $q, c>0$ such that if $\Var_{\mu_{h+ux}|_\Om}(\an{\si,x})> (1+\de) B^2$, then there exist $i,j$ with $x_i\le -c$ and $x_j\ge c$ such that $\mu_{h+ux}|_{\Om}(\si_k=1)\in [q,1-q]$ for $k\in \{i,j\}$:
\begin{align*}
    \set{u}{\Var_{\mu_{h+ux}|_\Om}(\an{\si,x})> (1+\de) B^2}
    &\subseteq 
    \bigcup_{\substack{i:x_i\le-c\\ j: x_i\ge c}}
    \set{u}{\mu_{h+ux}|_{\Om}(\si_i=1), \mu_{h+ux}|_{\Om}(\si_j=1)\in [q,1-q]}\\
    &\subseteq 
    \bigcup_{\substack{i:x_i\le-c\\ j: x_i\ge c}} \set{u}{|(h_i+ux_i)-(h_j+ux_j)|\le C_q},
\end{align*}
where the second inclusion follows from \pref{cor:tilt-far-marginal-far} for appropriate constant $C_q$. For given $x$, there are at most $\pf{C}{c}^2$ coordinates $i$ with $|x_i|\ge c$, so there are at most $\ba{\rc2\pf{C}{c}^2}^2$ such pairs of $(i,j)$, and the interval of $u$ where each inequality is true has measure at most $\fc{C_q}{c}$, so the measure of the set is at most $\pf{C^4}{4c^4}\fc{C_q}{c}$.
\end{prf}

\begin{lemma}[Variance bound for mostly log-concave distribution]\label{l:mostly-lc-var}
    Let $L, \al>0$. 
    Consider a probability measure $\mu$ with density $p$ on $\R$ such that $-(\log p)''(u)$ such that 
    \begin{enumerate}
        \item $-(\log p)''(u)\ge -L$ everywhere, and
        \item $-(\log p)''(u)\ge \al$ except on a set of Lebesgue measure $M$.
    \end{enumerate}
    Fix $A>0$. 
    Then $\mu$ has log-Sobolev constant at least $\fc{\al}{2} \exp\pa{-\fc{M^2(L+\al)^2}{\al}}$. 
    Thus, $\Var(\mu)\le \fc2\al \exp\pa{\fc{M^2(L+\al)^2}{\al}}$.
\end{lemma}
\begin{prf}
    We will show that $\mu$ has bounded ratio with a strongly log-concave measure. 
    Define an unnormalized probability density $\td q$ as follows. Let $\td q(0)=p(0)$ and define $g$ (which will be set to be $(\log \td q)'$) by
    \begin{align*}
        &&g(0) &= (\log p)'(0)\\
        \text{for }u>0, &&
        \dd{^+}{u}g(u) &= \begin{cases}
            -\fc{\al}2, &
            (\log p)'(u)>-\fc\al2 \text{ or } g(u) < (\log p)'(u)\\
            (\log p)''(u), & (\log p)'(u)\le -\fc\al2 \text{ and }g(u)=(\log p)'(u).
        \end{cases}\\
        \text{for }u<0, &&
        \dd{^-}{u}g(u) &= \begin{cases}
            -\fc{\al}2, & (\log p)'(u)>-\fc\al2 \text{ or } g(u) < (\log p)'(u)\\
            (\log p)''(u), & (\log p)'(u)\le -\fc\al2 \text{ and }g(u)=(\log p)'(u).
        \end{cases}
    \end{align*}
    Then we let 
$\log \td q(u) = \int_0^u g(t) dt$. Note that $\td q$ is $\fc{\al}2$-strongly log-concave.
    We now bound $\log \td q - \log p$. This is maximized when the area between the graphs of $g$ and $(\log p)'$ is maximized. Intuitively, this happens when the part when $(\log p)'\not\le -\fc{\al}2$ is a single interval and maximally violated, $(\log p)' = L$, and then this gap is shrunk as slowly as possible (with rate $\fc\al2$). We now show this formally.
    Let $I_+= \set{u>0}{(\log p)''(u)>-\al}$ and $I_-= \set{u<0}{(\log p)''(u)>-\al}$. These is the set where $(\log p)'$ is increasing, or not decreasing enough.  Let $M_{\pm}=\lm(I_{\pm})$. Note $M_++M_-\le M$. 
    Define also the set where $(\log p)'$ is decreasing but $g\ne (\log p)'$:
    \begin{align*}
        D_+ &= \set{u>0}{(\log p)''(u)\le -\al \text{ and }g(u) \ne (\log p)'(u)}.
    \end{align*}
    Note that $I_+,D_+$ form a partition of the set $\set{u>0}{g(u) \ne (\log p)'(u)}$. 
    Let $F_+^I, F_+^D$ be the cumulative distribution functions for $I_+$ and $D_+$, respectively. 
    Because $(\log p)''(u)-g'(u)\le L+\fc{\al}2$ on $I_+$, and $(\log p)''(u)-g'(u)\le 0$ elsewhere, 
    \begin{align*}
        (\log p)'(u)-g(u) &\le \pa{L+\fc \al 2} F_+^I(u).
    \end{align*}
    Because $\lim_{u\to \iy} (\log p)'(u)-g(u)=0$ and $(\log p)''(u) - g'(u)\le -\al - \pa{-\fc{\al}2} = -\fc\al 2$ on $D_+$,
    \begin{align*}
        (\log p)'(u)-g(u) &\le \fc{\al}2 (F_+^D(\iy)-F_+^D(u)).
    \end{align*}
    Note the following claim. If $f(x) \ge b $ and $\int_A f(x)dx \le B$ for $A\subset \R_{\ge 0}$, then letting $F(x) = \int_0^x f(t)\one_A(t)dt$, using the substitution $t=F^{-1}(u)$, 
    \begin{align*}
        \int_0^x F(t) dt \le \int_0^{F(x)} u\dd{F^{-1}}{u}(u)du \le \int_0^{F(x)} \fc ub du \le \fc{F(x)^2}{2b}\le \fc{B^2}{2b}.
    \end{align*}
    We apply this to $f(x) = [(\log p)''(x_0-u) - g'(x_0-u)] - [(\log p)''(x_0) - g'(x_0)]$ for $u\to \iy$, noting 
    \[-\int_0^\iy (\log p)''(u) - g'(u) 
        dF_-^I(u) = 
        \int_0^\iy (\log p)''(u) - g'(u) 
        dF_+^I(u) \le \pa{L+\fc\al2} M_+
    \]
    to conclude 
    \begin{align*}
        \int_0^\iy (\log p)'(u) - g(u) 
        dF_-^I(u) \le \rc 2\cdot \fc 2\al \pa{\pa{L+\fc\al2} M_+}^2
    \end{align*}
    Then
    \begin{align*}
        \int_0^\iy 
        (\log p)'(u) - g(u)  du &\le 
        \int_0^\iy (\log p)'(u) - g(u)
        dF_+^I(u) + 
        \int_0^\iy (\log p)'(u) - g(u) 
        dF_-^I(u)\\
        &\le \int_0^\iy \pa{L+\fc\al2} F_+^I(u)
        dF_+^I(u) + 
        \rc{\al}\pa{\pa{L+\fc\al2} M_+}^2\\
        &\le \rc 2 M_+^2\pa{L+\fc\al2} + \rc{\al}\pa{\pa{L+\fc\al2} M_+}^2\\
        &\le M_+^2 \pa{L+\fc\al2}\pa{\fc L\al+1} \le \fc{M_+^2}{\al}(L+\al)^2. 
    \end{align*}
    Similarly, we have $\int_{-\iy}^0 g(u) - (\log p)'(u) du \le \fc{M_-^2}{\al}(L+\al)^2$. Thus $(\sup - \inf)(\log \td q - \log p)\le \fc{M^2}{\al}(L+\al)^2$. Taking into account the normalization constant, letting $q(u) \propto \td q(u)$ be a probability density, $|\log  q - \log p|\le \fc{M^2}{\al}(L+\al)^2$. By Bakry-Emery, because $q$ is $\fc{\al}2$-strongly log-concave, it has log-Sobolev constant at least $\fc{2}{\al}$. By Holley-Stroock, $p$ has log-Sobolev constant at least $\fc{\al}{2} \exp\pa{-\fc{M^2(L+\al)^2}{\al}}$. This implies a Poincaré inequality with the same constant and hence the variance bound.
\end{prf}
\begin{lemma}[Entropy contraction for Ising model on wedges, II]\label{l:LS-from-AS-II}
Let $x_0\in \{\pm 1\}^n$. 
Given $\de>0$, $C$, $B<\rc{1+\de}$, there is $C_\de$ such that the following holds. Suppose $A$ is a PSD matrix such that $\opnorm{A}\le C^2$ and for all $|S|\le 4\ep n \wedge \fc{48m(\de)}{\de}$ ($m(\de)$ as in \Cref{cor:cov-hi-under-uncert}), $\opnorm{A_{S\times S}}\le B^2$. Then for any $h\in \R^n$, $n \cdot \rec(\mu_{A,h}|_{\ball{\ep n}{x_0}})$ is bounded below by a positive constant depending only on $\de$, $C$, and $B$.
\end{lemma}
\begin{proof}
We use the two-stage decomposition, \pref{lem:LS-2-stage}. 
Let $x\in \R^n$ be such that $xx^{\sT} \preceq A$, so $\ve{x}\le C$.  
It suffices to bound $\Var(p)$, where $p$ is given by \eqref{e:1d-p}, for $\lm=1$.
We have 
\[
\lm\ba{\set{u\in \R}{\Var_{\mu_{h+ux}|_{\ball{\ep n}{x_0}}}(\an{\si,x})> (1+\de) B^2}}\le K
\]  
    for some constant $K$ by \Cref{l:var-large-exist-sets}. 
    Note also that 
    \[
\Var_{\mu_{h+ux}|_{\ball{\ep n}{x_0}}}(\an{\si,x})
\le x^{\sT}\Cov_{\mu_{h+ux}|_{\ball{\ep n}{x_0}}}(\si)x\le 2C^2
    \]
    for all $u\in \R$ by \Cref{l:prod-dist-cov}.
    Then 
    \begin{align*}
    -\log p''(u)
    = 1 - \Var_{\mu_{ux+h}|_{\ball{\ep n}{x_0}}}(\an{x,\si})
    \ge \begin{cases}
        -(2C^2-1), &\text{everywhere}\\
        1-(1+\de)B>0, &\text{except on a set of Lebesgue measure $K$}.
    \end{cases}
\end{align*}
Therefore, by \Cref{l:mostly-lc-var}, $\Var_p(u)$ is bounded by a constant. 
The result then follows from \pref{lem:LS-2-stage}.
\end{proof}
\Cref{l:ls-sk-wedge-ii} then follows in the same way as \Cref{l:ls-sk-wedge}.
\begin{prf}[Proof of \Cref{l:ls-sk-wedge-ii}]
For $A\sim \GOE(n)$, consider $\be A + \ga I$, where $2\be < \ga < 1$ (e.g. $\ga=\be+\rc{2}$). With probability $\ge 1-e^{-\Om(n)}$, $\be A + \ga I$ is PSD and has operator norm bounded by $2$. 
Also with probability $1-e^{-\Om(n)}$, by \Cref{l:AS}, for $|S|\le 4\ep n$,
\[
\opnorm{(\be A+\ga I)_{S\times S}} \le \ga + \be c_1\sqrt{h(4\ep)} <1
\]
when $\ep$ is small enough.
The result than follows from \Cref{l:LS-from-AS-II}.
\end{prf}

\subsection{Concentration of $\an{\si,x_0}$}
\label{s:hi-prob-wedge}

The fact that the localized distribution is concentrated on a wedge is a generic property of the stochastic localization process.
\begin{lemma}\label{l:conc-x0}
Let $\mu$ be any distribution on the hypercube, and consider the stochastic localization process with $x_0\sim \mu$, $y_t = tx_0+B_t$.
Let $0<\ep<\rc 2$. 
\begin{enumerate}
    \item With probability $\ge  1-2e^{-\fc n2 D_{\KL}(\ep\| \Phi(-\sqrt t))}$, 
    \[
\mu_{y_t}(\ball{\ep n}{\sign(y_t)}) \ge 1-e^{-\fc n2 D_{\KL}(\ep\| \Phi(-\sqrt t))}.
    \]
    In particular, for $t \ge 2\log \pf{e^2}{\ep}$, with probability $\ge 1-2e^{-\fc{n\ep}2}$, 
\[
 \mu_{y_t}(\ball{\ep n}{\sign(y_t)}) \ge 1-e^{-\fc{n\ep}2} .
\]
\item 
With probability $\ge  1-2e^{-\fc n2 D_{\KL}(\fc\ep2\| \Phi(-\sqrt t))}$,
\[
 \mu_{y_t}(\ball{\ep n}{x_0}) \ge 1-e^{-\fc n2 D_{\KL}(\fc{\ep}2\| \Phi(-\sqrt t))} .
\]
In particular, for  $t\ge 2\log \pf{2e^2}{\ep}$, with probability $\ge 1-2e^{-\fc{n\ep}4}$, 
\[
 \mu_{y_t}(\ball {\ep n}{x_0}) \ge 1-e^{-\fc{n\ep}4} .
\]
\end{enumerate}
\end{lemma}
\begin{prf}

    Let $\Phi$ denote the cdf of the normal distribution. Note that for each $i$, 
    \[\P((x_0)_i\ne \sign((y_t)_i)) = \Phi(-\sqrt t)\] 
    and these events are independent.
    By a Chernoff bound,
    \begin{align}\label{e:chernoff-round}
    \P(d(x_0,\sign(y_t)) \ge \ep)\le e^{-nD_{\KL}(\ep\| \Phi(-\sqrt t))}. 
    \end{align}
        Because 
        \[\P(d(x_0,\sign(y_t)) > \ep) = 
\E[\P[d(\sign(y_t), x_0)> \ep|y_t]]
= \E[\mu_{y_t}(\ball{n\ep}{\sign(y_t)}^c)]
,\] we have by Markov's inequality that 
\[
\P\pa{ \mu_{y_t}(\ball{n\ep}{\sign(y_t)}) \ge 1-2e^{-\fc n2 D_{\KL}(\ep\| \Phi(-\sqrt t))} } \ge 1-e^{-\fc n2 D_{\KL}(\ep\| \Phi(-\sqrt t))}.
\]
Using $\Phi(u)\le e^{-\rc2u^2}$ for $u\le 0$, we have for $e^{-\rc2t}\le \ep$ that
\begin{align*}
D_{\KL}\pa{\ep\|\Phi(-\sqrt t)}
\ge D_{\KL}\pa{\ep\|e^{-\rc 2t}}
\ge \ep\log \fc{\ep}{e^{-t/2}}
+ \pa{1-\ep} \log \fc{1-\ep}{1-e^{-t/2}}
\ge \ep\pa{\log \ep + \fc t2 - 1}.
\end{align*}
Plugging in $t\ge \log \pf{2e^2}{\ep}$ gives the simplified bound.

For (2), consider $\hat x_0$ drawn from the localized distribution $\mu_{y_t}$. Note that this is a posterior distribution of $x_0$ given the observation $y_t$ (with prior $\mu$), so the distribution of $(x_0,y_t)$ and $(\hat x_0, y_t)$ are identical.

    Since $(x_0,y_t)\stackrel d=(\hat x_0, y_t)$, we also have that \eqref{e:chernoff-round} holds for $d(\hat x_0,\sign(y_t))$. Hence
    \[
\P(d(x_0,\hat x_0) \ge 2\ep) \le 
\P(d(x_0,y_t) \ge \ep) + \P(d(\hat x_0,y_t) \ge \ep)
\le 2e^{-nD_{\KL}(\ep\| \Phi(-\sqrt t))}
    \]
    Because $\P(d(x_0,\hat x_0) \ge 2\ep) = 
\E[\P[d(\hat x_0, x_0)\ge 2\ep|x_0]]$, we have by Markov's inequality that 
\[
\P\pa{ \mu_{y_t}(\ball{2\ep n}{x_0}) \ge 1-2e^{-\fc n2 D_{\KL}(\ep\| \Phi(-\sqrt t))} } \ge 1-e^{-\fc n2 D_{\KL}(\ep\| \Phi(-\sqrt t))}.
\]

Now replace $\ep$ with $\fc \ep2$. 
Plugging in $t\ge 2\log \pf{2e^2}{\ep}$ gives the final statement.
\end{prf}
\begin{prf}[Proof of \pref{t:ls-sk}]
This follows from combining \pref{l:conc-x0} and \pref{l:ls-sk-wedge-ii}.
\end{prf}

\subsection{Sampling and partition function calculation}
\label{s:samp-Z}
\begin{corollary}[Sampling]
\label{cor:wedge-sampling}
Let $\be<\rc 2$. 
    There are constants $\ep(\be), T, \rh_\be$ depending on $\be$ such that with probability $\ge 1-e^{-\Om(n)}$, for any $\be'<\be$, 
    \begin{enumerate}
        \item For any $y\in \R^n$, $\mu_{\be' A,y}|_{\ball{\ep(\be) n}{\sign(y)}}$ satisfies entropy contraction with constant $\ge \rc n \rh_\be$, and 
        \[
\tmix(\mu_{\be' A,y}|_{\ball{\ep(\be) n}{\sign(y)}}, \de_{\sign(y)}, \ep)
\le \ce{\fc{n}{\rh_\be}\log\pf{2n}{\ep^2}
}.
        \]
        \item 
        For $y_T = Tx_0 + B_T$ drawn according to SL, taking the polarized walk for 
        $\mu_{\be' A,y_T}|_{\ball{\ep(\be) n}{\sign(y_T)}}$ for 
        $\ge \ce{\fc{n}{\rh_\be}\log\pf{4n}{\ep^2}}$ steps starting from $\sign(y_T)$ gives a sample that is $\ep$-close in TV distance to $\mu_{\be' A,y_T}$.
    \end{enumerate}
\end{corollary}
\begin{prf}
For (1), We note that by \Cref{t:ls-sk}, with probability $\ge 1-e^{-\Om(n)}$, $\mu_{\be' A,y}|_{\ball{\ep(\be) n}{\sign(y)}}$ satisfies entropy contraction with constant $\rc n \rh_\be$. 
Let $\osc_\Om f := \max_\Om f- \min_\Om f$. 
As part of the high-probability event, we may assume $\opnorm{\be A}\le 1$. 
Note that 
\[
\dd{\delta_{\sign(y)}}{\mu_{\be' A, y}}
\le e^{\osc_{\si\in \{\pm 1\}^n} (\be' \an{x,Ax})} 2^n
\le e^{2n} 2^n
\]
The mixing time bound follows from \Cref{l:tmix}.

For (2), apply (1) with $\ep$ replaced by $\fc \ep2$ and $y=y_T$, and note further that by \Cref{t:ls-sk}, we can choose $T$ large enough so that with high probability, $\TV (\mu_{\be' A,y}|_{\ball{\ep(\be) n}{\sign(y)}}, \mu_{\be' A,y})\le \fc \ep 2$.
\end{prf}

\begin{algorithm}[!ht]
 \caption{Simulated annealing for partition function estimation} 
 \begin{algorithmic}
 \Require{Sampling oracles for $\td p_\ell$ (approximations to  $p_\ell\propto q_\ell$) for $1\le \ell\le M$ (distributions on $\Om$), for example, Glauber dynamics or the polarized walk; $Z_1 = \int_{\Om} q_1\,d\om$; number of samples $N$; number of trials $R$.}
 \Ensure{Estimate of $\int_{\Om} q_{\ell}\,d\om$ for each $1\le \ell\le M+1$.}
 \State Let $g_\ell(x) := \fc{q_{\ell+1}(x)}{q_\ell(x)}$.
 \For{$1\le r\le R$}
     \State Let $\hat Z^r_1 = Z_1$.
    \For{$1\le \ell\le M$}
         \State Obtain samples $x_1,\ldots, x_N\sim \td p_\ell$.
         \State Let $\hat Y_{\ell} = \rc N \sumo kN g_\ell(x_k)$.
         \State Let $\hat Z^r_{\ell+1} = \hat Z^r_\ell \hat Y_\ell$.
    \EndFor
 \EndFor 
 \For{$2\le \ell\le M+1$}
     Let $\hat Z_\ell$ be the median of $\set{\hat Z_\ell^r}{1\le r\le R}$.
\EndFor 
\end{algorithmic}
  \label{alg:sa}
 \end{algorithm}

\begin{lemma}[{Partition function estimation with simulated annealing, \cite[Lemma C.4]{koehler2022sampling}, cf. \cite{dyer1991computing}}]
\label{lem:sa}
Let $0<\ep<1$. Suppose that $p_\ell, 1\le\ell\le M+1$ are distributions on $\Om$, and that in \pref{alg:sa} we are given sampling oracles for $\td p_\ell$, $1\le \ell\le M$ such that the following hold for each $1\le \ell \le M$.
\begin{enumerate}
    \item (Variance bound) $\fc{\Var_{P_\ell}(g_\ell(x))}{(\E_{P_\ell} g_\ell(x))^2} \le \si^2$.
    \item (Bias bound) $\ab{\E_{P_\ell}g_\ell(x) - \E_{\td P_\ell}g_\ell(x)}\le \fc{\ep}{4M}$.
\end{enumerate}
Then taking $N\ge \fc{320\si^2 M}{\ep^2}$ and $R\ge 32 \log  \prc\de$, with probability $1-\de$, the output $\hat Z$ satisfies $\hat Z\in [e^{-\ep}, e^{\ep}] \cdot Z$.
    
\end{lemma}
A non-adaptive temperature schedule of length $O(n)$ is sufficient for partition function estimation. Note that a shorter schedule of length $O(\sqrt n \log  n \log  \log  n)$ is possible, and can be found in $n\poly\log(n)$ total queries to approximate sampling oracles 
at the different temperatures~\cite{vstefankovivc2009adaptive}, or constructed non-adaptively~\cite{liu2024work}, but we use a geometric non-adaptive schedule for simplicity.

\begin{corollary}
    \label{c:wedge-Z}
Let $Z_{\be A, y}|_{\Om} := \sum_{\si\in \Om} e^{\rc 2\be \an{\si,A\si} + \an{\si,y}}$.
Let $\ep = e^{-O(n)}$. 
    Under the same conditions as \pref{cor:wedge-sampling} (for fixed $\be<\rc 2$), we have that for \pref{alg:sa} with 
    unnormalized distributions 
    $q_\ell(\si) = e^{\rc 2\beta_\ell \an{\si, A\si} + \an{\si,y}}$,
    schedule $\beta_\ell = \fc{(\ell-1)\beta}{4n}$, $1\le \ell \le 4n+1$, $N=\Om\pf{n}{\ep^2}$, and $R=\Om\pa{\log\prc{\de}}$,
    with the polarized walk started at $\sign(y)$ with $\Om\pa{n\log\pf{n}{\ep}}$ steps, for $\de=e^{-O(n)}$, with probability $\ge 1-\de$, the following hold.
\begin{enumerate}
    \item For any $y\in \R^n$, 
    \[\hat Z_{\be'A, y}|_{\ball{\ep(\be) n}{\sign(y)}} \in [e^{-\ep}, e^{\ep}] \cdot  Z_{\be'A, y}|_{\ball{\ep(\be) n}{\sign(y)}}
    \subset [0,e^\ep] \cdot Z_{\be'A, y}
    \]
    \item For $y_T = Tx_0+B_T$ drawn according to SL, 
    \[\hat Z_{\be'A, y}|_{\ball{\ep(\be) n}{\sign(y)}} \in [e^{-\ep}, e^{\ep}] \cdot  Z_{\be'A, y}.\]
\end{enumerate}
\end{corollary}
\begin{prf} 
With probability $\ge 1-e^{-\Om(n)}$, $\opnorm{A}\le 4$. Then $g_\ell(\si) = e^{(\be_{\ell+1}-\be_\ell)\an{\si, A\si}}\in [e^{-1},e]$, so $\fc{\Var_{p_\ell}(g_\ell(\si))}{(\E_{p_\ell}g_\ell(\si))^2} \le e^4$. 
By \Cref{cor:wedge-sampling}, 
with probability $\ge 1-e^{-\Om(n)}$, 
the polarized walk started at $\sign(y)$ mixes to TV distance $\fc{\ep}{16n}$ within $O\pa{\fc{n}{\rh_\be}\log\pf{n}{\ep}}$ steps. Then by \pref{lem:sa}, with 
$N=\Om\pf{n}{\ep^2}$ samples and $R=\Om\pa{\log \prc{\de}}$ repetitions,
we obtain a $e^{\ep}$ factor approximation for (1).
For (2), apply (1) with $\fc{\ep}2$ to get $\hat Z_{\be'A, y}|_{\ball{\ep(\be) n}{\sign(y)}} \in [e^{-\ep/2}, e^{\ep/2}]\cdot Z_{\be'A, y}|_{\ball{\ep(\be) n}{\sign(y)}}$. 
It remains to note by \Cref{t:ls-sk}(1) that 
$Z_{\be'A, y}|_{\ball{\ep(\be) n}{\sign(y)}}\in [1-e^{-O(n)}, 1]\cdot Z_{\be'A, y}$.
\end{prf}

\section{Proof of main theorem}
\label{s:main-proof}
We prove the following formal version of \Cref{t:informal}.

\begin{theorem}\label{t:formal}
Given $\be<\rc 2$ and $\de>0$, 
with probability $1-\de$ over $A\sim \mathrm{GOE}(n)$, for any $\ep>e^{-cn}$ for appropriate constant $c$, 
\pref{alg:main} outputs a sample from a distribution $\hat \mu$ such that $\TV(\hat \mu, \mu_{\be A})\le \ep$ in time $\poly(n,\exp\pa{\rc \ep})$.
In particular, this gives a polynomial-time algorithm for sampling from the SK model for $\be<\rc 2$ with TV distance $o_n(1)$ in polynomial time.
\end{theorem}

We first prove the theorem assuming that in \pref{alg:asl-ta-je-dre}, $\td \mg_{t}$ is an exact solution of the TAP equation \eqref{e:TAP} with $\td y_t$ for each $t\in \{0,\ldots, T\}$. In \pref{s:mirror-tap}, we show that solving the TAP equation is efficient, and take into account the error involved.

\begin{prf}
We check the conditions of \Cref{l:je-rs} where the ideal process $y_t$ is \eqref{e:SL}, the approximate process is \eqref{e:asl-tap-je} with weights targeting $\rh_t$ given by \eqref{e:rho-t}. The algorithm computes the Euler-Maruyama discretization with some step size $h$ to be determined:
\begin{align*}
    \ty_{t+h} &= \ty_t + h \td{\mg}_t + \sqrt{h}\xi_t, \quad \xi_t \sim \calN(0,I_n), & \ty_0&=0\\
    \tw_{t+h} &= \tw_t + \fc h2 \ba{\Tr(\hQ(\hm(\td y_{t}))) + \ve{\hm(\td y_{t})}^2}. 
\end{align*}

\ppart{Part 1 (drift error)} This holds with $\ep(t) = O(1)$ by \Cref{l:m-error}. 

\ppart{Part 2 \& 3 (tails for log-weights)} Choose any $\lm>1$. By \ref{d:JE-warm-start} and \Cref{l:wt-tail} we have that the tail bounds hold with $C_1 = \pf{\CJET{\lm}}{\ep_1}^{\rc{\lm -1}}$ and $c_3 = \fc{\CJET{-1}\CJET{1}}{\ep_3}$.

\ppart{Part 4 (path error)} We use \Cref{l:sde-disc}. For this, we need a bound $\opnorm{D_y\hm} \le L$ and $\ve{\hm}_{\iy}\le M$, where here $\norm{\cdot}_{\infty}$ is the functional $L^{\infty}$ norm. Note 
\[\opnorm{D_y\hm} = \opnorm{\hQ(\hm(y))}\le \La\]
by \ref{d:Q-reg}, so we can take $L=\La = O(1)$. Note that $\ve{\hm}_{\iy}\le \sqrt N$,
so we can take $M=\sqrt n$. To obtain $\TV(\dist((\hy_t)_{t\in [0,T]}), \dist((\ty_t)_{t\in [0,T]}))\le \ep_4$ it suffices to take $h\le \fc{\sqrt{3/2}\ep_4}{\sqrt T LM}\wedge \fc{\ep_4^2}{2nTL^2}$, so it suffices for $h=O\pf{\ep_4^2}{n}$ with an appropriate constant.

\ppart{Part 5 (weight error)} We use \Cref{l:je-wt-disc}. For this, we need a bound $\ve{\gd \om}_2 \le L_\om$. 
First we bound $\ve{f(m)}_2$ where $f$ is defined in \eqref{e:Itomag}. Note that \ref{d:Q-reg} gives that $\hQ(m)^{-2}\preceq \La^2 D(m)^{-2}$ and considering the diagonal elements individually, $E_{\calD_n}\left[\hQ(\mg)^2\right]D(\mg)^2\preceq \La^2 \Id_n$.
Hence 
\begin{align*}
-\be^2 \La^2 \pa{1+\fc 2n}\Id_n&\preceq 
\left(E_{\calD_n}\left[\hQ(\mg)^2\right]D(\mg)^2 - \frac{\beta^2}{n}\left(\Tr\left[\hat{Q}^2(\mg)\right]\Id_n +2\hat{Q}^2(\mg)\right)\right)\preceq \La^2 \Id_n\\
\implies 
    \ve{f(m)}_2 &= \ve{\left(E_{\calD_n}\left[\hQ(\mg)^2\right]D(\mg)^2 - \frac{\beta^2}{n}\left(\Tr\left[\hat{Q}^2(\mg)\right]\Id_n +2\hat{Q}^2(\mg)\right)\right)\mg}_2\le \La^2 \sqrt n.
\end{align*}
By \eqref{e:d-om}, we have $\gd_{\hm} \om = \hm-f(\hm)$ so
\begin{align*}
    \ve{\gd_y\om }_2 &= \ve{(D_y\hm) (\hm - f(\hm))}_2\\
    &\le \La \cdot (\La^2+1)\sqrt n.
\end{align*}
Hence we can take $L_\om =\La(\La^2+1)\sqrt n = O(\sqrt n)$.
To obtain $\P\pa{|\tw_T((\ty_t)_{t\in [0,T]}) - \hw_T((\ty_t)_{t\in [0,T]})|\ge\fc{\ep_4}2}\le \fc{\ep_4}2$, it suffices to take $h\le \fc{\ep_4/2}{2L_\om TM} \wedge \fc{(\ep_4/2)^2}{2L_\om^2 T^2 n} \wedge \fc{1}{4L_\om \sqrt{TnC \log \prc{\ep_4/2}}}$ for appropriate constant $C$, so it suffices for $h= O\pf{\ep_4^2}{n^2}$ with an appropriate constant.

\ppart{Part 6 (ratio error)} The joint density of $(x, y_T)$ where $\si\sim \mu_{\be A}$ and $y_T = Tx_0+B_T$ is given by the SL process, relative to $\lm_{\mathrm{Unif}}\ot\lm_{\mathrm{Leb}}$ where $\lm_{\mathrm{Unif}}$ is the uniform measure on $\{\pm 1\}^n$ and $\lm_{\mathrm{Leb}}$ is Lebesgue measure on $\R^n$, is
\[
\mu(\si, y) \propto e^{\rc2 \be \ip{\si}{A\si}}e^{-\fc{\ve{y-T\si}^2}{2T}},
\]
so summing over $\si$ gives (letting $p_T=\dist(y_T)$)
\[
p_T(y) \propto \sum_{\si\in \{\pm 1\}^n} \mu(\si, y)
\propto \sum_{\si\in \{\pm 1\}^n} e^{\rc2 \be \ip{\si}{A\si}} e^{-\fc{\ve{y}^2}{2T} + \ip{y}{\si}}.
\]
Recalling that $\rh_T(y) \propto e^{\FT(\hat\mg(y),y) - \fc{\ve{y}^2}{2T}}$ from \eqref{e:rho-t}, 
\begin{align*}
    \dd{ p_T}{\rh_T}(y) \propto \fc{Z_{\be A, y}}{e^{\FT(\hat\mg(y),y)}} =:R_2(y)
\end{align*}
where $Z_{\be A, y} = \sum_{\si\in \{\pm 1\}^n} e^{\rc 2\be \an{\si,A\si} + \an{\si,y}}$. Hence, with 
\[
\td R_2(y) = \fc{\hat Z_{\be A, y}}{e^{\FT(\hat\mg(y),y)}}
\]
with $\hat Z_{\be A, y}$ calculated as in \pref{c:wedge-Z}, we have that $\fc{\td R_2(y)}{R_2(y)} = \fc{\hat Z_{\be A, y}}{Z_{\be A, y}} $ and the guarantees transfer:
\begin{align*}
    \forall y,\quad  \P\pa{\fc{\td R_2(y_T)}{R_2(y_T)}\nin [0, e^{\ep_4/2}] \Big| y_T=y} &\le \fc{\ep_4}2\\
    \P\pa{\fc{\td R_2(y_T)}{R_2(y_T)}\nin [e^{-\ep_4/2}, e^{\ep_4/2}]}&\le \ep
\end{align*}
by taking $\ep\mapsfrom \ep_4$.

\ppart{Putting everything together} Thus by \Cref{l:je-rs}, the output of Part 1 of \pref{alg:main} is a sample $\ty_T\sim \td  p_T$ with $\TV(\td p_T,  p_T)=O(\ep)$. This takes $O\pf{C_1^2L^2\log\prc{\ep}}{c_3^2}$ simulations, each of which takes time polynomial in $\rc{\ep_4}$ and $n$. Here $L=\exp\pa{\fc{3}{\ep}\pa{\fc 12 \int_0^T \ep(t)^2dt+1}} = \exp(O(1/\ep))$ and $\rc{\ep_4} = O\pf{L^2}{\ep} = \exp(O(1/\ep))$. By \pref{cor:wedge-sampling}, with probability $1-e^{-\Om(n)}$ over $y_T\sim  p_T$, the polarized walk for $\ce{\fc{n}{\rh_\be}\log\pf{4n}{\ep^2}}$ steps from $\sign(y)$ gives a sample $\hat \si$ that is $\ep$ in TV distance from $ \mu_{\be A, y_T}$. Letting $\td \mu$ be the distribution of $\hat\si$ for tilt $\ty_T$, we have by the chain rule for TV distance that
\begin{align*}
    \TV(\td \mu, \mu)
    &\le 
    \TV(\td  p_T,  p_T) + 
    \E_{y_T\sim  p_T} \TV (\calL(\hat\si | y_T), \calL(\si | y_T) ) = O(\ep) + e^{-\Om(n)} + \ep = O(\ep)\,,
\end{align*}
when $\ep=e^{-O(n)}$. Choosing $\ep = \Te{\prc{\log n}}$, we obtain TV distance $\ep=o(1)$ in polynomial time.
\end{prf}

\newpage
\addtocontents{toc}{\protect\setcounter{tocdepth}{-1}}

\section*{Acknowledgements}
\addtocontents{toc}{\protect\setcounter{tocdepth}{1}}

\iffocs{}{
HL and JSS are grateful to the organizers of the Rocky Mountain Summer Workshop on Combinatorics, Probability and Algorithms (2025), where the initial stages of this work were conducted.
This work and the workshop were supported by NSF grant 2309707.

JS and JSS thank David Jekel for pointing them to a particularly elegant strategy for the argument in \pref{sec:AIk+Bk}.}

\addtocontents{toc}{\protect\setcounter{tocdepth}{-1}}

\renewcommand{\baselinestretch}{0.928}\normalsize
{
    \small\hypersetup{urlcolor=Black}
    \bibliographystyle{alpha_beta_doi}
    \bibliography{main.bib}
}
\renewcommand{\baselinestretch}{1.0}\normalsize

\addtocontents{toc}{\protect\setcounter{tocdepth}{1}}

\newpage

\appendix
\normalsize

\section{Replica identities and higher-order remainder terms}\label{sec:app-a}
The goal of this section is three-fold: 
\begin{enumerate}[itemsep=0.3em]
    \item To collect a lemma which basically shows that the diagonals of the $\mathsf{GOE}(n)$ matrix plays no role in Gibbs averages by using the fact that $\sigma^2_i=1$, and
    \item Collect various rewrites of tracial terms that show up in the covariance estimates so that the cavity estimates can be applied in a straight-forward manner.
    \item Give $O(1)$ bounds for various higher-order moment terms that arise after applying Gaussian integration--by-parts to remainder terms in $\calE_A$, as well as a stand-alone lemma that controls the rank-1 term in $\calE_A$ using similar estimates.
\end{enumerate}

\subsection{Independence of the Gibbs covariance from diagonal entries} While $A \sim \mathsf{GOE}(n)$ has diagonal entries with variance $2/n$, it is easy to show that any Gibbs statistics for the planted SK model with SL tilt is ``independent'' of the Gibbs because of the fact that $\sigma^2_i=1$.

\begin{lemma}[Diagonal GOE entries do not affect Gibbs averages]\label{lem:diag-does-not-enter-gibbs}
Let $A$ be symmetric and consider the planted SK Hamiltonian
\[
\beta H_n(\sigma)
=\frac{\beta}{2}\langle \sigma,A\sigma\rangle
+\frac{\beta^2}{2n}\Big(\sum_{k=1}^n\sigma_k\Big)^2
+\langle y_t,\sigma\rangle,
\qquad \sigma\in\{-1,+1\}^n.
\]
Let $A^{\circ}:=A-\diag(A)$ be the matrix obtained by zeroing the diagonal of $A$.
Then, for every bounded observable $F:\{-1,+1\}^n\to\mathbb R$,
\begin{equation}\label{eq:gibbs-diag-invariant}
\langle F\rangle_{A,y_t,\beta}=\langle F\rangle_{A^{\circ},y_t,\beta}.
\end{equation}
In particular, the magnetization vector $m=\langle\sigma\rangle$ and covariance matrix
$P=\Cov(G_{A,y_t,\beta})$ depend on $A$ only through its off-diagonal entries:
\[
m=m(A^{\circ},y_t),\qquad P=P(A^{\circ},y_t).
\]
Consequently, for every $i \in [n]$ and every component $m_k$ and $P_{k\ell}$,
\[
\partial_{A_{ii}}m_k=0,\qquad \partial_{A_{ii}}P_{k\ell}=0,
\]
and the same holds for any matrix/function built from $m$ and $P$. 
\end{lemma}

\begin{proof}
Since $A$ is symmetric,
\[
\langle\sigma,A\sigma\rangle=\sum_{i,j=1}^n A_{ij}\sigma_i\sigma_j
=2\sum_{1\le i<j\le n}A_{ij}\sigma_i\sigma_j+\sum_{i=1}^n A_{ii}\sigma_i^2.
\]
But $\sigma_i^2\equiv 1$, so
\[
\frac{\beta}{2}\langle\sigma,A\sigma\rangle
=\beta\sum_{1\le i<j\le n}A_{ij}\sigma_i\sigma_j+\frac{\beta}{2}\sum_{i=1}^n A_{ii}.
\]
The last term $\frac{\beta}{2}\sum_i A_{ii}$ is independent of $\sigma$. Therefore, writing
\[
\Phi(\sigma;A^{\circ},y_t):=
\beta\sum_{i<j}A_{ij}\sigma_i\sigma_j
+\frac{\beta^2}{2n}\Big(\sum_{k=1}^n\sigma_k\Big)^2
+\langle y_t,\sigma\rangle,
\]
we have
\[
e^{\beta H_n(\sigma)}
=
\exp\Big(\frac{\beta}{2}\sum_{i=1}^n A_{ii}\Big)\cdot e^{\Phi(\sigma;A^{\circ},y_t)}.
\]
Therefore, the partition function factorizes as
\[
Z(A,y_t)=\sum_{\sigma}e^{\beta H_n(\sigma)}
=
\exp\Big(\frac{\beta}{2}\sum_{i=1}^n A_{ii}\Big)\sum_{\sigma}e^{\Phi(\sigma;A^{\circ},y_t)}
=
\exp\Big(\frac{\beta}{2}\sum_{i=1}^n A_{ii}\Big)\,Z(A^{\circ},y_t).
\]
This implies that for any bounded $F$,
\[
\langle F\rangle_{A,y_t,\beta}
=
\frac{1}{Z(A,y_t)}\sum_{\sigma}F(\sigma)e^{\beta H_n(\sigma)}
=
\frac{1}{Z(A^{\circ},y_t)}\sum_{\sigma}F(\sigma)e^{\Phi(\sigma;A^{\circ},y_t)}
=
\langle F\rangle_{A^{\circ},y_t,\beta},
\]
which proves \eqref{eq:gibbs-diag-invariant}. Since the RHS does not depend on $(A_{ii})$, all derivatives
$\partial_{A_{ii}}\langle F\rangle$ vanish identically. Applying this with $F(\sigma)=\sigma_k$ and $F(\sigma)=\sigma_k\sigma_\ell$
gives $\partial_{A_{ii}}m_k=0$ and $\partial_{A_{ii}}P_{k\ell}=0$. The same holds for any composition
$D(m),c(m),V(m),S(m)$.
\end{proof}

\subsection{Replica simplifications for cross terms via tracial identities and index algebra}
Below is a helper lemma that accounts for cancellations in the multi-overlaps that come from an expansion using replica-based identities for the $\E\left[\Tr[PA^2P]\right]$ term. This is a (minor) generalization of the computations in~\cite[Pg.~6 to 9]{el2024bounds} for the $\mathrm{III}$ term~\cite[Pg.~6]{el2024bounds} to account for the case where the covariance is (potentially) uncentered. 

\begin{lemma}[Exact $P$-pairing decomposition for $\partial_{ij}P_{ab}$]\label{lem:dP-pairing}
Fix $i<j$. For all $a,b\in\{1,\dots,n\}$,
\begin{equation}\label{eq:dP-pairing}
\partial_{ij}P_{ab}
=
\beta\Big(P_{ai}P_{bj}+P_{aj}P_{bi}\Big)
+\beta\Big(m_i\,T_{abj}+m_j\,T_{abi}+\Gamma_{abij}\Big).
\end{equation}
\end{lemma}

\begin{proof}
Start from the exact Gibbs derivative identity that follows by the division rule
\[
\partial_{ij}\langle F\rangle
=
\beta\Big(\langle F\,\sigma_i\sigma_j\rangle-\langle F\rangle\langle\sigma_i\sigma_j\rangle\Big)
=\beta\,\langle F\,(\sigma_i\sigma_j-\langle\sigma_i\sigma_j\rangle)\rangle.
\]
Applying this to $F=\tilde\sigma_a\tilde\sigma_b$ and recalling $P_{ab}=\langle\tilde\sigma_a\tilde\sigma_b\rangle$ gives
\begin{align}\label{eq:dP-start}
\partial_{ij}P_{ab}
&=\partial_{ij}\langle\sigma_a\sigma_b\rangle
-(\partial_{ij}m_a)m_b-m_a(\partial_{ij}m_b)\nonumber\\
&=\beta\Big(\langle\sigma_a\sigma_b\sigma_i\sigma_j\rangle-\langle\sigma_a\sigma_b\rangle\langle\sigma_i\sigma_j\rangle\Big)
-\beta m_b\Big(\langle\sigma_a\sigma_i\sigma_j\rangle-m_a\langle\sigma_i\sigma_j\rangle\Big)
-\beta m_a\Big(\langle\sigma_b\sigma_i\sigma_j\rangle-m_b\langle\sigma_i\sigma_j\rangle\Big).
\end{align}
Note that, expanding $\sigma_i\sigma_j=(m_i+\tilde\sigma_i)(m_j+\tilde\sigma_j)=m_im_j+m_i\tilde\sigma_j+m_j\tilde\sigma_i+\tilde\sigma_i\tilde\sigma_j$.
Also $\langle\sigma_i\sigma_j\rangle=m_im_j+P_{ij}$.
Inserting these into \eqref{eq:dP-start} cancels the $m_im_j$ term to finally give
\[
\partial_{ij}P_{ab}
=
\beta\Big(
m_i\langle\tilde\sigma_a\tilde\sigma_b\tilde\sigma_j\rangle
+m_j\langle\tilde\sigma_a\tilde\sigma_b\tilde\sigma_i\rangle
+\langle\tilde\sigma_a\tilde\sigma_b\tilde\sigma_i\tilde\sigma_j\rangle
-P_{ab}P_{ij}
\Big).
\]
Finally, note that the definition of $\Gamma_{abij}$ gives
\[
\langle\tilde\sigma_a\tilde\sigma_b\tilde\sigma_i\tilde\sigma_j\rangle-P_{ab}P_{ij} = P_{ai}P_{bj}+P_{aj}P_{bi}+\Gamma_{abij}\, ,
\]
and substituting this into the expression for $\partial_{ij} P_{ab}$ yields \eqref{eq:dP-pairing}.
\end{proof}


\begin{lemma}[Exact replica representations for tracial functions]\label{lem:trace-to-replica}
Let $G$ be the Gibbs measure, $P=\mathrm{Cov}(G)$ and $T=\diag(P)$.
For replicas $\sigma^1,\sigma^2,\sigma^3,\sigma^4\sim G^{\otimes 4}$, define
\[
a_i:=(\sigma_i^1-\sigma_i^2)(\sigma_i^3-\sigma_i^4),
\qquad
f:=\frac1n\sum_{i=1}^n a_i,
\qquad
\hat\sigma_i^\ell:=\frac{\sigma_i^\ell-m_i}{\sqrt{1-m_i^2}},
\qquad
\hat f:=\frac1n\sum_{i=1}^n (\hat\sigma_i^1-\hat\sigma_i^2)(\hat\sigma_i^3-\hat\sigma_i^4).
\]
Also, define the ``diagonally weighted'' rectangular sum
\[
f^{\mathrm{diag}}:=\frac1n\sum_{i=1}^n P_{ii}\,a_i.
\]
Then, the following tracial equalities hold conditioned on any choice of $\{A_{ij}\}_{i,j}$,
\begin{align}
\Tr[P^2] &= \frac{n^2}{4}\,\langle f^2\rangle,\label{eq:TrP2-f2}\\
\Tr[PDP] &= \frac{n^2}{4}\,\langle f\hat f\rangle,\label{eq:TrPDP-ff-hat} \\
\frac1n\sum_{i=1}^n \langle a_i^2\rangle &= \frac{4}{n}\sum_{i=1}^n P_{ii}^2. \label{eq:diag-a2} \\
\Tr[TP^2] &= \frac{n^2}{4}\,\langle f\,f^{\mathrm{diag}}\rangle, \label{eq:TrTP2-ffdiag}\\
\Tr[TPDP] &= \frac{n^2}{4}\,\langle f^{\mathrm{diag}}\,\hat f\rangle. \label{eq:TrTPDP-fdiagfhat}
\end{align}
\end{lemma}

\begin{proof}
For \eqref{eq:TrP2-f2}, note that
\[
P_{ij}=\langle\sigma_i\sigma_j\rangle-m_im_j
=\frac12\langle(\sigma_i^1-\sigma_i^2)(\sigma_j^1-\sigma_j^2)\rangle.
\]
This immediately gives
\[
P_{ij}^2
=\frac14\langle(\sigma_i^1-\sigma_i^2)(\sigma_j^1-\sigma_j^2)\rangle
\langle(\sigma_i^3-\sigma_i^4)(\sigma_j^3-\sigma_j^4)\rangle,
\]
and by the replica product identity,
\[
P_{ij}^2=\frac14\langle(\sigma_i^1-\sigma_i^2)(\sigma_i^3-\sigma_i^4)
(\sigma_j^1-\sigma_j^2)(\sigma_j^3-\sigma_j^4)\rangle.
\]
Summing over $i,j$ gives
\[
\Tr[P^2]=\sum_{i,j}P_{ij}^2
=\frac14\Big\langle\Big(\sum_i a_i\Big)\Big(\sum_j a_j\Big)\Big\rangle
=\frac{n^2}{4}\langle f^2\rangle.
\]

Using $D_{jj}=\frac1{1-m_j^2}$, one obtains
\[
\Tr[PDP]=\sum_{i,j}P_{ij}^2D_{jj}.
\]
Insert the same replica representation for $P_{ij}^2$ as above and pull $D_{jj}$ inside by rewriting
\[
D_{jj}(\sigma_j^1-\sigma_j^2)(\sigma_j^3-\sigma_j^4)
=
\frac{(\sigma_j^1-m_j)-(\sigma_j^2-m_j)}{\sqrt{1-m_j^2}}
\cdot
\frac{(\sigma_j^3-m_j)-(\sigma_j^4-m_j)}{\sqrt{1-m_j^2}}.
\]
This is exactly $(\hat\sigma_j^1-\hat\sigma_j^2)(\hat\sigma_j^3-\hat\sigma_j^4)$.
Therefore
\[
\Tr[PDP]
=\frac14\Big\langle\Big(\sum_i a_i\Big)\Big(\sum_j \hat a_j\Big)\Big\rangle
=\frac{n^2}{4}\langle f\hat f\rangle,
\]
which proves \eqref{eq:TrPDP-ff-hat}.

For \eqref{eq:diag-a2}, condition on disorder and use the independence of the replicas under the Gibbs to obtain
\[
\langle a_i^2\rangle
=\langle(\sigma_i^1-\sigma_i^2)^2\rangle\;\langle(\sigma_i^3-\sigma_i^4)^2\rangle.
\]
However, using the fact that $(\sigma_i^1-\sigma_i^2)^2=2(1-\sigma_i^1\sigma_i^2)$ we simplify further to
\[
\langle(\sigma_i^1-\sigma_i^2)^2\rangle
=2\big(1-\langle\sigma_i^1\sigma_i^2\rangle\big)
=2(1-m_i^2)=2P_{ii}.
\]
Since $\langle a_i^2\rangle=4P_{ii}^2$, summing and dividing by $n$ gives \eqref{eq:diag-a2}.

The expressions for~\eqref{eq:TrTP2-ffdiag} and~\eqref{eq:TrTPDP-fdiagfhat} follow via similar algebraic simplifications. 
\begin{align*}
\Tr[TP^2]
&=\sum_i P_{ii}(P^2)_{ii}
=\sum_{i,j}P_{ii}P_{ij}^2\\
&=\frac14\sum_{i,j}P_{ii}\left\langle(\sigma_i^1-\sigma_i^2)(\sigma_j^1-\sigma_j^2)\right\rangle^2\\
&=\frac14\sum_{i,j}P_{ii}\left\langle(\sigma_i^1-\sigma_i^2)(\sigma_i^3-\sigma_i^4)(\sigma_j^1-\sigma_j^2)(\sigma_j^3-\sigma_j^4)\right\rangle\\
&=\frac14\left\langle\left(\sum_i P_{ii}a_i\right)\left(\sum_j a_j\right)\right\rangle
=\frac14\,\langle (n f^{\mathrm{diag}})\,(n f)\rangle
=\frac{n^2}{4}\,\langle f\,f^{\mathrm{diag}}\rangle\, ,
\end{align*}
and
\[
\Tr[TPDP]=\sum_{i,j}P_{ii}P_{ij}^2\frac1{1-m_j^2}
=\frac14\Big\langle \Big(\sum_i P_{ii}a_i\Big)\Big(\sum_j\hat a_j\Big)\Big\rangle
=\frac{n^2}{4}\,\langle f^{\mathrm{diag}}\,\hat f\rangle\, ,
\]
where the proof is identical to the previous step, with one $a_j$ replaced by $\hat a_j$. \qedhere
\end{proof}


\begin{lemma}[Index trace identities for $P$ and $D$]\label{lem:trace-identities}
Let $m\in(-1,1)^n$ and let $P\in\mathbb R^{n\times n}$ be symmetric with
\[
P_{ij}=\langle \sigma_i\sigma_j\rangle-\langle\sigma_i\rangle\langle\sigma_j\rangle,
\qquad
P_{ii}=1-m_i^2,
\qquad
D:=\mathrm{diag}\Big(\frac1{1-m_1^2},\dots,\frac1{1-m_n^2}\Big)
=\mathrm{diag}\Big(\frac1{P_{11}},\dots,\frac1{P_{nn}}\Big).
\]
Then the following identities hold:

\begin{enumerate}[itemsep=0.35em]
\item[(i)] 
\[
\Tr[P^2]=\sum_{i,j=1}^n P_{ij}^2
=\sum_{i=1}^n P_{ii}^2+2\sum_{i<j}P_{ij}^2.
\]

\item[(ii)]
\[
\Tr[P]^2=\Big(\sum_{i=1}^n P_{ii}\Big)^2
=\sum_{i=1}^n P_{ii}^2+2\sum_{i<j}P_{ii}P_{jj}.
\]
Equivalently,
\[
\sum_{i<j}P_{ii}P_{jj}=\frac12\Big(\Tr(P)^2-\sum_i P_{ii}^2\Big).
\]

\item[(iii)]
\[
\Tr[PDP]=\sum_{i,j=1}^n P_{ij}D_{jj}P_{ji}
=\sum_{j=1}^n D_{jj}\sum_{i=1}^n P_{ij}^2
=\Tr(P)+\sum_{j=1}^n D_{jj}\sum_{i\neq j}P_{ij}^2.
\]
Therefore,
\begin{equation}\label{eq:TPDP-minus-TP}
\Tr[PDP]-\Tr[P]=\sum_{j=1}^n D_{jj}\sum_{i\neq j}P_{ij}^2
=\sum_{i<j}(D_{ii}+D_{jj})P_{ij}^2.
\end{equation}
\end{enumerate}
\end{lemma}

\begin{proof}
The proof for (i) and (ii) are immediate expansions. For (iii), since $D$ is diagonal,
\[
\Tr[PDP]=\sum_{i,j}P_{ij}D_{jj}P_{ji}=\sum_j D_{jj}\sum_i P_{ij}^2.
\]
Splitting the inner sum into $i=j$ and $i\neq j$ components gives that
\[
\Tr[PDP]=\sum_j D_{jj}P_{jj}^2+\sum_j D_{jj}\sum_{i\neq j}P_{ij}^2.
\]
However, since $D_{jj}P_{jj}^2=P_{jj}$, this simplifies to give
\[
\Tr[PDP]=\sum_j P_{jj}+\sum_j D_{jj}\sum_{i\neq j}P_{ij}^2=\Tr(P)+\sum_j D_{jj}\sum_{i\neq j}P_{ij}^2.
\]
Finally, rewrite the last double sum symmetrically yields
\[
\sum_j D_{jj}\sum_{i\neq j}P_{ij}^2=\sum_{i<j}\big(D_{ii}+D_{jj}\big)P_{ij}^2\,,
\]
which is the desired identity.
\end{proof}

\subsection{Bounds for the $\mathrm{HO}$ terms} In this subsection we will state a lemma that will tersely illustrate the computations showing that each $\mathrm{HO}_i$ term is $O(1)$. We do this using (rather repetitively) the same replica rewrites and summing over rectangular sums to turn every $\mathrm{HO}_i$ term to be equivalent to $\mathsf{span}\left\{\nu(fG), \nu(f_\alpha f_\beta f_\gamma f_\delta)\right\}$ or a term that decays even faster, and then applying some combination of~\pref{prop:D15} or H\"older's inequalities in conjunction with~\pref{thm:overlap-moment-concentration} and~\pref{lem:replicon-moments-from-mgf}.

\begin{lemma}[$O(1)$ bounds for the $\mathrm{HO}$ terms]\label{lem:bounds-for-ho-terms}
    Let $\mathrm{HO}_i$ for $i \in \{1,\dots,16\}$ be as defined in~\pref{lem:remainder-bound-ea}. Then, 
    \[
        \mathrm{HO}_i = O(1)\, ,
    \]
    for every $i$.
\end{lemma}
\begin{prf}
    We prove the bounds for every term in chronological order. Each proof (roughly) follows the same pattern -- use replica rewrites for the centered third-moment or fourth cumulant terms, perform the inner summation, symmetrize over the outer sums, and then use some combination of~\pref{thm:overlap-moment-concentration},~\pref{lem:replicon-moments-from-mgf} or~\pref{prop:D15} in conjunction with H\"older's inequality to deal with ``grouped sum'' terms, and deal with the diagonal terms using established bounds on trace powers of $P$.
    \ppart{$\mathrm{HO}_1$ term} Using replica rewrites for the third-centered moment and symmetrization of sums
    \allowdisplaybreaks
    \begin{align*}
        &\frac{2\beta^4}{n^2}\E\left[\sum_{i<j}m_i\sum_k(P^2)_{kj}T_{ikj}\right] = \frac{C(\beta)}{n}\sum_{i<j}\E\left[m_i\sum_k \an{\Delta^{12}_k\Delta^{34}_jf_{1234}\Delta^{56}_i\Delta^{57}_k\Delta^{56}_j}\right] \\
        &\qquad\qquad = C(\beta)n^2 \nu((R_{57}-R_{67})f_{3456}f_{1234}f_{1257}) - C'(\beta)\nu\left(f_{1234}f_{1257}\left(\sum_i \sigma^7_i(\Delta^{56}_i)^2\Delta^{34}_i\right)\right)\,.
    \end{align*}
    The first term is shown to be $O(1)$ by a $(4,4,4,4)$-H\"older's in conjunction with~\pref{thm:overlap-moment-concentration} and~\pref{lem:replicon-moments-from-mgf}, while the second term is $O(1)$ by observing that $|\sum_i\sigma^7_i(\Delta^{56}_i)^2\Delta^{34}_i| \le 80 n$ and using Cauchy--Schwarz in conjunction with~\pref{lem:replicon-moments-from-mgf}.

    \ppart{$\mathrm{HO}_2$ term}  Using replica rewrites with $T_{iki}=-2m_iP_{ik}$ and symmetrization of sums again gives
    \allowdisplaybreaks
    \begin{align*}
        \frac{2\beta^4}{n^2}\E\left[\sum_{i<j}m_j\sum_k (P^2)_{kj}T_{iki}\right] &= \frac{C(\beta)}{n^2}\E\left[-\mg^\sT P^3 \mg +\Tr[\diag(\mg^2)P^3]\right] \,.
    \end{align*}
    Now $|\mg^\sT P^3 \mg| \le \Tr[P^3]\norm{\mg}^2_2 \le n \Tr[P^3]$ which immediately implies that $\frac{C(\beta)}{n^2}\E \mg^\sT P^3 \mg = O(1)$ via an application of \eqref{eq:c-p3-bound}. For the second term, note that the proof of~\pref{lem:remainder-bound-ea} already demonstrates it is $O(1)$.

    \ppart{$\mathrm{HO}_3$ term} Using the fact that $\Gamma_{ikij} = -2P_{ki}P_{ij} -2m_iT_{ikj}$ the term simplifies down to
    \allowdisplaybreaks
    \begin{align*}
        \frac{2\beta^4}{n^2}\E\left[\sum_{i<j}\sum_k (P^2)_{kj}\Gamma_{ikij}\right] &= \frac{C(\beta)}{n^2}\sum_{i<j}\E\left[\sum_k(P^2)_{kj}P_{ki}P_{ij}\right] + \frac{C(\beta)}{n^2}\E\left[\sum_{i<j}m_i\sum_k (P^2)_{kj}T_{ikj}\right] \\
        &= \frac{C(\beta)}{n^2}\E\left[\Tr[P^4]\right] + \frac{C(\beta)}{n^2}\E\left[\sum_i P_{ii}(P^3)_{ii}\right] + \mathrm{HO}_1\,. 
    \end{align*}
    The first term is $O(1)$ immediately by using \eqref{eq:c-p4-bound}, while the third term is $O(1)$ by our previous bound. For the second term, observe that $0 \le P_{ii} \le 1$ and so $\sum_i P_{ii}(P^3)_{ii} \le \Tr[P^3]$, at which point $\frac{C(\beta)}{n^2}\E\left[\sum_i P_{ii}(P^3)_{ii}\right] = O(1)$ by \eqref{eq:c-p3-bound}.

    \ppart{$\mathrm{HO}_4$ term} Using the identities that $T_{k\ell j} = \an{\Delta^{12}_k\Delta^{13}_\ell \Delta^{12}_j}$, $P_{\ell j} = \an{\Delta^{45}_\ell\Delta^{45}_j}$ and $P_{ik} = \an{\Delta^{67}_i\Delta^{67}_k}$ we obtain
    \allowdisplaybreaks
    \begin{align*}
        \frac{2\beta^4}{n^2}\E\left[\sum_{i<j}m_i\sum_{k,\ell} P_{ik}P_{\ell j}T_{k\ell j}\right] &= \frac{C(\beta)}{n}\sum_{i<j}\sum_k\an{\sigma^8_i\Delta^{67}_i\Delta^{67}_k\Delta^{12}_k\Delta^{45}_j\Delta^{12}_jf_{1345}} \\
        &= C(\beta)n^2 \nu\left((R_{68}-R_{78})f_{1245}f_{1267}f_{1345}\right) - C(\beta)\nu(B^n f_{1345}f_{1267})\, ,
    \end{align*}
    where $B^n$ is some function of the spins that is absolutely bounded by $80 n$. Both these terms can be shown to be $O(1)$ as they are of exactly the same form as the terms that $\mathrm{HO}_1$ evaluates to.

    \ppart{$\mathrm{HO}_5$ term} For the fifth term, we use the fact that $\Gamma_{k\ell i j} = \an{\Delta^{12}_k\Delta^{13}_\ell\Delta^{14}_i\Delta^{15}_j} - P_{k\ell}P_{ij}-P_{ki}P_{\ell j} - P_{kj}P_{\ell i}$ and symmetrize sums to get an expression involving bulk deviations products, and three auxiliary trace power terms. More specifically,
    \allowdisplaybreaks
    \begin{align*}
        &\frac{2\beta^4}{n^2}\E\left[\sum_{i<j}\sum_{k,\ell} P_{ik}P_{\ell j}\Gamma_{k\ell ij}\right] = C(\beta)n^2\nu\left(f_{1267}f_{1389}f_{1589}f_{1467}\right) - C(\beta)\nu\left(B^n f_{1267}f_{1389}\right) \\
        &\qquad\qquad\qquad - \frac{C(\beta)}{n^2}\E\left[\sum_{i<j}P_{ij}(P^3)_{ij} - \sum_{i<j}(P^2)_{ii}(P^2)_{jj} - \sum_{i<j}(P^2)^2_{ij}\right]\,,
    \end{align*}
    where $B^n$ is again some function of the spins that is absolutely bounded by $100n$. The first two terms are bounded by $O(1)$ using the same strategy used on the prior terms, while the remaining terms can be symmetrized and written as
    \allowdisplaybreaks
    \begin{align*}
        \sum_{i<j} P_{ij}(P^3)_{ij} &= \frac{1}{2}\Tr[P^4] - \frac{1}{2}\sum_i P_{ii}(P^3)_{ii}\,, \\
        \sum_{i<j}(P^2)_{ii}(P^2)_{jj} &= \frac{1}{2}\Tr[P^2]^2 - \frac{1}{2}\sum_i (P^2_{ii})^2\,, \\
        \sum_{i<j} (P^2)^2_{ij} &= \frac{1}{2}\Tr[P^4] - \frac{1}{2}\sum_i (P^2)_{ii}^2\, ,
    \end{align*}
    which can be bounded by observing that $\frac{C(\beta)}{n^2}\E[\Tr[P^4]] = O(1)$ by \eqref{eq:c-p4-bound}, that $|\sum_i P_{ii}(P^3)_{ii}| \le \Tr[P^3]$ and so $\frac{C(\beta)}{n^2}\E[\Tr[P^3]] = O(1)$ by \eqref{eq:c-p3-bound}, that $\frac{C(\beta)}{n^2}\E\left[\Tr[P^2]^2\right] = \frac{n^2}{16}\nu(f_{1234}^2f_{5678}^2) = O(1)$, and that $\sum_i (P^2)^2_{ii} \le \Tr[P^4]$ which immediately yields $\frac{C(\beta)}{n^2}\sum_i (P^2)^2_{ii}\le O(1)$ by \eqref{eq:c-p4-bound}. Putting all the bounds together immediately gives $\mathrm{HO}_5 = O(1)$.

    \ppart{$\mathrm{HO}_6$, $\mathrm{HO}_7$ and $\mathrm{HO}_8$ terms} The proof for this is (essentially) an easier variant of the proof used to bound $R_2$ in~\pref{lem:remainder-bound-ea}. Observe that $0 \le P_{ii} + P_{jj} \le 2$, and then note that
    \allowdisplaybreaks
    \begin{align*}
        &\mathrm{HO}_6 + \mathrm{HO}_7 + \mathrm{HO}_8 = \frac{2\beta^4}{n^2}\sum_{i<j}\E\left[(P_{ii} + P_{jj})\left(m_i\sum_k P_{kj}T_{ikj} + m_j\sum_k P_{kj}T_{iki} + \sum_k P_{kj}\Gamma_{ikij}\right)\right] \\
        &= \frac{2\beta^4}{n^2}\sum_{i<j}\E\left[(P_{ii}+P_{jj})\left(\sum_k P_{kj}(m_iT_{ikj} + m_jT_{iki} + P_{kj}\Gamma_{ikij})\right)\right] \\
        &= \frac{2\beta^4}{n^2}\sum_{i<j}\E\left[(P_{ii}+P_{jj})\left(\sum_kP_{kj}(-m_iT_{ikj}-2m_im_jP_{ik} -2P_{ik}P_{ij})\right)\right] \\
        &= \frac{C(\beta)}{n^2}\sum_{i<j}\E\left[(P_{ii}+P_{jj})m_im_j(P^2)_{ij}\right] + \frac{C(\beta)}{n^2}\sum_{i<j}\E\left[(P_{ii}+ P_{jj})P_{ij}(P^2)_{ji}\right] + \frac{C(\beta)}{n^2}\sum_{i<j}\E\left[m_i(P_{ii}+P_{jj})\sum_kP_{jk}T_{ikj}\right]\,.
    \end{align*}
    We control each of these terms individually. For the first one, after symmetrization of sums, we obtain the following
    \allowdisplaybreaks
    \begin{align*}
        \frac{C(\beta)}{n^2}\sum_{i<j}\E\left[(P_{ii}+P_{jj})m_im_j(P^2)_{ij}\right] &= \frac{C(\beta)}{n^2}\E\left[\mg^\sT E_{\calD_n}[P]P^2\mg\right] - \frac{C(\beta)}{n^2}\E\left[\sum_i \mg^2_iP_{ii}(P^2)_{ii}\right]\,.
    \end{align*}
    Since $\sum_i \mg_i (P_{ii}(P^2)_{ii}) \le_{|\mg_i| \le 1} \sum_i P_{ii}(P^2)_{ii}\le_{|P_{ii}| \le 1} \sum_{i}(P^2)_{ii} = \Tr[P^2]$, we have $\frac{C(\beta)}{n^2}\E\left[\sum_i \mg_i(P_{ii}(P^2)_{ii})\right] = O(1)$ by invoking \eqref{eq:c-p2-bound}. Similarly, $\mg^\sT E_{\calD_n}[P]^{1/2}P^2E_{\calD_n}[P]^{1/2}\mg \le\norm{P^2}_F \norm{E_{\calD_n}[P]^{1/2}\mg}^2_2$. A $(1,\infty)$-H\"older's gives that $\norm{E_{\calD_n}[P]^{1/2}\mg}^2_2 \le_{\mg^2_i \le 1} \Tr[P] \le n$ and $\norm{P^2}_F = \sqrt{\Tr[P^4]}$. This immediately implies that $\frac{C(\beta)}{n^2}\E[\mg^\sT E_{\calD_n}[P]P^2\mg] \le \frac{C(\beta)}{n}\E\left[\sqrt{\Tr[P^4]}\right]$ which is $O(1)$ by Jensen's in conjunction with \eqref{eq:c-p4-bound}. \\
    For the second term, observe that
    \allowdisplaybreaks
    \begin{align*}
        \frac{C(\beta)}{n^2}\sum_{i<j}\E\left[(P_{ii}+ P_{jj})P_{ij}(P^2)_{ji}\right] &= \frac{C(\beta)}{n^2}\sum_{i,j}\E\left[P_{ii}P_{ij}(P^2)_{ij}\right] - \frac{C(\beta)}{n^2}\sum_i\E\left[(P_{ii})^2(P^2)_{ii}\right] \\
        &= \frac{C(\beta)}{n^2}\E\left[\sum_i P_{ii}(P^3)_{ii}\right] - \frac{C(\beta)}{n^2}\sum_i\E\left[(P_{ii})^2(P^2)_{ii}\right]\,.
    \end{align*}
    Using the fact that $0 \le P_{ii} \le 1$ and $P\succeq 0$ (again) gives that $\left|\sum_i P_{ii}(P^3)_{ii}\right| \le \Tr[P^3]$ which immediately yields $\frac{C(\beta)}{n^2}\sum_{i,j}\E\left[P_{ii}P_{ij}(P^2)_{ij}\right] = O(1)$ by using \eqref{eq:c-p3-bound}. For the second term, we use a similar strategy with the invocation that $(P_{ii})^2 \le 1$ in conjunction with \eqref{eq:c-p2-bound}. This gives $\frac{C(\beta)}{n^2}\sum_{i<j}\E\left[(P_{ii}+ P_{jj})P_{ij}(P^2)_{ji}\right] = O(1)$. \\
    For the last term, the following rewrite is immediately obvious using the replica identities and some algebra once we isolate the sum over the $k$ index. Specifically,
    \allowdisplaybreaks
    \begin{align*}
        \left|\sum_k P_{jk}T_{ikj}\right| &= \sum_k \an{(\sigma_i-\mg_i)(\sigma_j-\mg_j)\sum_k P_{jk}(\sigma_k-\mg_k)} \\
        &\le_{\text{CS}} \sqrt{\an{(\sigma_i-\mg_i)^2(\sigma_j-\mg_j)^2}^{1/2}\an{\left(\sum_k P_{kj}(\sigma_k-\mg_k)\right)\left(\sum_\ell P_{j \ell}(\sigma_\ell - \mg_\ell)\right)}^{1/2}} \\
        &\le_{|P_{ij}| \le 4,\,P_{k\ell}=(\sigma_k-\mg_k)(\sigma_\ell-\mg_\ell)} 4\sqrt{(P^3)_{jj}}\,. 
    \end{align*}
    Using the fact that $|m_i(P_{ii}+P_{jj})| \le 4$ and Jensen's inequality, we immediately obtain that
    \[
        \frac{C(\beta)}{n^2}\sum_{i<j}\E\left[m_i(P_{ii}+P_{jj})\sum_kP_{jk}T_{ikj}\right] \le \frac{C(\beta)}{n}\sum_j\sqrt{\E\left[(P^3)_{jj}\right]} \le_{\text{CS}} C(\beta)\sqrt{\frac{\E\left[\Tr[P^3]\right]}{n}} \le_{\text{\eqref{eq:c-p3-bound}}} O(1)\,. 
    \]

    \ppart{$\mathrm{HO}_{9}$ and $\mathrm{HO}_{10}$ terms}  The proofs for both these terms are very similar to the bound for $\mathrm{HO}_8$, so we keep it concise.  Since 
    \allowdisplaybreaks
    \begin{align*}
        \E\left[\sum_{ij}m_iP_{ij}\sum_kP_{ik}T_{ikj}\right] &\le_{\text{CS}} \sqrt{\E\left[\sum_{i,j}m_i^2(P_{ij})^2\right]}\sqrt{\E\left[\sum_{i,j}\left(\sum_k P_{ik}T_{ijk}\right)^2\right]} \\
        &\le_{|m_i^2|\le 1,\,\text{CS}}\sqrt{\E\left[\Tr[P^2]\right]}\sqrt{\E\left[\sum_{i,j}(\sum_k(P_{ik})^2)(\sum_k T^2_{ijk})\right]} \\
        &\le_{\text{\eqref{eq:c-p2-bound}},\,|T^{2}_{ijk}|\le32} C\,\sqrt{n}\sqrt{\E\left[\sum_{i,j}(P^2)_{ii}n\right]} \\
        &\le C\sqrt{n^{3/2}}\sqrt{\E\left[\Tr[P^2]\right]} \le_{\text{\eqref{eq:c-p2-bound}}} O(n^2)\,, 
    \end{align*}
    which immediately implies that $\mathrm{HO}_9 = O(1)$, since the diagonal term is easily shown to be $O(1/n)$ by using the fact that $T_{iki}=-2m_iP_{ki}$, and summing over $i=j$ gives a bound on the order of $\frac{C(\beta)}{n^2}\E\left[\Tr[P^2]\right] = O(1/n)$. \\
    For the $\mathrm{HO}_{10}$ term, observe that 
    \allowdisplaybreaks
    \begin{align*}
        \frac{C(\beta)}{n^2}\sum_{i<j}\E\left[m_jP_{ij}\sum_k P_{ki}T_{iki}\right] &= \frac{C'(\beta)}{n^2}\sum_{i<j}\E\left[m_mm_jP_{ij}(P^2)_{ii}\right] \\
        &\le_{|m_im_jP_{ij}| \le 1} \frac{C'(\beta)}{n}\E\left[\sum_i (P^2)_{ii}\right] =_{\text{\eqref{eq:c-p2-bound}}} O(1)\,. 
    \end{align*}

    \ppart{$\mathrm{HO}_{11}$ term} The $\mathrm{HO}_{11}$ term involves a fourth cumulant, and requires bounding sub-terms that follow from replica rewrites. Fortunately, a bound on the absolute value of the fourth cumulant (akin to the one used to bound $|T_{ikj}| \le C$) combines with a $\ell_1-\ell_2$ inequality to yield the desired bound.
    \allowdisplaybreaks
    \begin{align*}
         \left|\sum_{i<j}P_{ij}\sum_k P_{ik}\Gamma_{ikij}\right| &\le_{|\Gamma_{ikij}| \le C} C\sum_{i<j}\sum_k|P_{ij}||P_{ik}| \le C\sum_{i}\left(\sum_j|P_{ij}|\sum_k |P_{ik}|\right) \\
         &= C\sum_i \left(\sum_j |P_{ij}|\right)^2 \le_{\text{CS}} Cn\sum_i \sum_j (P_{ij})^2 =Cn\Tr[P^2]\, , 
    \end{align*}
    at which point \eqref{eq:c-p2-bound} immediately implies that $\frac{4\beta^4}{n^2}\sum_{i<j}\E\left[P_{ij}\sum_k P_{ik}\Gamma_{ikij}\right]$.

    \ppart{$\mathrm{HO}_{12}$ term} Using the fact that $T_{jii}=-2m_iP_{ij}$ and some trace algebra, the bound follows rather straightforwardly as
    \allowdisplaybreaks
    \begin{align*}
        \sum_{i<j}T_{jii}\sum_k m_k(P^2)_{kj} &= -2\sum_{i<j}\sum_k m_im_k(P^2)_{kj}P_{ij}\,.
    \end{align*}
    Since the joint law of $(A,y_t)$ (with fixed plant $x_0 = \mathbf{1}_n$) is invariant under the permutation that swaps $i$ and $j$, a symmetrization argument immediately implies that
    \[
        \sum_{i<j}\E\left[\sum_k m_im_kP_{ij}(P^2)_{kj}\right] = \frac{1}{2}\sum_{i\ne j}\E\left[\sum_k m_i m_k  P_{ij}(P^2)_{kj}\right]\,.
    \]
    At this point, simple linear algebra yields that
    \allowdisplaybreaks
    \begin{align*}
        \frac{C(\beta)}{n^2}\sum_{i<j}\E\left[\sum_k m_im_kP_{ij}(P^2)_{kj}\right] &= \frac{C(\beta)}{n^2}\E\left[\mg^\sT P^3 \mg - \mg^\sT\diag(P)P^2\mg\right]\,.
    \end{align*}
    Using the fact that $\norm{\mg}^2_2 \le n$ and $\opnorm{P^3} \le \Tr[P^3]$, we immediately obtain that the first term is $O(1)$ by \eqref{eq:c-p3-bound}, and the second term is small because 
    \allowdisplaybreaks
    \begin{align*}
        |\mg^\sT P^2\diag(P) \mg| &\le_{\text{CS}}\norm{\diag(P)m}_2\norm{P^2\mg}_2 \le_{\diag(P)\preceq \Id}\sqrt{n}\norm{P^2\mg}_2 \\
        &\le \sqrt{n}\sqrt{\opnorm{P}^4\norm{\mg}^2_2} \le_{\opnorm{P^4}\le(\Tr[P^3])^{4/3}} n\Tr[P^3]^{2/3}\,. 
    \end{align*}
    At this point, invoking Jensen's (with the concavity of $x^{2/3}$ for $x\ge 0$) along with \eqref{eq:c-p3-bound} immediately yields that $\frac{C(\beta)}{n^2}\E[\mg^\sT\diag(P)P^2\mg] = O(n^{-1/3})$.

    \ppart{$\mathrm{HO}_{13}$ term} The bound for this term follows by straightforward replica algebra on using the identities that $(P^2)_{jk}=\frac{n}{4}\an{f_{1234}\Delta^{12}\Delta^{34}_k}$ and $T_{jik}=\an{\Delta^{56}_j\Delta^{57}_i\Delta^{58}_k}$ along with symmetrization over the $i<j$ sum to obtain
    \allowdisplaybreaks 
    \begin{align*}
        \sum_{i<j}m_i\sum_k (P^2)_{jk}T_{jik} = C(\beta)n^4\an{(R_{59}-R_{79})f_{1256}f_{3458}f_{1234}}-C'(\beta)n^2\an{B^n_i f_{1234}f_{3458}}\,,
    \end{align*}
    where $B^n_i$ is a function absolutely bounded by $Cn$ after summing over the $i$ indices. At this point, taking $\E\left[\cdot\right]$ and using $(4,4,4,4)$-H\"older's for the first term (after separating it into its 2 constituent bulk centered terms) with \pref{thm:overlap-moment-concentration} and \pref{lem:replicon-moments-from-mgf} shows the first term is $O(1)$, and using Cauchy--Schwarz along with $|B^n_i| \le Cn$ and the same concentration shows the second term is also $O(1)$.

    \ppart{$\mathrm{HO}_{14}$ term} Note that $\Gamma_{jiki}=-2m_iT_{ijk}-2P_{ji}P_{ik}$. This immediately implies that
    \allowdisplaybreaks
    \begin{align*}
        \sum_{i<j}\sum_k (P^2)_{jk}\Gamma_{jiki} &= -2\sum_{i<j}\sum_k (P^2)_{kj}m_iT_{ijk} -2\sum_{i<j}\sum_k (P^2)_{kj}P_{ji}P_{ik}\,.
    \end{align*}
    Note that $\frac{C(\beta)}{n^2}\E\left[\sum_{i<j}\sum_k(P^2)_{kj}m_iT_{ijk}\right] = O(1)$ since the term reduces (up to an index re-shuffle) the $\mathrm{HO}_{13}$ term, and the exact same argument works. For the second term, observe that exchangeability of $i$ and $j$ under $\E$ and $\an{\cdot}$ implies that
    \allowdisplaybreaks
    \begin{align*}
        \E\left[\sum_{i<j}\sum_k (P^2)_{kj}P_{ji}P_{ik}\right] &= \frac{1}{2}\E\left[\sum_{i,j,k}(P^2)_{kj}P_{ji}P_{ik}\right] - \frac{1}{2}\E\left[\sum_{i,k}(P^2)_{ik}P_{ki}P_{ii}\right] \\
        &= \frac{1}{2}\E\left[\Tr[P^4]\right] - \frac{1}{2}\E\left[\sum_i P^{3}_{ii}P_{ii}\right]\,.
    \end{align*}
    Using \eqref{eq:c-p4-bound} immediately yields that $\frac{C(\beta)}{n^2}\E\left[\sum_{i,j,k}(P^2)_{kj}P_{ji}P_{ik}\right] = O(1)$, and the second term can be bounded by $O(1)$ by using the fact that $0 \le P_{ii} \le 1$ in conjunction with \eqref{eq:c-p3-bound}.

    \ppart{$\mathrm{HO}_{15}$ term} This term is bounded via replica rewrites and standard invocations of \eqref{eq:c-p2-bound} and \eqref{eq:c-p4-bound}. Using the fact that $T_{kak}=-2m_kP_{ka}$, $T_{aii}=-2m_iP_{ai}$ and $T_{kai}=T_{aik}$, we obtain that
    \allowdisplaybreaks
    \begin{align*}
        \mathrm{HO}_{15} &= \frac{C(\beta)}{n^2}\E\left[-4\sum_{i,k}m^3_im_k(P^2)_{ki} + \sum_{i,k}m_i^3\sum_aP_{ka}T_{aik} + \sum_{i,k}m_i^2m_k\sum_aP_{ai}T_{aik}\right] 
    \end{align*}
    Bounding the first term is fairly straight-forward using the fact that $|m_i| \le 1$ and an application of Cauchy--Schwarz. Specifically,
    \allowdisplaybreaks
    \begin{align*}
        \sum_{i,k}m^3_im_k(P^2)_{ki} &= \langle \diag(m^3 \ot m), P^2 \rangle \le \norm{P^2}_F\norm{\diag(m^3\ot m)}_F \\
        &\le_{|m_i| \le 1} n\sqrt{\Tr[P^4]}\,.
    \end{align*}
    At this point, taking expectation with respect to $A$ and invoking \eqref{eq:c-p4-bound} in conjunction with Jensen's inequality yields that $\frac{C(\beta)}{n^2}\E\left[\sum_{i,k}m^3_im_k(P^2)_{ki}\right] = O(1)$. \\
    For the second and third term, we will use a lemma already used in the bound for the $\mathrm{HO}_ + \mathrm{HO}_7 + \mathrm{HO}_8$ terms. Namely,
    \[
        |\sum_a P_{ka}T_{aik}| \le 4\sqrt{(P^3)_{kk}}\,\qquad\qquad |\sum_a P_{ai}T_{aik}| \le 4\sqrt{(P^3)_{ii}}\,.
    \]
    This immediately yields that
    \allowdisplaybreaks
    \begin{align*}
        \sum_{i,k}m^3_i \sum_a P_{ka}T_{aik} \le_{|m_i| \le 1}4n \sum_k \sqrt{(P^3)_{kk}} \le_{\text{CS}} 4n^{3/2}\sqrt{\Tr[P^3]}\,.
    \end{align*}
    The above bound implies, via an invocation of Jensen's with \eqref{eq:c-p3-bound}, that $\frac{C(\beta)}{n^2}\E\left[\sum_{i,k}m_i^3\sum_aP_{ka}T_{aik}\right] = O(1)$. An immediately analogous argument with the bound $|m^2_im_k| \le 4$ implies that $\frac{C(\beta)}{n^2}\E\left[\sum_{i,k}m_i^2m_k\sum_aP_{ai}T_{aik}\right] = O(1)$.

    \ppart{$\mathrm{HO}_{16}$ term} The bound for the final term follows by decomposing the fourth cumulants with repeated indices into centered third moments, corrected by powers of $P$. These are bounded by replica rewrites and trace powers of $P$. More specifically, since
    \[
        \Gamma_{kaki} = -2P_{ka}P_{ki}-2m_kT_{aik}\,,\qquad\qquad \Gamma_{aiki} = -2P_{ai}P_{ik} -2m_iT_{aik}\,,
    \]
    the following simplifications ensue for the sums over the index $a$
    \allowdisplaybreaks
    \begin{align*}
        \sum_a \Gamma_{kaki}P_{ai} &= -2P_{ki}\sum_a P_{ka}P_{ai} - 2m_k\sum_a P_{ai}T_{aik} = -2P_{ki}(P^2)_{ki} -2m_k\sum_a P_{ai}T_{aik}\,, \\
        \sum_a \Gamma_{aiki}P_{ka} &= -2P_{ki}\sum_a P_{ka}P_{ai} - 2m_i\sum_a P_{ka}T_{aik }= -2P_{ki}(P^2)_{ki} -2m_i\sum_a P_{ka}T_{aik}\,.
    \end{align*}
    By the bounds for the prior terms, $|\sum_a P_{ka}T_{aik}| \le 4\sqrt{(P^3)_{kk}}$. This immediately implies that 
    \allowdisplaybreaks
    \begin{align*}
        \left|\sum_{i,k} \mg^3_i\sum_a P_{ka}T_{aik}\right| \le_{|\mg^3_i| \le 1,\,\text{CS}} 4n^{3/2}\sqrt{\Tr[P^3]}\, ,  
    \end{align*}
    and, by the same bounding strategy invoking $|\mg_i^2\mg_k| \le 1$,
    \allowdisplaybreaks
    \begin{align*}
        \left|\sum_{i,k}\mg_i^2\mg_k\sum_a P_{ai}T_{aik}\right| \le 4n^{3/2}\sqrt{\Tr[P^3]}\,.
    \end{align*}
    At this point, taking expectation $\E$ and invoking \eqref{eq:c-p3-bound} with Jensen's immediately gives the third centered terms with scaling $\frac{C(\beta)}{n^2}$ are $O(1)$. \\
    To bound the ``product=of-P'' terms, observe that
    \allowdisplaybreaks
    \begin{align*}
        \left|\sum_{i,k}\mg^2_i(P^2)_{ki}P_{ki}\right| \le_{|\mg^2_i| \le 1} \Tr[P^3]\, , 
    \end{align*}
    at which point invoking \eqref{eq:c-p3-bound} immediately gives that $\frac{C(\beta)}{n^2}\E\left[\sum_{i,k}\mg^2_i(P^2)_{ki}P_{ik}\right] = O(1)$. \qedhere 
\end{prf}

We now state the final lemma in this section, which shows that the term $\E[\Tr[AP^2V]]$ in \pref{lem:EA-term} is $O(1)$. The idea is simply to use Gaussian integration-by-parts one, followed by \pref{lem:dP-pairing} and some elementary bounds on the simplified terms that follow (quite similar to the process in \pref{lem:bounds-for-ho-terms}).

\begin{lemma}[The rank-$1$ term in $\calE_A$ is $O(1)$]\label{lem:ea-rank-1-term}
    \[
        2\beta\E\left[\Tr[AP^2V]\right] = O(1)\, ,
    \]
    where $V = \frac{2\beta^2}{n}\mg\mg^\sT$.
\end{lemma}
\begin{prf}
    A lengthy (but straightforward) algebraic calculation using \eqref{eq:GOE-ibp} and \pref{lem:dP-pairing} with exact identities for ``repeated index'' centered third-moments an fourth cumulants leads to the following decomposition of the original term
    \allowdisplaybreaks
    \begin{align*}
        \beta\E\left[\Tr[AP^2V]\right] &= \frac{4\beta^3}{n}\E\left[\mg^\sT A (P^2\mg)\right] = \frac{4\beta^3}{n^2}\sum_{i<j}\E\left[\partial_{ij}\left(m_i(P^2\mg)_j + m_j(P^2\mg)_i\right)\right] \\
        &= \frac{4\beta^3}{n^2}\E\left[4\mg^\sT P^3 \mg + 2\Tr[P]m^\sT P^2\mg + \Tr[P^2]\mg^\sT P \mg + \norm{\mg}^2_2\Tr[P^3]\right] \\
        &- \frac{4\beta^3}{n^2}\E\left[4\sum_i \mg_iP_{ii}(P^2\mg)_i - 2\sum_i \mg_i(P\mg)_i(P^2)_{ii} - 2\sum_i \mg^2_i(P^3)_{ii}\right] \\
        &+ \frac{4\beta^3}{n^2}\E\left[\sum_{i<j}\left(T_{iij}(P^2\mg)_j + T_{jij}(P^2\mg)_i\right)\right] \\
        &+\frac{4\beta^3}{n^2}\E\left[\sum_{i<j} \left(m_i\sum_a\left(m_iT_{jaj} + m_jT_{jai} + \Gamma_{jaij}\right)(P\mg)_a\right)\right] \\
        &+ \frac{4\beta^3}{n^2}\E\left[\sum_{i<j} \left(m_j\sum_a\left(m_iT_{iaj} + m_jT_{iai} + \Gamma_{iaij}\right)(P\mg)_a\right)\right] \\
        &+ \frac{4\beta^3}{n^2}\E\left[\sum_{i<j} \left(m_i\sum_{a,k}P_{ja}\mg_k\left(m_iT_{akj} + m_jT_{aki} + \Gamma_{akij}\right)\right)\right] \\
        &+ \frac{4\beta^3}{n^2}\E\left[\sum_{i<j} \left(m_j\sum_{a,k}P_{ia}\mg_k\left(m_iT_{akj} + m_jT_{aki} + \Gamma_{akij}\right)\right)\right] \\
        &+ \frac{4\beta^3}{n^2}\E\left[\sum_{i<j} \left(m_i\sum_k(P^2)_{jk}T_{kij} + m_j\sum_k(P^2)_{ik}T_{kij}\right)\right]\,.
    \end{align*}
    We first bound the tracial terms. Note that
    \allowdisplaybreaks
    \begin{align*}
        \frac{C(\beta)}{n^2}\E\left[\mg^\sT P^3\mg\right] &\le_{\norm{\mg}^2_2 \le n} \frac{C(\beta)}{n}\E\left[\Tr[P^3]\right] =_{\text{\eqref{eq:c-p3-bound}}} O(1)\,, \\
        \frac{C(\beta)}{n^2}\E\left[\Tr[P]\mg^\sT P^2 \mg\right] &\le_{\Tr[P] \le n} \frac{C(\beta)}{n}\E\left[\mg^\sT P^2 \mg\right] =_{\text{Proof of~\pref{lem:bound-d-remainder}}} O(1) \, ,\\
        \frac{C(\beta)}{n^2}\E\left[\Tr[P^2]\mg^\sT P \mg\right] &\le_{\text{CS}}\sqrt{\E[\Tr[P^2]^2]}\sqrt{\E[(\mg^\sT P\mg)^2]}\le_{\text{$\mathrm{R}_{t2}$}}\frac{C(\beta)}{n}\sqrt{\mathsf{span}\left\{\nu\left(ff'\hat{R}_\alpha\hat{R}_\beta\hat{R}_\gamma\hat{R}_\delta\right)\right\}} = O(1)\,,\\
        \frac{C(\beta)}{n^2}\E\left[\norm{\mg}^2_2\Tr[P^3]\right] &\le_{\norm{\mg}^2_2 \le n} \frac{C(\beta)}{n}\E\left[\Tr[P^3]\right] =_{\text{\eqref{eq:c-p3-bound}}} O(1)\,,\\
        \frac{C(\beta)}{n^2}\E\left[|\sum_i m_iP_{ii}(P^2\mg)_i|\right] &= \frac{C(\beta)}{n^2}\E\left[(E_{\calD_n}[P]\mg)^\sT P^2 \mg\right] \le_{\text{CS},\,\norm{\diag(P)\mg}_2 \le \sqrt{n}} \le \frac{C(\beta)}{n^{3/2}} \E\left[\norm{P^2\mg}_2\right] \\
        &\le\frac{C(\beta)}{n}\E\left[\Tr[P^2]\right] =_{\text{\eqref{eq:c-p2-bound}}} O(1)\,, \\
        \frac{C(\beta)}{n^2}\E\left[\sum_i \mg_i(P\mg)_i(P^2)_{ii}\right] &\le_{|(P^2)_{ii}| \le Cn, |m_i| \le 1}\frac{C(\beta)}{n}\E\left[\sum_i (P^2)_{ii}\right] \le_{\text{\eqref{eq:c-p2-bound}}} O(1)\,, \\
        \frac{C(\beta)}{n^2}\E\left[\sum_i \mg^2_i (P^3)_{ii}\right] &\le_{|\mg^2_i| \le 1} \frac{C(\beta)}{n^2}\E\left[\sum_i (P^3)_{ii}\right] =_{\text{\eqref{eq:c-p3-bound}}} O(1)\,. 
     \end{align*}
     For the higher-order terms, we bound them line-by-line. For the first line, note that $T_{iij}=-2m_iP_{ij}$ and $T_{jij}=2-m_jP_{ij}$, which gives
     \allowdisplaybreaks
     \begin{align*}
         \frac{4\beta^3}{n^2}\E\left[\sum_{i<j}\left(T_{iij}(P^2\mg)_j + T_{jij}(P^2\mg)_i\right)\right] &= \frac{C(\beta)}{n^2}\E\left[\sum_{i<j} \mg^\sT P^3 \mg - \sum_i \mg_i P_{ii}(P^2\mg)_i\right] = O(1)\, ,
     \end{align*}
     which follows by the bounds developed above for the tracial terms. \\
     The second and third line have terms with the exact same structure (up to an index flip) and so we bound just one, since an exactly analogous argument applies to the other. We again have that $T_{jaj}=-2m_iP_{ia}$ and $\Gamma_{iaij} = m_iT_{iaj}-2P_{ia}P_{ij}$, which gives
     \[
        m_iT_{iaj} + m_jT_{iai} +\Gamma_{iaij} = 2\left(m_iT_{iaj}-2m_im_jP_{ia} - 2P_{ia}P_{ij}\right)\, , 
     \]
     and so
     \allowdisplaybreaks
     \begin{align*}
        \frac{4\beta^3}{n^2}\E\left[\sum_{i<j} \left(m_j\sum_a\left(m_iT_{iaj} + m_jT_{iai} + \Gamma_{iaij}\right)(P\mg)_a\right)\right]  &= \frac{C(\beta)}{n^2}\E\left[\sum_{i<j}\mg_i\mg_j\sum_a T_{iaj}(P\mg)_a\right] \\
        &- \frac{C(\beta)}{n^2}\E\left[\sum_{i<j}\mg_i\mg_j^2\sum_a P_{ia}(P\mg)_a\right] \\
        &-\frac{C(\beta)}{n^2}\E\left[\sum_{i<j}m_jP_{ij}\sum_a P_{ia}(P\mg)_a\right]\,,
     \end{align*}
     where the last two terms further simplify (after symmetrization) to
     \[
        \frac{C(\beta)}{n^2}\E\left[\sum_{i<j}\mg_jP_{ij}\sum_a P_{ia}(P\mg)_a\right] = \frac{C(\beta)}{n^2}\E\left[\mg^\sT P^3 \mg - \sum_i \mg_iP_{ii}(P^2\mg)_i\right]\,,
     \]
     and
     \[
        \frac{C(\beta)}{n^2}\E\left[\sum_{i<j}\mg_i\mg_j^2\sum_a P_{ia}(P\mg)_a\right] = -\frac{C(\beta)}{n^2}\E\left[\norm{\mg}^2_2\mg^\sT P^2 \mg - \mg^\sT P^2 \mg + \sum_i \mg_i P_{ii}(P^2\mg)_i\right]\, .
     \]
     Each of the quantities in the RHS of both terms are controlled to be $O(1)$ by the prior bounds in the proof on the tracial terms. It remains to bound the first term in the simplification, and this follows by using the standard technique of doing replica-rewrites and reducing to a span of terms which are controlled by $(4,4,4,4)$-H\"older's and \pref{prop:D15}. In this case,
     \[
        \sum_a T_{iaj}(P\mg)_a = \frac{n^2}{2}\an{\Delta^{12}_i\Delta^{14}_j f_{1356}(R_{57}-R_{67})}\, ,
     \]
     which immediately gives
     \[
        \frac{C(\beta)}{n^2}\E\left[\sum_{i<j}\mg_i\mg_j^2\sum_a P_{ia}(P\mg)_a\right] = C(\beta)\E\left[\an{f_{1356}(R_{57}-R_{67})\left(\sum_{i<j}\mg_i\mg_j\Delta^{12}_i\Delta^{14}_j\right)}\right]\,.
     \]
     Symmetrizing the sum over $i$ and $j$ yields that
     \allowdisplaybreaks
     \begin{align*}
          \frac{C(\beta)}{n^2}\E\left[\sum_{i<j}\mg_i\mg_j^2\sum_a P_{ia}(P\mg)_a\right] = \mathsf{span}\left\{C(\beta)n^2\nu\left(f_{1356}\Xi_{ab}\Xi_{cd}\Xi_{ce}\right), C(\beta)n\E\left|\an{f_{1356}\Xi_{ab}}\right|\right\}\, ,
     \end{align*}
     which is $O(1)$ by the appropriate invocations of H\"older and \pref{prop:D15}. \\
     We now move to bound the terms on the fourth and fifth line. Again, because of the symmetry in their structure, we bound only one as an exactly analogous argument works for the other. We split the main term into three components as
     \allowdisplaybreaks
     \begin{align*}
         \frac{4\beta^3}{n^2}\E\left[\sum_{i<j} \left(m_j\sum_a\left(m_iT_{iaj} + m_jT_{iai} + \Gamma_{iaij}\right)(P\mg)_a\right)\right] &= \frac{C(\beta)}{n^2}\E\left[\sum_{i<j}\sum_{a,k}\mg^2_i\mg_kP_{ja}T_{akj}\right] \\
         &- \frac{C(\beta)}{n^2}\E\left[\sum_{i<j}\sum_{a, k}\mg_i\mg_j\mg_kP_{ja}T_{aki}\right] \\
         &- \frac{C(\beta)}{n^2}\E\left[\sum_{i<j}\sum_{a,k} \mg_i\mg_k P_{ja}\Gamma_{akij}\right]\,.
     \end{align*}
     The first term can be written (using the same replica rewrites as before) as
     \[
        \frac{C(\beta)}{n^2}\E\left[\sum_{i<j}\sum_{a,k}\mg^2_i\mg_kP_{ja}T_{akj}\right] = C(\beta)\E\left[\sum_{i<j}\mg^2_i\an{\Delta^{56}_j\Delta^{14}_jf_{1256}(R_17)-R_37}\right]\,.
     \]
     Symmetrization over the indices $i$ and $j$ yields that
     \allowdisplaybreaks
     \begin{align*}
         \frac{C(\beta)}{n^2}\E\left[\sum_{i<j}\sum_{a,k}\mg^2_i\mg_kP_{ja}T_{akj}\right] = \mathsf{span}\left\{C(\beta)n^2\nu\left(f_{1456}f_{1256}\Xi_{ab}\right),\,C(\beta)n\nu\left(\left|f_{1256}\Xi_{ab}\right|\right)\right\}\, ,
     \end{align*}
     which is $O(1)$ by invocations of \pref{prop:D15} and a Cauchy--Schwarz followed by \pref{thm:overlap-moment-concentration} and \pref{lem:replicon-moments-from-mgf}. For the second term, we again symmetrize the sum over $i$ and $j$ to obtain that
     \[
        - \frac{C(\beta)}{n^2}\E\left[\sum_{i<j}\sum_{a, k}\mg_i\mg_j\mg_kP_{ja}T_{aki}\right] = \frac{C(\beta)}{n^2}\E\left[\sum_{a,k}(P\mg)_a\sum_i \mg_i\mg_kT_{aki}\right] - \frac{C(\beta)}{n^2}\E\left[\sum_{i,a,k}\mg_i^2\mg_k P_{ia}T_{aki}\right]\,.
     \]
     Note that the second term in the RHS has the same form (up to an index change) with the first term, and so the same replica-rewrite followed by \pref{prop:D15} and Cauchy--Schwarz yields it is $O(1)$. For the first term, we again take another replica--rewrite and group the sums to obtain 
     \[
        \frac{C(\beta)}{n^2}\E\left[\sum_{a,k}(P\mg)_a\sum_i \mg_i\mg_kT_{aki}\right] = \mathsf{span}\left\{C(\beta)n\nu\left(f_{1256}(R_{59}-R_{69})(R_{17}-R_{37})(R_{18}-R_{48})\right)\right\}\,,
     \]
     which is $O(1/n)$ by an immediate application of \pref{prop:D15}. The fourth cumulant term is easily shown to be rewritten using replica rewrites after the centering definition is applied as
     \allowdisplaybreaks
     \begin{align*}
         \frac{C(\beta)}{n^2}\E\left[\sum_{i<j}\sum_{a,k}\mg^2_i\mg_kP_{ja}T_{akj}\right] &= C(\beta)n^2\nu\left(f_{1567}f_{1267}(R_{18}-R_{48})(R_{19}-R_{39})\right) - \frac{C(\beta)}{n^2}\E\left[\sum_{i<j}\mg_iP_{ij}\sum_{k}\mg_k(P^2)_{kj}\right] \\
         &- \frac{C(\beta)}{n^2}\E\left[\sum_{i<j}\mg_i(P^2)_{ij}\sum_k\mg_kP_{kj}\right] - \frac{C(\beta)}{n^2}\E\left[\sum_{i<j}\mg_i(P^2)_{jj}\sum_k\mg_kP_{ik}\right]\,.
     \end{align*}
    The first term in the RHS is clearly $O(1)$ by an invocation of $(4,4,4,4)$-H\"older's in conjunction with \pref{thm:overlap-moment-concentration} and \pref{lem:replicon-moments-from-mgf}. For the remaining terms we use permutation-invariance under $\E$ to use $\E[\sum_{i<j}X_{ij}] = \frac{1}{2}\E[\sum_{i \ne j}X_{ij}]$ as in the proof of \pref{lem:bounds-for-ho-terms}. So, for the second term in the RHS, observe after symmetrization
    \[
        \E\left[\sum_{i<j}\mg_iP_{ij}\sum_k \mg_k(P^2)_{kj}\right] = \frac{1}{2}\E\left[\mg^\sT P^3 \mg - \sum_i \mg_i P_{ii}(P^2\mg)_i\right]\, , 
    \]
    where the RHS we have bounded before in the tracial terms, implying that $\frac{C(\beta)}{n^2}\E\left[\sum_{i<j}\mg_iP_{ij}\sum_{k}\mg_k(P^2)_{kj}\right] = O(1)$. Similarly, symmetrizing and doing some algebra gives the following for the third and fourth terms in the RHS
    \allowdisplaybreaks
    \begin{align*}
        \frac{C(\beta)}{n^2}\E\left[\sum_{i<j}\mg_i(P^2)_{ij}\sum_k\mg_kP_{kj}\right] &= \frac{C(\beta)}{n^2}\E\left[\mg^\sT P^3 \mg - \sum_i \mg_i (P^2)_{ii}(P\mg)_i\right]\,, \\
        \frac{C(\beta)}{n^2}\E\left[\sum_{i<j}\mg_i(P^2)_{jj}\sum_k\mg_kP_{ik}\right] &= \frac{C(\beta)}{n^2}\E\left[\Tr[P^2]\mg^\sT P \mg - \sum_i\mg_i(P\mg)_i(P^2)_{ii}\right]\,,
    \end{align*}
    and, once again, each of the terms in the RHS has been bounded before in the tracial terms, implying the two terms in the LHS are $O(1)$. \\
    Lastly, the higher-order terms in the last line have \emph{exactly} the same form as those in $\mathrm{HO}_1$. Therefore, both terms on the last line are $O(1)$ by \pref{lem:bounds-for-ho-terms}.
\end{prf}

\section{Resolvent analysis: Regularity and control of the diagonal}\label{app:deformed-wigner-resolvent}
This appendix supplies the random-matrix input used in \pref{sec:alg-properties} and develops the accompanying regularity estimates for deformed-Wigner resolvents. Its main output is a uniform comparison between the diagonal of the finite-\(n\) real-axis squared resolvent \(E_{\mathcal D_n}[D(a\Id_n-\beta A + D)^{-2}D]\) and its freely independent limit. We first establish deterministic operator-norm (\pref{sec:resolvent-op-norm-bounds}) and derivative (\pref{sec:resolvent-lipschitz-bounds}) bounds for resolvent products under a spectral-gap assumption. We then combine a chaining argument (\pref{sec:chaining-diagonal}), analytic subordination (\pref{sec:free-prob}), and Gaussian interpolation (\pref{sec:free-interpolation}) to upgrade this comparison to a uniform high-probability bound over the relevant diagonal deformations. Relative to the analogous deformed-Wigner analysis in \cite{jekel2024pha,jekel2025pha2}, the main new difficulty is that the observable lies on the real axis and can encounter singularities, so the argument requires a cutoff together with sharper edge and derivative estimates for squared resolvents.
In particular, a high-probability bound on the edges of the spectra of matrices of the form $A \ot \Id_k + B_k$ where $A$ and $B_k$ are both GOE is provided in \pref{sec:AIk+Bk}. Another (comparatively) minor technical issue stems from the fact that a naive $\eps$-net argument does not give the required rate of convergence, therefore, the use of Dudley entropy bounds is necessitated to shave-off a wasteful $\log(\cdot)$ factor. \\
The culmination is \pref{cor:diagonal-Q2-controlled-in-l2}, which gives the diagonal control with quantitative rate needed for the ASL--TAP versus PHD comparison.

\subsection{Operator norm bounds on resolvent products}
\label{sec:resolvent-op-norm-bounds}

The next subsections will be  analogues of various lemmata found in \iffocs{\cite{jekel2024pha}}{\cite{jekel2025pha2}}, but for the setting where $(a - \beta A + D) \succeq \gamma > 0$ and we can take the imaginary part $\to$ 0.

\begin{lemma}[Bounds on resolvent products]
\label{lem:resolvent-bounds}
Let $W, T$ be real symmetric matrices and $M$ be a Hermitian matrix of the same dimensions.
    Let
\[
    D \;:=\; W^{-1},
    \qquad
    X \;:=\; M + D,
    \qquad
    R \;:=\; X^{-1}.
\]
    Suppose that there is a $\gamma > 0$ so that $\min \Spec(|X|) \ge \gamma$ and $\|M\| \le C_1$.
    Then when $k \ge 1$ and $j \ge 2$,
\begin{align*}
    &\opnorm{R^k} \;\le\; \gamma^{-k},
    \qquad\qquad
    &&\opnorm{DR^k} \;\le\; (\gamma + C_1)\gamma^{-k},
\\
    &\opnorm{DR^jD} \;\le\; (\gamma + C_1)^2\gamma^{-j},
    \qquad\quad
    &&\schnorm{R^kDTDR^{j-k}} \le\;\; (\gamma + C_1)^2\gamma^{-j}\schnorm{T}.
\end{align*}
If $\Im[M] = b\Id \ne 0$, we also have the imaginary-part bounds
\begin{align*}
    &\opnorm{\frac{1}{b}\Im[R]} \;\le\; \gamma^{-2},
    \qquad\qquad
    \schnorm{\frac{1}{b}\Im[RDTDR]} \le\;\; 2(\gamma + C_1)^2\gamma^{-3}\schnorm{T},
\\
    &\schnorm{\frac{1}{b}\Im[R^2DTDR]} \le\;\; 3(\gamma + C_1)^2\gamma^{-4}\schnorm{T} + O(b^2).
\end{align*}
\end{lemma}
\begin{prf}
Since $X$ is Hermitian,
\[\opnorm{R^k} = \opnorm{X^{-k}} = \max\Spec(|X|^{-k}) = (\min\Spec(|X|))^{-k} \le \gamma^{-k}.\] 
Next, use $D = X - M$, triangle inequality, and $\norm{R^k} \le \gamma^{-k}$ to write
\[
    \opnorm{DR^k} \;=\; \opnorm{\bigl(X-M\bigr)R^k} \;=\; \opnorm{R^{k-1} - MR^k} \;\le\; \gamma^{-(k-1)} + C_1\gamma^{-k}.
\]
Similarly, when $j \ge 2$, and invoking the preceding bound on $\opnorm{DR^{j-1}}$ at the end,
\begin{align*}
    \opnorm{DR^jD} &=\;\; \opnorm{\bigl(X-M\bigr)R^jD} =\;\; \opnorm{R^{j-1}D - MR^jD}
    \\&\le\;\; \opnorm{DR^{j-1}} + \opnorm{M}\opnorm{DR^j} \le\;\; \opnorm{DR^{j-1}}(1 + C_1\gamma^{-1}).
\end{align*}
For the imaginary-part bounds, we have $X^*R^* = \Id$ and $R^*R$ is real, so
\[
\opnorm{\frac{1}{b}\Im[R]} =\; \opnorm{\frac{1}{b}\Im[X^*R^*R]} \le\; \opnorm{R^*R},
    \qquad\;
    \schnorm{\frac{1}{b}\Im[RDTDR]} =\; 2\schnorm{\frac{1}{b}\Im[R]DTD\,\Re[R]},
\]
and similarly for $\schnorm{\frac{1}{b}\Im[R^2DTDR]}$.
\end{prf}

\subsection{Derivative bounds and Lipschitz properties of real--axis resolvent squares}
\label{sec:resolvent-lipschitz-bounds}

 To show that various products involving resolvents are Lipschitz, we will bound their Fr\'echet derivatives.

 \begin{definition}
 \label{def:frechet}
     Let $\phi: U \to V$ be a function between two finite-dimensional normed vector spaces $U$ and $V$ with norms $\norm{\cdot}_U$ and $\norm{\cdot}_V$ respectively.
     If it exists, the Fr\'echet derivative of $\phi$ at a point $u \in U$ is the linear operator $d\phi(u): U \to V$ such that
     \[ \lim_{\norm{h}_U\to 0} \frac{\norm{\phi(u + h) - \phi(u) - d\phi(u)[h]}_V}{\norm{h}_U} = 0. \]

     We write $\frenorm{\cdot}$ for the operator norm between the two tangent spaces of a Fr\'echet derivative.
     \[ \frenorm{d\phi(u)} := \sup_{\norm{h}_U = 1}\norm{d\phi(u)[h]}_V.\]
    When $V = M_n(\R)$, then unless otherwise specified, we will take the norm on $V$ to be the normalized Schatten-2 norm $\norm{\cdot}_V \gets \schnorm{\cdot}$.
    We may write $\frenorm[\to \infty]{d\phi(u)}$ to indicate that the norm on $V$ is taken to be the matrix operator norm.

When applying the Poincar\'e inequality, we may consider the Hilbert-Schmidt norm
\[
\norm{d\phi(u)}_{\mathrm{HS}(U,V)}^2
:=
\sum_{\alpha}\|d\phi(u)[e_\alpha]\|_{V}^2,\]
where $\{e_{\alpha}\}$ is any orthonormal basis of $U$.
The point $u$ where we evaluate the Fr\'echet derivative $d\phi(u)$ will often be implicit---for example, $d\phi$ or $\norm{d\phi}$. 
 \end{definition}

 Fr\'echet derivatives of matrix-valued functions satisfy the non-commutative sum, product, power, inverse, and chain rules.
 Explicitly,
 \begin{itemize}
     \item $d(f+g) = df + dg$,
     \item $d(fg) = (df)g + f(dg)$,
     \item $d(f^k) = \sum_{i = 1}^{k}f^{i-1}(df)f^{k-i}$,
     \item $d(f^{-1}) = -f^{-1}(df)f^{-1}$,
     \item $d(f\circ g)(u)[h] = df(g(u))[dg(u)[h]]$.
 \end{itemize}
 A full treatment can be found in~\cite[Section A.5]{absil2008optimization}.

 \begin{fact}
    \label{fact:frechet-lipschitz}
     If $\phi$ is Fr\'echet differentiable in some open set $S$ and $\frenorm{d\phi(u)} \le L$ for all $u\in S$, then $\phi$ is $L$-Lipschitz on $S$.
 \end{fact}

\begin{lemma}[Bounds on derivatives of resolvent products]
    \label{lem:resolvent-frechet}
    Let $W$ be a real symmetric matrix and $M$ be a Hermitian matrix of the same dimensions that are functions of parameters that lie in an open subset of some vector space.
    Let $T$ be a constant matrix of the same dimensions.
    Let $D := W^{-1}$ and $X := M + D$ and $R := X^{-1}$.
    Suppose that there is a $\gamma > 0$ so that for all parameters of $W$ and $M$, $\min \Spec(|X|) \ge \gamma$ and $\opnorm{M} \le C_1$.
    Then, letting $C_2 := (\gamma + C_1)\gamma^{-1}$,
    \begin{enumerate}[label=(\alph*), ref=\thetheorem(\alph*)]
        \item \[\frenorm{dR} \le \gamma^{-2}\frenorm{dM} + C_2^2\frenorm{dW}.\]
        \label{lem:resolvent-frechet-dR}
        \item \[\frenorm{d(R^k)} \le k\gamma^{-(k+1)}\frenorm{dM} + kC_2^2\gamma^{-(k-1)}\frenorm{dW}.\]
        \label{lem:resolvent-frechet-dRk}
        \item \[\frenorm{d(DR)} \le C_2\gamma^{-1}\frenorm{dM} + C_1C_2^2\frenorm{dW}.\]
        \label{lem:resolvent-frechet-dDR}
        \item For $k \ge 1$,
        \[\frenorm{d(DR^k)} \le kC_2\gamma^{-k}\frenorm{dM} + (kC_2-1)C_2^2\gamma^{-(k-2)}\frenorm{dW}.\]
        \label{lem:resolvent-frechet-dDRk}
        \item For $k \ge 2$, \[\frenorm{d(DR^kD)} \le kC_2^2\gamma^{-(k-1)}\frenorm{dM} + (kC_2-2)C_2^3\gamma^{-(k-3)}\frenorm{dW}.\]
        \label{lem:resolvent-frechet-dDRkD}
        \item For $k \ge 2$ and $1 \le j \le k-1$,
        \begin{align*}
            \frenorm{d(R^jDTDR^{k-j})} \le  kC_2^2\gamma^{-(k-1)}\frenorm{dM}\opnorm{T} + (kC_2-2)C_2^3\gamma^{-(k-3)}\frenorm{dW}\opnorm{T}.
        \end{align*}
        \label{lem:resolvent-frechet-dRjDTDRkj}
        \item For $k \ge 2$ and $1 \le j \le k-1$, and letting $T$ be variable for this bound only,
        \begin{align*}
            \frenorm{d(R^jDTDR^{k-j})} \le{} &  kC_2^2\gamma^{-(k-1)}\frenorm[\to\infty]{dM}\schnorm{T} 
            \\&+ (kC_2-2)C_2^3\gamma^{-(k-3)}\frenorm[\to\infty]{dW}\schnorm{T} + C_2^2\gamma^{-(k-2)}\frenorm{dT}.
        \end{align*}
        \label{lem:resolvent-frechet-dRjDTDRkj-T}
        \item When $\Im[X]$ is invertible, \[\frenorm{d((\Im [X])^{-1}\Im[R])} \le 2\gamma^{-3}\frenorm{dM} + 2C_2^2\gamma^{-1}\frenorm{dW}.\]
        \label{lem:resolvent-frechet-dImR}
        \item When $\Im[X] = b\Id \ne 0$,
        \[\frenorm{d(b^{-1}\Im[RDTDR])} \le 6C_2^2\gamma^{-2}\frenorm{dM}\opnorm{T} + 2(3C_2-2)C_2^3\frenorm{dW}\opnorm{T}.\]
        \label{lem:resolvent-frechet-dImRDTDR}
    \end{enumerate}
\end{lemma}
\begin{proof}
    By \pref{lem:resolvent-bounds}, we have $\opnorm{R^k} \le \gamma^{-k}$ and $\opnorm{DR^{k+1}} \le C_2\gamma^{-k}$.
    Note also $dD = -DdWD$ and $dX = dM - D(dW)D$.

    \ppart{Part \ref{lem:resolvent-frechet-dR}} We compute the derivative
    \[dR = -R(dX)R = -R(dM)R + (RD)(dW)(DR).\]

    \ppart{Part \ref{lem:resolvent-frechet-dRk}}
    We use 
    \[d(R^k) = \sum_{j=0}^{k-1} R^j(dR)R^{k-j-1},\]
    then apply \pref{lem:resolvent-frechet-dR}.

    \ppart{Part \ref{lem:resolvent-frechet-dDR}}
    We use $DR = \Id - MR$ and hence $d(DR) = -(dM)R - M(dR)$. then apply \pref{lem:resolvent-frechet-dR} and $C_2 = 1 + C_1\gamma^{-1}$.

    \ppart{Part \ref{lem:resolvent-frechet-dDRk}}
    Compute the derivative 
    \[d(DR^k) = d(DR)R^{k-1} + (DR)d(R^{k-1}).\]
    Then apply \pref{lem:resolvent-frechet-dDR} and \pref{lem:resolvent-frechet-dRk} and, for the $\frenorm{dW}$ term, $C_1C_2^2\gamma^{-(k-1)} + (k-1)C_2^3\gamma^{-(k-2)} = (kC_2-1)C_2^2\gamma^{-(k-2)}$.

    \ppart{Part \ref{lem:resolvent-frechet-dDRkD}} 
    Compute the derivative
    \[d(DR^kD) = d(DR)(R^{k-1}D) + (DR)d(R^{k-1}D),\]
    then apply \pref{lem:resolvent-frechet-dDR} and \pref{lem:resolvent-frechet-dDRk} and collect terms.
    
    \ppart{Part \ref{lem:resolvent-frechet-dRjDTDRkj}} 
    Compute the derivative
    \[d(R^jDTDR^{k-j}) = d(R^jD)T(DR^{k-j}) + (R^jD)Td(DR^{k-j}).\]
    Then apply \pref{lem:resolvent-frechet-dDRk} and collect terms.

    \ppart{Part \ref{lem:resolvent-frechet-dRjDTDRkj-T}} 
    The same as \pref{lem:resolvent-frechet-dRjDTDRkj} but with the addition of another $(R^jD)(dT)(DR^{k-j})$ term, and also taking H\"older's inequality in a different way.

    \ppart{Part \ref{lem:resolvent-frechet-dImR}}
    $\Im[R] = \Im[X^*R^*R] = -\Im[X]R^*R$ since $X^*R^* = \Id$ and $R^*R$ is real, so 
    \[d((\Im[X])^{-1}\Im[R]) = d(R^*R) = d(R^*)R + R^*dR.\]
    Then apply \pref{lem:resolvent-frechet-dR} and $\opnorm{R} = \opnorm{R^*} \le \gamma^{-1}$.

    \ppart{Part \ref{lem:resolvent-frechet-dImRDTDR}}
    Compute the imaginary part
    \[b^{-1}\Im[RDTDR] = b^{-1}\Re[R]\,DTD\,\Im[R] + b^{-1}\Im[R]\,DTD\,\Re[R] = \Re[RD]\,TDRR^* + R^*RDT\,\Re[DR].\]
    Then perform a calculation similar to \pref{lem:resolvent-frechet-dRjDTDRkj} using the contractivity of $\Re[\cdot]$.
\end{proof}

\subsection{Chaining argument for diagonal of resolvent}
\label{sec:chaining-diagonal}
We perform a chaining argument to bound the centered fluctuation of the diagonal squared resolvent uniformly over the real parameter \(a\) and the diagonal deformation \(D\). After truncating \(A\) by an operator-norm projection (denoted \(\overline A\) below), the quantity of interest is
\[
\sup_{a\in[0,1]}\sup_{W\in\mathcal D_n((0,1])}
\schnorm{
E_{\mathcal D_n}\!\Big[D(a\Id_n-\beta\overline A + D)^{-2}D\Big]
-
\E_A\!\left[
E_{\mathcal D_n}\!\Big[D(a\Id_n-\beta\overline A + D)^{-2}D\Big]
\right]
},
\qquad D:=W^{-1}.
\]
By diagonal norm duality, this is equivalent to controlling the corresponding centered scalar process against diagonal test matrices \(T\) with \(\schnorm{T}=1\).
Chaining is used here to remove the logarithmic loss that would arise from a naive \(\varepsilon\)-net over this $(a,W,T)$ parameter family, obtaining a bound of scale $O(n^{-1/2})$.

First we introduce Talagrand's $\gamma_2(T,d)$, as given in \cite[Definition 2.7.3]{talagrand2014upper}.
It is a multiscale complexity parameter of the metric space $(T,d)$; informally, it measures how well $T$ can be successively approximated by increasingly fine discrete nets. It is the natural quantity controlling suprema of sub-Gaussian processes.

\begin{theorem}[Upper bound for the $\gamma_2$-functional,~{\cite[Exercise 2.7.6(b), Eq.~(2.40), and (2.37)]{talagrand2014upper}}]\label{thm:talagrand-gamma-2}
    Let $(T,d)$ be a separable metric space. Then, for some absolute constant $C > 0$,
    \[
        \gamma_2(T,d) \le C\int_{0}^{\mathsf{diam}(T,d)}\sqrt{\log N(T,d,\epsilon)}d\eps\, ,
    \]
    where $N(T,d,\eps)$ is the minimal $\eps$-covering number.
\end{theorem}

\begin{theorem}[Chaining bound for sub-Gaussian processes,~{\cite[Theorem 2.7.11]{talagrand2014upper}}]\label{thm:talagrand-subgaussian-chaining}
    Consider a separable metric space $(T,d)$ and stochastic process $\{X_t\}_{t\in T}$ with centering $\E[X_t] = 0$ for every $t \in T$. Assume the process satisfies the increment condition
    \[
        \forall u > 0,\, \P\left\{|X_s - X_t| > u \right\} \le 2 \exp\left(\frac{-u^2}{2d(s,t)^2}\right)\, ,
    \]
    Then
    \[
        \E\left[\sup_{t \in T}X_t\right] \le C\,\gamma_2(T,d)
    \]
    for some absolute constant $C$.
\end{theorem}

\begin{lemma}[Herbst concentration for normalized GOE]
\label{lem:herbst-goe}
Let $A\in M_n(\mathbb R)$ be a normalized GOE matrix, i.e.
$(A_{ij})_{1\le i<j\le n}$ are independent $\mathcal N(0,1/n)$,
$A_{ji}=A_{ij}$, and $(A_{ii})_{1\le i\le n}$ are independent $\mathcal N(0,2/n)$.
If $f:\mathrm{Sym}_n(\mathbb R)\to\mathbb R$ is $L$-Lipschitz with respect to $\schnorm{\cdot}$, then for all $\delta\ge 0$,
\[
\Pr\Big(|f(A)-\mathbb Ef(A)|\ge \delta\Big)
\;\le\;
2\exp\!\left(-\frac{n^2\,\delta^2}{4L^2}\right).
\]
\end{lemma}

\begin{proof}
Write $A=T(g)$ as a linear image of a standard Gaussian vector $g\in\mathbb R^{m}$
with $m=n(n+1)/2$ coordinates, so that
$T(g)_{i,j}=\frac{1}{\sqrt n}\,g_{\min(i,j),\max(j,i)}$ for $i \ne j$ and $T(g)_{i,i}=\sqrt{\frac{2}{n}}\,g_{i,i}$.
A direct computation then gives, for all $g,h$,
\[
\schnorm{T(g)-T(h)}\le \frac{\sqrt2}{n}\,\lpnorm{g-h}.
\]
Therefore $F(g):=f(T(g))$ is $(L\sqrt2/n)$-Lipschitz on $\mathbb R^m$ in Euclidean norm.
Standard Herbst Gaussian concentration~\cite[Theorem 2.3.5]{anderson2010introduction} yields
\[
\Pr\Big(|F(g)-\mathbb EF(g)|\ge \delta\Big)
\le
2\exp\!\left(-\frac{\delta^2}{2(L\sqrt2/n)^2}\right)
=
2\exp\!\left(-\frac{n^2\delta^2}{4L^2}\right).
\]
Since $F(g)=f(A)$, this is the claim.
\end{proof}

\begin{lemma}[Mixed Lipschitzness for chaining]
\label{lem:mixed-lipschitz-gradient}
Fix $\beta>0$ and $\gamma>0$.  Let $A\in M_n(\mathbb R)$ be symmetric and let
$a,\tilde a\in[0,1]$.  Let $W,\widetilde W\in\mathcal D_n((0,1])$ and set
$D:=W^{-1}$ and $\widetilde D:=\widetilde W^{-1}$.  Let $T,\widetilde T\in\mathcal D_n(\mathbb R)$
satisfy $\schnorm{T}\le 1$ and $\schnorm{\widetilde T}\le 1$.
Define
\[
X:=a\Id_n-\beta A + D,\qquad \widetilde X:=\tilde a\,\Id_n-\beta A + \widetilde D,
\]
and assume $\opnorm{\beta A} \le 1 - \gamma$.
For each triple $(a,W,T)$ define the scalar functional
\[
\Phi_{a,W,T}(A)
\;:=\;
\tr_n\left(
T\,E_{\mathcal D_n}\bigl[D\,(a-\beta A + D)^{-2}D\right]
\Bigr)
=
\tr_n\left(T\,D\,X^{-2}D\right),
\]
where the equality holds since $T$ is diagonal and $E_{\mathcal D_n}$ is trace-preserving.

Let $\nabla_A \Phi_{a,W,T}(A)\in M_n(\mathbb R)$ denote the gradient with respect to the
Hilbert structure $\langle U,V\rangle=\tr_n(UV)$ on $M_n(\mathbb R)$, i.e.
$d\Phi[H]=\langle \nabla_A\Phi,\,H\rangle$ for all symmetric directions $H$.
Then
\begin{equation}
\label{eq:mixed-lipschitz-bound}
\schnorm{\nabla_A \Phi_{a,W,T}(A)-\nabla_A \Phi_{\tilde a,\widetilde W,\widetilde T}(A)}
\ \le\
2\beta\left(
\frac{12}{\gamma^4}\,|a-\tilde a|
\;+\;
\frac{32}{\gamma^4}\,\opnorm{W-\widetilde W}
\;+\;
\frac{4}{\gamma^3}\,\schnorm{T-\widetilde T}
\right).
\end{equation}
In other words, the map $(a,W,T)\mapsto \nabla_A\Phi_{a,W,T}(A)$ is Lipschitz in $\schnorm{\cdot}$ with constants depending only on
$\beta,\gamma$ and the operator bound $\opnorm{A}$, and independent of $n$.
\end{lemma}

\begin{proof}
Since $dX=-\beta dA$ and $d(X^{-1})=X^{-1}(dX)X^{-1}$, we have
\[
d(X^{-2})=\beta\bigl(X^{-1}(dA)X^{-2}+X^{-2}(dA)X^{-1}\bigr).
\]
Therefore, plugging in $dA[H] = H$ for the gradient with respect to $A$,
\begin{align*}
d\Phi_{a,W,T}(A)[H]
&=\tr_n\!\bigl(TD\,(d(X^{-2}))\,D\bigr)\\
&=\beta\,\tr_n\!\Bigl(TD\bigl(X^{-1}HX^{-2}+X^{-2}HX^{-1}\bigr)D\Bigr)\\
&=\beta\,\tr_n\!\Bigl(H\bigl(X^{-2}DTDX^{-1}+X^{-1}DTDX^{-2}\bigr)\Bigr),
\end{align*}
using cyclicity of $\tr_n$.
Thus
\begin{equation}
\label{eq:grad-formula-mixed}
\nabla_A \Phi_{a,W,T}(A)
=
\beta\Bigl(
X^{-2}\,D T D\,X^{-1}
+
X^{-1}\,D T D\,X^{-2}
\Bigr).
\end{equation}
Then apply \pref{fact:frechet-lipschitz} and \pref{lem:resolvent-frechet-dRjDTDRkj-T} with $dW \gets dW$ and $dM \gets \Id_n da$ and $dT \gets dT$, noting that the assumptions that $\opnorm{\beta A} \le 1 - \gamma$ and $\opnorm{W} \le 1$ and $a \in [0,1]$ imply $\opnorm{a\Id_n - \beta A} \le 2-\gamma$ as well as $X\succeq \gamma \Id_n$.
Use the resulting Lipschitz constants to integrate pathwise, using triangle inequality on path segments varying one parameter at a time.
\end{proof}

\begin{theorem}[Chaining bound for centered diagonal resolvent squares]
\label{thm:chaining-square-no-sqrtlog}
Fix $0<\beta<\tfrac12$ and $\gamma>0$.
Let $A$ be a normalized GOE matrix.
Choose an operator--norm cutoff $C_{\mathrm{op}}>0$ such that
\begin{equation}
\label{eq:chaining-gap-cutoff}
\beta C_{\mathrm{op}}\le 1-\gamma.
\end{equation}

Let $\Pi_{C_{\mathrm{op}}}:\mathrm{Sym}_n(\mathbb R)\to\mathrm{Sym}_n(\mathbb R)$ be the metric projection
(with respect to $\schnorm{\cdot}$) onto the closed convex set $\{M:\opnorm{M}\le C_{\mathrm{op}}\}$.
Set
\[
\overline A:=\Pi_{C_{\mathrm{op}}}(A).
\]

For $a\in[0,1]$ and $W\in\mathcal D_n((0,1])$ let $D:=W^{-1}$ and define
\[
X(a,W):=a\Id_n-\beta\overline A + D.
\]
For $T\in\mathcal D_n$ with $\schnorm{T}=1$, define the scalar functional
\[
\Phi_{a,W,T}(\overline A):=\tr_n\!\Big(T\,D\,X(a,W)^{-2}D\Big),
\qquad
Z_{a,W,T}:=\Phi_{a,W,T}(\overline A)-\E_A\big[\Phi_{a,W,T}(\overline A)\big].
\]
Let the index set be
\[
\mathcal{J} \;:=\; [0,1] \;\times\; \mathcal D_n((0,1]) \;\times\; \{T\in\mathcal D_n:\schnorm{T}=1\}.
\]

Then there exist universal constants $C_{\mathrm{ch}},C_{\mathrm{tail}}>0$ such that
\begin{equation}
\label{eq:chaining-expectation}
\E_A\left[\sup_{(a,W,T)\in \mathcal{J}}|Z_{a,W,T}|\right]
\ \le\
C_{\mathrm{ch}}\,
\frac{\beta}{\sqrt n}\,
\frac{1}{\gamma^4},
\end{equation}
and, for every $p\in(0,1)$,
\begin{equation}
\label{eq:chaining-highprob}
\Pr\Bigg\{
\sup_{(a,W,T)\in \mathcal{J}}|Z_{a,W,T}|
\ \le\
C_{\mathrm{ch}}\,
\frac{\beta}{\sqrt n}\,
\frac{1}{\gamma^4}
\;+\;
C_{\mathrm{tail}}\,
\frac{\beta}{n}\,
\frac{1}{\gamma^3}\,
\sqrt{\log\frac{2}{p}}
\Bigg\}
\ \ge\ 1-p.
\end{equation}

By diagonal norm duality,
\[
\sup_{a\in[0,1]}\sup_{W\in\mathcal D_n((0,1])}\schnorm{E_{\mathcal D_n}\!\Big[D\,X(a,W)^{-2}D\Big] - \E\left[E_{\mathcal D_n}\!\Big[D\,X(a,W)^{-2}D\Big]\right]}
\;=\;
\sup_{(a,W,T)\in \mathcal{J}}|Z_{a,W,T}|,
\]
so \eqref{eq:chaining-expectation}--\eqref{eq:chaining-highprob} also bound the left-hand side.
\end{theorem}

\begin{proof}
\ppart{Increment metric}
Define the weights
\[
\kappa_a:=\frac{24}{\gamma^4},
\qquad
\kappa_W:=\frac{64}{\gamma^4},
\qquad
\kappa_{T}:=\frac{8}{\gamma^3}.
\]
Define on $\mathcal{J}$ the metric
\begin{equation}
\label{eq:rho-def}
\rho\bigl((a,W,T),(\tilde a,\widetilde W,\widetilde T)\bigr)
:=
\kappa_a|a-\tilde a|+\kappa_W\opnorm{W-\widetilde W}+\kappa_{T}\schnorm{T-\widetilde T}.
\end{equation}

Fix two indices $\theta:=(a,W,T)$ and $\tilde\theta:=(\tilde a,\widetilde W,\widetilde T)$.
Consider the difference
\[
G_{\theta,\tilde\theta}(M)
:=
\Phi_{\theta}(M)-\Phi_{\tilde\theta}(M),
\qquad M\in\mathrm{Sym}_n(\mathbb R).
\]
Due to \eqref{eq:chaining-gap-cutoff}, \pref{lem:mixed-lipschitz-gradient} applies with $A\gets M$ and yields
\[
\schnorm{\nabla_M \Phi_{\theta}(M)-\nabla_M \Phi_{\tilde\theta}(M)}
\ \le\ \beta\,\rho(\theta,\tilde\theta).
\]
Since $\nabla_M G_{\theta,\tilde\theta}=\nabla_M\Phi_{\theta}-\nabla_M\Phi_{\tilde\theta}$, we conclude that
\begin{equation}
\label{eq:G-Lip-on-ball}
\|G_{\theta,\tilde\theta}\|_{\mathrm{Lip}(\schnorm{\cdot};\,\opnorm{M}\le C_{\mathrm{op}})}
\ \le\
\beta\,\rho(\theta,\tilde\theta).
\end{equation}

\ppart{Sub-Gaussian increments}
Define the globally defined function on $M_n$
\[
\bar{G}_{\theta,\tilde\theta}(A) \;:=\; G_{\theta,\tilde\theta}(\Pi_{C_{\mathrm{op}}}(A))
\;=\;
\Phi_{\theta}(\overline A)-\Phi_{\tilde\theta}(\overline A).
\]
Because $\Pi_{C_{\mathrm{op}}}$ is $1$--Lipschitz in $\schnorm{\cdot}$ and maps into $\{\opnorm{M} \le C_{\mathrm{op}}\}$,
\eqref{eq:G-Lip-on-ball} implies
\[
\|\bar{G}_{\theta,\tilde\theta}\|_{\mathrm{Lip}(\schnorm{\cdot})}
\ \le\
\beta\,\rho(\theta,\tilde\theta).
\]
Now note
\[
Z_{\theta}-Z_{\tilde\theta}
\;=\;
\bar{G}_{\theta,\tilde\theta}(A)-\E[\bar{G}_{\theta,\tilde\theta}(A)].
\]
Applying the GOE Herbst inequality (\pref{lem:herbst-goe}) to $\bar{G}_{\theta,\tilde\theta}$
gives, for all $u>0$,
\[
\Pr\big(|Z_{\theta}-Z_{\tilde\theta}|>u\big)
\;\le\;
2\exp\!\left(-\frac{n^2u^2}{4\beta^2\rho(\theta,\tilde\theta)^2}\right)
\;\le\;
2\exp\!\left(-\frac{u^2}{2d(\theta,\tilde\theta)^2}\right),
\]
where we set the increment metric
\begin{equation}
\label{eq:d-def}
d(\theta,\tilde\theta)\;:=\;\frac{\sqrt{2}\beta}{n}\,\rho(\theta,\tilde\theta).
\end{equation}
Thus $(Z_\theta)_{\theta\in \mathcal{J}}$ has centered subgaussian increments with respect to $(\mathcal{J},d)$.

\ppart{Talagrand's chaining}
Since $\mathcal{J}$ is compact (take its closure to $[0,1] \times \mathcal D_n([0,1]) \times \{T\in\mathcal D_n:\schnorm{T}=1\}$ and extend $Z_{\theta}$ continuously since $Z_{\theta}$ is Lipschitz in $\theta$), it is separable.  By~\pref{thm:talagrand-subgaussian-chaining},
\[
\E\left[\sup_{\theta\in \mathcal{J}} Z_\theta\right]\le C\,\gamma_2(\mathcal{J},d),
\qquad
\E\left[\sup_{\theta\in \mathcal{J}} (-Z_\theta)\right]\le C\,\gamma_2(\mathcal{J},d),
\]
and hence
\begin{equation}
\label{eq:Esup-abs}
\E\left[\sup_{\theta\in \mathcal{J}} |Z_\theta|\right]\le 2\,C\,\gamma_2(\mathcal{J},d).
\end{equation}

\ppart{Entropy integral}
By~\pref{thm:talagrand-gamma-2},
\[
\gamma_2(\mathcal{J},d)\;\le\; C\int_0^{\mathrm{diam}(\mathcal{J},d)} \sqrt{\log N(\mathcal{J},d,\varepsilon)}\,d\varepsilon.
\]
Since $d=(\sqrt{2}\beta/n)\rho$, we have $\mathrm{diam}(\mathcal{J},d)=(\sqrt{2}\beta/n)\mathrm{diam}(\mathcal{J},\rho)$ and
$N(\mathcal{J},d,\varepsilon)=N\bigl(\mathcal{J},\rho,\varepsilon n/(\sqrt{2}\beta)\bigr)$. Therefore
\begin{equation}
    \label{eq:chaining-gamma2-bound-by-entropy-integral}
\gamma_2(\mathcal{J},d)
\;\le\;
C\frac{\sqrt{2}\beta}{n}\int_0^{\mathrm{diam}(\mathcal{J},\rho)}\sqrt{\log N(\mathcal{J},\rho,u)}\,du
.
\end{equation}

To bound $N(\mathcal{J},\rho,u)$, note that $\rho$ is a sum metric, so for every $u>0$,
\[
N(\mathcal{J},\rho,u) \;\le\;
N\Big([0,1],|\cdot|,\frac{u}{3\kappa_a}\Big)\cdot
N\Big(\mathcal{D}_n([0,1]),\opnorm{\cdot},\frac{u}{3\kappa_W}\Big)\cdot
N\Big(\diag(\sqrt{n}S^{n-1}),\schnorm{\cdot},\frac{u}{3\kappa_{T}}\Big),
\]
where $\diag(\sqrt{n}S^{n-1})$ denotes the set of diagonal matrices whose vector of diagonal entries is on the radius-$\sqrt{n}$ Euclidean sphere in $\mathbb R^n$, equivalently the set of diagonal matrices $T$ satisfying $\schnorm{T} = 1$. 
Using the standard volumetric bounds,
\[
N([0,1],|\cdot|,r)\;\le\; 1+\frac{1}{r},
\quad\;\;
N(\mathcal{D}_n([0,1]),\opnorm{\cdot},r)\;\le\; \left(\frac{3}{r}\right)^n,
\quad\;\;
N(S^{n-1},\schnorm{\cdot},r)\;\le\; \left(\frac{3}{r}\right)^n,
\]
we obtain for $u\in(0,1]$,
\[
\log N(\mathcal{J},\rho,u)\ \le\ (2n+1)\log\Big(\frac{C_1}{u}\Big),
\qquad
C_1:=9\max\{\kappa_a,\kappa_W,\kappa_{T}\}.
\]
Consequently,
\[
\int_0^{\mathrm{diam}(\mathcal{J},\rho)}\sqrt{\log N(\mathcal{J},\rho,u)}\,du
\ \le\
\int_0^{C_1}\sqrt{(2n+1)\log\Big(\frac{C_1}{u}\Big)}\,du
\ \le\
C_1\sqrt{2n+1}\int_0^\infty \sqrt{t}\,e^{-t}\,dt
=
C_1\sqrt{2n+1}\,\frac{\sqrt\pi}{2}.
\]
Substituting this into \pref{eq:chaining-gamma2-bound-by-entropy-integral} yields
\[
\gamma_2(\mathcal{J},d)
\ \le\
C'\,\frac{\sqrt{2}\beta}{n}\,C_1\sqrt{n}
\ \le\
C''\,\frac{\beta}{\sqrt n}\,\frac{1}{\gamma^4},
\]
where in the last step we used that $C_1\lesssim 1/\gamma^4$ (since $\kappa_W$ dominates).
Together with \eqref{eq:Esup-abs}, this proves \eqref{eq:chaining-expectation}.

\ppart{High-probability bound}
Define the supremum functional
\[
\mathcal Z(A):=\sup_{\theta\in \mathcal{J}} Z_\theta
=
\sup_{\theta\in \mathcal{J}}\Big(\Phi_\theta(\Pi_{C_{\mathrm{op}}}(A))-\E[\Phi_\theta(\Pi_{C_{\mathrm{op}}}(A))]\Big).
\]
For each fixed $\theta$, \pref{lem:resolvent-frechet-dImRDTDR} (or \pref{eq:grad-formula-mixed})
implies $\Phi_\theta(M)$ is $L_0$--Lipschitz in $M$ (in $\schnorm{\cdot}$) on $\{\opnorm{M}\le C_{\mathrm{op}}\}$ with
\[
L_0\ \le\ \frac{8\beta}{\gamma^3}.
\]
Composing with the $1$--Lipschitz map $\Pi_{C_{\mathrm{op}}}$ gives that $A\mapsto \Phi_\theta(\Pi_{C_{\mathrm{op}}}(A))$
is $L_0$--Lipschitz in $A$, uniformly in $\theta$. Hence $\mathcal Z$ is $L_0$--Lipschitz as a supremum of $L_0$--Lipschitz
functions. Applying \pref{lem:herbst-goe} to $\mathcal Z$ yields
\[
\Pr\big\{\mathcal Z-\E\mathcal [Z] > t\big\}\le \exp\!\left(-\frac{n^2t^2}{4L_0^2}\right).
\]
Taking $t= \frac{2L_0}{n}\sqrt{\log(2/p)}$ gives the one-sided bound, and then applying the same to $-\mathcal Z$ and using
$\sup|Z_\theta|=\max\{\sup Z_\theta,\sup(-Z_\theta)\}$ yields \eqref{eq:chaining-highprob}.
\end{proof}

\subsection{Free convolutions and analytic subordination}
\label{sec:free-prob}
To prepare ourselves for the free interpolation argument in \pref{sec:free-interpolation}, we briefly recall elementary definitions from free probability~\cite{nica2006lectures,mingo2017free} and state Belinschi's result for subordination of free sums~\cite{biane1998processes,belinschi2005complex}.

Here, the upper-half of the complex plane is denoted as $\mathbb{H}$, and its topological closure is $\bar{\mathbb{H}}$.

\begin{definition}[Non-commutative probability space]\label{def:non-comm-prob-space}
    A \emph{tracial non-commutative probability space} is a pair $(\mathcal{M},\tau)$ where $\mathcal{M}$ is a von Neumann algebra and $\tau$ is a faithful normal tracial state.
\end{definition}

Self-adjoint elements $X \in \mathcal{M}$ have a well-defined unique compactly supported measure $\mu_X$ such that
\[
\tau[f(X)] = \int_{\R} f\,d\mu_X,
\]
for polynomials $f$; $\mu_X$ is termed the \emph{(spectral) distribution of $X$}.
We also let $\Spec(X) := \supp(\mu_X)$.

\begin{definition}[Non-commutative $L^p$ space]\label{def:non-comm-lp-space}
    Given a tracial non-commutative probability space $(\mathcal{M},\tau)$ and $p \in [1,\infty)$, the $L^p$ norm $\Lpnorm[p]{X} := (\tau(|X|^p))^{1/p}$ where $|X| = (X^*X)^{1/2}$. $L^p(\mathcal{M},\tau)$ denotes the completion of $\mathcal{M}$ with respect to $\Lpnorm[p]{\cdot}$. 
\end{definition}  

These norms satisfy the non-commutative H\"older's inequality:  If $1/p = 1/p_1 + \dots + 1/p_k$, then
\[
\Lpnorm[p]{x_1 \dots x_k} \leq \Lpnorm[p_1]{x_1} \dots \Lpnorm[p_k]{x_k}\,.
\]
Furthermore,
\[
    \Lpnorm[\infty]{x} = \lim_{p \to \infty} \Lpnorm[p]{x}\,,
\]
with $\Lpnorm[\infty]{x}$ being the operator norm.
$L^2(\mathcal{M},\tau)$ is a Hilbert space with inner product $\la x,y \ra_\tau = \tau(x^*y)$.  Left multiplication by $x$ defines a representation of $\mathcal{M}$ on $L^2(\mathcal{M},\tau)$, i.e.\ a $*$-homomorphism $\mathcal{M} \to B(L^2(\mathcal{M},\tau))$.  

\begin{definition}[Non-commutative conditional expectation]\label{def:non-comm-conditional-expectation}
    Let $(\mathcal{M},\tau)$ be a tracial non-commutative probability space and $\mathcal{A} \subseteq \mathcal{M}$ be a von Neumann subalgebra of $\mathcal{M}$.  Then, $L^2(\mathcal{A},\tau|_{\mathcal{A}})$ is a subspace of $L^2(\mathcal{M},\tau)$.  The orthogonal projection $E_{\mathcal{A}}$ onto $L^2(\mathcal{A},\tau)$ restricts to a mapping $\mathcal{M} \to \mathcal{A}$, and is termed the \emph{canonical conditional expectation onto $\mathcal{A}$}.
\end{definition}    
We will use $\mathcal{M} = M_n(\mathbb{C})$ and $\mathcal{A} = \mathcal{D}_n$. Consequently, $E_{\mathcal{D}_n}(X)$ is the diagonal matrix obtained by zeroing out the off-diagonal entries of $X$. Note that $E_\calA$ is a contraction for every $L^p$ norm.

A tracial non-commutative probability space $(\mathcal{M},\tau)$ admits a non-commutative version of independence, known as \emph{free independence} \cite{voiculescu1985symmetries,voiculescu1986addition}.  
\begin{definition}[Free independence]\label{def:free-independence}
    Let $(\mathcal{M},\tau)$ be a tracial non-commutative probability space, and $\mathcal{A}_1$, \dots, $\mathcal{A}_d$ be $*$-subalgebras of $\calM$. Then, $\mathcal{A}_1$, \dots, $\mathcal{A}_d$ are \emph{freely independent} if whenever $i_1$, \dots, $i_k \in [d]$ with $i_1 \neq i_2 \neq i_3 \neq \dots \neq i_k$ and $X_j \in \mathcal{A}_{i_j}$, it is the case that
    \[
        \tau \left[ (X_1 - \tau(X_1)) \dots (X_k - \tau(X_k)) \right] = 0.
    \]
\end{definition}
Elements, tuples, or sets in $\mathcal{M}$ are said to be freely independent if the $*$-subalgebras that they generate are freely independent. If $X_1$, \dots, $X_d$ are freely independent self-adjoint operators, then the joint moments $\tau(X_{i_1} \dots X_{i_k})$ are uniquely determined by the individual moments $\tau(X_i^k)$.  


Voiculescu's \emph{asymptotic} freeness theory (roughly) implies that independent random matrices distributionally invariant under unitary (or orthogonal) conjugation become freely independent as $n \to \infty$ \cite{voiculescu1991limit-laws,voiculescu1998strengthened}. 

\begin{theorem}[Asymptotic freeness of GOE and deterministic matrix~{\cite[Theorem 5.4.2]{anderson2010introduction}}] \label{thm: asymptotic freeness}
Let $Z_n$ be a GOE random matrix, and let $D_n$ be a deterministic matrix with $\norm{D_n} \leq M$ for some constant $M$.  Assume that the empirical spectral distribution $\mu_{D_n}$ converges to some $\mu$ as $n \to \infty$.  Consider a tracial non-commutative probability space $\mathcal{M}$ generated by two freely independent self-adjoint elements $S$ and $D$, with spectral distributions $(1/2\pi) \one_{[-2,2]}(t) \sqrt{4 - t^2}\,dt$ and $\mu$ respectively.  Then, almost surely, for non-commutative polynomials $f$ in two variables,
\[
\lim_{n \to \infty} \tr_n[f(Z_n,D_n)] = \tau[f(S,D)].
\]
\end{theorem}

If $X$ and $Y$ have distributions $\mu$ and $\nu$ respectively, then the distribution of $X+Y$ is called the \emph{free convolution of $\mu$ and $\nu$} and is denoted $\mu \boxplus \nu$.  
Free convolution is computed using several complex-analytic functions induced by the measure $\mu$. 

\begin{definition}[Cauchy--Stieltjes Transform]
\label{def:cs-transform}
The \emph{Cauchy--Stieltjes transform} of a probability measure $\mu$ is the function
\[
G_\mu(z) = \int_{\R} \frac{1}{z - x}\,d\mu(x),
\]
defined for $z \in \mathbb{C} \setminus \supp(\mu)$. Similarly, for self-adjoint $X \in \mathcal{M}$,
\[
    G_X(z) = \tau[(z - X)^{-1}] \text{ for } z \in \mathbb{C} \setminus \Spec(X),
\]
which is easily seen to be the Cauchy--Stieltjes transform of the spectral distribution of $X$.
\end{definition}

The following are standard properties of the Cauchy--Stieltjes transform.

\begin{proposition}[Properties of Cauchy--Stieltjes transform~{\cite[\S 6,~Theorem~3]{akhiezer1963theory}}]
\label{prop:cs-transform}
Let $\mu \in \mathcal{P}(\R)$ and $z \in \C$.
\begin{enumerate}
    \itemsep0.3em
    \item $G_\mu(z)$ maps the upper half-plane into the lower half-plane and vice versa.
    \item $\displaystyle |G_\mu(z)| \leq \frac{1}{d(z,\supp(\mu))}$. 
    \item A point $a \in \R$ is \emph{not} in the support of $\mu$ if and only if there is an analytic function defined in a neighborhood $O$ of $a$ that agrees with $G_\mu$ on $O \setminus \R$.
    \item If $\mu$ is compactly supported, then
    \[
    \lim_{z \to \infty} z G_\mu(z) = 1.
    \]
    In particular, $G_\mu^{-1}$ is defined in a neighborhood of $0$ and $G_\mu^{-1}(z) - 1/z$ is analytic near $0$.
    \item 
For the semicircular measure $\mu_{\mathrm{sc}}$ of variance 1, we have
\begin{equation}\label{eq:Gsc-Stieltjes}
G_{\mu_{\mathrm{sc}}}(z)\ =\ \frac{z-\sqrt{z^2-4}}{2},
\qquad z\in\mathbb{C}\setminus[-2,2],
\end{equation}
where the branch of the square root is chosen to satisfy $G_{\mu_{\mathrm{sc}}}(z)\sim 1/z$ as $|z|\to\infty$. 
\end{enumerate}
\end{proposition}

\begin{theorem}[Analytic subordination for free sums~{\cite[Theorem 3.1]{biane1998processes}}] \label{thm: subordination}
Let $\mathcal{A}$ and $\mathcal{B}$ be freely independent subalgebras of $\mathcal{M}$.  Let $X \in \mathcal{A}$ and $Y \in \mathcal{B}$ be freely independent self-adjoint operators in a tracial von Neumann algebra $(\mathcal{M},\tau)$.
\begin{enumerate}
    \itemsep0.3em
    \item There exists a unique function $F: \mathbb{H} \to \mathbb{H}$ such that $F(z) = z + O(1)$ for sufficiently large $z$ and $G_{X+Y}(z) = G_X(F(z))$.
    \item Let $E_{\mathcal{A}}$ denote the canonical conditional expectation onto $\mathcal{A}$.  Then for $z \in \mathbb{H}$,
    \[
    E_{\mathcal{A}}[(z - X + Y)^{-1}] = (F(z) - X)^{-1}.
    \]
\end{enumerate}
\end{theorem}

\begin{lemma}[Properties of the subordination function for semicircular operators]
\label{lem:subordination-semicircular}
As in \cite{jekel2024pha}, let $\mathcal{B} = L^\infty[-2,2]$, let $\tau_{\mathcal{B}}(f) = \frac{1}{2\pi} \int_{-2}^2 f(x)\sqrt{4 - x^2}\,dx$, and let $S \in \mathcal{B}$ be the identity function, so that $S$ is a standard semicircular element. Let $(\mathcal{M}_n,\tau_n)$ be the free product of $(M_n(\mathbb{C}),\tr_n)$ with $(\mathcal{B},\tau_{\mathcal{B}})$, so that $\mathcal{M}_n$ is generated by its subalgebra $M_n(\mathbb{C})$ and a freely independent semicircular element $S$.

Fix $0<\beta<\tfrac12$ and let $D \in M_n(\mathbb{R})$ satisfy $D\succeq \Id_n$. 
Let $F_{\beta,-D}:\mathbb H\to\mathbb H$
be the subordination function of $G_{\beta S-D}$ subordinate to $G_{-D}$ from \pref{thm: subordination}, i.e.
$G_{\beta S-D}(z)=G_{-D}(F_{\beta,-D}(z))$.

Then the following properties are true of $F_{\beta,-D}$:
\begin{enumerate}[label=(\alph*), ref=\thetheorem(\alph*)]
    \item \label{lem:subordination-semicircular-continuation}
    $F_{\beta,-D}$ extends analytically onto the non-negative real ray, i.e. $F_{\beta,-D}(a)$ exists and is differentiable for $a \in \R$ satisfying $a \ge 0$.
    \item \label{lem:subordination-semicircular-0}
    \[F_{\beta,-D}\left(\beta^2\tr_n(D^{-1})\right) = 0.\]
    \item \label{lem:subordination-semicircular-derivative}
    \[F_{\beta,-D}'\left(\beta^2\tr_n(D^{-1})\right) = \frac{1}{1-\beta^2\tr_n(D^{-2})}.\]
\end{enumerate}
\end{lemma}
\begin{proof}
        Combining~\cite[Theorem 1.23 (3)]{belinschi2005complex} with~\cite[Theorem 3.1]{biane1998processes} yields that $F_{\beta,-D}: \bar{\mathbb{H}} \to \C$ is continuous for every $z \in \bar{\mathbb{H}}$.
        Let $\mu$ be the semicircular law of variance $\beta^2$ (and therefore the spectral law of the operator $\beta S$) and let $\nu$ be the spectral law of $-D$.
    Now, $\mu \boxplus \nu$ has compact support since both $\mu$ and $\nu$ do, and also $\mathsf{supp}\left(\mu \boxplus \nu\right) < 0$  since $\supp(\nu) \le -1$ and $\supp(\mu) < 2\beta$.
    Therefore by \pref{prop:cs-transform}, $G_{\mu \boxplus \nu}(z) = G_{\nu}(F_{\beta,-D}(z))$ is analytic for $z$ in a neighborhood of every non-negative real number. Consequently, the identity
    \begin{equation}
    \label{eq:subordination-inverse}
        z = F_{\beta,-D}(z) + \beta^2 G_{-D}(F_{\beta,-D}(z))
    \end{equation}
    in the proof of~\cite[Proposition 4.10]{jekel2024pha} extends continuously to $z \in \{a \in \R: a \ge 0\}$.
    
    For \pref{lem:subordination-semicircular-0}, setting $F_{\beta,-D}(a) = 0$ by its unique analytic continuation onto $z \in \bar{\mathbb{H}}\setminus\mathsf{supp}\left(\mu_{\beta^2}\boxplus \nu\right)$,
    \[
        \begin{aligned}
            a &=  0 + \beta^2 G_{-D}(0) = \beta^2\tr_n\left[\left( 0 + D\right)^{-1}\right] = \beta^2\tr_n[D^{-1}].
        \end{aligned}
    \]
    at which point we invoke invertibility via the identity theorem on the unique analytic continuation and assert that the above choice of $a$ is unique, since, as \pref{lem:subordination-semicircular-derivative} will show, $F_{\beta,-D}$ is strictly increasing and therefore injective on the real non-negative numbers\footnote{ In fact, the proof of~\cite[Proposition 1.22]{belinschi2005complex} directly invokes the analyticity of $F$ in a neighborhood of $z$ not in the support of the freely convolved measure.}.
    
    For \pref{lem:subordination-semicircular-derivative}, we implicitly differentiate \pref{eq:subordination-inverse} to find, letting $a := \beta^2\tr_n(D^{-1})$ so that $F_{\beta,-D}(a) = 0$,
    \[ 1 = F_{\beta,-D}'(z) + \beta^2 G_{-D}'(F_{\beta,-D}(z))F_{\beta,-D}'(z) = F_{\beta,-D}'(z)(1 + \beta^2 G_{-D}'(0))
     \]
     Dividing both sides by $1 + \beta^2G_{-D}'(0)$ and noting that $G_{D}'(0) = -\tr_n(D^{-2})$ completes the argument.
    \qedhere
\end{proof}

\subsection{Strong convergence for GOE tensor sum}
\label{sec:AIk+Bk}

We will soon embark on a free interpolation argument, using Gaussian integration by parts to characterize the interpolated resolvent
\[R_k(t) := (a\Id_{kn}-\beta W_k(t)+D \ot \Id_k)^{-1}, \qquad \text{with} \qquad W_k(t) := \sqrt{1-t}(A\ot \Id_k) + \sqrt{t}B_k.\]
Then $R_k(0)$ is simply $k$ direct-sum copies of the object of interest $\bar{Q}$ from \pref{sec:alg-properties}, $R_k(1)$ is an object that approaches an ideal free convolution as $k \to \infty$, and if time-derivatives of observables $\frac{d}{dt}\phi(R_k(t))$ of the interpolated resolvent are small, we conclude that $\phi(R_k(0))$ is close to $\phi(R_k(1))$ in expectation.

Before we take on this journey, we give a high-probability bound on the upper limit of the spectrum of $W_k(t)$ in \pref{cor:AIk+Bk-final-bound}, in order to ensure that $R_k(t)$ is defined and analytic throughout the entire interpolation.

Several recent advances in strong convergence in free probability~\cite{pisier2016strong,bandeira2023matrix,van2025strong,jekel2025strong} address control of the edge of the spectra of matrices similar to $W_k(t)$; however the specific tensor structure of $W_k(t)$ is incompatible with the application of these bounds (see \cite[Section 8]{bandeira2023matrix}).

Thus we take a two-stage approach, where we first condition on a specific instantiation $A_0$ of the GOE matrix $A$, show that if the spectrum of $A_0$ is close to an ideal semicircular spectrum then the extremal eigenvalues of $W_k(t)$
are well controlled (\pref{lem:strong-conv-to-a0+s} and \pref{lem:edge-uniform-close-2-Stieltjes}), then show that with high probability, $A_0$ is close enough to get that control (\pref{lem:GOE-free-edge}).

\begin{lemma}
    \label{lem:strong-conv-to-a0+s}
    Let $A_0 \in \R^{n \times n}$ be a fixed symmetric matrix.
    For each $k \in \N$, let $B_k \in \R^{(nk) \times (nk)}$ be a GOE matrix. 
    Then 
    \[ \sup\Spec\left(\sqrt{1-t}A_0 \otimes \Id_k + \sqrt{t}B_k\right) \xrightarrow[k\to\infty]{a.s.} \sup\Spec\left(\sqrt{1-t}A_0 + \sqrt{t}S\right), \]
    where $S$ is a freely independent semicircular operator of variance 1 as in \pref{lem:subordination-semicircular}.
\end{lemma}
\begin{proof}
    By \cite[Theorem 1]{lehner1999computing} (as stated in \cite[(3.1)]{van2025strong} or \cite[Lemma 2.4]{bandeira2023matrix}),
    \[ \lim_{k \to \infty} d_H(\Spec(\sqrt{1-t}A_0 \otimes \Id_k + \sqrt{t}B_k), \Spec(\sqrt{1-t}A_0 \otimes \Id + \sqrt{t}X)) = 0\quad\text{a.s.},\]
    where $d_H$ is the Hausdorff distance and $X$ is an $n \times n$ matrix whose $(i,j)$th entry is $S_{\min(i,j),\max(i,j)}/\sqrt{n}$ if $i \ne j$ and $\sqrt{2}S_{i,i}/\sqrt{n}$ if $i = j$, where $\{S_{i,j}\}_{1 \le i \le j \le n}$ is a family of freely independent semicircular operators of variance 1.

    By calculating traces, we can see that $X$ is itself also a freely independent ideal semicircular operator of variance 1, as $\tau(X^p)$ is a sum over monomials of the form $\tau(S_{i_1,i_2}S_{i_2,i_3}\dots S_{i_{p-1},i_p}S_{i_p,i_1})$, and the only non-zero monomials are those corresponding to non-crossing partitions, by the definition of free independence.
    Furthermore, since $\Spec(\sqrt{1-t}A_0 \otimes \Id) = \Spec(\sqrt{1-t}A_0)$, they have the same free additive convolutions with semicircular measures, and so $\Spec(\sqrt{1-t}A_0 \otimes \Id + \sqrt{t}X) = \Spec(\sqrt{1-t}A_0 + \sqrt{t}S)$.

    Then convergence of the spectra in Hausdorff distance implies convergence of the largest points in their supports.
\end{proof}

\begin{lemma}[Uniform right-edge control for $\sqrt{1-t}A_0 + \sqrt{t}S$]\label{lem:edge-uniform-close-2-Stieltjes}
Let $A=A^\sT \in M_n(\R)$ have spectral law $\mu_A$, and let
$S$ be a semicircular element of variance $1$ that is free from $A$, as in \pref{lem:subordination-semicircular}.
For $t\in[0,1]$ set
\[
X_t \;:=\;\sqrt{1-t}\,A+\sqrt t\,S,
\qquad
\mu_t \;:=\;\mu_{X_t},
\qquad
E_+(t)\;:=\;\sup\supp(\mu_t).
\]
Assume that for some $\delta\ge 0$ and $\varepsilon\ge 0$,
\begin{equation}\label{eq:assump-support-W1-Stieltjes}
\sup\supp(\mu_A) \le 2+\delta,
\qquad
W_1(\mu_A,\mu_{\mathrm{sc}})\le \varepsilon,
\end{equation}
where $\mu_{\mathrm{sc}}$ is the standard semicircle law. 
Then, for every $t\in[0,1]$,
\begin{equation}\label{eq:Eplus-uniform-bound-Stieltjes}
E_+(t)\ \le\ 2\;+\;\delta\;+\;2\,\varepsilon^{1/3}.
\end{equation}
\end{lemma}

\begin{proof}

For $z\in\mathbb{C}^+$, $\Im G_\mu(z)<0$.
Since $E_+(1)=2$ trivially at $t = 1$, fix $t\in[0,1)$. 

\ppart{Using the edge-variational formula}
Set
\[
s\ :=\ \frac{t}{1-t}\ \in[0,\infty),
\qquad
Y_s\ :=\ A+\sqrt{s}\,S.
\]
Then
\[
X_t=\sqrt{1-t}\,Y_s=\frac{1}{\sqrt{1+s}}\,Y_s.
\]
Let $\nu_s:=\mu_{Y_s}=\mu_A\boxplus \mu_{\mathrm{sc},s}$ where $\mu_{\mathrm{sc},s}$ is the semicircular law of variance $s$ and write
\[
r_s\ :=\ \sup\supp(\nu_s).
\]
By scaling of supports,
\begin{equation}\label{eq:Eplus-scale-Stieltjes}
E_+(t)\ =\ \sqrt{1-t}\,r_s\ =\ \frac{r_s}{\sqrt{1+s}}.
\end{equation}

A standard subordination identity for free convolution with a semicircle of variance $s$~\cite[Proof of Proposition 4.10]{jekel2024pha} gives
\[
G_{\nu_s}(z)\ =\ G_{\mu_A}\!\bigl(z-s\,G_{\nu_s}(z)\bigr),\qquad z\in\mathbb{C}^+.
\]
As in \cite[Lemma 4.4]{capitaine2009largest}
and \cite{biane1997free}, taking $u(z) \gets z-s\,G_{\nu_s}(z)$, we see that $z = u(z) + sG_{\mu_A}(u(z))$ and so
\begin{equation}\label{eq:rs-variational-Stieltjes}
r_s\ =\ \inf_{u>R_A}\Big\{u+s\,G_{\mu_A}(u)\Big\},
\qquad R_A:=\sup\supp(\mu_A).
\end{equation}
In particular, for any $u>R_A$,
\begin{equation}\label{eq:rs-testpoint-Stieltjes}
r_s\ \le\ u+s\,G_{\mu_A}(u).
\end{equation}

\ppart{Semicircular test point} 
Define, for $s\ge 0$,
\[
u_s\ :=\ \frac{s+2}{\sqrt{1+s}}
\]
A direct computation gives
\[
u_s^2-4=\frac{(s+2)^2-4(1+s)}{1+s}=\frac{s^2}{1+s},
\qquad\text{so}\qquad
\sqrt{u_s^2-4}=\frac{s}{\sqrt{1+s}},
\]
where the last equality uses that $u_s>2$ for $s>0$ and the branch choice $\sqrt{x^2-4}>0$ for $x>2$.
This choice of branch allows us to plug into \eqref{eq:Gsc-Stieltjes} to yield
\begin{equation}\label{eq:Gsc-us-Stieltjes}
G_{\mu_{\mathrm{sc}}}(u_s)
=\frac{u_s-\sqrt{u_s^2-4}}{2}
=\frac{\frac{s+2}{\sqrt{1+s}}-\frac{s}{\sqrt{1+s}}}{2}
=\frac{1}{\sqrt{1+s}}.
\end{equation}
Therefore,
\begin{equation}\label{eq:semicircle-min-value-Stieltjes}
u_s+s\,G_{\mu_{\mathrm{sc}}}(u_s)
=
\frac{s+2}{\sqrt{1+s}}+\frac{s}{\sqrt{1+s}}
=
2\sqrt{1+s},
\end{equation}
which is exactly the right edge of $\mu_{\mathrm{sc}}\boxplus\mu_{\mathrm{sc},s}=\mu_{\mathrm{sc},1+s}$.

\ppart{Comparing $G_{\mu_A}$ to $G_{\mu_{\mathrm{sc}}}$ via $W_1$ convergence}
Fix $\eta>0$ and set the test point
\[
u\ :=\ u_s+\delta+\eta.
\]
By \eqref{eq:assump-support-W1-Stieltjes}, $R_A\le 2+\delta\le u_s+\delta<u$, hence $u>R_A$ and
\eqref{eq:rs-testpoint-Stieltjes} applies. Decompose
\begin{align}
u+sG_{\mu_A}(u)
&=\bigl(u+sG_{\mu_{\mathrm{sc}}}(u)\bigr)
+s\bigl(G_{\mu_A}(u)-G_{\mu_{\mathrm{sc}}}(u)\bigr).
\label{eq:decompose-Stieltjes}
\end{align}

\emph{(a) Semicircle part.}
For $x>2$, one has $G'_{\mu_{\mathrm{sc}}}(x)<0$, hence the function
\[
h(x):=x+sG_{\mu_{\mathrm{sc}}}(x)
\]
satisfies $h'(x)=1+sG'_{\mu_{\mathrm{sc}}}(x)\le 1$ for all $x>2$.
Moreover, $h'(u_s)=0$, and for $x\ge u_s$ one has
$h'(x)\ge 0$, so $h$ is increasing on $[u_s,\infty)$ with slope at most $1$.
Therefore, by \pref{eq:semicircle-min-value-Stieltjes},
\begin{equation}
\label{eq:u+sGu-bound}
u+sG_{\mu_{\mathrm{sc}}}(u)=h(u)\le h(u_s)+(u-u_s)
=2\sqrt{1+s}+(\delta+\eta).
\end{equation}

\emph{(b) Difference term via $W_1$.}
For fixed $u>2+\delta$, the function $f_u(x):=(u-x)^{-1}$ is Lipschitz on
$\supp(\mu_A)\cup\supp(\mu_{\mathrm{sc}})\subset(-\infty,2+\delta]$ with
\[
\mathrm{Lip}(f_u)
\le\sup_{x\le 2+\delta}\frac{1}{(u-x)^2}
=\frac{1}{(u-(2+\delta))^2}.
\]
By Kantorovich--Rubinstein duality,
\[
|G_{\mu_A}(u)-G_{\mu_{\mathrm{sc}}}(u)|
=\left|\int f_u\,d\mu_A-\int f_u\,d\mu_{\mathrm{sc}}\right|
\le \mathrm{Lip}(f_u)\,W_1(\mu_A,\mu_{\mathrm{sc}})
\le \frac{\varepsilon}{(u-(2+\delta))^2}.
\]
But $u-(2+\delta)=u_s-2+\eta\ge \eta$, hence
\begin{equation}\label{eq:G-diff-Stieltjes}
|G_{\mu_A}(u)-G_{\mu_{\mathrm{sc}}}(u)|\ \le\ \frac{\varepsilon}{\eta^2}.
\end{equation}

\ppart{Optimize $\eta$}
Combining \pref{eq:rs-testpoint-Stieltjes}, \pref{eq:decompose-Stieltjes}, and the bounds \pref{eq:u+sGu-bound} and \pref{eq:G-diff-Stieltjes} gives
\[
r_s
\le
u+sG_{\mu_A}(u)
\le
2\sqrt{1+s}+\delta+\eta+\frac{s\varepsilon}{\eta^2}.
\]
Choose $\eta=(2s\varepsilon)^{1/3}$. Then
\[
\eta+\frac{s\varepsilon}{\eta^2}
=
(2s\varepsilon)^{1/3}+\frac{s\varepsilon}{(2s\varepsilon)^{2/3}}
=
\Bigl(2^{1/3}+2^{-2/3}\Bigr)(s\varepsilon)^{1/3}
\le 2\,(s\varepsilon)^{1/3},
\]
so
\[
r_s\ \le\ 2\sqrt{1+s}+\delta+2\,(s\varepsilon)^{1/3}.
\]
Divide by $\sqrt{1+s}$ and use \eqref{eq:Eplus-scale-Stieltjes}:
\[
E_+(t)=\frac{r_s}{\sqrt{1+s}}
\le
2+\frac{\delta}{\sqrt{1+s}}+2\,\frac{s^{1/3}}{\sqrt{1+s}}\,\varepsilon^{1/3}
\le
2+\delta+2\,\varepsilon^{1/3},
\]
since $s^{1/3}/\sqrt{1+s}\le 1$ for all $s\ge 0$.
This proves \eqref{eq:Eplus-uniform-bound-Stieltjes} for $t\in[0,1)$ and the $t=1$ case is trivial.
\end{proof}

\begin{lemma}[Edge control from GOE concentration of $A$]\label{lem:GOE-free-edge}
Let $A\in M_n(\R)$ be a normalized GOE matrix
and let $S$ be a semicircular operator of variance $1$ that is freely independent of $A$, as in \pref{lem:subordination-semicircular}.
For $t\in[0,1]$ define
\[
X_t:=\sqrt{1-t}\,A+\sqrt t\,S,\qquad \mu_t:=\mu_{X_t},\qquad E_+(t):=\sup\supp(\mu_t).
\]
Let $\mu_{\mathrm{sc}}$ be the standard semicircle law and  $\mu_A$ denote the empirical eigenspectrum of $A.$

Then for every $\xi>1$ there exists $L=L(\xi)>0$ and constants $c,C,C'>0$ such that, with probability
at least $1-C\exp\!\big(-c(\log n)^\xi\big)$, the following holds simultaneously for all $t\in[0,1]$:
\begin{equation}\label{eq:edge-bound-GOE}
E_+(t)\ \le\ 2\;+\;C\,(\log n)^L\,n^{-2/3}\;+\;C\,(\log n)^{L/3}\,n^{-1/3}
\ \le\ 2\;+\;C'(\log n)^{L'}\,n^{-1/3}.
\end{equation}
\end{lemma}

\begin{proof}
We apply ~\pref{lem:edge-uniform-close-2-Stieltjes}.

\ppart{High-probability rigidity event for GOE}
Let $\lambda_1\le\cdots\le \lambda_n$ be the eigenvalues of $A$ and let
$\gamma_j$ be the $j$th classical location under $\mu_{\mathrm{sc}}$ (i.e.\ $j/n=\mu_{\mathrm{sc}}((-\infty,\gamma_j])$).
By eigenvalue rigidity for generalized Wigner matrices~\cite{erdHos2012rigidity}, for any $\xi>1$ there exists
$L=L(\xi)$ such that with probability at least $1-C\exp(-c(\log n)^\xi)$,
\begin{equation}\label{eq:rigidity}
|\lambda_j-\gamma_j|
\ \le\
(\log n)^L\,\bigl[\min(j,n-j+1)\bigr]^{-1/3}\,n^{-2/3}
\qquad\text{for all }j=1,\dots,n.
\end{equation}
We work on this event.
Since $\gamma_n=2$, \eqref{eq:rigidity} implies
\[
\lambda_{\max}(A)=\lambda_n \le 2+(\log n)^L n^{-2/3}.
\]
Hence
\[
\sup\supp(\mu_A) \le 2+\delta_n,
\qquad \delta_n:=(\log n)^L n^{-2/3}.
\]

\ppart{Bound $W_1(\mu_A,\mu_{\mathrm{sc}})$ from rigidity}
Couple $\mu_A$ to $\mu_{\mathrm{sc}}$ by matching the atom at $\lambda_j$ with the quantile point $\gamma_j$.
This coupling yields
\[
W_1(\mu_A,\mu_{\mathrm{sc}})
\ \le\ \frac{1}{n}\sum_{j=1}^n |\lambda_j-\gamma_j|
\ \le\ \frac{(\log n)^L n^{-2/3}}{n}\sum_{j=1}^n \bigl[\min(j,n-j+1)\bigr]^{-1/3}.
\]
Using the estimate
\[
\sum_{j=1}^n \bigl[\min(j,n-j+1)\bigr]^{-1/3}
\le 2\sum_{j=1}^{\lfloor n/2\rfloor} j^{-1/3}
\le C\,n^{2/3},
\]
we get
\[
W_1(\mu_A,\mu_{\mathrm{sc}})\ \le\ C\,\frac{(\log n)^L}{n}.
\]
Set $\varepsilon_n:=C(\log n)^L/n$.
Then \pref{lem:edge-uniform-close-2-Stieltjes} gives, deterministically (conditional on $A$),
\[
E_+(t)\ \le\ 2+\delta_n+2\,\varepsilon_n^{1/3}
\qquad\text{for all }t\in[0,1].
\]
Substituting $\delta_n=(\log n)^L n^{-2/3}$ and $\varepsilon_n=C(\log n)^L/n$ yields \eqref{eq:edge-bound-GOE}.
\end{proof}

\begin{corollary}
    \label{cor:AIk+Bk-final-bound}
    Let $A \in \R^{n \times n}$ be a normalized GOE random matrix.
    For each $k \in \N$, let $B_k \in \R^{(nk) \times (nk)}$ be a normalized GOE matrix. 
    Then
    \[ \sup_{t \in [0,1]}\left[\opnorm{\sqrt{1-t}A \otimes \Id_k + \sqrt{t}B_k}\right] \;\xrightarrow[k\to\infty]{a.s.}\; \zeta_A, \]
    where $\zeta_A$ is a scalar function of $A$ such that for every $\xi > 1$, there is an $L$ and constants $c,C,C' > 0$ such that
    \[ \Pr_{A}\left\{\zeta_A \;>\; 2\,+\,C'(\log n)^{L}\,n^{-1/3}\right\} \;\le\; C\exp(-c(\log n)^{\xi}). \]
\end{corollary}
\begin{prf}
For each $k\in\mathbb{N}$, define the random function on $[0,1]$
\[
f_k(t)\;:=\;\lambda_{\max}\!\Big(\sqrt{1-t}\,A\otimes \Id_k+\sqrt{t}\,B_k\Big)
\;=\;\sup\Spec\!\Big(\sqrt{1-t}\,A\otimes \Id_k+\sqrt{t}\,B_k\Big),
\]
and define the deterministic limit function
\[
f(t)\;:=\;E_+(t)\;=\;\sup\Spec\!\big(\sqrt{1-t}\,A+\sqrt{t}\,S\big).
\]

\ppart{Almost sure equicontinuity of $(f_k)_k$}
For $t,t'\in[0,1]$, Weyl's inequality gives
\[
|f_k(t)-f_k(t')|
\;\le\;
\opnorm{
\big(\sqrt{1-t}-\sqrt{1-t'}\big)\,A\otimes \Id_k
+\big(\sqrt{t}-\sqrt{t'}\big)\,B_k}.
\]
Hence, using $\opnorm{A\otimes \Id_k}=\opnorm{A}$,
\begin{equation}\label{eq:fk-holder}
|f_k(t)-f_k(t')|
\;\le\;
\opnorm{A}\,|\sqrt{1-t}-\sqrt{1-t'}|
+\opnorm{B_k}\,|\sqrt{t}-\sqrt{t'}|
\;\le\;(\opnorm{A}+\opnorm{B_k})\sqrt{|t-t'|},
\end{equation}
where we used $|\sqrt a-\sqrt b|\le \sqrt{|a-b|}$ for $a,b\in[0,1]$.

By the Bai--Yin law~\cite{bai1993limit},
\[
\opnorm{B_k}\ \xrightarrow[k\to\infty]{a.s.}\ 2.
\]
In particular, on the almost-sure event
\[
\Omega_1:=\left\{\sup_{k\ge 1}\opnorm{B_k}<\infty\right\},
\]
the constant $L:=\opnorm{A}+\sup_{k\ge 1}\opnorm{B_k}$ is finite and \eqref{eq:fk-holder} implies that
$\{f_k\}_{k\ge 1}$ is uniformly (H\"older-$1/2$) equicontinuous on $[0,1]$.
Moreover, for each draw of $A, B_1, B_2, \dots$ falling within $\Omega_1$ the family is uniformly bounded since
\[
|f_k(t)| \;\le\; \opnorm{\sqrt{1-t}\,A\otimes \Id_k+\sqrt t\,B_k}
\;\le\; \opnorm{A}+\sup_{k\ge 1}\opnorm{B_k}.
\]

\ppart{Almost sure pointwise convergence on a dense set}
Let $\mathbb{Q}\cap[0,1]$ be the rationals in $[0,1]$. By \pref{lem:strong-conv-to-a0+s},
for each fixed $t\in[0,1]$ we have $f_k(t)\to f(t)$ almost surely.
Taking a countable intersection over $t\in\mathbb{Q}\cap[0,1]$, we obtain an event $\Omega_2$ of probability $1$
such that under $\Omega_2$,
\begin{equation}\label{eq:pointwise-rationals}
f_k(t)\ \xrightarrow[k\to\infty]{}\ f(t)\qquad\text{for every }t\in\mathbb{Q}\cap[0,1].
\end{equation}

\ppart{Arzel\`a--Ascoli and uniform convergence}
Work on the almost sure event $\Omega:=\Omega_1\cap\Omega_2$ and condition on a particular draw of $B_1, B_2, \dots$.
Then the family $\{f_k\}_{k\ge 1}$ is uniformly bounded and equicontinuous on
the compact interval $[0,1]$, hence by the Arzel\`a--Ascoli theorem, from any subsequence $(f_{k_m})_{m\ge 1}$ we can extract a
further subsequence $(f_{k_{m_\ell}})_{\ell\ge 1}$ converging uniformly to some continuous function
$g$.

But uniform convergence implies pointwise convergence, so for every rational $t\in\mathbb{Q}\cap[0,1]$ we have
\[
g(t)=\lim_{\ell\to\infty} f_{k_{m_\ell}}(t).
\]
On the other hand, \eqref{eq:pointwise-rationals} gives
\[
\lim_{\ell\to\infty} f_{k_{m_\ell}}(t)=f(t)\qquad\forall\,t\in\mathbb{Q}\cap[0,1].
\]
Hence $g(t)=f(t)$ for all rational $t\in[0,1]$.
The function $f$ is continuous since a faithful trace implies $|\sup\Spec(X)-\sup\Spec(Y)|\le \Lpnorm[\infty]{X-Y}$ and then the same calculation as in \pref{eq:fk-holder} applies.
Since $g$ was also continuous by uniform convergence,
$g$ and $g$ must coincide on all of $[0,1]$, i.e.\ $g \equiv f$.

Therefore \emph{every} subsequence of $f_k$ has a further subsequence converging uniformly to $f$.
This forces the whole sequence to converge uniformly:
\begin{equation}\label{eq:uniform-conv}
\left|\sup_{t\in[0,1]} f_k(t)-\sup_{t\in[0,1]} f(t)\right| \le \sup_{t\in[0,1]}|f_k(t)-f(t)|\ \xrightarrow[k\to\infty]{}\ 0.
\end{equation}

\ppart{Passing to the limit}
We now allow $B_1, B_2, \dots$ to be random again.

By \pref{lem:GOE-free-edge}, for every $\xi>1$ there exists $L=L(\xi)>0$ and constants $c,C,C'>0$ such that, with probability
at least $1-C\exp\!\big(-c(\log n)^\xi\big)$ in $A$,
\[ \sup_{t \in [0,1]} f(t) \;=\; \sup_{t \in [0,1]}\left[\sup\Spec\!\big(\sqrt{1-t}\,A+\sqrt{t}\,S\big)\right] \;\le\; 2 + C'(\log n)^{L}n^{-1/3}. \]
Therefore, on that event in $A$, by \pref{eq:uniform-conv} we have
\[ \sup_{t\in[0,1]}\sup\Spec\!\Big(\sqrt{1-t}\,A\otimes \Id_k+\sqrt t\,B_k\Big) \;=\; \sup_{t\in[0,1]} f_k(t) \;\xrightarrow[k\to\infty]{a.s.}\; \sup_{t\in[0,1]} f(t). \]
Due to symmetry, and applying a union bound under the addition of the same event on $-A$, the same convergence holds for 
    \[ \sup_{t \in [0,1]}\left[\opnorm{\sqrt{1-t}A \otimes \Id_k + \sqrt{t}B_k}\right] \;\xrightarrow[k\to\infty]{a.s.}\; \sup_{t \in [0,1]} f(t) \;=:\; \zeta_A \;\le\; 2 + C'(\log n)^{L}n^{-1/3}.\qedhere \]
\end{prf}

\subsection{Free interpolation}
\label{sec:free-interpolation}

\begin{proposition}[Free interpolation bound for real--axis resolvent squares in the $\gamma$--regime]
\label{prop:gamma-regime-square-expectation-estimate}
Fix parameters $a,\beta$ satisfying $0\le a \le1$ and $0\le\beta<\frac12$, and set
\[
\gamma \;:=\; (1-2\beta)/2\in(0,1).
\]
Let $D\in M_n(\R)$ be a deterministic diagonal matrix with $D\succeq \Id_n$ and let $A\in M_n(\R)$ be a normalized GOE matrix.


Let $\xi>1$ be fixed and let $L,C',c,C>0$ be the constants from \pref{cor:AIk+Bk-final-bound} (applied with this $\xi$).
Define deterministic sequences
\[
\eta_n \,:\, =C'(\log n)^{L}\,n^{-1/3},
\qquad
\delta_n \,:=\, \frac{1}{\log n},
\qquad
\theta_n \,:=\, 2+\eta_n+\delta_n.
\]
Let
\begin{equation}
\label{eq:Omega-gamma}
\Omega_\gamma \;:=\; \{\zeta_A\le 2+\eta_n\},
\end{equation}
where $\zeta_A$ is as in \pref{cor:AIk+Bk-final-bound}. Then $\Pr\Omega_\gamma^c \;\le\; C e^{-c(\log n)^\xi}$.
Assume that $n$ is large enough that $n^2\Pr\Omega_\gamma^c$ is bounded by some positive rational function of $\beta$ and $1/\gamma$ and also
\begin{equation}
\label{eq:n-large-cutoff-conditions}
\beta\bigl((\theta_n+\delta_n)-2\bigr) \;\le\; \gamma.
\end{equation}

Let $F := F_{\beta,-D}$ be the complex analytic subordination function in \pref{lem:subordination-semicircular}.
Then $F$ extends analytically to $a$, and there exists a universal constant $C_0>0$ such that
\begin{equation}
\label{eq:gamma-square-target-bound}
\schnorm{
\E_A\!\left[E_{\mathcal D_n}\bigl[D(a\Id_n-\beta A+D)^{-2}D\bigr]\,\one_{\Omega_\gamma}\right]
\,-\,
F'(a)\,D\bigl(F(a)\Id_n+D\bigr)^{-2}D
}
 \le\
\frac{C_0\,\beta^4}{n^{3/2}\gamma^6}
\;+\;
\frac{C_0\,\beta^2}{n\gamma^4}.
\end{equation}
\end{proposition}

\begin{proof}
The proof follows the same broad strategy as in \cite[Proposition 4.20]{jekel2024pha}\iffocs{}{ and \cite[Section 4]{jekel2025pha2}}, including setting up a Gaussian interpolation in \pref{eq:hk-def}, applying Gaussian integration by parts to the derivative of the interpolant in \pref{lem:gamma-square-ibp-B} and \pref{lem:gamma-square-ibp-A}, using the Gaussian Poincar\'e inequality in \pref{lem:gamma-square-poincare-main}, and taking a limit that approaches the ideal freely independent case in \pref{eq:finite-b-estimate-with-chi-infty}.

To handle the square of the resolvent, we use the identity
$\lim_{b\to 0} -(1/b)\Im[(bi + X)^{-1}] = X^{-2}$ where $X$ is real and when both sides are defined.
Because this involves taking a limit as $b\to 0$, we cannot depend on the regularizing effect of the imaginary component this time.
Therefore, the main new addition to this argument is multiplication by a cutoff function \pref{eq:chi-k-def} to ensure that the interpolant is still defined in the rare event that a large fluctuation in $\opnorm{A}$ causes the bare resolvent to hit a singularity.

Fix $k\in\N$ and write $N:=nk$ and let $B_k\in M_{N}(\R)$ be a normalized independent GOE matrix.. Define
\[
W_k(t):=\sqrt{1-t}\,(A\otimes \Id_k)+\sqrt t\,B_k,
\qquad t\in[0,1],
\qquad
D_k:=D\otimes \Id_k.
\]
For $0 \le a \le 1$, define
\[
X_k(t):=a\Id_{N}-\beta W_k(t)+D_k.
\]
For $z=a+ib\in\C$ with $b>0$, define the resolvent
\[
R_k(z,t):= \bigl(z \Id_{N}-\beta W_k(t)+D_k\bigr)^{-1}.
\]

\ppart{Cutoff function}
Define the path maximum
\[
\Xi_k:=\sup_{t\in[0,1]}\opnorm{W_k(t)},
\]
and let $\phi\in C^\infty(\R;[0,1])$ satisfy $\phi(u)=1$ for $u\le0$ and $\phi(u)=0$ for $u\ge1$.
Set the cutoff
\begin{equation}
\label{eq:chi-k-def}
\chi_k:=\phi\!\left(\frac{\Xi_k-\theta_n}{\delta_n}\right)\in[0,1].
\end{equation}

This cutoff is used to define the interpolation in \pref{eq:hk-def}.
Multiplying by the cutoff ensures that the interpolation observable is defined, bounded, and differentiable for all realizations of $A$ and $B_k$ (as recorded in \pref{lem:resolvent-bounds-cutoff-support} and \pref{lem:resolvent-lipschitz-cutoff-support}), enabling the application of Gaussian integration by parts to obtain the post-IBP derivative \pref{eq:hprime-post-ibp}.
This incurs error terms involving the derivative of the cutoff function, which will be negligible since the cutoff derivative is supported only on rare events in $A$ and $B_k$, as argued in \pref{lem:cutoff-gradient-product} and \pref{lem:gamma-square-cutoff-ibp-error}.

\begin{fact}
    \label{fact:cutoff-gradient-max}
We may take $\sup_u |\phi'(u)| < 2$ and so 
\begin{align*}
|\chi_k(A,B_k) - \chi_k(A',B_k')| &\le 2\delta_n^{-1}\left|\Xi_k(A,B_k) - \Xi_k(A',B_k')\right| \le 2\delta_n^{-1}(\opnorm{A - A'} + \opnorm{B_k - B_k'})
\\ &\le 2\delta_n^{-1}\left(\sqrt{n}\schnorm{A - A'} + \sqrt{N}\schnorm{B_k - B_k'}\right),
\end{align*}
meaning $\chi_k$ is $2\sqrt{N}/\delta_n$-Lipschitz as a function of $B_k$ and $2\sqrt{n}/\delta_n$-Lipschitz as a function of $A$ in $\schnorm{\cdot}$ distance.
\end{fact}

\begin{lemma}[Resolvent bounds on the support of the cutoff]
\label{lem:resolvent-bounds-cutoff-support}
Let $R_k := R_k(z,t)$.
Let $T \in M_N(\R)$ be any constant symmetric matrix.
    On $\{\chi_k \ne 0\}$ and when $\Re[z] =: a \in [0,1]$ and $\Im[z] =: b$ satisfies $|b| < 1$, we have $\min\Spec(|z \Id_{N}-\beta W_k(t)+D_k|) \ge \gamma$ and $\opnorm{z \Id_{N}-\beta W_k(t)} \le 3-\gamma$.

As a result, $z \mapsto R_k(z, t)$ is analytic when $\Re[z] \ge 0$ and $\chi_k \ne 0$, and the following bounds hold over $\Re[z] \ge 0$:

\begin{enumerate}[label=(\alph*), ref=\thetheorem(\alph*)]
\item 
\label{lem:free-interp-bounds-R}
\begin{equation*}
\opnorm{\chi_kR_k}\le \gamma^{-1},
\qquad
\opnorm{\frac{1}{b}\Im [\chi_kR_k]}
\le \gamma^{-2}.
\end{equation*}
\item 
\label{lem:free-interp-bounds-DRRD}
\begin{equation*}
\opnorm{\chi_kD_kR_k^2D_k}\le 9\gamma^{-2}.
\end{equation*}
\item
\label{lem:free-interp-bounds-RDTDR}
\begin{equation*}
\schnorm{\chi_k R_k D_k T D_k R_k} \le 9\gamma^{-2}\schnorm{T},
\qquad
\schnorm{\frac{1}{b}\Im[\chi_kR_kD_kTD_kR_k]}
\le 18\gamma^{-3}\schnorm{T}.
\end{equation*}
\item
\label{lem:free-interp-bounds-R^2DTDR}
\begin{equation*}
\schnorm{\chi_kR_k^2D_kTD_kR_k} \le 9\gamma^{-3}\schnorm{T},
\qquad
\schnorm{\frac{1}{b}\Im[\chi_kR_k^2D_kTD_kR_k]}
\le 27\gamma^{-4}\norm{T}_2.
\end{equation*}
\end{enumerate}
\end{lemma}
\begin{proof}
    Since $\phi(u)=0$ for $u\ge1$, on $\{\chi_k>0\}$ one has $\Xi_k\le \theta_n+\delta_n$, hence
\[
\opnorm{W_k(t)} \le \theta_n+\delta_n\qquad\forall\,t\in[0,1].
\]
Using $D_k\succeq \Id_N$ and the previous inequality yields, on $\{\chi_k>0\}$,
\[
X_k(t)\succeq \bigl(a+1-\beta(\theta_n+\delta_n)\bigr)\,\Id_N.
\]
By \eqref{eq:n-large-cutoff-conditions},
\[
a+1-\beta(\theta_n+\delta_n)\ge 1-2\beta - \gamma \ge \gamma,
\]
so on $\{\chi_k\ne0\}$,
\begin{equation}
\label{eq:X-gap}
X_k(t)\succeq \gamma \Id_N\qquad\forall\,t\in[0,1].
\end{equation}
Consequently, due to \pref{lem:resolvent-bounds}, the bounds hold on $\{\chi_k\ne0\}$.
\end{proof}

\begin{lemma}[Resolvent Lipschitzness on the support of the cutoff]
\label{lem:resolvent-lipschitz-cutoff-support}
Let $R_k := R_k(z,t)$.
Let $T \in M_N(\R)$ be any constant symmetric matrix.

    On $\{\chi_k \ne 0\}$ and $\Re[z] = a \in [0,1]$, we have the following Lipschitz constants with respect to changes in $B$, where the underlying matrix norms for the Lipschitz constants are the normalized Schatten-2 norm $\|\cdot\|_2$:
\begin{equation*}
\lipnorm{B_k \mapsto R_k} \le \beta\gamma^{-2}\sqrt{t},
\qquad
\lipnorm{B_k \mapsto \frac{1}{b}\Im[R_k]}
\le 2\beta\gamma^{-3}\sqrt{t}.
\end{equation*}
\begin{equation*}
\lipnorm{B_k \mapsto R_kD_kTD_kR_k} \le 18\beta\gamma^{-3}\sqrt{t}\opnorm{T},
\qquad
\lipnorm{B_k \mapsto \frac{1}{b}\Im[R_kD_kTD_kR_k]}
\le 54\beta\gamma^{-4}\sqrt{t}\opnorm{T}.
\end{equation*}
\end{lemma}
\begin{prf}
Invoke \pref{lem:resolvent-frechet-dR}, \pref{lem:resolvent-frechet-dImR}, \pref{lem:resolvent-frechet-dRjDTDRkj}, and \pref{lem:resolvent-frechet-dImRDTDR} with $dW \gets 0$ and $dM \gets -\beta\sqrt{t} dB_k$ and \pref{fact:frechet-lipschitz} and use the bounds from \pref{lem:resolvent-bounds-cutoff-support}.
\end{prf}

\ppart{Step 2: the observable (with $\chi_k^2$) and endpoints}

Fix a diagonal test matrix $T\in\mathcal D_n$ and set $T_k:=T\otimes \Id_k$ and
\[
K_k:=D_kT_kD_k=(DTD)\otimes \Id_k.
\]
Define the scaled imaginary part
\[
H_k(z,t):=
-\frac{1}{b}\Im\bigl[D_k\,R_k(z,t)\,D_k\bigr]
=
D_k\,(X_k(t)^2+b^2)^{-1}\,D_k,
\]
and define the scalar interpolation observable
\begin{equation}
\label{eq:hk-def}
h_k(t):=\E\,\tr_{N}\!\bigl[\chi_k^2\,H_k(z,t)\,T_k\bigr]
=
-\frac{1}{b}\Im\,\E\,\tr_{N}\!\bigl[\chi_k^2\,D_k\,R_k(z,t)\,D_k\,T_k\bigr].
\end{equation}

At $t=0$, $W_k(0)=A\otimes \Id_k$ and $R_k(z,0)=R_n(z)\otimes \Id_k$ with $R_n(z):=(z\Id_n-\beta A+D)^{-1}$, hence
\begin{equation}
\label{eq:h0-identification}
h_k(0)=\E\Bigl[\chi_k^2\;\tr_n\!\left(-\frac{1}{b}\Im\bigl[D\,R_n(z)\,D\bigr]\,T\right)\Bigr]
=
\E\Bigl[\chi_k^2\;\tr_n\!\left(E_{\mathcal D_n}\!\left[-\frac{1}{b}\Im\bigl[D\,R_n(z)\,D\bigr]\right]\,T\right)\Bigr].
\end{equation}

\ppart{Differentiating the interpolant}

\begin{lemma}[Resolvent derivative along the interpolation]
\label{lem:gamma-square-resolvent-derivative}
For $t\in(0,1)$, on $\{\chi_k \ne 0\}$,
\[
\frac{d}{dt}R_k(z,t)
=
\frac{\beta}{2}t^{-1/2}\,R_k(z,t)\,B_k\,R_k(z,t)
-
\frac{\beta}{2}(1-t)^{-1/2}\,R_k(z,t)\,(A\otimes \Id_k)\,R_k(z,t).
\]
\end{lemma}

\begin{proof}
Differentiate $W_k(t)=\sqrt{1-t}(A\otimes \Id_k)+\sqrt t\,B_k$ to obtain
\[
W_k'(t)= -\frac{1}{2}(1-t)^{-1/2}(A\otimes \Id_k)+\frac{1}{2}t^{-1/2}B_k.
\]
Since $R_k(z,t)=(z\Id_N-\beta W_k(t)+D_k)^{-1}$, the identity $R' = R\,(\beta W')\,R$ yields the claim.
\end{proof}

Differentiating \eqref{eq:hk-def} and using \pref{lem:gamma-square-resolvent-derivative} gives, with $R:=R_k(z,t)$,
\begin{equation}
\label{eq:hprime-pre-ibp}
\begin{aligned}
h_k'(t)
&=
-\frac{1}{b}\Im\left[\E\,\tr_{N}\!\bigl[\chi_k^2\,D_k\,R_k'(z,t)\,D_k\,T_k\bigr]\right]
\\
&=
-\frac{\beta}{2}t^{-1/2}\,\frac{1}{b}\Im\left[\E\,\tr_{N}\!\bigl[\chi_k^2\,B_k\,R\,K_k\,R\bigr]\right]
\;+\;
\frac{\beta}{2}(1-t)^{-1/2}\,\frac{1}{b}\Im\left[\E\,\tr_{N}\!\bigl[\chi_k^2\,(A\otimes \Id_k)\,R\,K_k\,R\bigr]\right].
\end{aligned}
\end{equation}

\ppart{Gaussian integration by parts (IBP)}

\begin{lemma}[IBP for the $B_k$--term]
\label{lem:gamma-square-ibp-B}
For $t\in(0,1)$,
\begin{equation}
\label{eq:gamma-square-ibp-B}
\frac{\beta}{2}t^{-1/2}\,\E\,\tr_{N}\!\bigl[\chi_k^2\,B_k\,R\,K_k\,R\bigr]
=
\beta^2\,\E\Bigl[
\tr_{N}[\chi_k R]\;\tr_{N}[\chi_k\,R\,K_k\,R]
+\frac{1}{N}\tr_{N}\!\bigl[\chi_k^2\,R^3\,K_k\bigr]
\Bigr]
+\mathrm{Err}_{B,\chi}(t),
\end{equation}
where
\begin{equation}
\label{eq:ErrBchi-def}
\mathrm{Err}_{B,\chi}(t)
:=
\frac{\beta}{2}t^{-1/2}\cdot \frac{1}{N^2}\sum_{i,j=1}^{N}
\E\left[
\frac{\partial (\chi_k^2)}{\partial (B_k)_{ij}}\;
\bigl[R\,K_k\,R\bigr]_{ji}
\right].
\end{equation}
\end{lemma}

\begin{proof}
Write
\[
\tr_N\!\bigl[\chi_k^2\,B_k\,R\,K_k\,R\bigr]
=\frac{1}{N}\sum_{i,j=1}^{N} (B_k)_{ij}\,\bigl[\chi_k^2\,R\,K_k\,R\bigr]_{ji},
\]
and apply Gaussian integration by parts entrywise. For the derivative of $R$ with respect to $(B_k)_{ij}$, note that
$W_k(t)$ depends on $B_k$ through the symmetric term $\sqrt t\,B_k$, hence (in matrix form)
\[
\frac{\partial}{\partial (B_k)_{ij}}\bigl(\beta W_k(t)\bigr)=\beta\sqrt t\,(E_{ij}+E_{ji}),
\]
so
\[
\frac{\partial R}{\partial (B_k)_{ij}}=R\bigl(\beta\sqrt t\,(E_{ij}+E_{ji})\bigr)R.
\]
Using the product rule
\[
\partial_{ij}\!\bigl(\chi_k^2\,R\,K_k\,R\bigr)
=(\partial_{ij}\chi_k^2)\,RKR+\chi_k^2\,\partial_{ij}(RKR),
\]
the term with $\partial_{ij}\chi_k^2$ gives \eqref{eq:ErrBchi-def}. The term with $\chi_k^2\,\partial_{ij}(RKR)$ yields two Wick/tensor contractions, giving $\chi_k^2\,\tr_N[R]\tr_N[RKR]$ and $\frac{\chi_k^2}{N}\tr_N[R^3K]$.
Finally, $\chi_k^2\,\tr_N[R]\tr_N[RKR]=\tr_N[\chi_kR]\tr_N[\chi_kRKR]$,
and similarly for the cubic term, yielding \eqref{eq:gamma-square-ibp-B}.
\end{proof}

\begin{lemma}[IBP for the $A\otimes \Id_k$--term]
\label{lem:gamma-square-ibp-A}
For $t\in(0,1)$,
\begin{equation}
\label{eq:gamma-square-ibp-A}
\begin{aligned}
&\frac{\beta}{2}(1-t)^{-1/2}\,\E\,\tr_{N}\!\bigl[\chi_k^2\,(A\otimes \Id_k)\,R\,K_k\,R\bigr]
\\&\qquad=
\beta^2\,\E\,\tr_k\!\Bigl[(\tr_n\otimes\id)[\chi_k R]\;(\tr_n\otimes\id)[\chi_k\,R\,K_k\,R]\Bigr]
+\frac{\beta^2}{n}\,\E\,\tr_{N}\!\bigl[(\mathsf T\otimes\id)[\chi_k R]\;\chi_k\,R\,K_k\,R\bigr]
+\mathrm{Err}_{A,\chi}(t),
\end{aligned}
\end{equation}
where
\begin{equation}
\label{eq:ErrAchi-def}
\mathrm{Err}_{A,\chi}(t)
:=
\frac{\beta}{2}(1-t)^{-1/2}\cdot\frac{1}{n^2}\sum_{i,j=1}^{n}
\E\left[
\frac{\partial (\chi_k^2)}{\partial A_{ij}}\;
\tr_{N}\!\bigl[(E_{ij}\otimes \Id_k)\,R\,K_k\,R\bigr]
\right].
\end{equation}
\end{lemma}

\begin{proof}
This is the same contraction computation as in the standard partial trace/transpose IBP identity.
The dependence on $A$ is through $\sqrt{1-t}(A\otimes \Id_k)$, hence
\[
\frac{\partial}{\partial A_{ij}}\bigl(\beta W_k(t)\bigr)=\beta\sqrt{1-t}\,\bigl((E_{ij}+E_{ji})\otimes \Id_k\bigr),
\qquad
\frac{\partial R}{\partial A_{ij}}=R\bigl(\beta\sqrt{1-t}\,((E_{ij}+E_{ji})\otimes \Id_k)\bigr)R.
\]
Applying Gaussian integration by parts for the entries of $A$ and expanding with the product rule as in the proof of
\pref{lem:gamma-square-ibp-B} yields
\eqref{eq:gamma-square-ibp-A}; the terms where the derivative lands on $\chi_k^2$ collect into \eqref{eq:ErrAchi-def}.
\end{proof}

\ppart{Step 5: decomposition of $h_k'(t)$}

Substituting \eqref{eq:gamma-square-ibp-B} and \eqref{eq:gamma-square-ibp-A} into \eqref{eq:hprime-pre-ibp} gives
\begin{equation}
\label{eq:hprime-post-ibp}
\begin{aligned}
h_k'(t)
&=
-\beta^2\,\frac{1}{b}\Im\,\E\Bigl[
\tr_{N}[\chi_k R]\;\tr_{N}[\chi_k\,R\,K_k\,R]
-
\tr_k\!\Bigl[(\tr_n\otimes\id)[\chi_k R]\;(\tr_n\otimes\id)[\chi_k\,R\,K_k\,R]\Bigr]
\Bigr]
\\
&\quad
-\beta^2\,\frac{1}{b}\Im\Bigl[\frac{1}{N}\E\,\tr_{N}[\chi_k^2\,R^3K_k]\Bigr]
+\frac{\beta^2}{n}\,\frac{1}{b}\Im\Bigl[\E\,\tr_{N}\bigl[(\mathsf T\otimes\id)[\chi_k R]\;\chi_k\,R\,K_k\,R\bigr]\Bigr]
+\mathrm{Err}_{\chi}(t),
\end{aligned}
\end{equation}
where
\[
\mathrm{Err}_{\chi}(t):=
-\frac{1}{b}\Im\bigl[\mathrm{Err}_{B,\chi}(t)\bigr]
+\frac{1}{b}\Im\bigl[\mathrm{Err}_{A,\chi}(t)\bigr].
\]

\ppart{Bounds for the two single-trace terms}

\begin{lemma}[Bounds for the $N^{-1}$ and $n^{-1}$ terms]
\label{lem:gamma-square-correction-terms}
For every $k\in\N$ and every $t\in[0,1]$,
\[
\left|\frac{1}{b}\Im\left[\frac{1}{N}\E\,\tr_{N}\bigl[\chi_k^2\,R^3K_k\bigr]\right]\right|
\ \le\ \frac{27}{N}\cdot\frac{\schnorm{T}}{\gamma^4} + O(b^2)
\]
and
\[
\left|\frac{1}{b}\Im\left[\frac{1}{n}\E\,\tr_{N}\bigl[(\mathsf T\otimes\id)[\chi_k R]\;\chi_k\,R\,K_k\,R\bigr]\right]\right|
\ \le\ \frac{27}{n}\cdot\frac{\schnorm{T}}{\gamma^4},
\]
where $N:=nk$, $z=a+ib$ with $b>0$, $K_k:=D_kT_kD_k$ and $C>0$ is universal.
\end{lemma}

\begin{proof}
For the $N^{-1}$ term, we apply \pref{lem:free-interp-bounds-R^2DTDR} to find $\schnorm{b^{-1}\Im[\chi_k^2R^2K_kR]} \le 27\gamma^{-4}\schnorm{T}$.
Then $|\tr_n[\cdot]\,| \le \schnorm{\cdot}$ and trace cyclicity yield the first bound.

For the $n^{-1}$ term, let
\[
U:=(\mathsf T\otimes\id)[\chi_kR],
\qquad
V:=\chi_k\,R\,K_k\,R.
\]
Using $\Im(UV)=(\Re U)(\Im V)+(\Im U)(\Re V)$ and dividing by $b$ gives
\[
\frac{1}{b}\Im\,\tr_N[UV]
=
\tr_N\!\Bigl[(\Re U)\,\frac{1}{b}\Im V\Bigr]
+
\tr_N\!\Bigl[\Bigl(\frac{1}{b}\Im U\Bigr)\,\Re V\Bigr].
\]
Since $\mathsf T\otimes\id$ is an isometry for $\schnorm{\cdot}$ and $\Re$ is a contraction, we apply \pref{lem:free-interp-bounds-R} and \pref{lem:free-interp-bounds-RDTDR} to find
\begin{align*}
&\schnorm{\Re U} \le \schnorm{\chi_kR} \le \opnorm{\chi_kR}\le \gamma^{-1},
\quad
&&\schnorm{\tfrac{1}{b}\Im U}\le \opnorm{\tfrac{1}{b}\Im[\chi_k R]}\le \gamma^{-2},
\\
&\schnorm{\Re V} \le \schnorm{\chi_kRD_kT_kD_kR} \le \frac{9}{\gamma^2}\schnorm{T},
\;\;
&&\schnorm{\tfrac{1}{b}\Im V} = \schnorm{\tfrac{1}{b}\Im[\chi_kRD_kT_kD_kR]} \le \frac{18}{\gamma^3}\schnorm{T}. 
\end{align*}
Combining the displayed estimates yields
\[
\left|\frac{1}{b}\Im\,\tr_N[UV]\right|
\le
\frac{1}{\gamma}\cdot\frac{18}{\gamma^3}\schnorm{T}
+
\frac{1}{\gamma^2}\cdot\frac{9}{\gamma^2}\schnorm{T}
\le \frac{27}{\gamma^4}\schnorm{T}.
\]
Multiplying by $1/n$ and taking expectation gives the second inequality.
\end{proof}

\ppart{Poincar\'e estimate for the covariance term}
Next we apply the Gaussian Poincar\'e inequality to bound the remaining terms.
First, we give a version of the inequality for vector-valued functions of GOE matrices.

\newcommand{\Sym}{\operatorname{Sym}}
\begin{corollary}[Vector Poincar\'e inequality for GOE]
\label{cor:vector-poincare-goe-gradient}
Let $N\in\N$, let $\Sym_N$ be the set of symmetric matrices in $M_N(\R)$ equipped with the metric induced by $\schnorm{\cdot}$, and let $B\in \Sym_N$ be a normalized GOE matrix.
Let $(\calH,\langle\cdot,\cdot\rangle_\calH)$ be a finite-dimensional real Hilbert space.
Let $f,g: \Sym_N \to \calH$ be locally Lipschitz with respect to $\schnorm{\cdot}$ on $\Sym_N$.
Then
\begin{equation}
\label{eq:vector-poincare-goe-gradient-correct}
\E\Bigl|\bigl\langle f(B)-\E f(B),\,g(B)-\E g(B)\bigr\rangle_{\calH}\Bigr|
\ \le\
\frac{2}{N^2}\,
\Bigl(\E\|df(B)\|_{\mathrm{HS}(\Sym_N,\calH)}^2\Bigr)^{1/2}\,
\Bigl(\E\|dg(B)\|_{\mathrm{HS}(\Sym_N,\calH)}^2\Bigr)^{1/2},
\end{equation}
where $\|df(B)\|_{\mathrm{HS}(\Sym_N,\calH)}$ is as given in \pref{def:frechet}, so that
\begin{equation}
\label{eq:intrinsic-gradient-norm-coordinates}
\|df(B)\|_{\mathrm{HS}(\Sym_N,\calH)}^2
=
\sum_{i=1}^N
\left\|
\frac{\partial f(B)}{\partial(B_{ii}/\sqrt N)}
\right\|_{\calH}^2
+
\sum_{1\le i<j\le N}
\left\|
\frac{\partial f(B)}{\partial(\sqrt2\,B_{ij}/\sqrt N)}
\right\|_{\calH}^2.
\end{equation}
\end{corollary}
\begin{proof}
Let $X:=f(B)-\E f(B)$ and $Y:=g(B)-\E g(B)$.
By Cauchy--Schwarz,
\begin{equation}
\label{eq:vector-poincare-cs}
    \E|\langle X,Y\rangle_{\calH}|
\le
\bigl(\E\|X\|_{\calH}^2\bigr)^{1/2}\bigl(\E\|Y\|_{\calH}^2\bigr)^{1/2}.
\end{equation}
Fix an orthonormal basis $(w_k)_{k=1}^{\dim \calH}$ of $\calH$ and define scalar maps
$f_k:=\langle f,w_k\rangle_{\calH}$ and $g_k:=\langle g,w_k\rangle_{\calH}$.
Then
\begin{equation}
\label{eq:vector-poincare-variance-decomp}
\E\|X\|_{\calH}^2=\sum_{k=1}^{\dim \calH} \Var\bigl(f_k(B)\bigr),
\qquad
\E\|Y\|_{\calH}^2=\sum_{k=1}^{\dim \calH} \Var\bigl(g_k(B)\bigr).
\end{equation}
Let
\[
e_{ii}:=\sqrt N\,E_{ii},
\qquad
e_{ij}:=\sqrt{\frac N2}\,(E_{ij}+E_{ji}),
\quad 1\le i<j\le N,
\]
where $E_{ij}$ is the matrix whose $(i,j)$th entry is 1 and all other entries are 0.
Then \(\{e_{ii}\}_{i=1}^N\cup\{e_{ij}\}_{1\le i<j\le N}\) is an orthonormal basis of
\((\Sym_N,\schnorm{\cdot})\), which is a Euclidean space of dimension \(N(N+1)/2\), and
\[
\langle B,e_{ii}\rangle_{L^2}=\frac{B_{ii}}{\sqrt N},
\qquad
\langle B,e_{ij}\rangle_{L^2}=\frac{\sqrt2\,B_{ij}}{\sqrt N}
\quad (i<j).
\]
Since $B$ is a normalized GOE matrix with $\Var(B_{ij}) = 1 + \delta_{i,j}$, these coordinates are independent centered Gaussians with variance \(2/N^2\).

So by the Gaussian Poincar\'e inequality \pref{thm:gaussian-poincare} applied on each scalar component $f_k$,
\[
\Var\bigl(f_k(B)\bigr)
\le
\frac{2}{N^2}\,
\E\left[\norm{df_k(B)}_{\mathrm{HS}(\Sym_N,\R)}^2\right].
\]
Combine with \pref{eq:vector-poincare-cs} and \pref{eq:vector-poincare-variance-decomp} and repeat the same computation for $g$ to complete the proof.
\end{proof}

Next, it is convenient to argue that the cutoff function has a negligible contribution to the Poincar\'e bound.

\begin{lemma}[Cutoff gradient second moment bound]
\label{lem:cutoff-gradient-product}
Let $N\in\N$, let $\Sym_N$ be the set of real symmetric $N \times N$ matrices, and let $(\Sym_N,\schnorm{\cdot})$ be the normalized Schatten-2 norm metric on $\Sym_N$. Let $\calH$ be a real Hilbert space.
Let $f: \Sym_N\to\calH$ be $L_f$-Lipschitz on $S$ with respect to the stated metric and let $\chi:\Sym_N\to\R$ be globally $L_{\chi}$-Lipschitz.
Fix open measurable sets $\mathcal L\subseteq S\subseteq \bbR^N$ and assume:
\begin{enumerate}
\item $0\le \chi\le 1$ everywhere and $\chi=0$ on $S^{c}$,
\item $d\chi=0$ a.e.\ on $\mathcal L^{c}$ (equivalently, $\supp(d\chi)\subseteq\mathcal L$),
\item $\sup_{x\in S}\|f(x)\|_{\calH}\le J$ for some $J<\infty$.
\end{enumerate}
Then $h:=\chi f$ is locally Lipschitz, and for every random variable $X$ with values in $M_N(\R)$,
\begin{equation}
\label{eq:cutoff-gradient-product-bound-L}
\E\bigl[\|d(\chi f)(X)\|_{\mathrm{HS}(\Sym_N,\calH)}^2\bigr]
\ \le\
2(\dim \calH)L_f^2 \;+\; 2\,J^2\,L_{\chi}^2\,\Pr\{X\in\mathcal L\},
\end{equation}
where $\norm{\cdot}_{\mathrm{HS}(\Sym_N,\calH)}$ is as in \pref{eq:intrinsic-gradient-norm-coordinates}.
\end{lemma}

\begin{proof}
By Rademacher's theorem, $f$ and $\chi$ are differentiable a.e.
On the set of differentiability, for each coordinate $\ell\in\{1,\dots,N\}$ the product rule gives
\[
\partial_{\ell}(\chi f)= (\partial_{\ell}\chi)\,f+\chi\,(\partial_{\ell}f).
\]
Taking squared norms and using $(a+b)^2\le 2a^2+2b^2$ yields
\[
\|d(\chi f)\|_{\mathrm{HS}(\Sym_N,\calH)}^2
\le
2\,\chi^2\,\|df\|_{\mathrm{HS}(\Sym_N,\calH)}^2
+
2\,\|f\|_{\calH}^2\,\|d\chi\|_{\mathrm{HS}(\Sym_N,\R)}^2.
\]
Since $0\le\chi\le 1$ and $\chi=0$ on $S^c$, the first term is bounded by
\[2\cdot\one_S\,\|df\|_{\mathrm{HS}(\Sym_N,\calH)}^2 \;\le\; 2\cdot\one_S\,(\dim \calH)\frenorm{df}^2 \;\le\; 2(\dim \calH)L_f^2.\]
Since $d\chi=0$ a.e.\ on $\mathcal L^c$ and $\|f\|_{\calH}\le J$ on $S\supseteq\mathcal L$,
the second term is bounded by $2J^2L_{\chi}^2\,\one_{\mathcal L}$.
Taking expectations yields \eqref{eq:cutoff-gradient-product-bound-L}.
\end{proof}

Now define $M_k(\C)$--valued maps (as functions of $B_k$, with $A$ held fixed)
\[
F(B_k):=(\tr_n\otimes\id)[\chi_k\,R_k(z,t)^*],
\qquad
G(B_k):=(\tr_n\otimes\id)[\chi_k\,R_k(z,t)\,K_k\,R_k(z,t)].
\]
Then $\tr_k[F(B_k)^*]=\tr_N[\chi_k R]$, $\tr_k[G(B_k)]=\tr_N[\chi_k R K_k R]$, and
\[
\tr_k[F(B_k)^*G(B_k)]
=
\tr_k\!\Bigl[(\tr_n\otimes\id)[\chi_k R]\;(\tr_n\otimes\id)[\chi_k R K_k R]\Bigr].
\]

\begin{lemma}[Main term bound via Poincar\'e inequality with cutoff]
\label{lem:gamma-square-poincare-main}
Let $\mathcal L_{k,n}$ be the event that $\chi_k \in (0,1)$, equivalently $\theta_n < \Xi_k \le \theta_n + \delta_n$.
Then
\begin{equation}
\label{eq:poincare-main-bound-concrete}
\left|
\E\left[
\frac{1}{b}\Im\Bigl(\tr_k[F(B_k)^*]\;\tr_k[G(B_k)]-\tr_k[F(B_k)^*G(B_k)]\Bigr)
\right]
\right|
\ \le\
\frac{C\,\beta^2}{n^{3/2}\gamma^6}\,\schnorm{T}
\;+\;
\frac{C}{n^{1/2}\delta_n^2\gamma^4k}\,\Pr\mathcal L_{k,n}\,\schnorm{T},
\end{equation}
with a universal constant $C>0$.
\end{lemma}

\begin{proof}
Write $F^*=F_R+iF_I$ and $G=G_R+iG_I$. The algebraic identity
\begin{equation}
\label{eq:Im-covariance-splitting-again}
\frac{1}{b}\Im\Bigl(\tr_k[F^*]\;\tr_k[G]-\tr_k[F^*G]\Bigr)
=
\Bigl(\tr_k[F_R]\;\tr_k\!\bigl[\tfrac{1}{b}G_I\bigr]-\tr_k\!\bigl[F_R\,\tfrac{1}{b}G_I\bigr]\Bigr)
+
\Bigl(\tr_k\!\bigl[\tfrac{1}{b}F_I\bigr]\;\tr_k[G_R]-\tr_k\!\bigl[\tfrac{1}{b}F_I\,G_R\bigr]\Bigr)
\end{equation}
reduces the claim to bounding the expectations of the two displayed differences.

\ppart{Conditioning and orthogonal invariance}
Condition on $A$ and write $\E_B$ for expectation over $B_k$.
For any $O\in O_k$, conjugation by $\Id_n\otimes O$ preserves the law of $B_k$ and leaves
$A\otimes \Id_k$, $D_k$, $T_k$, and $\chi_k$ invariant.
Since $(\tr_n\otimes\id)$ is equivariant under this conjugation, both $F(B_k)$ and $G(B_k)$ satisfy
\[
F\bigl((\Id_n\otimes O)B_k(\Id_n\otimes O)^{\sT}\bigr)=O\,F(B_k)\,O^{\sT},
\qquad
G\bigl((\Id_n\otimes O)B_k(\Id_n\otimes O)^{\sT}\bigr)=O\,G(B_k)\,O^{\sT}.
\]
Therefore $\E_B[F]$ and $\E_B[G]$ commute with all $O\in O_k$, hence are scalar multiples of $\Id_k$.
In particular, $\E_B[F_R]$ and $\E_B[(1/b)G_I]$ are scalar multiples of $\Id_k$.

\ppart{Rewrite as inner product of centered $k\times k$ matrices}
For the first bracket, set $U:=F_R$ and $V:=\tfrac1b G_I$, and define $\widetilde V:=\tr_k(V)\Id_k-V$.
Since $\E_B[V]$ is a scalar multiple of $\Id_k$, one has $\E_B[\widetilde V]=0$ and
\[
\tr_k(U)\tr_k(V)-\tr_k(UV)=\tr_k(U\widetilde V)=\langle U,\widetilde V\rangle_{\tr_k}.
\]
Hence
\[
\E_B\bigl[\tr_k(U)\tr_k(V)-\tr_k(UV)\bigr]
=
\E_B\left\langle U-\E_B U,\ \widetilde V\right\rangle_{\tr_k}.
\]
The second bracket is treated identically with $(U,V)=(\tfrac1b F_I,\,G_R)$.

\ppart{Apply Gaussian Poincar\'e inequality}
Work in the Hilbert space $W:=\mathrm{Sym}_k(\R)$ with inner product $\langle X,Y\rangle_W:=\tr_k(XY)$ implying $\norm{\cdot}_W = \schnorm{\cdot}$ and the norm on gradients $\norm{\cdot}_{\mathrm{HS}(\Sym_N,\calH)}$ given by \pref{eq:intrinsic-gradient-norm-coordinates}. 
Henceforth, we shorten $\norm{\cdot}_{\mathrm{HS}(\Sym_N,W)}$ to just $\norm{\cdot}_{\mathrm{HS}}$.
Apply \pref{cor:vector-poincare-goe-gradient} to the $W$--valued maps $B_k\mapsto U(B_k)$ and $B_k\mapsto \widetilde V(B_k)$.
Using $\|d\widetilde V\|_{\mathrm{HS}}\le 2\|dV\|_{\mathrm{HS}}$ (pointwise, because $\widetilde V=\tr_k(V)\Id_k-V$ and $\tr_k$ is contractive in $\schnorm{\cdot}$) gives
\[
\Bigl|\E_B\bigl[\tr_k(U)\tr_k(V)-\tr_k(UV)\bigr]\Bigr|
\ \le\
\frac{4}{N^2}\,
\Bigl(\E_B\|dU\|_{\mathrm{HS}}^2\Bigr)^{1/2}\,
\Bigl(\E_B\|dV\|_{\mathrm{HS}}^2\Bigr)^{1/2}.
\]
The same estimate holds for the second bracket in \eqref{eq:Im-covariance-splitting-again}.
Averaging over $A$ yields
\begin{equation}
\label{eq:poincare-main-bound-cutoff-product}
\begin{aligned}
&\left|
\E\left[
\frac{1}{b}\Im\Bigl(\tr_k[F(B_k)^*]\;\tr_k[G(B_k)]-\tr_k[F(B_k)^*G(B_k)]\Bigr)
\right]
\right|
\
\\&\qquad\qquad\qquad\le
\frac{4}{N^2}\,\left(
\Bigl(\E\|dF_R\|_{\mathrm{HS}}^2\Bigr)^{1/2}\,
\Bigl(\E\bigl\|d\bigl(\tfrac{1}{b}G_I\bigr)\bigr\|_{\mathrm{HS}}^2\Bigr)^{1/2}+\Bigl(\E\bigl\|d\bigl(\tfrac{1}{b}F_I\bigr)\bigr\|_{\mathrm{HS}}^2\Bigr)^{1/2}\,
\Bigl(\E\|dG_R\|_{\mathrm{HS}}^2\Bigr)^{1/2}\right).
\end{aligned}
\end{equation}

\ppart{Bound the gradient factors using \pref{lem:cutoff-gradient-product}}
Each of $F_R$, $\tfrac1b F_I$, $G_R$, $\tfrac1b G_I$ is of the form $\chi_k$ times a $\chi_k$--free base map.
Define the base maps
\[
f_R(B_k):=\Re\bigl[(\tr_n\otimes\id)R(B_k)^*\bigr],\quad
f_I(B_k):=\frac{1}{b}\Im\bigl[(\tr_n\otimes\id)R(B_k)^*\bigr],
\]
\[
g_R(B_k):=\Re\bigl[(\tr_n\otimes\id)\bigl(R(B_k)K_kR(B_k)\bigr)\bigr],\quad
g_I(B_k):=\frac{1}{b}\Im\bigl[(\tr_n\otimes\id)\bigl(R(B_k)K_kR(B_k)\bigr)\bigr],
\]
so that $F_R=\chi_k f_R$, $\tfrac1bF_I=\chi_k f_I$, $G_R=\chi_k g_R$, $\tfrac1bG_I=\chi_k g_I$.

\pref{lem:cutoff-gradient-product} applies with the sets $S:=\{\chi_k>0\}$ and $\mathcal L:=\mathcal L_{k,n}$
because $\chi_k=0$ on $S^c$ and $\supp(\nabla\chi_k)\subseteq \mathcal L_{k,n}$ by construction.

Since $\tr_n\otimes\id$ and $\Re$ are contractions in $\schnorm{\cdot}$ and $\schnorm{\cdot} \le \opnorm{\cdot}$, we apply \pref{lem:free-interp-bounds-R}, \pref{lem:free-interp-bounds-RDTDR}, and \pref{lem:resolvent-lipschitz-cutoff-support} with $\sqrt{t} \le 1$ to find
\begin{equation}
\label{eq:base-map-LJ-bounds}
\begin{aligned}
&\|f_R\|_{\mathrm{Lip}(S)}\le \frac{\beta}{\gamma^2},\qquad &&\sup_S\|f_R\|_W\le \frac{1}{\gamma},\\
&\|f_I\|_{\mathrm{Lip}(S)}\le \frac{2\beta}{\gamma^3},\qquad &&\sup_S\|f_I\|_W\le \frac{1}{\gamma^2},\\
&\|g_R\|_{\mathrm{Lip}(S)}\le \frac{18\beta}{\gamma^3}\opnorm{T},\qquad &&\sup_S\|g_R\|_W\le \frac{9}{\gamma^2}\opnorm{T},\\
&\|g_I\|_{\mathrm{Lip}(S)}\le \frac{54\beta}{\gamma^4}\opnorm{T},\qquad &&\sup_S\|g_I\|_W\le \frac{18}{\gamma^3}\opnorm{T}.
\end{aligned}
\end{equation}

\pref{lem:cutoff-gradient-product} applied to $F_R=\chi_k f_R$ and \pref{fact:cutoff-gradient-max} then give
\[
\E\|dF_R\|_{\mathrm{HS}}^2
\le
2k^2\|f_R\|_{\mathrm{Lip}(S)}^2
+
8\left(\sup_S\|f_R\|_W\right)^2\,N\delta_n^{-2}\,\Pr\mathcal L_{k,n},
\]
and similarly for $\tfrac1bF_I$, $G_R$, and $\tfrac1bG_I$.
Letting
\[ \zeta \;:=\; C\left(\frac{\beta^2 k^2}{\gamma^4}+\frac{N}{\delta_n^2\gamma^2}\,\Pr\mathcal L_{k,n}\right)\]
for some large enough $C$, we have
\[\E\|dF_R\|_{\mathrm{HS}}^2 \le \zeta
, \qquad \E\bigl\|d\left(\tfrac1bG_I\right)\bigr\|_{\mathrm{HS}}^2 \le \frac{\opnorm{T}^2\zeta}{\gamma^4}
, \qquad \E\bigl\|d\left(\tfrac1bF_I\right)\bigr\|_{\mathrm{HS}}^2 \le \frac{\zeta}{\gamma^2}
, \qquad \E\|dG_R\|_{\mathrm{HS}}^2 \le \frac{\opnorm{T}^2\zeta}{\gamma^4}.
\]
Applying the AM-GM inequality and substituting into \eqref{eq:poincare-main-bound-cutoff-product} gives
\[
\left|\E\left[\frac{1}{b}\Im\Bigl(\tr_k[F^*]\tr_k[G]-\tr_k[F^*G]\Bigr)\right]\right|
\le
\frac{C\opnorm{T}}{n^2}\left(\frac{\beta^2}{\gamma^6}+\frac{n}{\delta_n^2\gamma^4k}\,\Pr\mathcal L_{k,n}\right).
\]
Finally, use $\opnorm{T} \le n^{1/2}\schnorm{T}$ to obtain \eqref{eq:poincare-main-bound-concrete}.
\end{proof}

\ppart{Bounding the explicit IBP cutoff-derivative errors}

\begin{lemma}[IBP cutoff-derivative error bound]
\label{lem:gamma-square-cutoff-ibp-error}
For every $k\in\N$,
\[
\int_0^1 \E\left[\bigl|\mathrm{Err}_{\chi}(t)\bigr|\right]dt
\ \le\
C\left(\frac{\beta}{n\delta_n\gamma^3}+\frac{\beta}{N\delta_n\gamma^3}\right)\,
\Pr\{\theta_n\le \Xi_k\le \theta_n+\delta_n\}\,\schnorm{T},
\]
with a universal constant $C>0$.
\end{lemma}

\begin{proof}
Use \eqref{eq:ErrBchi-def}--\eqref{eq:ErrAchi-def}.
For the $-\frac{1}{b}\Im\bigl[\mathrm{Err}_{B,\chi}(t)\bigr]$ part of $\mathrm{Err}_{\chi}(t)$, with the gradient taken with respect to $\schnorm{\cdot}$,
\begin{align*}
\frac{1}{N^2}\sum_{i,j=1}^{N}\frac{\partial (\chi_k^2)}{\partial (B_k)_{ij}}\,(R K_k R)_{ji}
&=
\frac{1}{N^{5/2}}\sum_{i,j=1}^{N}\frac{\partial (\chi_k^2)}{\partial (B_k/\sqrt{N})_{ij}}\,(R K_k R)_{ji}
\\&=\frac{1}{N^{3/2}}\tr_N\!\bigl[(\nabla_{B_k}\chi_k^2)\,(R K_k R)\bigr]
\\&=\frac{2}{N^{3/2}}\tr_N\!\bigl[(\nabla_{B_k}\chi_k)\,(\chi_k R K_k R)\bigr],
\end{align*}
Then Cauchy--Schwarz gives
\[
\left|\frac{2}{N^{3/2}}\tr_N\!\bigl[(\nabla_{B_k}\chi_k)\,(\chi_k R K_k R)\bigr]\right|
\le
\frac{2}{N^{3/2}}\,\schnorm{\nabla_{B_k}\chi_k}\,\schnorm{\chi_k R K_k R}.
\]
On $\{\chi_k>0\}$, $\schnorm{\frac{1}{b}\Im[R K_k R]}\le 18\gamma^{-3}\schnorm{T}$ by \pref{lem:free-interp-bounds-RDTDR}.
Moreover, $\chi_k^2$ is supported in $[0,1]$ while
$\nabla\chi_k$ is supported on the cutoff layer and satisfies $\schnorm{\nabla\chi_k}\le 2\sqrt{N}\delta_n^{-1}$ there by \pref{fact:cutoff-gradient-max}, so that the error term is non-zero only in the event when $\theta_n\le \Xi_k\le \theta_n+\delta_n$.
Applying the same argument to $\frac{1}{b}\Im\bigl[\mathrm{Err}_{A,\chi}(t)\bigr]$ and
integrating the prefactors $t^{-1/2}$ and $(1-t)^{-1/2}$ over $[0,1]$ yields
the stated bound.
\end{proof}

\ppart{The $k\to\infty$ limit to the freely independent case}
Integrate \eqref{eq:hprime-post-ibp} from $t=0$ to $t=1$ and apply Lemmas
\ref{lem:gamma-square-correction-terms}, \ref{lem:gamma-square-poincare-main}, and \ref{lem:gamma-square-cutoff-ibp-error}.
This yields, for every $k$,
\begin{equation}
\label{eq:hk-difference-bound-with-cutoff}
|h_k(1)-h_k(0)|
\ \le\
\left(
\frac{C\beta^4}{n^{3/2}\gamma^6}
+\frac{C\beta^2}{n\gamma^4}
+\frac{C\beta^2}{N\gamma^4}
\right)\schnorm{T}
\;+\;
\mathrm{Tail}_{k,n}\,\schnorm{T}
\end{equation}
where
\begin{equation}
\label{eq:Tailkn-def}
\mathrm{Tail}_{k,n}
:=
C\left(\frac{\beta^2}{n^{1/2}\delta_n^2\gamma^4k}+\frac{\beta}{n\delta_n\gamma^3}+\frac{\beta}{N\delta_n\gamma^3}\right)\,
\Pr\{\theta_n\le \Xi_k\le \theta_n+\delta_n\}.
\end{equation}

By \pref{cor:AIk+Bk-final-bound}, under the coupling of $(B_k)$ one has almost surely
$\Xi_k\to \zeta_A$ as $k\to\infty$. Define the limiting cutoff
\begin{equation}
\label{eq:chi-infty-def}
\chi_\infty:=\phi\!\left(\frac{\zeta_A-\theta_n}{\delta_n}\right)\in[0,1].
\end{equation}
Then $\chi_k^2\to \chi_\infty^2$ almost surely and, by dominated convergence since $0 \le \chi_k \le 1$, also in $L^1$.

At $t=0$, the resolvent depends only on $A$, so letting $k\to\infty$ in \eqref{eq:h0-identification} gives
\[
\lim_{k\to\infty} h_k(0)
=
\E\left[\chi_\infty^2\,\tr_n\!\left(-\frac{1}{b}\Im\bigl[D(z\Id_n-\beta A+D)^{-1}D\bigr]\,T\right)\right].
\]
At $t=1$, the resolvent depends only on $B_k$ (and $D,T$), while $\chi_\infty$ depends only on $A$ and is independent of
$B_k$; \pref{thm: asymptotic freeness} and \pref{thm: subordination} then yield, with $F$ the subordination function satisfying $G_{\beta S - D}(z) = G_{-D}(F(z))$ as defined in the statement of this proposition,
\[
\lim_{k\to\infty} h_k(1)
\;=\;
\E[\chi_\infty^2]\;
\tr_n\!\left(-\frac{1}{b}\Im\bigl[D(F(z)\Id_n+D)^{-1}D\bigr]\,T\right).
\]
Letting $k\to\infty$ in \eqref{eq:hk-difference-bound-with-cutoff} and using
$\Pr\{\theta_n\le \Xi_k\le \theta_n+\delta_n\}\to \Pr\{\theta_n\le \zeta_A\le \theta_n+\delta_n\}$ yields
\begin{equation}
\label{eq:finite-b-estimate-with-chi-infty}
\begin{aligned}
&\left|
\E\Bigl[\chi_\infty^2\;\tr_n\!\left(-\frac{1}{b}\Im\bigl[D(z\Id_n-\beta A+D)^{-1}D\bigr]\,T\right)\Bigr]
-
\E[\chi_\infty^2]\;
\tr_n\!\left(-\frac{1}{b}\Im\bigl[D(F(z)\Id_n+D)^{-1}D\bigr]\,T\right)
\right|
\\&\qquad\qquad\qquad\qquad\qquad\qquad\qquad\qquad\qquad\qquad\le
\left(
\frac{C\beta^4}{n^{3/2}\gamma^6}
+\frac{C\beta^2}{n\gamma^4}
\right)\schnorm{T}
+\mathrm{Tail}_{\infty,n}\schnorm{T},
\end{aligned}
\end{equation}
where $\mathrm{Tail}_{\infty,n}$ is \eqref{eq:Tailkn-def} with $\Xi_k$ replaced by $\zeta_A$ and with the terms having $k$ or $N=nk$ in the denominator removed.

\ppart{The $b\downarrow 0$ limit}

On $\{\chi_\infty>0\}$ one has $\zeta_A\le \theta_n+\delta_n$, hence $\lambda_{\max}(A)\le \zeta_A\le \theta_n+\delta_n$.
By \eqref{eq:n-large-cutoff-conditions} this implies
\[
a\Id_n-\beta A+D\succeq \gamma \Id_n \qquad\text{on }\{\chi_\infty>0\}.
\]
Therefore,
\begin{align*}
-\frac{1}{b}\Im\bigl[D(z\Id_n-\beta A+D)^{-1}D\bigr]
\;\;\;=\;\;\;\,
&\;\;D\bigl((a\Id_n-\beta A+D)^2+b^2\bigr)^{-1}D
\\\ \xrightarrow[b\downarrow0]{}\ 
&\;\;D(a\Id_n-\beta A+D)^{-2}D
\quad\;\;\text{in operator norm on }\{\chi_\infty>0\}.
\end{align*}
Since $F$ is analytic at $a$ with $F(a)\in\R$ due to \pref{lem:subordination-semicircular-continuation} and the fact that $a\Id_n-2\beta +D\succeq \gamma \Id_n$, one has
\[
-\frac{1}{b}\Im\bigl[D\left(F(a+ib)\Id_n+D\right)^{-1}D\bigr]\;\;\xrightarrow[b\downarrow0]{}\;\; F'(a)\;D\left(F(a)\Id_n+D\right)^{-2}D.
\]
Letting $b\downarrow0$ in \eqref{eq:finite-b-estimate-with-chi-infty} gives
\begin{equation}
\label{eq:real-axis-estimate-with-chi-infty}
\begin{aligned}
&\left|
\E\Bigl[\chi_\infty^2\;\tr_n\!\bigl(D(a\Id_n-\beta A+D)^{-2}D\,T\bigr)\Bigr]
-
\E[\chi_\infty^2]\;\tr_n\!\bigl(F'(a)\,D(F(a)\Id_n+D)^{-2}D\,T\bigr)
\right|
\\&\qquad\qquad\qquad\qquad\qquad\qquad\qquad\qquad\qquad\;\;\le\;\;
\left(
\frac{C\beta^4}{n^{3/2}\gamma^6}
+\frac{C\beta^2}{n\gamma^4}
\right)\schnorm{T}
+\mathrm{Tail}_{\infty,n}\schnorm{T}.
\end{aligned}
\end{equation}

\ppart{Conditioning on $\Omega_\gamma$}

Since $\theta_n=2+\eta_n+\delta_n$ and $\phi(u)=1$ for $u\le0$, on $\Omega_\gamma=\{\zeta_A\le 2+\eta_n\}$ one has
$\zeta_A\le \theta_n-\delta_n$, hence $\chi_\infty=1$ and therefore $\chi_\infty^2=1$ on $\Omega_\gamma$.

Let
\[
\Delta(A) \;:=\; E_{\mathcal D_n}\bigl[D(a\Id_n-\beta A+D)^{-2}D\bigr] - F'(a)\,D(F(a)\Id_n+D)^{-2}D
\quad\in\;\;\mathcal D_n. \]
By norm duality on diagonal matrices and taking the supremum over $\schnorm{T}\le 1$ in \eqref{eq:real-axis-estimate-with-chi-infty},
\[
\schnorm{\E[\chi_{\infty}^2\,\Delta(A)]} \;=\; \sup\left\{|\tr_n[\Delta(A) T]|:\ T\in\mathcal D_n,\ \schnorm{T}\le 1\right\}
\;\le\;
\frac{C\beta^4}{n^{3/2}\gamma^6}
+\frac{C\beta^2}{n\gamma^4}
+\mathrm{Tail}_{\infty,n}
.
\]
Since $\chi_\infty^2=1$ on $\Omega_\gamma$, we have the decompositions
\begin{equation*}
\label{eq:chi-decompositions}
\chi_\infty^2=\one_{\Omega_\gamma}+\chi_\infty^2\,\one_{\Omega_\gamma^c},
\qquad
\Delta(A)\,\one_{\Omega_\gamma} =
\chi_\infty^2\,\Delta(A) - \chi_\infty^2\,\Delta(A)\,\one_{\Omega_\gamma^c}.
\end{equation*}
So combining the last two displayed equations with the triangle inequality,
\begin{equation}
    \label{eq:free-interp-Delta-decomp}
\schnorm{\E[\Delta(A)\,\one_{\Omega_\gamma}} \le \frac{C\beta^4}{n^{3/2}\gamma^6}
+\frac{C\beta^2}{n\gamma^4}
+\mathrm{Tail}_{\infty,n}
+\schnorm{\E\left[\chi_\infty^2\,\Delta(A)\,\one_{\Omega_\gamma^c}\right]}
.
\end{equation}
By triangle inequality and \pref{lem:free-interp-bounds-DRRD} applied at $t \gets 0$ and $z \gets a + bi$ with $b \to 0$,
\[
\schnorm{\E\left[\chi_\infty^2\,\Delta(A)\,\one_{\Omega_\gamma^c}\right]} \le 9\gamma^{-2}\Pr\Omega_{\gamma}^c +  \schnorm{\E\left[\chi_\infty^2\,\one_{\Omega_\gamma^c}\right] F'(a)\,D(F(a)\Id_n+D)^{-2}D}
\]
Now, by \pref{thm: subordination} and \pref{lem:subordination-semicircular-continuation},  $F'(a)\,D(F(a)\Id_n+D)^{-2} = \lim_{b \to 0} -b^{-1}\Im[D(F(a+bi)\Id_n + D)^{-1}D]$ has the same spectrum as $\lim_{b \to 0}-b^{-1}\Im[D((a+bi)\Id_n - \beta S + D)^{-1}D] = D(a\Id_n - \beta S + D)^{-2}D$, where $S$ is an ideal semicircular operator of variance 1.
Because $\beta < 1/2 - \gamma$ and $D \succeq 1$, $a\Id_n - \beta S + D > \gamma$ is invertible in the tracial non-commutative probability space and $D(a\Id_n - \beta S + D)^{-1} = \Id_n - (a\Id_n - \beta S)(a\Id_n - \beta S + D)^{-1}$ has spectrum bounded by $1+3\gamma^{-1}$, so $F'(a)\,D(F(a)\Id_n+D)^{-2}D$ has spectrum bounded by $F'(a)(1+3\gamma^{-1})^2$.
By \pref{lem:subordination-semicircular-derivative}, $0 \le F'(a) \le 1/(1-\beta^2) \le 4/3$ and $1 \le \gamma^{-1}/2$, so
\begin{equation}
\label{eq:free-interp-bound-on-omega-c-Delta}
    \schnorm{\E\left[\chi_\infty^2\,\Delta(A)\,\one_{\Omega_\gamma^c}\right]} \le 26\gamma^{-2}\P\Omega_{\gamma}^c
    .
\end{equation}

Combining \pref{eq:free-interp-Delta-decomp} with \pref{eq:free-interp-bound-on-omega-c-Delta} and observing that 
$\theta_n\le \Xi_k\le \theta_n+\delta_n$ implies $\Omega_{\gamma}$, 
\[
\schnorm{\E[\Delta(A)\,\one_{\Omega_\gamma}]} \le \frac{C\beta^4}{n^{3/2}\gamma^6}
+\frac{C\beta^2}{n\gamma^4}
+\left(\frac{\beta}{n\delta_n\gamma^3} + 26\gamma^{-2}\right)\Pr\Omega_{\gamma}^c
.\]
Finally, since $\Pr\Omega_\gamma^c\le Ce^{-c(\log n)^\xi}$, the term involving $\Pr\Omega_\gamma^c$ can be absorbed into the universal constant.
This yields \eqref{eq:gamma-square-target-bound}.
\end{proof}

\begin{corollary}[Control on diagonal of squared resolvent]
    \label{cor:diagonal-Q2-controlled-in-l2}
     Let $A \in M_n(\R)$ be a normalized GOE matrix.
     Let $0 \le \beta < 1/2$ and $\gamma := (1 - 2\beta)/2$.
     Let $F := F_{\beta, -D}$ be the subordination function so that $G_{\beta S - D} = G_{-D}\circ F$ as in \pref{lem:subordination-semicircular}.
     Under an event that occurs with probability at least $1 - e^{-O((\log n)^{\xi})}$ in the randomness of $A$ for some $\xi > 1$ and large enough $n$, we have uniformly for all $D \succeq \Id_n$ and all $a \in [0, 1]$,
    \begin{equation}
    \label{eq:diagonal-DQ2D-controlled-in-l2-a}
        \schnorm{E_{\calD_n}\left[D(a -\beta A + D)^{-2}D\right] - F'(a)D(F(a) + D)^{-2}D} \le \frac{C\beta}{\gamma^4\sqrt{n}},
    \end{equation}
    for some universal constant $C$.
    In particular, when we choose $a \gets \beta^2\tr_n(D^{-1})$,
    \begin{equation}
    \label{eq:diagonal-DQ2D-controlled-in-l2}
        \schnorm{E_{\calD_n}\left[D(\beta^2\tr_n(D^{-1}) -\beta A + D)^{-2}D\right] - \frac{\Id_n}{1 - \beta^2\tr_n(D^{-2})}} \le \frac{C\beta}{\gamma^4\sqrt{n}}.
    \end{equation} 
\end{corollary}
\begin{proof}
    Let $\bar{A}$ be the projection of $A$ onto the set of matrices with operator norm at most $(1-\gamma)/\beta$, as in    \pref{thm:chaining-square-no-sqrtlog}.
    Let $Y(a,D) := E_{\calD_n}[D\,(a -\beta A + D)^{-2}D]$ and let $\bar{Y}(a,D) = E_{\calD_n}[D\,(a -\beta \bar{A} + D)^{-2}D]$.
    
    \[\schnorm{\bar{Y}(a,D) - \E_A\left[\bar{Y}(a,D)\right]} \le \frac{C_1\beta}{\gamma^4\sqrt{n}} \]
    for all $a \in [0,1]$ and all $D \succeq \Id_n$.
    By \pref{prop:gamma-regime-square-expectation-estimate}, 
    \[ \schnorm{\E_A\left[Y(a,D) \one_{\Omega_\gamma} \right] -   F'(a)D(F(a) + D)^{-2}D} \;\le\;\frac{C_0\,\beta^2}{n\gamma^4}
\;+\; \frac{C_0\,\beta^4}{n^{3/2}\gamma^6}  \;<\; o\left(\frac{1}{\sqrt{n}}\right)
 \]
    for some universal constant $C_0$, where $\Omega_{\gamma}=\{\zeta_A\le 2+\eta_n\}$ for $\zeta_A$ and $\eta_n$ as in \pref{cor:AIk+Bk-final-bound}, noting that $\Omega_{\gamma}$ occurs with probability at least $1 - C e^{-c(\log n)^\xi}$ for some $\xi > 1$.

    Since $1 = \one_{\Omega_{\gamma}} + \one_{\Omega_{\gamma}^c}$, we have
    \[\E_A\left[\bar{Y}(a,D)\right] = \E_A\left[\bar{Y}(a,D)\one_{\Omega_{\gamma}}\right] + \E_A\left[\bar{Y}(a,D)\one_{\Omega_{\gamma}^c}\right]\]

    And finally, by \pref{lem:resolvent-bounds} and $\schnorm{\cdot} \le \opnorm{\cdot}$, 
    \[ \schnorm{\E_A\left[\bar{Y}(a,D)\one_{\Omega_{\gamma}^c}\right]} \;\le\; 4\gamma^{-2}\Pr\Omega_{\gamma}^c \;<\; o\left(\frac{1}{\sqrt{n}}\right). \]
    
    The conclusion \pref{eq:diagonal-DQ2D-controlled-in-l2-a} follows by triangle inequality and union bound with the observation that $\Omega_{\gamma}$ entails $\{\bar{A} = A\}$ and hence $\bar{Y}(a,D)\one_{\Omega_{\gamma}} = Y(a,D)\one_{\Omega_{\gamma}}$.
    
     For \pref{eq:diagonal-DQ2D-controlled-in-l2}, the choice $a = \beta^2\tr_n(D^{-1})$ implies $F(a) = 0$ by \pref{lem:subordination-semicircular-0}
     and $F'(a) = 1/(1-\beta^2\tr_n(D^{-2}))$ by \pref{lem:subordination-semicircular-derivative}.
\end{proof}

\section{Quantitative local limit theorem for Bernoulli random variables}

\label{s:llt}
\begin{theorem}
\label{t:llt}
    Let $X_i\sim \mathsf{Bernoulli}(p_i)$, $1\le i\le n$ be independent Bernoulli random variables. Let $\si_n^2 = \sumo in p_i(1-p_i)$. 
    Let $p_n(x)$ be the probability mass function of $S_n = \sumo in X_i$ and 
    $\ga(x) = \rc{\sqrt{2\pi}}e^{-\fc{x^2}{2}}$. Then
\begin{align*}
    \sup_{x\in \Z}
    \ab{\si p_n(x) - \ga\pf{x-\sumo in p_i}{\si}} \le f(\si)
\end{align*}
for some function $f$ such that $\lim_{\si\to \iy} f(\si)=0$. For $\si\ge2$, we can take $f(\si) = \fc{c\log^2(\si)}{\si}$ for some constant $c$.
\end{theorem}
\begin{prf}
    Let $Y_i = X_i - p_i$, $\ph_i$ be the characteristic function of $Y_i$, and $\ph$ be the characteristic function of $T_n = \rc{\si}\sumo in Y_i$. Note 
    $\ph(t) = \prodo in \ph_i\pf{t}{\si}$. 
    Let $q_n(x)$ be the probability mass function of $T_n$. We equivalently bound $\ab{\si q_n(x) - \ga(x)}$ for $x\in -\sumo in p_i + \rc{\si}\Z$.
    As in \cite[\S 3.5]{durrett2019probability}, we have by the inversion formula that 
    \begin{align*}
        \si q_n(x) &= \rc{2\pi} \int_{-\pi\si}^{\pi\si} e^{-itx} \prodo in\ph_i \pf{t}{\si}dt\\
        \ga(x) &= \rc{2\pi} \int_{-\iy}^{\iy} e^{-itx} e^{-\fc{t^2}2}.
    \end{align*}
We bound 
\begin{multline*}
    \ab{\si q_n(x) - \ga(x)}
    \le \rc{2\pi} 
    \ub{\int_{-a}^a
    \ab{\prodo in \ph_i\pf{t}{\si} - e^{-\fc{t^2}2}}dt}{\pat{I}}
    + \rc{2\pi} 
    \ub{\int_{[-\pi\si,\pi\si] \bs [-a,a]}
    \ab{\prodo in \ph_i\pf{t}{\si} 
    } dt}{\pat{II}}
    + \rc{2\pi}\ub{\int_{[-a,a]^c
    } e^{\fc{-t^2}{2}}dt}{\pat{III}}.
\end{multline*}
Note (III)$\to 0$ as $a\to \iy$. Using the tail bound $\rc{2\pi}\pat{III}\le \fc{2}{\sqrt{2\pi}a}e^{-\rc 2a^2}$, it suffices to take $a=\Te\pa{\sqrt{\log\prc\ep}}$ for sufficiently large constants to make this  $<\ep/4$.

    For (I), by \cite[Lemma 3.3.19]{durrett2019probability},
    \begin{align*}
        \ab{\ph_i(t) - 
        \ba{1-p_i(1-p_i) \fc{t^2}2}}
        &\le \E \min\bc{\fc{|tY_i|^3}{6}, |tY_i|^2}.
    \end{align*}
    By \cite[Lemma 3.4.3]{durrett2019probability},
    \begin{align}
    \allowdisplaybreaks
    \nonumber
        &\ab{\prodo in \ph_i\pf{t}{\si} - \prodo in 
        \exp\pa{\fc{p_i(1-p_i)}{\si^2}\fc{t^2}2}}\\
        \nonumber
        &\le \sumo in 
        \ab{\ph_i\pf{t}{\si} - \exp\pa{\fc{p_i(1-p_i)}{\si^2}\fc{t^2}2}}\\
        \nonumber
        &\le 
        \sumo in \ab{\ph_i\pf{t}{\si} - \ba{1-\fc{p_i(1-p_i)}{\si^2}\fc{t^2}2}}
        + \sumo in 
        \ab{\ba{1-\fc{p_i(1-p_i)}{\si^2} \fc{t^2}2}- \exp\pa{\fc{p_i(1-p_i)}{\si^2}\fc{t^2}2}}\\
        &\le 
        \sumo in 
        \min\bc{\fc{|tY_i|^3}{6\si^3}, \fc{|tY_i|^2}{\si^2}} + \sumo in O\pf{(p_i(1-p_i))^2t^4}{\si^4}. 
    \label{e:diff-prod}
    \end{align}
    Let $a>0$ be a constant. 
    Note $\sumo in \fc{(p_i(1-p_i))^2t^4}{\si^4}\le \max_{1\le i\le n} \fc{p_i(1-p_i)}{\si^2}\sumo in \fc{p_i(1-p_i)}{\si^2}t^4 
    \le \fc{a^4}{\si^2} \to 0$ as $\si\to \iy$. 
    Using $|X_i|\le 1$ and $\sumo in \E |X_i|^2 = \sumo in p_i(1-p_i)=\si^2$, for any constant $a$, we have 
    \begin{align*}
        \sumo in \int_{-a}^a 
        \E \min\bc{\fc{|tY_i|^3}{6\si^3}, \fc{|tY_i|^2}{\si^2}}
        &\le 2\int_0^a \fc{t^3}{6\si} = \fc{a^4}{12\si}\to 0
    \end{align*}
    as $\si\to \iy$. Quantitatively, we can make both terms in \eqref{e:diff-prod} $< \fc{\ep}4$ by taking $\si= \Om\pf{a^4}{\ep} = \Om\pf{\log\prc{\ep}^2}{\ep}$.

    For (II), we bound $\ph(t)$ on $[a, \pi\si]$.
    Consider the function
    \begin{align*}
        g(p,t) &= \fc{1-\ab{pe^{it(1-p)}+ (1-p) e^{-itp}}}{p(1-p)t^2}
    \end{align*}
    on $[0,1]\times [0,\pi]$, defined by continuity on the boundary. We claim that this function is continuous and nonzero on the boundary; since it is $>0$ in the interior, it is bounded below by a positive constant $c$. For this, we need to check the limit in 3 situations.
    \begin{enumerate}
        \item $p\to 0$, $t\ne 0$ fixed: We have the Taylor expansion 
        \[
pe^{it(1-p)}+ (1-p) e^{-itp} = pe^{it}(1-O(p)) + (1-p) (1-itp+O(p^2)) = 1+p(e^{it}-1-it) + O(p^2)
        \]
        so as $p\to 0$, $g(p,t) = \fc{1-\ab{1 - p\Re(1+it - e^{it})  + O(p^2)}}{pt^2(1+o(1))} \to \fc{\Re(1+it - e^{it})}{t^2}$, where the convergence is uniform in a neighborhood around $t$, as the constants in the O-notation are uniform in $t\in [0,\pi]$. 
        \item $t\to 0$, $p\not\in \{0,1\}$ fixed: We have the Taylor expansion
        \[
pe^{it(1-p)}+ (1-p) e^{-itp} = 1-p(1-p)\fc{t^2}2+O(t^3),
        \]
        so as $t\to 0$, $g(p,t) = \fc{1-\ab{ 1-p(1-p)\fc{t^2}2+O(t^3)}}{p(1-p)t^2} \to\rc 2$, again locally uniformly in $t$.
        \item 
        $(p,t)\to 0$: We have the Taylor expansion
    \[
pe^{it(1-p)}+ (1-p) e^{-itp} = 1 - \fc{pt^2}{2} + O(pt^3) + O(p^2t^2) + O(p^3).
    \]
    Hence, as $(p,t)\to (0,0)$, $g(p,t) = \fc{1-\ab{1-\fc{pt^2}{2}(1+o(1))}}{pt^2(1+o(1))}\to \rc 2$, showing continuity at $(0,0)$ and by symmetry $(1,0)$.
    \end{enumerate}
Hence we have
\allowdisplaybreaks
\begin{align*}
    \pat{II}
    \le 
    \int_{[-\pi\si,\pi\si] \bs [-a,a]} 
    \prodo in
    \pa{1-c\fc{p_i(1-p_i)}{\si^2}t^2}dt
    &\le \int_{[-\pi\si,\pi\si] \bs [-a,a]}
    \prodo in
    \exp\pa{-c\fc{p_i(1-p_i)}{\si^2}t^2}dt\\
    &=\int_{[-\pi\si,\pi\si] \bs [-a,a]} \exp(-ct^2)dt
\end{align*}
Choosing $a=\Te\pa{\sqrt{\log\prc\ep}}$ also makes  $\rc{2\pi}\pat{II}<\fc{\ep}4$. 
    
Then for large enough $\si = \Om\pf{\log\prc{\ep}^2}{\ep}$, the integral over $[-a,a]$ of each of the terms in \eqref{e:diff-prod} is $<\fc\ep4$. Solving for $\ep$, assuming $\si\ge 2$, this holds for $\ep = O\pf{\log(\si)^2}{\si}$ with appropriate constant. This finishes the proof.
\end{prf}

\section{Mirror descent for solving the TAP equation}

\label{s:mirror-tap}

Let $\psi$ be a strictly convex mirror map defined on a closed convex set $Q$.
The mirror descent algorithm for optimizing a function $F:Q\to \R$,
initialized at $v_0$ with step size $\eta$ is defined as follows. Given $v_k$, let $v_{k+1}$ be such that 
\[
\gd \psi(v_{k+1}) = \gd \psi(v_k) - \eta \gd F(v_k).
\]
Letting $x_{k} = \gd \psi(v_k)$, we can rewrite this in the dual space as 
\[
x_{k+1} = x_k - \eta (\gd F)((\gd\psi)^{-1}(x_k)).
\]
We will apply this to $\psi(\mg) = - H(\mg) := -\sumo in h(\mg_i)$ and $Q=[-1,1]^n$, where $h$ is defined in \eqref{e:h}. We have $\gd \psi(\mg)  = \rc 2 \ln \pf{1+\mg}{1-\mg} = \tanh^{-1}(\mg)$, so we have $\gd\psi:Q^{\circ}\to \R^n$ is bijective. We fix $y$ and take $F(\mg)=\FT(\mg, y)$. This exactly gives the mirror descent iteration in \pref{alg:asl-ta-je-dre}. 

\textbf{Remark}. Note that the same algorithm is used in \cite{el2022sampling}, where it is referred to as natural gradient descent (NGD) and analyzed locally in a ball. We refer to it as mirror descent as the updates are in the dual space \cite{gunasekar2021mirrorless}. 

Convergence guarantees for mirror descent rely on the following definitions.
\begin{definition}
The Bregman divergence of $\psi:Q\to \R$ is $D_\psi:Q\times Q\to \R$ defined as 
\[
D_\psi (x,y) = \psi(x) - \psi(y) - \ip{\gd \psi(y)}{x-y}.
\]
    We say that $F:Q\to \R$ is $c$-strongly convex with respect to $D_\psi$ if for all $x,y\in Q$,
\[
F(x) \ge F(y) + \ip{\gd F(y)}{x-y} + cD_\psi(x,y).
\]
We say that $F:Q\to \R$ is $C$-smooth with respect to $D_\psi$ if for all $x,y\in Q$.
\[
F(x) \le F(y) + \ip{\gd F(y)}{x-y} + CD_\psi(x,y).
\]
\end{definition}
We can check strong convexity with respect to $D_\psi$ by considering the Hessian.
\begin{proposition}[{\cite[Proposition 1.1]{lu2018relatively}}]
\label{p:rel-cvx}
    Suppose that $F$ is twice differentiable. Then
    \begin{enumerate}
        \item $F$ is $c$-strongly convex with respect to $D_\psi$ iff $\gd^2 F\succeq c\gd^2 \psi$.
        \item $F$ is $C$-smooth with respect to $D_\psi$ iff $\gd^2 F\preceq C\gd^2 \psi$.
    \end{enumerate}
\end{proposition}
For $\psi = -H$, we have $\gd^2 \psi(\mg) = D(\mg)$.
For $\FT$, we have under \ref{d:Q-reg} that 
\begin{align*}
    \La^{-1}D(\mg) \preceq 
    \gd_{\mg}^2 \FT(\mg,y) &= \hQ(\mg)^{-1} \preceq \lm^{-1} D(\mg).
\end{align*}
Hence by \pref{p:rel-cvx}, $\FT$ is $\La^{-1}$-strongly convex and $\lm^{-1}$-smooth with respect to $\psi$.
\begin{theorem}[Mirror descent for relatively strongly convex and smooth functions, {\cite[(30)]{lu2018relatively}}]
\label{t:mirror-descent}
    Suppose that $F$ is $c$-strongly convex and $C$-smooth with respect to $D_\psi$. Let $v^*$ be the minimizer of $F$. Let $v_0,v_1,\ldots$ be the iterates of mirror descent with step size $\eta\le \rc C$. Then
    \[
    D_\psi (v_k, v^*) 
    \le (1-\eta c)^k D_\psi(v_0,v^*).
    \]
\end{theorem}
To obtain bounds on Euclidean distance, we note the following.
\begin{lemma}[{\cite[Lemma A.3]{el2022sampling}}]
\label{l:DH-bound}
    We have for all $m,w\in (-1,1)^n$
    \[
\fc{\ve{m-w}^2}2
\le D_{-H}(m,w) \le 
10n\pa{1+\fc{\ve{\tanh^{-1}(w)}_2}{\sqrt n}}\wedge
\ve{\tanh^{-1}(m) - \tanh^{-1}(w)}_2^2.
    \]
\end{lemma}
To use \pref{t:mirror-descent}, we take   $w=\mg^*$ where $\mg^*$ is the minimizer of $\FT$, which means we need to $\mg^*$ to not be too close to the boundary.
\begin{lemma}
\label{l:m-not-too-close-boundary}
    If $\be\le 1$, the solution $\mg^*$ to \eqref{e:TAP-expanded} satisfies
    \[
\ve{\tanh^{-1}(\mg^*)}_{\iy} 
\le \opnorm{\be A} \sqrt n + \ve{y}_\iy + 1.
    \]
\end{lemma}
\begin{prf}
    Let $y\in \R^n$. Considering the $i$th coordinate of \eqref{e:TAP-expanded}, we must have for $\mg=\mg^*$ that 
    \begin{align*}
\tanh^{-1}(\mg_i)
&\le \ab{(-\be A \mg)_i + y_i + \be^2\pa{1-\rc n \ve{\mg}_2^2}\mg_i}
\le \opnorm{\be A} \sqrt n + |y_i| + 1.
    \end{align*}
    The conclusion follows.
\end{prf}

Now we consider the $\td y_t$ that arise in \pref{alg:asl-ta-je-dre}.
An inductive comparison shows that no matter what the estimates $\td \mg_t$ are, for each $i\in [n]$, $(\td y_t)_i - t$ is stochastically dominated by 
$B_t$, and that $(\td y_t)_i + t$ stochastically dominates $B_t$, trajectory-wise.
Fixing $T$, for $\ep<\fc 12$, we hence have with probability $1-\ep$ that $\ve{\ty_t}_{\iy} = O\pa{ \sqrt{\ln \pf{n}{\ep}} }$. Taking $\ep=e^{-n}$, we have that with probability $1-e^{-n}$ that $\ve{\ty_t} = O(\sqrt n)$ for all $0\le t\le T$. 
Let $\mg_t^*$ solve \eqref{e:TAP-expanded}.
By \pref{l:m-not-too-close-boundary}, additionally under the event \eqref{eq:event-opnorm}, $\ve{\tanh^{-1}(\mg_t^*)} = O(\sqrt n)$.
Suppose in \pref{a:asl-ta-je-dre} that $\td \mg_t^{(0)} = \td x_t^{(0)}$; note we are running mirror descent on $\td \mg_t^{(s)} = \tanh (\td x_t^{(s)})$. 
Taking step size $\eta = \lm$, we have by \pref{t:mirror-descent} and \pref{l:DH-bound} that 
\[
\ve{\td \mg_t - \mg_t^*}^2 \le 
2D_{-H}(\td \mg_t, \mg_t^*) \le 
2\pa{1-\fc{\lm}{\La}}^S D_{-H}(\td \mg_t^{(0)}, \mg_t^*)
\le \pa{1-\fc{\lm}{\La}}^S\cdot 20n \pa{1+\fc{\tanh^{-1}(\mg_t^*)}{\sqrt n}}.
\]
Note that this works with an arbitrary initialization $\td \mg_t^{(0)}$; in \pref{a:asl-ta-je-dre}, we initialize it with the value obtained by discretizing ASL-TAP in the dual space.
Thus, taking $S = \Om\pa{\ln\pf{n}{\ep_1}}$ suffices for $\ve{\td \mg_t - \mg_t^*} \le \ep_1$. 
Taking a union bound over all times \pref{alg:asl-ta-je-dre} is called, the probability that this is violated in any step in any call is $e^{-n}\poly(n) e^{O(1/\ep)} = e^{-\Om(n)}$.
For $t\in \{0,\eta,\ldots, T\}$, let $\td y_t^*$ be the iterates if the TAP equation were solved exactly. A coupling argument shows inductively that $\ve{\td y_t - \td y_t^*}\le \int_0^T \ep_1 e^{(T-t)\La }dt = \fc{\ep_1}{\La} (e^{\La T}-1) = O(\ep_1)$, using $D_y\hat\mg(y) = \hQ(y)\preceq \La D(m)^{-1}$. 
We can couple the entire \pref{alg:main} when exact $\td y_T^*$ are used and when the $\td y_T$ are computed from mirror descent. Each step of Glauber or rejection sampling incurs error $\ep_1\poly(n)$, from the differences in the Hamiltonian between the two tilts, and Lipschitzness of $\log \rh_t$ when $\ve{y}_2 = O(\sqrt n)$ (using \eqref{e:gd-y-FTAP}). Now take $\ep_1=\fc{\ep}{\poly(n, e^{O(1/\ep)})}$ so that this contributes error at most $\ep$. It suffices to take the number of steps for each mirror descent to be $S=O\pa{\rc{\ep} + \ln n}$. 

\end{document}